\documentclass{article}

\usepackage{fullpage}

\usepackage{enumerate}
\usepackage{amssymb}
\usepackage{amsfonts}
\usepackage{amsmath}
\usepackage{amsthm}
\usepackage{mathrsfs}
\usepackage{esint} 

\usepackage[colorlinks=true]{hyperref}

\usepackage{tikz}
\usetikzlibrary{calc}
\usetikzlibrary{decorations.markings}

\newtheorem{thm}{Theorem}[section]
\newtheorem{cor}[thm]{Corollary}
\newtheorem{lem}[thm]{Lemma}
\newtheorem{prop}[thm]{Proposition}

\theoremstyle{remark}
\newtheorem{remark}[thm]{Remark}

\newcommand{\dist}{\mathrm{dist}}
\newcommand{\re}{\mathrm{Re}}
\newcommand{\im}{\mathrm{Im}}
\newcommand{\Aut}{\mathrm{Aut}}
\DeclareMathOperator{\arcsinh}{\textnormal{arcsinh}}

\newcommand{\CC}{\mathbb{C}}
\newcommand{\RR}{\mathbb{R}}
\newcommand{\ZZ}{\mathbb{Z}}
\newcommand{\QQ}{\mathbb{Q}}
\newcommand{\FF}{\mathbb{F}}
\newcommand{\PP}{\mathbb{P}}

 \newcommand{\CCC}{{\mathcal {C}}}

\newcommand{\CI}{{\mathcal {I}}}

\newcommand{\CL}{{\mathcal {L}}}

\newcommand{\CO}{{\mathcal {O}}}
\newcommand{\CP}{{\mathcal {P}}}

\newcommand{\CX}{{\mathcal {X}}}

\newcommand{\lra}{\longrightarrow}
\newcommand{\an}{{\mathrm{an}}}

\title{Quantitativity on the Number of Rational Points in the Mordell Conjecture}
\author{Jiawei Yu, Xinyi Yuan, Shengxuan Zhou}

\begin{document}

\maketitle

\begin{abstract}
In this paper, we prove an explicit upper bound on the number of rational points on a smooth projective curve of genus at least two over a number field. This gives explicit constants in the uniform Mordell conjecture proposed by Mazur and 
 proved by Vojta, Dimitrov--Gao--Habegger, and K\"uhne. 
 The main body of this paper consists of two parts: Part I for arithmetic estimates and Part II for analytic estimates.  
\end{abstract}

\tableofcontents

\section{Introduction}

The finiteness of rational points on algebraic curves of genus at least two, conjectured by Mordell \cite{Mor22} and proved by Faltings \cite{Fal}, stands as one of the central achievements of modern number theory. Following this landmark result, Vojta \cite{Vojta_Ar} provided a new proof using techniques from Diophantine approximation, and
Lawrence--Venkatesh \cite{LawrenceVenkatesh} gave a third proof based on $p$-adic Hodge theory.

A natural refinement of Faltings’ theorem was suggested by Mazur \cite[p. 234]{Mazur1986}, who asked whether the number of rational points can be bounded uniformly in terms of the genus of the curve and the Mordell--Weil rank of its Jacobian variety. 
Building on Vojta’s method, Dimitrov--Gao--Habegger \cite{DGH} and K\"uhne \cite{Kuhne} solves Mazur's uniform Mordell conjecture. 
Their theorem has the following form. 

\begin{thm}[Vojta \cite{Vojta_Ar}, Dimitrov--Gao--Habegger \cite{DGH}, K\"uhne \cite{Kuhne}]
\label{thm_DGHKuniformMordellsubgroup}
Let $K$ be a field of characteristic $0$, and let $C$ be a curve of genus $g \geq 2$ over $K$. 
Denote by $J$ the Jacobian variety of $C$. 
Let $\Lambda \subset J(K)$ be a subgroup of finite rank. Fix a point $P_0 \in C(K)$, and view $C(K)-P_0$ as a subset of $J(C)(K)$ via the Abel--Jacobi embedding. Then
\begin{equation*}
\#((C(K)-P_0)\cap \Lambda ) \leq c_1(g)c_2(g)^{\mathrm{rk}(\Lambda)},
\end{equation*}
where $c_1(g)$ and $c_2(g)$ are positive constants depending only on the genus $g$.
\end{thm}

Throughout this paper, a \emph{curve} means a smooth, projective, and geometrically connected curve over a base field. 
In the theorem,  $\Lambda$ is not assumed to be finitely generated, and the rank 
$\mathrm{rk}(\Lambda)=\dim_\QQ (\Lambda\otimes_\ZZ\QQ)$. 
By a specialization argument, it is easy to reduce the theorem to the case where $K$ is a number field.

The theorem is proved by a combination of the following bounds. 
\begin{enumerate}
\item (Large points)
Vojta's proof of the Mordell conjecture in \cite{Vojta_Ar} actually gives a deep inequality concerning distribution of rational points of large heights, which particularly implies an upper bound on the number of rational points of large heights. Note that this is sufficient for the Mordell conjecture by the Northcott property of heights.
\item (Small points)
The works \cite{DGH, Kuhne} prove a uniform Bogomolov conjecture, which  
 gives an upper bound on the number of algebraic points of small heights. 
By a simple argument of sphere packing, the uniform Bogomolov conjecture also bounds the number of all rational points complementary to (1). 
\end{enumerate}

Previously, the original Bogomolov conjecture for a single curve was proved by Ullmo \cite{Ull} by the equidistribution theorem of Szpiro--Ullmo--Zhang \cite{SUZ}. 
A less well-known proof of the Bogomolov conjecture, which establishes the positivity of the admissible volumes of curves over number fields, is obtained by 
the works of Zhang \cite{Zhang_Admissible, Zhang_phi}, Cinkir \cite{cinkir}, and de Jong \cite{dejong_NeronTate}. 
Before the works \cite{DGH, Kuhne}, DeMarco--Krieger--Ye \cite{DKY} proved the uniform Bogomolov conjecture for bi-elliptic curves of genus 2.

Based on the theory of adelic line bundles of Yuan--Zhang \cite{YZ}, recently Yuan \cite{Yuan_Bogomolov} gave a different proof of the uniform Bogomolov conjecture and thus the uniform Mordell conjecture. The key ingredient of Yuan's approach is a family version of the positivity of the admissible volume by Zhang \cite{Zhang_Admissible, Zhang_phi}, Cinkir \cite{cinkir}, and de Jong \cite{dejong_NeronTate}, while the key ingredient of K\"uhne's approach is a family version of the equidistribution argument of Ullmo \cite{Ull}. 

\subsection{Main theorem: quantitative Mordell conjecture}

It is natural to ask whether one can obtain a quantitative version of the uniform Mordell conjecture, i.e. to determine explicit constants $c_1(g)$ and $c_2(g)$ 
in Theorem \ref{thm_DGHKuniformMordellsubgroup}. 
There are serious difficulties to extract explicit constants in the above approaches. 
In fact, the approach of \cite{DGH, Kuhne} is based on an equidistribution argument, the approach of Yuan \cite{Yuan_Bogomolov}  is based on a compactness argument at archimedean places, and both approaches are based on arithmetic of chains of subvarieties of the moduli space of curves. 

However, in the function field case, Looper--Silverman--Wilms \cite{LSW} established the uniform Bogomolov conjecture with surprisingly explicit constants. 
Their proof is through estimates on single curves instead of moduli spaces, and still based on the positivity of admissible volumes of 
Zhang \cite{Zhang_Admissible, Zhang_phi}, Cinkir \cite{cinkir}, and de Jong \cite{dejong_NeronTate}. 
Building on the result of \cite{LSW} and adapting Vojta’s proof of the Mordell conjecture over function fields, Yu \cite{Yu} achieved $c_1(g)=16g^2 + 32g + 129$ and $c_2(g)=40g$ in Theorem \ref{thm_DGHKuniformMordellsubgroup} for non-isotrivial curves over function fields of characteristic $0$.

The work \cite{LSW} does not cover number fields due to many obstacles at archimedean places. In this paper, we manage to overcome all these archimedean obstacles and prove a quantitative Bogomolov conjecture over number fields. We also establish a quantitative Vojta inequality (with explicit constants) over number fields. 
Consequently, we prove the following our main theorem, which is an explicit version of Theorem \ref{thm_DGHKuniformMordellsubgroup}.  

\begin{thm}[Theorem \ref{mordell_lang}]
\label{thm_quantitativeMordellsubgroup}
Let $K$ be a field of characteristic $0$, $C$ a curve of genus $g\geq 2$ over $K$, and $J$ the Jacobian variety of $C$. For any line bundle $\alpha$  of degree $1$ on $C$ and any subgroup $\Lambda\subset J(K)$ of finite rank,
$$\#((C(K)-\alpha)\cap\Lambda)\le 10^{13}g^8\cdot\min\left\{1+\frac{5}{4\sqrt g},1+\frac{3\log g}{g}\right\}^{\mathrm{rk}(\Lambda)}.$$
\end{thm}

If $K$ is a number field, then $J(K)$ is finitely generated by the Mordell--Weil theorem and we can take $\Lambda=J(K)$. The gives the following quantitative version of the original Mordell conjecture.

\begin{thm}
\label{thm_quantitativeMordellfullJacobian}
Let $C$ be a curve of genus $g\geq 2$ over a number field $K$, and $J$ the Jacobian variety of $C$. Then
$$\#C(K)\le 10^{13}g^8\cdot\min\left\{1+\frac{5}{4\sqrt g},1+\frac{3\log g}{g}\right\}^{\mathrm{rk}\, J(K)}.$$
\end{thm}

In the theorems, the Vojta constant 
$$\displaystyle c_2(g)=\min\left\{1+\frac{5}{4\sqrt g},1+\frac{3\log g}{g}\right\},$$
 which essentially comes from Vojta's proof of the Mordell conjecture, tends to 1 as $g\to\infty$. This provides an affirmative answer to a question of Gao and Habegger (cf. \cite[Que. 1.19]{Habegger_ICM}).  
 
On the other hand, recall that the Manin--Mumford conjecture proved by Raynaud 
\cite{Ray83a,Ray83b} asserts that $(C(K)-\alpha)\cap J(K)_{\rm tor}$ is finite.
Here $J(K)_{\rm tor}$ denotes the subgroup of torsion elements of $J(K)$. 
In Theorem \ref{thm_quantitativeMordellsubgroup}, 
the constant $c_1(g)=10^{13} g^8$ bounds the order of $(C(K)-\alpha)\cap J(K)_{\rm tor}$ by taking $\Lambda=J(K)_{\rm tor}$, so we call $c_1(g)$ the Manin--Mumford constant. 
Our proof of the theorem actually gives the following quantitative Manin--Mumford conjecture with better constants. 

\begin{thm}[Theorem \ref{manin_mumford}]
\label{manin_mumford intro}
Let $K$ be a field of characteristic $0$, $C$ a curve of genus $g\geq 2$ over $K$, and $J$ the Jacobian variety of $C$. Then for any line bundle $\alpha$ of degree $1$ on $C$,
$$\#\big((C(K)-\alpha)\cap J(K)_\mathrm{tor}\big)\le 3.2\cdot 10^{11}g^{\frac{17}3}.$$
\end{thm}

In Theorem \ref{thm_quantitativeMordellsubgroup}, 
if $C$ is not defined over a number field, we are in the function field case and have the following theorem with better constants.

\begin{thm}[Theorem \ref{geometric_mordell_lang}]
\label{geometric_mordell_lang intro}
Let $K$ be a field of characteristic $0$, $C$ a curve of genus $g\geq 2$ over $K$, and $J$ the Jacobian variety of $C$. Assume that $C$ is not isotrivial over $\QQ$ in that  $C_{\overline K}$ cannot be descended to $\overline\QQ$. Then for any line bundle $\alpha$ of degree $1$ on $C$  and any subgroup $\Lambda\subset J(K)$ of finite rank,
$$\#((C-\alpha)\cap\Lambda)\le1.8\cdot10^6g^{3}\cdot\left(1+\frac5{4\sqrt g}\right)^{\mathrm{rk}(\Lambda)},$$
$$\#((C-\alpha)\cap\Lambda)\le2.5\cdot10^4g^{8}\cdot\left(1+\frac{3\log g}{g}\right)^{\mathrm{rk}(\Lambda)}.$$
\end{thm}

Our proof of this function field version applies a variant of the main result of 
Looper--Silverman--Wilms \cite{LSW}, but the other parts of the proof are similar to that of Theorem \ref{thm_quantitativeMordellsubgroup}. 
Note that our Vojta constant here is better than that obtained by Yu \cite{Yu}, which comes from a refined version of the Vojta inequality and a more delicate method of point counting.

\subsection{Some  consequences}

Let us consider some consequences of  Theorem \ref{thm_quantitativeMordellfullJacobian}, which give different bounds of the number of rational points.  

\subsubsection*{First consequence: beyond Chabauty range}

For the first consequence, recall that in the Chabauty range 
$\mathrm{rk}\, J(K)<g$, 
Coleman \cite{Coleman} applied the Chabauty method give an upper bound of 
$\#C(K)$ which is almost linear in $g$. Our theorem gives a clean bound beyond the Chabauty range. 
In fact, by the inequality
$$\left(1+\frac{3\log g}{g}\right)^g\le g^3,$$
Theorem \ref{thm_quantitativeMordellfullJacobian} gives 
$$\# C(K)\le 10^{13}g^{3\frac{\mathrm{rk}\, J(K)}{g}+8}.$$
This is polynomial in $g$ if $\displaystyle\frac{1}{g}\mathrm{rk}\, J(K)$ is bounded. 

\subsubsection*{Second consequence: bounds by bad reductions}

Our second consequence is the following ``more explicit" bound on the number of rational points.

\begin{thm} \label{bound by bad reduction}
Let $C$ be a curve of genus $g\geq 2$ over a number field $K$. 
Denote by $d=[K:\mathbb{Q}]$ the degree of $K$ over $\QQ$, by $\Delta_{K/\mathbb{Q}}$ the discriminant of $K$ over $\mathbb{Q}$, and by $N_{C/K}^0$ the product of the norms of the prime ideals of $O_K$ where $C$ has bad reduction. 
Then
$$\#C(K)\le c_1(g,d)\cdot (N_{C/K}^0)^{c_2(g,d)}\cdot|\Delta_{K/\mathbb{Q}}|^{c_3(g,d)},$$
where
\begin{align*}
c_1(g,d)=&\ 10^{13}g^8\cdot\min\left\{1+\frac{5}{4\sqrt g},1+\frac{3\log g}{g}\right\}^{2g^3d^32^{8g^2}} \
<\  10^{13}g^{6g^2d^32^{8g^2}+8},  \\
c_2(g,d)=&\ 4g^3d^22^{8g^2}\cdot\log_4\min\left\{1+\frac{5}{4\sqrt g},1+\frac{3\log g}{g}\right\}
\ <\ 8.7 \cdot g^2d^2 \cdot 2^{8g^2} \log g,  \\
c_3(g,d)=&\ gd2^{8g^2}\cdot\log_4\min\left\{1+\frac{5}{4\sqrt g},1+\frac{3\log g}{g}\right\}
\ <\ 2.2 \cdot d \cdot 2^{8g^2} \log g.
\end{align*}
\end{thm}

The theorem is a consequence of Theorem \ref{thm_quantitativeMordellfullJacobian}
and the following bound of R\'emond \cite[Proposition 5.1]{Rem10} on the Mordell--Weil rank:
$$\mathrm{rk}\, J(K)\le\frac{gd2^{8g^2}}{\log 4}(4dg^2\log N_{J/K}^0+\log|\Delta_{K/\mathbb{Q}}|+g^2d^2\log 16)-1.$$
Here $N_{J/K}^0$ is the product of the norms of the prime ideals of $O_K$ where $J$ has bad reduction, which is a factor of $N_{C/K}^0$. 

Recall that the admissible volume $\bar \omega_{C/K,a}^2$ is the self-intersection number of Zhang's admissible canonical bundle; see \S\ref{section_ad} for a review. 
In Theorem \ref{bound by bad reduction}, we can replace 
$N_{C/K}^0$ by $e^{88 \bar \omega_{C/K,a}^2}$, since we have
$$
\bar \omega_{C/K,a}^2 \geq \frac{\max\{g-1,2\}}{2g+1}\varphi(C)
\geq \frac{1}{88} \log N_{C/K}^0. 
$$
Here the first inequality is from Theorem \ref{phi_inequality}, a consequence of works of Zhang \cite{Zhang_phi}, de Jong \cite{dejong_NeronTate}, and Wilms \cite{Wilms}, 
 and the second inequality is by Cinkir \cite[Theorem 2.11]{cinkir} and Cinkir \cite{cinkir15}. 

We can further change the expression in terms of the admissible volume $\bar \omega_{C/K,a}^2$ by that in terms of the Faltings height $h_{\rm Fal}(J)$, since we have 
$$
\frac{1}{[K:\QQ]}
\bar \omega_{C/K,a}^2 
\leq 12 h_{\rm Fal}(J)+ 6g \log(2\pi^2).
$$
See \cite[Theorem 4.14]{Yuan_Bogomolov} for the inequality, a consequence of the arithmetic Noether formula of Faltings \cite{Faltings}, 
an inequality of Zhang \cite{Zhang_Admissible}, and a lower bound of the delta invariant by Wilms \cite{Wil17}.

Interested readers can also compare our results with the results of
R\'emond \cite[Theorem 1.1, Theorem 1.2]{Rem10} for rational points on planar curves. Our current results strengthen the ones in the loc. cit. in some sense.

Finally, for integral points, Corvaja--Zannier \cite{CZ03} obtained explicit upper bounds on the number of integral points on a hyperbolic curve with at least three points at infinity, which is based on a sharp explicit version of Schmidt’s subspace Theorem and does not involve the Mordell--Weil rank.

\subsubsection*{Third consequence: hyperelliptic curves}

If $C$ is a hyperelliptic curve, we have a better bound for $\mathrm{rk}\, J(K)$, and thus a better bound for $\#C(K)$ as follows.
\begin{thm}[Theorem \ref{hyperelliptic}]
\label{hyperelliptic_intro}
Let $K$ be a number field of degree $d$ over $\mathbb{Q}$. Let $f(x)\in O_K[x]$ be a monic and square-free polynomial.  Let $\Delta_f\in O_K$ be the discriminant of $f$.
\begin{enumerate}
\item[(1)]If $\deg f=2g+1$ for some integer $g\ge 2$, then
$$\#\{(x,y)\in K^2:y^2=f(x)\}\le
10^{13}\cdot 2^dg^{9d\log_2(2gd+d)+18d+8}(2gd+d)^{\frac32d}\cdot|N_{K/\mathbb{Q}}\Delta_{f}|^{4\log_2g+\frac12}\cdot|\Delta_{K/\mathbb{Q}}|^{3\log_2g+\frac12}.$$
\item[(2)]If $\deg f=2g+2$ for some integer $g\ge 2$, then
$$\#\{(x,y)\in K^2:y^2=f(x)\}\le
10^{13}\cdot2^{6d+1}g^{9d\log_2(2gd+2d)+18d+8}(2gd+2d)^{\frac92d}
\cdot|N_{K/\mathbb{Q}}\Delta_{f}|^{4\log_2g+\frac52}
\cdot|\Delta_{K/\mathbb{Q}}|^{3\log_2g+\frac32}.$$
\end{enumerate}
\end{thm}
Here the discriminant $\Delta_f\in O_K$ is the usual one defined by 
$$
\Delta_f=\prod_{1\leq i<j \leq \deg f}(x_i-x_j)^2,$$
where $x_1,\dots, x_{\deg f}$ are the roots of $f$ in $\overline K$.

\subsubsection*{Fourth consequence: average for hyperelliptic curves}

For the last consequence, 
recall that Bhargava--Gross \cite{Bhargava-Gross} proved that the average size of the $2$-Selmer groups of hyperelliptic curves with marked rational Weierstrass points over $\QQ$ is $3$.
If the rational point is non-Weierstrass,  Shankar--Wang \cite{SW} proved that the average size is at most $6$. 
By the easy inequality
$$
\left(1+\frac{5}{4\sqrt g}\right)^{\mathrm{rk}\, J(K)}
< 2^{\mathrm{rk}\, J(K)}
\leq \#\mathrm{Sel}_2(J),
$$
we have the following consequence of Theorem \ref{thm_quantitativeMordellfullJacobian}. 

\begin{thm}[Theorem \ref{average}]
\label{average intro}
\begin{enumerate}
    \item For hyperelliptic curves of genus $g\geq 2$ over $\mathbb{Q}$ with a marked rational Weierstrass point, the average number of rational points is at most $3\cdot10^{13}g^8.$
    \item For hyperelliptic curves of genus $g\geq 2$ over $\mathbb{Q}$ with a marked rational non-Weierstrass point, the average of number of rational points is at most $6\cdot10^{13}g^8.$
\end{enumerate}
\end{thm}

\subsection{Quantitative Vojta inequality}

Our proof of Theorem \ref{thm_quantitativeMordellsubgroup} is based on a quantitative Vojta inequality and a quantitative Bogomolov conjecture with some delicate arguments in sphere packing. 
We will first introduce these two quantitative theorems and then give some ingredients of our sphere packing. 

We start with some notation for heights and arithmetic numbers, which are introduced in \S \ref{sec prelim} with details. 
Let $C$ be a curve of genus $g\geq 2$ over a number field $K$, and let $J$ be the Jacobian variety of $C$ over $K$. 
Taking the base change to a finite extension of $K$ if necessary, we can assume that there is a line bundle $\alpha_0$ on $C$ such that $(2g-2)\alpha_0$ is isomorphic to the canonical sheaf $\omega=\omega_{C/K}$.
Here we usually write tensor products of line bundle additively. 
Via $\alpha_0$, we have an Abel--Jacobi map 
$$i_{\alpha_0}:C\lra J,\quad x\longmapsto (x)-\alpha_0. $$
This gives a theta divisor $\theta$ on $J$, defined as the image of $C^{g-1}\to J$. 
The theta divisor gives a N\'eron--Tate height function $\hat h: J(\overline K)\to \RR$, which extends to a positive definite quadratic form on 
$J(\overline K)_\RR=J(\overline K)\otimes_\ZZ\RR$. 
This gives a metric on $J(\overline K)_\RR$ by 
$$|x|=\sqrt{[K:\QQ]\hat h(x)}, \quad x\in J(\overline K)_\RR.$$
Denote by $\angle(x,y)$ the angle between two vectors $x,y\in J(\overline K)_\RR$ under this metric. 
The notations $|x|$ and $\angle(x,y)$ for $x,y\in C(\overline K)$ are understood via the map $i_{\alpha_0}:C\to J$. 

Let $\bar\omega_a=\bar\omega_{C/K,a}$ be Zhang's admissible canonical line bundle on $C$. The Bogomolov conjecture proved by Ullmo is equivalent to the strictly positive of the arithmetic self-intersection number $\bar\omega_a^2$. 
The number $\bar\omega_a^2$ is a canonical arithmetic invariant of $C/K$, and it is called the \emph{admissible volume} of $C/K$. 
We refer to \S\ref{section_ad} for more details on this. 

In this paper, we prove the following quantitative version of Vojta's inequality.

\begin{thm}[Theorem \ref{vojta}]
\label{thm_vojta_introduction}
Let $C$ be a curve of genus $g \geq 2$ over a number field $K$. Let $P_1,P_2\in C(K)$ be distinct rational points. Assume
$$|P_1|\ge1.2\cdot10^9g^{\frac{7}{3}}\sqrt{\bar{\omega}_a^2}$$
and
$${|P_2|}{}\ge10^5g^{\frac52}|P_1|.$$
Then
$$\cos\angle(P_1,P_2)\le\sqrt{\frac{1.01}{g}}.$$
\end{thm}

The theorem with implicit constants is proved by Vojta \cite{Vojta_Ar} in order to prove the Mordell conjecture. 
The most technical part of the proof of Vojta \cite{Vojta_Ar} is an application of the arithmetic Riemann--Roch theorem of Gillet--Soul\'e \cite{GS92} to construct a small section. 
This technical part is replaced by an application of the more elementary Siegel's lemma by Bombieri \cite{Bom90}, and replace by an application of Yuan's arithmetic Siu inequality by Yuan \cite{Yuan_Vojta}. 
We refer to \cite[\S1]{Yuan_Vojta} for a brief history of these approaches. 

Note that R\'emond \cite[Thm. 1.1]{Rem00b} proved a Vojta inequality for subvarieties of abelian varieties, most of whose constants are also explicitly written down. 

Our proof of Theorem \ref{thm_vojta_introduction} is based on the framework of Yuan \cite{Yuan_Vojta} on Vojta's method, and our additional idea is to use Zhang's admissible adelic line bundles to have explicit constants. 
Moreover, we also need some explicit inequalities to compare the Arakelov metric and the hyperbolic metric of a compact Riemann surface. 

Note that the bound $\displaystyle\sqrt{\frac{1.01}{g}}$ in the theorem is very sharp by the current method. In fact, all the bounds we can derive from Vojta's method should be strictly bigger than $\displaystyle\sqrt{\frac{1}{g}}$. This can also be seen from the statement of \cite[Theorem 2.1]{Yuan_Vojta}. 

As a companion of the quantitative Vojta inequality, 
we have the following quantitative version of Mumford's inequality in \cite{Mum65} (cf. \cite[\S5.7]{Serre}).

\begin{thm}[Theorem \ref{mumford}]
\label{mumford_introduction}
Let $C$ be a curve of genus $g\geq 2$ over a number field $K$. Let $P_1,P_2\in C(\overline K)$ be distinct algebraic points. Assume
$$|P_1|\ge 10^{9}g^{\frac{7}{3}}\sqrt{\bar\omega_a^2}$$
and
$$|P_1|\le{|P_2|}{}\le1.15|P_1|.$$
Then
$$\cos\angle(P_1,P_2)\le{\frac{1.01}{g}}.$$
\end{thm}

As in Mumford's original inequality, the proof of our theorem uses some techniques of those of Theorem \ref{thm_vojta_introduction}, but it is much easier than that of Theorem \ref{thm_vojta_introduction}.

\subsection{Quantitative Bogomolov conjecture}

The following theorem is our quantitative version of the Bogomolov conjecture.

\begin{thm}[Theorem \ref{bogomolov}]
\label{bogomolov_introduction}
Let $K$ be a number field, $C$ a curve of genus $g\geq 2$ over $K$, and $J$ the Jacobian variety of $C$. 
\begin{enumerate}
    \item For any $0\le r<\sqrt\frac1{8g}$ and $x\in J(\overline K)_\mathbb R$, 
    $$\#\{P\in C(\overline K):|P-x|\le r\sqrt{\bar\omega_a^2}\}<\frac{3.2\cdot 10^{11}g^{\frac{17}3}}{1-8gr^2}\left(1+10^{-6}\log\left(\frac{1}{1-8gr^2}\right)\right).$$
    \item Let $\kappa>1$ and $\theta\in(0,\pi/2)$ be real numbers satisfying
    $$g\cos^2\theta>\frac{(\kappa+1)^2}{4\kappa}+\frac{1}{3.2\cdot 10^{11} g^{\frac{11}{3}}}.$$
    Then for any nonzero $x\in J(\overline K)_\mathbb R$,
    $$\#\{P\in C(\overline K):|x|\le|P|\le\kappa|x|,\, \angle(P,x)\le\theta\}
    \leq 
    3.2\cdot 10^{11}g^{\frac{17}3}.$$
\end{enumerate}
\end{thm}

The first inequality is a number field version of \cite[Theorem 1.3]{LSW}, and it has a much larger coefficient $3.2\cdot 10^{11}$ due to the larger coefficients from archimedean places. 
The second inequality is an effective variant of the inequality coming from the extra term in   
\cite[Theorem 1.1]{Yuan_Bogomolov}. 
In fact, an elementary calculation gives
$$
    \{P\in C(\overline K):|x|\le|P|\le\kappa|x|,\, \angle(P,x)\le\theta\}\  
       \subseteq\  
  \{P\in C(\overline K):|P-x|\le (\kappa^2+1-2\kappa\cos\theta)^{\frac12}|x| \}.$$
Then an upper bound of the order of the right-hand side gives the effect of the extra term in \cite[Theorem 1.1]{Yuan_Bogomolov}. 
See also Remark \ref{rmk extra term}. 

Let us describe our idea to prove Theorem \ref{bogomolov_introduction}. 
As mentioned above, the global part of our proof follows closely the proof \cite[Theorem 5.2]{LSW}, where the key ingredient is a quadratic identity, and the local part of our proof consists of a few technical estimates of invariants of compact Riemann surfaces. 
We will first describe the global part, and come back to the local part later.

Assume that $P_1,\dots, P_n$ are $n$ distinct points of $C({\overline{K}})$ with small distances $|P_i -x|$. We hope to obtain an upper bound for $n$. By base change, we can assume that $P_1,\dots, P_n$ are defined over $K$. Extend $P_i$ to the admissible adelic divisor $\overline{\mathcal{O}}(P_i)_a$ on $C$ (cf. \S\ref{section_ad}). 
By the arithmetic Hodge index theorem (cf. Theorem \ref{hodge index}) and the arithmetic adjunction formula (cf. Theorem \ref{adjunction}), we have
a quadratic identity
$$\frac{4(n+g-1)(g-1)}{n-1}\sum_{1\le i<j\le n}|P_i-P_j|^2=\frac{4ng(g-1)}{n-1}\sum_{1\le i<j\le n}h_{\overline{\mathcal{O}}(\Delta)_a}(P_i,P_j)+{\frac{n^2}2\bar\omega_a^2}+\left|(2g-2)\sum_{i=1}^nP_i-{n\omega}\right|^2.$$
Here the height
$$
h_{\overline{\mathcal{O}}(\Delta)_a}(P_i,P_j)
 =\sum_{v\in M_K}  \epsilon_{v} G_{a,v}(P_i, P_j),
$$
where $G_{a,v}:(C^2\setminus \Delta)(\overline K_v)\to \RR$ is Zhang's admissible Green function at $v$ as in \S\ref{section_ad}.
See Lemma \ref{quadratic} for a proof of this quadratic identity. 
This quadratic identity is extracted from the proof of \cite{LSW}. 

To prove the first inequality of Theorem \ref{bogomolov_introduction}, assume that $|P_i -x|^2<t$ for a small positive number $t$, and we want to have an upper bound of $n$ in terms of $t$. 
We first have
$$
|P_i-P_j|^2 \leq (|P_i-x|+|P_j-x|)^2 
\leq (\sqrt t+\sqrt t)^2=4t. 
$$
We will prove a key bound of the form
$$
\sum_{1\leq i<j\leq n} h_{\overline{\mathcal{O}}(\Delta)_a}(P_i,P_j) \geq -O((n\log n)\bar{\omega}_a^2).
$$
It follows that the inequality becomes 
$$
8(g-1)n^2 t
\geq  \left|(2g-2) \sum_{i=1}^n P_i -n\omega\right|^2+\frac{n^2}{2} \bar{\omega}_a^2-O((n\log n)\bar{\omega}_a^2)
\geq  \frac{n^2}{2} \bar{\omega}_a^2-O((n\log n)\bar{\omega}_a^2).
$$
This gives an upper bound of $n$ 
as long as $\displaystyle 2t<\frac{1}{8(g-1)} \bar{\omega}_a^2$. 
By some refined method, we eventually get the first inequality of Theorem \ref{bogomolov_introduction}.

For the second inequality of Theorem \ref{bogomolov_introduction}, the conditions
$|x|\le|P|\le\kappa|x|$ and $\angle(P,x)\le\theta$ give
$$
|P_i-P_j|^2=|P_i|^2+|P_j|^2-{2}|P_i||P_j|\cos\angle(P_i,P_j)\\
\le 2\kappa |x|^2-2|x|^2\cos(2\theta)  .$$
The projection of $P_i$ to the direction of $x$ has length $|P_i|\cos\angle(P_i,x)$. It follows that 
\begin{align*}
\left|(2g-2)\sum_{i=1}^nP_i-{n\omega}\right|^2
\ge(2g-2)^2\left(\sum_i|P_i|\cos\angle(P_i,x)\right)^2
\ge (2g-2)^2n^2|x|^2\cos^2\theta .
\end{align*}
Then the quadratic identity gives 
$$
 2 |x|^2 (\kappa-\cos(2\theta))  n^2
 \geq 
 (2g-2)^2n^2|x|^2\cos^2\theta
+  \frac{n^2}{2} \bar{\omega}_a^2-O((n\log n)\bar{\omega}_a^2).$$
This gives an upper bound of $n$
as long as $2 (\kappa-\cos(2\theta)) < (2g-2)^2 \cos^2\theta$. 
A refinement of this argument gives the second inequality of Theorem \ref{bogomolov_introduction}.

It remains to prove an explicit bound of the form
$$
\sum_{1\leq i<j\leq n} h_{\overline{\mathcal{O}}(\Delta)_a}(P_i,P_j) \geq -O((n\log n)\bar{\omega}_a^2).
$$
We need Zhang's global $\varphi$-invariant
$$
\varphi(C) = \sum_{v\in M_K} \epsilon_{v} \varphi_v(C),
$$
where the local $\varphi$-invariant
$$\varphi_v(C)=-\int_{(C^2)_{K_v}^{\mathrm{an}}}G_{a,v}\, c_1(\overline{\mathcal{O}}(\Delta)_{a,v})^{\wedge 2}\geq 0.$$
This is recalled in \S\ref{section_ad}. 
By the works of
Zhang \cite{Zhang_Admissible, Zhang_phi}, Cinkir \cite{cinkir}, de Jong \cite{dejong_NeronTate}, and Wilms \cite{Wilms} on the original Bogomolov conjecture, we have an inequality 
$$
\bar\omega_{a}^2 \geq \frac{g-1}{2g+1} \varphi(C). 
$$
We refer to Theorem \ref{phi_inequality} for more details. 

It is reduced to prove an explicit bound of the form
$$
\sum_{1\leq i<j\leq n} h_{\overline{\mathcal{O}}(\Delta)_a}(P_i,P_j) \geq -O((n\log n)  \varphi(C) ).
$$
Note that $h_{\overline{\mathcal{O}}(\Delta)_a}(P_i,P_j)$
and $\varphi(C)$ are both sums of local components over places $v$ of $K$. 
The problem is reduced to find a suitable positive constant $c(g)$ independent of $C$ and $v$ such that 
$$
\mathrm{FE}_{v}(n,C)+ c(g)\,(n\log n)\, \varphi_v(C) \geq 0.
$$
Here the \emph{Faltings--Elkies} invariant
$$
\mathrm{FE}_{v}(n,C)=\inf_{P_1,\dots, P_n\in C(K_v)} \sum_{1\leq i<j\leq n} G_{a,v}(P_i, P_j).
$$
Note that the problem is purely local for $v$. 

If $v$ is a non-archimedean place, then both $\mathrm{FE}_{v}(n,C)$ and $\varphi_v(C)$ can be computed in terms of the reduction dual graph of $C$. 
In this case, a lower bound of $\varphi_v(C)$ is given by Cinkir \cite{cinkir}, and a suitable lower bound of $\mathrm{FE}_{v}(n,C)$ is given by \cite{LSW}. 
Their combination gives $c(g)$ at $v$. 

If $v$ is archimedean, the term $\mathrm{FE}_{v}(n,C)$ is an invariant of the compact Riemann surface $C_v(\CC)$. Lower bounds of $\mathrm{FE}_{v}(n,C)$ was previously obtained by Faltings \cite{Fal} and Elkies (cf. \cite[\S VI.5]{lang1}), but their bounds are not uniform as the Riemann surface varies. 
Our explicit bound in this case is as follows. 

\begin{thm}[Theorem \ref{thmfaltingselkieszhangphiinv part I},  Theorem \ref{thmfaltingselkieszhangphiinv}]
\label{thmfaltingselkieszhangphiinv_introduction}
Let $C$ be a curve of genus $g\geq 2$ over a number field $K$. 
Let $v$ be an archimedean place of $K$. 
Then we have
$$ \mathrm{FE}_v(C,n) \geq - \left(4\cdot 10^3 g^{\frac{1}{3}} n\log n+ 1.32\cdot 10^{10} g^{\frac{11}{3}} n \right) \cdot \varphi_v(C).$$
\end{thm}
The proof of this explicit bound is the most technical part of our analytic part of this paper. We will come back to this later.

\subsection{Sphere packing}

Now we describe our argument of sphere packing to prove 
Theorem \ref{thm_quantitativeMordellsubgroup} by  
Theorem \ref{thm_vojta_introduction}, Theorem \ref{mumford_introduction}, and Theorem \ref{bogomolov_introduction}. 
It is easy reduce to reduce the problem the case $\alpha=\alpha_0$. 
We divide the set $C\cap \Lambda=(C(K)-\alpha_0)\cap\Lambda$ into the following three parts:
\begin{align*}
(C\cap\Lambda)_{\mathrm{small}}=&\ \left\{P\in C\cap\Lambda:|P|\le\sqrt\frac{\bar\omega_a^2}{16g}\right\},  \\
(C\cap\Lambda)_{\mathrm{medium}}=&\ \left\{P\in C\cap\Lambda:\sqrt\frac{\bar\omega_a^2}{16g}<|P|\le1.2\cdot10^9g^{\frac{7}{3}}\sqrt{\bar\omega_a^2}\right\},  \\
(C\cap\Lambda)_{\mathrm{large}}=&\ \left\{P\in C\cap\Lambda:|P|>1.2\cdot10^9g^{\frac{7}{3}}\sqrt{\bar\omega_a^2}\right\}.
\end{align*}
This is similar to the partition by Dimitrov--Gao--Habegger \cite{DGH}, but we will apply many refined methods to obtain good constants in the bounds. 

Our cut-off limits for $|x|$ in the three sets are chosen from the constants in Theorem \ref{thm_vojta_introduction} and Theorem \ref{bogomolov_introduction}. 
For example,  the set $(C\cap\Lambda)_{\mathrm{small}}$ of small points comes from Theorem \ref{bogomolov_introduction}(1) with $r=\sqrt{\frac{1}{16g}}$.
This gives the bound
$$\#(C\cap\Lambda)_{\mathrm{small}} \leq  6.5\cdot 10^{11}g^\frac{17}{3}.$$
It remains to bound the medium points and the large points. 

Denote $n=\mathrm{rk}(\Lambda)$ and assume $n\geq 1$.
The space $V=\Lambda_\RR$ under the norm $|\cdot|$ is isometric to the Euclidean space $\RR^n$.
For any nonzero vector $u\in V$ and any angle $\theta\in (0,\pi)$, denote the cone
$$
\mathrm{cone}(u,\theta)=\{x\in V: \angle(u,x)\leq \theta\}
$$
centered at $u$ and of angular radius $\theta$ in $V$, and denote the spherical cap 
$$\mathrm{cap}(u,\theta)=\{x\in V:|x|=1,\ \angle(u,x)\le\theta\}$$ 
centered at $u$ and of angular radius $\theta$ in $V$. 
Note that $\mathrm{cap}(u,\theta)$ lies in the sphere
$$S_V=\{x\in V:|x|=1\}.$$

\subsubsection*{Medium points}

Now we describe our idea to bound the set 
$(C\cap\Lambda)_{\mathrm{medium}}$.
For any $0<r<r'$, denote 
$$(C\cap\Lambda)_{r,r'}=\{P\in C\cap\Lambda:r<|P|\le r'\}.$$
Then we have 
$$
(C\cap\Lambda)_{\mathrm{medium}}=(C\cap\Lambda)_{r_1, r_2}
$$
with 
$$
r_1=\sqrt\frac{\bar\omega_a^2}{16g},\qquad r_2=1.2\cdot10^9g^{\frac{7}{3}}\sqrt{\bar\omega_a^2}.
$$
Cover $(C\cap\Lambda)_{\mathrm{medium}}$ by the sets 
$$(C\cap\Lambda)_{2^{i-1}r_1, 2^i r_1},\quad i=1,\dots, \lceil\log_2(r_2/r_1)\rceil.$$
It suffices to bound 
 $\#(C\cap\Lambda)_{r,2 r}$ for every $r>0$. 

Apply Theorem \ref{bogomolov_introduction}(2) to
$$\theta=\arccos\sqrt{\frac{1.13}{g}}, \quad \kappa=2.$$ 
As a consequence, under the natural map
$$
\coprod_{P\in(C\cap\Lambda)_{r,2 r}} \mathrm{cone}(P,\theta)
\longrightarrow V,
$$
every nonzero point has at most $3.2\cdot10^{11}g^{\frac{17}3}$ preimages. 
Take intersections with the sphere $S$ and compare the volumes. We have 
$$\#(C\cap\Lambda)_{r,2r}\cdot \mathrm{vol}(\mathrm{cap}(P,\theta))\leq 3.2\cdot10^{11}g^{\frac{17}3}\mathrm{vol}(S_V).$$
The remaining part is just some volume estimates in Euclidean spaces.
Our final result is 
$$\#(C\cap\Lambda)_{\mathrm{medium}}
\leq 1.64\cdot10^{13}g^7\left(1+\frac{5}{2\sqrt2g}\right)^{n-1}.$$

\subsubsection*{Large points}
Now we describe our idea to bound $C(K)_{\rm large}$, which requires more delicate results in sphere packing. 

We first introduce an important counting function in sphere packing. 
For any positive integer $n$ and any angle
$\theta\in(0,\pi)$, denote by $A(n,\theta)$ the maximal number of nonzero points 
$u_1,\dots, u_N\in \mathbb R^n$ such that any two of them form an angle at least 
$\theta$. 
The condition is also equivalent to that the interiors of the spherical caps
$\mathrm{cap}(u_1,\theta/2),\dots, \mathrm{cap}(u_N,\theta/2)$
are disjoint from each other. 
It follows that 
$A(n,\theta)$ is also the maximal number of spherical caps of angular radius $\theta/2$ in $\RR^n$ such that any two such spherical caps have no common interior points.

Now we return to the bound of large points. 
For any $\theta\in (0,\pi)$, there is a subset $S$ of $(C\cap\Lambda)_\mathrm{large}$
satisfying the following properties.
\begin{enumerate}
    \item For any $P\in(C\cap\Lambda)_\mathrm{large}$, there is a $Q\in S$ such that $\angle(P,Q)\le\theta$ and $|Q|\leq |P|$.
    \item For distinct $Q_1,Q_2\in S$, we have $\angle(Q_1,Q_2)>\theta$. As a consequence, $\# S\le A(n,\theta)$.
\end{enumerate}
This set can be constructed inductively.

For any $\kappa>1$, we have 
\begin{align*}
(C\cap\Lambda)_\mathrm{large}
=&\ \bigcup_{Q\in S}\{P\in C\cap\Lambda:|P|\ge|Q|,\ \angle(P,Q)<\theta\}\\
=&\ \bigcup_{Q\in S}\bigcup_{i=1}^{\infty}\{P\in C\cap\Lambda:\kappa^{i-1}|Q|\le|P|\le\kappa^{i}|Q|,\ \angle(P,Q)<\theta\}\\
=&\ \bigcup_{Q\in S}\bigcup_{i=1}^{\lceil\log_{\kappa}(10^5g^{\frac52})\rceil}\{P\in C\cap\Lambda:\kappa^{i-1}|Q|\le|P|\le\kappa^{i}|Q|,\ \angle(P,Q)<\theta\}.
\end{align*}
Here the third equality is a consequence of the quantitative Vojta inequality in Theorem \ref{thm_vojta_introduction}
by taking $\displaystyle\theta=\arccos\left(\sqrt\frac{1.01}{g}\right)$. 

Take $\kappa=1.15$. For each $Q\in S$ and $i$, suppose we have two distinct points
$$P_1, P_2\in\{P\in C(K)\cap\Lambda:\kappa^{i-1}|Q|\le|P|\le\kappa^i|Q|,\ \angle(P,Q)<\theta\}.$$
Denote by $P_1', P_2'$ the orthogonal projections of $P_1, P_2$ to the linear space 
$Q^\perp=\{x\in V:\langle x,Q\rangle=0\}$. 
A simple estimate gives $\langle P'_1,P'_2\rangle<0$ and hence $\angle(P'_1,P'_2)>\frac\pi2$. 
This strong property implies 
that 
$$\#\{P\in C\cap\Lambda:\kappa^{i-1}|Q|\le|P|\le\kappa^i|Q|,\ \angle(P,Q)<\theta\}\le 
A(n-1,\frac{\pi}{2})\le
n.$$
As a consequence, we have the following result. 
$$\#(C\cap\Lambda)_\mathrm{large}\le \lceil\log_{1.15}(10^5g^{\frac52})\rceil nA(n,\theta).$$
We refer to 
Proposition \ref{intermediate large point} for more details. 

Then the problem is reduced to give a good upper bound of $A(n,\theta)$. 
By volume comparison, we have a bound
$$
A(n,\theta) \leq \frac{\mathrm{vol}(S^{n-1})}{\mathrm{vol}(\mathrm{cap}(P,\theta/2))}.
$$
This bound grows roughly as $(\sqrt 2)^{n}$ when $\theta$ is close to $\pi/2$, and is too large for our purpose of making the Vojta constant close to 1. 

We have two other methods to seek stronger upper bounds of $A(n,\theta)$. 
The first method is due to  
Rankin \cite[Thm. 2]{Rankin}, which eventually leads to 
$$\#(C\cap\Lambda)_\mathrm{large}
\le3.4\cdot10^6g^2\left(1+\frac{5}{4\sqrt{g}}\right)^{n-1}.$$
The second method is due to
Kabatjanski\u i--Leven\v ste\u in \cite{KL}, 
which is more delicate and eventually leads to 
$$
\#(C\cap\Lambda)_\mathrm{large} \leq  2.4\cdot 10^{4} g^8 
\left(1+\frac{3\log g}{g}\right)^{n-1}
$$
for $g\geq 142$. 
This explains the proof of Theorem \ref{thm_quantitativeMordellsubgroup}, 
and also explains the origin of the Vojta constant 
$$c_2(g)=\min\left\{1+\frac{5}{4\sqrt g},1+\frac{3\log g}{g}\right\}. $$

The main body of this paper consists of two parts: Part I for arithmetic estimates and Part II for analytic estimates.  
Up to now, we have sketched the main ideas of Part I. 
On the other hand,
Part II proves three estimates required for Part I, which can be found from the first paragraph of Part I or that of Part II.
In the remaining part of this section, we sketch the main idea to prove Theorem \ref{thmfaltingselkieszhangphiinv_introduction}, which is the most difficult one of the three. The proofs of the other two are based on similar strategies.

\subsection{Systoles of hyperbolic surfaces}

Now we start to sketch our proof of 
Theorem \ref{thmfaltingselkieszhangphiinv_introduction}.
For this, we restrict our attention to smooth projective curves $C$ of genus $g\geq 2$ over $\mathbb{C}$, i.e., compact Riemann surfaces of genus $g\geq 2$. 
Then 
$$G_a=G_{\mathrm{Ar}}:C^2\setminus \Delta\longrightarrow \RR$$
 is just the Arakelov Green function on $C$. 
We have Zhang's $\varphi$-invariant
$$\varphi(C)=-\int_{C^2}G_{\mathrm{Ar}}\, c_1( \Delta, G_{\mathrm{Ar}})^{\wedge 2}$$
and  the {Faltings--Elkies} invariant
$$
\mathrm{FE}(n,C)=\inf_{P_1,\dots, P_n\in C} \sum_{1\leq i<j\leq n} G_{\mathrm{Ar}}(P_i, P_j).
$$

\

Since no direct connection between $\varphi(C)$ and $\mathrm{FE}(C)$ is currently known, our estimation strategy proceeds in two steps: a lower bound for $\varphi (C)$ and a lower bound for $\mathrm{FE} (n,C)$. The key idea is to relate both invariants to a certain geometric quantity, which will serve as an intermediate bridge. We now introduce this quantity.

Let $C$ be a smooth projective curve of genus $g\geq 2$ over $\mathbb{C}$. Then $C$ carries a unique K\"ahler metric $\mu_{\mathrm{KE}}$ of constant curvature $-1$. Since $\dim C =1 $, we identify a K\"ahler metric $\mu$ with its associated volume form and measure. Let $\mathrm{sys}(C)$ denote the 
\emph{systole} of $(C,\mu_{\mathrm{KE}})$, i.e., the length of the shortest non-trivial closed geodesic. 
We will use the quantity 
$$\max \left\{\mathrm{sys}(C)^{-1},1\right\}$$ 
as the intermediate quantity. 

Although it is not seriously needed in our estimates, the behavior of the systole when varying $C$ helps us understand our final results. 
Consider the coarse moduli scheme $\mathcal{M}_{g,\mathbb{C}}$ of curves of genus $g$ over $\mathbb{C}$. 
We have a function
$$\mathrm{sys} : \mathcal{M}_{g,\mathbb{C}} (\mathbb{C})\longrightarrow (0,\infty) ,\quad C\longmapsto \mathrm{sys}(C). $$
The function is continuous on $\mathcal{M}_{g,\mathbb{C}} (\mathbb{C})$ under the Euclidean topology. See also the Teichm\"uller-theoretic construction of $\mathcal{M}_{g,\mathbb{C}}$ in \cite[Section 2.C]{harrismorrison1}.
As $C$ approaches the boundary of the moduli space, 
the systole $\mathrm{sys}(C)$ converges to 0. 

The following estimate shows that, up to constants depending only on $g$, we can bound the $\varphi$-invariant $\varphi(C)$ in terms of the invariant
$\max \left\{\mathrm{sys}(C)^{-1},1\right\}$.

\begin{thm}[combining Theorem \ref{thmglobalestimatezhangphiinvariant} and Corollary \ref{coro_phiinv_upperbound}]
\label{thmglobalestimatezhangphiinvariant_introduction}
Let $C$ be a compact Riemann surface of genus $g \geq 2$. Then 
$$  
 10^{-7} g^{-\frac{5}{3}} \cdot \max \left\{\mathrm{sys}(C)^{-1},1\right\}
  \leq \varphi (C) \leq
10^6 g^5\cdot \max \left\{\mathrm{sys}(C)^{-1},1\right\}.$$
\end{thm}

To understand the theorem abstractly, consider the function 
$$\varphi: \mathcal{M}_{g,\mathbb{C}} (\mathbb{C})\longrightarrow (0,\infty),\quad C\longmapsto \varphi(C).$$
It is smooth and tends to infinity as $C$ approaches the boundary of the moduli space (with logarithmic growth). These facts follow from the interpretation of $\varphi$ as an adelic divisor
by \cite[\S3.2.3]{Yuan_Bogomolov}, and the proof of 
\cite[Theorem 3.10]{Yuan_Bogomolov} considering the boundary behavior. 
Then the loc. cit. concludes that $\varphi (C)\geq  c_0(g)$ for some implicit constant $c_0(g)>0$ by compactness. 
Now we can view the theorem as a comparison of two continuous functions on the moduli space. 
It is worth noting that the bound $\varphi (C)\geq  10^{-7} g^{-\frac{5}{3}}$
from the theorem gives an effective expression for $c_0(g)$. 

The upper bound in Theorem \ref{thmglobalestimatezhangphiinvariant_introduction} is comparatively straightforward. It follows from a comparison between the systole and the eigenvalues of the Laplacian, together with the standard estimate comparing the $L^\infty$ norm and $L^2$ norm of holomorphic sections. This part of the argument is similar in spirit to the classical $L^2$ estimate for the Poisson equation. We also note that the upper bound in Theorem \ref{thmglobalestimatezhangphiinvariant_introduction} will not be used later in the paper.

The lower bound in Theorem \ref{thmglobalestimatezhangphiinvariant_introduction}, however, plays a central role in our analysis. The argument is inspired by Tian’s peak section method as in \cite{tg1, tg5, tg6, tg3, sxz1}, the localization of Bergman kernels due to Ma--Marinescu \cite{mamari1}, and the relationship between loops and harmonic $1$-forms established by Buser--Sarnak \cite{busa1}.

We also need the following result, which bounds the Faltings--Elkies invariant in terms of the systole. 
\begin{prop}[Proposition \ref{propfaltingselkiessystole}]
\label{propfaltingselkiessystole_introduction}
Let $C$ be a compact Riemann surface of genus $g \geq 2$. Then
$$ \mathrm{FE}(n,C) \geq -\frac{1}{4}  n\log n - 400\pi^3(g-1)^3\cdot 
\max \left\{ \mathrm{sys}(C)^{-1},1 \right\}\cdot n.$$
\end{prop}

By combining Theorem \ref{thmglobalestimatezhangphiinvariant_introduction} with Proposition \ref{propfaltingselkiessystole_introduction}, we can obtain a result in the form of Theorem \ref{thmfaltingselkieszhangphiinv_introduction}.
By introducing some additional refinements, we can achieve the better constants in Theorem \ref{thmfaltingselkieszhangphiinv_introduction}.
In the following, we sketch the proofs of Theorem \ref{thmglobalestimatezhangphiinvariant_introduction} and Proposition \ref{propfaltingselkiessystole_introduction}.

\subsection{Estimate of the \texorpdfstring{$\varphi$}{Lg}-invariant}

We now sketch the proof of the lower bound in Theorem \ref{thmglobalestimatezhangphiinvariant_introduction}. For clarity, we only outline the main ideas leading to an estimate of the form $\varphi (C)\geq c(g) \cdot  \max \left\{\mathrm{sys} (C)^{-1},1\right\} $ for some constant $c(g)>0$. A complete proof with an improved value of $c(g)$ will be provided later.

We begin by introducing a lower bound for Zhang’s $\varphi$-invariant that is easier to analyze. The lower bound is derived using a potential-theoretic expression of Zhang’s $\varphi$-invariant (cf. Proposition \ref{propzhangphiinvariantpotential}), which we view as an analytic reformulation of Zhang's original formula (cf. Proposition \ref{propzhangphiinvariantcalculation}).

Let $\{ \alpha_j \}_{j=1}^g $ be an $L^2$-orthonormal basis of $\Gamma (C,\omega_C)$ with respect to the $L^2$ inner product 
$$\langle\alpha,\alpha'\rangle=\frac{i}{2}\int_C \alpha\wedge\bar{\alpha}',$$
 where $\omega_C$ denotes the canonical bundle of $C$. The Arakelov K\"ahler metric is defined by 
$$\mu_{\mathrm{Ar}} = \frac{i}{2g} \sum_{k=1}^g \alpha_k \wedge \bar{\alpha}_k ,$$
and is independent of the choice of basis. By the standard theory of Poisson equations, for each $j$, there exists a unique smooth $\mathbb{R}$-valued function $\mathfrak{u}_{j} $ on $C$ satisfying  
$$dd^c \mathfrak{u}_{j} = i\alpha_j \wedge\bar{\alpha}_j - 2 \mu_{\mathrm{Ar}} ,\quad 
\int_C \mathfrak{u}_{j} \mu_{\mathrm{Ar}} =0,$$
where $d^c = \frac{1}{2\pi i} (\partial -\bar{\partial}) $. Then Proposition \ref{propzhangphiinvariantcalculation} yields the lower bound
$$ \varphi (C) \geq \frac{1}{2}\sum_{j=1}^g \int_{C} d\mathfrak{u}_{j} \wedge d^c \bar{\mathfrak{u}}_{j} .$$
Thus, a suitable choice of $L^2$ orthonormal basis for $\Gamma (C,\omega_C)$ leads directly to explicit lower bounds for the $\varphi$-invariant.

Let $\epsilon=\epsilon (g)\in (0,1)$ be a small constant, and let $\gamma_0 \subset C$ be the shortest non-trivial closed geodesic on $(C,\mu_{\mathrm{KE}})$. We say that $\gamma_0$ is separating if its homology class $  [\gamma_0 ] \in H_1 (C;\mathbb{Z}) $ vanishes. Otherwise, $\gamma_0$ is non-separating. By the classification theorem of compact surfaces, a simple closed geodesic is separating if and only if it divides $C$ into two connected components. See also \cite[Chapter I, 41A]{alhsait1}.

We then decompose the moduli space as
$$ \mathcal{M}_{g,\mathbb{C}} (\mathbb{C}) = \mathcal{M}_{g,\mathrm{thick}}(\epsilon) \cup \mathcal{M}_{g,\mathrm{sep}}(\epsilon)\cup \mathcal{M}_{g,\mathrm{nsep}}(\epsilon) ,$$
where
\begin{align*}
\mathcal{M}_{g,\mathrm{thick}}(\epsilon) &= \{ C\in \mathcal{M}_{g,\mathbb{C}} (\mathbb{C}) : \mathrm{sys} (C) \geq \epsilon \} ,\\
\mathcal{M}_{g,\mathrm{sep}}(\epsilon) &= \{ C\in \mathcal{M}_{g,\mathbb{C}} (\mathbb{C}) : \mathrm{sys} (C) \leq \epsilon ,\textrm{ and $\gamma_0\subset C$ is separating} \} ,\\
\mathcal{M}_{g,\mathrm{nsep}}(\epsilon) &= \{ C\in \mathcal{M}_{g,\mathbb{C}} (\mathbb{C}) : \mathrm{sys} (C) \leq \epsilon ,\textrm{ and $\gamma_0\subset C$ is non-separating} \} .
\end{align*}
We call these the thick, separating, and non-separating parts, respectively. As before, we do not distinguish between a curve $C$ and the corresponding complex point in $\mathcal{M}_{g,\mathbb{C}}$. This decomposition reduces the proof to establishing the existence of positive constants $\epsilon_1 (g),\epsilon_2 (g),\epsilon_3 (g)$ such that 
$$\varphi \geq \epsilon_1 \textrm{ on $\mathcal{M}_{g,\mathrm{thick}}(\epsilon)$},\quad \varphi \geq \frac{\epsilon_2}{\mathrm{sys}(C)} \textrm{ on $\mathcal{M}_{g,\mathrm{sep}}(\epsilon)$}, \quad \varphi \geq \frac{\epsilon_3}{\mathrm{sys}(C)} \textrm{ on $\mathcal{M}_{g,\mathrm{nsep}}(\epsilon)$} .$$

We now analyze each case separately. In what follows, we sketch the argument and highlight the main ideas, suppressing certain technical refinements that will be carried out in the full proof later in the paper.

\subsubsection*{Thick part}

Assume that $C\in\mathcal{M}_{g,\mathrm{thick}}(\epsilon)$. Since $C$ is compact, the function $\frac{\mu_{\mathrm{Ar}}}{\mu_{\mathrm{KE}}}$ attains its maximum at some point $x\in C$. Moreover, a classical argument (see also Lemma \ref{lembasicpropertiesBergmankernelfunction}) shows that there exists a section $\alpha\in\Gamma (C,\omega_C)$ such that
$$\frac{i}{2}\int_C \alpha\wedge\bar{\alpha}=1 ,\quad \textrm{ and } \quad \frac{i}{2}\alpha\wedge\bar{\alpha}(x)=\sup_{C} \left|\frac{g\mu_{\mathrm{Ar}}}{\mu_{\mathrm{KE}}}\right|\mu_{\mathrm{KE}} (x) .$$
By a unitary transformation, we can assume that the $L^2$-orthonormal basis $\{ \alpha_j \}_{j=1}^g $ satisfies $\alpha_1=\alpha$.

Let $(\mathbb{D},\mu_{\mathbb{D}})$ be the Poincar\'e disk. Then there exists a universal covering map $\rho:\mathbb{D}\to C$ such that $\rho(0)=x$ and $\rho^* \mu_{\mathrm{KE}} = \frac{4\mu_{\mathrm{Euc}}}{(1-|z|^2)^2} =\mu_{\mathbb{D}} $, where $z$ is the standard coordinate on $\mathbb{D}$, and $\mu_{\mathrm{Euc}} = \frac{idz\wedge d\bar{z}}{2} $ is the standard Euclidean K\"ahler metric on $\mathbb{C}$. See also Theorem \ref{thmpoincaremodel}. 

Let $f$ be the holomorphic function on $\mathbb{D}$ such that $\rho^*\alpha =f(z)dz$. Let $\tilde{\mathfrak{u}}_1=\rho^*\mathfrak{u}_1 $ and define the radial average $\tilde{\mathfrak{u}}_{1,\mathrm{radial}} (z)=\frac{1}{2\pi} \int_0^{2\pi} \tilde{\mathfrak{u}}_{1} (ze^{i\theta}) d\theta$. By definition,
$$ |f(0)|^2 = \frac{i\rho^*(\alpha\wedge\bar{\alpha})(0)}{2\mu_{\mathrm{Euc}} (0)} = \frac{2i\alpha\wedge \bar{\alpha} (x)}{\mu_{\mathrm{KE}} (x)} = 4g \sup_{C} \left|\frac{\mu_{\mathrm{Ar}}}{\mu_{\mathrm{KE}}}\right| \geq 4g\frac{\mu_{\mathrm{Ar}}(C)}{\mu_{\mathrm{KE}}(C)} = \frac{4g}{4\pi (g-1)} >\frac{1}{\pi} .$$
Moreover,
$$ dd^c \tilde{\mathfrak{u}}_{1} = i|f|^2 dz\wedge d\bar{z} - 2\rho^* \mu_{\mathrm{Ar}} \geq 2|f|^2 \mu_{\mathrm{Euc}} - \frac{|f(0)|^2}{4g} \rho^*\mu_{\mathrm{KE}} \geq \left(2|f|^2 - \frac{|f(0)|^2}{g(1-|z|^2)^2}\right) \mu_{\mathrm{Euc}} .$$
Since $f$ is holomorphic, the function $|f|^2$ is subharmonic and satisfies the mean value inequality. Taking the radial average yields
$$dd^c \tilde{\mathfrak{u}}_{1,\mathrm{radial}} (z) \geq \left(2-\frac{1}{g(1-|z|^2)^2}\right) |f(0)|^2 \mu_{\mathrm{Euc}} > |f(0)|^2 \mu_{\mathrm{Euc}} > \frac{1}{\pi} \mu_{\mathrm{Euc}} ,\quad \forall |z|\leq \frac{1}{10}.$$
By the Laplacian inequality for the radial function $\tilde{\mathfrak{u}}_{1,\mathrm{radial}}$, we obtain
\begin{align*}
\int_{\mathbb{D}_{\frac{1}{10}}} d\tilde{\mathfrak{u}}_{1,\mathrm{radial}} \wedge d^c \bar{\tilde{\mathfrak{u}}}_{1,\mathrm{radial}}&\ \geq   \int_{\mathbb{D}_{\frac{
1}{10}}} \frac{1}{4}d|z|^2 \wedge d^c |z|^2 = \frac{1}{40000} .
\end{align*}

By the thick-thin decomposition and the Cauchy--Schwarz inequality, we obtain
$$ \varphi (C) > \frac{1}{2} \int_{C} d\mathfrak{u}_{1} \wedge d^c \bar{\mathfrak{u}}_{1} > \frac{5\cdot \mathrm{sys}(C)}{ 2+10\cdot \mathrm{sys}(C)} \int_{\mathbb{D}_{\frac{1}{10}}} d\tilde{\mathfrak{u}}_{1} \wedge d^c \bar{\tilde{\mathfrak{u}}}_{1} \geq \frac{5\cdot \mathrm{sys}(C)}{ 2+10\cdot \mathrm{sys}(C)}\int_{\mathbb{D}_{\frac{1}{10}}} d\tilde{\mathfrak{u}}_{1,\mathrm{radial}} \wedge d^c \bar{\tilde{\mathfrak{u}}}_{1,\mathrm{radial}} .$$
where the second inequality follows from Lemma \ref{lemlocalstructureorbitcollartheorem}. Hence $\varphi (C)>10^{-6}\epsilon$.

\subsubsection*{Separating part}

Assume that $C\in\mathcal{M}_{g,\mathrm{sep}}(\epsilon)$. Since $\gamma_0 $ is separating, $C\setminus \gamma_0$ has two connected components, which we denote by $C_1$ and $C_2$. There exists a smooth $\mathbb{R}$-valued function $\eta$ on $C$ such that $0\leq \eta\leq 1$,
$$ \int_C d\eta\wedge d^c\eta \leq \mathrm{sys} (C) ,$$
and
$$ \sum_{j=1}^g \left( i\int_{C_1} |\eta|\alpha_j \wedge \bar{\alpha}_j + i\int_{C_2} |\eta -1|\alpha_j \wedge \bar{\alpha}_j \right) \leq 10^{-10}g^{-10} .$$
The construction of such a cutoff function $\eta$ is standard (see \cite[Theorem 8.1.3]{bu1}) and relies on the fact that, when $\gamma_0$ is separating, the Arakelov measure $\mu_{\mathrm{Ar}}$ is uniformly small in a neighborhood of $\gamma_0$. See also Proposition \ref{propexplicitestimateskahlercollar}.

Using this function $\eta$ and the Cauchy--Schwarz inequality, we obtain
\begin{align*}
\varphi (C) &\ \geq \frac{1}{2}\sum_{j=1}^g \int_{C} d\mathfrak{u}_{j} \wedge d^c \bar{\mathfrak{u}}_{j} \geq \frac{\sum_{j=1}^g \left|\int_{C} d\eta \wedge d^c \mathfrak{u}_{j}\right|^2}{2\int_C d\eta\wedge d^c\eta} \\
&\ \geq \frac{1}{2\cdot \mathrm{sys}(C)} \sum_{j=1}^g \left|\int_{C} \eta dd^c \mathfrak{u}_{j}\right|^2 = \frac{1}{2\cdot \mathrm{sys}(C)} \sum_{j=1}^g \left|\int_{C} \eta (i\alpha_j\wedge \bar{\alpha}_j -2\mu_{\mathrm{Ar}})\right|^2 .
\end{align*}

Thus, to obtain a uniform lower bound for $\varphi(C)$, it suffices to construct a holomorphic $1$-form $\alpha\in\Gamma(C,\omega_C)$ with $\frac{i}{2}\int_C \alpha\wedge \bar{\alpha}=1$, such that the quantity $ \left|\int_{C} \eta (i\alpha \wedge \bar{\alpha}  -2\mu_{\mathrm{Ar}})\right| $ is bounded below by a positive constant depending only on $g$.

Since $\gamma_0$ is a simple closed geodesic, the Gauss--Bonnet formula implies that both $C_1$ and $C_2$ have positive genus. Choose a smooth loop $\mathbf{c}_0\subset C_1 $ whose homology class $[\mathbf{c}_0]\in H_1 (C;\mathbb{Z}) $ does not vanish. By Poincar\'e duality, let $\beta$ be the harmonic $1$-form representing the class dual to $[\mathbf{c}_0]$.

Since $\gamma_0$ is a separating closed short geodesic, there exists a smooth closed $1$-form $\beta'$ cohomologous to $\beta$, supported in $C_1$, such that
$$\int_C \beta' \wedge \star\beta'\leq \int_{C_1} \beta \wedge \star\beta + \epsilon'_2 (g,\epsilon)\int_C \beta \wedge \star\beta,$$
where $\epsilon'_2$ is an explicit function depending only on $g$ and $\epsilon$, with $\epsilon'_2(g,\epsilon)\to 0$ as $\epsilon\to 0$ for each fixed $g\geq2$. Note that the cohomology class is determined by its integrals along cycles; in particular, this can be verified using suitable loops $\mathbf{c}_1$. See Figure 1. Here we again use the fact that the Arakelov measure $\mu_{\mathrm{Ar}}$
is uniformly small in a neighborhood of $\gamma_0$. By the $L^2$-minimality of harmonic forms (see also Lemma \ref{lemexplicitpeaksectionalongacycle}), we obtain 
$$\int_{C_2} \beta \wedge \star\beta \leq\epsilon'_2 (g,\epsilon)\int_C \beta \wedge \star\beta ,$$
where $\star$ denotes the Hodge $\star$-operator. 

\begin{center}
\begin{tikzpicture}[scale=1.2]

\def\Rxx{6}      
\def\Ryy{3}     
\def\gapgraph{1.5}    
\def\pinchgraph{0.1}  
\def\hxgraph{0.8}      
\def\hygraph{0.4}     
\def\hshortgeo{0.05}     
\def\hround{1}     
\def\hroundx{0.8}     
\def\hroundy{0.6}     
\def\hxgap{0.1}     

\draw[thick]
   (-\Rxx,0) .. controls (-\Rxx,\Ryy) and (-\gapgraph,\pinchgraph) .. (-\hxgap,\pinchgraph);

\draw[thick]
   (\hxgap,\pinchgraph) .. controls (\gapgraph,\pinchgraph) and (\Rxx,\Ryy) .. (\Rxx,0);

\draw[thick]
   (-\Rxx,0) .. controls (-\Rxx,-\Ryy) and (-\gapgraph,-\pinchgraph) .. (-\hxgap,-\pinchgraph);
   
\draw[thick]
   (\hxgap,-\pinchgraph) .. controls (\gapgraph,-\pinchgraph) and (\Rxx,-\Ryy) .. (\Rxx,0);

  \draw[thick] (\Rxx-3.2,0.27) arc[start angle=180,end angle=360,x radius=1.5*\hxgraph,y radius=1.4*\hygraph];
  \draw[thick] (\Rxx-3.04,0) arc[start angle=180,end angle=0,x radius=1.3*\hxgraph,y radius=1.2*\hygraph];

  \draw[thick] (-\Rxx+1,0.65) arc[start angle=180,end angle=360,x radius=\hxgraph,y radius=\hygraph];
  \draw[thick] (-\Rxx+1.16,0.42) arc[start angle=180,end angle=0,x radius=0.8*\hxgraph,y radius=0.8*\hygraph];

\draw[thick] (-\Rxx+1.2,-0.25) arc[start angle=180,end angle=360,x radius=\hxgraph,y radius=\hygraph];
\draw[thick] (-\Rxx+1.36,-0.48) arc[start angle=180,end angle=0,x radius=0.8*\hxgraph,y radius=0.8*\hygraph];

\draw[thick] (-\hxgap,-\pinchgraph) arc[start angle=270,end angle=450,x radius=\hshortgeo,y radius=2*\hshortgeo];
\draw[dashed] (-\hxgap,-\pinchgraph) arc[start angle=270,end angle=-90,x radius=\hshortgeo,y radius=2*\hshortgeo];
\node at (0.1,-0.3) {$\gamma_0$};

\draw[thick] (\hxgap,-\pinchgraph) arc[start angle=270,end angle=450,x radius=\hshortgeo,y radius=2*\hshortgeo];
\draw[thick] (\hxgap,-\pinchgraph) arc[start angle=270,end angle=-90,x radius=\hshortgeo,y radius=2*\hshortgeo];

\draw[thick] (-\Rxx+1.9,-0.15) arc[start angle=270,end angle=450,x radius=0.1,y radius=0.2];
\draw[dashed] (-\Rxx+1.9,-0.15) arc[start angle=270,end angle=-90,x radius=0.1,y radius=0.2];
\node at (-\Rxx+2.2,0) {$\mathbf{c}_0$};

\draw[densely dashed,->] (-\Rxx+\hroundx,\hroundy) arc[start angle=180, end angle=360, x radius=\hround, y radius=0.5*\hround];
\draw[densely dashed,->] (-\Rxx+\hroundx+2*\hround,\hroundy) arc[start angle=0, end angle=180, x radius=\hround, y radius=0.5*\hround];
\node at (-\Rxx+3,0.6) {$\mathbf{c}_1$};

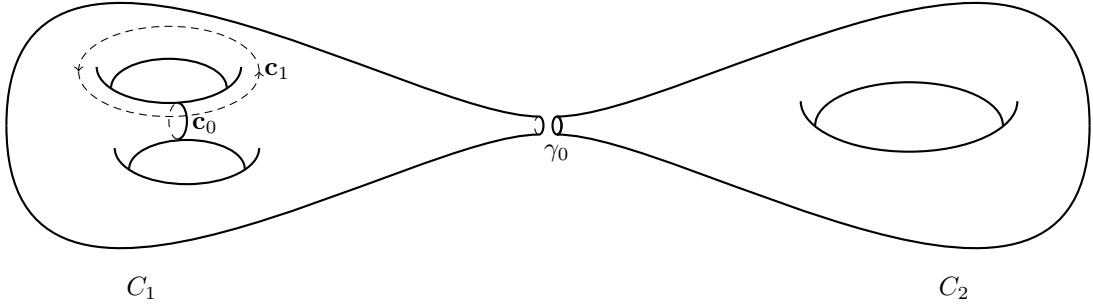
\node at (0,-\Ryy-0.1) {Figure 1: The separating short geodesic $\gamma_0$ and the loops $\mathbf{c}_0, \mathbf{c}_1$};
\node at (-4.5,-1.8) {$C_1$};
\node at (4.5,-1.8) {$C_2$};

\end{tikzpicture}
\end{center}

Set $\alpha = \frac{1}{\sqrt{\int_C\beta\wedge\star\beta}}(\beta +i\star\beta)$. Then $\frac{i}{2}\int_C\alpha\wedge \bar{\alpha}=1$, and $ \frac{i}{2}\int_{C_2} \alpha\wedge \bar{\alpha} \leq\epsilon'_2 $, which implies
$$ \mu_{\mathrm{Ar}} (C_j)\geq (1-\epsilon'_2)g^{-1} ,\quad j=1,2 .$$

Using the second inequality involving $\eta$, we obtain
\begin{align*}
 \left|\int_{C} \eta (i\alpha_1\wedge \bar{\alpha}_1 -2\mu_{\mathrm{Ar}})\right|  &\ \geq \left|\int_{C_2}  (2\mu_{\mathrm{Ar}}-i\alpha_1\wedge \bar{\alpha}_1)\right| - 10^{-8}g^{-8} \\
 &\ \geq 2\mu_{\mathrm{Ar}} (C_2) - 10^{-8}g^{-8}-2\epsilon'_2 \geq 2g^{-1} - 10^{-8}g^{-8}-2\epsilon'_2 ,
\end{align*}
which yields a strictly positive lower bound depending only on $g$, and thus the existence of $\epsilon_2$ follows.

\subsubsection*{Non-separating part}

Assume that $C\in\mathcal{M}_{g,\mathrm{nsep}}(\epsilon)$. By the same argument as in the thick part, there exists an explicit constant $c'(g)>0$ depending only on $g$, such that
$$ \varphi (C) \geq c'(g) \cdot \mathrm{sys} (C) \cdot \sup_{C} \left|\frac{\mu_{\mathrm{Ar}}}{\mu_{\mathrm{KE}}}\right|^2 .$$
It remains to prove that the quantity $\mathrm{sys}(C)\cdot \sup\left|\frac{\mu_{\mathrm{Ar}}}{\mu_{\mathrm{KE}}}\right| $ admits a strictly positive lower bound depending only on $g$. Let $\gamma_0$ be the shortest non-separating closed geodesic on $C$. By definition of $\mu_{\mathrm{Ar}}$, for $\alpha\in\Gamma (C,\omega_C)$ satisfying $\frac{i}{2}\int_{C}\alpha\wedge\bar{\alpha}=1$, we have
$$ \sup_{C} \left|\frac{\mu_{\mathrm{Ar}}}{\mu_{\mathrm{KE}}}\right| \geq \sup_{C} \left|\frac{i\alpha\wedge \bar{\alpha}}{2g\cdot\mu_{\mathrm{KE}}}\right| \geq \frac{2\pi^2}{2g\cdot\mathrm{sys}(C)^2} \cdot \left| \frac{1}{2\pi}\int_{\gamma_0}\alpha \right|^2 .$$
It therefore suffices to construct a holomorphic $1$-form $\alpha$ satisfying $\frac{i}{2}\int_{C}\alpha\wedge\bar{\alpha}=1$, such that the quantity $\frac{1}{\mathrm{sys}(C)} \cdot \left| \frac{1}{2\pi}\int_{\gamma_0}\alpha \right|^2$ has a strictly positive lower bound depending only on $g$.

Let $\beta$ be the harmonic form representing the Poincar\'e dual of $\gamma_0 $. By Lemma \ref{lemexplicitpeaksectionalongacycle}, the desired lower bound can be reduced to a lower bound on the quantity
$$ \frac{1}{\mathrm{sys}(C)} \cdot \int_{C} \beta\wedge\star\beta .$$
The idea is that if this quantity were too small, one could construct a loop $\mathbf{c}$ having a nonzero intersection number with $\gamma_0$ but such that the absolute value of $ \int_{\mathbf{c}}\beta $ is less than $1$, contradicting the defining property of $\beta$ as the Poincar\'e dual of $\gamma_0$. Let $\mathbf{c}_0$ be the shortest geodesic loop such that $[\gamma_0]\cap[\mathbf{c}_0]\neq 0$. See Figure 2.

\begin{center}
\begin{tikzpicture}[scale=1.2]

\def\Rxx{6}      
\def\Ryy{3}     
\def\gapgraph{1.5}    
\def\pinchgraph{0.5}  
\def\hxgraph{0.8}      
\def\hygraph{0.4}     
\def\hshortgeo{0.1}     
\def\hround{1.2}     
\def\hroundx{0.6}     
\def\hroundy{0.6}     
\def\radiussmalldisks{0.2}     

\draw[thick]
   (-\Rxx,0) .. controls (-\Rxx,\Ryy) and (-\gapgraph,\pinchgraph) .. (0,\pinchgraph)
   .. controls (\gapgraph,\pinchgraph) and (\Rxx,\Ryy) .. (\Rxx,0);

\draw[thick]
   (-\Rxx,0) .. controls (-\Rxx,-\Ryy) and (-\gapgraph,-\pinchgraph) .. (0,-\pinchgraph)
   .. controls (\gapgraph,-\pinchgraph) and (\Rxx,-\Ryy) .. (\Rxx,0);

  \draw[thick] (\Rxx-3.2,0.27) arc[start angle=180,end angle=360,x radius=1.5*\hxgraph,y radius=1.4*\hygraph];
  \draw[thick] (\Rxx-3.04,0) arc[start angle=180,end angle=0,x radius=1.3*\hxgraph,y radius=1.2*\hygraph];

  \draw[thick] (-\Rxx+1,0.65) arc[start angle=180,end angle=360,x radius=\hxgraph,y radius=\hygraph];
  \draw[thick] (-\Rxx+1.16,0.42) arc[start angle=180,end angle=0,x radius=0.8*\hxgraph,y radius=0.8*\hygraph];

\draw[thick] (-\Rxx+1.2,-0.25) arc[start angle=180,end angle=360,x radius=\hxgraph,y radius=\hygraph];
\draw[thick] (-\Rxx+1.36,-0.48) arc[start angle=180,end angle=0,x radius=0.8*\hxgraph,y radius=0.8*\hygraph];

\draw[thick] (-\Rxx+1.9,-0.15) arc[start angle=270,end angle=450,x radius=\hshortgeo,y radius=2*\hshortgeo];
\draw[dashed] (-\Rxx+1.9,-0.15) arc[start angle=270,end angle=-90,x radius=\hshortgeo,y radius=2*\hshortgeo];
\node[text = blue] at (-\Rxx+1.97,-0.37) {$\mathscr{C}(\gamma_0)$};

\draw[densely dashed,->] (-\Rxx+\hroundx,\hroundy) arc[start angle=180, end angle=360, x radius=\hround, y radius=0.5*\hround];
\draw[densely dashed,->] (-\Rxx+\hroundx+2*\hround,\hroundy) arc[start angle=0, end angle=180, x radius=\hround, y radius=0.5*\hround];
\node at (-\Rxx+3.2,0.7) {$\mathbf{c}_0$};

\draw[thick] (-\Rxx+1.5,-0.265) arc[start angle=270,end angle=450,x radius=0.134,y radius=0.268];
\draw[dashed] (-\Rxx+1.5,0.271) arc[start angle=90,end angle=270,x radius=0.134,y radius=0.268];

\draw[thick] (-\Rxx+2.34,-0.19) arc[start angle=270,end angle=450,x radius=0.134,y radius=0.268];
\draw[dashed] (-\Rxx+2.34,-0.19) arc[start angle=270,end angle=-90,x radius=0.134,y radius=0.268];

\draw[fill=blue, opacity=0.1]
    (-\Rxx+2.34,0.34) arc[start angle=310,end angle=245,x radius=\hxgraph,y radius=\hygraph]
    -- (-\Rxx+1.5,0.271) arc[start angle=90,end angle=270,x radius=0.134,y radius=0.268]
    -- (-\Rxx+1.5,-0.28)
    arc[start angle=140,end angle=57,x radius=0.8*\hxgraph,y radius=0.8*\hygraph]
    -- (-\Rxx+2.34,-0.19) arc[start angle=270,end angle=450,x radius=0.134,y radius=0.268]
    -- cycle;

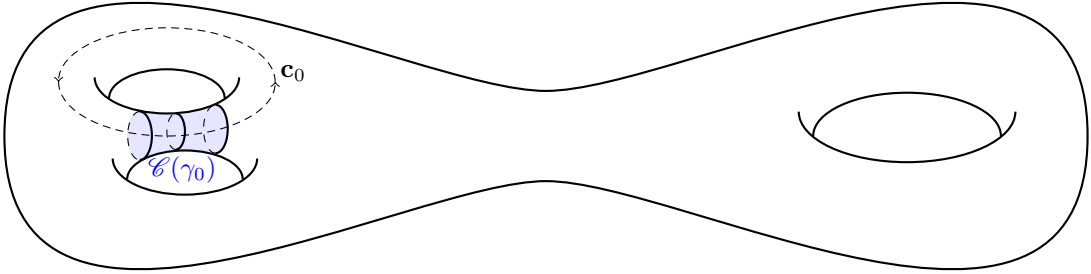
\node at (0,-2.5) {Figure 2: The collar $\mathscr{C}(\gamma_0)$ and the loop $\mathbf{c}_0$};

\end{tikzpicture}
\end{center}

By combining the collar theorem with the Cheeger--Colding's segment inequality (cf. Theorem \ref{thmcheegercoldingsegmentinequality}),
we obtain a piecewise smooth loop $\mathbf{c}_{\gamma_0}$ homotopic to $\mathbf{c}_0$ such that
$$ \left|\int_{\mathbf{c}_{\gamma_0}} \beta \right| \leq \frac{c''(g)}{\sqrt{\mathrm{sys}(C)}} \cdot \sqrt{\int_C \beta\wedge\star\beta} ,$$
where $c''(g)$ is a constant depending only on $g$. Note that our previous smallness condition for $\beta$ concerns integrals over a real surface, whereas here we are interested in an integral taken along a real curve. By the Cheeger--Colding's segment inequality, which can be viewed as an analogue of Fubini’s theorem, we obtain the desired estimate by integrating along a suitable perturbation curve $\mathbf{c}_{\gamma_0}$ of $\mathbf{c}_0$. See Figure 3.

Since $\mathbf{c}_{\gamma_0}$ is close to $\mathbf{c}_{0}$, it has the same nonzero intersection number with $\gamma_0$ as $\mathbf{c}_0$. Therefore,
$$  \int_C \beta\wedge\star\beta \geq \frac{ \mathrm{sys}(C) }{c''(g)^2} \left|\int_{\mathbf{c}_{\gamma_0}} \beta \right|^2 \geq \frac{ \mathrm{sys}(C) }{c''(g)^2} .$$
Hence there exists a constant $\epsilon_3$.

Combining the estimates for $\epsilon_1$, $\epsilon_2$, and $\epsilon_3$, we obtain the lower bound of $\varphi (C)$ in Theorem \ref{thmglobalestimatezhangphiinvariant_introduction}. 

\begin{center}
\begin{tikzpicture}[scale=4]
\def\Rxx{1}      
\def\Ryy{3}     
\def\hxgraph{0.8}      
\def\hygraph{0.4}     
\def\hshortgeo{0.1}     
\def\hround{1.2}     
\def\hroundx{0.6}     
\def\hroundy{0.6}     
\def\radiussmalldisks{0.2}     
  \draw[thick] (-\Rxx+1,0.65) arc[start angle=180,end angle=360,x radius=\hxgraph,y radius=\hygraph];
  \draw[thick] (-\Rxx+1.16,0.42) arc[start angle=180,end angle=0,x radius=0.8*\hxgraph,y radius=0.8*\hygraph];
\draw[thick] (-\Rxx+1.2,-0.25) arc[start angle=180,end angle=360,x radius=\hxgraph,y radius=\hygraph];
\draw[thick] (-\Rxx+1.36,-0.48) arc[start angle=180,end angle=0,x radius=0.8*\hxgraph,y radius=0.8*\hygraph];
\draw[thick] (-\Rxx+1.9,-0.15) arc[start angle=270,end angle=450,x radius=\hshortgeo,y radius=2*\hshortgeo];
\draw[dashed] (-\Rxx+1.9,-0.15) arc[start angle=270,end angle=-90,x radius=\hshortgeo,y radius=2*\hshortgeo];
\node[text = blue] at (-\Rxx+1.95,-0.25) {$\mathscr{C}(\gamma_0)$};
\draw[densely dashed,->] (-\Rxx+\hroundx,\hroundy) arc[start angle=180, end angle=360, x radius=\hround, y radius=0.5*\hround];
\draw[densely dashed,->] (-\Rxx+\hroundx+2*\hround,\hroundy) arc[start angle=0, end angle=180, x radius=\hround, y radius=0.5*\hround];
\node at (-\Rxx+3.1,0.7) {$\mathbf{c}_0$};
\draw[thick] (-\Rxx+1.5,-0.265) arc[start angle=270,end angle=450,x radius=0.134,y radius=0.268];
\draw[dashed] (-\Rxx+1.5,0.271) arc[start angle=90,end angle=270,x radius=0.134,y radius=0.268];
\draw[thick] (-\Rxx+2.34,-0.205) arc[start angle=270,end angle=450,x radius=0.138,y radius=0.276];
\draw[dashed] (-\Rxx+2.34,-0.205) arc[start angle=270,end angle=-90,x radius=0.138,y radius=0.276];
\draw[fill=blue, opacity=0.1]
    (-\Rxx+2.33,0.345) arc[start angle=310,end angle=245,x radius=\hxgraph,y radius=\hygraph]
    -- (-\Rxx+1.5,0.271) arc[start angle=90,end angle=280,x radius=0.134,y radius=0.268]
    -- (-\Rxx+1.5,-0.275)
    arc[start angle=140,end angle=57,x radius=0.8*\hxgraph,y radius=0.8*\hygraph]
    -- (-\Rxx+2.34,-0.205) arc[start angle=270,end angle=450,x radius=0.138,y radius=0.276]
    -- cycle;
\draw[fill=red, opacity=0.1] (-4.3+6-\Rxx,0.03) circle (\radiussmalldisks);
\draw[fill=red, opacity=0.1] (-4.7+6-\Rxx,0.06) circle (\radiussmalldisks);
\draw[fill=red, opacity=0.1] (-5.05+6-\Rxx,0.17) circle (\radiussmalldisks);
\draw[fill=red, opacity=0.1] (-5.35+6-\Rxx,0.4) circle (\radiussmalldisks);
\draw[fill=red, opacity=0.1] (-4.3+6-\Rxx,2*\hroundy-0.01) circle (\radiussmalldisks);
\draw[fill=red, opacity=0.1] (-4.7+6-\Rxx,2*\hroundy-0.06) circle (\radiussmalldisks);
\draw[fill=red, opacity=0.1] (-5.05+6-\Rxx,2*\hroundy-0.17) circle (\radiussmalldisks);
\draw[fill=red, opacity=0.1] (-5.35+6-\Rxx,2*\hroundy-0.4) circle (\radiussmalldisks);
\draw[fill=red, opacity=0.1] (-3.8+6-\Rxx,0.06) circle (\radiussmalldisks);
\draw[fill=red, opacity=0.1] (-3.45+6-\Rxx,0.12) circle (\radiussmalldisks);
\draw[fill=red, opacity=0.1] (-3.15+6-\Rxx,0.32) circle (\radiussmalldisks);
\draw[fill=red, opacity=0.1] (-3+6-\Rxx,0.6) circle (\radiussmalldisks);
\draw[fill=red, opacity=0.1] (-3.9+6-\Rxx,2*\hroundy-0.01) circle (\radiussmalldisks);
\draw[fill=red, opacity=0.1] (-3.55+6-\Rxx,2*\hroundy-0.1) circle (\radiussmalldisks);
\draw[fill=red, opacity=0.1] (-3.2+6-\Rxx,2*\hroundy-0.27) circle (\radiussmalldisks);
\def\mypoints{(-3.88+6-\Rxx,0.11) -- (-3.35+6-\Rxx,0.18) -- (-3.15+6-\Rxx,0.32) -- (-3+6-\Rxx,0.6) -- (-3.2+6-\Rxx,2*\hroundy-0.27) -- (-3.59+6-\Rxx,2*\hroundy-0.21) -- (-3.9+6-\Rxx,2*\hroundy-0.01) -- (-4.3+6-\Rxx,2*\hroundy+0.1) -- (-4.7+6-\Rxx,2*\hroundy-0.09) -- (-5.15+6-\Rxx,2*\hroundy-0.15) -- (-5.25+6-\Rxx,2*\hroundy-0.4) --  (-5.45+6-\Rxx,0.4) -- (-5.05+6-\Rxx,0.17) -- (-4.7+6-\Rxx,0.16) -- (-4.3+6-\Rxx,0.09)}
\draw[densely dotted,green,
      postaction={decorate},
      decoration={markings, mark=between positions 0.05 and 1 step 0.1 with {\arrow{>}},},
     ] \mypoints;

\draw[densely dotted,green,
      postaction={decorate},
      decoration={markings, mark=at position 0.5 with {\arrow{>}}}
     ] (-4.3+6-\Rxx,0.09) -- (-3.88+6-\Rxx,0.11);

\fill[green]   (-3.88+6-\Rxx,0.11) circle (0.3pt); 
\fill[green]  (-3.35+6-\Rxx,0.18) circle (0.3pt);  
\fill[green]  (-3.15+6-\Rxx,0.32) circle (0.3pt); 
\fill[green]  (-3+6-\Rxx,0.6) circle (0.3pt);  
\fill[green] (-3.2+6-\Rxx,2*\hroundy-0.27) circle (0.3pt); 
\fill[green]   (-3.59+6-\Rxx,2*\hroundy-0.21) circle (0.3pt); 
\fill[green]  (-3.9+6-\Rxx,2*\hroundy-0.01) circle (0.3pt);
\fill[green]  (-4.3+6-\Rxx,2*\hroundy+0.1) circle (0.3pt); 
\fill[green]  (-4.7+6-\Rxx,2*\hroundy-0.09) circle (0.3pt);  
\fill[green] (-5.15+6-\Rxx,2*\hroundy-0.15) circle (0.3pt);
\fill[green]   (-5.25+6-\Rxx,2*\hroundy-0.4) circle (0.3pt); 
\fill[green]  (-5.45+6-\Rxx,0.4) circle (0.3pt);  
\fill[green]  (-5.05+6-\Rxx,0.17) circle (0.3pt); 
\fill[green]  (-4.7+6-\Rxx,0.16) circle (0.3pt);  
\fill[green] (-4.3+6-\Rxx,0.09) circle (0.3pt);


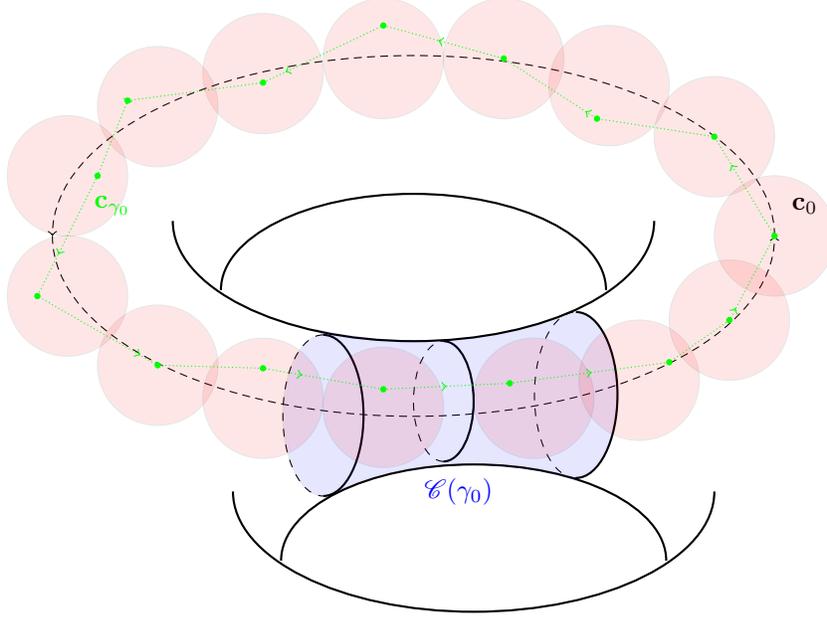
\node at (0.8,-1) {Figure 3: Construction of the piecewise smooth loop $\mathbf{c}_{\gamma_0}$};
\node[text=green] at (-0.2,0.7) {$\mathbf{c}_{\gamma_0}$};
\end{tikzpicture}
\end{center}

\subsection{Estimate of the Faltings--Elkies invariant}

Now we sketch our proof of Proposition \ref{propfaltingselkiessystole_introduction}. 
Inspired by the works of Autissier \cite[Proposition 3.2.5]{autissier1} and Faltings--Elkies (cf. \cite[Theorem 5.1]{lang1}), our key idea to prove Proposition \ref{propfaltingselkiessystole_introduction} is to recast the Arakelov Green function  into one involving hyperbolic geometric data, achieved by convolving the hyperbolic heat kernel with the K\"ahler potential of the Arakelov K\"ahler form.

By the definition of $\mathrm{FE}(n,C)$, it suffices to show that for pairwise distinct points $P_1,\ldots,P_n$ on $C$, we have
$$ \sum_{1\leq j< k\leq n} G_{\mathrm{Ar}} (P_j,P_k) \geq
-\frac{401}{4}  n\log n -\frac{401}{4}  n\log n - \frac{400\pi^3(g-1)^3\max\left\{ 1,\mathrm{sys}(C) \right\}}{\mathrm{sys}(C)} n.$$
Recall that $G_{\mathrm{Ar}}$ is the Green function associated with the Arakelov measure $\mu_{\mathrm{Ar}}$.

Our strategy is to reduce this to an estimate involving the hyperbolic measures $\mu_{\mathrm{KE}}$ and $\mu_{\mathrm{hyp}}=\frac{1}{4\pi (g-1)}\mu_{\mathrm{KE}}$. Note that $\mu_{\mathrm{hyp}}(C)=\mu_{\mathrm{Ar}}(C)=1$. As a quick fact (cf. Theorem \ref{thmdifferencearakelovgreenhyperbolicgreen}), there exists a unique $\mathbb{R}$-valued smooth function $\psi_{\mathrm{Ar}}$ on $C$, such that
$$G_{\mathrm{Ar}} (x,y) - G_{\mathrm{hyp}} (x,y) = \psi_{\mathrm{Ar}} (x) +\psi_{\mathrm{Ar}} (y) ,$$
where $G_{\mathrm{hyp}}$ is the Green function associated with the hyperbolic measure $\mu_{\mathrm{hyp}}$.

Following the standard argument of Autissier \cite[Proposition 3.2.5]{autissier1} and Faltings-Elkies \cite[Theorem 5.1]{lang1}, one applies the hyperbolic heat equation 
$$\frac{\partial u}{\partial t} -\Delta_{\mathrm{hyp}} u=0 $$
 to deduce that for any $t>0$,
$$ \sum_{1\leq j< k\leq n} G_{\mathrm{Ar}} (P_j,P_k) \geq - \frac{1}{2} \sum_{j=1}^n G_{\mathrm{hyp}} (P_j,P_j;t ) -\frac{n(n-1)t}{2} \left( 1 +\sup_C \left| \frac{\mu_{\mathrm{Ar}} }{\mu_{\mathrm{hyp}}} -1 \right| \right) ,$$
where $G_{\mathrm{hyp}}(x,y;t)$ denotes the solution to the hyperbolic heat equation with initial data $\displaystyle \lim_{t\to 0^+ }G_{\mathrm{hyp}} (x,y;t )=G_{\mathrm{hyp}} (x,y )$ in the $L^2$ sense.

Set $t=\frac{1}{40n(g-1)}$. By a local argument, the term $-\frac{n(n-1)t}{2} \left( 1 +\sup_C \left| \frac{\mu_{\mathrm{Ar}} }{\mu_{\mathrm{hyp}}} -1 \right| \right) $ admits a lower bound of the form $-c'(g)n\cdot \mathrm{sys} (C)^{-1}$, where $c'(g)>0$ is a constant depending only on $g$. Thus it remains to estimate the quantity $\displaystyle - n\sup_{x\in C} G_{\mathrm{hyp}} (x,x;t ) $.

We expand $G_{\mathrm{hyp}} (x,x;t )$ in terms of the eigenfunctions $\phi_{\mathrm{hyp},0}, \phi_{\mathrm{hyp},1}, \phi_{\mathrm{hyp},2},\cdots $, where $\phi_{\mathrm{hyp},l}$ is the eigenfunction of the hyperbolic Laplacian corresponding to the eigenvalue $\lambda_{\mathrm{hyp},l}$, with $\lambda_{\mathrm{hyp},0}=0$ and $\phi_{\mathrm{hyp},0}=1$. Then
$$ G_{\mathrm{hyp}} (x,x;t ) = \sum_{0<\lambda_{\mathrm{hyp},l} \leq \frac{g-1}{5}} \frac{e^{-t\lambda_{\mathrm{hyp},l}} |\phi_{\mathrm{hyp} , l }(x) |^2}{ \lambda_{\mathrm{hyp},l}} + \sum_{ \lambda_{\mathrm{hyp},l}>\frac{g-1}{5} } \frac{e^{-t\lambda_{\mathrm{hyp},l}} |\phi_{\mathrm{hyp} , l }(x) |^2}{ \lambda_{\mathrm{hyp},l}} ,$$
and the desired estimate splits accordingly into two parts.

For the first term, by a theorem of Otal--Rosas \cite[Th\'eor\`eme 2]{otro1}, if 
$\lambda_{\mathrm{hyp},l} \leq \frac{g-1}{5} $, then $l \leq 2g-3$. Hence one can control this term via an estimate of the form 
$$\frac{1}{\lambda_{\mathrm{hyp},1}}\leq \frac{c''(g)}{\mathrm{sys}(C)},$$
where $c''(g)>0$ is a constant depending only on $g$. Such an estimate was proved by Schoen--Wolpert--Yau \cite{schwolyau1}; it was later improved by Wu--Xue \cite{wuxue1}. In fact, the estimates of Schoen--Wolpert--Yau and Wu--Xue are stronger: the right-hand side can be taken to be the reciprocal of the minimal total length of a separating multicurve dividing $C$ into two components. Moreover, Wu--Xue's estimate is optimal up to an absolute constant. We shall also provide a new proof using conformal perturbation and Cheeger's inequality (see Proposition \ref{propfirsteigenvalueestimate}).

For the second term, we use the hyperbolic heat kernel $K_{\mathrm{hyp}}(x,y;t)$. By definition,
$$ G_{\mathrm{hyp}}(x,x;2t)-G_{\mathrm{hyp}}(x,x;t)=\int_{t}^{2t} K_{\mathrm{hyp}}(x,y;\varsigma)d\varsigma ,$$
and hence the integral $\int_{t}^{2t} K_{\mathrm{hyp}}(x,y;\varsigma)d\varsigma$ provides an upper bound for the second term. Therefore, the required estimate for the second term follows from an upper bound for the hyperbolic heat kernel.

\subsection{Notation and terminology}

For a real number $x$, denote by $\lfloor x\rfloor$ the greatest integer less than or equal to $x$, and denote by $\lceil x\rceil$ the least integer greater than or equal to $x$. 

For an abelian group $M$, we denote $M_\QQ=M\otimes_\ZZ\QQ$, $M_\RR=M\otimes_\ZZ\RR$, and $M_\CC=M\otimes_\ZZ\CC$. 
The \emph{rank} of $M$ is $\mathrm{rk}(M)=\dim_\QQ (M_\QQ).$

By a \emph{variety} over a field, we mean an integral scheme separated of finite type over the field. By a \emph{curve} over a field, we mean a smooth, projective, and geometrically connected variety of dimension 1 over the field. 

We usually write tensor products of line bundles additively. For example, $aL-bM$ means $L^{\otimes a}\otimes M^{\otimes (-b)}$ for line bundles $L,M$ and integers $a,b$.

\subsubsection*{Acknowledgments}
This article is motivated by the previous works on the uniform Mordell conjecture by Vojta, Dimitrov--Gao--Habegger, K\"uhne, Looper--Silverman--Wilms, and Yuan. 
The authors would like to thank Pascal Autissier for sharing his ideas to improve several estimates in an early draft of this paper. 
The authors would also like to thank Ziyang Gao, Yuxin Ge, Henri Guenancia, Philipp Habegger, George Marinescu, Gang Tian, Rafael von Kanel, Yuanqi Wang, Yunhui Wu, Junyi Xie for their helpful conversations and comments in the preparation of this article. 

The second author would also like to thank the support by the Xplorer Prize from the New Cornerstone Science Foundation, the support by the grants NO. 12250004 and NO. 12321001 from the National Science Foundation of China, and the support by the China--Russia Mathematics Center.

The third author is supported by LabEx CIMI. He also thanks the Beijing International Center for Mathematical Research, Peking University, for its hospitality, where part of this work was carried out.

\newpage
\part{Arithmetic Estimates}

This part consists of arithmetic arguments of this paper. 
The main goal of this section is to prove the quantitative Mordell conjecture in Theorem \ref{thm_quantitativeMordellsubgroup} (or equivalently Theorem \ref{mordell_lang}). 
The treatment of this part assumes three analytic results: 
Theorem \ref{thmfaltingselkieszhangphiinv part I},
 Theorem \ref{propestimateArakelovofdiagonal_compareTwoMetrics part I}, and 
Theorem \ref{thmglobalestimatezhangphiinvariant part I}. 
These three analytic results will be proved in Part II of this paper.

\section{Preliminaries on arithmetic intersection theory}
\label{sec prelim}

This section is review some preliminary terminology and results on heights and arithmetic intersection theory. 

\subsection{N\'eron--Tate heights}\label{sec height}

By a \emph{variety} over a field, we mean an integral scheme separated of finite type over the field. 
Let us briefly recall the definition of Weil heights on projective varieties and N\'eron--Tate heights on abelian varieties. For a detailed introduction, we refer to \cite{Serre, HS, BG}. 

\subsubsection*{Weil heights}
Let $K$ be a number field. Denote by $M_K$ the set of places of $K$. 
Denote by $M_{K,\infty}$ (resp. $M_{K,f}$) the set of archimedean places 
(resp. non-archimedean places) of $K$. 
For any place $v$ of $K$, normalize the absolute value $|\cdot|_v$ on the completion $K_v$ as follows.
\begin{enumerate}[(1)]
\item If $v$ is an archimedean place, take $|\cdot|_v$ to be the usual absolute value on $K_v\cong{\mathbb{R}},{\mathbb C}$;
\item If $v$ is a non-archimedean place, set $|a|_v=\big(\#(O_{K_v}/aO_{K_v})\big)^{-1}$ for any $a\in O_{K_v}$. 
\end{enumerate}
The valuations satisfy the product formula 
$$\prod_{v\in M_K} |a|_v^{\epsilon_v} =1, \quad \forall \ a\in K^\times. $$
Here $\epsilon_v=1$ for all real or non-archimedean places $v$; $\epsilon_v=2$ for all complex places $v$.

Let ${\mathbb P}^n$ be the projective space of dimension $n$ over $K$.
The \emph{standard height function $h: {\mathbb P}^n(\overline K)\to {\mathbb{R}}$} is defined to be
$$h(x_0, x_1, \cdots, x_n)=\frac{1}{[K':{\mathbb Q}]}\sum_{w\in M_{K'}} \epsilon_w \log \max\{|x_0|_w, |x_1|_w,\cdots, |x_m|_w \},$$
where $K'$ is a finite extension of $K$ containing all the coordinates $x_i$. It is independent of the choice of the homogeneous coordinate by the product formula.

Let $X$ be a projective variety over $K$, and $L$ an ample line bundle on $X$. 
Let $i: X\to {\mathbb P}^n$ be any morphism such that $i^*\mathcal O(1)\cong L^{\otimes e}$ for some positive integer $e\geq 1$. 
We obtain a height function 
$$h_{L,i}= \frac1e h\circ i : X(\overline K)\longrightarrow {\mathbb{R}}$$ 
as the composition of 
$i: X(\overline K)\to {\mathbb P}^n(\overline K)$ and $\frac1e h: {\mathbb P}^n(\overline K)\to {\mathbb{R}}$. 
It depends on the choices of $d$ and $i$.

More generally, let $L$ be any line bundle on $X$. We can always write $L= A_1\otimes A_2^{\otimes(-1)}$ for two ample line bundles $A_1$ and $A_2$ on $X$. For $k=1,2$, let $i_k: X\to {\mathbb P}^{n_k}$ be any two morphisms such that $i_k^*\mathcal O(1)=A_k^{\otimes e_k}$ for some positive integer $e_k$. We obtain a height function 
$$h_{L,i_1,i_2}= h_{A_1,i_1}-h_{A_2,i_2}: X(\overline K)\longrightarrow {\mathbb{R}}.$$ 
It depends on the choices of $(A_1, A_2, i_1,i_2)$.
However, the following result asserts that it is unique up to bounded functions.

\begin{thm}[Weil's height machine] \label{Weil}
The above construction $L\mapsto h_{L,i_1,i_2}$ gives a group homomorphism
$$
{\mathcal H}: \mathrm{Pic}(X) \longrightarrow  
\frac{\{ \mathrm{functions\ }\phi: X(\overline K)\to {\mathbb{R}} \}}
{ \{ \mathrm{bounded\ functions\ }\phi: X(\overline K)\to {\mathbb{R}} \}}.
$$
\end{thm}
A function $h_L: X(\overline K)\to {\mathbb{R}} $ in the class ${\mathcal H}(L)$ in the theorem is called a \emph{Weil height function} associated to $L$.

\subsubsection*{N\'eron--Tate heights on abelian varieties}
Let $A$ be an abelian variety over a number field $K$. Denote by $[m]:A\to A$ the multiplication by an integer $m$. 
Let $L$ be a symmetric and ample line bundle on $A$. 
Here $L$ is called \emph{symmetric} if $[-1]^*L\cong L$, which implies $[m]^*L\cong L^{\otimes (m^2)}$ for all integers $m$. 
Let $h_L: A(\overline K)\to {\mathbb{R}} $ be any Weil height function associated to $L$.
The \emph{N\'eron--Tate height function} (or the \emph{canonical height function})
$$\hat h_{L}: A(\overline K)\longrightarrow {\mathbb{R}}$$
associated to $L$ is defined by Tate's limit
$$\hat h_{L}(x) =\lim_{n\rightarrow \infty} \frac{1}{4^n}h_{L}(2^n x), \quad 
x\in A(\overline K).$$
An easy argument proves the convergence and independence of the choice of the Weil height function $h_L$. Moreover, $\hat h_L:A(\overline K)\to {\mathbb{R}}$ is also a Weil height function associated to $L$. 
We further have the following quadraticity and positivity properties.

\begin{thm}\label{canonical}
\begin{enumerate}[(1)]
\item The height $\hat h_L(x)\geq 0$ for any $x\in A(\overline K)$, and the equality holds if and only if $x$ is torsion.
\item
The function $\hat{h}_L: A(\overline{K})\to{\mathbb{R}}$ is {quadratic} in the sense that
it satisfies the parallelogram rule 
$$
\hat{h}_L(x+y)+\hat{h}_L(x-y)=2\left(\hat{h}_L(x)+\hat{h}_L(y)\right),\quad \forall x,y\in A(\overline{K}).
$$
Moreover, we have
$$
\hat{h}_L(mx)=m^2\hat{h}_L(x),\quad \forall m\in{\mathbb Z},\quad \forall x\in A(\overline{K}).
$$
\item
The quadratic form $\hat{h}_L: A(\overline{K})\otimes{\mathbb{R}} \to{\mathbb{R}}$, induced by the height function $\hat{h}_L: A(\overline{K})\to{\mathbb{R}}$ via the quadraticity, is positive definite on the real vector space $A(\overline{K})\otimes{\mathbb{R}}$.
\end{enumerate}
\end{thm}

Denote $A({\overline K})_{\mathbb{R}}=A({\overline K})\otimes{\mathbb{R}}$. Fix a symmetric and ample line bundle $L$ on $A$ and its associated N\'eron--Tate height $\hat{h}(\cdot)=\hat{h}_L(\cdot)$.
For convenience, we take the normalization 
$$
\hat h(x)_K=[K:{\mathbb Q}]\cdot \hat h(x), \quad x\in A({\overline K})_{\mathbb{R}}. 
$$
The factor $[K:{\mathbb Q}]$ will bring convenience when comparing heights with arithmetic intersection numbers. Define the N\'eron--Tate height pairing
$$\langle{x,y}\rangle: A(\overline{K})_{\mathbb{R}}\times A(\overline{K})_{\mathbb{R}}\longrightarrow {\mathbb{R}}$$
by 
$$ \langle{x,y}\rangle= \frac12 \left(\hat{h}(x+y)_K-\hat{h}(x)_K-\hat{h}(y)_K\right), \quad x, y\in A({\overline K})_{\mathbb{R}}.$$
The associated norm of the pairing is given by 
$$|x| := \sqrt{\hat{h}(x)_K}, \quad x\in A({\overline K})_{\mathbb{R}}.$$
The angle $\angle(x,y)$ between $x,y\in A({\overline K})_{\mathbb{R}}$ is defined as
$$\angle(x,y):=\arccos\left(\frac{\langle{x,y}\rangle}{|x||y|}\right).$$
The left-hand side is set to be 0 if $x=0$ or $y=0$.

\subsubsection*{N\'eron--Tate heights on curves}

By a \emph{curve} over a field, we mean a smooth and geometrically connected projective variety of dimension 1 over the field. 

Let $C$ be a curve of genus $g\geq 2$ over a number field $K$. Denote by $J$ the Jacobian variety of $C$ over $K$. 
As the multiplication
$$[2g-2]\colon J({{\overline K}})\longrightarrow J({{\overline K}})$$
is surjective, there is a line bundle $\alpha_0$ on $C_{{\overline K}}$ such that $(2g-2)\alpha_0$ is isomorphic to the canonical sheaf $\omega$. 
Replacing $K$ by a finite extension if necessary, we can assume that $\alpha_0$ is actually a line bundle on $C$. 
This process does not change our definition of heights. Nor does the choice of $\alpha_0$.

Consider the Abel--Jacobi embedding
$$
i_{\alpha_0}: C\longrightarrow J, \quad x\longmapsto x-\alpha_0.
$$
Recall that the theta divisor on $J$ is given by
$$
\theta_{\alpha_0}= \underbrace{i_{\alpha_0}(C)+\dots + i_{\alpha_0}(C)}_{g-1 \text{ copies}}.$$
It is well-known that $\theta_{\alpha_0}$ is ample and gives a principal polarization of $J$. 
By \cite[p. 74, eq. (1)]{Serre}, $\theta_{\alpha_0}$ is symmetric in the sense that $[-1]^*\theta_{\alpha_0}$ is linearly equivalent to $\theta_{\alpha_0}$. It defines a N\'eron--Tate height function $\hat{h}=\hat h_{\theta_{\alpha_0}}$ on $J({\overline K})_{\mathbb{R}}$.

We apply $\hat h(\cdot)$, $\hat h(\cdot)_K$, $|\cdot|$, $\langle{\cdot,\cdot}\rangle$ and $\angle(\cdot,\cdot)$ to $C({\overline K})$ via the embedding $i_{\alpha_0}: C\to J$. 
For example, for $x\in C({\overline K})$, we have
$$\hat h(x)=\hat h(i_{\alpha_0}(x))=\hat h(x-\alpha_0)=\hat h_{\theta_{\alpha_0}}(x-\alpha_0), \quad
|x|^2=\hat h(x)_K=[K:{\mathbb Q}]\hat h(x).$$

\subsection{Intersection theory of hermitian line bundles}

Let us introduce some terminology of Arakelov geometry developed by Arakelov \cite{Ara}, Deligne \cite{Deligne}, and Gillet--Soul\'e \cite{GS}. 

Let ${\mathcal{X}}$ be an arithmetic variety of dimension $n+1$, i.e. a projective and flat integral scheme over ${\mathrm{Spec}}({\mathbb Z})$ of absolute dimension $n+1$. 
A \emph{hermitian line bundle} on ${\mathcal{X}}$ is a pair ${\overline{\mathcal{L}}}=({\mathcal{L}},{\|\cdot\|})$, where ${\mathcal{L}}$ is a line bundle on ${\mathcal{X}}$, and ${\|\cdot\|}$ is a smooth hermitian metric of ${\mathcal{L}}({\mathbb C})$ on the complex variety ${\mathcal{X}}({\mathbb C})$ invariant under the complex conjugation. 
If ${\mathcal{X}}({\mathbb C})$ is singular, the smoothness of the metric means that the metric is locally equal to the pull-back of a smooth metric via a closed embedding into a smooth complex manifold. Denote by $\widehat{\mathrm{Pic}}({\mathcal{X}})$ the group of isometry classes of hermitian line bundles on ${\mathcal{X}}$.

If $\dim {\mathcal{X}}=1$ and $\CX$ is normal, then ${\mathcal{X}}={\mathrm{Spec}}(O_K)$ for some number field $K$, and ${\mathcal{L}}$ is an $O_K$-module of rank $1$ and the hermitian metric ${\|\cdot\|}$ is a collection $(\|\cdot\|_\sigma)_{\sigma:K\to\mathbb{C}}$ of metrics on the complex line ${\mathcal{L}}_\sigma({\mathbb C})$. Let $s$ be any nonzero element of ${\mathcal{L}}$. The \emph{arithmetic degree} of ${\mathcal{L}}$ is defined as
$$\widehat{\mathrm{deg}}({\mathcal{L}})=\log\#({\mathcal{L}}/s\mathcal L)-\sum_{\sigma:K\to{\mathbb C}}\log\|s\|_\sigma.$$
It is independent of the choice of $s$ by the product formula. The arithmetic degree map
$$
\widehat{\mathrm{deg}}: \widehat{\mathrm{Pic}}({\mathcal{X}}) \longrightarrow {\mathbb R}
$$
is additive.

In general, let ${\overline{\mathcal L}}_1, {\overline{\mathcal L}}_2, \cdots, {\overline{\mathcal L}}_{n+1}$ be $n+1$ hermitian line bundles on ${\mathcal{X}}$, and let $s_{n+1}$ be any nonzero rational section of ${\mathcal{L}}_{n+1}$ on ${\mathcal{X}}$. 
The \emph{arithmetic intersection number} is defined inductively by 
\begin{eqnarray*}
{\overline{\mathcal L}}_1\cdot {\overline{\mathcal L}}_2 \cdots {\overline{\mathcal L}}_{n+1}
= {\overline{\mathcal L}}_1\cdot {\overline{\mathcal L}}_2 \cdots {\overline{\mathcal L}}_{n} \cdot {\mathrm{wdiv}}(s_{n+1})
- \int_{{\mathcal{X}}({\mathbb C})} \log \|s_{n+1}\| c_1({\overline{\mathcal L}}_1)\wedge\cdots \wedge c_1({\overline{\mathcal L}}_n).
\end{eqnarray*}
The right-hand side depends on the Weil divisor ${\mathrm{wdiv}}(s_{n+1})$ linearly, so it suffices to explain the case that ${\mathcal D}={\mathrm{wdiv}}(s_{n+1})$ is irreducible (and reduced). 

If ${\mathcal D}$ is horizontal in the sense that it is flat over ${\mathbb Z}$, then
$$
{\overline{\mathcal L}}_1\cdot {\overline{\mathcal L}}_2 \cdots {\overline{\mathcal L}}_{n} \cdot {\mathcal D}= {\overline{\mathcal L}}_1|_{\mathcal D}\cdot {\overline{\mathcal L}}_2|_{\mathcal D}\cdots {\overline{\mathcal L}}_{n}|_{\mathcal D}
$$
is an arithmetic intersection on ${\mathcal D}$ defined by the induction hypothesis. 
If ${\mathcal D}$ is a vertical divisor in the sense that it is a variety over $\mathbb F_p$ for some prime $p$, then 
$$
{\overline{\mathcal L}}_1\cdot {\overline{\mathcal L}}_2 \cdots {\overline{\mathcal L}}_{n} \cdot {\mathcal D}= ({\mathcal{L}}_1|_{\mathcal D}\cdot {\mathcal{L}}_2|_{\mathcal D}\cdots {\mathcal{L}}_{n}|_{\mathcal D})\log p.
$$
Here the intersection is the usual intersection on the projective variety $\mathcal D$ over $\mathbb F_p$.
The definition does not depend on the choice of the rational section $s_{n+1}$. 
It gives a symmetric and multi-linear intersection pairing
$$
\widehat{\mathrm{Pic}}({\mathcal{X}})^{n+1} \longrightarrow {\mathbb R}. 
$$
When $n=0$, it is just the arithmetic degree map.

Let ${\mathcal{X}}$ be an arithmetic variety and ${\overline{\mathcal{L}}}=({\mathcal{L}},{\|\cdot\|})$ be a {hermitian line bundle} on ${\mathcal{X}}$. We say that ${\overline{\mathcal{L}}}$ is \emph{nef} if its hermitian metric is semi-positive and the intersection number ${\overline{\mathcal{L}}}\cdot \mathcal{C}\geq0$ for any 1-dimensional closed integral subscheme $\mathcal{C}$ of ${\mathcal{X}}$.

\subsection{Metrized line bundles on Berkovich analytic spaces}

We review the theory of analytic spaces developed by Berkovich \cite{Berkovich}. It allows us to define metrics at non-archimedean places more flexibly.

\subsubsection*{Berkovich analytic spaces}
Let $(K,|\cdot|_K)$ be a complete valuation field. 
For any variety $X$ over $K$, its \emph{analytification} $X^\mathrm{an}$ is the set of pairs $x=(\overline{x},|\cdot|_x)$ consisting of the following data:
\begin{enumerate}
    \item a schematic point $\overline{x}\in X$, called the \emph{center} of $x$;
    \item a multiplicative norm $|\cdot|_x$ on the residue field $K(\overline{x})$ of $\overline{x}$ which restricts to $|\cdot|_K$ on $K$. Denote the completion of $(K(\overline{x}), |\cdot|_x)$ by $H(x)$.
\end{enumerate}
For any pair $(U,f)$ with a regular function $f$ on an open subvariety $U$ of $X$, we have a function 
$|f|:U^\an\to \RR$ sending $x\in U^\an$ to $|f(x)|=|f|_x$. 
The analytification $X^\mathrm{an}$ is endowed with the weakest topology such that 
$U^\an$ is open in $X^\an$ and $|f|:U^\an\to \RR$ is continuous for every pair $(U,f)$.

As $X$ is separated and connected, the space $X^\an$ is Hausdorff and path-connected. If further $X$ is projective, then $X^\mathrm{an}$ is compact. 

Let $K'$ be a finite normal extension of $K$. The analytification $X_{K'}^\mathrm{an}$ of the base change $X_{K'}$ has an action by $\mathrm{Aut}(K'/K)$, which gives an isomorphism $X^{\mathrm{an}}\cong X_{K'}^{\mathrm{an}}/\mathrm{Aut}(K'/K)$. 

If $K=\mathbb C$, the space $X^\mathrm{an}$ is usual complex analytic space $X(\mathbb{C})$ with the Euclidean topology.
If $K=\mathbb{R}$, then $X^{\mathrm{an}}$ is the quotient of $X(\mathbb{C})$ by the complex conjugation.
If $K$ is non-archimedean, $X^\an$ is the \emph{Berkovich analytic space} associated to $X$. 

Let $X,Y$ be two varieties over $K$ and $f:Y\to X$ a morphism. Then there is a natural continuous map $f^{\mathrm{an}}:Y^\mathrm{an}\to X^\mathrm{an}$. In particular, any closed point $x\in X$ induces a point $x^\mathrm{an}\in X^\mathrm{an}$, called a \emph{classical point} (or a \emph{rigid point}). The map
$$\{\mathrm{closed\ points\ of\ }X\}\lra X^\mathrm{an}$$
is injective with dense image.

\subsubsection*{Metrized line bundles}
Let $(K,|\cdot|_K)$ be a complete valuation field which is either archimedean or discrete non-archimedean. 

Let $X$ be a projective variety over $K$ and $L$ a line bundle on $X$. A \emph{metric} of $L$ on $X$ (or on $X^\an$) is a collection $\|\cdot\|=(\|\cdot\|_x)_{x\in X^\mathrm{an}}$ of norms $\|\cdot\|_x$ on the fiber $L(x)=L_{H(x)}$, which is required to be \emph{continuous} in the following sense. 
For any pair $(U,s)$ with a regular section $s$ of $L$ on an open subvariety $U$ of $X$, we have a function 
$\|s\|:U^\an\to \RR$ sending $x\in U^\an$ to $\|s(x)\|=\|s\|_x$. 
The collection is \emph{continuous} if $\|s\|:U^\an\to \RR$ is continuous for every pair $(U,s)$. 
The pair $(L,\|\cdot\|)$ is called a \emph{metrized line bundle} on $X^\an$. 

Let $(K,|\cdot|_K)$ be a discrete non-archimedean field.
Let $X$ a projective variety over $K$. An integral model $\mathcal{X}$ of $X$ over $O_K$ is a projective and flat integral scheme $\mathcal{X}$  over $O_K$ with an isomorphism $\mathcal{X}\otimes_{O_K}K\cong X$. Let $k$ be the residue field of $K$ and $\mathcal{X}_{k}$ the special fiber of $\mathcal X$. There is a reduction map
$$\mathrm{red}:X^\mathrm{an}\lra\mathcal{X}_{k}$$
defined as follows. For $x\in X^{\mathrm{an}}$, let $R(x)=\{a\in H(x):|a|_x\le1\}$ be the valuation ring of $H(x)$. By the valuative criterion, the morphism $\mathrm{Spec}(H(x))\to\mathcal{X}$ extends to a unique morphism $\mathrm{Spec}(R(x))\to\mathcal{X}$. Then $\mathrm{red}(x)$ is defined as the image of the closed point of $\mathrm{Spec}(R(x))$ in $\CX$, which lies in the special fiber $\CX_k$.

Let $L$ be a line bundle on $X$. Let $({\mathcal{X}},{\mathcal{L}}')$ be an integral model of $(X,L^{\otimes e})$ for some positive integer $e$, i.e. $\mathcal{X}$ is an integral model of $X$ over $O_K$ and ${\mathcal{L}}'$ is a line bundle on ${\mathcal{X}}$ which restricts to $L^{\otimes e}$ on $X$. This integral model induces a \emph{model metric} on $L$ in the following way. For $x\in X^\mathrm{an}$, there is a morphism $\mathrm{Spec}(R(x))\to\mathcal{X}$ as above. Then $\mathcal L_{R(x)}'$ is an $R(x)$-submodule of rank $1$ in $L_{H(x)}^{\otimes e}$. There is a unique norm on $L_{H(x)}$ such that $\mathcal{L}_{R(x)}'$ is the unit ball. More precisely, for $s\in L_{H(x)}$,
$$\|s\|_x:=\inf\{|a|_x^{\frac1e}:a\in H(x)^\times,a^{-1}{s}^{\otimes e}\in\mathcal L_{R(x)}'\}.$$
The model metric is called  \emph{semi-positive} if the line bundle $\CL'$ is nef on (every irreducible component of) the special fiber of $\CX$. 

For either archimedean or discrete non-archimedean $K$, let $X$ be a projective variety of dimension $n$ over $K$. 
Let $(L_1,\|\cdot\|_1),\dots,(L_n,\|\cdot\|_n)$ be \emph{semi-positive} metrized line bundles on $X$, to be explained in the following. Then there is a \emph{Monge--Amp\`ere measure}
$$c_1(L_1,\|\cdot\|_1)\wedge\cdots\wedge c_1(L_n,\|\cdot\|_n)$$
on $X^\mathrm{an}$, which is also explained in the following.
\begin{enumerate}
    \item When $K=\mathbb{C}$, the first Chern current $c_1(L,\|\cdot\|)$ of a metrized line bundle $(L,\|\cdot\|)$ is the $(1,1)$-current defined as
    $$c_1(L,\|\cdot\|)=\frac{1}{\pi i}\partial\overline\partial\log\|s\|+\delta_{\mathrm{wdiv}(s)}.$$
    Here $s$ is any nonzero rational section of $L$. The metric (or the metrized line bundle) is \emph{semi-positive} if for any analytic curve $Y$ in $X$, $c_1(L|_Y,\|\cdot\|)$ is a non-negative measure. The Monge--Amp\`ere measure is defined as the wedge product by the theory of Bedford-Taylor \cite{BT}.
    \item When $K=\mathbb{R}$, everything is defined by descending from $\mathbb{C}$. To be precise, a metric (or a metrized line bundle) is \emph{semi-positive} if its base change to $\mathbb C$ is semi-positive, and the Monge--Amp\`ere measure is the push-forward from $\mathbb C$.
    \item When $K$ is discrete non-archimedean, a metric $\|\cdot\|$ on a line bundle $L$ (or the corresponding metrized line bundle) is \emph{semi-positive} if it is a uniform limit of semi-positive model metrics $\|\cdot\|_m$ on $L$; i.e. the sequence of functions ${\|\cdot\|_m}/{\|\cdot\|}$ on $X^\mathrm{an}$ converges uniformly to $1$. In this case, the Monge--Amp\`ere measure will be defined as the weak limit of those for the semi-positive model metrics, so it suffices to define it for model metrics.
Hence, for $i=1,2,\dots,n$, assume that $(L_i,\|\cdot\|_i)$ are induced by integral models $(\mathcal{X}_i,\mathcal L'_i)$ of $(X,L_i^{\otimes e_i})$. We may assume all $\mathcal{X}_i$ are the same by taking a model $\mathcal{X}$ dominating them. The pull-back of $\mathcal{L}_i$ to $\mathcal{X}$ induces the same model metric on $L_i$. Assume that $X$ is normal. Then we may take ${\mathcal{X}}$ to be normal. For any irreducible component $V$ of the special fiber of $\CX$, there is a unique point $\eta_V\in X^\mathrm{an}$, called the \emph{divisorial point} corresponding to $V$, such that the reduction $\mathrm{red}(\eta_V)$ is the generic point of $V$. Denote by $\delta_{\eta_V}$ the Dirac measure supported at $\eta_V$. In this case, Chambert-Loir \cite{CL} defined the Monge--Amp\`ere measure as
    $$c_1(L_1,\|\cdot\|_1)\wedge\cdots\wedge c_1(L_n,\|\cdot\|_n)=\sum_V(\mathcal L_1|_V\cdots\mathcal L_n|_V)\delta_{\eta_V}.$$
    For general $X$, the Monge--Amp\`ere measure is the push-forward from its normalization. 
\end{enumerate}
A metric (or the corresponding metrized line bundle) is \emph{integrable} if it is the quotient of two semi-positive metrics (of different line bundles). Let $(L_1,\|\cdot\|_1),\dots,(L_d,\|\cdot\|_d)$ be integrable metrized line bundles on $X$. The Monge--Amp\`ere measure
$$c_1(L_1,\|\cdot\|_1)\wedge\cdots\wedge c_1(L_d,\|\cdot\|_d)$$
is defined by multi-linearity.

\subsection{Adelic line bundles on projective varieties}

In this subsection, we review the theory of adelic line bundles on projective varieties introduced by Zhang \cite{Zhang_Adelic}. We only consider the adelic metrics which extend continuously to the Berkovich analytic space at all places, or equivalently the ``m-continuous" metrics defined in \cite[\S A.5.1]{YZ}.

Let $K$ be a number field, $X$ a projective variety over $K$, and $L$ a line bundle on $X$. For any place $v\in M_K$, a $K_v$-metric of $L$ on $X$ is a metric of $L_{K_v}$ on $X_{K_v}^\mathrm{an}$. An \emph{adelic metric} on $L$ is a \emph{coherent} collection $(\|\cdot\|_v)_{v\in M_K}$ of $K_v$-metrics $\|\cdot\|_v$ of $L$ on $X$ over all places $v$ of $K$. 
That the collection $(\|\cdot\|_v)_v$ is \emph{coherent} means that, there exists a finite set $S$ of non-archimedean places of $K$ and a (projective and flat) model 
$({\mathcal{X}},{\mathcal{L}})$ of $(X,L)$ over ${\mathrm{Spec}}(O_K)\setminus S$, such that the $K_v$-metric $\|\cdot\|_v$ is the model metric induced by $({\mathcal{X}}_{O_{K_v}},{\mathcal{L}}_{O_{K_v}})$ for all $v\in{\mathrm{Spec}}(O_K)\setminus S$. In the above situation, we write $\overline L= (L, (\|\cdot\|_v)_v)$ and call it an \emph{adelic line bundle on} $X$.

Let $({\mathcal{X}}, {\overline{\mathcal{L}}}')$ be a \emph{arithmetic model} of $(X,L^{\otimes e})$ for some positive integer $e$, i.e., ${\mathcal{X}}$ is an arithmetic variety over $O_K$, and ${\overline{\mathcal{L}}}'=({\mathcal{L}}', \|\cdot\|)$ is a hermitian line bundle on ${\mathcal{X}}$, such that the generic fiber $({\mathcal{X}}_K, {\mathcal{L}}'_K)=(X,L^{\otimes e})$. Then $({\mathcal{X}}, {\overline{\mathcal{L}}}')$ induces an adelic metric $(\|\cdot\|_v)_{v}$ on
$L$. The metrics at non-archimedean places are model metrics, and the metrics at archimedean places are just given by the root of the hermitian metric. 
Such an adelic metric (resp. adelic line bundle) is called a \emph{model adelic metric} (resp. \emph{model adelic line bundle}). 
The model adelic metric (resp. model adelic line bundle) is called \emph{nef} if the hermitian line bundle ${\overline{\mathcal{L}}}'$ is nef on ${\mathcal{X}}$. 

An adelic line bundle $\overline L= (L, (\|\cdot\|_v)_v)$ on $X$ is called \emph{nef} if the adelic metric $(\|\cdot\|_v)_v$ is a uniform limit of nef model adelic metrics on $L$. Namely, there exists a sequence $\{(\|\cdot\|_{m,v})_v\}_m$ of nef model adelic metrics on $L$, and a finite set $S$ of non-archimedean places of $K$, such that $\|\cdot\|_{m,v}=\|\cdot\|_{v}$ for any $v\in{\mathrm{Spec}}(O_K)\setminus 
S$ and any $m$, and such that $\|\cdot\|_{m,v}/\|\cdot\|_{v}$ converges uniformly to 1 at all places $v$. 
An adelic line bundle is called \textit{integrable} if it is isometric to the tensor quotient of two nef adelic line bundles.
The metric at every place of a nef (resp. integrable) adelic metric is nef (resp. integrable). Denote by $\widehat{\mathrm{Pic}}(X)_{\mathrm {int}}$ the group of isometry classes of integrable adelic line bundles on $X$.

Let $X$ be a projective variety of dimension $n$ over $K$. 
By a limit process, the intersection pairing of hermitian line bundles on models of $X$ extends to a symmetric and multi-linear intersection pairing
$$
\widehat{\mathrm{Pic}}(X)_{\rm int}^{n+1} \longrightarrow \mathbb R, \quad 
(\overline L_1, \overline L_2, \dots, \overline L_{n+1})
\longmapsto \overline L_1\cdot \overline L_2 \cdots \overline L_{n+1}.
$$

For any closed subvariety $Y$ of dimension $d$ in $X$, and any integrable adelic line bundles $\overline L_1, \dots,  \overline L_{d+1}$ on $X$, 
denote 
$$
\overline L_1\cdot \overline L_2 \cdots \overline L_{d+1} \cdot Y= \overline L_1|_Y\cdot \overline L_2|_Y \cdots \overline L_{d+1}|_Y.
$$
By linear combination, the definition extends to Chow cycles of dimension $d$ on $X$. 

If $X=\mathrm{Spec}(K)$, then the line bundle $L_1$ on $X$ is just a vector space over $K$ of dimension one.
We simply have 
$$
\widehat{\deg}(\overline L_1)= -\sum_{v\in M_K} \epsilon_v\log \|s\|_v. 
$$
Here $s$ is any nonzero element of $L_1$, and the degree is independent of the choice of $s$ by the product formula. 

Return to the high-dimensional case. We have an induction formula of Chambert-Loir--Thuillier (cf. \cite[Thm. 4.1]{CLT}).
Let $s_{n+1}$ be a nonzero rational section of $L_{n+1}$ on $X$. 
Then their induction formula gives
$$
 \overline{L}_1\cdots \overline{L}_n\cdot \overline{L}_{n+1} 
 =
 \overline{L}_1\cdots \overline{L}_n\cdot \mathrm{wdiv}(s_{n+1})-\sum_{v\in M_K} \epsilon_v
\int_{X_{K_v}^\mathrm{an}}\big( \log \|s_{n+1}\|_{n+1,v} \big) c_1(\overline{L}_1)_v\wedge\cdots \wedge c_1(\overline{L}_n)_v. 
$$
Here $c_1(\overline{L}_{1,v})\wedge\cdots \wedge c_1(\overline{L}_{n,v})$ denotes the Monge--Amp\`ere measure on $X_{K_v}^\mathrm{an}$ induced by the metrics of $\overline{L}_1,\dots, \overline{L}_{n}$ at $v$.

Let $\overline L$ be an adelic line bundle on $X$. The set of \emph{effective sections} of $\overline L$ is
$$\widehat H^0(X,\overline L)=\{s\in H^0(X,L):\|s(x)\|_{\overline L, v,\sup}\le1,\ \forall v\in M_K\}.$$
Here the supremum norm 
$$\|s(x)\|_{\overline L, v,\sup}=\sup_{x\in X_{K_v}^\mathrm{an}}\|s(x)\|_v.$$
Then $\widehat H^0(X,\overline L)$ is a finite set. We  denote
$$\widehat h^0(X,\overline L)=\log\#\widehat H^0(X,\overline L).$$
The volume of $\overline L$ is
$$\widehat{\mathrm{vol}}(\overline L)=\limsup_{n\to\infty}\frac{(d+1)!}{n^{d+1}}\widehat h^0(X,\overline L).$$
Here the limit always exists by the works of Chen \cite{Che08} and Yuan \cite{Yua09}, but we do not need this fact here. 
We have the following arithmetic Siu inequality by proved Yuan \cite{Yuan_big}.
\begin{thm}[arithmetic Siu inequality]
\label{Siu_inequality}
Let $\overline L_1,\overline L_2$ be two nef adelic line bundles on $X$. Then
$$\widehat{\mathrm{vol}}(\overline L_1-\overline L_2)\ge \overline L_1^{d+1}-(d+1)\overline L_1^d\cdot\overline L_2.$$
\end{thm}

For any adelic line bundle $\overline{L}$ on $X$, we have an associated \emph{height function} 
$$
h_{\overline{L}}: X(\overline K)\longrightarrow \mathbb R
$$
defined by 
$$
h_{\overline{L}}(x)=\frac{1}{\deg(\tilde x)} \overline{L}\cdot \tilde x, \quad x\in X(\overline K). 
$$
Here $\tilde x$ denotes the closed point of $X$ corresponding to $x$, and $\deg(\tilde x)=[K(\tilde x):K]$ denotes the degree of the residue field of $\tilde x$ over $K$. 
It is known that $[K:\mathbb Q]^{-1}h_{\overline{L}}$ is a Weil height function on $X(\overline K)$ associated to $L$. If $s$ is an effective section of $\overline L$, then
$$h_{\overline L}(x)\ge0,\quad\forall x\in(X\setminus|\mathrm{div}(s)|)(\overline K).$$

\subsection{Admissible adelic line bundles}
\label{section_ad}

We review Zhang's theory of admissible adelic line bundles on curves over number fields. We refer to \cite{Zhang_Admissible} and \cite[Appendix A]{Yuan_Bogomolov} for more details.

\subsubsection*{Local setting: metrics}

Let $K$ be a complete value field which is either archimedean or discrete non-archimedean. Let $C$ be a curve of genus $g\geq 2$ over $K$. Let $\omega$ be the canonical line bundle on $C$. Let $\Delta:C\to C^2$ be the diagonal. We also denote by $\Delta$ its image in $C^2$.

A metric $\|\cdot\|$ of $\mathcal{O}(\Delta)$ on $C^2$ induces the following data:
\begin{enumerate}
\item a {Green function} $G=-\log\|1\|:(C^2\setminus\Delta)^{\mathrm{an}}\to\mathbb{R}$;
\item a metric of $\omega$ on $C$ via the natural isomorphism $\Delta^*\mathcal{O}(-\Delta)\cong \omega$;
\item a metric of $\mathcal{O}(x)$ on $C_{K'}$ for any finite extension $K'/K$ and any$x\in C(K')=C_{K'}(K')$ via the natural isomorphism $(x,\mathrm{id})^*\mathcal{O}(\Delta)=\mathcal{O}(x)$. Here
$$(x,\mathrm{id}):{\mathrm{Spec}}(K')\times_{{\mathrm{Spec}}(K)}C\longrightarrow C\times_{{\mathrm{Spec}}(K)}C.$$
\end{enumerate}
There exists a unique metric $\|\cdot\|_{\Delta,a}$ of $\mathcal{O}(\Delta)$ on $C^2$, called the \emph{admissible metric}, satisfying the following properties:
\begin{enumerate}
\item It is symmetric; i.e. it is invariant under the swap map
$$C^2\longrightarrow C^2,\quad(x,y)\longmapsto(y,x).$$
As a corollary, the \emph{admissible Green function} $G_a=-\log\|1\|_{\Delta,a}$ is also symmetric.
\item The equality
$$(2g-2)c_1(\mathcal{O}(x),\|\cdot\|_a)=c_1(\omega,\|\cdot\|_a)_{K'}$$
of the Monge--Amp\`ere measure on $C_{K'}^\mathrm{an}$ holds.
\item The integral
$$\int_{C_{K'}^\mathrm{an}}G_a(x,\cdot)c_1(O(x),\|\cdot\|_a)=0.$$
\end{enumerate}
We also call the induced metric on $\omega$ and $\mathcal{O}(x)$ 
\emph{the admissible metric}.

Following Zhang \cite{Zhang_phi}, we have the \emph{local $\varphi$-invariant}
$$\varphi(C):=-\int_{(C^2)^{\mathrm{an}}}G_{a}\, c_1(\overline{\mathcal{O}}(\Delta)_{a})^{\wedge 2}.$$
We will see that this invariant is always non-negative.

If $K=\mathbb{C}$, the admissible metrics $\|\cdot\|_a$ on 
$\omega$ and $\|\cdot\|_{\Delta, a}$ on 
$\CO(\Delta)$
recover the original \emph{Arakelov metrics} introduced by Arakelov \cite{Ara}. 
They are also denoted by $\|\cdot\|_{\mathrm{Ar}}$ and 
$\|\cdot\|_{\Delta,\mathrm{Ar}}$. 
The \emph{Arakelov measure}
$$\mu_{\mathrm{Ar}}=\frac{1}{2g-2}c_1(\omega,\|\cdot\|_{\mathrm{Ar}})$$
on $C^\an=C(\CC)$
has the following description. Let $\alpha_1,\dots,\alpha_g$ be a basis of $H^0(C(\mathbb{C}),\omega)$ orthonormal with respect to the hermitian pairing
$$\langle{\alpha,\beta}\rangle=\frac i2\int_{C(\mathbb{C})}\alpha\wedge\overline\beta.$$
Then
$$\mu_{\mathrm{Ar}}=\frac{i}{2g}\sum_{j=1}^g\alpha_j\wedge\overline\alpha_j.$$

In the archimedean case, Zhang \cite{Zhang_phi} obtained a spectral expansion of the local $\varphi$-invariant, which implies the positivity $\varphi(C)>0$. 
See also \S\ref{sectionzhangphi1} for more details.

\subsubsection*{Local setting: graph theory}
In the above local setting, 
assume that $K$ is a non-archimedean field with a discrete valuation. 
Assume that $C$ has semi-stable reduction over $O_K$, which can be obtained by passing to a finite extension of $K$. Let $\mathcal{C}$ be the minimal regular model of $C$ over $O_{K}$ and $\omega_{\mathcal{C}/O_{K}}$ its relative dualizing sheaf. The reduction graph $\Gamma=\Gamma(C)$ of $C$ is the graph determined by the following conditions.
\begin{enumerate}
    \item The vertices of $\Gamma$ correspond to the irreducible components of the special fiber $\mathcal{C}_{k}$ of $\mathcal C$.
    \item The edges of $\Gamma$ correspond to the nodes of $\mathcal{C}_{k}$.
    \item Two (possibly equal) vertices of $\Gamma$ are connected by an edge if and only if the irreducible components of $\CCC_k$ represented by the vertices intersect at a node.
\end{enumerate}
There is a natural injective map $\Gamma\to C^\mathrm{an}$, and we identify $\Gamma$ with its image in $C^\mathrm{an}$. The vertex corresponding to a component $V$ is its divisorial point $\eta_V$. There is a retraction map $r:C^\mathrm{an}\to\Gamma$. For a node $\eta_V \in\Gamma$ corresponding to a component $V$, $r^{-1}(\eta_V)$ consists of those points $x\in C^{\mathrm{an}}$ whose reduction $\mathrm{red}(x)$ lies in the smooth locus of $V$. For any piecewise smooth function $f:\Gamma\to\mathbb R$, define the measure
$$\Delta f=-f''dx-\sum_P\left(\sum_{\vec{v}} d_{\vec{v}}f(P)\right)\delta_P,$$
where
\begin{enumerate}
    \item each edge is identified with the closed interval $[0,1]$, on which $f''$ is the second derivative function of $f$ and $dx$ is the Lebesgue measure;
    \item the first summation is over all points $P\in \Gamma$ and $\delta_P$ is the Dirac measure at $P$, and the summation has only finitely many nonzero terms;
    \item the second summation is over all tangent directions $\vec v$ at $P$ and $d_{\vec v}f(P)$ is the directional derivative.
\end{enumerate}
The \emph{canonical divisor} $K_\Gamma$ is the formal linear combination
$$K_\Gamma=\sum_{V}\deg(\omega_{\mathcal{C}/O_K}|_{V})\eta_V,$$
where the summation is over all components $V$ of the special fiber, and $\eta_V$ is the corresponding node. There is a unique pair of a piecewise smooth function $G_\Gamma:\Gamma^2\to\mathbb{R}$ and a measure $\mu_\Gamma$ on $\Gamma$ satisfying 
$$\Delta g_\Gamma(x,\cdot)=\delta_x-\mu_\Gamma,\quad \int_{\Gamma}g_\Gamma(x,\cdot)\mu_\Gamma=0,$$
and such that
$$g_\Gamma(K_\Gamma,x)+g_\Gamma(x,x)$$
is a constant for $x\in\Gamma$. By abuse of notation, we also denote by $G
_\Gamma:(C^2)^\mathrm{an}\to\mathbb{R}$ its pull-back via the composition $(C^2)^\mathrm{an}\to (C^\mathrm{an})^2 \stackrel{(r,r)}{\to} \Gamma^2$.

Let $\|\cdot\|_{\omega_{\mathcal{C}/O_K}}$ be the model metric of $\omega$ on $C$ induced by $\omega_{\mathcal{C}/O_K}$. Then
$$-\log\|\cdot\|_a(x)=-\log\|\cdot\|_{\omega_{\mathcal{C}/O_K}}(x)
+{G_\Gamma(x,x)}\log |\varpi|$$
for $x\in C^\mathrm{an}$. Here $\varpi$ is a generator of maximal ideal of $O_K$. 

For a finite extension $K'$ of $K$, let $\CCC'$ be the minimal regular model of $C_{K'}$ over $O_{K'}$. Then $\CCC'\times_{{\mathrm{Spec}}(O_{K'})}\CCC'$ is singular exactly at the products of nodes, and its blow-up $\CX'$ along these singular points is regular. Let $\overline\Delta_{{K'}}$ be the closure of the diagonal of $C_{K'}^2$ in $\CX'$. Its induced model metric on $\mathcal{O}(\Delta_{K'})$ descends to a metric $\|\cdot\|^{(K')}$ on $\mathcal{O}(\Delta)$, and
$$-\log\|\cdot\|_{\Delta,a}(x_1,x_2)=-\log\left(\lim_{K'}\|\cdot\|^{(K')}\right)(x_1,x_2)+G_\Gamma(x_1,x_2)\log |\varpi|^{-1}$$
for distinct points $(x_1,x_2)\in(C^\mathrm{an})^2$ (cf. \cite[\S3.5]{Zhang_phi}). 
As a consequence, we have 
$$G_a(x_1,x_2)=-\log\|1\|_{\Delta,a}(x_1,x_2)\geq G_\Gamma(x_1,x_2)\log |\varpi|^{-1},$$
where $G_a$ denotes the admissible Green function.

For the case $K'=K$, denote $\|\cdot\|_{\overline\Delta}=\|\cdot\|^{(K)}$. If $r(x_1)$ and $r(x_2)$ are vertices in the graph $\Gamma$, then $\mathrm{red}(x_1,x_2)$ lies in the smooth locus of the special fiber. The sequence of divisors $\overline\Delta_{{K'}}$ is stable on this locus. It follows that
$$-\log\|\cdot\|_{\Delta,a}(x_1,x_2)=-\log\|\cdot\|_{\overline{\Delta}}(x_1,x_2)-G_\Gamma(x_1,x_2)\log |\varpi|.$$
Moreover, if $C$ has good reduction, then the admissible metrics on $\omega$ and $\mathcal{O}(\Delta)$ are the model metrics induced by the integral models $\omega_{\mathcal{C}/O_K}$ and $\mathcal{O}(\overline{\Delta})$ respectively.

Finally, define the \emph{$\varphi$-invariant} of the graph $\Gamma$ by 
$$\varphi(\Gamma)=-\frac14\delta(\Gamma)+\frac14\int_\Gamma g_\Gamma(x,x)((10g+2)\mu_\Gamma-\delta_{K_\Gamma}),$$
where $\delta(\Gamma)$ is the number of edges in $\Gamma$. 
The relation with Zhang's local $\varphi$-invariant
$$\varphi(C)=-\int_{(C^2)^{\mathrm{an}}}G_{a}\, c_1(\overline{\mathcal{O}}(\Delta)_{a})^{\wedge 2}$$
is simply 
$$
\varphi(C)=  \varphi(\Gamma) \log|\varpi|^{-1}. 
$$
By Cinkir \cite{cinkir}, we always have $\varphi(C)\geq 0$, where the equality holds if and only if $C$ has good reduction over $O_K$.

\subsubsection*{Global setting}

Let $K$ be a number field. Let $C$ be a curve of genus $g\geq 2$ over $K$. Let $\omega$ be the canonical line bundle on $C$. Let $\Delta:C\to C^2$ be the diagonal. We also denote by $\Delta$ its image in $C^2$.

With the admissible metrics at all places of $K$, we have \emph{admissible adelic line bundles} $\overline{\mathcal{O}}(\Delta)_a,\bar\omega_a$ and $\mathcal{O}(x)_a$ with underlying line bundles $\overline{\mathcal{O}}(\Delta),\bar\omega$ and $\mathcal{O}(x)$, where $x\in C(K')$ is a point for a finite extension $K'$ of $K$. All of them are integrable. Moreover, $\bar\omega_a$ is nef. 

The self-intersection number $\bar\omega_a^2$, called the \emph{admissible volume} of $C/K$, is strictly positive.
This is equivalent to the Bogomolov conjecture proved by Ullmo \cite{Ull}, and we will have a quantitative version form of this positivity in the following.

Following
Zhang \cite{Zhang_phi}, the \emph{global $\varphi$-invariant} of $C$ over $K$ is defined by
$$\varphi(C)=\sum_{v\in M_K}\epsilon_v\varphi(C_{K_v}),$$
where the local $\varphi$-invariant
$$\varphi(C_{K_v})=-\int_{(C^2)_{K_v}^{\mathrm{an}}}G_{a,v}c_1(\overline{\mathcal{O}}(\Delta)_{a})_v^{\wedge 2}.$$
By the above positivity of local $\varphi$-invariants by Zhang and Cinkir, we always have $\varphi(C)>0$. 
The $\varphi$-invariant gives a lower bound of the admissible volume.

\begin{thm}[Zhang \cite{Zhang_phi}, Cinkir \cite{cinkir}, de Jong \cite{dejong_NeronTate}, Wilms \cite{Wilms}]
\label{phi_inequality}
Let $C$ be a curve of genus $g\geq 2$ over a number field $K$. Then
$$\bar{\omega}_a^2\ge\frac{g-1}{2g+1}\varphi(C)>0.$$
Moreover, for $g=2$, $\bar{\omega}_a^2\ge\frac{2}{5}\varphi(C)$; for $g=3$, $\bar{\omega}_a^2\ge\frac13\varphi(C)$; for $g=4$, $\bar{\omega}_a^2\ge\frac{38}{109}\varphi(C)$.
\end{thm}

By Zhang \cite{Zhang_phi}, the difference
$$d(C)=\frac{2g+1}{2(g-1)}\bar{\omega}_a^2-\varphi(C)$$
is equal to the Beilinson--Bloch height of the Gross--Schoen cycle on $C^3$, so it is conjectured to be non-negative. 
By the recent work of Gao--Zhang \cite{GZ24}, $d(C)$ satisfies some type of the Northcott property on the moduli space of curves, and thus $d(C)>0$ for ``almost all'' $C$. For such curves, we have a stronger version of Theorem \ref{phi_inequality}.

The following arithmetic Hodge index theorem is a variant of the model case of the original theorem of Faltings \cite{Faltings} and Hriljac \cite{Hriljac}. 

\begin{thm}[arithmetic Hodge index theorem]
\label{hodge index}
Let $D$ and $E$ be divisors of degree 0 on $C$. Then
$$
\overline{\mathcal{O}}(D)_a\cdot \overline{\mathcal{O}}(E)_a = -2 \langle{D,E}\rangle. 
$$
\end{thm}

Here $\mathcal{O}(D)_a$ is an arbitary admissible extension of $\mathcal{O}(D)$. For example, if $\mathcal O(D)$ is a linear combination of $\omega$ and $\mathcal{O}(x)$ for $x\in C(K)$, then the corresponding combination of $\bar\omega_a$ and $\overline{\mathcal{O}}(x)_a$ is an admissible extension. The right-hand side is the N\'eron--Tate height pairing on the Jacobian variety of $C$.

The following admissible version of the arithmetic adjunction formula of Arakelov \cite{Ara} is a major motivation of the definition of the admissible metrics. 

\begin{thm}[arithmetic adjunction formula]
\label{adjunction}
Let $P$ be a rational point of $C$ over $K$. Then
$$
\overline{\mathcal{O}}(P)_a\cdot \overline{\mathcal{O}}(P)_a = -\bar\omega_a\cdot \overline{\mathcal{O}}(P)_a
= -\bar\omega_a\cdot P, 
$$
\end{thm}

The fundamental theorems give the following clean expressions, which demonstrate the importance of the above theory of admissible adelic line bundles. 

\begin{thm}\label{admissible heights}
\begin{enumerate}
\item 
Let $Q\in C(K)$ be a rational point. Then 
$$h_{\bar{\omega}_a}(Q)=\frac{2g-2}{g}|Q|^2+\frac{\bar\omega_a^2}{4g(g-1)}.$$
\item
Let $P_1,P_2\in C(K)$ be distinct rational points. Then 
$$h_{\overline{\mathcal O}(\Delta)_a}(P_1, P_2)
=\overline{\mathcal{O}}(P_1)_a\cdot \overline{\mathcal{O}}(P_2)_a=\frac{1}{g}(|P_1|^2+|P_2|^2)-2\langle P_1,P_2\rangle-\frac{\bar\omega_a^2}{4g(g-1)}.$$
\end{enumerate}

\end{thm}
\begin{proof}
By Theorem \ref{hodge index} and Theorem \ref{adjunction}, we have
\begin{align*}
-2|Q|^2=&\ \left(\overline{\mathcal{O}}(Q)_a-\frac{\bar\omega_a}{2g-2}\right)^2\\
=&\ \overline{\mathcal{O}}(Q)_a^2-\frac{(\overline{\mathcal{O}}(Q)_a\cdot{\bar\omega_a})}{g-1}+\frac{\bar\omega_a^2}{(2g-2)^2}
=\ -\frac{g}{g-1}h_{\bar\omega_a}(Q)+\frac{\bar\omega_a^2}{(2g-2)^2}.
\end{align*}
This gives (1).
 
For (2), by the induction formula,
\begin{align*}
\overline{\mathcal{O}}(P_1)_a\cdot \overline{\mathcal{O}}(P_2)_a=&\
\overline{\mathcal{O}}(P_1)_a\cdot P_2+\sum_{v\in M_K}\epsilon_v\int_{C_{K_v}^\mathrm{an}}G_{a,v}(P_2,\cdot)c_1(\overline{\mathcal{O}}(P_1)_a)\\
=&\ \overline{\mathcal{O}}(P_1)_a\cdot P_2\\
=&\ \sum_{v\in M_K}\epsilon_vG_{a,v}(P_1,P_2)\\
=&\ \overline{\mathcal{O}}(\Delta)_a\cdot (P_1,P_2). 
\end{align*}
Now Theorem \ref{hodge index} and Theorem \ref{adjunction} imply
\begin{align*}
-2\langle P_1,P_2\rangle=&\left(\overline{\mathcal{O}}(P_1)_a-\frac{\bar\omega_a}{2g-2}\right)\cdot\left(\overline{\mathcal{O}}(P_2)_a-\frac{\bar\omega_a}{2g-2}\right)\\
=&\overline{\mathcal{O}}(P_1)_a\cdot\overline{\mathcal{O}}(P_2)_a-\frac{1}{2g-2}\bar\omega_a\cdot\left(\overline{\mathcal{O}}(P_1)_a+\overline{\mathcal{O}}(P_1)_a\right)+\frac{\bar\omega_a^2}{(2g-2)^2}\\
=&h_{\overline{\mathcal O}(\Delta)_a}(P)-\frac{1}{2g-2}(h_{\bar\omega_a}(P_1)+h_{\bar\omega_a}(P_2))+\frac{\bar\omega_a^2}{(2g-2)^2}.
\end{align*}
Apply (1).
\end{proof}

\subsection{Function fields}
Let $k$ be an algebraically closed field, $B$ a projective curve over $k$ and $K=k(B)$ its function field. We call $K$ a \emph{function field (of one variable)} over $k$. The above theory has an analog over $K$. We sketch it to finish this section.

The set $M_K$ of places of $K$ is naturally bijective to the set of closed points of $B$. 
For any place $v\in M_K$, normalize the absolute value by 
$$|a|_v=e^{-\mathrm{ord}_v(a)},\quad\forall a\in K^\times.$$
Denote by $K_v$ its completion. These absolute values satisfy the product formula
$$\prod_{v\in M_K} |a|_v=1, \quad \forall \ a\in K^\times. $$

Let ${\mathbb P}^n$ be the projective space over $K$.
The \emph{standard height function} is defined as
$$h(x_0, x_1, \cdots, x_n)_K=\frac{1}{[K':K]}\sum_{w\in M_{K'}}-\log \max\{|x_0|_w, |x_1|_w,\cdots, |x_m|_w \}.$$
There is no subfield $\mathbb{Q}$ and the standard height function depends on the base field, as indicated by the subscript. Theorem \ref{Weil} is still valid.

Let $A$ be an abelian variety over $K$ and $L$ a symmetric and ample line bundle on $A$. The N\'eron--Tate height $\hat{h}_L$ is still defined by Tate's limit
$$\hat h_{L}(x) =\lim_{n\rightarrow \infty} \frac{1}{4^n}h_{L}(2^n x), \quad 
x\in A(\overline K).$$
It is quadratic as in Theorem \ref{canonical}(2). However, it is positive only on the quotient space $A(\overline K)/f(\mathrm{tr}_{K/k}(A)(k))$. Here \emph{Chow's $K/k$-trace} $\mathrm{tr}_{K/k}(A)$ is an abelian variety over $k$ with a homomorphism $f:\mathrm{tr}_{K/k}(A)\otimes_kK\to A$ satisfying the universal property that for any abelian variety $A'$ over $k$ with a homomorphism $g:A'\otimes_kK\to A$, $g$ factors through $f$ uniquely. See \cite{Conrad} for its existence.

Let $X$ be a projective variety over $K$. The adelic line bundles on $X$ are defined simply by omitting archimedean data from the definition over number fields. To be precise, let $L$ be a line bundle on $X$. An {adelic metric} on $L$ is a collection $(\|\cdot\|_v)_{v\in M_K}$ of $K_v$-metrics $\|\cdot\|_v$ on $L$ satisfying the coherent condition that there exist an open subscheme $U\subset B$ and a (projective and flat) model 
$({\mathcal{X}},{\mathcal{L}})$ of $(X,L)$ over $U$, such that the $K_v$-metrics $\|\cdot\|_v$ is induced by $({\mathcal{X}}_{O_{K_v}},{\mathcal{L}}_{O_{K_v}})$ for all $v\in U$. 
The intersection number of nef adelic line bundles is the limit of the classical intersection numbers on projective varieties (cf. \cite{Fulton}). The space of effective sections
$$\widehat H^0(X,\overline L)=\{s\in H^0(X,L):\|s(x)\|_{\overline L,v,\sup}\le1,\ \forall v\in M_K\}$$
is a finite-dimensional vector space over $k$. Define
$$\widehat h^0(X,\overline L)=\dim_k\widehat H^0(X,\overline L),$$
$$\widehat{\mathrm{vol}}(\overline L)=\lim_{n\to\infty}\frac{(d+1)!}{n^{d+1}}\widehat h^0(X,\overline L).$$
The geometric version of Theorem \ref{Siu_inequality}, which is a limit version of the original 
Siu inequality, is still valid.

Let $C$ be a curve of genus $g\geq 2$ over $K$. The admissible adelic line bundles are defined in the same way as over number fields, except that the factor $\log |\varpi|=-1$ for quantities from the reduction graph. For example,
$$-\log\|\cdot\|_a(x)=-\log\|\cdot\|_{\omega_{\mathcal{C}/O_{K_v}}}(x)-{G_\Gamma(x,x)}$$
for $x\in C_{K_v}^\mathrm{an}$,
$$-\log\|\cdot\|_a(x_1,x_2)=-\log\|\cdot\|_{\overline\Delta}(x_1,x_2)+G_\Gamma(x_1,x_2)$$
for $(x_1,x_2)\in(C^2)_{K_v}^\mathrm{an}$ satisfying that $r(x_1)$ and $r(x_2)$ are vertices in $\Gamma$.

 Zhang's global $\varphi$-invariant is
$$\varphi(C)=\sum_{v\in M_K}
\varphi(C_{K_v}).$$
Theorem \ref{phi_inequality} still holds in the current setting. 
A major difference is that the admissible volume $\bar\omega_a^2$ is strictly positive if and only if $C$ is non-isotrivial. Here $C$ is \emph{isotrivial} if $C\otimes_K\overline K\cong C_0\otimes_k\overline{ K}$ for some curve $C_0$ over $k$. If $k$ is of characteristic $0$, then Theorem \ref{phi_inequality} can be strengthened to
$$\bar\omega_a^2\ge\frac{2g-2}{2g+1}\varphi(C).$$

\section{Quantitative Vojta inequality}

Recall that the norm $|Q|=\sqrt{\hat h(Q)_K}$ and the angle $\angle(P_1,P_2)$ for algebraic points $Q, P_1, P_2$ on a curve $C$ over a number field $K$ are defined in terms of the N\'eron--Tate heights
as in \S\ref{sec height}. 
The main goal of this section is to prove the following quantitative version of Vojta's inequality. 

\begin{thm}[Theorem \ref{thm_vojta_introduction}]
\label{vojta}
Let $C$ be a curve of genus $g\geq 2$ over a number field $K$. Let $P_1,P_2\in C(K)$ be distinct rational points. Assume
$$|P_1|\ge1.2\cdot10^9g^{\frac{7}{3}}\sqrt{\bar{\omega}_a^2}$$
and
$${|P_2|}{}\ge10^5g^{\frac52}|P_1|.$$
Then
$$\cos\angle(P_1,P_2)\le\sqrt{\frac{1.01}{g}}.$$
\end{thm}

The theorem with implicit constants is due to Vojta \cite{Vojta_Ar}.
Our proof of the quantitative version is based on the framework of 
Yuan \cite{Yuan_Vojta}, and our additional idea is to use admissible adelic line bundles to track the constants.

Let us sketch our idea of the proof before going into details. Let $p_1,p_2:C^2\to C$ be two projections. Consider the (adelic) \emph{Vojta line bundle}
$$\overline{L}=d_1p_1^*\bar\omega_a+d_2p_2^*\bar\omega_a+d((2g-2)\overline{\mathcal{O}}(\Delta)_a-p_1^*\bar\omega_a-p_2^*\bar\omega_a)$$
on $C^2$. Here $d_1,d_2$ and $d$ are positive integers.
By Theorem \ref{admissible heights}, we have the following result. 

\begin{lem}\label{height identity}
Let $P_1,P_2\in C(K)$ be distinct rational points. Then
$$h_{\overline L}(P_1, P_2)=\frac{2g-2}{g}(d_1|P_1|^2+d_2|P_2|^2)-(4g-4)d\langle P_1,P_2\rangle+\frac{d_1+d_2-2gd}{4g(g-1)}\bar\omega_a^2.$$
\end{lem}

Now we take  $(d_1,d_2, d)$ with growth
$$d_1\approx\sqrt{g+\frac{1}{10^5g}}\frac{|P_2|}{|P_1|}d,\quad d_2\approx\sqrt{g+\frac{1}{10^5g}}\frac{|P_1|}{|P_2|}d.$$
This gives
$$h_{\overline L}(P_1,P_2)\approx(4g-4)d|P_1||P_2|\left(\sqrt{\frac{1}{g}+\frac{1}{10^5g^3}}-\cos\angle(P_1,P_2)\right).$$

Then the key is to obtain a suitable explicit lower bound of $h_{\overline L}(P_1,P_2)$, which gives an upper bound of $\cos\angle(P_1,P_2)$ immediately.
Such a lower bound is obtained as follows. First, we construct a small section $s$ of $nL$ for some $n>0$ by applying Yuan's arithmetic Siu inequality. Here ``small" means that
$$-\sum_{v\in M_K}\epsilon_v\log \|s\|_{v,\sup}$$
has a suitable lower bound. If $s$ is non-vanishing at $P=(P_1,P_2)$, this is a lower bound of $h_{n\overline L}(P)$. In general, we apply Vojta's idea to take the ``derivatives" of $s$ to eliminate the zero at $P$ and estimate its norm. The height is bounded below by the supremum norms of $s$ and its multiplicities at $P$. Finally, Dyson's lemma controls the multiplicities. This finishes the proof.

\subsection{Global sparsity}
The goal of this section is to present the following technical estimate, which will be used twice in this paper. 

\begin{thm}[global sparsity] \label{sparsity}
Let $C$ be a curve of genus $g\geq 2$ over a number field $K$. Let $P_1,\dots, P_n\in C(\overline K)$ be pairwise distinct algebraic points.
Then 
$$\sum_{1\le i<j\le n}h_{\overline{\mathcal{O}}(\Delta)_a}(P_i,P_j)\ge-\left(1.2\cdot 10^4 g^{\frac{1}{3}} n\log n+ 3.96\cdot 10^{10} g^{\frac{11}{3}} n \right)\bar{\omega}_a^2.$$
\end{thm}

Recall that the admissible adelic line bundle $\overline\CO(\Delta)_a$ corresponding to the diagonal $\Delta$ of $C^2$ is introduced in \S\ref{section_ad}. 
We may view $h_{\overline{\mathcal{O}}(\Delta)_a}(P_i,P_j)$ as a (logarithmic) global distance between $P_i,P_j$, as it is a summation of local Green functions. Then the theorem asserts that the global distances are bounded below on average, which can also be viewed as a lower bound of some global Faltings--Elkies invariant of $C$. 
We will prove the theorem in the following, and also take the opportunity to introduce our bounds on the (local) Faltings--Elkies invariant in terms of the (local) $\varphi$-invariant.

From \S\ref{section_ad}, we have an admissible Green function 
$G_{a,v}:(C^2\setminus \Delta)_{K_v}^\an \to \RR$ and a local $\varphi$-invariant 
$\varphi(C_{K_v})$ at every place $v$ of $K$. 
For $n\geq2$, we have the Faltings--Elkies invariant
$$
\mathrm{FE}_v(C,n)= \inf_{P_1,\cdots ,P_n\in C(K_v)}  \sum_{1\leq j< k\leq n} G_{a,v} (P_j,P_k).
$$
Let us first recall the following result from Looper--Silverman--Wilms \cite{LSW}, which is based estimates of Baker--Rumely \cite{BR} and Cinkir \cite{cinkir} on the reduction graph.  

\begin{thm} \label{thmfaltingselkieszhangphiinv non-archi}
Let $C$ be a curve of genus $g\geq 2$ over a number field $K$. 
Let $v$ be a non-archimedean place of $K$. 
Then we have
$$ \mathrm{FE}_v(C,n) \geq -\frac{15}{8} n\, \varphi (C_{K_v}).$$
\end{thm}

\begin{proof}
We can assume that $C$ has semistable reduction by base change. We have  
$$G_{a,v} (P_j,P_k)\geq G_{\Gamma(C_{K_v})} (P_j,P_k)\log N_v$$
from \S\ref{section_ad}. Then the result follows from \cite[Proposition 2.1]{LSW}.  
\end{proof}

The following archimedean version of the theorem will be proved in Part II of this paper, which is actually the most difficult estimate in Part II.

\begin{thm}[Theorem \ref{thmfaltingselkieszhangphiinv}]
\label{thmfaltingselkieszhangphiinv part I}
Let $C$ be a curve of genus $g\geq 2$ over a number field $K$. 
Let $v$ be an archimedean place of $K$. 
Then we have
$$ \mathrm{FE}_v(C,n) \geq - \left(4\cdot 10^3 g^{\frac{1}{3}} n\log n+ 1.32\cdot 10^{10} g^{\frac{11}{3}} n \right) \cdot \varphi (C_{K_v}).$$
\end{thm}

Now we are ready to prove Theorem \ref{sparsity}. 

\begin{proof}[Proof of Theorem \ref{sparsity}]
There is a finite extension $K'/K$ such that $P_1,P_2,\dots,P_n\in C(K')$. Taking the base change of $C$ to $K'$ multiplies each term in the equality by $[K':K]$. So we may assume $P_1,P_2,\dots,P_n\in C(K)$.

By Theorem \ref{thmfaltingselkieszhangphiinv non-archi} and Theorem \ref{thmfaltingselkieszhangphiinv part I},
\begin{align*}
\sum_{1\le i<j\le n}h_{\overline{\mathcal{O}}(\Delta)_a}(P_i,P_j)=&\ \sum_{v\in M_K}\sum_{1\le i<j\le n}\epsilon_vG_{a,v}(P_i,P_j)\\
\ge&\ -\sum_{v\mid\infty}\left(4\cdot 10^3 g^{\frac{1}{3}} n\log n+ 1.32\cdot 10^{10} g^{\frac{11}{3}} n \right)\epsilon_v\varphi(C_{K_v})-\sum_{v\nmid\infty}\frac{15}{8}n\varphi(C_{K_v})\\
\ge&\ -\left(4\cdot 10^3 g^{\frac{1}{3}} n\log n+ 1.32\cdot 10^{10} g^{\frac{11}{3}} n \right)\varphi(C).
\end{align*}
Now the result follows from Wilms' estimate 
$\varphi(C)\leq 3 \bar\omega_a^2$ in Theorem \ref{phi_inequality}. 
\end{proof}

\subsection{Quantitative Mumford inequality}

Before Faltings' proof of the Mordell conjecture, Mumford proved an inequality on angles between two points under the assumption that the heights of the points are large and close. 
Mumford's inequality was not sufficient for the Mordell conjecture, but actually inspired Vojta's inequality in history. We refer to \cite[\S3, Cor. 1]{Mum65} and
\cite[\S5.7]{Serre} for Mumford's original work. 
In this subsection, we derive the following quantitative version of Mumford's inequality.  

\begin{thm}[Theorem \ref{mumford_introduction}]
\label{mumford}
Let $C$ be a curve of genus $g\geq 2$ over a number field $K$. Let $P_1,P_2\in C(\overline K)$ be distinct algebraic points. Assume
$$|P_1|\ge 10^{9}g^{\frac{7}{3}}\sqrt{\bar\omega_a^2}$$
and
$$|P_1|\le{|P_2|}{}\le1.15|P_1|.$$
Then
$$\cos\angle(P_1,P_2)\le{\frac{1.01}{g}}.$$
\end{thm}

\begin{proof}
By Theorem \ref{sparsity} for $n=2$, 
$$
h_{\overline{\mathcal O}(\Delta)_a}(P_1, P_2)\ge
 -\left(4.71\cdot 10^6 g^{\frac{1}{3}}\cdot (2\log2)+3.96\cdot 10^{10}g^{\frac{11}{3}}\cdot 2\right) \bar\omega_a^2
 \ge-10^{11}g^{\frac{11}{3}}\bar\omega_a^2.$$
By Theorem \ref{admissible heights}, we further have
$$\frac{1}{g}(|P_1|^2+|P_2|^2)-2\langle P_1,P_2\rangle\ge-10^{11}g^{\frac{11}{3}}\bar\omega_a^2.$$
By the assumption $|P_1|\ge 10^{9}g^{\frac{7}{3}}\sqrt{\bar\omega_a^2}$, we have
$$\frac{1+10^{-7}}{g}|P_1|^2+\frac{1}{g}|P_2|^2-2\langle P_1,P_2\rangle\ge\frac{1}{10^7g}|P_1|^2-10^{11}g^{\frac{11}{3}}\bar\omega_a^2\ge0.$$
By the assumption $|P_1|\le{|P_2|}{}\le1.15|P_1|$, we have
$$\cos\angle(P_1,P_2)=\frac{\langle P_1,P_2\rangle}{|P_1||P_2|}\le\frac{1+10^{-7}}{2g}\frac{|P_1|}{|P_2|}+\frac{1}{2g}\frac{|P_2|}{|P_1|}\le\frac{1.01}{g}.$$
Here the last inequality follows from the fact that the function $(1+10^{-7}) x^{-1}+ x$ for $x\in [1, 1.15]$ takes the biggest value at $x=1.15.$
This finishes the proof. 
\end{proof}

\subsection{Existence of small section}

The goal of this subsection is to find a suitable small section of our version of Vojta's line bundle. This will be our first step to prove the quantitative Vojta inequality in Theorem \ref{vojta}.

Let $C$ be a curve of genus $g\geq 2$ over a number field $K$. 
Let $P_1,P_2\in C(K)$ be distinct rational points as in Theorem \ref{vojta}. 
Denote $P=(P_1,P_2)\in C^2(K)$. 
As above, consider the adelic \emph{Vojta line bundle}
$$\overline{L}=d_1p_1^*\bar\omega_a+d_2p_2^*\bar\omega_a+d((2g-2)\overline{\mathcal{O}}(\Delta)_a-p_1^*\bar\omega_a-p_2^*\bar\omega_a)$$
on $C^2$. Here $d_1,d_2, d$ are positive integers to be determined later. 
The underlying line bundle on $C^2$ is
$${L}=d_1p_1^*\omega+d_2p_2^*\omega+d((2g-2){\mathcal{O}}(\Delta)-p_1^*\omega-p_2^*\omega)$$
with $\omega=\omega_{C/K}$.
The main result of this subsection is as follows. 

\begin{prop}[small section]
\label{small}
If $d_1\ge gd$, $d\ge d_2$, and $d_1d_2>gd^2$, then there exists a positive integer $n$ such that $nL$ has a nonzero section $s$ on $C^2$ satisfying
$$-\sum_{v\in M_K}\epsilon_v\log\|s\|_{n\overline L, v,\sup}\ge
-\frac{ng(g+1)d_1d^2}{4(g-1)(d_1d_2-gd^2)}\bar\omega_a^2.$$
\end{prop}

\begin{proof}
The proof is  almost the same as \cite[Lemma 3.2]{Yu}, except that we need to apply Yuan's arithmetic version of Siu's inequality. For completeness we include it here.
Denote
$$c=\frac{g(g+1)d_1d^2}{4(g-1)(d_1d_2-gd^2)} \bar\omega_a^2.$$
Take an adelic line bundle ${\mathcal{O}}(c)$ on $\mathrm{Spec}(K)$ of degree $c$. Denote $\overline{L}(c)=\overline{L}+\pi^*{\mathcal{O}}(c)$. Here $\pi:C\to\mathrm{Spec}(K)$ is the structure morphism. Then for any nonzero section $s$ of $nL$,
$$\sum_{v\in M_K}\epsilon_v\log\|s\|_{n(\overline L(c)),v,\sup}=\sum_{v\in M_K}\epsilon_v\log \|s\|_{n\overline L,v,\sup}-nc.$$
So the lemma is equivalent to that $n(\overline L(c))$ has a nonzero section $s$ with 
$\log\|s\|_{n(\overline L(c)),v,\sup}\leq 1$ for every place $v$. 
It suffices to prove $\widehat{\mathrm{vol}}(\overline{L}(c))>0$.

Write 
$$\overline L=\overline L_1-\overline L_2$$
with 
$$\overline{L}_1=(d_1-d)p_1^*{\bar\omega_a}+d(g-1)p_2^*{\bar\omega_a}+d(2g-2)\overline{\mathcal{O}}({\Delta})_a,$$
and
$$\overline{L}_2=(gd-d_2)p_2^*{\bar\omega_a}.$$
By \cite[Theorem 2.10(2)]{Yuan_Bogomolov}, the adelic line bundle
$$p_1^*\bar\omega_a+p_2^*\bar\omega_a+2\overline{\mathcal{O}}(\Delta)_a$$
is the pull-back of a nef adelic line bundle from the Jacobian variety $J$, so it is nef. 
Then $\overline L_1$ is nef by the assumption $d_1\ge gd$. 
It follows that $\overline{L}_1(c)$ and $\overline{L}_2$ are nef. 
Apply Theorem \ref{Siu_inequality} to $\overline L(c)=\overline L_1(c)-\overline L_2$. We have
\begin{align*}
\widehat{\mathrm{vol}}(\overline{L}(c))
\ge&\ \overline{L}_1(c)^3-3\overline{L}_1(c)^2\overline{L}_2\\
=&\ (\overline{L}_1(c)-\overline{L}_2)^3-3\overline{L}_1(c)\overline{L}_2^2+\overline{L}_2^3\\
=&\ \overline{L}(c)^3-6(g-1)d_1(gd-d_2)^2\bar{\omega}_a^2\\
=&\ \overline L^3+3cL^2-6(g-1)d_1(gd-d_2)^2\bar{\omega}_a^2\\
=&\ \overline L^3+24(g-1)^2(d_1d_2-gd^2)c-6(g-1)d_1(gd-d_2)^2\bar{\omega}_a^2\\
=&\ \overline L^3+6(g-1)d_1(gd^2+2gd_2d-d_2^2)\bar\omega_a^2.
\end{align*}
Here $\overline L^3$ is a polynomial in $d_1,d_2,d$ of degree $3$. The coefficients are listed below. They can be computed as in the proof of \cite[Theorem 3.6]{Yuan_Bogomolov} or that of \cite[Theorem 1.4]{Wilms}.
\begin{enumerate}
    \item For $i=1,2$, the coefficient of $d_i^3$ is $(p_i^*\bar{\omega}_a)^3=0$.
    \item For $i=1,2$, let $j=3-i$. The coefficient of $d_i^2d_j$ is $3(p_i^*\bar{\omega}_a)^2(p_{j}^*\bar{\omega}_a)=6(g-1)\bar{\omega}_a^2.$
    \item For $i=1,2$, the coefficient of $d_i^2d$ is $3((2g-2)\overline{\mathcal{O}}({\Delta})_a-p_1^*{\bar\omega_a}-p_2^*{\bar\omega_a})(p_i^*\bar{\omega}_a)^2=0.$
    \item The ciefficient if $d_1d_2d$ is $6((2g-2)\overline{\mathcal{O}}({\Delta})_a-p_1^*{\bar\omega_a}-p_2^*{\bar\omega_a})(p_1^*\bar{\omega}_a)(p_2^*\bar{\omega}_a)=-12(g-1)\bar{\omega}_a^2.$
    \item For $i=1,2$, the coefficient of $d_id^2$ is $3((2g-2)\overline{\mathcal{O}}({\Delta})_a-p_1^*{\bar\omega_a}-p_2^*{\bar\omega_a})^2(p_i^*\bar{\omega}_a)=-6(g-1)(2g-1)\bar{\omega}_a^2.$
    \item The coefficient of $d^3$ is
    \begin{align*}
    &((2g-2)\overline{\mathcal{O}}({\Delta})_a-p_1^*{\bar\omega_a}-p_2^*{\bar\omega_a})^3\\
    =&\ (2g-2)(4g^2+4g-2)\bar{\omega}_a^2-(2g-2)^3\varphi(C)\\
    \ge&\ (2g-2)(-4g^2+8g+2)\bar{\omega}_a^2\\
    \ge&\ -12(g-1)g(g-2)\bar{\omega}_a^2.
    \end{align*}
\end{enumerate}
Here the first inequality in (6) comes from Theorem \ref{phi_inequality}. 
Putting these together, we have
$$\widehat{\mathrm{vol}}(\overline{L}(c))\ge6(g-1)\Big(d_1^2d_2+(2g-2)d_1dd_2-(g-1)d_1d^2-(2g-1)d_2d^2-(2g-4)gd^3\Big){\bar\omega_a^2}.$$
By assumption, we have 
$d_1^2d_2>gd_1d^2$ and $d_1dd_2\ge gd^3\ge gd^2d_2$.
It follows that 
    \begin{align*}
  \widehat{\mathrm{vol}}(\overline{L}(c)) 
  >&\ 6(g-1)\Big((2g-2)d_1dd_2-(2g-1)d_2d^2-(2g-4)gd^3\Big){\bar\omega_a^2}\\
  \geq &\ 6(g-1)\Big((2g-2)gd^3-(2g-1)d_2d^2-(2g-4)gd^3\Big){\bar\omega_a^2}\\
  >&\ 0.  
    \end{align*}
   This finishes the proof. 
  \end{proof}

\subsection{Lower bound of the height}

Let $C$ be a curve of genus $g\geq 2$ over a number field $K$. 
Let $P_1,P_2\in C(K)$ be distinct rational points. 
Recall that on $C^2$, have the {adelic line bundle}
$$\overline{L}=d_1p_1^*\bar\omega_a+d_2p_2^*\bar\omega_a+d((2g-2)\overline{\mathcal{O}}(\Delta)_a-p_1^*\bar\omega_a-p_2^*\bar\omega_a)$$
with underlying line bundle
$${L}=d_1p_1^*\omega+d_2p_2^*\omega+d((2g-2){\mathcal{O}}(\Delta)-p_1^*\omega-p_2^*\omega).$$
Here $d, d_1, d_2$ are positive integers.
Suppose that $s$ is a nonzero section of $nL$ on $C^2$ for some positive integer $n$. The goal here is to give a lower bound of the height $h_{n\overline L}(P)$ at the rational point $P=(P_1,P_2)$ of $C^2$ in terms of supremum norms of $s$.

Take a local section $s_0$ of $L$ on a neighborhood of $P$ in $C^2$ with $s_0(P)\neq 0$. Then ${s}/{s_0}$ is a regular function on a neighborhood of $P$. Take a local coordinate $z_i$ of $C$ at $P_i$,  i.e. a generator of the maximal ideal of the local ring $\CO_{C, P_i}$. 
The completion of the local ring of $C^2$ at $P$ is $K[[z_1,z_2]]$. So we have a formal expansion
$$\frac s{s_0}=\sum_{i_1,i_2\geq 0} a_{i_1,i_2}z_1^{i_1}z_2^{i_2}, \quad a_{i_1,i_2}\in K.$$

We endow the set of pairs $(i_1,i_2)$ with the product partial order, i.e. $(i_1,i_2)\le(i_1',i_2')$ if and only if both $i_1\le i_1'$ and $i_2\le i_2'$. A pair $(e_1,e_2)$ is called \emph{admissible} if $a_{e_1,e_2}\ne 0$ and $a_{i_1,i_2}=0$ for all $(i_1,i_2)<(e_1,e_2)$. We claim that the admissibility is independent of the choices of $s_0$ and $(z_1,z_2)$. Indeed, for another local section $s_0'$ of $L$ non-vanishing at $P$, assume
$$\frac s{s_0'}=\sum_{i_1,i_2} a_{i_1,i_2}'z_1^{i_1}z_2^{i_2}.$$
There is an expansion
$$\frac{s_0}{s_0'}=\sum_{j_1,j_2} b_{j_1,j_2}z_1^{j_1}z_2^{j_2}.$$
So
$$a_{i_1,i_2}'=\sum_{(j_1,j_2)\le(i_1,i_2)}a_{i_1-j_1,i_2-j_2}b_{j_1,j_2}.$$
Therefore, $a_{i_1,i_2}'=0$ for $(i_1,i_2)<(e_1,e_2)$, and
$$a_{e_1,e_2}'=a_{e_1,e_2}b_{0,0}=a_{e_1,e_2}\left(\frac{s_0}{s_0'}\right)(P).$$
On the other hand, for another local coordinate $z_1'$ at $P_1$ and
$$\frac s{s_0}=\sum_{i_1,i_2} a_{i_1,i_2}''(z_1')^{i_1}z_2^{i_2},$$
we have
$$z_1=\sum_{j\ge1} b_j(z_1')^j.$$
Here $b_1\ne 0$. So $a_{i_1,i_2}''-a_{i_1,i_2}b_1^{i_1}$ is a linear combination of $a_{j_1,j_2}$ for $(j_1,j_2)<(i_1,i_2)$. Hence, $a_{i_1,i_2}'=0$ for $(i_1,i_2)<(e_1,e_2)$, and
$$a_{e_1,e_2}''=a_{e_1,e_2}b_1^{e_1}=a_{e_1,e_2}\left(\frac{dz_1}{dz_1'}\right)^{e_1}(P).$$
The situation is similar for $z_2$.

Fix an admissible pair $(e_1,e_2)$ for $s$, we have a nonzero section
$$s^\circ=a_{e_1,e_2}s_0\otimes (dz_1)^{\otimes e_1}(dz_2)^{\otimes e_2}$$
of the line bundle
$$L^\circ=(n L+e_1p_1^*\omega+e_2p_2^*\omega)|_P$$
on $P\cong \mathrm{Spec}(K)$. 
This section depends on $s$, but  is independent of the choices  of $s_0$ and $(z_1,z_2)$ by the previous calculations.
It is natural to take the adelic line bundle
$$\overline L^\circ=(n\overline L+e_1p_1^*\bar\omega_{a}+e_2p_2^*\bar\omega_{a})|_P$$
on $P$.
We have the following estimate of the norm of $s^\circ$ under the metric of this adelic line bundle at each place.

\begin{prop}
\label{height_bound}
Let $s$ be a nonzero section of $nL$ and $(e_1,e_2)$ be an admissible pair for $s$. 
\begin{enumerate}
    \item At a non-archimedean place $v$,
    $$-\log\|s^\circ\|_{v}\ge-\log\|s\|_{n\overline L,v,\sup}-\frac{15}{4}(e_1+e_2)\varphi(C_{K_v}).$$
    \item At an archimedean place $v$,
    \begin{align*}
    -\log\|s^\circ\|_{v}\ge-\log\|s\|_{n\overline L,v\sup }-2.1\cdot10^{14}g^{\frac{14}{3}}(2nd_1+2nd_2+4gnd+e_1+e_2)\varphi(C_{K_v}).
    \end{align*}
\end{enumerate}
\end{prop}

We prove two parts of the proposition separately in the following.

\subsubsection*{Non-archimedean case}
Let us first prove Proposition \ref{height_bound}(1). 
We may assume $n=1$ by replacing $(d_1,d_2,d)$ by $(nd_1,nd_2,nd)$, and assume that $C$ has semi-stable reduction at $v$ by passing to a finite extension of $K$.

Let $\mathcal{C}\to\mathrm{Spec}(O_{K_v})$ be the minimal regular model of $C_{K_v}$ over $O_{K_v}$. By the valuative criterion, $P_i\in C(K)$ (resp. $P\in C^2(K)$) extends to $\mathcal P_i\in\mathcal{C}(O_{K_v})$ (resp. $\mathcal{P}\in\mathcal{C}^2(O_{K_v})$). Moreover, the closed points of $\mathcal{P}_
i$ and $\mathcal P$ lie in the smooth locus of the special fiber (cf. \cite[\S9.1, Cor. 1.32]{Liu}). Let $V_i$ (resp. $V$) be the irreducible component of the special fiber of $\mathcal C$ (resp. $\mathcal C^2$) intersecting $\mathcal P_i$ (resp. $\mathcal P$). Let $\eta_i\in C_{K_v}^{\mathrm{an}}$ (resp. $\eta\in(C^2)_{K_v}^{\mathrm{an}}$) be the corresponding divisorial point, i.e. $|f|_{\eta_i}=N_v^{-\mathrm{ord}_{V_i}(f)}$ (resp. $|f|_\eta=N_v^{-\mathrm{ord}_V(f)}$) for any rational function $f$ on $\mathcal C$ (resp. $\mathcal C^2$). It suffices to prove
$$-\log\|s^\circ(P)\|_{\overline L^\circ}\ge-\log\|s(\eta)\|_{\overline L}-\frac{15}4(e_1+e_2)\varphi(C_v).$$
Here we indicate the source of the metrics in the subscripts, and omit the dependence on $v$.

As recalled in \S\ref{section_ad}, the blowing-up $\mathcal{X}$ of $\mathcal C^2$ along all singular points is regular. Denote by $p_i:\mathcal{X}\to\mathcal{C}$ the composition of the blowing-up $\mathcal X\to\mathcal C^2$ and the $i$-th projection $\mathcal{C}^2\to \mathcal{C}$. Let $\overline\Delta$ be the closure of the diagonal $\Delta_{K_v}\subset (C^2)_{K_v}$ in $\mathcal X$. The strict transform $V'$ of $V$ is a vertical prime divisor on $\mathcal X$.

Define the line bundle
$$\mathcal{L}=d_1p_1^*\omega_{\mathcal{C}/O_v}+d_2p_2^*\omega_{\mathcal{C}/O_{K_v}}+d((2g-2)\mathcal{O}(\overline\Delta)-p_1^*\omega_{\mathcal{C}/O_v}-p_2^*\omega_{\mathcal{C}/O_v})$$
on $\mathcal{X}$. 
Then $(\CX,\CL)$ is an integral model of $((C^2)_{K_v}, L_{K_v})$, and thus induces a metric $\|\cdot\|_\CL$ of $L_{K_v}$. 

Similarly, define the line bundle 
$$\mathcal{L}^\circ=(\mathcal{L}+e_1p_1^*\omega_{\mathcal{C}/O_{K_v}}+e_2p_2^*\omega_{\mathcal{C}/O_v})|_{\mathcal P}$$
on $\CP\cong\mathrm{Spec}(O_{K_v})$. 
Then $(\CP,\CL^\circ)$ is an integral model of $(P_{K_v}, (L^\circ)_{K_v})$, and thus induces a metric $\|\cdot\|_{\CL^\circ}$ of $(L^\circ)_{K_v}$.

 Let $r:C_{K_v}^\mathrm{an}\to\Gamma$ be the retraction map to the reduction graph $\Gamma\subset C_{K_v}^\mathrm{an}$, and $G_\Gamma:(C^2)_{K_v}^\mathrm{an}\to\mathbb R$ the pull-back of the Green function of $\Gamma$. 
 As recalled in \S\ref{section_ad}, we have
$$-\log{\|\cdot\|_{\bar\omega_a}(x)}=-\log{\|\cdot\|_{\omega_{\mathcal C/O_v}}(x)}-G_\Gamma(x,x)\log N_v$$
for $x\in C_v^\mathrm{an}$, and
$$-\log{\|\cdot\|_{\overline{\mathcal{O}}(\Delta)_a}(x_1,x_2)}=-\log{\|\cdot\|_{\mathcal{O}(\overline{\Delta})}(x_1,x_2)}+G_\Gamma(x_1,x_2)\log N_v$$
for $(x_1,x_2)\in(C^2)_{K_v}^\mathrm{an}$ such that $r(x_1),r(x_2)\in\Gamma$ are vertices of the graph. Since $r(P_i)=\eta_i$, we have
$$-\log{\|\cdot\|_{\overline L}(\eta)}=-\log{\|\cdot\|_{\mathcal L}(\eta)}-((d_1-d)G_\Gamma(\eta_1,\eta_1)+(d_2-d)G_\Gamma(\eta_2,\eta_2)-(2g-2)dG_\Gamma(\eta_1,\eta_2))\log N_v,$$
$$-\log{\|\cdot\|_{\overline L^\circ}(P)}=-\log{\|\cdot\|_{\mathcal L^\circ}(P)}-((d_1-d+e_1)G_\Gamma(\eta_1,\eta_1)+(d_2-d+e_2)G_\Gamma(\eta_2,\eta_2)-(2g-2)dG_\Gamma(\eta_1,\eta_2))\log N_v,$$
and thus
$$-\log{\|s^\circ(P)\|_{\overline L^\circ}}=-\log{\|s^\circ(P)\|_{\mathcal L^\circ}}-\log{\|s(\eta)\|_{\overline L}}+\log{\|s(\eta)\|_{\mathcal L}}-(e_1G_\Gamma(\eta_1,\eta_1)+e_2G_\Gamma(\eta_2,\eta_2))\log N_v.$$
By \cite[Lemma 2.3]{LSW},
$$G_\Gamma(x,y) \le\frac{15}{4}\varphi(\Gamma),\quad \forall x,y\in C_{K_v}^\an.$$
Thus it suffices to prove 
$$\|s^\circ(P)\|_{\mathcal L^\circ}\leq \|s(\eta)\|_{\mathcal L}.$$

View the section $s$ of $L$ as a rational section of $\CL$ on $\CX$. Then the corresponding divisor $\mathrm{div}_{\CL}(s)$ on $\CX$ is a vertical divisor. 
Denote by $m$ the multiplicity of $V'$ in $\mathrm{div}_{\CL}(s)$.
This gives $\|s(\eta)\|_{\mathcal L}=N_v^{-m}.$
Let $\varpi$ be a generator of the maximal ideal of $O_{K_v}$.
As the multiplicity of $V'$ in $\mathrm{div}_{\CL}(\varpi^{-m}s)$ is 0, we see that $\varpi^{-m}s$, as a rational section of $\CL$ on $\CX$, is regular at every point $x\in V'$ which does not lie in any other irreducible component of the special fiber of $\CX$. 
In particular,  $\varpi^{-m}s$ is a regular section of $\CL$ on a neighborhood of $\CP$.

Recall that the section $s^\circ$ of $L^\circ$ is defined as 
$$s^\circ= a_{e_1,e_2}s_0\otimes (d z_1)^{\otimes e_1}(d z_2)^{\otimes e_2}$$
where $s_0$ is local generator of $L$ at $P$,  $z_1$ (resp. $z_2$) is a local coordinate of $C$ at $P_1$ (resp. $P_2$), and  
$$\frac{s}{s_0}=\sum_{i_1,i_2\geq 0}{a}_{i_1,i_2} z_1^{i_1} z_2^{i_2}, \quad {a}_{i_1,i_2}\in {K}.$$
The definition of $s^\circ$ is independent of the choice of $(s_0, z_1, z_2)$. 
Now we choose $(s_0, z_1, z_2)$ carefully as follows. 
First, we choose $s_0$ to be a local generator of $\CL$ on an open neighborhood of $\CP$. 
Second, for $i=1,2$, as the ideal sheaf $\CI_i$ of $\CP_i$ on $\mathcal C$ is a line bundle, we can choose $z_i$ to be a local generator of $\CI_i$ on an open neighborhood of $\CP_i$. 
As $\mathcal P_i$ lies in the smooth locus of $\mathcal C$ over $O_{K_v}$, the completion of the local ring of $\mathcal{X}$ at the closed point of $\mathcal{P}$ along the ideal sheaf of $\mathcal P$
is isomorphic to $O_{K_v}[[\tilde z_1,\tilde z_2]]$. 

As $\varpi^{-m}s$ is a regular section of $\CL$ over a neighborhood of $\CP$ in $\CX$, the expansion 
$$\frac{\varpi^{-m} s}{s_0}=\sum_{i_1,i_2\geq 0} \varpi^{-m} {a}_{i_1,i_2} z_1^{i_1} z_2^{i_2}$$
has integral coefficients $\varpi^{-m} {a}_{i_1,i_2}\in O_{K_v}.$
It follows that 
$$\varpi^{-m} s^\circ=\varpi^{-m} a_{e_1,e_2}s_0\otimes (d  z_1)^{\otimes e_1}(d  z_2)^{\otimes e_2}$$
is a regular section of the line bundle
$\CL^\circ$ on $\CP$. 
As a consequence,  
$$\|s^\circ(P)\|_{\mathcal L}\leq N_v^{-m}=\|s(\eta)\|_{\mathcal L}.$$
This proves Proposition \ref{height_bound}(1).

\subsubsection*{Archimedean case}

To prove Proposition \ref{height_bound}(2), we start with some analytic terminology. 
We refer to \S\ref{subsec hyperbolic} for more details. 
Denote the discs 
$$\mathbb D=\{z\in\mathbb C:|z|<1\},\quad 
\mathbb{D}_r=\{z\in\mathbb{D}:|z|<r\}.$$
The sheaf $\Omega_{\mathbb D}^1$ of holomorphic 1-forms on $\mathbb D$ is endowed with the hyperbolic metric 
$\|\cdot\|_{{\mathrm{hyp},\mathbb{D}}}$ 
normalized by 
$$\Vert dz\Vert_{\mathrm{hyp} , \mathbb{D}} = 1- |z|^2.$$
The metric is invariant under biholomorphic automorphisms on $\mathbb D$. 
Let $C$ be a (connected) compact Riemann surface of genus $g\geq 2$.
Let $\rho_0:\mathbb{D}\to C$ be a universal covering map. 
Via the canonical identity $\rho_0^* \omega_{C}\cong \Omega_{\mathbb D}^1$, 
the metric $\|\cdot\|_{{\mathrm{hyp},\mathbb{D}}}$ descends to a hyperbolic metric 
$\|\cdot\|_{{\mathrm{hyp}}}$ on $\omega_{C}$. 
  
We need the following theorem from Part II of this paper, which particularly compares the hyperbolic metric $\|\cdot\|_{{\mathrm{hyp}}}$ with the Arakelov metric 
$\|\cdot\|_{{\mathrm{Ar}}}=\|\cdot\|_{a}$
on $\omega_{C}$. 

\begin{thm}[Theorem \ref{propestimateArakelovofdiagonal_compareTwoMetrics}]
\label{propestimateArakelovofdiagonal_compareTwoMetrics part I}
Let $C$ be a compact Riemann surface of genus $g\geq 2$. 
Then the Arakelov metric $\Vert\cdot\Vert_{\mathrm{Ar}}$ and the hyperbolic metric $\Vert\cdot \Vert_{\mathrm{hyp}}$ on $\omega_C$ satisfy
$$ \left| \log \frac{\Vert\cdot\Vert_{\mathrm{Ar}}}{\Vert\cdot\Vert_{\mathrm{hyp}}} \right|  < 2\cdot 10^{14}g^{\frac{14}{3}} \cdot \varphi (C).$$

Let $\Delta$ be the diagonal of $C\times C$.  
For any $(x_1,x_2)\in C\times C$, let $\rho=(\rho_1,\rho_2):\mathbb{D}\times \mathbb{D}\to C\times C$ be a universal covering map with $\rho_1(0)=x_1$ and $\rho_2(0)=x_2$. Then there exists a local holomorphic section $s$ of $(\rho^*\mathcal{O}(\Delta),\rho^*\Vert \cdot\Vert_{\Delta})$ on $\mathbb{D}_{\frac{1}{20}}\times\mathbb{D}_{\frac{1}{20}}$, such that
$$ \big| \log \Vert s\Vert_{\Delta, \mathrm{Ar}} \big|  < 2\cdot 10^{14}g^{\frac{14}{3}} \cdot \varphi (C).$$
\end{thm}

We also need the following neat theorem from Part II, which is also a part of Theorem \ref{thmglobalestimatezhangphiinvariant_introduction}. It gives the explicit constant in \cite[Theorem 3.10]{Yuan_Bogomolov}. 

\begin{thm}[Theorem \ref{thmglobalestimatezhangphiinvariant}]
\label{thmglobalestimatezhangphiinvariant part I}
Let $C$ be a compact Riemann surface of genus $g\geq 2$. 
Then 
$$ \varphi (C) \geq 10^{-7} g^{-\frac{5}{3}}.$$
\end{thm}

Now we are ready to prove Proposition \ref{height_bound}(2). 
Let $C, K, \overline L$ and $P=(P_1, P_2)$ be as in the proposition. 
As in the non-archimedean case, we may assume $n=1$ by replacing $(d_1,d_2,d)$ by $(nd_1,nd_2,nd)$ and assume $K_v=\mathbb{C}$ by replacing $K$ by a finite extension.

Let $\rho=(\rho_1,\rho_2):\mathbb{D}^2\to (C(\mathbb C))^2$ be a holomorphic universal covering map such that $\rho(0,0)=P$. The pull-back $\rho^*\overline L$ is a metrized line bundle on $\mathbb{D}^2$ and $\rho^*s$ is a section of $\rho^*L$. We can also define a section $(\rho^*{s})^\circ$ of
$$(\rho^*\overline{L})^\circ=(\rho^*\overline{L}+e_1p_1^*\rho_1^*\bar\omega_\mathrm{a}+e_2p_2^*\rho_2^*\bar\omega_\mathrm{a})|_{(0,0)}$$
with the analytic data. More precisely, we can take a holomorphic non-vanishing local section $s_0$ of $\rho^*{L}$ near $(0,0)$. Then $f={\rho^*{s}}/{s_0}$ is a holomorphic function near $(0,0)$ and has a power series expansion
$$f(z_1,z_2)=\sum_{i_1,i_\geq 0} a_{i_1,i_2}z_1^{i_1}z_2^{i_2}.$$
Set
$$(\rho^*s)^\circ=a_{e_1,e_2}s_0(dz_1)^{\otimes e_1}(dz_2)^{\otimes e_2}.$$
It is compatible with the previous definition in the sense that there is a canonical isometry $(\rho^*\overline{L})^\circ=\rho^*(\overline{L}^\circ)$ such that $(\rho^*{s})^\circ=\rho^*({s}^\circ)$. So
\begin{align*}
-\log\|s^\circ(P)\|_{\overline{L}^\circ}
=&\ -\log\|\rho^*s^\circ(0,0)\|_{\rho^*\overline L^\circ}\\
=&\ -\log|a_{e_1,e_2}|-\log\|s_0(0,0)\|_{\rho^*\overline L}-e_1\log\|(dz_1)(0)\|_{\rho_1^*\bar\omega_a}-e_2\log\|(dz_2)(0)\|_{\rho_2^*\bar\omega_a}.
\end{align*}
By the first part of Theorem  \ref{propestimateArakelovofdiagonal_compareTwoMetrics part I},
$$
-\log\|(dz_i)(0)\|_{\rho_i^*\bar\omega_a}
\ge -\big|\log\|(dz_i)(0)\|_{\mathrm{hyp},\mathbb{D}}\big|-2\cdot10^{14}g^{\frac{14}{3}}\varphi(C(\mathbb C))
=-2\cdot10^{14}g^{\frac{14}{3}}\varphi(C(\mathbb C)).$$
Then we only need to prove that there is an $s_0$ satisfying
$$-\log|a_{e_1,e_2}|-\log\|s_0(0,0)\|_{\rho^*\overline L}\ge-\log \|s\|_{\overline L,\sup}-(2.1\cdot10^{14}g^{\frac{14}{3}}(2d_1+2d_2+4gd)+10^{13}g^{\frac{14}{3}}(e_1+e_2))\varphi(C(\mathbb{C})).$$

Assume that $s_0$ is a holomorphic section of $\rho^*L$ on the polydisc ${\mathbb{D}_r^2}$ for some $r\in (0,1)$ such that
$$M(s_0,r)=\sup_{{\mathbb{D}_r^2}} \left|\log\|s_0\|_{\rho^*\overline{L}}\right|$$
is finite. 
Then $s_0$ does not vanish at any point of $\mathbb{D}_r$. 
Thus $f={\rho^*s}/{s_0}$ is a holomorphic function on $\mathbb{D}_r^2$. By the Cauchy integral formula, for $0<t<r$,
$$a_{e_1,e_2}=\frac{1}{e_1!e_2!}\frac{\partial^{e_1+e_2}f(0,0)}{\partial z_1^{e_1}\partial z_2^{e_2}}=\int_{(\partial\mathbb{D}_t)^2}\frac{f(z_1,z_2)}{z_1^{e_1}z_2^{e_2}}d\mu.$$
Here $d\mu$ is the probability Haar measure on $(\partial\mathbb{D}_t)^2$. So
\begin{align*}
|a_{e_1,e_2}|&\le\sup_{(\partial\mathbb{D}_t)^2}\left|\frac{f(z_1,z_2)}{z_1^{e_1}z_2^{e_2}}\right|\\
&={t^{-e_1-e_2}}\sup_{(\partial\mathbb{D}_t)^2}\frac{\|\rho^*s\|_{\rho^*\overline L}}{\|s_0\|_{\rho^*\overline L}}\\
&\le{t^{-e_1-e_2}}e^{M(s_0,r)}\sup_{(\partial\mathbb{D}_t)^2}{\|\rho^*s\|_{\rho^*\overline L}}\\
&\leq {t^{-e_1-e_2}}e^{M(s_0,r)}{\sup_{C(\mathbb C)^2}\|s\|_{\overline L}}.
\end{align*}
We get
$$-\log|a_{e_1,e_2}|-\log\|s_0(0,0)\|_{\rho^*\overline L}\ge-\log
\|s\|_{\overline L,\sup}+(e_1+e_2)\log(t)-2M(s_0,r).$$
Set $t\to r$. Then
$$-\log|a_{e_1,e_2}|-\log\|s_0(0,0)\|_{\rho^*\overline L}\ge-\log\|s\|_{\overline L,\sup}+(e_1+e_2)\log(r)-2M(s_0,r).$$
If $(s_0,r)$ satisfies
$$\log(r)\ge-10^{13}g^{\frac{14}{3}}(e_1+e_2)\varphi(C(\mathbb{C}))$$ and 
$$M(s_0,r)\le 2.1\cdot10^{14}g^{\frac{14}{3}}(d_1+d_2+2gd)\varphi(C(\mathbb{C})),$$
 then we can finish the proof.
In the following, we prove that such a pair $(s_0,r)$ exists. 

By the second part of Theorem  \ref{propestimateArakelovofdiagonal_compareTwoMetrics part I}, there is a section $s_1$ of $\rho^*\mathcal{O}(\Delta)$ on $(\mathbb{D}_{\frac1{20}})^2$ such that
$$\sup_{{(\mathbb{D}_\frac{1}{20})^2}}\left| \log\|s_1\|_{\rho^*\overline{\mathcal{O}}(\Delta)_a}\right|\le2\cdot10^{14}g^{\frac{14}{3}}\varphi(C(\mathbb{C})).$$
Then
$$s_0=s_1^{\otimes(2g-2)d}(dz_1)^{\otimes d_1-d}(dz_2)^{\otimes d_2-d}$$
is a section of
$$\rho^*L=(2g-2)d\rho^*\mathcal{O}(\Delta)+(d_1-d)p_1^*\rho_1^*\omega+(d_2-d)p_2^*\rho_2^*\omega$$
on $(\mathbb{D}_{\frac1{20}})^2$. 
We will check that $(s_0, 1/20)$ satisfies the requirement. 

Start with the trivial bound
$$\left|\log\|s_0\|_{\rho^*\overline{L}}\right|
\le(2g-2)d\left|\log\|s_1\|_{\rho^*\overline{\mathcal{O}}(\Delta)_a}\right|
+(d_1+d)\big|\log\|dz_1\|_{\rho_1^*\bar{\omega}_a}\big|+(d_2+d)\big|\log\|dz_2\|_{\rho_2^*\bar{\omega}_a}\big|.$$
By Theorem  \ref{propestimateArakelovofdiagonal_compareTwoMetrics part I} again,
$$\sup_{{\mathbb{D}_\frac{1}{20}}}\left|\log\|dz_i\|_{\rho_i^*\bar{\omega}_a}\right|\le2\cdot10^{14}g^{\frac{14}{3}}\varphi(C(\mathbb{C}))+\sup_{{\mathbb{D}_\frac{1}{20}}}\big|\log\|dz_i\|_{{\mathrm{hyp},\mathbb{D}}}\big|
=2\cdot10^{14}g^{\frac{14}{3}}\varphi(C(\mathbb{C}))+\log(400/399).$$
This gives 
\begin{align*}
M(s_0, 1/20)=&\ \sup_{(\mathbb{D}_\frac{1}{20})^2} \left|\log\|s_0\|_{\rho^*\overline{L}}\right| \\
\le &\ (2g-2)d\cdot \left(2\cdot10^{14}g^{\frac{14}{3}}\varphi(C(\mathbb{C}))\right)
+(d_1+d_2+2d)\cdot 
\left(2\cdot10^{14}g^{\frac{14}{3}}\varphi(C(\mathbb{C}))+\log(400/399)\right)\\
< &\ (2g-2)d\cdot \left(2.1\cdot10^{14}g^{\frac{14}{3}}\varphi(C(\mathbb{C})\right)
+(d_1+d_2+2d)\cdot 
\left(2.1\cdot10^{14}g^{\frac{14}{3}}\varphi(C(\mathbb{C})) \right)\\
=&\ (d_1+d_2+2gd)\cdot 
\left(2.1\cdot10^{14}g^{\frac{14}{3}}\varphi(C(\mathbb{C})) \right).
\end{align*}
Here the last inequality follows from Theorem \ref{thmglobalestimatezhangphiinvariant part I}. 
On the other hand, 
Theorem \ref{thmglobalestimatezhangphiinvariant part I} also implies 
$$\log\left(\frac{1}{20}\right)\ge-10^{13}\cdot g^{\frac{14}{3}}\varphi(C(\mathbb{C})).$$
The finishes the proof of Proposition \ref{height_bound}(2).

\subsection{Completion of the proof}

Now we can finish the proof of Theorem \ref{vojta}. 
Let $P_1,P_2\in C(K)$ be distinct rational points as in the theorem. 
Denote $P=(P_1,P_2)\in C^2(K)$. 
Recall the adelic {line bundle}
$$\overline{L}=d_1p_1^*\bar\omega_a+d_2p_2^*\bar\omega_a+d((2g-2)\overline{\mathcal{O}}(\Delta)_a-p_1^*\bar\omega_a-p_2^*\bar\omega_a)$$
on $C^2$. 
Here $d$ is a large positive integer, and we further  take
$$d_1=\left\lfloor\sqrt{g+\frac{1}{10^5g}}\frac{\vert P_2\vert}{\vert P_1\vert}d\right\rfloor,\quad d_2=\left\lfloor\sqrt{g+\frac{1}{10^5g}}\frac{\vert P_1\vert}{\vert P_2\vert}d\right\rfloor.$$
Here $\lfloor x\rfloor$ is the largest integer which is less than or equal to $x$. 
For sufficiently large $d$, the inequalities among $(d,d_1,d_2)$ in Proposition \ref{small} is satisfied. Then the proposition implies that there exists a positive integer $n$ such that $nL$ has a section $s$ with
$$-\sum_{v\in M_K}\epsilon_v\log\|s\|_{v,\sup}\ge-\frac{ng(g+1)d_1d^2}{4(g-1)(d_1d_2-gd^2)}\bar\omega_a^2.$$
By Proposition \ref{height_bound},
\begin{align*}
h_{\overline L^\circ}(P)=&\ -\sum_{v\in M_K}\epsilon_v\log\|s^\circ(P)\|_{\overline L^\circ,v}\\
\ge&\ -\sum_{v|\infty}\left(\epsilon_v\log \|s\|_{v,\sup}+2.1\cdot10^{14}g^{\frac{14}{3}}(2nd_1+2nd_2+4gnd+e_1+e_2)\epsilon_v\varphi(C_{K_v})\right)\\
&-\sum_{v\nmid\infty}\left(\epsilon_v\log \|s\|_{v,\sup}+\frac{15}{4}(e_1+e_2)\varphi(C_{K_v})\log N_v\right)\\
\ge&\ -\sum_{v\in M_K}\epsilon_v\log\|s\|_{v,\sup}-2.1\cdot10^{14}g^{\frac{14}{3}}(2nd_1+2nd_2+4gnd+e_1+e_2)\varphi(C)\\
\ge&\ -\left(\frac{ng(g+1)d_1d^2}{4(g-1)(d_1d_2-gd^2)}+6.3\cdot10^{14}g^{\frac{14}{3}}(2nd_1+2nd_2+4gnd+e_1+e_2)\right)\bar{\omega}_a^2.
\end{align*}
Here the last inequality follows from Theorem \ref{phi_inequality}.

On the other hand, by definition
$$\overline L^\circ=\left(n \overline L+ e_1p_1^*\bar{\omega}_a+e_2p_2^*\bar{\omega}_a \right)|_P.$$
By Theorem \ref{admissible heights} and Lemma \ref{height identity}, we have 
\begin{align*}
h_{\overline L^\circ}(P)
=&\ h_{n \overline L}(P)+ e_1 h_{\bar{\omega}_a}(P_1)+ e_2 h_{\bar{\omega}_a}(P_2)\\
=&\ 
\frac{2g-2}{g}((nd_1+e_1)|P_1|^2+(nd_2+e_2)|P_2|^2)-(4g-4)nd\langle P_1,P_2\rangle+\frac{nd_1+nd_2-2gnd+e_1+e_2}{4g(g-1)}\bar\omega_a^2.
\end{align*}
Combine it with the above lower bound. We have
\begin{align*}
&\frac{2g-2}{g}\left((nd_1+e_1)|P_1|^2+(nd_2+e_2)|P_2|^2\right)-(4g-4)nd\langle P_1,P_2\rangle\\
\ge&\ - \left( \frac{ng(g+1)d_1d^2}{4(g-1)(d_1d_2-gd^2)}
+6.3\cdot10^{14}g^{\frac{14}{3}}(2nd_1+2nd_2+4gnd+e_1+e_2)
+\frac{nd_1+nd_2-2gnd+e_1+e_2}{4g(g-1)} \right)\bar\omega_a^2\\
\ge&\ - \left( \frac{ng(g+1)d_1d^2}{4(g-1)(d_1d_2-gd^2)}
+6.31\cdot10^{14}g^{\frac{14}{3}}(2nd_1+2nd_2+4gnd+e_1+e_2)
 \right)\bar\omega_a^2. 
\end{align*}
It follows that 
\begin{align*}
&\frac{2g-2}{g}\left((d_1+\frac{e_1}{n})|P_1|^2+(d_2+\frac{e_2}{n})|P_2|^2\right)-(4g-4)d\langle P_1,P_2\rangle\\
\ge&\ -\left(\frac{g(g+1)d_1d^2}{4(g-1)(d_1d_2-gd^2)}+6.31\cdot10^{14}g^{\frac{14}{3}}(2d_1+2d_2+4gd+\frac{e_1}{n}+\frac{e_2}{n}) \right)\bar{\omega}_a^2.
\end{align*}

\subsubsection*{Dyson's lemma}
As in Vojta's original proof, we use the following Dyson's lemma to bound $e_1$ and $e_2$.

\begin{thm}[Vojta \cite{Vojta_Dyson}]
\label{dyson}
Let $K$ be a field of characteristic $0$ and $C$ a curve of genus $g\geq 2$ over $K$. 
Let $s$ be a nonzero section of a line bundle $L$ on $C^2$. 
Assume $L$ has degree $m_1>0$ (resp. $m_2>0$) on the fiber of the projection map $p_{2}:C^2\to C$ (resp. $p_{1}:C^2\to C$). Denote a function
$$V(t)=\int_{0\le x\le 1,\ 0\le y\le 1,\ x+y\le t}dxdy, \quad  \ t\geq 0.$$
Then for any $P\in C^2(K)$ there is an admissible pair $(e_1,e_2)$ such that
$$V\left(\frac{e_1}{m_1}+\frac{e_2}{m_2}\right)\le\frac{L^2}{2m_1m_2}+\frac{m_2}{2m_1}(2g-1).$$
\end{thm}

Apply Dyson's lemma to our section $s$ of $n L$. 
In this case, $m_i=(2g-2)nd_i$. 
Then Dyson's lemma gives an admissible pair $(e_1,e_2)$  such that 
$$V\left(\frac{e_1}{(2g-2)nd_1}+\frac{e_2}{(2g-2)nd_2}\right)\le\frac{d_1d_2-gd^2}{d_1d_2}+\frac{d_2}{2d_1}(2g-1).$$
Recall  
$$d_1=\left\lfloor\sqrt{g+\frac{1}{10^5g}}\frac{\vert P_2\vert}{\vert P_1\vert}d\right\rfloor,\quad d_2=\left\lfloor\sqrt{g+\frac{1}{10^5g}}\frac{\vert P_1\vert}{\vert P_2\vert}d\right\rfloor.$$
We have 
$$
\lim_{d\to\infty} \left( \frac{d_1d_2-gd^2}{d_1d_2}+\frac{d_2}{2d_1}(2g-1) \right)
=\frac{1}{10^5g^2+1}+\frac{(2g-1)|P_1|^2}{2|P_2|^2}
< \frac{1}{10^5g^2+1}+\frac{2g-1}{2\cdot 10^{10}g^5}
< \frac{1}{10^5g^2}.
$$
Therefore, for sufficiently large $d$, we have 
$$V\left(\frac{e_1}{(2g-2)nd_1}+\frac{e_2}{(2g-2)nd_2}\right)< \frac{1}{10^5g^2}.$$
By definition, the function $V(t)$ is increasing and $\displaystyle V(t)=\frac{t^2}2$ for $0\le t\le1$. It follows that
$$\frac12\left(\frac{e_1}{(2g-2)nd_1}+\frac{e_2}{(2g-2)nd_2}\right)^2\le\frac{1}{10^5g^2},$$
and thus
$$\frac{e_1}{nd_1}+\frac{e_2}{nd_2}\le\frac{g-1}{(50\sqrt{5})g}.$$
As
$$
d_1\leq \sqrt{g+\frac{1}{10^5g}}\frac{\vert P_2\vert}{\vert P_1\vert}d,
\quad d_2\leq \sqrt{g+\frac{1}{10^5g}}\frac{\vert P_1\vert}{\vert P_2\vert}d,
$$
we further have
$$\frac{e_1}{n} \frac{\vert P_1\vert}{\vert P_2\vert}+\frac{e_2}{n} \frac{\vert P_2\vert}{\vert P_1\vert}
\le  \frac{g-1}{(50\sqrt{5})g}\cdot d\sqrt{g+\frac{1}{10^5g}}
\leq  \frac{1}{50\sqrt{5}}\cdot (d_1+1)\frac{\vert P_1\vert}{\vert P_2\vert}.$$
Here the second inequality holds for sufficiently large $d$. 
Therefore, we finally get
$$\frac{e_1|P_1|^2+e_2|P_2|^2}n\le\frac{(d_1+1)|P_1|^2}{50\sqrt5},$$
and
$$\frac{e_1+e_2}n\le\frac{d_1+1}{50\sqrt5}\le0.01(d_1+1).$$

Return to our height inequality right before Theorem \ref{dyson}. By the bounds for $(e_1,e_2)$, we further have
\begin{align*}
&\frac{2g-2}{g}\left(d_1\left(1+\frac{1}{50\sqrt5}\right)|P_1|^2+\frac{|P_1|^2}{50\sqrt5}+d_2|P_2|^2\right)-(4g-4)d\langle P_1,P_2\rangle\\
\ge&\ -\left(\frac{g(g+1)d_1d^2}{4(g-1)(d_1d_2-gd^2)}+6.31\cdot10^{14}g^{\frac{14}{3}}(2.01d_1+2d_2+4gd+0.01)\right)\bar{\omega}_a^2.
\end{align*}
Divide the inequality by $d$ and set $d\to\infty$. We have
\begin{align*}
&\frac{2g-2}{g}\left(2+\frac{1}{50\sqrt5}\right)\sqrt{g+\frac{1}{10^5g}}|P_1||P_2|-(4g-4)\langle P_1,P_2\rangle\\
\ge&\ -\left(\frac{10^5g^2(g+1)}{4(g-1)}\sqrt{g+\frac{1}{10^5g}}\frac{|P_2|}{|P_1|}+6.31\cdot10^{14}g^{\frac{14}{3}}\left(2.01\sqrt{g+\frac{1}{10^5g}}\frac{|P_2|}{|P_1|}+2\sqrt{g+\frac{1}{10^5g}}\frac{|P_1|}{|P_2|}+4g\right)\right)\bar{\omega}_a^2\\
\geq &\ -6.31\cdot10^{14}g^{\frac{14}{3}}\left(2.011\sqrt{g+\frac{1}{10^5g}}\frac{|P_2|}{|P_1|}+2\sqrt{g+\frac{1}{10^5g}}\frac{|P_1|}{|P_2|}+4g\right)\bar{\omega}_a^2\\
=&\ -6.31\cdot10^{14}g^{\frac{14}{3}}\left(2.011\sqrt{g+\frac{1}{10^5g}}
+2\sqrt{g+\frac{1}{10^5g}}\frac{|P_1|^2}{|P_2|^2}+4g\frac{|P_1|}{|P_2|} \right)
\cdot \frac{|P_2|}{|P_1|}\bar{\omega}_a^2.
\end{align*}
By the assumptions in Theorem \ref{vojta},  
$$
\frac{|P_1|}{|P_2|} \leq 10^{-5}g^{-\frac52}, \quad 
\bar{\omega}_a^2 \leq (1.44)^{-1}\cdot10^{-18}g^{-\frac{14}{3}} |P_1|^2.
$$
Then the above inequality implies
\begin{align*}
&\frac{2g-2}{g}\left(2+\frac{1}{50\sqrt5}\right)\sqrt{g+\frac{1}{10^5g}}|P_1||P_2|-(4g-4)\langle P_1,P_2\rangle\\
\geq&\ -6.31\cdot10^{14}g^{\frac{14}{3}} \cdot 
\left(2.012\sqrt{g+\frac{1}{10^5g}} \right)  \cdot (1.44)^{-1}\cdot10^{-18}g^{-\frac{14}{3}}\cdot |P_1|\cdot |P_2| \\
\ge&\ -\left(8.82\cdot 10^{-4} \sqrt{g+\frac{1}{10^5g}} \right)|P_1|\cdot |P_2|.
\end{align*}
This further implies
\begin{align*}
 \frac{\langle P_1,P_2\rangle}{|P_1||P_2|} 
 \le&\ \frac{1}{4g-4} 
\left( \frac{2g-2}{g} \big(2+\frac{1}{50\sqrt5} \big)\sqrt{g+\frac{1}{10^5g}}+
8.82\cdot 10^{-4} \sqrt{g+\frac{1}{10^5g}}
 \right)\\
  \le&\ \frac{1}{4g-4} 
\cdot \frac{2g-2}{g}  \left(2+\frac{1}{50\sqrt5}   +
8.82\cdot 10^{-4}  \right) \sqrt{g+\frac{1}{10^5g}} \\
  \le&\  \frac{1.00492}{g}     \sqrt{g+\frac{1}{10^5g}}.
\end{align*}
Note that 
$$
\sqrt{g+\frac{1}{10^5g}}= \sqrt{g} \sqrt{1+\frac{1}{10^5g^2}}
\leq \sqrt{g} \left(1+\frac{1}{2\cdot 10^5g^2}\right)
\leq   (1+10^{-5} )\sqrt{g}. 
$$
We have
$$\cos\angle(P_1,P_2)
=\frac{\langle P_1,P_2\rangle}{|P_1||P_2|}
\le 
\frac{1.00492}{g}\cdot (1+10^{-5} )\sqrt{g}
\le\sqrt{\frac{1.01}{g}}.$$
This finishes the proof of Theorem \ref{vojta}.

\subsection{Functions fields}

As a geometric analogue of Theorem \ref{vojta}, we have the following theorem over function fields of characteristic $0$. 
It has better constants than the arithmetic case due to absence of archimedean places. 
On the other hand, comparing it with Yu \cite[Theorem 3.1]{Yu}, the current theorem gives more accurate results by assuming stronger conditions. 

\begin{thm}
\label{vojta_function_field}
Let $K$ be a function field of one variable over a field $k$ of characteristic $0$. 
Let $C$ be a curve of genus $g\geq 2$ over $K$. 
Let $P_1,P_2\in C(K)$ be distinct rational points. Assume
$$|P_1|\ge10^4g\sqrt{\bar{\omega}_a^2}$$
and
$${|P_2|}{}\ge10^5g^{\frac52}|P_1|>0.$$
Then
$$\cos\angle(P_1,P_2)\le\sqrt{\frac{1.01}{g}}.$$
\end{thm}

\begin{proof}
We sketch the proof here. 
Still consider the adelic {line bundle}
$$\overline{L}=d_1p_1^*\bar\omega_a+d_2p_2^*\bar\omega_a+d((2g-2)\overline{\mathcal{O}}(\Delta)_a-p_1^*\bar\omega_a-p_2^*\bar\omega_a)$$
on $C^2$. 
Here $d$ is a large positive integer, and we still take
$$d_1=\left\lceil\sqrt{g+\frac{1}{10^5g}}\frac{\vert P_2\vert}{\vert P_1\vert}d\right\rceil,\quad d_2=\left\lceil\sqrt{g+\frac{1}{10^5g}}\frac{\vert P_1\vert}{\vert P_2\vert}d\right\rceil.$$
We still have a small section $s$ of $nL$ as in Proposition \ref{small}. 
This still gives a section $s^\circ$ of 
$$\overline L^\circ=(n \overline L+e_1p_1^*\bar\omega_a+e_2p_2^*\bar\omega_a)|_P.$$

Since all places are non-archimedean, 
we have
\begin{align*}
h_{\overline L^\circ}(P)=&\ -\sum_{v\in M_K}\log\|s^\circ(P)\|_{v}\\
\ge&\ -\sum_{v\in M_K}\log \|s\|_{n\overline L,v, \sup}-\frac{15}{4}(e_1+e_2)\varphi(C)\\
\ge&\ -\frac{ng(g+1)d_1d^2}{4(g-1)(d_1d_2-gd^2)}\bar\omega_a^2
-\frac{75}{8}(e_1+e_2)\bar{\omega}_a^2.
\end{align*}
It follows that
\begin{align*}
&\frac{2g-2}{g}((nd_1+e_1)|P_1|^2+(nd_2+e_2)|P_2|^2)-(4g-4)nd\langle P_1,P_2\rangle\\
\ge&\ -\frac{ng(g+1)d_1d^2}{4(g-1)(d_1d_2-gd^2)}\bar\omega_a^2-\left(\frac{nd_1+nd_2+e_1+e_2}{4g(g-1)}+\frac{75}{8}(e_1+e_2)\right)\bar\omega_a^2.
\end{align*}
We have still assumed $|P_2|\ge10^5g^\frac{5}{2}|P_1|$. 
Apply Dyson's lemma.
For sufficiently large $d$, we still have
$$\frac{e_1|P_1|^2+e_2|P_2|^2}n\le\frac{d_1|P_1|^2}{50\sqrt5},$$
and
$$\frac{e_1+e_2}n\le\frac{d_1}{50\sqrt5}\le0.01d_1.$$
It follows that 
\begin{align*}
&\frac{2g-2}{g}\left(2+\frac{1}{50\sqrt5}\right)\sqrt{g+\frac{1}{10^5g}}|P_1||P_2|-(4g-4)\langle P_1,P_2\rangle\\
\ge&\ -\left(\frac{10^5g^2(g+1)}{4(g-1)}\sqrt{g+\frac{1}{10^5g}}\frac{|P_2|}{|P_1|}+\left(\frac{1.01}{4g(g-1)}+\frac{3}{32}\right)\sqrt{g+\frac{1}{10^5g}}\frac{|P_2|}{|P_1|}+\frac{1}{4g(g-1)}\sqrt{g+\frac{1}{10^5g}}\frac{|P_1|}{|P_2|}\right)\bar{\omega}_a^2\\
\ge&\ -\left(\frac{10^5g^2(g+1)}{4(g-1)}+\left(\frac{1.01}{4g(g-1)}+\frac{3}{32}\right)+\frac{10^{-10} g^{-5}}{4g(g-1)}
 \right)\sqrt{g+\frac{1}{10^5g}} \frac{|P_2|}{|P_1|} \bar{\omega}_a^2\\
 \ge&\ -\left(\frac{10^5g^2(g+1)}{4(g-1)}+0.23
 \right)\sqrt{g+\frac{1}{10^5g}} \frac{|P_2|}{|P_1|} \bar{\omega}_a^2\\
 \ge&\ -\left(\frac{10^5g^2(g+1)}{4(g-1)}+0.23
 \right)\sqrt{g+\frac{1}{10^5g}}\cdot 10^{-8} g^{-2} |P_1|\cdot |P_2|\\
  \ge&\ -6\cdot 10^{-4}\sqrt{g+\frac{1}{10^5g}}\cdot   |P_1|\cdot |P_2| .
\end{align*}
Here the second inequality has used $|P_2|\ge10^5g^\frac{5}{2}|P_1|$, 
and the fourth inequality has used $|P_1|\ge10^4g\sqrt{\bar\omega_a^2}$. 
The remaining part of the proof is the same as that of Theorem \ref{vojta}. 
\end{proof}

We also have the following geometric analogue of the quantitative Mumford inequality in Theorem \ref{mumford}. 

\begin{thm}
\label{mumford_function_field}
Let $K$ be a function field of one variable over a field $k$ of characteristic $0$. 
Let $C$ be a curve of genus $g\geq 2$ over a function field $K$ of characteristic $0$. Let $P_1,P_2\in C(K)$ be distinct rational points. Assume
$$|P_1|\ge10^{4}g^{\frac{7}{3}}\sqrt{\bar\omega_a^2}.$$
and
$$|P_1|\le{|P_2|}{}\le1.15|P_1|.$$
Then
$$\cos\angle(P_1,P_2)\le{\frac{1.01}{g}}.$$
\end{thm}

\begin{proof}
The proof is similar to Theorem \ref{mumford}. 
In fact, Theorem \ref{thmfaltingselkieszhangphiinv non-archi} and Theorem \ref{phi_inequality} imply 
\begin{align*}
h_{\overline{\mathcal O}(\Delta)_a}(P)=\sum_{v\in M_K} 
G_{a,v}(P_1,P_2)
\ge
-\frac{15}{4}\sum_{v\in M_K}\varphi(C_{K_v})
\ge -\frac{45}{4}\, \bar\omega_a^2.
\end{align*}
By Theorem \ref{admissible heights}, this gives
$$
\frac{1}{g}(|P_1|^2+|P_2|^2)-2\langle P_1,P_2\rangle
\geq -\frac{45}{4}\, \bar\omega_a^2.$$
By  $|P_1|\ge10^4g\sqrt{\bar{\omega}_a^2}$, we have
$$\frac{1+10^{-7}}{g}|P_1|^2+\frac{1}{g}|P_2|^2-2\langle P_1,P_2\rangle\ge0.$$
The remaining part of the proof is the same as that of Theorem \ref{mumford}. 
\end{proof}

\section{Quantitative Bogomolov conjecture}

Recall that the norm $|Q|=\sqrt{\hat h(Q)_K}$ and the angle $\angle(P_1,P_2)$ for algebraic points $Q, P_1, P_2$ on a curve $C$ over a number field $K$ are defined in terms of the N\'eron--Tate heights as in \S\ref{sec height}. 
The goal of this section is to prove the following quantitative 
 Bogomolov conjecture. 

\begin{thm}[Theorem \ref{bogomolov_introduction}]
\label{bogomolov}
Let $K$ be a number field, $C$ a curve of genus $g\geq 2$ over $K$, and $J$ the Jacobian variety of $C$. 
\begin{enumerate}
    \item For any $0\le r<\sqrt\frac1{8g}$ and $x\in J(\overline K)_\mathbb R$, 
    $$\#\{P\in C(\overline K):|P-x|\le r\sqrt{\bar\omega_a^2}\}<\frac{3.2\cdot 10^{11}g^{\frac{17}3}}{1-8gr^2}\left(1+10^{-6}\log\left(\frac{1}{1-8gr^2}\right)\right).$$
    \item Let $\kappa>1$ and $\theta\in(0,\pi/2)$ be real numbers satisfying
    $$g\cos^2\theta>\frac{(\kappa+1)^2}{4\kappa}+\frac{1}{3.2\cdot 10^{11} g^{\frac{11}{3}}}.$$
    Then for any nonzero $x\in J(\overline K)_\mathbb R$,
    $$\#\{P\in C(\overline K):|x|\le|P|\le\kappa|x|,\, \angle(P,x)\le\theta\}
    \leq 
    3.2\cdot 10^{11}g^{\frac{17}3}.$$
\end{enumerate}
\end{thm}

As mentioned in the introduction, part (1) of the theorem with inexplicit constants is the uniform Bogomolov conjecture proved by  Dimitrov--Gao--Habegger \cite{DGH} and K\"uhne \cite{Kuhne}, and by Yuan \cite{Yuan_Bogomolov} using a different method.
Part (2) of the theorem was inspired by the extra term in the main theorem of \cite{Yuan_Bogomolov} .

A geometric version of part (1) was originally proved by Looper--Silverman--Wilms \cite{LSW}. Our current arithmetic version follows the framework of \cite{LSW} by applying lower bound of the Faltings--Elkies invariant by the $\varphi$-invariant
 from Part II of this paper.

Recall that the Manin--Mumford conjecture was proved by Raynaud \cite{Ray83a,Ray83b}. 
Taking $r=0$ in Theorem \ref{bogomolov}, we obtain 
the following quantitative version of the Manin--Mumford conjecture.

\begin{thm}[quantitative Manin--Mumford]
\label{manin_mumford}
Let $K$ be a field of characteristic $0$, $C$ a curve of genus $g\geq 2$ over $K$ and $J$ the Jacobian variety of $C$. For any line bundle $\alpha$ of degree $1$ on $C$,
$$\#\big((C(K)-\alpha)\cap J(K)_\mathrm{tor}\big)\le 3.2\cdot 10^{11}g^{\frac{17}3}.$$
\end{thm}
\begin{proof}
By a reduction argument as in \cite[Lemma 3.1]{DGH2},we can assume that $K$ is a number field.
Then the theorem is just the case $r=0$ of Theorem \ref{bogomolov}. 
\end{proof}

\subsection{The first part}

The idea to prove Theorem \ref{bogomolov} is to obtain a suitable lower bound of $\sum_{1\leq i>j\leq n} |P_i-P_j|^2$ to bound the number $n$. The following quadratic identity is extracted from the proof of Looper--Silverman--Wilms \cite{LSW}.

\begin{lem}[quadratic identity]
\label{quadratic}
For $P_1,P_2,\dots,P_n\in C(\overline K)$,
$$\frac{4(n+g-1)(g-1)}{n-1}\sum_{1\le i<j\le n}|P_i-P_j|^2=\frac{4ng(g-1)}{n-1}\sum_{1\le i<j\le n}h_{\overline{\mathcal{O}}(\Delta)_a}(P_i,P_j)+{\frac{n^2}2\bar\omega_a^2}+\left|(2g-2)\sum_{i=1}^nP_i-{n\omega}\right|^2.$$
\end{lem}

\begin{proof}
There is a finite extension $K'/K$ such that $P_1,P_2,\dots,P_n\in C(K')$. Taking the base change of $C$ to $K'$ multiplies each term in the target equality by $[K':K]$. So we may assume $P_1,P_2,\dots,P_n\in C(K)$.

By the arithmetic Hodge index theorem (Theorem \ref{hodge index}) and the arithmetic adjunction formula (Theorem \ref{adjunction}), we have
\begin{align*}
2|P_i-P_j|^2=&\ -(\overline{\mathcal{O}}(P_i)_a-\overline{\mathcal{O}}(P_j)_a)^2\\
=&\ -\overline{\mathcal{O}}(P_i)_a^2-\overline{\mathcal{O}}(P_j)_a^2+2\overline{\mathcal{O}}(P_i)_a\cdot\overline{\mathcal{O}}(P_j)_a\\
=&\ \bar\omega_a\cdot\overline{\mathcal{O}}(P_i)_a+\bar\omega_a\cdot\overline{\mathcal{O}}(P_j)_a+2h_{\overline{\mathcal{O}}(\Delta)_a}(P_i,P_j).
\end{align*}
Here we also used the first equality of Theorem \ref{admissible heights}(2). 
This gives
$$\sum_{1\le i<j\le n}|P_i-P_j|^2=\frac{n-1}2\sum_{i=1}^n\bar\omega_a\cdot\overline{\mathcal{O}}(P_i)_a+\sum_{1\le i<j\le n}h_{\overline{\mathcal{O}}(\Delta)_a}(P_i,P_j).$$
Similarly,
\begin{align*}
2\left|(2g-2)\sum_{i=1}^nP_i-{n\omega}\right|^2=&\ -\left((2g-2)\sum_{i=1}^n\overline{\mathcal{O}}(P_i)_a-{n\bar\omega_a}\right)^2\\
=&\ (2g-2)(2n+2g-2)\sum_{i=1}^n\bar\omega_a\cdot\overline{\mathcal{O}}(P_i)_a-2(2g-2)^2\sum_{1\le i<j\le n}h_{\overline{\mathcal{O}}(\Delta)_a}(P_i,P_j)-n^2\bar\omega_a^2.
\end{align*}
Canceling the term $\sum_i\bar\omega_a\cdot\overline{\mathcal{O}}(P_i)_a$ gives the equality. 
\end{proof}

Now we are ready to prove the first part of Theorem \ref{bogomolov}. 
In the following, the idea to apply the generalized parallelogram rule to get better estimates is suggested to the authors by Pascal Autissier. 

\begin{proof}[Proof of Theorem \ref{bogomolov}(1)]
Let $P_1,P_2\dots, P_n\in C(\overline K)$ be distinct points satisfying  
$|P_i-x|^2\le r^2{\bar\omega_a^2}$. 
This easily implies 
$$
|P_i-P_j|^2\le (|P_i-x|+|P_j-x|)^2\le 4r^2{\bar\omega_a^2}
$$ 
and thus 
$$\sum_{1\le i<j\le n}|P_i-P_j|^2\leq 2n(n-1)r^2{\bar\omega_a^2}.$$

In fact, the generalized parallelogram rule 
$$\sum_{1\le i<j\le n}|P_i-P_j|^2=n\sum_{i=1}^n|P_i-x|^2-\left|\sum_{i=1}^nP_i-nx\right|^2$$
implies the stronger bound
$$\sum_{1\le i<j\le n}|P_i-P_j|^2\le n^2r^2\bar\omega_a^2.$$

By Lemma \ref{quadratic}, and omitting the non-negative term
$$\left|(2g-2)\sum_{i=1}^nP_i-{n\omega}\right|^2\ge0,$$
we get a quadratic inequality
$$\left(\frac{4n^2(n+g-1)(g-1)}{n-1}r^2-\frac{n^2}2\right)\bar{\omega}_a^2\ge\frac{4ng(g-1)}{n-1}\sum_{1\le i<j\le n}h_{\overline{\mathcal{O}}(\Delta)_a}(P_i,P_j).$$

We need to prove 
    $$n<\frac{3.2\cdot 10^{11}g^{\frac{17}3}}{1-8gr^2}\left(1+10^{-6}\log\left(\frac{1}{1-8gr^2}\right)\right).$$ 
If $n< g^2$, we have nothing to prove. 
Assume that $n\geq g^2$ in the following, which gives
$$\frac{(n+g-1)(g-1)}{n-1}< g.$$ 
Then the quadratic inequality implies
$$\frac{(8gr^2-1)n^2}2\bar{\omega}_a^2\geq\frac{4ng(g-1)}{n-1}\sum_{1\le i<j\le n}h_{\overline{\mathcal{O}}(\Delta)_a}(P_i,P_j).$$

By the bounds of the global sparsity in Theorem \ref{sparsity},
$$\sum_{1\le i<j\le n}h_{\overline{\mathcal{O}}(\Delta)_a}(P_i,P_j)\ge-\left(1.2\cdot 10^4 g^{\frac{1}{3}} n\log n+ 3.96\cdot 10^{10} g^{\frac{11}{3}} n \right)\bar{\omega}_a^2.$$
It follows that
$$\frac{(8gr^2-1)n^2}2\bar{\omega}_a^2\ge-\frac{4ng(g-1)}{n-1}
\left(1.2\cdot 10^4 g^{\frac{1}{3}} n\log n+ 3.96\cdot 10^{10} g^{\frac{11}{3}} n \right)\bar{\omega}_a^2.$$
Since $\bar{\omega}_a^2>0$, we get
$$1-8gr^2\le\frac{8g(g-1)}{n-1}\left(1.2\cdot 10^4 g^{\frac{1}{3}} \log n+ 3.96\cdot 10^{10} g^{\frac{11}{3}}  \right)< 10^5g^\frac{7}{3}\frac{\log n}{n-1}+3.17\cdot10^{11}g^\frac{17}{3}\frac{1}{n-1}.$$

Denote $A=(1-8gr^2)^{-1}$. Then we have 
$$ 10^5g^\frac{7}{3}A\frac{\log n}{n-1}+3.17\cdot10^{11}g^\frac{17}{3}A\frac{1}{n-1}>1.$$
We need to prove that this implies 
$$n<3.2\cdot 10^{11}g^{\frac{17}3} A  (1+10^{-6}\log A).$$ 
This is elementary, but we give a proof to check the constants. 
In fact, for $n>4$, we have 
$$\frac{1}{n-1}<\frac{\log n}{n-1}< \frac{1}{\sqrt{n}}.$$
This gives
$$
1 < 
 10^5g^\frac{7}{3}A\frac{1}{\sqrt {n}}+3.17\cdot10^{11}g^\frac{17}{3}A \frac{1}{\sqrt {n}}
< 3.18\cdot10^{11}g^\frac{17}{3}A\frac{1}{\sqrt n}.
$$
Then
$$
{\sqrt{n}} <   3.18\cdot10^{11}g^\frac{17}{3} A
$$
and 
$$
\log n < 2 \log (3.18\cdot10^{11}g^\frac{17}{3} A).
$$
Putting back to the original inequality, we have 
\begin{align*}
10^5g^\frac{7}{3}\cdot {2 \log (3.18\cdot10^{11}g^\frac{17}{3} A)}\cdot A
\frac{1}{n-1}+3.17\cdot10^{11}g^\frac{17}{3}A\frac{1}{n-1}> 1.
\end{align*}
It follows that 
\begin{align*}
n-1
< &\ 
10^5g^\frac{7}{3}\cdot 2 \log (3.18\cdot10^{11}g^\frac{17}{3} A) \cdot A+3.17\cdot10^{11}g^\frac{17}{3} A\\
< &\ 
10^5g^\frac{7}{3}\cdot 2 A \log A+3.18\cdot10^{11}g^\frac{17}{3} A\\
< &\ 
3.18\cdot10^{11}g^\frac{17}{3} A(1+10^{-6}\log A). 
\end{align*}
This finishes the proof of Theorem \ref{bogomolov}(1).
\end{proof}

\subsection{The second part}

To prove the second part of Theorem \ref{bogomolov}, we need the following inverse version of the classical Cauchy--Schwarz inequality.

\begin{lem}
\label{lemma_reversecauchyschwarz}
For $\kappa>1$ and $a_1,a_2,\dots,a_n\in[1,\kappa]$,
$$\sum_{i=1}^na_i^2\le\frac{(\kappa+1)^2}{4\kappa n}\left(\sum_{i=1}^na_i\right)^2.$$
\end{lem}

\begin{proof}
For $i=1,2,\dots,n$,
$$(a_i-1)(a_i-\kappa)\le0.$$
Their summation gives
$$\sum_{i=1}^na_i^2-(\kappa+1)\sum_{i=1}^na_i+\kappa n\le0.$$
So
$$(\kappa+1)^2\left(\sum_{i=1}^na_i\right)^2\ge\left(\sum_{i=1}^na_i^2+\kappa n\right)^2\ge4\kappa n\sum_{i=1}^na_i^2.$$
\end{proof}

\begin{proof}[Proof of Theorem \ref{bogomolov}(2)]    
The proof of still based on the quadratic identity. 
Let $P_1,P_2\dots, P_n\in C(\overline K)$ be distinct points satisfying   
$$|x|\le|P_i|\le\kappa|x|, \quad \angle(P_i,x)\le\theta.$$
Before applying the quadratic identity, we claim that 
$$\left|(2g-2)\sum_{i=1}^nP_i-{n\omega}\right|^2
\geq 
\frac{(2g-2)^2\cos^2\theta}{\gamma-\cos^2\theta}
\sum_{1\le i<j\le n}|P_i-P_j|^2.$$
Here we denote $\displaystyle\gamma=\frac{(\kappa+1)^2}{4\kappa}$ for simplicity. 

We first prove the claim, which is also written as 
$$(\gamma-\cos^2\theta)\left|\sum_{i=1}^nP_i\right|^2
\geq 
(\cos^2\theta)
\sum_{1\le i<j\le n}|P_i-P_j|^2.$$
By the generalized parallelogram law
$$\sum_{1\le i<j\le n}|P_i-P_j|^2=n\sum_{i=1}^n|P_i|^2-\left|\sum_{i=1}^nP_i\right|^2,$$
 it suffices to prove 
$$ \gamma \left|\sum_{i=1}^nP_i\right|^2
\geq 
n(\cos^2\theta)\sum_{i=1}^n|P_i|^2.$$
We  have 
$$\left|\sum_{i=1}^nP_i\right|^2\\
\ge \left(\sum_{i=1}^n|P_i|\cos\angle(P_i,x)\right)^2\\
\ge (\cos^2\theta)\left(\sum_{i=1}^n|P_i|\right)^2.
$$
Here the first inequality holds as the projection of $P_i$ to the direction of $x$ is 
$|P_i|\cos\angle(P_i,x)$ times $x/|x|$.  
Then the claim is reduced to 
$${\gamma}\left(\sum_{i=1}^n|P_i|\right)^2\geq n\sum_{i=1}^n|P_i|^2.$$
This follows from Lemma \ref{lemma_reversecauchyschwarz} by $|x|\le|P_i|\le\kappa|x|$. 
The claim is proved. 

To prove
Theorem \ref{bogomolov}(2), we are going to have two cases. If 
$$
\frac{(2g-2)^2\cos^2\theta}{\gamma-\cos^2\theta}
\geq \frac{4(n+g-1)(g-1)}{n-1},$$
then the claim implies
$$\left|(2g-2)\sum_{i=1}^nP_i-{n\omega}\right|^2\ge\frac{4(n+g-1)(g-1)}{n-1}\sum_{1\le i<j\le n}|P_i-P_j|^2.$$
In this case,  Lemma \ref{quadratic} implies
$$-\frac{n^2}2\bar{\omega}_a^2\ge\frac{4ng(g-1)}{n-1}\sum_{1\le i<j\le n}h_{\overline{\mathcal{O}}(\Delta)_a}(P_i,P_j),$$
Then the proof of Theorem \ref{bogomolov}(1) for $r=0$ gives exactly
$$n\le3.2\times10^{11}g^{\frac{17}{3}}.$$

On the other hand, if 
$$
\frac{(2g-2)^2\cos^2\theta}{\gamma-\cos^2\theta}
< \frac{4(n+g-1)(g-1)}{n-1},$$
then 
$$
\frac{(g-1)\cos^2\theta}{\gamma-\cos^2\theta}
< \frac{g}{n-1}+1,$$
and thus 
$$n< \frac{g(\gamma-\cos^2\theta)}{g\cos^2\theta-\gamma}+1.$$
The assumption
$$g\cos^2\theta>\gamma
    +(3.2\cdot 10^{11} g^{\frac{11}{3}})^{-1}$$
further gives
$$n< 
\frac{g(g\cos^2\theta-\cos^2\theta)}{(3.2\cdot 10^{11} g^{\frac{11}{3}})^{-1}}+1
< g(g-1)\cdot 3.2\cdot 10^{11} g^{\frac{11}{3}}+1
< 3.2\cdot 10^{11} g^{\frac{17}{3}}.$$
This proves Theorem \ref{bogomolov}(2).
\end{proof}

\begin{remark}\label{rmk extra term}
We can also prove the following explicit result in the form of \cite[Theorem 1.1]{Yuan_Bogomolov}. For any line bundle $\alpha$ of degree $1$ on $C$,
$$\#\left\{P\in C(\overline K):\hat{h}(P-\alpha)\le\frac{1}{32g}\bar\omega_a^2+\frac{g-2}{2g-2}\hat{h}\left(\alpha-\frac{\omega}{2g-2}\right)\right\}< 6.5\cdot 10^{11}g^{\frac{17}3}.$$
Since it is not used in the sequel, we omit the details.
\end{remark}

\subsection{Function fields}

Now we present the following analogue of Theorem \ref{bogomolov} over  function fields. Part (1) of the theorem is essentially due to Looper--Silverman--Wilms (cf. \cite[Theorem 1.3]{LSW}). Note that the canonical height of the loc. cit. is twice of ours, so our statement is slightly stronger than that of the loc. cit in that sense. This improvement is due to the application of the generalized parallelogram rule suggested by Autissier as in the proof of the number field case. Part (2) can be proved in the same way as number fields.

\begin{thm}
\label{geometric_bogomolov}
Let $K$ be a function field of one variable over a field $k$, $C$ be a non-isotrivial curve of genus $g\geq 2$ over $K$, and $J$ be the Jacobian variety of $C$. 
\begin{enumerate}
    \item For $0\le r<\sqrt\frac1{8g}$, and $x\in J(\overline K)_\mathbb R$, 
    $$\#\{P\in C(\overline K):|P-x|\le r\sqrt{\bar\omega_a^2}\}\le\frac{1}{1-8gr^2}\frac{16g^4+36g^2-26g-2}{(g-1)^2}+1.$$
    \item Let $\kappa>1$ and $\theta\in(0,\pi/2)$. Assume
    $$g\cos^2\theta>\frac{(\kappa+1)^2}{4\kappa}+\frac{1}{16}.$$
    Then for any $x\in J(\overline K)_\mathbb R$,
    $$\#\{P\in C(\overline K):|x|\le|P|\le\kappa|x|,\angle(P,x)\le\theta\}\le  16g^2+32g+124.$$
\end{enumerate}
\end{thm}

\section{Counting points}

The goal of this section is to prove the following main theorem of this paper.

\begin{thm}[Theorem \ref{thm_quantitativeMordellsubgroup}]
\label{mordell_lang}
Let $K$ be a field of characteristic $0$, $C$ a curve of genus $g\geq 2$ over $K$, and $J$ the Jacobian variety of $C$. For any line bundle $\alpha$  of degree $1$ on $C$ and any subgroup $\Lambda\subset J(K)$ of finite rank,
$$\#((C(K)-\alpha)\cap\Lambda)\le 10^{13}g^8\cdot\min\left\{1+\frac{5}{4\sqrt g},1+\frac{3\log g}{g}\right\}^{\mathrm{rk}(\Lambda)}.$$
\end{thm}

As we mentioned in the introduction, this gives explicit constants in the uniform Mordell conjecture proved by Vojta \cite{Vojta_Ar}, Dimitrov--Gao--Habegger \cite{DGH}, and K\"uhne \cite{Kuhne}. 
It is not surprising that our quantitative versions of the Vojta inequality, the Mumford inequality, and the uniform Mordell conjecture imply the current theorem with explicit constants. Nonetheless, we still need to make lots of efforts to decrease the Vojta constant to the current form.

\subsection{Different ranges of heights}

Now we sketch our idea to prove Theorem \ref{mordell_lang}. 
By a reduction argument as in \cite[Lemma 3.1]{DGH2}, we may assume that $\Lambda$ is finitely generated and $K$ is a number field. 
The case $\mathrm{rk}(\Lambda)=0$ is just the quantitative Manin--Mumford theorem in  Theorem \ref{manin_mumford}. 
We assume $\mathrm{rk}(\Lambda)\ge1$ in the following.

After a finite field extension of $K$, we may assume that there is a line bundle $\alpha_0$ of degree $1$ on $C$ with an isomorphism $(2g-2)\alpha_0\cong\omega$. Then we only need to prove the stronger bound
$$\#((C(K)-\alpha)\cap\Lambda)\le 10^{13}g^8\cdot
\min\left\{1+\frac{5}{4\sqrt g},1+\frac{3\log g}{g}\right\}^{\mathrm{rk}(\Lambda)-1},$$
under the additional assumption that $\alpha=\alpha_0$. Indeed, let $\Lambda'$ be the subgroup of $J(K)$ generated by $\Lambda$ and $\alpha-\alpha_0$. Then $\mathrm{rk}(\Lambda')\le\mathrm{rk}(\Lambda)+1$ and
$$(C(K)-\alpha)\cap\Lambda\ \subset\ (C(K)-\alpha)\cap\Lambda'=(C(K)-\alpha_0)\cap\Lambda'.$$
So it suffices to prove
$$\#((C(K)-\alpha_0)\cap\Lambda')\le 10^{13}g^8\cdot\min\left\{1+\frac{5}{4\sqrt g},1+\frac{3\log g}{g}\right\}^{\mathrm{rk}(\Lambda')-1}.$$

Denote $n=\mathrm{rk}(\Lambda)$. We abbreviate $C-\alpha_0$ as $C$, which viewed as a subvariety of $J$. By abuse of notation, we further abbreviate $C(K)\cap\Lambda$ as $C\cap\Lambda$.
We divide the set $C\cap\Lambda$ into the following three parts:
\begin{align*}
(C\cap\Lambda)_{\mathrm{small}}=&\ \left\{P\in C\cap\Lambda:|P|\le\sqrt\frac{\bar\omega_a^2}{16g}\right\},  \\
(C\cap\Lambda)_{\mathrm{medium}}=&\ \left\{P\in C\cap\Lambda:\sqrt\frac{\bar\omega_a^2}{16g}<|P|\le1.2\cdot10^9g^{\frac{7}{3}}\sqrt{\bar\omega_a^2}\right\},  \\
(C\cap\Lambda)_{\mathrm{large}}=&\ \left\{P\in C\cap\Lambda:|P|>1.2\cdot10^9g^{\frac{7}{3}}\sqrt{\bar\omega_a^2}\right\}.
\end{align*}
Our bounds for them are 
\begin{align*}
\#(C\cap\Lambda)_{\mathrm{small}} \leq &\ 6.5\cdot 10^{11}g^\frac{17}{3},  \\
\#(C\cap\Lambda)_{\mathrm{medium}}
\leq&\ 1.64\cdot10^{13}g^7\left(1+\frac{5}{2\sqrt2g}\right)^{n-1},  \\
\#(C\cap\Lambda)_{\mathrm{large}}
\leq &\  3.4\cdot10^6g^2\left(1+\frac{5}{4\sqrt{g}}\right)^{n-1}.
\end{align*}
Their sum gives
$$\#(C\cap\Lambda)\le2\cdot10^{13}g^7\left(1+\frac{5}{4\sqrt g}\right)^{n-1}\le 10^{13}g^8\left(1+\frac{5}{4\sqrt g}\right)^{n-1}.$$

Moreover, for $g\ge142$ and $n\ge1$, we have a refined bound 
$$\#(C\cap\Lambda)_\mathrm{large}\le2.4\cdot10^4g^8\left(1+\frac{3\log g}{g}\right)^{n-1}.$$
Replacing this one in the above bound of large points, we have for $g\ge142$, 
$$\#(C\cap\Lambda)\le 10^{13}g^8\left(1+\frac{3\log g}{g}\right)^{n-1}.$$
Note that for integers $g\geq2$, 
$$
1+\frac{3\log g}{g}\leq 1+\frac{5}{4\sqrt g} \quad \Longleftrightarrow \quad
g\geq 142.
$$
Therefore, these two cases imply that for $g\geq2$, 
$$\#(C\cap\Lambda)\le 10^{13}g^8\cdot
\min\left\{1+\frac{5}{4\sqrt g},1+\frac{3\log g}{g}\right\}^{n-1}.$$

Then the task is to prove the bounds for the subsets of different range of heights.  
Note that the bound for $\#(C\cap\Lambda)_{\mathrm{small}}$ is a consequence of Theorem \ref{bogomolov}(1) with $r=\sqrt{\frac{1}{16g}}$.
In the following, we are going to prove the bound for medium points and the two bounds for large points. 

\subsection{Medium points}
\label{subsec medium}

For any $0<r<r'$, let
$$(C\cap\Lambda)_{r,r'}=\{P\in C\cap\Lambda:r<|P|\le r'\}$$
be the set of points in a spherical shell of inner radius $r$ and outer radius $r'$. We will estimate $\#(C\cap\Lambda)_{r,r'}$. The medium point set is the special case that
$$r=\sqrt\frac{\bar\omega_a^2}{16g},\quad r'=1.2\cdot10^9g^{\frac{7}{3}}\sqrt{\bar\omega_a^2}.$$
We only need to bound $\#(C\cap\Lambda)_{r,\kappa r}$ for some $\kappa>1$ and all $r>0$ because of the covering
$$(C\cap\Lambda)_{r,r'}\ \subset\ \bigcup_{i=1}^{\left\lceil\log_\kappa(\frac{r'}{r})\right\rceil}(C\cap\Lambda)_{\kappa^{i-1}r,\kappa^i r}.$$

Denote $n=\mathrm{rk}(\Lambda)\geq 1$.
The space $\Lambda_\RR$ under the norm $|\cdot|$ is isometric to the Euclidean space $\RR^n$.
Denote the sphere
$$S^{n-1}=\{x\in\Lambda_\mathbb{R}:|x|=1\}.$$
Enow $S^{n-1}$ with the standard Euclidean measure $\mathrm{vol}$. 

Take 
$$\theta=\arccos\sqrt{\frac{1.13}{g}}, \quad \kappa=2.$$ 
Then $(\theta, \kappa)$ satisfies the assumption of Theorem \ref{bogomolov}(2). For any $P\in(C\cap\Lambda)_{r,\kappa r}$, let
$$\mathrm{cap}(P,\theta)=\{x\in\Lambda_\mathbb{R}:|x|=1,\angle(x,P)\le\theta\}$$
be the spherical cap centered at $\displaystyle\frac{P}{|P|}$ of angular radius $\theta$ on the sphere $S^{n-1}$.

For any $s\in S^{n-1}$, Theorem \ref{bogomolov}(2) for $x=rs$ implies that $s$ belongs to at most $3.2\cdot10^{11}g^{\frac{17}3}$ spherical caps 
$\mathrm{cap}(P,\theta)$ for $P\in(C\cap\Lambda)_{r,\kappa r}$. 
In other words, under the natural map
$$
\coprod_{P\in(C\cap\Lambda)_{r,\kappa r}} \mathrm{cap}(P,\theta)
\longrightarrow S^{n-1},
$$
every point has at most $3.2\cdot10^{11}g^{\frac{17}3}$ preimages. 
Compare the volumes. We have 
$$\#(C\cap\Lambda)_{r,\kappa r}\cdot \mathrm{vol}(\mathrm{cap}(P,\theta))\leq 3.2\cdot10^{11}g^{\frac{17}3}\mathrm{vol}(S^{n-1}).$$
The total measure is
$$\mathrm{vol}(S^{n-1})=\frac{2\pi^{\frac n2}}{\Gamma\left(\frac n2\right)}.$$
When $n\ge2$, 
$$\mathrm{vol}(\mathrm{cap}(P,\theta))
=\int_{0}^\theta\sin^{n-2}(t)\frac{2\pi^{\frac {n-1}2}}{\Gamma\left(\frac {n-1}2\right)}dt
> \frac{2\pi^{\frac {n-1}2}}{\Gamma\left(\frac {n-1}2\right)}\int_{0}^\theta\sin^{n-2}(t)\cos tdt
=\frac{2\pi^{\frac {n-1}2}}{\Gamma\left(\frac {n-1}2\right)}\frac{\sin^{n-1}\theta}{n-1}
=\frac{\pi^{\frac {n-1}2}\sin^{n-1}\theta}{\Gamma(\frac{n+1}{2})}.$$
This is also true for $n=1$. So
$$\#(C\cap\Lambda)_{r,\kappa r}\le3.2\cdot10^{11}g^{\frac{17}3}\frac{2\sqrt\pi\Gamma(\frac{n+1}{2})}{\Gamma\left(\frac n2\right)}\frac{1}{\sin^{n-1}\theta}\le 3.2\cdot10^{11}g^{\frac{17}3}\frac{\sqrt{2\pi n}}{\sin^{n-1}\theta}.$$
Here the last inequality follows from Wendel's classical inequality (cf. \cite{Wen48}), 
$$
\Gamma(x+s) \leq x^{s} \Gamma(x),\quad x>0, \ 0< s<1. 
$$

It follows that 
$$\#(C\cap\Lambda)_{r,\kappa r}
\le 3.2\cdot10^{11}g^{\frac{17}3} \cdot\sqrt{2\pi n}
\left(1-\frac{1.13}{g}\right)^{-\frac{n-1}{2}}.$$
To simplify the bound further, we will use the following basic fact several times. 

\begin{lem}
\label{calculus}
For any $n,t>0$ and $a>1$,
$$n^t\le\left(\frac{t}{e\log a}\right)^ta^n.$$
\end{lem}
\begin{proof}
It is an easy consequence of
$$\frac{n\log a}{t}\le e^{\frac{n\log a}{t}-1}=\frac{a^{\frac nt}}{e},\quad\forall n,t>0,\ a>1.$$
\end{proof}

Take $t=\frac12$ and
$$a=\left(1+\frac{5}{2\sqrt2g}\right)\sqrt{1-\frac{1.13}{g}}.$$
It satisfies
$$e^{\frac{1}{eg}}\le a\le\sqrt2.$$
Here to prove the first inequality, just note that $\log a-\frac{1}{eg}$ is a concave function of $\frac1g$.
Then Lemma \ref{calculus} implies
$$
\sqrt{n}\leq \left(\frac{1}{2e\log a}\right)^{\frac12} a^n
\leq \sqrt{g} a^{n-1}
= \sqrt{g}\left(1+\frac{5}{2\sqrt2g}\right)^{n-1} \left(1-\frac{1.13}{g}\right)^{\frac{n-1}{2}}.
$$
It follows that
$$
\sqrt{2\pi n} \left(1-\frac{1.13}{g}\right)^{-\frac{n-1}{2}}
\le\sqrt{2\pi g}\left(1+\frac{5}{2\sqrt2g}\right)^{n-1}.$$

Therefore,
$$\#(C\cap\Lambda)_{r,\kappa r}\le3.2\cdot10^{11}g^{\frac{17}3}\cdot\sqrt{2\pi g}\left(1+\frac{5}{2\sqrt2g}\right)^{n-1}\le 8.1\cdot10^{11}g^{\frac{37}6}\left(1+\frac{5}{2\sqrt2g}\right)^{n-1}.$$
So
$$\#(C\cap\Lambda)_{r,r'}\le8.1\cdot10^{11}g^{\frac{37}6}\left(1+\frac{5}{2\sqrt2g}\right)^{n-1}\left\lceil\log_2(\frac{r'}{r})\right\rceil.$$
In particular,
$$\#(C\cap\Lambda)_{\mathrm{medium}}\le8.1\cdot10^{11}g^{\frac{37}6}\left(1+\frac{5}{2\sqrt2g}\right)^{n-1}\lceil\log_2(1.2\cdot10^9g^{\frac{7}{3}}\sqrt{16g})\rceil\le1.64\cdot10^{13}g^7\left(1+\frac{5}{2\sqrt2g}\right)^{n-1}.$$
This gives the bound for medium points.

\subsection{Large points: first bound}
\label{sec large points}

The goal here is prove  
$$\#(C\cap\Lambda)_\mathrm{large}
\le3.4\cdot10^6g^2\left(1+\frac{5}{4\sqrt{g}}\right)^{n-1}.$$
This is the first bound for large points.

\subsubsection*{A bound from sphere packing}

We first recall a result in spherical packing. For any positive integer $n$ and 
$\theta\in(0,\pi)$, denote by $A(n,\theta)$ the maximal number of nonzero points 
$P_1,\dots, P_N\in \mathbb R^n$ such that any two of them form an angle at least 
$\theta$. 
At each point $P_i$, we have the spherical cap of angular radius $\theta/2$  given by
$$\mathrm{cap}(P_i,\theta/2)=\{x\in\RR^n: |x|=1, \angle(x,P_i)\le\theta/2\}.$$
Then the interiors of these $N$ spherical caps are disjoint from each other. 
It follows that 
$A(n,\theta)$ is also the maximal number of spherical caps of angular radius $\theta/2$ in the standard unit sphere in $\RR^n$ such that any two such spherical caps have no common interior points.  

The volume comparison from the previous subsection gives a bound
$$
A(n,\theta) \leq \frac{\mathrm{vol}(S^{n-1})}{\mathrm{vol}(\mathrm{cap}(P,\theta/2))}
\leq \frac{\sqrt{2\pi n}}{\sin^{n-1}\frac{\theta}{2}}.
$$
In our application, $\theta$ is close to $\pi/2$, and thus the bound grows roughly as 
$(\sqrt 2)^{n}$, which is too large for our purpose of making the Vojta constant (i.e. the base) close to 1. Therefore, we are going to seek stronger upper bounds of $A(n,\theta)$.

The following lemma estimate of $A(n,\theta)$ is essentially due to Rankin. 
\begin{lem}\label{rankin}
Let $n$ be a positive integer.  
If $\theta\in(0,\pi)$,  then 
$$A(n,\theta)\le \frac{(1+\sqrt{\cos\theta})n^{\frac32}}{\sqrt{\cos\theta}(1-\cos\theta)^{\frac{n-1}2}}.$$
If $\theta\in (\pi/2, \pi)$, then $A(n,\theta)\le n+1$. 
\end{lem}
\begin{proof}
The result holds for $n=1$ because $A(1,\theta)=2$.
Assume $n\geq2$ in the following.  
We are going to apply results of Rankin \cite{Rankin}. 
Note that our $\theta/2$ corresponds to $\alpha$ in the loc. cit. 
If $\theta\in (\pi/2, \pi)$,  we have $N(n, \theta)=n+1$ by \cite[Theorem 1(iii)]{Rankin}.

Now we consider general
$\theta\in (0,\pi/2]$ in the following. 
By Rankin \cite[Theorem 2]{Rankin}, we have the following bound
$$A(n,\theta)\le\frac{\pi^{\frac{1}{2}}\Gamma\left(\frac{n-1}{2}\right)\sin\beta\tan\beta}{2\Gamma\left(\frac{n}{2}\right)\int_{0}^\beta(\sin x)^{n-2}(\cos x-\cos\beta)dx}, \quad 
n\ge2.$$
Here the angle
$$\beta=\arcsin\left(\sqrt 2\sin  \frac{\theta}2 \right)=\arccos(\sqrt{\cos\theta}).$$

Let us first simplify the bound, which also weakens the bound slightly. Let $t=\sin x\in[0,\sin\beta]$. Then
$$\int_{0}^\beta(\sin x)^{n-2}(\cos x-\cos\beta)dx=\int_0^{\sin\beta}
\left(t^{n-2}-\frac{t^{n-2}}{\sqrt{1-t^2}}\cos\beta \right)dt.$$
By the convexity of $\displaystyle\frac1{\sqrt{x}}$, we have
$$\frac{1}{\sqrt{1-t^2}}\le\frac{1-\cos\beta}{\sin^2\beta\cos\beta}t^2+1.$$
One can also check that the inequality by noting that it holds for $t^2=0$ and $t^2=\sin^2\beta$, and the left-hand side is convex in $t^2$ and the right-hand side is linear in $t^2$. 

It follows that
$$
\int_{0}^\beta(\sin x)^{n-2}(\cos x-\cos\beta)dx
\ge (1-\cos\beta)\int_0^{\sin\beta} \left(t^{n-2}-\frac{t^n}{\sin^2\beta}\right)dt
= \frac{2(1-\cos\beta)\sin^{n-1}\beta}{n^2-1}.
$$
We remark that this result catches the exponential part in the decay for large $n$, comparing the opposite bound 
$$
\int_{0}^\beta(\sin x)^{n-2}(\cos x-\cos\beta)dx
< \int_{0}^\beta(\sin x)^{n-2}\cos x dx
= \frac{\sin^{n-1}\beta}{n-1}.
$$
As a consequence, 
$$A(n,\theta)\le\frac{(n^2-1)\pi^{\frac{1}{2}}\Gamma\left(\frac{n-1}{2}\right)\tan\beta}{4\Gamma\left(\frac{n}{2}\right) (1-\cos\beta)\sin^{n-2}\beta}.$$

By Wendel's classical inequality (cf. \cite{Wen48}) again, 
$$\frac{(n^2-1)\pi^{\frac{1}{2}}\Gamma\left(\frac{n-1}{2}\right)}{4\Gamma\left(\frac{n}{2}\right)}
=\frac{(n+1)\pi^{\frac{1}{2}}\Gamma\left(\frac{n+1}{2}\right)}{2\Gamma\left(\frac{n}{2}\right)}
\leq
\frac{(n+1)\pi^{\frac{1}{2}}n ^{\frac{1}{2}}}{2\sqrt2}
\le n^{\frac{3}{2}}.$$
Then we have 
$$A(n,\theta)\le\frac{n^{\frac{3}{2}}\tan\beta}{ (1-\cos\beta)\sin^{n-2}\beta}
\le\frac{(1+\cos\beta)n^{\frac32}}{\sin^{n-1}\beta\cos\beta}=\frac{(1+\sqrt{\cos\theta})n^{\frac32}}{\sqrt{\cos\theta}(1-\cos\theta)^{\frac{n-1}2}}.$$
\end{proof}

\subsubsection*{Counting points}

Now we return to the bound of large points. 
For any $\theta$, there is a subset $S$ of $(C\cap\Lambda)_\mathrm{large}$
satisfying the following properties.
\begin{enumerate}
    \item For any $P\in(C\cap\Lambda)_\mathrm{large}$, there is a $Q\in S$ such that $\angle(P,Q)\le\theta$ and $|Q|\leq |P|$.
    \item For distinct $Q_1,Q_2\in S$, we have $\angle(Q_1,Q_2)>\theta$. As a consequence, $\# S\le A(n,\theta)$.
\end{enumerate}
This set can be constructed inductively as follows. By the Northcott property, we can take $Q_1\in(C\cap\Lambda)_\mathrm{large}$ such that $|Q_1|$ is minimal.  After taking $Q_1,\dots,Q_n$, let $Q_{n+1}\in(C\cap\Lambda)_\mathrm{large}$ be a point satisfying (2) with minimal $|Q_{n+1}|$. This process ends in $A(n,\theta)$ steps. So for any $P\in(C\cap\Lambda)_\mathrm{large}$, there is a $Q_i\in S$ such that $\angle(P,Q_i)\le\theta$. Let $i$ be the minimal such $i$. Then the minimality of $|Q_i|$ at the $i$-th step implies that $|P|\ge|Q_i|$; otherwise, this step would choose $P$. This verifies (1).

For any $\kappa>1$, we have the covering
\begin{align*}
(C\cap\Lambda)_\mathrm{large}
=&\ \bigcup_{Q\in S}\{P\in C\cap\Lambda:|P|\ge|Q|,\ \angle(P,Q)<\theta\}\\
=&\ \bigcup_{Q\in S}\bigcup_{i=1}^{\infty}\{P\in C\cap\Lambda:\kappa^{i-1}|Q|\le|P|\le\kappa^{i}|Q|,\ \angle(P,Q)<\theta\}.
\end{align*}

Take $\theta=\arccos\left(\sqrt\frac{1.01}{g}\right)$. By the quantitative Vojta inequality in Theorem \ref{vojta}, if $P\in(C\cap\Lambda)_\mathrm{large}$ satisfies $\angle(P,Q)<\theta$ for some $Q\in S$, then $|P|\le10^5g^{\frac52}|Q|$. So
$$(C\cap\Lambda)_\mathrm{large}=\bigcup_{Q\in S}\bigcup_{i=1}^{\lceil\log_{\kappa}(10^5g^{\frac52})\rceil}\{P\in C(K)\cap\Lambda:\kappa^{i-1}|Q|\le|P|\le\kappa^{i}|Q|,\ \angle(P,Q)<\theta\}$$

Take $\kappa=1.15$. For each $Q\in S$ and $i$, suppose we have two distinct points
$$P_1, P_2\in\{P\in C(K)\cap\Lambda:\kappa^{i-1}|Q|\le|P|\le\kappa^i|Q|,\ \angle(P,Q)<\theta\}.$$
Then by the quantitative Mumford inequality in Theorem \ref{mumford}, 
$$\angle(P_1,P_2)\ge\arccos\left(\frac{1.01}g\right).$$
 Take the orthogonal decomposition
$$\frac{P_i}{|P_i|}=a_i\frac{Q}{|Q|}+x_i.$$
Here $a_i\in\mathbb{R}$ and $x_i\in Q^\perp=\{x\in\Lambda_\mathbb R:\langle x,Q\rangle=0\}$. So
$$a_i=\left\langle\frac{P_i}{|P_i|},\frac{Q}{|Q|}\right\rangle=\cos\angle(P_i,Q)>\sqrt\frac{1.01}{g}$$
and
$$a_1a_2+\langle x_1,x_2\rangle=\left\langle\frac{P_1}{|P_1|},\frac{P_2}{|P_2|}\right\rangle=\cos\angle(P_1,P_2)\le\frac{1.01}{g}.$$
We get $\langle x_1,x_2\rangle<0$ and hence $\angle(x_1,x_2)>\frac\pi2$. 

Apply Lemma \ref{rankin} to the 
orthogonal projection of 
$$
\{P\in C\cap\Lambda:\kappa^{i-1}|Q|\le|P|\le\kappa^i|Q|,\ \angle(P,Q)<\theta\}.$$
to the space $Q^\perp$. 
We have 
$$\#\{P\in C\cap\Lambda:\kappa^{i-1}|Q|\le|P|\le\kappa^i|Q|,\ \angle(P,Q)<\theta\}\le 
A(n-1,\frac{\pi}{2})\le
n.$$
As a consequence, we have obtained the following result. 

\begin{prop} \label{intermediate large point}
We have
$$\#(C\cap\Lambda)_\mathrm{large}\le \lceil\log_{1.15}(10^5g^{\frac52})\rceil nA(n,\theta).$$
\end{prop}

By Lemma \ref{rankin} again, the proposition gives
\begin{align*}
\#(C\cap\Lambda)_\mathrm{large} 
\le&\ \lceil\log_{1.15}(10^5g^{\frac52})\rceil \frac{(1+\sqrt{\cos\theta})n^{\frac52}}{\sqrt{\cos\theta}(1-\cos\theta)^{\frac{n-1}2}}\\
=&\ \lceil\log_{1.15}(10^5g^{\frac52})\rceil \,
n^{\frac52}
\left(\sqrt[4]{\frac{g}{1.01}}+1\right) 
\left(1-\sqrt{\frac{1.01}g}\right)^{-\frac{n-1}2}.
\end{align*}

Apply Lemma \ref{calculus} to $t=\frac{5}{2}$ and
$$a=\left({1+\frac{5}{4\sqrt{g}}}\right){\sqrt{{1}-{\sqrt{\frac{1.01}g}}}}.$$
Note that
$$e^{\frac{1}{60\sqrt g}}\le a\le1.2.$$
The lemma gives 
$$n^{\frac{5}{2}}\le\left(\frac{5}{2e\log a}\right)^{\frac{5}{2}} a^n
\le 1.2 \left(\frac{150\sqrt g}{e}\right)^{\frac{5}{2}} a^{n-1}.$$
Then we get
$$n^{\frac52}
\left(\sqrt[4]{\frac{g}{1.01}}+1\right) 
\left(1-\sqrt{\frac{1.01}g}\right)^{-\frac{n-1}2}\le1.2\left(\frac{150\sqrt g}{e}\right)^{\frac52}\left(\sqrt[4]{\frac{g}{1.01}}+1\right)\cdot\left(1+\frac{5}{4\sqrt{g}}\right)^{n-1}\le 5\cdot10^4g^\frac{3}{2}\left(1+\frac{5}{4\sqrt{g}}\right)^{n-1}.$$
Therefore,
$$\#(C\cap\Lambda)_\mathrm{large}\le\lceil\log_{1.15}(10^5g^{\frac52})\rceil\cdot5\cdot10^4g^\frac{3}{2}\left(1+\frac{5}{4\sqrt{g}}\right)^{n-1}\le3.4\cdot10^6g^2\left(1+\frac{5}{4\sqrt{g}}\right)^{n-1}.$$
This proves the first bound of large points.

\subsection{Large points: refined bound}

The goal of this subsection is to prove that 
 for $g\ge142$ and $n\ge1$,
$$\#(C\cap\Lambda)_\mathrm{large}\le 2.4\cdot10^4g^8\left(1+\frac{3\log g}{g}\right)^{n-1}.$$
In fact, this is implied by the following bound. 

\begin{prop} \label{large point refined}
Assume that $g\ge142$ and $n\geq 2g$. Then 
$$\#(C\cap\Lambda)_\mathrm{large}\le 2.4\cdot10^4g\left(1+\frac{3\log g}{g}\right)^{n-1}.$$
\end{prop}

It is easy to extend the proposition to 
 $g\ge 142$ and all $n\geq 1$ with a weaker bound. 
Note that if $n<2g$, we may take a subgroup $\Lambda'$ of rank $2g$ in $J(K)$ containing $\Lambda$ by adding $2g-n$ linearly independent generators, which is possible after replacing $K$ with a finite extension. Then
$$\#(C\cap\Lambda)_\mathrm{large}\le\#(C\cap\Lambda')_\mathrm{large}
\le2.4 \cdot 10^4 g^2  \left( 1+\frac{3\log g}{g}\right)^{2g-1}
\le2.4 \cdot 10^4 g^8.$$
Here the last inequality uses
$$\left(1+\frac{3\log g}{g}\right)^g=\left(\left(1+\frac{3\log g}{g}\right)^\frac{g}{3\log g}\right)^{3\log g}\le e^{3\log g}=g^3.$$
Combining with the proposition, we have the bound for all $g\ge142$ and $n\ge1$.

Then the goal of this subsection is to prove Proposition \ref{large point refined}.
The proof is obtained by improving the bound of $A(n,\theta)$ in Lemma \ref{rankin}.

\subsubsection*{Refined bound from sphere packing}

For $g\geq2$, take $\theta=\arccos\sqrt\frac{1.01}g$ as in the previous subsection. Our goal is to derive a sharp upper bound for $A(n,\theta)$. 

The best known asymptotic bound for $A(n,\theta)$ is due to 
Kabatjanski\u i and Leven\v ste\u in
in \cite{KL,Lev} (cf. \cite{SZ}). For any non-negative integer $m$, the differential equation
$$(1-x^2)\frac{d^2f}{dx^2}-(n-1)x\frac{df}{dx}+m(m+n-2)f=0$$
has a unique polynomial solution up to scaling. Denote its largest root by $x_{(n,m)}$. Then
$$-1<x_{(n,m)}<x_{(n,m+1)}<1.$$
The key is the following result from
Kabatjanski\u i--Leven\v ste\u in
 \cite[\S6 (52)]{KL}.
\begin{thm}\label{thm KL}
If $\cos\theta\le x_{(n,m)}$, then
$$A(n,\theta)\le\frac{4\tbinom{m+n-2}{n-2}}{1-x_{(n,m+1)}}.$$
\end{thm}

By taking
$$m\approx\frac{1-\sin\theta}{2\sin\theta}n,$$
they get $x_{(n,m)}\to\cos\theta$ and
$$\limsup_{n\to\infty}\frac{1}{n}\log A(n,\theta)\le\frac{1+\sin\theta}{2\sin\theta}\log\frac{1+\sin\theta}{2\sin\theta}-\frac{1-\sin\theta}{2\sin\theta}\log\frac{1-\sin\theta}{2\sin\theta}.$$
However, this asymptotic bound is not sufficient for our purpose. We need an explicit estimate of $x_{(n,m)}$. 

In the following, we take 
$$\beta=\arccos\sqrt{\frac{10.7}g}, \quad\
m=\left\lceil\frac{1-\sin\beta}{2\sin\beta}n\right\rceil.$$
Since $x_{(n,m)}\to\cos\beta>\cos\theta$, we have $x_{(n,m)}\ge\cos\theta$ 
for $n$ large enough, and thus we can apply the above bound. We will make this explicit, which gives the following result.

\begin{prop} \label{rankin refined}
For $g\ge142$ and $n\geq 2g$,
\begin{align*}
A(n,\theta)\le 
6\binom{m+n-2}{n-2}
\leq {6e}  \left( \frac{1+\sin\beta}{1-\sin\beta}\right)^{\frac{1-\sin\beta}{2\sin\beta}n+1}
 \left( \frac{1+\sin\beta}{2\sin\beta}\right)^{n-2}.
\end{align*}
\end{prop}

To apply Theorem \ref{thm KL} to prove the proposition, 
our main tool to estimate the zeroes is the following classical Sturm comparison theorem.

\begin{thm} \label{sturm}
For $i=1,2$, let $r_i(x)$ be a nonzero continuous function on an open interval, and let $F_i(x)$ be a nonzero solution of the differential equation 
$$\frac{d^2F}{dx^2}+r_i(x)F=0.$$
Assume that $r_1(x)\leq r_2(x)$. 
Then for any two roots $a<b$ of $F_1(x)$, there is a root of $F_2(x)$ lying in $[a,b]$.
\end{thm}

Now we prove Proposition \ref{rankin refined}. The proof
consists of a few steps. Let $f_{(n,m)} (x)$ be a polynomial solution of the previous differential equation with the maximal root $x_{(n,m)}$. Then 
$$F_{(n,m)}(x)=(1-x^2)^{\frac{n-1}{4}}f_{(n,m)}(x)$$
 is a solution of the differential equation
$$\frac{d^2F}{dx^2}+\frac{p_{(n,m)}-q_{(n,m)}x^2}{4(1-x^2)^2}F=0,$$
where
$$q_{(n,m)}=(2m+n-1)(2m+n-3),$$
$$p_{(n,m)}=q_{(n,m)}-(n-1)(n-5).$$
The largest root of $F_{(n,m)}$ in the interval $(-1,1)$ is still $x_{(n,m)}$. 

Let us first collect some trivial bounds. 
Note that $n\geq 1.99g$ implies
$$m\ge\frac{1-\sin\beta}{2\sin\beta}n\ge\frac{\cos^2\beta}{4}n>5.$$
As $n\ge6$, we have $p_{(n,m)}<q_{(n,m)}$. 
By
$$ \frac{1-\sin\beta}{2\sin\beta}n\leq m \leq \frac{1-\sin\beta}{2\sin\beta}n+1,$$
we have an easy bound
$$
\frac{n}{\sin \beta}\leq 2m+n
\leq \frac{n}{\sin \beta}+2.
$$

Our first bound on the root is the following quick result. 
\begin{lem} \label{lem trivial}
For $n\ge6$, we have
$x_{(n,m)}\le\sqrt{\frac{p_{(n,m)}}{q_{(n,m)}}}$.
\end{lem}

\begin{proof}
We first note that for $x\in(\sqrt{\frac{p_{(n,m)}}{q_{(n,m)}}},1)$, we have
$$\frac{p_{(n,m)}-q_{(n,m)}x^2}{4(1-x^2)^2}<0.$$
Then the lemma follows from Theorem \ref{sturm} by comparing $F_{(n,m)}$
with solutions of 
$$\frac{d^2F}{dx^2}=0. $$

Alternatively, we can prove the lemma as follows.
For the sake of contradiction, assume that  $u=x_{(n,m)}$ lies in $(\sqrt{\frac{p_{(n,m)}}{q_{(n,m)}}},1)$.
Note that $1$ is a root of $F_{(n,m)}$. 
Then the value of $F_{(n,m)}$ has a fixed sign on the interval $(u,1)$. 
If $F_{(n,m)}>0$ on $(u,1)$, then the differential equation gives 
$F_{(n,m)}''>0$ on $(u,1)$. Then $F_{(n,m)}$ is convex on $(u,1)$, and thus 
$$F_{(n,m)}\leq \max\{F_{(n,m)}(u), F_{(n,m)}(1)\}=0$$
 on $(u,1)$. This is a contradiction. 
A similar contradiction can be obtained if $F_{(n,m)}<0$  on $(u,1)$. This proves the lemma.
\end{proof}

For the other direction, we will prove $x_{(n,m)}\ge\cos\theta$ in the following. 
To apply Theorem \ref{sturm}, we need the following elementary result. 

\begin{lem}\label{sturm condition}
Assume that $n\geq 2g$. 
For $x\in [\cos\theta,\frac{1}{3}(2\cos\beta+\cos\theta)]$, 
$$\frac{p_{(n,m)}-q_{(n,m)}x^2}{4(1-x^2)^2}
>\frac{9\pi^2}{4({\cos\beta-\cos\theta})^2}.$$ 
\end{lem}
\begin{proof}
We only need to prove
$$p_{(n,m)}-q_{(n,m)}x^2>\frac{9\pi^2}{({\cos\beta-\cos\theta})^2}$$
on $[\cos\theta,\frac{1}{3}(2\cos\beta+\cos\theta)]$. Since ${p_{(n,m)}-q_{(n,m)}x^2}$ is decreasing on this interval, it is equivalent to
$${p_{(n,m)}-q_{(n,m)}(\frac{2\cos\beta+\cos\theta}{3})^2}
> \frac{9\pi^2}{({\cos\beta-\cos\theta})^2}.$$
On the one hand,
$$\frac{p_{(n,m)}}{q_{(n,m)}}=1-\frac{(n-1)(n-5)}{(2m+n-1)(2m+n-3)}\ge1-\left(\frac{n}{2m+n}\right)^2\ge\cos^2\beta.$$
We also have $q_{(n,m)}\ge n^2$. 
So
$${p_{(n,m)}-q_{(n,m)}(\frac{2\cos\beta+\cos\theta}{3})^2}\ge
{n^2\left(\cos^2\beta-(\frac{2\cos\beta+\cos\theta}{3})^2\right)}=
\left(10.7-(\frac{2\sqrt{10.7}+\sqrt{1.01}}{3})^2\right)\frac{n^2}{g}.$$
On the other hand,
$$\frac{9\pi^2}{( {\cos\beta-\cos\theta} )^2}=\frac{9 \pi^2}{({\sqrt{10.7}-\sqrt{1.01}})^2}g.$$
As $n\geq 2g$, it suffices to prove
$$\left(10.7-(\frac{2\sqrt{10.7}+\sqrt{1.01}}{3})^2\right) \cdot2^2>\frac{9\pi^2}{({\sqrt{10.7}-\sqrt{1.01}})^2}.$$
This can be verified numerically.
\end{proof}

Once we have the lemma, consider the differential equation 
$$\frac{d^2F}{dx^2}+\frac{9\pi^2}{4(\cos\beta-\cos\theta)^2}F=0,$$
which has a solution
$$F^*(x)=\sin \frac{3\pi(x-\cos\theta)}{ {2(\cos\beta-\cos\theta)} }.$$
Then $\cos\theta$ and $\frac{1}{3}(2\cos\beta+\cos\theta)$ are roots of $F^*$.
By Theorem \ref{sturm} and Lemma \ref{sturm condition}, 
the function $F_{(n,m)}$ has a root in $[\cos\theta,\frac{1}{3}(2\cos\beta+\cos\theta)]$.
As a consequence, $x_{(n,m)}\ge\cos\theta$. 
Then Theorem \ref{thm KL} implies that 
$$A(n,\theta)\le\frac{4\tbinom{m+n-2}{n-2}}{1-x_{(n,m+1)}}.$$

By Lemma \ref{lem trivial}, we have
$$x_{(n,m+1)}^2\le{\frac{p_{(n,m+1)}}{q_{(n,m+1)}}}
=1-\frac{(n-1)(n-5)}{(2m+n+1)(2m+n-1)}.$$
For $g\geq 142$ and $n\geq 2g$, we have
$$2m+n\leq 
\frac{n}{\sin \beta}+2
=\frac{n}{\sqrt{1-\frac{10.7}{g}}}+2
\leq 1.04 n+2. 
$$
We also have
$$
\frac{(n-1)(n-5)}{(2m+n+1)(2m+n-1)}
\geq \frac{(n-1)(n-5)}{(1.04 n+3)(1.04 n+1)}
\geq 0.89. 
$$
It follows that 
$$x_{(n,m+1)}^2\le1-0.89<\frac19, \quad 
x_{(n,m+1)}<\frac13.$$
Therefore, the bound implies
$$
A(n,\theta)\le6\binom{m+n-2}{n-2}.$$
This proves the first inequality of Proposition \ref{rankin refined}. 

The second inequality of Proposition \ref{rankin refined} follows from the following lemma.

\begin{lem}\label{binomial}
We have
$$\binom{m+n-2}{n-2}\leq e  \left( \frac{1+\sin\beta}{1-\sin\beta}\right)^{\frac{1-\sin\beta}{2\sin\beta}n+1}
 \left( \frac{1+\sin\beta}{2\sin\beta}\right)^{n-2}.$$
\end{lem}
\begin{proof}
Write
\begin{align*}
\log \binom{m+n-2}{n-2}
= \log (m+n-2)!-\log m!-\log (n-2)!.
\end{align*}
By Robbins' effective version of Stirling's formula, for all positive integers $k$, 
$$
\sqrt{2\pi k} \big(\frac{k}{e}\big)^k  e^{\frac{1}{12k+1}} < 
k!< \sqrt{2\pi k} \big(\frac{k}{e}\big)^k  e^{\frac{1}{12k}}.
$$
This gives an upper bound of $\log (m+n-2)!$, and a lower bound for each of $\log m!$ and $\log (n-2)!$. 
As a consequence, we obtain 
\begin{align*}
& \log \binom{m+n-2}{n-2}\\
\leq&\ -\frac{1}{2}\log(2\pi)+\frac{1}{2}\log \frac{m+n-2}{m(n-2)}+(m+n-2)\log (m+n-2)
-m\log m-(n-2)\log (n-2)\\
\leq &\ -\frac{1}{2}\log(6\pi)
+m\log \frac{m+n-2}{m}+(n-2)\log \frac{m+n-2}{n-2}.
\end{align*}
As 
$$\displaystyle \frac{1-\sin\beta}{2\sin\beta}n\leq m \leq \frac{1-\sin\beta}{2\sin\beta}n+1,$$
we have
\begin{align*}
& \log \binom{m+n-2}{n-2}\\
<&\ -\frac{1}{2}\log(6\pi)
+\left(\frac{1-\sin\beta}{2\sin\beta}n+1\right)\log \frac{\frac{1-\sin\beta}{2\sin\beta}n+n-2}{\frac{1-\sin\beta}{2\sin\beta}n}+(n-2)\log \frac{\frac{1-\sin\beta}{2\sin\beta}n+n-1}{n-2}\\
\leq&\ -\frac{1}{2}\log(6\pi)
+\left(\frac{1-\sin\beta}{2\sin\beta}n+1\right)\log \frac{1+\sin\beta}{1-\sin\beta}+(n-2)\log \frac{1+\sin\beta}{2\sin\beta}
+(n-2)\log \frac{n}{n-2}\\
\leq&\ 
\left(\frac{1-\sin\beta}{2\sin\beta}n+1\right)\log \frac{1+\sin\beta}{1-\sin\beta}+(n-2)\log \frac{1+\sin\beta}{2\sin\beta}
+1.
\end{align*}
\end{proof}

\subsubsection*{Counting points}

Now it is easy to prove Proposition \ref{large point refined}, which asserts that 
 for $g\ge142$ and $n\ge 2g$,
$$\#(C\cap\Lambda)_\mathrm{large}\le2.4\cdot 10^4g\left(1+\frac{3\log g}{g}\right)^{n-1}.$$
Assume that $g\ge142$ and $n\ge 2g$ in the following.

By Proposition \ref{intermediate large point}, 
$$\#(C\cap\Lambda)_\mathrm{large}\le \lceil\log_{1.15}(10^5g^{\frac52})\rceil nA(n,\theta).$$
By Proposition \ref{rankin refined}, this further gives
$$\#(C\cap\Lambda)_\mathrm{large}\le 
\lceil\log_{1.15}(10^5g^{\frac52})\rceil \cdot 6en  
\left( \frac{1+\sin\beta}{1-\sin\beta}\right)^{\frac{1-\sin\beta}{2\sin\beta}n+1}
 \left( \frac{1+\sin\beta}{2\sin\beta}\right)^{n-2}.$$
 Here $\beta$ is defined by 
 $$\cos\beta=\sqrt{\frac{10.7}{g}}.$$
Note that 
$$\left( \frac{1+\sin\beta}{1-\sin\beta}\right)
 \left( \frac{1+\sin\beta}{2\sin\beta}\right)^{-2}=\frac{4\sin^2\beta}{\cos^2\beta}\leq \frac{4}{\cos^2\beta}=\frac{4g}{10.7}.$$
We have 
$$\#(C\cap\Lambda)_\mathrm{large}\le 
\lceil\log_{1.15}(10^5g^{\frac52})\rceil \cdot 6.1gn  
B^n,$$
where
$$
B= \left( \frac{1+\sin\beta}{1-\sin\beta}\right)^{\frac{1-\sin\beta}{2\sin\beta}}
\frac{1+\sin\beta}{2\sin\beta}.
$$

We have the following elementary estimate. 
\begin{lem}
For $g\geq 142$, 
$$
\log B< \frac{2.84\log g}{g}, 
$$
and
$$
B
<\left(1+\frac{3\log g}{g}\right) \left(1+\frac{0.01\log g}{g}\right)^{-1}.
$$
\end{lem}
\begin{proof}
Denote 
$$
s=\frac{1-\sin\beta}{2\sin\beta} >0.
$$
Note that for $g\geq 142$, 
$$
\cos\beta=\sqrt{\frac{10.7}{g}},\quad \
\sin\beta=\sqrt{1-\frac{10.7}{g}} > 0.961.
$$
Then we have
$$
s=\frac{1-\sin\beta}{2\sin\beta}
= \frac{\cos^2\beta}{2\sin\beta(1+\sin\beta)}
= \frac{10.7}{2\sin\beta(1+\sin\beta)g}
<\frac{2.84}{g}.
$$

By definition, 
$$
B= \left( \frac{1+\sin\beta}{1-\sin\beta}\right)^{\frac{1-\sin\beta}{2\sin\beta}}
\frac{1+\sin\beta}{2\sin\beta}
=\left(1+\frac{1}{s}\right)^{s}(s+1).
$$
Denote
$$
b=\log B= s\log \left(1+\frac{1}{s}\right)+\log (s+1)= (s+1)\log (s+1)-s \log s.
$$
View $b$ as a function of $s>0$. 
Then we have 
$$
\frac{db}{ds}=\log (s+1)-\log s >0.
$$
It follows that $b$ is increasing in $s$. 
Therefore, 
\begin{align*}
b< &\ \left(1+\frac{2.84}{g}\right) \log \left(1+\frac{2.84}{g}\right) -\frac{2.84}{g} \log \frac{2.84}{g}\\
= &\ \left(1+\frac{2.84}{g}\right) \log \left(1+\frac{2.84}{g}\right) -\frac{2.84}{g} \log {2.84}+\frac{2.84}{g} \log g\\
< &\ \left(1+\frac{2.84}{g}-\log {2.84}\right) \log \left(1+\frac{2.84}{g}\right) +\frac{2.84}{g} \log g\\
<& \frac{2.84}{g}\log g.
\end{align*}
This gives the bound for $B$ by taking quotients of
$$
e^{\frac{2.85}{g}\log g} 
<1+\frac{3\log g}{g}, \quad
e^{\frac{0.01}{g}\log g}
> 1+\frac{0.01}{g}\log g. 
$$
\end{proof}

Therefore, we have
\begin{align*}
\#(C\cap\Lambda)_\mathrm{large}
< \lceil\log_{1.15}(10^5g^{\frac52})\rceil \cdot 6.1gn  
\left( 1+\frac{3\log g}{g}\right)^n\left(1+\frac{0.01\log g}{g}\right)^{-n}.
\end{align*}
By the binomial bound
$$
\left(1+\frac{0.01\log g}{g}\right)^{n}
\geq \frac{0.01\log g}{g}n,
$$
we have 
\begin{align*}
\#(C\cap\Lambda)_\mathrm{large}
< &\ \frac{\log_{1.15}(10^5g^{\frac52})+1}{\log g} \cdot 610 g^2  
\left( 1+\frac{3\log g}{g}\right)^n \\
< &\ 34.8 \cdot 610 g^2  
\left( 1+\frac{3\log g}{g}\right)^n \\
< &\ 2.13 \cdot 10^4 g^2  
\left( 1+\frac{3\log g}{g}\right)^n\\
< &\ 2.4 \cdot 10^4 g^2  
\left( 1+\frac{3\log g}{g}\right)^{n-1}.
\end{align*}
This proves Proposition \ref{large point refined}, and thus finishes the proof of Theorem \ref{mordell_lang}.

\subsection{Function fields}

In the end, we have the following counterpart of Theorem \ref{mordell_lang}  over function fields of characteristic $0$. As expected, it is proved by a similar method, and it has better constants.  

\begin{thm}[Theorem \ref{geometric_mordell_lang intro}]
\label{geometric_mordell_lang}
Let $K$ be a field of characteristic $0$, $C$ a curve of genus $g\geq 2$ over $K$, and $J$ the Jacobian variety of $C$. Assume that $C$ is not isotrivial over $\QQ$ in that  $C_{\overline K}$ cannot be descended to $\overline\QQ$. Then for any line bundle $\alpha$  of degree $1$ on $C$ and any subgroup $\Lambda\subset J(K)$ of finite rank,
$$\#((C-\alpha)\cap\Lambda)\le1.8\cdot10^6g^{3}\cdot\left(1+\frac5{4\sqrt g}\right)^{\mathrm{rk}(\Lambda)},$$
$$\#((C-\alpha)\cap\Lambda)\le2.5\cdot10^4g^{8}\cdot\left(1+\frac{3\log g}{g}\right)^{\mathrm{rk}(\Lambda)}.$$
\end{thm}

\begin{proof}
By a standard argument as in \cite[\S4]{Yu}, we may assume $K$ is a function field of one variable over a field $k$, and that $C$ is non-isotrivial over $k$. As in the proof of Theorem \ref{mordell_lang}, we assume $\alpha=\alpha_0$ and divide $C\cap\Lambda$ into the following three parts:
\begin{align*}
(C\cap\Lambda)_{\mathrm{small}}=&\ \left\{P\in C\cap\Lambda:|P|\le\sqrt\frac{\bar\omega_a^2}{16g}\right\},  \\
(C\cap\Lambda)_{\mathrm{medium}}=&\ \left\{P\in C\cap\Lambda:\sqrt\frac{\bar\omega_a^2}{16g}<|P|\le10^4g\sqrt{\bar\omega_a^2}\right\},  \\
(C\cap\Lambda)_{\mathrm{large}}=&\ \left\{P\in C\cap\Lambda:|P|>10^4g\sqrt{\bar\omega_a^2}\right\}.
\end{align*}
Here the boundary between medium points and large points is adjusted to align the 	condition in Theorem \ref{vojta_function_field}. It suffices to prove that
\begin{align*}
\#(C\cap\Lambda)_{\mathrm{small}} \leq &\ 200g^2,  \\
\#(C\cap\Lambda)_{\mathrm{medium}}
\leq&\ 10^{4}g^3\left(1+\frac{5}{2\sqrt2g}\right)^{n-1},  \\
\#(C\cap\Lambda)_{\mathrm{large}}
\leq &\  3.4\cdot10^6g^2\left(1+\frac{5}{4\sqrt{g}}\right)^{n-1},
\end{align*}
and that for $g\ge142$ and $n\ge1$,
$$\#(C\cap\Lambda)_\mathrm{large}\le2.4\cdot10^4g^8\left(1+\frac{3\log g}{g}\right)^{n-1}.$$
For small points, we apply Theorem \ref{geometric_bogomolov}(1) with $r=\sqrt{\frac{1}{16g}}$ and get
$$\#(C\cap\Lambda)_{\mathrm{small}} \leq 200g^2.$$
For medium points, take $$\theta=\arccos\sqrt{\frac{1.13}{g}}, \quad \kappa=1.67$$
in Theorem \ref{geometric_bogomolov}(2). Then the method of \S\ref{subsec medium} gives
\begin{align*}
\#(C\cap\Lambda)_{r,\kappa r}&\ \le(16g^2+32g+124)\frac{2\sqrt\pi\Gamma(\frac{n+1}{2})}{\Gamma\left(\frac n2\right)}\frac{1}{\sin^{n-1}\theta}\\
&\ \le(16g^2+32g+124)\sqrt{2\pi g}\left(1+\frac{5}{2\sqrt2g}\right)^{n-1}\\
&\ \le158g^{\frac52}\left(1+\frac{5}{2\sqrt2g}\right)^{n-1}.
\end{align*}
So
$$\#(C\cap\Lambda)_{\mathrm{medium}} \leq \lceil\log_{1.67}(10^4g\sqrt{16g})\rceil158g^{\frac52}\left(1+\frac{5}{2\sqrt2g}\right)^{n-1}\le10^4g^{3}\left(1+\frac{5}{2\sqrt2g}\right)^{n-1}.$$
Finally, since the angles and ratios in Theorem \ref{vojta_function_field} and \ref{mumford_function_field} are exactly the same as their counterparts in Theorem \ref{vojta} and Theorem \ref{mumford}, the bounds of large points can be obtained similarly.
\end{proof}

\section{Further bounds for hyperelliptic curves}

Recall that a hyperelliptic curve over a field $k$ is a curve $C$ of genus $g\geq2$ over $k$ with a finite separable morphism $C\to \PP^1_k$ of degree $2$. 
There are exactly $2g+2$ Weierstrass points of $C$, i.e. points of $C(\overline k)$ ramified over $\PP^1_k$.
In this section, we consider some consequences of our quantitative Mordell conjecture for hyperelliptic curves.

\subsection{Bounds in terms of discriminants}

The goal here is to prove the following refinement of Theorem \ref{bound by bad reduction} for hyperelliptic curves.
  
\begin{thm}[Theorem \ref{hyperelliptic_intro}]
\label{hyperelliptic}
Let $K$ be a number field of degree $d$ over $\mathbb{Q}$. Let $f(x)\in O_K[x]$ be a monic and square-free polynomial.  Let $\Delta_f\in O_K$ be the discriminant of $f$.
\begin{enumerate}
\item[(1)]If $\deg f=2g+1$ for some integer $g\ge 2$, then
$$\#\{(x,y)\in K^2:y^2=f(x)\}\le
10^{13}\cdot 2^dg^{9d\log_2(2gd+d)+18d+8}(2gd+d)^{\frac32d}\cdot|N_{K/\mathbb{Q}}\Delta_{f}|^{4\log_2g+\frac12}\cdot|\Delta_{K/\mathbb{Q}}|^{3\log_2g+\frac12}.$$
\item[(2)]If $\deg f=2g+2$ for some integer $g\ge 2$, then
$$\#\{(x,y)\in K^2:y^2=f(x)\}\le
10^{13}\cdot2^{6d+1}g^{9d\log_2(2gd+2d)+18d+8}(2gd+2d)^{\frac92d}
\cdot|N_{K/\mathbb{Q}}\Delta_{f}|^{4\log_2g+\frac52}
\cdot|\Delta_{K/\mathbb{Q}}|^{3\log_2g+\frac32}.$$
\end{enumerate}
\end{thm}

To prove the theorem, we need parts (3)-(4) of the following theorem to bound the Mordell--Weil ranks of  hyperelliptic curves. 

\begin{thm} \label{descent}
Let $C$ be a hyperelliptic curve of genus $g\ge2$ over a number field $K$. Let $J$ be the Jacobian variety of $C$. 
\begin{enumerate}
\item 
Let $L_1,\dots,L_m$ be the residue fields of the Galois orbits of the Weierstrass points of $C$. 
Let $S$ be the subset of $M_K$ consisting of all archimedean places, all places above $2$, and all non-archimedean places where $C$ has bad reduction. 
Then
$$\mathrm{rk}\, J(K)\le1+(2g+3)\#S+\sum_{i=1}^m\dim_{\mathbb{F}_2}\mathrm{Cl}_S(L_i)[2]+\dim_{\mathbb{F}_2}\mathrm{Cl}_S(K)[2].$$
Here  $\mathrm{Cl}_S(L)$ denotes the ideal class group of the $S$-integer ring $O_{L,S}$ for any finite extension $L$ of $K$.  
\item
In part (1), if $[L_i:K]$ is odd for some $i$, then
$$\mathrm{rk}\, J(K)\le2g\#S+\sum_{i=1}^m\dim_{\mathbb{F}_2}\mathrm{Cl}_S(L_i)[2]-2\dim_{\mathbb{F}_2}\mathrm{Cl}_S(K)[2].$$

\item 
Assume that $C$ contains an affine curve $\mathrm{Spec}\, K[x,y]/(y^2-f(x))$ as an open subvariety, where $f(x)\in O_K[x]$ is a monic and square-free polynomial with discriminant $\Delta_f\in O_K$ and with $\deg f=2g+1$. 
Then
$$
\mathrm{rk}\, J(K)\le
(\frac43g+\frac12)\log_2|N_{K/\mathbb{Q}}\Delta_f|+(g+\frac12)\log_2|\Delta_{K/\mathbb{Q}}|+(3g+\frac32)d\log_2(2gd+d)+(6g+1)d.$$

\item 
Assume that $C$ contains an affine curve $\mathrm{Spec}\, K[x,y]/(y^2-f(x))$ as an open subvariety, where $f(x)\in O_K[x]$ is a monic and square-free polynomial with discriminant $\Delta_f\in O_K$ and with $\deg f=2g+2$. 
Then
$$
\mathrm{rk}\, J(K)\le
(\frac43g+\frac52)\log_2|N_{K/\mathbb{Q}}\Delta_f|+(g+\frac32)\log_2|\Delta_{K/\mathbb{Q}}|+(3g+\frac92)d\log_2(2gd+2d)+(6 g+6)d+1.$$
\end{enumerate}
\end{thm}

If all Weierstrass points are rational, then $L_1=\cdots=L_{2g+2}=K$ and the second bound becomes
$$\mathrm{rk}\, J(K)\le2g\#S+2g\dim_{\mathbb{F}_2}\mathrm{Cl}_S(K)[2],$$
which is just \cite[Theorem 2]{OT89}. If $K=\mathbb{Q}$ and $C$ has a rational Weierstrass point, the second bound can be derived from \cite[Lemma 4.10]{Sto01}. All these results are based on 2-descents, and a similar method also leads to our general case. Since the bound is implicit in the loc. cit., we will present a proof later. Before that, let us first prove Theorem \ref{hyperelliptic_intro}.

\begin{proof}[Proof of Theorem \ref{hyperelliptic}]
We first prove (1). 
Apply Theorem \ref{descent}(3) to Theorem \ref{mordell_lang}. Since
$$\left(1+\frac{3\log g}{g}\right)^{g}\le g^3,$$
we have
\begin{align*}
\# C(K)\le&\ 10^{13}g^8\cdot
\left(1+\frac{5}{4\sqrt{g}}\right)^{\frac12\log_2|N_{K/\mathbb{Q}}\Delta_{f}|+\frac12\log_2|\Delta_{K/\mathbb{Q}}|+\frac32d\log_2(2gd+d)+d}\\
&\quad \cdot\left(1+\frac{3\log g}{g}\right)^{g\left(\frac43\log_2|N_{K/\mathbb{Q}}\Delta_f|+\log_2|\Delta_{K/\mathbb{Q}}|+3d\log_2(2gd+d)+6d\right)}\\
\le&\ 10^{13}g^8\cdot 2^{\frac12\log_2|N_{K/\mathbb{Q}}\Delta_{f}|+\frac12\log_2|\Delta_{K/\mathbb{Q}}|+\frac32d\log_2(2gd+d)+d}\\
&\quad\cdot g^{4\log_2|N_{K/\mathbb{Q}}\Delta_f|+3\log_2|\Delta_{K/\mathbb{Q}}|+9d\log_2(2gd+d)+18d}\\
=&\ 10^{13}\cdot 2^dg^{9d\log_2(2gd+d)+18d+8}(2gd+d)^{\frac32d}\cdot|N_{K/\mathbb{Q}}\Delta_{f}|^{4\log_2g+\frac12}\cdot|\Delta_{K/\mathbb{Q}}|^{3\log_2g+\frac12}.
\end{align*}
This proves (1).

To prove (2), apply Theorem \ref{descent}(4) to Theorem \ref{mordell_lang}. 
We have
\begin{align*}
\# C(K)\le&\ 10^{13}g^8\cdot 2^{\frac52\log_2|N_{K/\mathbb{Q}}\Delta_{f}|+\frac32\log_2|\Delta_{K/\mathbb{Q}}|+\frac{9}{2}d\log_2(2gd+2d)+6d+1}\\
&\quad\cdot g^{4\log_2|N_{K/\mathbb{Q}}\Delta_f|+3\log_2|\Delta_{K/\mathbb{Q}}|+9d\log_2(2gd+2d)+18d}\\
=&\ 10^{13}\cdot2^{6d+1}g^{9d\log_2(2gd+2d)+18d+8}(2gd+2d)^{\frac92d}
\cdot|N_{K/\mathbb{Q}}\Delta_{f}|^{4\log_2g+\frac52}
\cdot|\Delta_{K/\mathbb{Q}}|^{3\log_2g+\frac32}.
\end{align*}
This proves (2).
\end{proof}

To prove Theorem \ref{descent}, we need the following classical fact, for which we include a proof for lack of a good reference. 

\begin{lem} \label{torsion}
Let $k$ be an algebraically closed field with $\mathrm{char}(k)\neq2$. 
Let $C$ be a hyperelliptic curve of genus $g\geq 2$ over $k$. 
Let $P_0, \dots, P_{2g+1}$ be Weierstrass points of $C$. 
Then as an $\FF_2$-vector space, $\mathrm{Pic}(C)[2]$ is generated by 
$\{P_j-P_0: j=1,\dots, 2g+1\}$ with the only relation 
$\sum_{j=1}^{2g+1}(P_j-P_0)=0$. 
\end{lem}
\begin{proof}
Denote by $\mathrm{Div}_W^0(C)$ the $\ZZ$-module of divisors of $C$ supported at Weierstrass points and of degree 0, i.e. formal $\ZZ$-linear combinations of $P_j-P_0$ with $j=1,\dots, 2g+1$. 
Consider the natural map
$$
\Phi: \mathrm{Div}_W^0(C) \lra \mathrm{Pic}^0(C). 
$$
Denote by $\mathcal P$ the subgroup of $\mathrm{Div}_W^0(C)$ generated by $2\, \mathrm{Div}_W^0(C)$ and 
$$
D_0=(P_1+\cdots+P_{2g+1})-(2g+1)P_0. 
$$
We claim that $\mathcal P=\ker(\Phi)$. 

We first prove $\mathcal P\subseteq \ker(\Phi)$.
Denote by $\pi:C\to \PP^1$ the finite morphism of degree 2. 
As $2P_j$ is linearly equivalent to $\pi^*\CO(1)$, we see that
$2\, \mathrm{Div}_W^0(C)$ lies in $\ker(\Phi)$. 

As the extension $k(C)/k(\PP^1)$ has degree 2, we can write 
$k(C)=k(\PP^1)(\sqrt{f})$ for some non-constant $f\in k(\PP^1)^\times$. We can further assume that $P_0$ is over $\infty\in\mathbb{P}^1$ and $f$ is a square-free polynomial. Then $\mathrm{div}(\sqrt f)=D_0$. 
As $\mathrm{div}(\sqrt f)$ and $2\, \mathrm{Div}_W^0(C)$ are contained in $\ker(\Phi)$, we see that $\mathcal P$ is contained in $\ker(\Phi)$.

Now we prove $\ker(\Phi) \subseteq \mathcal P$.
For any $D\in \ker(\Phi)$, write $D=\mathrm{div}(\alpha)$ for some $\alpha\in k(C)^\times$.
Denote by $\sigma:C\to C$ the non-trivial automorphism of $C$ over $\PP^1$.
Then 
$$\mathrm{div}(\sigma^*\alpha)=\sigma^*D=D=\mathrm{div}(\alpha),\quad\
\mathrm{div}(\sigma^*\alpha/\alpha)=0.
$$
It follows that $\sigma^*\alpha/\alpha=c$ for some $c\in k^\times$. 
Apply a further action by $\sigma^*$. We have $\alpha/\sigma^*\alpha=c$.
It follows that $c^2=1$.
This gives $\sigma^*(\alpha^2)=\alpha^2$, and thus 
$\alpha^2\in k(\PP^1)^\times$. 
Comparing with the above expression 
$k(C)=k(\PP^1)(\sqrt{f})$, we have $\alpha= h \sqrt f$ or $\alpha= h$ for some $h\in k(\PP^1)^\times$. 
Then $D=\mathrm{div}(\alpha)$ lies in $\mathcal P$.
This proves  $\ker(\Phi) \subseteq \mathcal P$.

Once $\ker(\Phi)= \mathcal P$, we have an injection 
$$
\Phi: \mathrm{Div}_W^0(C)_{\FF_2}/ (\FF_2 D_0) \lra \mathrm{Pic}^0(C)[2]
$$
of $\FF_2$-vector spaces. 
It is an isomorphism since both sides have dimension $2g$. 
\end{proof}

Now we prove the bounds of the Mordell--Weil rank in Theorem \ref{descent}. 

\begin{proof}[Proof of Theorem \ref{descent}(1)(2)]
We first prove (1). 
Let $K_S$ be the maximal extension of $K$ which is unramified outside $S$. Let $G=\mathrm{Gal}(K_S/K)$ be the Galois group. The Kummer map
$$J({K})/2J({K})\longrightarrow H^1(G,J[2](K_S))$$
is injective, and thus
$$\mathrm{rk\ }J({K})\le\dim H^1(G,J[2](K_S)).$$
It suffices to prove the bound for $\dim H^1(G,J[2](K_S))$.

Let $P_j\in C(K_S)$ for $j=0,1,\dots,2g+1$ be the Weierstrass points of $C$. 
Denote by $\mathrm{Div}_W(C_{K_S})$ the $\ZZ$-module of divisors of $C_{K_S}$ supported at Weierstrass points.
Denote by $\mathrm{Div}_W^0(C_{K_S})$ 
the subgroup of $\mathrm{Div}_W(C_{K_S})$ of divisors of degree 0. 
By tensor products, we also have $\mathrm{Div}_W(C_{K_S})_{\FF_2}$ and $\mathrm{Div}_W^0(C_{K_S})_{\FF_2}$. 

By Lemma \ref{torsion}, the natural map 
$$
\mathrm{Div}_W^0(C_{K_S}) \lra \mathrm{Pic}^0(C_{\overline K})
$$
has an image $\mathrm{Pic}^0(C_{\overline K})[2]=J[2](\overline K)$, and thus induces a surjection 
$$
\mathrm{Div}_W^0(C_{K_S})_{\FF_2} \lra J[2](K_S).
$$
The kernel of this map is generated by 
$$
D_1=\sum_{i=0}^{2g+1}  P_i, 
$$
Here $D_1$ is viewed as an element of
$$
\mathrm{Div}_W^0(C_{K_S})_{\FF_2}
=\ker(\deg: \mathrm{Div}_W(C_{K_S})_{\FF_2} \lra \FF_2),
$$
since $\deg(D_1)=2g+2$ is 0 in $\FF_2$. 
Hence, we have a $G$-equivariant exact sequence
$$
0\lra \FF_2 \longrightarrow \mathrm{Div}_W^0(C_{K_S})_{\FF_2} \stackrel{}{\longrightarrow} J[2](K_S) \longrightarrow 0. 
$$
The definition also gives a $G$-equivariant exact sequence
$$
0 \longrightarrow \mathrm{Div}_W^0(C_{K_S})_{\FF_2}
\longrightarrow \mathrm{Div}_W(C_{K_S})_{\FF_2} \stackrel{\deg}{\longrightarrow} 
\FF_2 \longrightarrow 0. 
$$

Fix an embedding of $L_i$ into $K_S$ for each $i$. Denote by $G_i=\mathrm{Gal}(K_S/L_i)$ the Galois group of $K_S$ over $L_i$. Then
$$\mathrm{Div}_W(C_{K_S})_{\FF_2}=\bigoplus_{i=1}^m\left(\bigoplus_{K(P_j)=L_i}\FF_2 P_j\right)\cong\bigoplus_{i=1}^m\mathrm{Ind}_{G_i}^G\FF_2.$$
By taking the group cohomology of the two short exact sequences, we get
\begin{align*}
\dim H^1(G,J[2](K_S))\le&\ \dim H^1(G,\mathrm{Div}_W^0(C_{K_S})_{\FF_2})+\dim H^2(G,\FF_2)\\
\le&\ \dim H^0(G,\FF_2)+\dim H^1(G,\mathrm{Div}_W(C_{K_S})_{\FF_2})+\dim H^2(G,\FF_2)\\
=&\ 1+\sum_{i=1}^m\dim H^1(G_i,\FF_2)+\dim H^2(G,\FF_2).
\end{align*}

Let $O_{S}$ be the integral closure of $O_{K,S}$ in $K_S$. We identify $\FF_2$ with $\{\pm1\}$ via the unique isomorphism. Then the short exact sequence
$$0\longrightarrow\FF_2\longrightarrow O_S^\times\xrightarrow{[2]} O_S^\times\longrightarrow 0$$
implies
$$\dim H^1(G,\FF_2)=\dim H^0(G,O_S^\times)\otimes\FF_2+\dim H^1(G,O_S^\times)[2],$$
$$\dim H^2(G,\FF_2)=\dim H^1(G,O_S^\times)\otimes\FF_2+\dim H^2(G,O_S^\times)[2].$$
By \cite[Propostion 8.3.11]{NSW},
$$H^0(G,O_S^\times)=O_{K,S}^\times,$$
$$H^1(G,O_S^\times)=\mathrm{Cl}_S(K),$$
$$H^2(G,O_S^\times)[2]=\mathrm{ker}(\bigoplus_{\substack{v\in S\\ K_v\ne\mathbb{C}}}\FF_2\stackrel{\Sigma}\longrightarrow\FF_2).$$
So we have
$$\dim H^1(G,\FF_2)=\dim O_{K,S}^\times/(O_{K,S}^\times)^2+\dim \mathrm{Cl}_S(K)[2]$$
and
$$\dim H^2(G,\FF_2)\le\dim \mathrm{Cl}_S(K)/2\mathrm{Cl}_S(K)+\#S
=\dim \mathrm{Cl}_S(K)[2]+\#S.$$
Here the last equality holds since $\mathrm{Cl}_S(K)$ is finite.

Note that $K_S$ is also the maximal extension of $L_i$ which is unramified outside $S$. Thus,
\begin{align*}
\dim H^1(G_i,\FF_2)=&\ \dim O_{L_i,S}^\times/(O_{L_i,S}^\times)^2+\dim \mathrm{Cl}_S(L_i)[2]\\
\le&\ [L_i:K]\dim O_{K,S}^\times/(O_{K,S}^\times)^2+\dim \mathrm{Cl}_S(L_i)[2],
\end{align*}
Here the inequality follows from the structure of the $S$-unit groups. 

In summary,
\begin{align*}
\dim H^1(G,J[2](K_S))\le&\ 1+\sum_{i=1}^m\left([L_i:K]\dim O_{K,S}^\times/(O_{K,S}^\times)^2+\dim \mathrm{Cl}_S(L_i)[2]\right)+\dim \mathrm{Cl}_S(K)[2]+\#S\\
=&\ 1+(2g+2)\dim O_{K,S}^\times/(O_{K,S}^\times)^2+\sum_{i=1}^m\dim \mathrm{Cl}_S(L_i)[2]+\dim \mathrm{Cl}_S(K)[2]+\#S\\
\le&\ 1+(2g+3)\#S+\sum_{i=1}^m\dim \mathrm{Cl}_S(L_i)[2]+\dim \mathrm{Cl}_S(K)[2].
\end{align*}
This finishes the proof of Theorem \ref{descent}(1).

For Theorem \ref{descent}(2), by assumption, $[L_{i}:K]$ is odd for some $i$, 
so we assume that $[K(P_0):K]$ is odd. We claim that in this case, both short exact sequences
$$
0\lra \FF_2 \longrightarrow \mathrm{Div}_W^0(C_{K_S})_{\FF_2} \stackrel{}{\longrightarrow} J[2](K_S) \longrightarrow 0 
$$
and
$$
0 \longrightarrow \mathrm{Div}_W^0(C_{K_S})_{\FF_2}
\longrightarrow \mathrm{Div}_W(C_{K_S})_{\FF_2} \stackrel{\deg}{\longrightarrow} 
\FF_2 \longrightarrow 0
$$
are $G$-split. Indeed, to split the first short exact sequence, we take the projection
$$\mathrm{Div}_W^0(C_{K_S})_{\FF_2} \longrightarrow \FF_2,
\quad\sum_{j=0}^{2g+1}a_j P_j \longmapsto 
\sum_{P_j \in \mathrm{Gal}(K_S/K) P_0}a_j.$$
To split the second short exact sequence, we take the injection 
$$\FF_2 \longrightarrow \mathrm{Div}_W(C_{K_S})_{\FF_2},\quad a\longmapsto\sum_{P_j \in \mathrm{Gal}(K_S/K) P_0} a P_j.$$

With the splittings, the estimate on the dimensions can be refined as
\begin{align*}
 \dim H^1(G,J[2](K_S))
=&\ \dim H^1(G,\mathrm{Div}_W^0(C_{K_S})_{\FF_2})-\dim H^1(G,\FF_2)\\
=&\ \dim H^1(G,\mathrm{Div}_W(C_{K_S})_{\FF_2})-2\dim H^1(G,\FF_2)\\
=&\ \sum_{j=1}^m\dim H^1(G_i,\FF_2)-2\dim H^1(G,\FF_2).
\end{align*}
The above computation further gives
\begin{align*}
&\ \dim H^1(G,J[2](K_S))\\
=&\ \sum_{j=1}^m\left(\dim O_{L_i,S}^\times/(O_{L_i,S}^\times)^2+\dim \mathrm{Cl}_S(L_i)[2]\right)-2\left(\dim O_{K,S}^\times/(O_{K,S}^\times)^2+\dim \mathrm{Cl}_S(K)[2]\right)\\
\le&\ \sum_{j=1}^m\left([L_i:K]\dim O_{K,S}^\times/(O_{K,S}^\times)^2+\dim \mathrm{Cl}_S(L_i)[2]\right)-2\left(\dim O_{K,S}^\times/(O_{K,S}^\times)^2+\dim \mathrm{Cl}_S(K)[2]\right)\\
=&\ 2g\dim O_{K,S}^\times/(O_{K,S}^\times)^2+\sum_{j=1}^m\dim \mathrm{Cl}_S(L_i)[2]-2\dim \mathrm{Cl}_S(K)[2]\\
\le&\ 2g\#S+\sum_{j=1}^m\dim \mathrm{Cl}_S(L_i)[2]-2\dim \mathrm{Cl}_S(K)[2].
\end{align*}
This proves Theorem \ref{descent}(2). 
\end{proof}

To prove Theorem \ref{descent}(3)(4), we need to bound the number of $2$-torsions of ideal class groups of a number field $K$. 
The Minkowski bound with a basic counting argument gives 
$$
\#\mathrm{Cl}(K)\leq |\Delta_{K/\QQ}|^{d}.  
$$
Here $d=[K:\QQ]$, and $\Delta_{K/\mathbb{Q}}$ is the discriminant of $K$ over $\mathbb{Q}$. 
This implies a rough bound of $\#\mathrm{Cl}(K)[2]$. 
For a more delicate bound, \cite[Theorem 1.1]{BSTTTZ20} makes use of the information ``2-torsion'' and proves
$$\#\mathrm{Cl}(K)[2]=O_{d,\epsilon}(|\Delta_{K/\QQ}|^{\frac12-\frac{1}{2d}+\epsilon}).$$
We will apply the following effective version of \cite[Theorem 1.1]{BSTTTZ20}. 

\begin{thm}\label{torsion_class_number}
Let $K$ be a number field of degree $d=[K:\mathbb{Q}]$. Then
$$\# \mathrm{Cl}(K)[2]\le 2^d d^{\frac32 d}
|\Delta_{K/\mathbb{Q}}|^{\frac 12}.
$$
\end{thm}

\begin{proof}
Our proof follows the proof of \cite[Theorem 1.1]{BSTTTZ20}, but we pay extra attention to get explicit constants. 

Define a norm $\|\cdot\|$ on $K$ by
$$\|x\|^2=\sum_{v\in M_{K,\infty}}\epsilon_v |x|_v^2, \quad x\in K.$$
Here $\epsilon_v=1$ for real places $v$; $\epsilon_v=2$ for  complex places $v$.
The norm extends to a norm on $V=K\otimes_{\mathbb{Q}}\mathbb{R}$. Denote by
$$B(r)=\{x\in V:\|x\|\le r\}$$
the ball of radius $r>0$ in $V$. Let $\mathrm{vol}(\cdot)$ be the unique Haar measure on $V$ satisfying $\mathrm{vol}(B(1))=1.$ The integer ring $O_K$ is a lattice in $V$. 
It is classical to have
$$\mathrm{vol}(V/O_K)
=\pi^{-\frac d2}\Gamma(\frac{d}{2}+1){|\Delta_{K/\mathbb{Q}}|}^{\frac12}.$$

By \cite[Theorem 1.5]{BSTTTZ20}, every ideal class in $\mathrm{Cl}(K)[2]$
is represented by an integral ideal $I$ of $O_K$ with $I^2=xO_K$ for some $x\in O_K$ satisfying that $|x|_v\le|\Delta_{K/\mathbb{Q}}|^{\frac1d}$ for all archimedean places $v$ of $K$. Such $x$ lies in $B(r_1)$ with
$$r_1={d}^\frac12|\Delta_{K/\mathbb{Q}}|^{\frac1d}.$$
It follows that
$$\#\mathrm{Cl}(K)[2]\le\#(O_K\cap B(r_1)).$$
It is reduced to estimate the right-hand side. 

Take a basis $e_0,e_1,\dots,e_{d-1}\in O_K$ of $V$ achieving the successive minima of $O_K$, that is, $e_i$ satisfies
$$\|e_i\|=\min\{\|x\|:x\in O_K\setminus\mathrm{Span}_{\mathbb{Q}}\{e_0,e_1,\dots,e_{i-1}\}\}.$$
Note that $e_0,e_1,\dots,e_{d-1}$ do not necessarily form a $\ZZ$-basis of $O_K$. 
Since for any $x\in O_K\setminus\{0\}$.
$$\|x\|^2\ge d\left(\prod_{v\in M_{K,\infty}}|x|_v^{2\epsilon_v}\right)^\frac{1}{d}=d\cdot|N_{K/\mathbb{Q}}(x)|^\frac{2}{d}\ge d=\|1\|^2,$$
we may take $e_0=1$. By the proof of \cite[Theorem 3.1]{BSTTTZ20}, there exists a permutation $\sigma$ of $\{1,\dots,d-2\}$ such that for $i=1,\dots,d-2$,
$$e_ie_{\sigma(i)}\notin\mathrm{Span}_{\mathbb{Q}}\{e_0,e_1,\dots,e_{d-2}\}.$$
So
$$\|e_{d-1}\|\le\|e_ie_{\sigma(i)}\|\le\|e_i\|\cdot\|e_{\sigma(i)}\|,$$
and thus
$$\|e_{d-1}\|^{d}\le\prod_{i=1}^{d-1}\|e_i\|^2\le\frac{4^{d}\mathrm{vol}(V/O_K)^2}{d}.$$
Here the last inequality follows from Minkowski's second theorem, which asserts that
$$\prod_{i=0}^{d-1}\|e_i\|\le2^d{\mathrm{vol}(V/O_K)}.$$
As a consequence,
$$\|e_0\|\le\|e_1\|\le\cdots\le\|e_{d-1}\|\le\frac{4\mathrm{vol}(V/O_K)^\frac{2}{d}}{d^{\frac{1}{d}}}
=4\pi^{-1}\Gamma(\frac{d}{2}+1)^{\frac 2d} d^{-\frac1d}{|\Delta_{K/\mathbb{Q}}|}^{\frac1d}.$$

By the above bound, the subset
$$D=\left\{\sum_{i=0}^{d-1}a_ie_i:-\frac12< a_i<\frac12\right\}$$
 of $V$
is contained in $B(r_2)$ with 
$$r_2=2\pi^{-1}\Gamma(\frac{d}{2}+1)^{\frac 2d} d^{1-\frac1d}{|\Delta_{K/\mathbb{Q}}|}^{\frac1d}.$$
Consider the natural map
$$\coprod_{x\in O_K\cap B(r_1)}(x+D)\longrightarrow B(r_1+r_2).$$
Note that every point in $B(r_1+r_2)$ has at most $[O_K:\sum_{i=0}^{d-1} \mathbb Ze_i]$ preimages. Comparing volumes gives
$$\#(O_K\cap B(r_1))\le\frac{\mathrm{vol}(B(r_1+r_2))}{\mathrm{vol}(D)}\cdot[O_K:\sum_{i=0}^{d-1} \mathbb Ze_i]=\frac{(r_1+r_2)^d}{\mathrm{vol}(V/O_K)}
=\frac{\left({d}^\frac12+
2\pi^{-1}\Gamma(\frac{d}{2}+1)^{\frac 2d} d^{1-\frac1d}\right)^d}{\pi^{-\frac d2}\Gamma(\frac{d}{2}+1)} {|\Delta_{K/\mathbb{Q}}|}^{\frac12}.$$
The $d$-th root of the coefficient on the right-hand side is just 
$$
\frac{ {d}^\frac12+
2\pi^{-1}\Gamma(\frac{d}{2}+1)^{\frac 2d} d^{1-\frac1d}}{\pi^{-\frac 12}\Gamma(\frac{d}{2}+1)^{\frac 1d}}  
=\frac{(d\pi)^\frac12}{A}
+\frac{2d^{1-\frac1d}}{\pi^{\frac 12} }A
$$
with 
$$
A=\Gamma(\frac{d}{2}+1)^{\frac1d}. 
$$
We need to prove that it is at most $2d^{\frac32}$.
We may assume $d\ge 2$. Then the above inequality is a direct consequence of
$$1\le A \le\big(\frac d2\big)^{\frac 12}.$$
\end{proof}

Now we are ready to prove the remaining parts of Theorem \ref{descent}, which can be viewed as effective versions of \cite[Theorem 1.3(b)]{BSTTTZ20}. 

\begin{proof}[Proof of Theorem \ref{descent}(3)(4)]
By Theorem \ref{torsion_class_number}, for any finite extension $L/K$,
$$\dim\mathrm{Cl}_S(L)[2]\le\dim\mathrm{Cl}(L)[2]\le\frac{1}{2}\log_2|\Delta_{L/\mathbb{Q}}|+\frac32{[L:\QQ]\log_2[L:\QQ]}+[L:\QQ].$$
Here the first inequality follows from the surjection $\mathrm{Cl}(L)\to \mathrm{Cl}_S(L)$. We also need the trivial bound
$$\#S\le 2d+\log_3|N_{K/\mathbb{Q}}\Delta_f|\le 2d+\frac23\log_2|N_{K/\mathbb{Q}}\Delta_f|.$$
This follows from the fact that every non-archimedean place $v\in S$ divides $2\Delta_f$. 

Now we prove part (3) of the theorem.
The point at infinity is a rational Weierstrass point. So we may assume $L_1=K$. Denote $d_i=[L_i:\mathbb{Q}]$. 
For $i=1,\dots,2g+2$, $L_i$ is isomorphic to $K[T]/(f_i(T))$ for some factor $f_i(x)\in O_K[x]$ of $f(x)$. 
As $f(x)=f_1(x)\cdots f_{2g+2}(x)$, we have $\Delta_{f_1}\cdots \Delta_{f_{2g+2}}$
divides $\Delta_f$. 
We also have $\Delta_{L_i/K}\mid\Delta_{f_i}$, and thus
$$\log_2|\Delta_{L_i/\mathbb{Q}}|
=\log_2|N_{K/\mathbb{Q}}\Delta_{L_i/K}|+\frac{d_i}{d}\log_2|\Delta_{K/\mathbb{Q}}|\le\log_2|N_{K/\mathbb{Q}}\Delta_{f_i}|+\frac{d_i}{d}\log_2|\Delta_{K/\mathbb{Q}}|.$$
We get
$$\dim\mathrm{Cl}_S(L_i)[2]\le \frac12\log_2|N_{K/\mathbb{Q}}\Delta_{f_i}|+
\frac{d_i}{2d}\log_2|\Delta_{K/\mathbb{Q}}|+\frac32d_i\log_2d_i+d_i.$$
Therefore, part (2) of the theorem implies that
\begin{align*}
\mathrm{rk}\, J(K)\le&\ 2g(2d+\frac{2}{3}\log_2|N_{K/\mathbb{Q}}\Delta_f|)
+\sum_{i=2}^m\left(\frac12\log_2|N_{K/\mathbb{Q}}\Delta_{f_i}|+
\frac{d_i}{2d}\log_2|\Delta_{K/\mathbb{Q}}|+\frac32d_i\log_2d_i+d_i\right)\\
\le&\ 4gd+(\frac43g+\frac12)\log_2|N_{K/\mathbb{Q}}\Delta_f|+(g+\frac12)\log_2|\Delta_{K/\mathbb{Q}}|+(3g+\frac32)d\log_2(2gd+d)+(2g+1)d\\
=&\ (\frac43g+\frac12)\log_2|N_{K/\mathbb{Q}}\Delta_f|+(g+\frac12)\log_2|\Delta_{K/\mathbb{Q}}|+(3g+\frac32)d\log_2(2gd+d)+(6g+1)d.
\end{align*}
This proves Theorem \ref{descent}(3). 

The proof of part (4) is similar to that of (3), except that we apply part (1) instead of part (2). Then we have
\begin{align*}
\mathrm{rk}\, J(K)
\le&\ 1+(2g+3)(2d+\frac{2}{3}\log_2|N_{K/\mathbb{Q}}\Delta_f|)+\sum_{i=1}^m
\left(\frac12\log_2|N_{K/\mathbb{Q}}\Delta_{f_i}|+\frac{d_i}{2d}\log_2|\Delta_{K/\mathbb{Q}}|+\frac32d_i\log_2d_i+d_i\right)\\
&\quad+ \frac12\log_2|\Delta_{K/\mathbb{Q}}|+\frac32d\log_2d+d\\
\leq &\ 1+(4g+6)d+(\frac43g+\frac52)\log_2|N_{K/\mathbb{Q}}\Delta_f|+(g+\frac32)\log_2|\Delta_{K/\mathbb{Q}}|+(3g+3)d\log_2(2gd+2d)\\
&\quad+\frac32d\log_2d+(2g+3)d\\
\le &\ 1+(6 g+6)d+(\frac43g+\frac52)\log_2|N_{K/\mathbb{Q}}\Delta_f|+(g+\frac32)\log_2|\Delta_{K/\mathbb{Q}}|+(3g+\frac92)d\log_2(2gd+2d).
\end{align*}
This proves Theorem \ref{descent}(4) and thus the whole theorem. 
\end{proof}

\subsection{Average bound}
\label{miscellaneous}

Let $C$ be a hyperelliptic curve of genus $g\geq 2$ over $\mathbb Q$ with a marked rational point $P\in C(\QQ)$. Let $J_C$ the Jacobian variety of $C$ over $\QQ$. If $P$ is a Weierstrass point, then the affine curve $C\setminus \{P\}$ has an equation
$$y^2=x^{2g+1}+c_2x^{2g-1}+c_3x^{2g-2}+\cdots+c_{2g+1}, \quad c_i\in \QQ.$$
The point $P$ is  infinity under this equation. Moreover, these coefficients are unique if we require the equation to be \emph{minimal} over $\mathbb Z$, i.e. $c_i\in\mathbb Z$ and there is no prime $p$ such that each $c_i$ is divided by $p^{2i}$. Assume that the equation is minimal, and define the height of $(C,P)$ as
$$h(C,P)=\log\max\{1, |c_2|^\frac{1}2,\dots, |c_{2g+1}|^\frac{1}{2g+1}\}.$$

By Bhargava--Gross \cite{Bhargava-Gross}, the average size of the $2$-Selmer groups $\mathrm{Sel}_2(J_C)$ is $3$ in that
$$\lim_{A\to\infty}\frac{\sum\limits_{h(C,P)<A}\#\mathrm{Sel}_2(J_C)}{\sum\limits_{h(C,P)<A}1}=3.$$
The Kummer map
$$J_C(\QQ)/2J_C(\QQ)\longrightarrow\mathrm{Sel}_2(J_C)$$
is injective. So the result implies
$$\limsup_{A\to\infty}\frac{\sum\limits_{h(C,P)<A} 2^{\mathrm{rk}\, J_C(\QQ)}}{\sum\limits_{h(C,P)<A}1}\le3.$$
Note that $1+\frac{5}{4\sqrt{g}}\le2$ for $g\geq 2$. 
By Theorem \ref{mordell_lang}, the average number of rational points is at most
$$\limsup_{A\to\infty}\frac{\sum\limits_{h(C,P)<A}\#C(\QQ)}{\sum\limits_{h(C,P)<A}1}\le 10^{13}g^8\cdot \limsup_{A\to\infty}\frac{\sum\limits_{h(C,P)<A}2^{\mathrm{rk}\, J_C(\QQ)}}{\sum\limits_{h(C,P)<A}1}\le3\cdot10^{13}g^8.$$

If $P$ is not a Weierstrass point, denote by $P'$ the image of $P$ under the hyperelliptic involution. 
Then the affine curve $C\setminus\{P, P'\}$ has an equation
$$y^2=x^{2g+2}+c_2x^{2g}+c_3x^{2g-1}+\cdots+c_{2g+2}, \quad c_i\in \QQ.$$
The points $P, P'$ are at infinity. 
Assume that the equation is minimal over $\mathbb Z$ as above.  
This determines $c_i\in\ZZ$ up to a sign. 
Define the height of $(C,P)$ as
$$h(C,P)=\log\max\{1, |c_2|^\frac{1}2,\dots, |c_{2g+2}|^\frac{1}{2g+2}\}.$$
In this case, Shankar--Wang \cite{SW} proved that the average size of the $2$-Selmer group is at most $6$. Then the above argument gives
$$\limsup_{A\to\infty}\frac{\sum\limits_{h(C,P)<A}\#C(\QQ)}{\sum\limits_{h(C,P)<A}1}\le
10^{13}g^8\cdot \limsup_{A\to\infty}\frac{\sum\limits_{h(C,P)<A}2^{\mathrm{rk}\, J_C(\QQ)}}{\sum\limits_{h(C,P)<A}1}\le6\cdot10^{13}g^8.$$ 

In summary, we have the following theorem. 

\begin{thm}[Theorem \ref{average intro}]
\label{average}
\begin{enumerate}
    \item For hyperelliptic curves of genus $g\geq 2$ over $\mathbb{Q}$ with a marked rational Weierstrass point, the average number of rational points is at most $3\cdot10^{13}g^8.$
    \item For hyperelliptic curves of genus $g\geq 2$ over $\mathbb{Q}$ with a marked rational non-Weierstrass point, the average of number of rational points is at most $6\cdot10^{13}g^8.$
\end{enumerate}
\end{thm}

We remark that the average relies on the order given by the height $h(C,P)$. This height can be viewed as some positive combination of a height of $C$ with a height of $P$. If $P$ is a Weierstrass point, the height of $P$ should be bounded above by a multiple of the height of $C$, and thus $h(C,P)$ should be bounded above by a multiple of the height of $C$. In this case, we might be able to convert the result in terms of the average by the height of $C$. 
If $P$ is not a Weierstrass point, the problem is much harder, since it is related to the effective Mordell conjecture.

\newpage
\part{Analytic Estimates}

This part consists of complex analytic arguments of this paper. 
In Part I, we have proved our main theorem 
(Theorem \ref{thm_quantitativeMordellsubgroup} or Theorem \ref{mordell_lang})
assuming three analytic results: 
Theorem \ref{thmfaltingselkieszhangphiinv part I},
 Theorem \ref{propestimateArakelovofdiagonal_compareTwoMetrics part I}, and Theorem \ref{thmglobalestimatezhangphiinvariant part I}. 
The main goal of this part is to prove these three results. 
In fact, their stronger froms are respectively proved in 
Theorem \ref{thmfaltingselkieszhangphiinv}, Theorem \ref{propestimateArakelovofdiagonal_compareTwoMetrics}, and Theorem \ref{thmglobalestimatezhangphiinvariant}.

\section{Preliminaries on hyperbolic geometry and Tian's peak sections}
\label{sectionlocalstructurehyperbolicsurfaces}
In this section, we collect several preliminaries needed for analyzing the geometric structure of complex curves of genus $g\geq 2$. We begin by recalling some local Riemannian and K\"ahler geometric properties of hyperbolic metrics. Then we review the basic setup of Tian’s peak section method following \cite{tg1, tg5, tg6, tg3}. See also \cite[Definition 5.1.1]{mamari1}. Following the approach in \cite[Section 4]{sxz1}, we construct peak sections via bounded linear functionals on the space of sections equipped with the $L^2$ inner product. For simplicity, we restrict our attention to complex $1$-dimensional manifolds and their cotangent bundles, while allowing general smooth sections, not necessarily holomorphic.

Let $C$ be a curve of genus $g\geq 2$ over $\mathbb{C}$, and denote by $\omega_C$ its canonical bundle. Let $\mu_{\mathrm{hyp}}$ denote the unique K\"ahler metric on $C $ with constant curvature $4\pi (1-g) $ and volume $1$, and let $\mu_{\mathrm{KE}} = 4\pi (g-1) \mu_{\mathrm{hyp}} $ denote the unique K\"ahler metric on $C $ with constant curvature $-1 $. Since $\dim C =1 $, we identify a K\"ahler metric $\mu$ with its associated volume form and measure.

\subsection{Local structure of hyperbolic metrics}
\label{subsec hyperbolic}

Away from short geodesics, the local structure of a hyperbolic metric is modeled on the Poincar\'e disk. In a neighborhood of a short geodesic, it is modeled on a hyperbolic cylinder, called a collar. We describe these two local models in detail below.

\subsubsection*{Poincar\'e models and uniformization}

Let $C$ be a compact Riemann surface of genus $g\geq 2$, equipped with a Riemannian metric $ds^2$ of constant curvature $-1$. Then the universal cover of $(C, ds^2)$ is isometric to the hyperbolic plane with constant curvature $-1$.

\begin{thm}
\label{thmanalysispartspaceform}
Let $\rho_{\widetilde{C}} : \widetilde{C} \to C$ be the universal covering map, and denote by $\widetilde{ds^2} = \rho_{\widetilde{C}}^* ds^2$ the pull-back of $ds^2$. Then the Riemannian manifold $(\widetilde{C},\widetilde{ds^2})$ is isometric to the real plane $\mathbb{R}^2$ equipped with the Riemannian metric $\widetilde{ds^2} = dr^2 + \sinh^2 r d\theta^2 $, where $(r,\theta)$ is the polar coordinate of $\mathbb{R}^2$.
\end{thm}

\begin{proof}
This follows directly from the uniqueness of simply connected space forms of a given constant curvature, up to isometry. See \cite[Example 1.4.6, Theorem 5.6.7]{pp1}.
\end{proof}

The metric introduced above comes from a warped product metric and is not naturally associated with a K\"ahler form. We now describe the K\"ahler structure on the universal covering space. Let $\mathbb{D}$ denote the unit disc
$$ \mathbb{D} = \{ z\in\mathbb{C} : |z|<1 \} ,$$
and let $\mathcal{H}$ denote the upper half plane
$$ \mathcal{H} = \{ \tau\in\mathbb{C} : \mathrm{Im} (\tau) > 0 \} .$$
The K\"ahler structure on the universal covering space of $C $ can be described as follows.

\begin{thm}[Poincar\'e model]
\label{thmpoincaremodel}
Let $\mu_{\mathrm{KE}} $ be the unique K\"ahler metric on $C $ with constant curvature $-1 $. Let $\rho_{\widetilde{C}} : \widetilde{C} \to C$ be the universal covering map, and $\widetilde{\mu_{\mathrm{KE}}} = \rho_{\widetilde{C}}^* \mu_{\mathrm{KE}}$ be the pull-back of $\mu_{\mathrm{KE}}$. Then the K\"ahler curve $(\widetilde{C},\widetilde{\mu_{\mathrm{KE}}})$ is isomorphic to each of the following two biholomorphically and isometrically equivalent models.
\begin{itemize}
    \item The Poincar\'e unit disk: $\left( \mathbb{D} , \mu_{\mathbb{D}} \right)=\left( \mathbb{D} , \frac{4\mu_{\mathrm{Euc}}}{(1-|z|^2)^2} \right) $.
    \item The Poincar\'e upper half-plane: $\left( \mathcal{H} , \mu_{\mathcal{H}} \right)=\left( \mathcal{H} , \frac{\mu_{\mathrm{Euc}}}{(\mathrm{Im}(\tau))^{2}} \right)  $.
\end{itemize}
Here $\mu_{\mathrm{Euc}} = \frac{idz\wedge d\bar{z}}{2} $ is the standard Euclidean K\"ahler metric on $\mathbb{C}$.
\end{thm}

\begin{proof}
See \cite[Subsection 4.3.3, Theorem 5.6.7]{pp1}. 
\end{proof}

Since these two Poincar\'e models are explicit, we now present an isomorphism between them. Define
$$f: \mathbb{D} \longrightarrow \mathcal{H} ,\;\;\; z\longmapsto \frac{z+i}{iz+1} .$$
One can verify that $f$ is an isomorphism from $\mathbb{D} $ onto $ \mathcal{H}$.

We next describe the hermitian metric on the canonical bundle of $C$. 
Let $\omega_{C}$ be the canonical bundle of $C$. There exists a unique hermitian metric $\Vert\cdot\Vert_{\mathrm{hyp},C}$ on $\omega_C$ such that for any open subset $U \subset C$ in the classical topology and any holomorphic $1$-form $\alpha$ on $U$,
$$i\alpha \wedge \bar{\alpha} = \frac{ \Vert \alpha \Vert_{\mathrm{hyp},C}^2}{2} \mu_{\mathrm{KE}} .$$

Let $\rho_{\mathbb{D}} : \mathbb{D} \to C$ and $\rho_{\mathcal{H}} : \mathcal{H} \to C$ be the universal covering maps from the unit disk and the upper half-plane models, respectively. Then $ \Omega_{\mathbb{D}}^1 = \rho_{\mathbb{D}}^* \omega_{C}$ and $ \Omega_{\mathcal{H}}^1 = \rho_{\mathcal{H}}^* \omega_{C}$ are the line bundles of holomorphic $1$-forms on the unit disk and the upper half-plane models, respectively.

Since $\rho_{\mathbb{D}}$ and $\rho_{\mathcal{H}}$ are universal covering maps, they induce hermitian metrics $\Vert\cdot\Vert_{\mathrm{hyp},\mathbb{D}}$ and $\Vert\cdot\Vert_{\mathrm{hyp},\mathcal{H}}$ such that for any open subset $U \subset C$ in the classical topology and any holomorphic $1$-form $\alpha$ on $U$,
$$ \Vert \alpha \Vert_{\mathrm{hyp},C}^2 (\rho_{\mathbb{D}} (z)) = \Vert \rho_{\mathbb{D}}^* \alpha \Vert_{\mathrm{hyp},\mathbb{D} }^2 (z),\;\; \Vert \alpha \Vert_{\mathrm{hyp},C}^2 (\rho_{\mathcal{H}} (\tau )) = \Vert \rho_{\mathcal{H}}^* \alpha \Vert_{\mathrm{hyp},\mathcal{H} }^2 (\tau ) ,$$
where $ z\in \rho_{\mathbb{D}}^{-1} (U) $, and $\tau\in \rho_{\mathcal{H}}^{-1} (U) $. By the definition of $\Vert \cdot \Vert_{\mathrm{hyp},C} $, we have
\begin{eqnarray*}
    idz \wedge d\bar{z} & = & \frac{ \Vert dz \Vert_{\mathrm{hyp},\mathbb{D}}^2}{2} \rho^*_{\mathbb{D}}\mu_{\mathrm{KE}} = \frac{2\Vert dz \Vert_{\mathrm{hyp},\mathbb{D}}^2 \mu_{\mathrm{Euc}}}{(1-|z|^2)^2} = \frac{i\Vert dz \Vert_{\mathrm{hyp},\mathbb{D}}^2 }{(1-|z|^2)^2} dz \wedge d\bar{z} ,\\
    id\tau \wedge d\bar{\tau} & = & \frac{ \Vert d\tau \Vert_{\mathrm{hyp},\mathcal{H}}^2}{2} \rho^*_{\mathcal{H}}\mu_{\mathrm{KE}} = \frac{\Vert d\tau \Vert_{\mathrm{hyp},\mathcal{H}}^2\mu_{\mathrm{Euc}}}{2(\mathrm{Im}(\tau))^{2}} = \frac{i \Vert d\tau \Vert_{\mathrm{hyp},\mathcal{H}}^2 }{4(\mathrm{Im}(\tau))^{2}} d\tau \wedge d\bar{\tau} .
\end{eqnarray*}
It follows that
$$ \Vert dz\Vert_{\mathrm{hyp} , \mathbb{D}} = 1- |z|^2 ,\;\; \Vert d\tau \Vert_{\mathrm{hyp} , \mathcal{H}} = 2 \mathrm{Im}(\tau) ,$$
and
$$\Vert f^* \alpha\Vert_{\mathrm{hyp} , \mathbb{D}} (z) = \Vert\alpha\Vert_{\mathrm{hyp} , \mathcal{H}} (f(z)) ,\;\;\;\forall z\in \mathbb{D}, \alpha\in\Gamma (\mathcal{H} ,\Omega^1_{\mathcal{H}}) .$$

As the hyperbolic hermitian metric is invariant under automorphisms, we shall write $\Vert\cdot\Vert_{\mathrm{hyp}}$ without specifying the space whenever the context is clear.

\begin{remark}
Let $\mu_{\mathrm{hyp}}^{\mathrm{BK}}$ be the K\"ahler metric on $C$ with constant bisectional curvature $-1$. Then $\mu_{\mathrm{hyp}}^{\mathrm{BK}} = \frac{1}{4} \mu_{\mathrm{KE}} $, and
$$\frac{i}{2}\alpha \wedge \bar{\alpha} = \Vert \alpha \Vert_{\mathrm{hyp},C}^2 \mu_{\mathrm{hyp}}^{\mathrm{BK}} .$$
See also \cite[Section 1.2]{tgbook1} for further details on bisectional curvature.
\end{remark}

\subsubsection*{Thick-thin decomposition of hyperbolic surfaces}

It is a classical result that any compact hyperbolic surface admits a thick-thin decomposition: collars are associated with sufficiently short geodesics, and the geometry of the complement is uniformly bounded. We now state the collar theorem describing this structure.

\begin{thm}[Collar Theorem]
\label{thmcollarthm}
Let $C$ be a compact Riemann surface of genus $g\geq 2$, equipped with a Riemannian metric $ds^2$ of constant curvature $-1$, and let $\gamma_1 ,\cdots ,\gamma_m $ be pairwise disjoint simple closed geodesics on $C$. Then the following hold:
\begin{enumerate}[(1)]
    \item $m\leq 3g-3$.
    \item There exist simple closed geodesics $\gamma_{m+1} ,\cdots ,\gamma_{3g-3} $ which, together with $\gamma_1 ,\cdots ,\gamma_m$, decompose $C$ into pairs of pants, namely, connected surfaces of genus zero with three boundary components.
    \item The collars
    $$ \mathscr{C} (\gamma_j ) = \{ p\in C \; |\; \mathrm{dist}_{ds^2} (p,\gamma_j ) \leq w(\gamma_j ) \} $$
    of widths
    $$ w(\gamma_j ) = \arcsinh \left( \frac{1}{\sinh  \frac{\ell (\gamma_j )}{2} } \right) $$
    are pairwise disjoint for $j=1,\cdots ,3g-3$, where $\ell ( \gamma )$ is the length of $\gamma$, and $\mathrm{dist}_{ds^2}$ is the distance on $C$ induced by the Riemannian metric $ds^2$.
    \item Each $\mathscr{C} (\gamma_j )$ is isometric to the cylinder $[-w(\gamma_j ) , w(\gamma_j ) ] \times (\mathbb{R}/\mathbb{Z}) $ with the Riemannian metric $ds^2 = dr^2 + \ell^2 (\gamma_j ) \cosh^2 r dt^2 $ .
\end{enumerate}

Moreover, let $\gamma'_1 ,\cdots ,\gamma'_{m'} $ be the set of all simple closed geodesics of length $\leq 2\arcsinh 1$ on $C$. Then the following hold:
\begin{enumerate}[(1)]
\addtocounter{enumi}{4}
    \item $m' \leq 3g-3$.
    \item $\gamma'_1 ,\cdots ,\gamma'_{m'} $ are pairwise disjoint.
    \item The injectivity radius $\mathrm{inj}_{ds^2} (x) > \arcsinh 1 $ for all $x\in C\setminus \left( \cup_{j=1}^{m'} \mathscr{C} (\gamma'_j) \right) $.
    \item The injectivity radius $\mathrm{inj}_{ds^2} (x)$ satisfies
    $$\sinh \mathrm{inj}_{ds^2} (x) = \sinh \left( \frac{\ell (\gamma'_j )}{2} \right) \cosh \mathrm{d}_{\gamma'_j} (x ) = \cosh \left( \frac{\ell (\gamma'_j )}{2} \right) \cosh \mathrm{d}_{\partial\mathscr{C} (\gamma'_j)} ( x ) -\sinh \mathrm{d}_{\partial\mathscr{C} (\gamma'_j)} (x ) ,$$ 
    for all $x\in  \mathscr{C} (\gamma'_j) $, where $\mathrm{d}_{\gamma'_j} (x) = \mathrm{dist}_{ds^2} (x,\gamma'_j ) $, and $\mathrm{d}_{\partial\mathscr{C} (\gamma'_j)} (x) = \mathrm{dist}_{ds^2} (x,\partial \mathscr{C} (\gamma'_j) ) $. 
\end{enumerate}
\end{thm}

\begin{remark}
The injectivity radius $\mathrm{inj}_{ds^2} (x)$ at a point $x\in C$ is, by definition, the supremum of radii $r>0$ such that the metric ball $B_{r}(\tilde{x})$ in the universal covering $(\widetilde{C},\widetilde{ds^2})$ of $(C,ds^2)$ projects isometrically onto its image in $C$, where $\tilde{x}$ is a lift of $x$. In our setting, since $(C,ds^2)=(C,\mu_{\mathrm{KE}})$ has negative sectional curvature, $\mathrm{inj}_{ds^2} (x)$ is equivalently given by half the length of the shortest nontrivial geodesic loop based at $x$. See also \cite[Lemma 5.7.12]{pp1} for more details.
\end{remark}

\begin{proof}
The proof is based on explicit calculations of the lengths of right-angled hexagons in the hyperbolic plane, which arise naturally from the pants decomposition of the surface. Using these formulas, one can construct the collars and verify their properties. For more details, see \cite[Theorem 4.1.1 and 4.1.6]{bu1}.
\end{proof}

As a direct consequence of the Collar Theorem, we obtain the following corollary, which concerns the relation of another closed geodesic to a short one.

\begin{cor}
\label{coranothergeodesicscrossorawayfromshortgeodesics}
Let $C$ be a compact Riemann surface of genus $g\geq 2$, equipped with a Riemannian metric $ds^2$ of constant curvature $-1$. Let $\gamma$ be a simple closed geodesic on $C$ of length $\ell (\gamma) \leq 2 \arcsinh 1$, and denote by $\mathscr{C}(\gamma)$ its associated collar. Let $\mathbf{c}:\mathbb{S}^1 \cong \mathbb{R}/\mathbb{Z} \to C$ be a simple closed geodesic distinct from $\gamma$. Then:
\begin{enumerate}[(1)]
    \item If $\mathbf{c}\cap \gamma\neq\varnothing$, then $\sinh (\frac{\ell (\gamma)}{2}) \cdot \sinh (\frac{\ell (\mathbf{c})}{2}) > 1$, and the function $\mathbf{P}_{[-w(\gamma ) , w(\gamma ) ]} \circ \mathbf{c} $ defined by composing $\mathbf{c}$ with the natural projection $\mathbf{P}_{[-w(\gamma ) , w(\gamma ) ]}:\mathscr{C}(\gamma)\to [-w(\gamma ), w(\gamma ) ]$ is strictly monotone on any connected component $\mathbf{I} \subset \mathbb{R}/\mathbb{Z}$ of $\mathbf{c}^{-1}(\mathscr{C}(\gamma))$ that contains a point of $\mathbf{c}\cap \gamma$.
    \item If $\mathbf{c}\cap \gamma=\varnothing $, then for any $x\in\mathbf{c}\cap \mathscr{C}(\gamma)$, the injectivity radius satisfies $\mathrm{inj}_{ds^2} (x) > \arcsinh \frac{\sqrt{2}}{2} $.
\end{enumerate}
\end{cor}

\begin{remark}
By the Collar Theorem, $\mathbf{I} \neq \mathbb{R}/\mathbb{Z}$; hence it can be naturally identified with a closed interval $\widetilde{\mathbf{I}}$ in $\mathbb{R}$. In addition, the same argument shows that if $\mathbf{c}:[0,1]\to\mathscr{C}(\gamma)$ is a geodesic segment such that $\mathbf{c}(\{0,1\})=\mathbf{c}([0,1])\cap\partial\mathscr{C}(\gamma)$, then the estimates in Corollary \ref{coranothergeodesicscrossorawayfromshortgeodesics} remain valid.
\end{remark}

\begin{proof}
We begin with (1). By Theorem \ref{thmcollarthm}, we have 
$$\sinh \left( \frac{\ell (\gamma)}{2} \right) \sinh \left( \frac{\ell (\mathbf{c})}{2} \right) > \sinh \left( \frac{\ell (\gamma)}{2} \right) \sinh (w(\gamma)) = 1,$$ 
see also \cite[Corollary 4.1.2]{bu1}. To prove the monotonicity, we consider the universal covering map $\rho_{\mathcal{H}}:\mathcal{H} \to C$ and the lifts $\tilde{\mathbf{c}}:\mathbb{R}\to \widetilde{C}\cong\mathcal{H}$ and $\tilde{\gamma}:\mathbb{R}\to \widetilde{C}\cong\mathcal{H}$ chosen so that $\tilde{\gamma}\cap\tilde{\mathbf{c}} \neq\varnothing $, and $\rho_{\mathcal{H}}(\tilde{\gamma}\cap\tilde{\mathbf{c}}) \in \mathbf{c}(\mathbf{I}) $. Note that $\tilde{\gamma}$ and $\tilde{\mathbf{c}}$ are geodesics in $\mathcal{H}$; in particular, $\#(\tilde{\gamma}\cap\tilde{\mathbf{c}})\leq 1$. Assume that $\tilde{\gamma}(0)=\tilde{\mathbf{c}}(0)$, and that $\mathbf{I}$ is identified with a closed interval $\widetilde{\mathbf{I}}\subset\mathbb{R}$ containing $0$. Let $\theta_{\tilde{\gamma},\tilde{\mathbf{c}}}$ denote the angle between $\tilde{\gamma}$ and $\tilde{\mathbf{c}}$ at $\tilde{\gamma}(0)$. Since $\mathbf{c}\neq \gamma$, $\sin \theta_{\tilde{\gamma},\tilde{\mathbf{c}}} \neq 0$. By the hyperbolic triangle formula (see \cite[Theorem 2.2.2]{bu1}), 
$$\sinh \dist_{\mu_{\mathcal{H}}} (\tilde{\mathbf{c}} (t) ,\tilde{\gamma}) \cdot |\sin \theta_{\tilde{\gamma},\tilde{\mathbf{c}}}| = \sinh (|t|\ell(\mathbf{c})),$$
where $\mu_{\mathcal{H}}=\rho_{\mathcal{H}}^*\mu_{\mathrm{KE}}$. By the Collar Theorem, for $t\in\widetilde{\mathbf{I}}$, we have $\sinh \dist_{\mu_{\mathrm{KE}}} (\mathbf{c} (t) ,\gamma)= \sinh \dist_{\mu_{\mathcal{H}}} (\tilde{\mathbf{c}} (t) ,\tilde{\gamma})$, and the monotonicity follows. This completes the proof of (1).

We now turn to (2). Without loss of generality, assume that $x\in \mathscr{C}(\gamma)^{\mathrm{int}}$, the interior of $\mathscr{C}(\gamma)$. Let $\mathbf{c}_x$ denote the connected component of $\mathbf{c}\cap \mathscr{C}(\gamma)$ containing $x$. Then $\mathbf{c}_x\cap\partial\mathscr{C}(\gamma)$ consists of two points, which we denote by $y$ and $y'$. Hence $\mathbf{c}_x$ is a geodesic connecting $y$ and $y'$. Since $\mathbf{c}\cap\gamma = \varnothing$, the points $y$ and $y'$ lie in the same connected component of $\partial\mathscr{C}(\gamma)$. By the description of $\mathscr{C}(\gamma)$ (Theorem \ref{thmcollarthm}), the natural projection $[-w(\gamma),w(\gamma)]\times (\mathbb{R}/\mathbb{Z})\to \mathbb{R}/\mathbb{Z}$ induces a map $\mathbf{P}_{\gamma}:\mathscr{C}(\gamma)\to\gamma $. Since $\mathbf{c}$ is simple, the restriction $\mathbf{P}_{\gamma}|_{\mathbf{c}_x}$ is injective. By applying the hyperbolic trirectangle formula (see \cite[Theorem 2.3.1 and (4.1.7)]{bu1}), we obtain
\begin{eqnarray*}
\sinh \left(\mathrm{dist}_{ds^2} (x,\gamma) \right) & \geq & \frac{\sinh (\mathrm{dist}_{ds^2} (y,\gamma))}{\sqrt{1+\sinh^2 (\frac{\ell (\mathbf{P}_{\gamma}(\mathbf{c}_x))}{2}) \cosh^2 (\mathrm{dist}_{ds^2} (y,\gamma))}} \\
& \geq & \frac{1}{\sinh \frac{\ell (\gamma)}{2} } \cdot \frac{1}{\sqrt{1+\sinh^2 (\frac{\ell (\gamma)}{2})\cdot (1+\frac{1}{\sinh^2 ( \frac{\ell (\gamma)}{2}) })}} \\
& = & \frac{1}{\sinh \frac{\ell (\gamma)}{2} } \cdot \frac{1}{\sqrt{1+\cosh^2 (\frac{\ell (\gamma)}{2}) }} .
\end{eqnarray*}

Applying the hyperbolic trirectangle formula once more,
$$ \sinh \left(\mathrm{inj}_{ds^2} (x) \right) \geq \sinh \left( \frac{\ell (\gamma)}{2} \right) \cosh \left(\mathrm{dist}_{ds^2} (x,\gamma) \right)= \frac{\cosh^2 \left( \frac{\ell (\gamma)}{2} \right)}{\sqrt{1+\cosh^2 \left( \frac{\ell (\gamma)}{2} \right)}} >\frac{\sqrt{2}}{2} . $$
This completes the proof.
\end{proof}

Let $\Gamma \leq \Aut(\mathcal{H}) \cong \mathrm{PSL}(2,\mathbb{R})$ be the deck transformation group associated with the uniformization of a compact Riemann surface $C$ of genus $g\geq 2$. Then $\mathcal{H}/\Gamma \cong C$ as a K\"ahler manifold. By the Collar Theorem above, we obtain the following description of the local structure of the orbit $\Gamma z$ for $z\in \mathcal{H}$.

\begin{lem}
\label{lemlocalstructureorbitcollartheorem}
Let $z\in\mathcal{H}$, and let $\sigma_z\in\Gamma $ be an element satisfying
$$\mathrm{dist}_{\mathcal{H}} ( z,\sigma_z z) =\inf_{\sigma\in\Gamma\setminus \left\{ \mathrm{id}_{\mathcal{H}} \right\}}  \mathrm{dist}_{\mathcal{H}} ( z,\sigma z) ,$$
where $\mathrm{dist}_{\mathcal{H}}$ is the distance on $\mathcal{H}$ induced by the K\"ahler metric $\mu_{\mathcal{H}}=\frac{idz\wedge d\bar{z}}{2(\mathrm{Im}(z))^{2}}$. For $r>0$ and $y\in\mathcal{H}$, denote $ B_{r} (y) = \left\{ \tau\in\mathcal{H} \mid \mathrm{dist}_{\mathcal{H}} (y,\tau) <r \right\} $. Then:
\begin{enumerate}[(1)]
    \item $B_{2\arcsinh1} (z)\cap \Gamma z \subset \left\{ \sigma_z^k z :\; k \in \mathbb{Z} \right\} $.
    \item If $B_{2\arcsinh1} (z)\cap \Gamma z\neq \{z\}$, then there exists a simple closed geodesic $\gamma$ of length at most $ 2\arcsinh 1$, such that $\rho_{\mathcal{H}}(z)\in\mathscr{C}(\gamma)$. Moreover, for all $k\geq 0$,
    $$\mathrm{dist}_{\mathcal{H}} ( z,\sigma_z^{k+1} z) \geq \mathrm{dist}_{\mathcal{H}} ( z,\sigma_z^k z) + \ell(\gamma) .$$
    \item Let $\rho_{\mathcal{H}} :\mathcal{H} \to C $ be the natural projection. If $\mathrm{inj}_{\mu_{\mathrm{KE}}} (\rho_{\mathcal{H}}(z)) = r_z \leq\arcsinh 1 $ and $\mathrm{dist}_{\mathcal{H}} (z,y) \leq \log 2 $, then 
    $$\# \left( B_{r_z} (\tau)\cap \Gamma y \right) \leq 2,\;\;\forall \tau \in\mathcal{H} .$$
\end{enumerate}
\end{lem}

\begin{remark}
In the unit disk model, the hyperbolic ball $B_{\arcsinh1} (0)$ of radius $\arcsinh 1$ centered at $0$ coincides with the Euclidean disk $\mathbb{D}_{\sqrt{2}-1}(0)$.
\end{remark}

\begin{proof}
We first prove (1). By Theorem \ref{thmcollarthm}, if $B_{2\arcsinh1} (z)\cap \Gamma z \neq \{z\}$, then $\rho_{\mathcal{H}} (z) \in \mathscr{C} (\gamma) $ for some collar $\mathscr{C} (\gamma) \subset C $ associated with a short simple closed geodesic $\gamma\in C$ of length $\leq 2\arcsinh 1$, and $\sigma_z$ is the element given by $\gamma$. Moreover, if a homotopic nontrivial closed loop $\mathbf{c}:\mathbb{S}^1 =\mathbb{R}/\mathbb{Z} \to C$ satisfies $\mathbf{c} (0)=\rho_{\mathcal{H}}(z)$ and $\ell (\mathbf{c})\leq 2\arcsinh 1$, then $\mathbf{c} \subset \mathscr{C} (\gamma) $. It follows that $B_{2\arcsinh1} (z)\cap \Gamma z \subset\left\{ \sigma_z^k z :\; k \in \mathbb{Z} \right\} .$

We now prove (2). The existence of $\gamma$ follows from the arguments in (1). Let $\tilde{\gamma}$ be the lift of $\gamma$ and let $w\in\tilde{\gamma} $ be the point on $\tilde{\gamma}$ closest to $z$. Then for any $k\in\mathbb{Z}$, $\mathrm{dist}_{\mathcal{H}} ( w,\sigma_z^k w) = |k|\ell (\gamma ) $. By the hyperbolic trirectangle formula \cite[Theorem 2.3.1 and (4.1.7)]{bu1},
\begin{equation*}
 \sinh \left( \frac{1}{2} \mathrm{dist}_{\mathcal{H}} ( z,\sigma_z^k z) \right) = \sinh \left( \frac{|k|\ell (\gamma )}{2} \right) \cdot \cosh \left( \mathrm{dist}_{\mathcal{H}} ( z, w) \right) ,\;\;\; \forall k\in\mathbb{Z}.
\end{equation*}
For any $k\geq 0$, define $u_k (t) = 2\arcsinh \left(t\sinh  \frac{|k|\ell (\gamma )}{2} \right) $. Then
$$ \frac{d (u_{k+1}-u_k )}{dt}(t)
= \frac{2}{\sqrt{t^2+\big( \sinh\big(\tfrac{(k+1)\ell(\gamma)}{2}\big) \big)^{-2}}}
  - \frac{2}{\sqrt{t^2+\big( \sinh\big(\tfrac{k\ell(\gamma)}{2}\big) \big)^{-2}}}
> 0 .$$
Since $\mathrm{dist}_{\mathcal{H}} ( z,\sigma_z^k z) = u_k (\cosh \left( \mathrm{dist}_{\mathcal{H}} ( z, w) \right)) $, we can conclude that
\begin{equation*}
 \mathrm{dist}_{\mathcal{H}} ( z,\sigma_z^{k+1} z) - \mathrm{dist}_{\mathcal{H}} ( z,\sigma_z^k z)  \geq u_{k+1} \left( 1 \right) - u_k \left( 1 \right) =  \ell (\gamma ) .
\end{equation*}

Finally, we prove (3). By Theorem \ref{thmcollarthm} again, if $\# \left( B_{r_z} (\tau)\cap \Gamma y \right) \geq 2 $ for some $\tau\in\mathcal{H}$, then $\mathrm{inj}_{\mu_{\mathrm{KE}}} (\rho_{\mathcal{H}}(y)) \leq r_z $. Consequently, $\rho_{\mathcal{H}} (y) \in \mathscr{C} (\gamma) $ for some collar $\mathscr{C} (\gamma) \subset C $ associated with a short simple closed geodesic $\gamma$ of length $\leq 2\arcsinh 1$. Let $\sigma_y\in\Gamma$ be the deck transformation corresponding to $\gamma$. Let $\tilde{\gamma}$ be a lift of $\gamma$ and let $w_y,w_z\in\tilde{\gamma} $ be the points on $\tilde{\gamma}$ closest to $y$ and $z$, respectively. By the hyperbolic trirectangle formula again,
\begin{equation*}
  \sinh \left( \frac{\ell (\gamma )}{2} \right) \cdot \cosh \left( \mathrm{dist}_{\mathcal{H}} ( z, w_z) \right) = \sinh \left( \frac{1}{2} \mathrm{dist}_{\mathcal{H}} ( z,\sigma_y z) \right) \geq \sinh \mathrm{inj}_{\mu_{\mathrm{KE}}} (\rho_{\mathcal{H}}(z)) = \sinh r_z ,
\end{equation*}
and for any $k\in\mathbb{Z}$,
\begin{equation*}
 \sinh \left( \frac{1}{2} \mathrm{dist}_{\mathcal{H}} ( z,\sigma_y^k z) \right) = \sinh \left( \frac{|k|\ell (\gamma )}{2} \right) \cdot \cosh \left( \mathrm{dist}_{\mathcal{H}} ( y, w_y) \right)  .
\end{equation*}
Since $\mathrm{dist}_{\mathcal{H}} (z,y) \leq \log 2 $, we have
$$ \cosh \left( \mathrm{dist}_{\mathcal{H}} ( y, w_y) \right) \geq \max\{ 1,\cosh \left( \mathrm{dist}_{\mathcal{H}} ( z, w_z)-\log 2 \right) \} > \frac{1}{2} \cosh \left( \mathrm{dist}_{\mathcal{H}} ( z, w_z) \right) .$$
Therefore, for any $k\geq 2$, 
\begin{equation*}
 \sinh \left( \frac{1}{2} \mathrm{dist}_{\mathcal{H}} ( y,\sigma_y^k y) \right) > \frac{1}{2}\sinh \left( \frac{k\ell (\gamma )}{2} \right) \cdot \cosh \left( \mathrm{dist}_{\mathcal{H}} ( z, w_z) \right) > \sinh \left( \frac{1}{2} \mathrm{dist}_{\mathcal{H}} ( z,\sigma_y z) \right) > \sinh r_z .
\end{equation*}
This shows that any ball of radius $r_z$ in $\mathcal{H}$ intersects the orbit $\Gamma y$ in at most two points, i.e.,
$$\# \left( B_{r_z} (\tau)\cap \Gamma y \right) \leq 2,\;\;\forall \tau \in\mathcal{H} .$$
This completes the proof.
\end{proof}

\subsubsection*{Metric structure on K\"ahler collars}

For convenience, we restate the metric on collars by using the complex coordinates. It is easy to see that the K\"ahler structures compatible with the Riemannian metric on $\mathscr{C} (\gamma)$ are equivalent.

\begin{lem}[K\"ahler Collar]
\label{lemkahlercollar}
Let $C$ be a compact Riemann surface of genus $g\geq 2$, equipped with a K\"ahler metric $\mu_{\mathrm{KE}}$ of constant curvature $-1$, and let $\gamma $ be a simple closed geodesic of length $\ell(\gamma)\leq 2\arcsinh 1$ on $C$. Then the collar $\mathscr{C} (\gamma)$ is isometric to $\overline{\mathbb{D}_{e^{-\nu (\gamma)}} (0) }\setminus {\mathbb{D}}_{e^{ \nu (\gamma) - \frac{2\pi^2}{\ell(\gamma)}}} (0) $ equipped with the K\"ahler metric
$$ \mu_{\mathrm{KE}} |_{\mathscr{C} (\gamma)} = \frac{ \ell^2 (\gamma) \mu_{\mathrm{Euc}} }{ 4\pi^2 |z|^2 \sin^2 \left( \frac{\ell (\gamma)}{2\pi} \log |z| \right) } ,$$
where $\mu_{\mathrm{Euc}} = \frac{idz\wedge d\bar{z}}{2} $ is the standard Euclidean metric on $ \mathbb{C} $, and $\nu (\gamma) \in (0, \frac{\pi}{2\ell (\gamma)}) $ satisfies $\nu (\gamma) \ell (\gamma) = {2\pi} \arcsin \left( \tanh \frac{\ell (\gamma)}{2} \right) $. Moreover, the image of the short simple closed geodesic $\gamma $ is the circle $|z|=e^{-\frac{\pi^2}{\ell (\gamma)}}$.
\end{lem}

\begin{proof}
The argument is the same as in \cite[Subsection 4.3.3]{pp1}. Since the hyperbolic metric on $\mathscr{C} (\gamma)$ can be represented as a warped product, we only need to find a diffeomorphism $(r ,\theta) \mapsto ( \varsigma (r) ,\theta) $ and a smooth function $\Upsilon (\varsigma)$ such that
$$ \Upsilon^2 (\varsigma) ( d\varsigma^2 +\varsigma^2 d\theta^2 ) = dr^2 + \ell^2 (\gamma ) \cosh^2 r dt^2 .$$
Note that $d\varsigma = \varsigma' (r) dr $ and $2\pi dt= d\theta$. This yields the following system of equations:
$$ \varsigma' (r) \Upsilon(\varsigma) = 1 ,\;\;\; 2\pi \varsigma (r) = \varsigma' (r) \ell (\gamma ) \cosh r . $$
One can verify that the following pair of functions is a particular solution to the system above:
$$ \varsigma (r) = e^{4\pi \cdot \frac{\arctan e^r -\pi/2}{\ell (\gamma )}} ,\;\;\; \Upsilon ( \varsigma) = \frac{-\ell (\gamma)}{ 2\pi \varsigma \sin \left( \frac{\ell (\gamma)}{2\pi } \log \varsigma \right) } .$$
Clearly, $z=\varsigma e^{i\theta}$ is just the holomorphic coordinate we need.
\end{proof}

Now we present some properties concerning the Riemannian and K\"ahler structure on the K\"ahler collar $(\mathscr{C} (\gamma) ,\mu_{\mathrm{KE}} )$.

\begin{prop}
\label{propbasicpropertiesriemanniankahlercollars}
Let $C$ be a compact Riemann surface of genus $g\geq 2$, equipped with a K\"ahler metric $\mu_{\mathrm{KE}}$ of constant curvature $-1$, and let $\gamma $ be a simple closed geodesic of length $\ell (\gamma)\leq 2\arcsinh 1$ on $C$. Let $\mathscr{C} (\gamma)$ and $\nu (\gamma)$ be defined as in Lemma \ref{lemkahlercollar}. Then the following properties hold.
\begin{enumerate}[(1)]
    \item $\frac{\pi^2}{\ell (\gamma) } - \nu (\gamma) $ and $\nu (\gamma)$ are decreasing functions of $\ell (\gamma)$. Moreover, $ \frac{\pi^2}{4\arcsinh 1} \leq \nu (\gamma) <\pi $.
    \item For any $\varsigma_0 \in (e^{ - \frac{\pi^2}{\ell (\gamma )} } , e^{-\nu (\gamma)} ) $ and $z_0 \in \overline{\mathbb{D}_{\varsigma_0} (0)} \setminus {\mathbb{D}}_{e^{ - \frac{\pi^2}{\ell (\gamma )} }} (0) $, the distance
    $$ \mathrm{dist}_{\mu_{\mathrm{KE}}} (z_0 , \partial \overline{\mathbb{D}_{\varsigma_0} (0)} ) = \log \left( \frac{ \tan \left( -\frac{\ell (\gamma) \log |z_0 | }{4\pi} \right)}{\tan \left(- \frac{\ell (\gamma) \log \varsigma_0 }{4\pi} \right)} \right) \geq \log \left( \frac{\log |z_0|}{\log \varsigma_0} \right) .$$
    \item For any $z_0 \in \overline{\mathbb{D}_{e^{-\nu (\gamma)}} (0)} \setminus {\mathbb{D}}_{e^{ - 2 \nu (\gamma) }} (0) $, the injectivity radius
    $$ \mathrm{inj}_{\mu_{\mathrm{KE}}} (z_0 ) \geq \arcsinh \left( -\frac{\nu (\gamma)}{ \log |z_0 |} \right) .$$
\end{enumerate}
\end{prop}

\begin{proof}
We first prove (1). By definition, $\nu (\gamma) = \frac{2\pi \arcsin \left( \tanh  \frac{\ell (\gamma)}{2} \right)}{\ell (\gamma)} $ is a function of $\ell (\gamma)$. To show that $\nu (\gamma)$ is decreasing, set $u_1 (t) = \frac{\arcsin (\tanh t)}{t}$, $\forall t>0$. Then
$$ \frac{d u_1 (t)}{dt} = \frac{1}{t \cosh t} - \frac{\arcsin (\tanh t)}{t^2} = \frac{ \frac{t}{\cosh t} - \arcsin (\tanh t) }{t^2} ,$$
and
$$ \frac{d}{dt} \left( \frac{t}{\cosh t} - \arcsin (\tanh t) \right) = -\frac{t\sinh t}{\cosh^2 t} <0 .$$
Hence $\frac{t}{\cosh t} - \arcsin (\tanh t) $ is strictly decreasing for $t>0$, and moreover it is negative. It follows that $\frac{du_1}{dt} < 0$ when $t>0$. Since $\nu (\gamma) = \pi u_1 (\frac{\ell (\gamma)}{2} ) $, we conclude that $\nu (\gamma )$ is strictly decreasing. Then
$$ \pi = \lim_{t\to 0^+ } \pi u_1 (t) > \pi u_1 (\frac{\ell (\gamma)}{2}) = \nu (\gamma) \geq \pi u_1 (\arcsinh 1) = \frac{\pi^2}{4\arcsinh 1} .$$
Similarly, for any $t>0$, we compute 
$$\frac{d }{dt} \left(  \frac{\pi}{2t}  -u_1 (t) \right) = -\frac{\pi}{2t^2} - \frac{ 1 }{t \cosh t} + \frac{\arcsin (\tanh t)}{t^2} < \frac{-\pi+\pi}{2t^2} =0 .$$ 
It follows that $\frac{\pi^2}{\ell (\gamma) } - \nu (\gamma) = \pi \left(  \frac{\pi}{2} \cdot \frac{1}{\frac{\ell (\gamma)}{2}} -u_1 (\frac{\ell (\gamma)}{2} ) \right) $ is also decreasing in $\ell (\gamma)$. This proves (1).

We now turn to (2). Since the metric $\mu_{\mathrm{KE}}$ is invariant under the transformations 
$z\mapsto e^{-2\pi^2 \ell^{-1} (\gamma ) } z^{-1}$ and $z\mapsto e^{i\theta } z$, $\forall \theta \in [0,2\pi )$, it follows that
\begin{eqnarray*}
\mathrm{dist}_{\mu_{\mathrm{KE}}} (z_0 , \partial \overline{\mathbb{D}_{\rho_0} (0)} ) & = & \mathrm{dist}_{\mu_{\mathrm{KE}}} (z_0 , \rho_0 |z_0|^{-1} z_0 ) .
\end{eqnarray*}

Set $\varsigma = |z|$ and $z=\varsigma e^{i\theta}$. Then
$$\mu_{\mathrm{KE}} |_{\mathscr{C} (\gamma)} = \frac{ \ell^2 (\gamma) (d\varsigma^2 + \varsigma^2 d\theta^2 ) }{ 4\pi^2 \varsigma^2 \sin^2 ( \frac{\ell (\gamma)}{2\pi} \log \varsigma ) } .$$ 
Since $\tan$ is convex and increasing on $(0,\frac{\pi}{2})$, we obtain
\begin{eqnarray*}
\mathrm{dist}_{\mu_{\mathrm{KE}}} (z_0 , \partial \overline{\mathbb{D}_{\varsigma_0} (0)} ) & = & \int_{|z_0 |}^{\varsigma_0} \frac{ -\ell (\gamma) d\varsigma }{ 2\pi \varsigma \sin ( \frac{\ell (\gamma)}{2\pi} \log \varsigma ) } = \log \left( \frac{ \tan \left( \frac{\ell (\gamma) \log |z_0 | }{-4\pi} \right)}{\tan \left( \frac{\ell (\gamma) \log \varsigma_0 }{-4\pi} \right)} \right)\\
& \geq & \log \left( \frac{\frac{\ell (\gamma) \log |z_0 | }{-4\pi} -0}{\frac{\ell (\gamma) \log \varsigma_0 }{-4\pi} -0} \right)  = \log \left( \frac{\log |z_0|}{\log \varsigma_0} \right) .
\end{eqnarray*}
This establishes (2).

Finally, we consider (3). For $t\in \left( 0,\frac{\ell (\gamma) \nu (\gamma)}{2\pi} \right)$, define 
$$u_2 (t) = \cosh t \cosh \left( \log\left(\frac{\tan \frac{\vartheta t}{2} }{\tan \frac{ t}{2}}\right) \right) - \sinh \left( \log\left(\frac{\tan \frac{\vartheta t}{2} }{\tan \frac{ t}{2}}\right) \right) ,$$
where $\vartheta = -\frac{\log |z_0|}{\nu (\gamma)} \in \left[ 1,2 \right) $. By Theorem \ref{thmcollarthm} and the previous properties, we have 
\begin{eqnarray*}
\sinh \left( \mathrm{inj}_{\mu_{\mathrm{KE}}} (z_0 ) \right) & = &  \cosh \left(\log \left( \frac{ \tan \left( \frac{\vartheta \ell (\gamma) \nu (\gamma ) }{4\pi} \right)}{\tan \left( \frac{\ell (\gamma) \nu (\gamma ) }{4\pi} \right)} \right) \right) \cosh \frac{\ell(\gamma)}{2} \\
& & -\sinh \left(\log \left( \frac{ \tan \left( \frac{\vartheta \ell (\gamma) \nu (\gamma ) }{4\pi} \right)}{\tan \left( \frac{\ell (\gamma) \nu (\gamma ) }{4\pi} \right)} \right) \right) \\
& > & u_2 \left( \frac{\ell(\gamma) \nu (\gamma)}{2\pi} \right) .
\end{eqnarray*}

Since
$$\lim_{t\to 0^+} u_2 (t) = \cosh 0 \cosh (\log (\vartheta)) - \sinh (\log (\vartheta)) = \vartheta^{-1} ,$$
it is sufficient to prove that $u_2 (t)$ is increasing for $t\in \left( 0,\frac{\ell (\gamma) \nu (\gamma)}{2\pi} \right)$. For $\vartheta \in \left[ 1,2 \right)$ and $t\in \left( 0,\frac{\ell (\gamma) \nu (\gamma)}{2\pi} \right)$, one checks that
$$ 0< \frac{\sin (\vartheta t)}{\vartheta } = t\int_{0}^{1} \cos (\vartheta t\varsigma)d\varsigma \leq t\int_{0}^{1} \cos ( t\varsigma)d\varsigma = \sin t .$$ A direct computation shows
\begin{eqnarray*}
\frac{d u_2 (t)}{dt} & = & \cosh \left( \log\left(\frac{\tan \frac{\vartheta t}{2}}{\tan\frac{t}{2}}\right) \right) \left( \sinh (t) - \frac{\vartheta}{\sin (\vartheta t)} + \frac{1}{\sin t} \right) \\
& & + \cosh (t) \sinh \left( \log\left(\frac{\tan \frac{\vartheta t}{2}}{\tan\frac{t}{2}}\right) \right) \left( \frac{\vartheta}{\sin (\vartheta t)} - \frac{1}{\sin t} \right) \\
& \geq & \cosh \left( \log\left(\frac{\tan \frac{\vartheta t}{2}}{\tan\frac{t}{2}}\right)\right) \left( \sinh (t) - \frac{\vartheta}{\sin (\vartheta t)} + \frac{1}{\sin t} \right) . 
\end{eqnarray*}

Setting $u_3 (t) = t - \frac{\vartheta}{\sin (\vartheta t)} + \frac{1}{\sin t} $, the inequality $u_3 (t)\geq 0$ implies $\frac{d u_2 (t)}{dt} \geq 0 $. By the Taylor expansion $\sin t = t-\frac{t^3}{6} + o (t^4)$ as $t\to 0^+$, we have
$$ \lim_{t\to 0^+} u_3 (t) = \lim_{t\to 0^+}\left( \frac{t-\sin t}{t\sin t} \right) - \lim_{t\to 0^+}\left( \frac{\vartheta t - \sin (\vartheta t)}{t\sin (\vartheta t)} \right) = 0 .$$
A further computation shows $\frac{d u_3 (t)}{dt} = 1+\frac{\vartheta^2 \cos (\vartheta t)}{\sin^2 (\vartheta t)} - \frac{\cos t}{\sin^2 t}$, and
\begin{eqnarray*}
\frac{d^2 u_3 (t)}{dt^2} & = &  \frac{\vartheta^3 }{\sin (\vartheta t)} - \frac{2\vartheta^3 }{\sin^3 (\vartheta t)} - \frac{1}{\sin t} + \frac{2 }{\sin^3 t} .
\end{eqnarray*}

Setting $u_4 (t) = t^3 \left( \frac{2}{\sin^3 t} - \frac{1}{\sin t} \right)$, $ t\in (0,\pi/2) $. If $\frac{\tan t}{t}\geq 2$, a direct calculation yields 
$$ \frac{d u_4 (t)}{dt} = \frac{t^3 \cos t}{\sin^2 t } \left( \frac{3\tan t}{t} - 5 + \frac{6}{\tan^2 t} \left( \frac{\tan t}{t} -1 \right) \right) > \frac{t^3 \cos t}{\sin^2 t } >0 .$$
Now assume $\frac{\tan t}{t}\in [1,2]$. In this case, the function $\frac{t^2}{\tan^2 t} \cdot \left( \frac{\tan t}{t} -1 \right) $ is increasing with respect to $\frac{\tan t}{t}$. Since $\tan t > t + \frac{t^3}{3}+ \frac{2t^5}{15}$ holds for all $\forall t\in (0,\pi/2) $, it follows that
\begin{eqnarray*}
\frac{d u_4 (t)}{dt} & = & \frac{t^3 \cos t}{\sin^2 t } \left( \frac{3\tan t}{t} - 5 + \frac{6}{ t^2} \cdot \frac{t^2}{\tan^2 t} \cdot \left( \frac{\tan t}{t} -1 \right) \right) \\
& \geq & \frac{t^3 \cos t}{\sin^2 t } \left(t^{2}+\frac{2t^{4}}{5}-2+\frac{2+\frac{4t^{2}}{5}}{\left(1+\frac{t^{2}}{3}+\frac{2t^{4}}{15}\right)^{2}}\right) \\
& = & \frac{t^3 \cos t}{\sin^2 t } \cdot \frac{t^2 (525 + 350 t^2 + 525 t^4 + 230 t^6 + 60 t^8 + 8 t^{10})}{5(15 + 5 t^2 + 2 t^4)^2} > 0 .
\end{eqnarray*}
Hence, $\frac{d^2 u_3 (t)}{dt^2} =  t^{-3} \left( u_4 (t) - u_4 (\vartheta t) \right) <0$, $\forall t\in \left( 0,\frac{\pi }{2\vartheta} \right)$. Therefore, for any $t\in \left( 0,\frac{\pi }{4} \right) $,
\begin{eqnarray*}
u_3 (t) & > & \min\left\{ \lim_{t\to 0^+} u_3 (t) , \lim_{t\to \left( \frac{\pi }{4} \right)^-} u_3 (t) \right\} = \min\left\{ 0 , \frac{\pi }{4} -\frac{\vartheta}{\sin \frac{\pi \vartheta}{4} } + \frac{1}{\sin \frac{\pi }{4}} \right\} \\
& \geq & \min\left\{ 0 , \frac{\pi }{4} -2 + \frac{1}{\sin \frac{\pi }{4}} \right\} =0 .
\end{eqnarray*}

Since $\nu (\gamma) \ell (\gamma) = {2\pi} \arcsin \left( \tanh  \frac{\ell (\gamma)}{2} \right) $, we can conclude that
$$ \frac{\ell (\gamma) \nu (\gamma) }{2\pi} = \arcsin \left( \tanh  \frac{\ell (\gamma)}{2} \right) \leq \arcsin \left( \tanh \left( \arcsinh 1 \right) \right) = \frac{\pi}{4} .$$
Then $\frac{du_2 (t)}{dt} >0 $, $\forall t\in \left( 0,\frac{\ell (\gamma) \nu (\gamma)}{2\pi} \right)$, and (3) is proven.
\end{proof}

\subsection{\texorpdfstring{$L^2$}{Lg} inner product and Hodge decomposition}

To prepare for the construction of peak sections, we recall the notion of the $L^2$ inner product on differential forms and the corresponding Hodge decomposition.

Let $C$ be a connected compact smooth complex curve equipped with the classical topology, and $\mu$ be a K\"ahler metric on $C$. Denote by $\mathscr{A}_{C, \mathbb{R}}^k $ the sheaf of smooth $k$-forms with coefficients in $\mathbb{R}$, and set $\mathscr{A}_{C, \mathbb{C}}^k = \mathscr{A}_{C, \mathbb{R}}^k \otimes_{\mathbb{R}} \mathbb{C} $, where $k\in \{0,1,2\}$. One can verify that $\mathscr{A}_{C, \mathbb{R}}^0 $ naturally forms a sheaf of rings and is isomorphic to the sheaf of smooth $\mathbb{R}$-valued functions on $C$. Consequently, for any $k\in\{0,1,2\}$ and $K\in \{\mathbb{R},\mathbb{C}\}$, $\mathscr{A}_{C, K}^k $ inherits a natural $\mathscr{A}_{C, \mathbb{R}}^0 $-module structure. Hence the sheaf of smooth differential forms with coefficients in $K$, $\mathscr{A}_{C, K} = \bigoplus_{k=0}^{2} \mathscr{A}_{C, K}^k $, is also an $\mathscr{A}_{C, \mathbb{R}}^0 $-module. Furthermore, for any classical open subset $U\subset C$, the canonical map 
$$\alpha\mapsto \alpha \otimes_{\mathbb{R}} 1 ,\;\forall \alpha \in \mathscr{A}_{C, \mathbb{R}} (U) ,$$ 
defines an $\mathscr{A}_{C, \mathbb{R}}^0 (U) $-homomorphism $\mathscr{A}_{C, \mathbb{R}} \to \mathscr{A}_{C, \mathbb{C}}  $. In this way, the sheaf $\mathscr{A}_{C, \mathbb{R}} $ naturally regarded as a $\mathscr{A}_{C, \mathbb{R}}^0 $-submodule of $\mathscr{A}_{C, \mathbb{C}} $.

We now introduce the Hodge $\star $-operator associated with the K\"ahler metric $\mu$ on $C$. Let $x\in C$, and let $U_x $ be a classical open neighborhood of $x$ with local coordinate $z:U_x\to \mathbb{C}$. Denote by $u\in  \mathscr{A}_{C, \mathbb{C}}^0 (U_x ) $ a smooth $\mathbb{C}$-valued function on $U_x $. Then the Hodge $\star$-operator 
$$\star \in \mathrm{Hom}_{\mathscr{A}_{C, \mathbb{R}}^0} \left(  \mathscr{A}_{C, \mathbb{C}} , \mathscr{A}_{C, \mathbb{C}} \right) $$ 
is an isomorphism, and locally it is determined by the following relations
$$ \star(u)=\bar{u}\mu ,\;\; \star(u\mu )=\bar{u},\;\; \star(udz)=i\bar{u}d\bar{z},\;\; \star(ud\bar{z}) = -i\bar{u}dz ,$$
where $\mu$ is the K\"ahler form on $C$. One can verify that the restrictions
$$\star|_{\mathscr{A}_{C, \mathbb{R}}} \in \mathrm{Hom}_{\mathscr{A}_{C, \mathbb{R}}^0} \left(  \mathscr{A}_{C, \mathbb{R}} , \mathscr{A}_{C, \mathbb{R}} \right),\;\; \star^2|_{\mathscr{A}^k_{C, \mathbb{C}}} = (-1)^{k} \mathrm{Id}_{\mathscr{A}^k_{C, \mathbb{C}}} , $$
where $\mathrm{Id}_{\mathscr{A}^k_{C, \mathbb{C}}} :\mathscr{A}^k_{C, \mathbb{C}} \to \mathscr{A}^k_{C, \mathbb{C}} $ is the identity morphism.  See also \cite[Section 1.2]{dh1} for more details about Hodge $\star$-operators.

In particular, the restriction $\star|_{\mathscr{A}^k_{C, \mathbb{C}}} : \mathscr{A}^k_{C, \mathbb{C}} \to \mathscr{A}^{1-k}_{C, \mathbb{C}} $ is conjugate-linear, and for $k=1$, the composition of the Hodge $\star$-operator with complex conjugation defines a $\mathbb{C}$-linear operator
$$ \bar{\star} :\mathscr{A}^1_{C, \mathbb{C}} (C) \to \mathscr{A}^1_{C, \mathbb{C}} (C) ,\;\; \alpha\mapsto \overline{\star\alpha} , $$
with eigenvalues $-i$ and $i$. The corresponding eigenspaces yield an $L^2$-orthogonal decomposition:
\begin{equation*}
  \mathscr{A}^1_{C, \mathbb{C}} (C) = \mathscr{A}^{1,0}_{C, \mathbb{C}} (C) \oplus \mathscr{A}^{0,1}_{C, \mathbb{C}} (C) .
\end{equation*}
For any $\alpha \in \mathscr{A}^1_{C, \mathbb{C}}(C)$, let $\alpha^{1,0}$ and $\alpha^{0,1}$ denote the projections of $\alpha$ onto $\mathscr{A}^{1,0}_{C, \mathbb{C}}(C)$ and $\mathscr{A}^{0,1}_{C, \mathbb{C}}(C)$, respectively. The differential operators $\partial$ and $\bar{\partial}$ are then defined by
$$ \partial u= (du)^{1,0} ,\; \bar{\partial} u = (du)^{0,1} ,\; \partial \alpha = d\alpha^{0,1} ,\; \bar{\partial} \alpha = d\alpha^{1,0} , $$
for any $u\in \mathscr{A}^0_{C, \mathbb{C}} (C) $, and $ \alpha\in \mathscr{A}^1_{C, \mathbb{C}} (C)$.

For any open subset $U\subset C$ (in the classical topology) and any $\alpha_U \in \mathscr{A}_{C, \mathbb{C}} (U) $, there exists a smooth nonnegative function $\eta_{\alpha_{U}} : U\to [0,\infty) \subset \mathbb{R} $ such that 
$$\alpha_U \wedge \star \alpha_U = \eta_{\alpha_{U}} \mu \geq 0.$$
This allows us to define a Hermitian inner product on $\mathscr{A}_{C, \mathbb{C}} (U)$ by
$$ \langle \alpha_U ,\beta_U \rangle_{L^2 ;U} = \frac{1}{2} \int_{U} \alpha_U \wedge \star \beta_U ,\;\; \forall  \alpha_U,\beta_U \in \mathscr{A}_{C, \mathbb{C}} (U) .$$
We denote by $\Vert\cdot\Vert_{L^2 ;U}$ the associated $L^2$-norm. Clearly, $\langle \alpha_U ,\beta_U \rangle_{L^2 ,U} \in\mathbb{R} $ if $\alpha_U,\beta_U \in \mathscr{A}_{C, \mathbb{R}} (U) $, so the restriction of the Hermitian inner product to $\mathscr{A}_{C, \mathbb{R}} (U) \times \mathscr{A}_{C, \mathbb{R}} (U) $ is a real inner product.

Next, we recall the Hodge decomposition. As before, we restrict to the case of Riemann surfaces for simplicity.

\begin{thm}[Hodge decomposition]
\label{thmhodgedecomposition}
Let $C$ be a connected compact smooth complex curve equipped with the classical topology, and $\mu$ be a K\"ahler metric on $C$. Then for $k\in \{ 0,1,2 \} $ and $K \in \{ \mathbb{R} ,\mathbb{C} \} $, the space of smooth $k$-forms with coefficients in $K$, $\mathscr{A}^k_{C, K} (C) $, has the following $L^2$-orthogonal decomposition:
\begin{equation*}
 \mathscr{A}^k_{C, K} (C) = d\mathscr{A}^{k-1}_{C, K} (C) \oplus \mathscr{H}^k_{C, K} (C) \oplus \star d\mathscr{A}^{1-k}_{C, K} (C) ,
\end{equation*}
where $\mathscr{A}^j_{C, K} =0$ for $j\notin \{0,1,2\}$, and $\mathscr{H}^k_{C, K}$ denotes the sheaf of harmonic $k$-forms with coefficients in $K$ defined by
$$ \mathscr{H}^k_{C, K} (U) = \left\{ \alpha\in \mathscr{A}^k_{C, K} (U) :d\alpha =d\star \alpha =0 \right\} ,$$
for any classical open subset $U\subset C$. In particular, $\dim_K \mathscr{H}^k_{C, K} (C) <\infty $. Moreover, in the case $K=\mathbb{C}$ and $k=1$, the composition of the Hodge $\star$-operator with complex conjugation defines a $\mathbb{C}$-linear operator
$$ \bar{\star} :\mathscr{H}^1_{C, \mathbb{C}} (C) \to \mathscr{H}^1_{C, \mathbb{C}} (C) ,\;\; \alpha\mapsto \overline{\star\alpha} , $$
with eigenvalues $-i$ and $i$. The corresponding eigenspaces yield an $L^2$-orthogonal decomposition:
\begin{equation*}
  \mathscr{H}^1_{C, \mathbb{C}} (C) = \Gamma (C,\omega_C) \oplus \overline{\Gamma (C,\omega_C)} ,
\end{equation*}
where 
$$\Gamma (C,\omega_C)  = \mathscr{H}^1_{C, \mathbb{C}} (C) \cap \mathscr{A}^{1,0}_{C, \mathbb{C}} (C) = \{ \alpha \in \mathscr{H}^1_{C, \mathbb{C}} (C) : \bar{\star} (\alpha) = -i\alpha \} $$ 
is the space of holomorphic $1$-forms on $C$.
\end{thm}

\begin{proof}
See \cite[Theorem 6.8]{fwar1} and \cite[Theorem 3.2.8]{dh1}.
\end{proof}

Let $\mathscr{Z}_{C,K}^k$ denote the sheaf of closed $k$-forms with coefficients in $K$, defined by
$$ \mathscr{Z}^k_{C, K} (U) = \left\{ \alpha\in \mathscr{A}^k_{C, K} (U) :d\alpha  =0 \right\} ,$$
for any classical open subset $U\subset C$ and $K\in\{\mathbb{R},\mathbb{C}\} $. The de Rham cohomology with coefficients in $K$ is then given by
$$ H^k_{\mathrm{dR}} (U;K) = \mathscr{Z}^k_{C, K} (U) /d\mathscr{A}^{k-1}_{C, K} (U) .$$

As a consequence of Theorem \ref{thmhodgedecomposition}, the space of harmonic $k$-forms is naturally isomorphic to the $k$-th de Rham cohomology group.

\begin{cor}
\label{corisomorphismharmonicdhcohomology}
Let $C$ be a connected compact smooth complex curve equipped with the classical topology, and $\mu$ be a K\"ahler metric on $C$. Then the embedding $\mathscr{H}^k_{C,K} (C) \to \mathscr{Z}_{C,K}^k (C)$ gives an isomorphism $\mathscr{H}^k_{C,K} (C) \cong H^k_{\mathrm{dR}} (C;K) $.
\end{cor}

\begin{proof}
For any $\alpha \in \mathscr{Z}^k_{C, K} (C) $ and $\beta \in \mathscr{A}^{1-k}_{C, K} (C) $, we have
\begin{equation*}
    \langle \alpha ,\star d \beta \rangle_{L^2 ;C} = \int_{C} \alpha \wedge \star^2 d\beta = (-1)^k \int_{C} \alpha \wedge d\beta = - \int_{C} d \alpha \wedge \beta =0 .
\end{equation*}
Hence $\mathscr{Z}^k_{C, K} (C)$ and $\star d \mathscr{A}^{1-k}_{C, K} (C) $ are $L^2$-orthogonal. By Theorem \ref{thmhodgedecomposition}, we have
\begin{equation*}
 \mathscr{Z}^k_{C, K} (C) = d\mathscr{A}^{k-1}_{C, K} (C) \oplus \mathscr{H}^k_{C, K} (C). 
\end{equation*}
It follows that $\mathscr{H}^k_{C,K} (C) \cong \mathscr{Z}^k_{C, K} (C) /d\mathscr{A}^{k-1}_{C, K} (C) = H^k_{\mathrm{dR}}  (C;K) $.
\end{proof}

Now we focus on the case $k = 1$. By the definition of the Hodge $\star$-operator, the restriction $\star|_{\mathscr{A}^1_{C,\mathbb{C}}(C)}$ depends only on the complex structure of $C$, and not on the particular choice of K\"ahler metric. As a real vector space, $\mathscr{H}^1_{C,\mathbb{C}}(C)$ admits two natural $\mathbb{R}$-linear decompositions. One is given by separating real and imaginary parts:
$$\mathscr{H}^1_{C,\mathbb{C}}(C) = \mathscr{H}^1_{C,\mathbb{R}}(C) \oplus i\mathscr{H}^1_{C,\mathbb{R}}(C),$$
and the other is the decomposition into eigenspaces of the complex-linear involution $\bar{\star}$, with eigenvalues $-i$ and $i$:
$$\mathscr{H}^1_{C,\mathbb{C}}(C) = \Gamma(C, \omega_C) \oplus \overline{\Gamma(C, \omega_C)}. $$
In the following, we compare these two decompositions and clarify their relation.

\begin{cor}
\label{cortwodecompositionharmonic1form}
Let $C$ be a connected compact smooth complex curve equipped with the classical topology, and $\mu$ be a K\"ahler metric on $C$. Then the $\mathbb{R}$-linear map 
\begin{eqnarray*}
    \mathscr{H}^1_{C,\mathbb{R}}(C) & \longrightarrow & \;\;\;\; \mathscr{H}^1_{C,\mathbb{C}}(C) ,\\
    \alpha \;\;\;\;\; &\longmapsto& \;\;\;\; \alpha +i\star\alpha ,
\end{eqnarray*}
gives an isomorphism $\mathscr{H}^1_{C,\mathbb{R}} (C) \cong \Gamma(C, \omega_C) $ between $\mathbb{R}$-vector spaces.
\end{cor}

\begin{proof}
Let $\alpha \in \mathscr{H}^1_{C,\mathbb{R}}(C)$. Since $\star$ preserves harmonic forms and $\alpha$ is real-valued, the form $\alpha + i \star \alpha$ lies in $\mathscr{H}^1_{C,\mathbb{C}}(C)$. Moreover, we compute:
$$ \bar{\star} ( \alpha +i\star\alpha ) = \bar{\star} \alpha +i \bar{\star} \star\alpha = \star \alpha + i\star^2 \alpha = \star \alpha -i\alpha =-i ( \alpha +i\star\alpha ) , $$
so $\alpha + i \star \alpha$ is an eigenvector of $\bar{\star}$ with eigenvalue $-i$, hence lies in $\Gamma(C, \omega_C)$. This map is clearly $\mathbb{R}$-linear and injective. To show surjectivity, we compare dimensions:
$$\dim_{\mathbb{R}} \mathscr{H}^1_{C,\mathbb{R}} (C) = \dim_{\mathbb{C}} \mathscr{H}^1_{C,\mathbb{C}} (C) = 2 \dim_{\mathbb{C}} \Gamma(C, \omega_C) = \dim_{\mathbb{R}} \Gamma(C, \omega_C) <\infty .$$
Hence the map gives an isomorphism $\mathscr{H}^1_{C,\mathbb{R}} (C) \cong \Gamma(C, \omega_C) $.
\end{proof}

\subsection{Holomorphic \texorpdfstring{$1$}{Lg}-forms on K\"ahler collars}

Motivated by the local $L^2$ analysis developed for complex hyperbolic cusps in \cite{aumamar1, aumamar2, sxz1}, we now carry out a similar study for K\"ahler collars. In this subsection, we give an explicit computation of the local $L^2$ inner product on a K\"ahler collar, which serves as an important example for later use. Such a local description of $L^2$-integrable holomorphic $1$-forms is essential for the geometric localization arguments.

Let $\mathscr{C} (\gamma)^{\mathrm{int}}$ denote the interior of $\mathscr{C} (\gamma)$, and let $\Omega^{1}_{\mathscr{C} (\gamma)^{\mathrm{int}}} $ denote the line bundle of holomorphic $1$-forms on $\mathscr{C} (\gamma)^{\mathrm{int}}$. We next construct an $L^2$-orthonormal basis of $\Gamma^{\mathrm{hol}}_{L^2}   (\mathscr{C} (\gamma)^{\mathrm{int}} ,\Omega^{1}_{\mathscr{C} (\gamma)^{\mathrm{int}}})$, the Hilbert space of $L^2$-integrable holomorphic $1$-forms on $\mathscr{C} (\gamma)^{\mathrm{int}}$, equipped with the local $L^2$-inner product
$$ \langle\alpha ,\beta\rangle_{L^2 ;\mathscr{C} (\gamma)^{\mathrm{int}}} = \frac{i}{2} \int_{\mathscr{C} (\gamma)^{\mathrm{int}} } \alpha \wedge \bar{\beta} .$$

\begin{lem}
\label{lemlocalorthonormalbasis}
Let $C$ be a compact Riemann surface of genus $g\geq 2$, equipped with a K\"ahler metric $\mu_{\mathrm{KE}}$ of constant curvature $-1$, and let $\gamma $ be a simple closed geodesic of length $\ell (\gamma)\leq 2\arcsinh 1$ on $C$. Let $\mathscr{C} (\gamma)$ and $\nu (\gamma)$ be defined as in Lemma \ref{lemkahlercollar}. Then $\{ a_k z^k \}_{k\in\mathbb{Z}}$ forms an $L^2$ orthonormal basis of $\Gamma^{\mathrm{hol}}_{L^2}   (\mathscr{C} (\gamma)^{\mathrm{int}} ,\Omega^{1}_{\mathscr{C} (\gamma)^{\mathrm{int}}})$, where
\begin{eqnarray*}
|a_k |^2 = \left\{
\begin{aligned}
\frac{1}{4\pi \left( \frac{\pi^2}{\ell (\gamma)} -\nu (\gamma) \right)}, \;\;\;\;\;\; & \; k=-1 ;\\
 \frac{(k+1)e^{2(k+1) \nu (\gamma ) }}{\pi\left( 1 - e^{4(k+1) (\nu (\gamma ) - \frac{\pi^2}{\ell (\gamma)} ) } \right)}, \; & \; k\neq -1 . 
\end{aligned}
\right.
\end{eqnarray*}
\end{lem}

\begin{proof}
From the Laurent series expansion, it follows that
$$\Gamma^{\mathrm{hol}}_{L^2}   (\mathscr{C} (\gamma)^{\mathrm{int}} ,\Omega^{1}_{\mathscr{C} (\gamma)^{\mathrm{int}}}) = \overline{\oplus_{k\in\mathbb{Z}} \mathbb{C} z^k dz } ,$$ 
where the closure is taken with respect to the $L^2$-inner product defined above. A straightforward computation yields
\begin{eqnarray*}
   \langle z^k dz ,z^k dz\rangle_{L^2 ;\mathscr{C} (\gamma)^{\mathrm{int}}} & = & \frac{i}{2} \int_{\mathscr{C} (\gamma) } |z|^{2k} dz\wedge d\bar{z} = 2\pi \int_{e^{\nu (\gamma ) -\frac{2\pi^2}{\ell (\gamma)} }}^{e^{-\nu (\gamma )}} t^{2k+1} dt  \\
  & = & \left\{
\begin{aligned}
\frac{4\pi^3 }{\ell (\gamma ) } - 4\pi \nu (\gamma ), \;\;\;\;\;\;\;\;\;\;\;\;\;\;\;\;\;\;\;\;\;\; & \; k=-1 ;\\
{ \frac{ \pi}{k+1} } \left(e^{-2(k+1) \nu (\gamma ) } - e^{2(k+1) (\nu (\gamma ) - \frac{2\pi^2}{\ell (\gamma)} ) } \right), \; & \; k\neq -1 , 
\end{aligned}
\right.
\end{eqnarray*}
and the lemma follows.
\end{proof}

We now present some explicit estimates for $ \Gamma^{\mathrm{hol}}_{L^2}   (\mathscr{C} (\gamma)^{\mathrm{int}} ,\Omega^{1}_{\mathscr{C} (\gamma)^{\mathrm{int}}}) $.

\begin{prop}
\label{propexplicitestimateskahlercollar}
Let $C$ be a compact Riemann surface of genus $g\geq 2$, equipped with a K\"ahler metric $\mu_{\mathrm{KE}}$ of constant curvature $-1$, and let $\gamma $ be a simple closed geodesic of length $\ell (\gamma)\leq 2\arcsinh 1$ on $C$. Let $\mathscr{C} (\gamma)$ and $\nu (\gamma)$ be defined as in Lemma \ref{lemkahlercollar}, and let $\Gamma^{\mathrm{hol}}_{L^2} (\mathscr{C} (\gamma)^{\mathrm{int}} ,\Omega^{1}_{\mathscr{C} (\gamma)^{\mathrm{int}}})$ be the Hilbert space defined above. 

Let $\alpha\in \Gamma^{\mathrm{hol}}_{L^2}   (\mathscr{C} (\gamma)^{\mathrm{int}} ,\Omega^{1}_{\mathscr{C} (\gamma)^{\mathrm{int}}}) $ satisfy $\Vert\alpha \Vert_{L^2 ; \mathscr{C} (\gamma)^{\mathrm{int}} } = 1 $. Then the following properties hold.
\begin{enumerate}[(1)]
    \item Assume that $\int_\gamma \alpha =0$. Then for any $z \in \mathscr{C} (\gamma)^{\mathrm{int}} $, we have the pointwise estimate
    $$ i\alpha \wedge \bar{\alpha} \leq \frac{2e^{2\nu (\gamma)}}{\pi (1-e^{-3\pi }) } \left( \frac{1}{\left( 1-|z|^2 e^{2\nu (\gamma)} \right)^2} + \frac{|z|^{-4} e^{-\frac{4\pi^2}{\ell (\gamma)}} }{\left(1-|z|^{-2} e^{2\left( \nu (\gamma) - \frac{2\pi^2}{\ell (\gamma)} \right)} \right)^2} \right) \mu_{\mathrm{Euc}} .$$
    \item Assume that $\int_\gamma \alpha =0$. Then for any two points $z_0 ,z_1\in \mathscr{C} (\gamma)^{\mathrm{int}} $ and any piecewise smooth curve $\mathbf{c}: [0,1] \to \mathscr{C} (\gamma)^{\mathrm{int}} $ with $\mathbf{c} (0) = z_0 $ and $\mathbf{c} (1) =z_1 $, the integral $\int_{\mathbf{c}} \alpha = \int_{0}^{1} \mathbf{c}^*\alpha $ depends only on $z_0$ and $z_1$. 
    
    Furthermore, if $|z_0| \leq |z_1|$, we have
    $$ \left| \int_{\mathbf{c}} \alpha \right|^2 \leq \frac{e^{2\nu (\gamma)}}{\pi (1-e^{-3\pi }) } \left( \frac{|z_1 - z_0|^2}{\left( 1-|z_1|^2 e^{2\nu (\gamma)} \right)^2} + \frac{|z_1^{-1} - z_0^{-1} |^{2} e^{-\frac{4\pi^2}{\ell (\gamma)}} }{\left(1-|z_0|^{-2} e^{2\left( \nu (\gamma) - \frac{2\pi^2}{\ell (\gamma)} \right)} \right)^2} \right) $$
    \item For any two points $z_0 ,z_1\in \mathscr{C} (\gamma)^{\mathrm{int}} $, there exists a piecewise smooth curve $\mathbf{c}: [0,1] \to \mathscr{C} (\gamma)^{\mathrm{int}} $ such that $\mathbf{c} (0) = z_0 $, $\mathbf{c} (1) =z_1 $, and
    \begin{eqnarray*}
  \left| \int_{\mathbf{c}} \re(\alpha) \right| & \leq & \left(  \frac{-\log \left( \left( 1-e^{2\nu (\gamma) } |z_1 |^2 \right) \left( 1-e^{2 ( \nu (\gamma ) - \frac{2\pi^2}{\ell (\gamma)} ) } |z_0 |^{-2} \right) \right)}{\pi (1-e^{-3\pi} ) \left( 1+\left|\frac{z_0}{z_1}\right| \right)^{-2}} \right)^{\frac{1}{2}} \\
  & & + \frac{1}{2\pi}\left| \int_{\gamma} \im( \alpha) \right| \left| \log \left|\frac{z_1}{z_0} \right| \right| + \frac{1}{2 } \left| \int_{\gamma} \re( \alpha) \right| ,
\end{eqnarray*}
where $\alpha = \re(\alpha) + i\im (\alpha) $ is the decomposition of $\alpha$ into its real part and imagine part.
\end{enumerate}
\end{prop}

\begin{proof}
Let $\{a_j z^j dz\}_{j\in\mathbb{Z}}$ be the $L^2$-orthonormal basis of $\Gamma^{\mathrm{hol}}_{L^2}   (\mathscr{C} (\gamma)^{\mathrm{int}} ,\Omega^{1}_{\mathscr{C} (\gamma)^{\mathrm{int}}})$ constructed in Lemma \ref{lemlocalorthonormalbasis}. By the Laurent series expansion, there exist constants $b_j\in\mathbb{C}$ such that $\sum\limits_{j\in\mathbb{Z}} |b_j |^2 =1 $, and $\alpha = \sum\limits_{j\in\mathbb{Z}} a_j b_j z^j dz $.

We now consider estimate (1). Since $b_{-1} = 0$, the Cauchy–Schwarz inequality yields
\begin{eqnarray*}
    i\alpha\wedge \bar{\alpha} & = & 2 \left| \sum_{j\in\mathbb{Z}} a_j b_j z^j \right|^2 \mu_{\mathrm{Euc}} \leq 2\left( \sum_{j\neq-1} \left| b_j \right|^2 \right) \left( \sum_{j\neq -1} \left| a_j z^j \right|^2 \right) \mu_{\mathrm{Euc}} \\
    & = & 2 \left( \sum_{j=0}^\infty \left| a_j\right|^2 |z|^{2j}  + \sum_{j=0}^\infty \left| a_{-2-j} \right|^2 |z|^{-4-2j} \right) \mu_{\mathrm{Euc}} .
\end{eqnarray*}

Since Proposition \ref{propbasicpropertiesriemanniankahlercollars} implies $\nu(\gamma) < \pi$, we see that
$$ \frac{4}{\pi}\left( \frac{\pi^2}{\ell (\gamma)} - \nu (\gamma ) \right) > 4 \left( \frac{\pi}{2\arcsinh 1} -1 \right) \thickapprox 3.1288559\cdots > 3 . $$
By Lemma \ref{lemlocalorthonormalbasis}, it follows that
\begin{eqnarray*}
    \sum_{j=0}^\infty \left| a_j\right|^2 |z|^{2j} & = & \sum_{j=0}^\infty \frac{(j+1) e^{2(j+1) \nu (\gamma)} |z|^{2j}}{\pi \left( 1 - e^{4(j+1) (\nu (\gamma ) - \frac{\pi^2}{\ell (\gamma)} ) } \right)}  \\
    & < & \frac{e^{2\nu (\gamma)}}{\pi \left( 1 - e^{4 (\nu (\gamma ) - \frac{\pi^2}{\ell (\gamma)} ) } \right) } \sum_{j=0}^\infty (j+1) \left( e^{2 \nu (\gamma)} |z|^{2} \right)^j \\
    & \leq & \frac{e^{2\nu (\gamma)}}{\pi \left( 1 - e^{-3\pi } \right) \left( 1 - e^{2 \nu (\gamma)} |z|^{2} \right)^2 } ,
\end{eqnarray*}
and similarly,
\begin{equation*}
    \sum_{j=0}^\infty \left| a_{-2-j} \right|^2 |z|^{-4-2j} \leq \frac{ |z|^{-4} e^{2 (\nu (\gamma ) - \frac{2\pi^2}{\ell (\gamma)} )}}{\pi \left( 1 - e^{-3\pi } \right) \left( 1 - e^{2 (\nu (\gamma ) - \frac{2\pi^2}{\ell (\gamma)} )} |z|^{-2} \right)^2 } .
\end{equation*}
This completes the proof of (1).

We next turn to estimate (2). The invariance of the integral $\int_{\mathbf{c}} \alpha $ follows from the Poincar\'e lemma \cite[Proposition 4.7]{bot1} and Poincar\'e duality \cite[(5.4)]{bot1} for de Rham cohomology. Since our setting is explicit, however, we provide a constructive proof and estimate its value.

Set $u (z) = \sum\limits_{j\neq -1} \frac{a_j b_j z^{j+1}}{j+1} $. Then $u$ is holomorphic on $\mathscr{C} (\gamma)^{\mathrm{int}}$, and satisfies 
$$ du = \sum_{j\in\mathbb{Z}} a_j b_j z^j dz =\alpha .$$
For any piecewise smooth real curve $\mathbf{c}: [0,1] \to \mathscr{C} (\gamma)^{\mathrm{int}} $ with $\mathbf{c} (0) = z_0$ and $\mathbf{c} (1) =z_1$, we have $ \int_{\mathbf{c}} \alpha = u(z_1) -u(z_0) $, so the integral depends only on $z_0$ and $z_1$. By the Cauchy-Schwarz inequality, 
\begin{eqnarray*}
  \left| \int_{\mathbf{c}} \alpha \right|^2 & = & \left| \sum_{j\neq -1} \frac{a_j b_j (z_1^{j+1} - z_0^{j+1} ) }{j+1} \right|^2 \leq \left( \sum_{j\neq -1} |b_j|^2 \right) \left( \sum_{j\neq -1} \frac{|a_j|^2 | z_1^{j+1} - z_0^{j+1} |^2 }{(j+1)^2} \right) \\
  & \leq & \sum_{j=0}^{\infty} \frac{|a_j|^2 | z_1^{j+1} - z_0^{j+1} |^2 }{(j+1)^2} + \sum_{j=0}^{\infty} \frac{|a_{-2-j}|^2 | z_1^{-j-1} - z_0^{-j-1} |^2 }{(j+1)^2} .
\end{eqnarray*}
As in the argument for (i), Lemma \ref{lemlocalorthonormalbasis} yields:
\begin{eqnarray*}
  \sum_{j=0}^{\infty} \frac{|a_j|^2 | z_1^{j+1} - z_0^{j+1} |^2 }{(j+1)^2} & = & |z_1 - z_0|^2 \sum_{j=0}^{\infty} \frac{|a_j|^2 | \sum_{k=0}^j z_1^{k} z_0^{j-k} |^2 }{(j+1)^2} \\
  & \leq & |z_1 - z_0|^2 \sum_{j=0}^{\infty} |a_j|^2 |z_1|^{2j}  \\
  & < & \frac{ |z_1 - z_0|^2 e^{2\nu (\gamma)}}{\pi \left( 1 - e^{-3\pi } \right) \left( 1 - e^{2 \nu (\gamma)} |z_1|^{2} \right)^2 } ,
\end{eqnarray*}
and similarly,
\begin{equation*}
  \sum_{j=0}^{\infty} \frac{|a_{-2-j}|^2 | z_1^{-j-1} - z_0^{-j-1} |^2 }{(j+1)^2} < \frac{\left| z_0^{-1} -z_1^{-1} \right|^2 e^{2 (\nu (\gamma ) - \frac{2\pi^2}{\ell (\gamma)} )} }{\pi (1-e^{-3\pi} ) \left( 1 - e^{2 (\nu (\gamma ) - \frac{2\pi^2}{\ell (\gamma)} )} |z_0|^{-2} \right)^2} , 
\end{equation*}
which gives (2).

Finally, we address estimate (3). Without loss of generality, assume $z_0 =|z_0|$, and $z_1 = |z_1|e^{i\theta_1}$ with $\theta_1\in [-\pi ,\pi]$. Let $\mathbf{c} :[0,1] \to \mathscr{C}(\gamma)^{\mathrm{int}}$ be the piecewise smooth curve given by
\begin{equation*}
\mathbf{c} (t) = \left\{
\begin{aligned}
2t |z_1| + (1-2t)|z_0|, \;\;\;\;& \;\;\;\; t\in [0,\frac{1}{2}] ,\\
|z_1|e^{i(2t-1)\theta_1 },  \;\;\;\;\;\;\;\;\; & \;\;\;\; t\in [\frac{1}{2} ,1] , \\
\end{aligned}
\right.  
\end{equation*}
It remains to estimate the integral $|\int_{\mathbf{c}} \alpha |$. Set $\alpha'=\alpha - a_{-1}b_{-1}z^{-1} dz$. Then $\int_{\gamma} \alpha' =0$, and $\int_{\gamma} \alpha = a_{-1}b_{-1} \int_{\gamma} z^{-1} dz = 2\pi a_{-1} b_{-1} i $. Write $a_{-1}'= a_{-1}b_{-1}$, the Cauchy–Schwarz inequality gives
\begin{eqnarray*}
    \left| \int_{\mathbf{c}} \alpha' \right| & \leq & \left( \sum_{j\neq -1} |b_j|^2 \right)^{\frac{1}{2}} \left( \sum_{j\neq -1} \frac{|a_j|^2 | z_1^{j+1} - z_0^{j+1} |^2 }{(j+1)^2} \right)^{\frac{1}{2}} \leq \left( \sum_{j\neq -1} \frac{|a_j|^2 | z_1^{j+1} - z_0^{j+1} |^2 }{(j+1)^2} \right)^{\frac{1}{2}} .
\end{eqnarray*}

Moreover, a direct computation yields
$$ \re\left(a_{-1}' \int_{\mathbf{c}} z^{-1} dz \right)  
  = \re\left( a_{-1}' \right) \log \left|\frac{z_1}{z_0}\right| -  \im\left( a_{-1}' \right) \theta_1 .$$

Hence
\begin{eqnarray*}
  \left| \int_{\mathbf{c}} \re (\alpha) \right| & \leq &  \left| \int_{\mathbf{c}} \alpha' \right| +  \left| \int_{\mathbf{c}} \re (\alpha-\alpha' ) \right| \\
  & \leq & \left( \sum_{j\neq -1} \frac{|a_j|^2 | z_1^{j+1} - z_0^{j+1} |^2 }{(j+1)^2} \right)^{\frac{1}{2}} + \left| \re\left( a_{-1}' \right) \log \left|\frac{z_1}{z_0}\right|\right| + \left| \im\left( a_{-1}' \right) \theta_1 \right|,
\end{eqnarray*}
where $z=\varsigma e^{i\theta}$ is the polar coordinate.

As in (1), Lemma \ref{lemlocalorthonormalbasis} implies the desired estimates for the sums, leading to
\begin{eqnarray*}
  \sum_{j=0}^{\infty} \frac{|a_j|^2 | z_1^{j+1} - z_0^{j+1} |^2 }{(j+1)^2} & = & \sum_{j=0}^{\infty} \left| 1- \frac{z_0^{j+1}}{z_1^{j+1}} \right|^2 \frac{|a_j|^2 | z_1|^{2(j+1)}  }{(j+1)^2}  \\
  & \leq & \frac{\left( 1+\left|\frac{z_0}{z_1}\right| \right)^2}{\pi (1-e^{-3\pi} )} \sum_{j=0}^{\infty} \frac{ e^{2(j+1) \nu (\gamma ) } | z_1|^{2(j+1)} }{ j+1 } \\
  & = & \frac{-\log \left( 1-e^{2\nu (\gamma) } |z_1 |^2 \right) }{\pi (1-e^{-3\pi} )} \left( 1+\left|\frac{z_0}{z_1}\right| \right)^2 ,
\end{eqnarray*}
and similarly,
\begin{equation*}
  \sum_{j=0}^{\infty} \frac{|a_{-2-j}|^2 | z_1^{-j-1} - z_0^{-j-1} |^2 }{(j+1)^2} \leq \frac{-\log \left( 1-e^{2 ( \nu (\gamma ) - \frac{2\pi^2}{\ell (\gamma)} ) } |z_0 |^{-2} \right) }{\pi (1-e^{-3\pi} )} \left( 1+\left|\frac{z_0}{z_1}\right| \right)^2 .
\end{equation*}

It follows that
\begin{eqnarray*}
  \left| \int_{\mathbf{c}} \re(\alpha) \right| & \leq & \left( \sum_{j\neq -1} \frac{|a_j|^2 | z_1^{j+1} - z_0^{j+1} |^2 }{(j+1)^2} \right)^{\frac{1}{2}}  + \left| \re\left( a_{-1}' \right) \log \left|\frac{z_1}{z_0}\right|\right| + \left| \im\left( a_{-1}' \right) \theta_1 \right| \\
  & \leq & \left(  \frac{-\log \left( \left( 1-e^{2\nu (\gamma) } |z_1 |^2 \right) \left( 1-e^{2 ( \nu (\gamma ) - \frac{2\pi^2}{\ell (\gamma)} ) } |z_0 |^{-2} \right) \right)}{\pi (1-e^{-3\pi} ) \left( 1+\left|\frac{z_0}{z_1}\right| \right)^{-2}} \right)^{\frac{1}{2}} \\
  & & + \frac{1}{2\pi}\left| \int_{\gamma} \im( \alpha) \right| \left| \log \left|\frac{z_1}{z_0} \right| \right| + \frac{1}{2 } \left| \int_{\gamma} \re( \alpha) \right| ,
\end{eqnarray*}
as claimed.
\end{proof}

\subsection{Tian's peak sections}

We are now in a position to introduce Tian's peak section for Riemann surfaces.

\subsubsection*{Basic properties of Tian's peak sections}

Let $C$ be a connected compact smooth complex curve equipped with the classical topology, and let $\mu$ be a K\"ahler metric on $C$. Denote by $\mathscr{A}_{C, K} (C) $ the space of smooth differential forms on $C$ with coefficients in $K\in \{\mathbb{R},\mathbb{C}\}$, endowed with the $L^2$-inner product $\langle \cdot ,\cdot \rangle_{L^2 ;C} $. Let $\mathbf{X} \subset \mathscr{A}_{C, K} (C) $ be a nonzero $K$-linear subspace, and let $\mathbf{T}_{\mathbf{X}}:\mathbf{X} \to K $ be a bounded $K$-linear functional. 

An element $\alpha\in\mathbf{X}$ is called the peak section of $\mathbf{T}_{\mathbf{X}} $ if and only if it satisfies $\Vert \alpha \Vert_{L^2 ;C}=1$ and attains the operator norm of the functional: $\mathbf{T}_{\mathbf{X}} (\alpha) = \Vert \mathbf{T}_{\mathbf{X}} \Vert $.

We next present some fundamental properties of peak sections. While these follow from classical functional analysis, we include brief sketches of the proofs.

\begin{prop}
\label{propfundamentalpropertiespeaksection}
Let $C, \mu ,K,\mathscr{A}_{C, K} (C) , \mathbf{X} $ and $\mathbf{T}_{\mathbf{X}} $ be as above. Then the following properties hold.
\begin{enumerate}[(1)]
    \item If $\mathbf{T}_{\mathbf{X}} =0 $, then $\alpha\in\mathbf{X}$ is a peak section of $\mathbf{T}_{\mathbf{X}}$ if and only if $\Vert \alpha \Vert_{L^2 ;C}=1$.
    \item If $\mathbf{T}_{\mathbf{X}} \neq 0 $, then
    $$ \inf \left\{ \Vert \alpha \Vert_{L^2 ;C} :\alpha\in\mathbf{X} ,\;\; |\mathbf{T}_{\mathbf{X}} (\alpha)| = \Vert\mathbf{T}_{\mathbf{X}} \Vert \right\} = 1 . $$
    \item If $\mathbf{T}_{\mathbf{X}} \neq 0 $, then $\mathbf{T}_{\mathbf{X}} $ admits at most one peak section.
    \item If $\mathbf{T}_{\mathbf{X}} \neq 0 $ and $\alpha\in \mathbf{X} $ is the peak section of $\mathbf{T}_{\mathbf{X}}$, then $\mathbf{X}$ admits the $L^2$-orthogonal decomposition:
    $$ \mathbf{X} = K\alpha \oplus \ker \mathbf{T}_{\mathbf{X}} . $$
\end{enumerate}
\end{prop}

\begin{proof}
The propositions (1) and (2) follow directly from the definition of $\Vert\mathbf{T}_{\mathbf{X}} \Vert$. We thus begin with property (3).

Let $\alpha ,\alpha' \in \mathbf{X}$ be two peak sections of $\mathbf{T}_{\mathbf{X}} $. It suffices to show that $\alpha =\alpha' $. By definition, we have $\mathbf{T}_{\mathbf{X}} (\frac{1}{2} (\alpha +\alpha') ) = \frac{1}{2} ( \mathbf{T}_{\mathbf{X}} ( \alpha ) + \mathbf{T}_{\mathbf{X}} ( \alpha' ) )=\Vert \mathbf{T}_{\mathbf{X}} \Vert $. Hence, the $L^2$-norm satisfies $\left\Vert \frac{1}{2} (\alpha +\alpha') \right\Vert_{L^2 ;C} \geq 1 $. On the other hand, we compute:
\begin{eqnarray*}
 \left\Vert \frac{1}{2} (\alpha +\alpha') \right\Vert_{L^2 ;C}^2 & = & \frac{1}{2} \left\Vert \alpha \right\Vert^2_{L^2 ;C} +\frac{1}{2} \left\Vert \alpha' \right\Vert^2_{L^2 ;C} -\left\Vert \frac{1}{2} (\alpha -\alpha') \right\Vert^2_{L^2 ;C} \\
 & = & 1-\left\Vert \frac{1}{2} (\alpha -\alpha') \right\Vert^2_{L^2 ;C} .
\end{eqnarray*}
It follows that $\left\Vert \frac{1}{2} (\alpha +\alpha') \right\Vert_{L^2 ;C} =1$ and $\left\Vert \frac{1}{2} (\alpha -\alpha') \right\Vert_{L^2 ;C} =0 $. Therefore, $\alpha =\alpha'$.

We now turn to property (4).

By definition, $\mathbf{T}_{\mathbf{X}} ( \alpha ) = \Vert\mathbf{T}_{\mathbf{X}} \Vert >0 $, and hence $\mathbf{X}$ admits the decomposition:
$$ \mathbf{X} = K\alpha \oplus \ker \mathbf{T}_{\mathbf{X}} . $$ 
It remains to show that this is an $L^2$-orthogonal decomposition. 

Let $\beta\in \ker \mathbf{T}_{\mathbf{X}} \setminus \{0\} $. Then,
$$\mathbf{T}_{\mathbf{X}} \left( \alpha -\frac{\langle \alpha ,\beta \rangle_{L^2 ;C} }{\langle \beta ,\beta \rangle_{L^2 ;C}} \beta \right) = \mathbf{T}_{\mathbf{X}} \left( \alpha \right) = \Vert\mathbf{T}_{\mathbf{X}} \Vert .$$
By property (ii), we must then have $\left\Vert \alpha -\frac{\langle \alpha ,\beta \rangle_{L^2 ;C} }{\langle \beta ,\beta \rangle_{L^2 ;C}} \right\Vert_{L^2 ;C} \geq 1$. However, a straightforward calculation shows:
\begin{eqnarray*}
\left\Vert \alpha -\frac{\langle \alpha ,\beta \rangle_{L^2 ;C} }{\langle \beta ,\beta \rangle_{L^2 ;C}} \beta \right\Vert^2_{L^2 ;C} & = & \left\Vert \alpha \right\Vert^2_{L^2 ;C} -  \frac{\left| \langle \alpha ,\beta \rangle_{L^2 ;C} \right|^2 }{\langle \beta ,\beta \rangle_{L^2 ;C}} \\
& = & 1 - \frac{\left| \langle \alpha ,\beta \rangle_{L^2 ;C} \right|^2 }{\langle \beta ,\beta \rangle_{L^2 ;C}} .
\end{eqnarray*}
Therefore, we must have $\langle \alpha ,\beta \rangle_{L^2 ;C} =0 $. This shows that the decomposition is indeed $L^2$-orthogonal, which completes the proof.
\end{proof}

In general, it is not easy to determine whether a given linear functional is bounded. Moreover, since the space $\mathscr{A}_{C, K} (C)$ is not complete, a given bounded functional does not necessarily admit a peak section. However, if the linear functional can be reduced to a finite-dimensional linear subspace, then its boundedness and the existence of a peak section follow automatically.

\begin{prop}
\label{propexistancepeaksection}
Let $C, \mu ,K,\mathscr{A}_{C, K} (C) $ and $\mathbf{X} $ be as above. Let $\mathbf{T}_{\mathbf{X}}:\mathbf{X} \to K $ be a $K$-linear functional (not necessarily bounded). Assume that $\mathbf{X} $ admits an $L^2$-orthogonal decomposition $\mathbf{X} = \mathbf{X}_1\oplus \mathbf{X}_2 $ such that $\dim_{K} \mathbf{X}_1 <\infty $ and the restriction $\mathbf{T}_{\mathbf{X}_2} := \mathbf{T}_{\mathbf{X}} |_{\mathbf{X}_2} =0 $. Then $\mathbf{T}_{\mathbf{X}}$ is a bounded $K$-linear functional, and there exists a peak section $\alpha \in \mathbf{X}_1 \subset \mathbf{X} $ of $\mathbf{T}_{\mathbf{X}} $.
\end{prop}

\begin{proof}
Without loss of generality, we can assume that $\mathbf{T}_{\mathbf{X}} \neq 0 $. Since $\dim_{K} \mathbf{X}_1 <\infty $, the restriction $\mathbf{T}_{\mathbf{X}_1} := \mathbf{T}_{\mathbf{X}} |_{\mathbf{X}_1} $ is a bounded $K$-linear functional on $\mathbf{X}_1$. By the Riesz representation theorem, there exists a unique element $\alpha\in \mathbf{X}_1 $ such that $\Vert\alpha\Vert_{L^2 ,C} =1$ and $\mathbf{T}_{\mathbf{X}_1} (\alpha) = \Vert \mathbf{T}_{\mathbf{X}_1} \Vert $. Hence, $\alpha$ is the peak section of $\mathbf{T}_{\mathbf{X}_1} $. By Proposition \ref{propfundamentalpropertiespeaksection}, the space $\mathbf{X}_1$ admits the $L^2$-orthogonal decomposition $ \mathbf{X}_1 = K\alpha \oplus \ker \mathbf{T}_{\mathbf{X}_1} $. Since the restriction $\mathbf{T}_{\mathbf{X}_2} =0 $, we have $\ker \mathbf{T}_{\mathbf{X} }= \ker \mathbf{T}_{\mathbf{X}_1} \oplus \mathbf{X}_2 $, and thus $\mathbf{X}$ admits the $L^2$-orthogonal decomposition $ \mathbf{X} = K\alpha \oplus \ker \mathbf{T}_{\mathbf{X} } $. It follows that $\Vert \mathbf{T}_{\mathbf{X} } \Vert = \Vert \mathbf{T}_{\mathbf{X}_1} \Vert = \mathbf{T}_{\mathbf{X} } (\alpha ) $, and hence $\alpha$ is the peak section of $\mathbf{T}_{\mathbf{X} } $.
\end{proof}

As a corollary, we see that in certain special cases relevant to our setting, a peak section always exists.

\begin{cor}
\label{corpeaksectiontwocases}
Let $C$ be a connected compact smooth complex curve equipped with the classical topology. Denote by $\mathscr{A}_{C, K} (C) $ the space of smooth differential forms on $C$ with coefficients in $K\in \{\mathbb{R},\mathbb{C}\}$, endowed with the $L^2$-inner product $\langle \cdot ,\cdot \rangle_{L^2 ,C} $. Let $\mathbf{X} \subset \mathscr{A}_{C, K} (C) $ be a nonzero $K$-linear subspace, and let $\mathbf{T}_{\mathbf{X}}:\mathbf{X} \to K $ be a $K$-linear functional (not necessarily bounded). 

Then, in each of the following two settings, $\mathbf{T}_{\mathbf{X}}$ is a bounded $K$-linear functional, and there exists a peak section $\alpha \in \mathbf{X} $ of $\mathbf{T}_{\mathbf{X}} $:
\begin{enumerate}[(1)]
    \item The field $K=\mathbb{C} $, the space $\mathbf{X} \subset \Gamma (C,\omega_C) \subset \mathscr{A}^1_{C, \mathbb{C}} (C) $, and there exists a point $x\in C$ together with a holomorphic local coordinate $z:U_x \to \mathbb{C} $ defined on an open neighborhood $U_x\subset C$ of $x$, such that
    $$\alpha (x) =\mathbf{T}_{\mathbf{X}} (\alpha) dz ,\;\;\; \forall \alpha\in \mathbf{X} , $$
    where $\alpha (x)$ is the evaluation of $\alpha$ at $x$.
    \item The space $\mathbf{X} \subset \mathscr{Z}^1_{C, \mathbb{C}} (C) $, and there exists a form $\beta\in \mathscr{Z}^1_{C, \mathbb{C}} (C) $ such that 
    $$ \mathbf{T}_{\mathbf{X}} (\alpha) =\int_C \alpha \wedge \beta ,\;\; \forall \alpha\in \mathbf{X} .$$
\end{enumerate}
\end{cor}

\begin{proof}
In the first setting, since $C$ is compact, we have $\dim_{\mathbb{C}} \mathbf{X} \leq \dim_{\mathbb{C}} \Gamma(C, \omega_C) < \infty$. Thus, Proposition \ref{propexistancepeaksection} applies, yielding the boundedness of $\mathbf{T}_{\mathbf{X}}$ and the existence of a peak section.

Now consider the second setting. By Corollary \ref{corisomorphismharmonicdhcohomology}, we have the $L^2$-orthogonal decomposition
\begin{equation*}
 \mathscr{Z}^1_{C, K} (C) = d\mathscr{A}^{0}_{C, K} (C) \oplus \mathscr{H}^1_{C, K} (C),
\end{equation*}
and the assumption on $\mathbf{T}_{\mathbf{X}}$ implies that $\mathbf{T}_{\mathbf{X}}|_{d\mathscr{A}^0_{C, K}(C)} = 0$. Since $\dim_{K}\mathscr{H}^1_{C, K}(C) <\infty $ by Theorem \ref{thmhodgedecomposition}, Proposition \ref{propexistancepeaksection} applies again, and the result follows.
\end{proof}

\subsubsection*{Two cases of peak sections: at a point and along a cycle}

The next step is to further analyze the two settings described in Corollary \ref{corpeaksectiontwocases}.

We now turn to a more detailed study of the first setting in Corollary \ref{corpeaksectiontwocases}.

Let $C$ be a connected compact smooth complex curve equipped with the classical topology, and let $\mu$ be a K\"ahler metric on $C$. Denote by $\Gamma (C,\omega_C) $ the space of holomorphic $1$-forms on $C$, endowed with the $L^2$-inner product
$$ \langle \alpha , \beta \rangle_{L^2 ;C} = \frac{1}{2} \int_{C} \alpha \wedge \star \beta = \frac{1}{2} \int_C \alpha \wedge \bar{\beta} ,\;\; \forall \alpha ,\beta \in \Gamma (C,\omega_C) . $$

We say that $\alpha$ is a peak section in $\mathbf{X}\subset\Gamma (C,\omega_C) $ at a point $x\in C$, if there exists a holomorphic local coordinate $z:U_x \to \mathbb{C} $ defined on an open neighborhood $U_x\subset C$ of $x$, such that $\beta (x) =\mathbf{T}_{\mathbf{X}} (\beta) dz $, $ \forall \beta\in \mathbf{X} $, where $\beta (x)$ denotes the evaluation of $\beta$ at the point $x$. The relation above is equivalent to $ \mathbf{T}_{\mathbf{X}} (\beta) = \beta (x) (\frac{\partial}{\partial z} |_{x}) $, $\forall \beta\in \mathbf{X}$. The set of peak sections in $\mathbf{X}$ at a given point $x\in C$ can be described as follows.

\begin{lem}
\label{lempeaksectionatapoint}
Let $C, \Gamma (C,\omega_C) $ and $\mathbf{X} $ be as above. Let $\alpha_0 \in \mathbf{X} $ satisfy $\Vert\alpha_0 \Vert_{L^2 ;C}=1$ and $\alpha_0 (x) \neq 0 $ for some $x\in C$. Then $\alpha_0 $ is a peak section in $\mathbf{X}$ at $x$ if and only if $\alpha_0 $ is $L^2$-orthogonal to the subspace of $\mathbf{X}$ consisting of sections vanishing at $x$:
$$ \mathbf{X} (x) =\{ \alpha \in \mathbf{X} : \alpha (x) = 0 \} .$$
Moreover, if $\alpha_0 $ and $\alpha_1 $ are both peak sections in $\mathbf{X}$ at $x$, then there exists a constant $\theta\in [0,2\pi )$ such that $\alpha_1 =e^{i\theta} \alpha_0 $.
\end{lem}

\begin{proof}
We begin with the ``only if'' direction. Suppose that $\alpha_0 $ is a peak section in $\mathbf{X}$ at $x$. Then there exists a holomorphic local coordinate $z:U_{x} \to \mathbb{C} $ defined on an open neighborhood $U_x\subset C$ of $x$ such that $\alpha_0$ is the peak section corresponding to the functional $\mathbf{T}_{\mathbf{X} ,0} : \mathbf{X}\to\mathbb{C} $ defined by $ \mathbf{T}_{\mathbf{X},0} (\alpha) = \alpha (x) (\frac{\partial}{\partial z} |_{x}) $, $ \forall \alpha\in \mathbf{X} $. It is clear that $\ker \mathbf{T}_{\mathbf{X},0} = \mathbf{X} (x) $. By Proposition \ref{propfundamentalpropertiespeaksection}, the peak section $\alpha_0$ must be $L^2$-orthogonal to $ \ker \mathbf{T}_{\mathbf{X},0} = \mathbf{X} (x) $, which proves the ``only if'' direction.

Now we turn to the ``if'' direction. Let $z:U_{x} \to \mathbb{C} $ be a holomorphic local coordinate on a neighborhood $U_x\subset C$ of $x$. Since $\alpha_0(x) \neq 0$, we may write $\alpha_0 (x) = a_0 dz $ for some $a_0 \in\mathbb{C}^* $. Define a new local coordinate $w=a_0 z $ on $U_x $, so that $dw (x) = a_0 dz(x) $. Now consider the $\mathbb{C}$-linear functional $ \mathbf{T}_{\mathbf{X},1} : \mathbf{X} \to \mathbb{C} $ defined by requiring $ \alpha (x) = \mathbf{T}_{\mathbf{X},1} (\alpha) dw $, $\forall \alpha\in\mathbf{X} $. Then $\mathbf{T}_{\mathbf{X},1} (\alpha_0 ) =1 $, and again $ \ker \mathbf{T}_{\mathbf{X},0} = \mathbf{X} (x) $. By assumption, $\mathbf{X}=\mathbb{C}\alpha_0 \oplus \mathbf{X} (x) $ is an $L^2$-orthogonal decomposition. Then, by Proposition \ref{propexistancepeaksection}, there exists a peak section $\alpha'_0 \in \mathbb{C}\alpha_0 $ of $ \mathbf{T}_{\mathbf{X},1}$. By definition, we have $\Vert\alpha'_0\Vert_{L^2 ;C}=1=\Vert\alpha_0\Vert_{L^2 ;C}$, and $ \mathbf{T}_{\mathbf{X},1} (\alpha_0 ), \mathbf{T}_{\mathbf{X},1} (\alpha'_0 ) \in (0,\infty) \subset \mathbb{R} $. It follows that $\alpha_0 =\alpha'_0$ is the peak section of $ \mathbf{T}_{\mathbf{X},1}$. This proves the ``if'' direction.

Finally, suppose that both $\alpha_0$ and $\alpha_1$ are peak sections in $\mathbf{X}$ at $x$. Then, by a similar argument, $\alpha_1 \in\mathbb{C}\alpha_0 $ and $\Vert\alpha_0\Vert_{L^2 ;C}=\Vert\alpha_1\Vert_{L^2 ;C} =1 $. Hence, $\alpha_1 =e^{i\theta} \alpha_0 $ for some $\theta\in [0,2\pi )$.
\end{proof}

By the above lemma, if there exists an element $\alpha\in\mathbf{X} $ that does not vanish at the given point $x$, then all peak sections at $x\in C$ are equivalent up to a complex rotation $e^{i\theta}$. However, to compare peak sections at different points, we need to use the K\"ahler metric $\mu$. 

Let $\{\alpha_j \}_{j=1}^{\dim_{\mathbb{C}} \mathbf{X} }$ be an $L^2$-orthonormal basis of $\mathbf{X} $. Then there exists a non-negative, $\mathbb{R}$-valued function $\mathbf{b}_{ \mathbf{X},\mu }$ on $C$ such that
$$ \frac{i}{2} \sum_{j=1}^{\dim_{\mathbb{C}} \mathbf{X}} \alpha_j \wedge \bar{\alpha}_j = \mathbf{b}_{ \mathbf{X} ,\mu } \cdot \mu .$$
It is clear that $\mathbf{b}_{ \mathbf{X},\mu }$ is independent of the choice of the $L^2$-orthonormal basis. This function $\mathbf{b}_{ \mathbf{X},\mu }$ is called the Bergman kernel function associated with $(  \mathbf{X},\mu )$. A point $x\in C$ is called a peak point of $(  \mathbf{X},\mu )$, if it satisfies $\mathbf{b}_{ \mathbf{X},\mu } (x) =\Vert\mathbf{b}_{ \mathbf{X},\mu }\Vert_{L^\infty } $.

We now record some basic properties of the Bergman kernel function.

\begin{lem}
\label{lembasicpropertiesBergmankernelfunction}
Let $C,\mu , \Gamma (C,\omega_C) , \mathbf{X} $ and $\mathbf{b}_{ \mathbf{X} ,\mu } $ be as above. Then the following properties hold.
\begin{enumerate}[(1)]
    \item Let $\alpha \in \mathbf{X} $ be a peak section at $x\in C$. Then at $x$, 
    $$\frac{i}{2} \alpha (x) \wedge \bar{\alpha} (x) = \mathbf{b}_{ \mathbf{X} ,\mu } (x) \mu (x) .$$
    \item The $L^\infty$-norm $\Vert \mathbf{b}_{ \mathbf{X} ,\mu } \Vert_{L^\infty} = \sup\limits_{C} |\mathbf{b}_{ \mathbf{X} ,\mu }| $ satisfies
    $$ \Vert \mathbf{b}_{ \mathbf{X} ,\mu } \Vert_{L^\infty} \geq \frac{\dim_{\mathbb{C}} \mathbf{X}}{\int_C \mu } .$$
\end{enumerate}
\end{lem}

\begin{proof}
Let $x \in C$, and let $\alpha \in \mathbf{X}$ be a peak section at $x$. Without loss of generality, we assume $\alpha(x) \neq 0$. Choose an $L^2$-orthonormal basis $\{\alpha_j\}_{j=1}^{\dim_{\mathbb{C}} \mathbf{X}}$ of $\mathbf{X}$ such that $\alpha_1 = \alpha$. Then by Lemma \ref{lempeaksectionatapoint}, we have $\alpha_j(x) = 0$ for all $j \geq 2$. Therefore, at $x$,
$$\frac{i}{2} \alpha(x) \wedge \bar{\alpha}(x) = \frac{i}{2} \sum_{j=1}^{\dim_{\mathbb{C}} \mathbf{X}} \alpha_j \wedge \bar{\alpha}_j = \mathbf{b}_{ \mathbf{X},\mu }(x) \mu(x),$$
which proves (1).

To prove (2), note that by the definition of $\mathbf{b}_{ \mathbf{X},\mu }$, we have
$$ \Vert \mathbf{b}_{ \mathbf{X} ,\mu } \Vert_{L^\infty} \cdot \int_C \mu \geq \int_C \mathbf{b}_{ \mathbf{X} ,\mu } \cdot \mu = \int_C \frac{i}{2} \sum_{j=1}^{\dim_{\mathbb{C}} \mathbf{X}} \alpha_j \wedge \bar{\alpha}_j = \dim_{\mathbb{C}} \mathbf{X} ,$$
which gives (2).
\end{proof}

Our attention now shifts to the second setting described in Corollary \ref{corpeaksectiontwocases}.

Let $C$ be a connected compact smooth complex curve equipped with the classical topology, and let $\mu$ be a K\"ahler metric on $C$. Denote by $\mathscr{A}^1_{C, K} (C) $ the space of smooth $1$-forms on $C$, endowed with the $L^2$-inner product
$$ \langle \alpha , \beta \rangle_{L^2 ;C} = \frac{1}{2} \int_{C} \alpha \wedge \star \beta ,\;\; \forall \alpha ,\beta \in \mathscr{A}^1_{C, K} (C) , $$ 
where $K=\mathbb{R}$ or $\mathbb{C}$. 

Let $\mathbf{X} \subset \mathscr{Z}^1_K(C) \subset \mathscr{A}^1_K(C)$ be a subspace of closed $1$-forms, where $K = \mathbb{R}$ or $\mathbb{C}$. Given a singular $1$-cycle $\sigma \in Z_1(C; K)$, we define a linear functional $\mathbf{T}_{\mathbf{X},\sigma} : \mathbf{X} \to K$ by
$$
\mathbf{T}_{\mathbf{X},\sigma}(\alpha) = \sum_{j=1}^N a_j \int_{\gamma_j'} \alpha = \int_C \alpha \wedge \beta,
$$
where $\sigma = \sum_{j=1}^N a_j \gamma_j'$ is a $K$-linear combination of smooth oriented loops representing the homology class $[\sigma] \in H_1(C;K)$, and $[\beta] \in H^1_{\mathrm{dR}} (C; K)$ is its Poincar\'e dual. The functional $\mathbf{T}_{\mathbf{X},\sigma}$ is well-defined on $\mathbf{X}$, since both expressions depend only on the cohomology class of $\alpha$ and the homology class of $\sigma$.

We say that a form $\alpha \in \mathbf{X}$ is a peak section along the cycle $\sigma$ if it is the peak section of the functional $\mathbf{T}_{\mathbf{X},\sigma}$. Moreover, if there exists $\alpha\in\mathbf{X}$ such that $\int_{C} \alpha \wedge \beta \neq 0 $, then the peak section in $\mathbf{X} $ along the cycle $\sigma $ is unique.

We now consider the explicit form of the peak sections in certain natural subspaces of closed differential forms.

\begin{lem}
\label{lemexplicitpeaksectionalongacycle}
Let $C $ and $\mu $ be as above. Let $\sigma \in Z_1 (C;K) $ be a cycle, and $[\beta]\in H^1_{\mathrm{dR}} (C;K)$ its Poincar\'e dual, where $K=\mathbb{R}$ or $\mathbb{C}$. Let $\beta_{\mathscr{H}}$ be the harmonic $1$-form corresponding to $[\beta] $. Assume that $[\beta]\neq 0$. Let $\alpha $ be the peak section in $\mathbf{X}\subset \mathscr{Z}^1_{C, K} (C)$ along the cycle $\sigma $. Then the following properties hold.
\begin{enumerate}[(1)]
    \item Assume that $\mathscr{H}^1_{C, K} (C) \subset \mathbf{X}$. Then $\alpha=-\frac{1}{\Vert \beta_{\mathscr{H}} \Vert_{L^2;C}} \star \beta_{\mathscr{H}} $.
    \item Assume that $K= \mathbb{R}$ and $\mathbf{X} =\Gamma (C,\omega_C) $. Then 
    $$\alpha=-\frac{\sqrt{2}}{2\Vert \beta_{\mathscr{H}} \Vert_{L^2;C}} \star \beta_{\mathscr{H}} + \frac{\sqrt{2} i}{2\Vert \beta_{\mathscr{H}} \Vert_{L^2;C}} \beta_{\mathscr{H}} .$$
\end{enumerate}
\end{lem}

\begin{proof}
We first prove (1).

By the Cauchy-Schwarz inequality, for any $\alpha'\in \mathbf{X}$, we have
\begin{eqnarray*}
\int_C \alpha'\wedge \beta & = & \int_C \alpha'\wedge \beta_{\mathscr{H}} = 2\langle \alpha' ,-\star \beta_{\mathscr{H}} \rangle_{L^2;C} \\
& \leq & 2\Vert \alpha' \Vert_{L^2;C} \cdot \Vert\star \beta_{\mathscr{H}}  \Vert_{L^2;C} = 2\Vert \alpha' \Vert_{L^2;C} \cdot \Vert \beta_{\mathscr{H}} \Vert_{L^2;C} ,
\end{eqnarray*}
with equality if and only if $\alpha' = -a (\star \beta_{\mathscr{H}}) $ for some $a\in (0,\infty)$. Since $\star \beta_{\mathscr{H}}\in \mathscr{H}^1_{C, K} (C) \subset \mathbf{X} $, it follows that $\alpha=-\frac{1}{\Vert \beta_{\mathscr{H}} \Vert_{L^2;C}} \star \beta_{\mathscr{H}} $. Note that $\Vert \alpha \Vert_{L^2} =1 $. This gives (1).

We now turn to (2).

By the Cauchy-Schwarz inequality again, for any $\alpha''\in \mathbf{X}$, we have
\begin{eqnarray*}
\int_C \alpha''\wedge \beta & = & \int_C \re (\alpha'')\wedge \beta_{\mathscr{H}} = 2\langle \re (\alpha'') ,-\star \beta_{\mathscr{H}} \rangle_{L^2} \\
& \leq & 2\Vert \re (\alpha'') \Vert_{L^2;C} \cdot \Vert\star \beta_{\mathscr{H}}  \Vert_{L^2;C} = 2\Vert \re (\alpha'') \Vert_{L^2;C} \cdot \Vert \beta_{\mathscr{H}} \Vert_{L^2;C} \\
& = & \Vert \alpha'' \Vert_{L^2;C} \cdot \Vert \beta_{\mathscr{H}} \Vert_{L^2;C} ,
\end{eqnarray*}
with equality if and only if $\re (\alpha'') = -a (\star \beta_{\mathscr{H}}) $ for some $a\in (0,\infty)$. Since $\Vert \alpha \Vert_{L^2;C} =1 $, we have
$$\alpha=-\frac{\sqrt{2}}{2\Vert \beta_{\mathscr{H}} \Vert_{L^2;C}} \star \beta_{\mathscr{H}} + \frac{\sqrt{2} i}{2\Vert \beta_{\mathscr{H}} \Vert_{L^2;C}} \beta_{\mathscr{H}} .$$
This completes the proof of (2).
\end{proof}

\section{Hyperbolic eigenvalues and eigenfunctions}

Let $C$ be a curve of genus $g \geq 2$ over $\mathbb{C}$, and let $\omega_C$ denote its canonical bundle. Let $\mu_{\mathrm{hyp}}$ be the unique K\"ahler metric on $C $ with constant curvature $4\pi (1-g) $ and volume $1$, and let $\mu_{\mathrm{KE}} $ be the unique K\"ahler metric on $C$ with constant curvature $-1 $. In particular, $\mu_{\mathrm{KE}} = 4\pi (g-1) \mu_{\mathrm{hyp}} $. Since $\dim C =1 $, we identify a K\"ahler metric $\mu$ with its associated volume form and measure.

For any K\"ahler metric $\mu$, the associated Laplacian is given by
$$\Delta_{\mu} = \frac{dd^c}{\mu} = \frac{i}{\pi} \frac{\partial\bar{\partial}}{\mu} = \frac{1}{2\pi} \star_{\mu} d \star_{\mu} d ,$$ 
where $d^c = \frac{1}{2\pi i} (\partial -\bar{\partial}) $, and $\star_{\mu} $ is the Hodge $\star$-operator associated with $(C,\mu)$. 

The spectrum of $\Delta_{\mu}$ consists of eigenvalues
$$ 0=\lambda_{\mu ,0 } < \lambda_{\mu ,1 } \leq \lambda_{\mu ,2 } \leq \cdots ,$$
listed with multiplicities. Correspondingly, there exists an $L^2$-orthonormal basis $\{ \phi_{\mu , l } \}_{l=0}^\infty $ of the $L^2$-completion of $\mathscr{A}^{0}_{C,\mathbb{R}} (C) $, $L^2 (C,\mu ;\mathbb{R} )$, such that 
$$\Delta_{\mu} \phi_{\mu , l } + \lambda_{\mu ,l } \phi_{\mu , l } =0 .$$

For $(x,y)\in C\times C$ and $t>0$, we introduce the kernels
\begin{eqnarray*}
G_{\mu} (x,y;t) & = & \sum_{l=1}^\infty \frac{e^{-t\lambda_{\mu,l} }}{\lambda_{\mu,l}} \phi_{\mu , l } (x) \phi_{\mu , l }(y) ,\\
K_{\mu} (x,y;t) & = & \sum_{l=1}^\infty e^{-t\lambda_{\mu,l} } \phi_{\mu , l } (x) \phi_{\mu , l }(y) .
\end{eqnarray*}
By construction, for any fixed $y_0 \in C$, the functions $(x,t) \mapsto G_{\mu}(x,y_0;t)$ and $(x,t) \mapsto K_{\mu}(x,y_0;t)$ satisfy the heat equation $\frac{\partial u (x,t)}{\partial t} = \Delta_{\mu} u(x,t) $. Moreover, we have the weak limits (in the sense of distributions)
$$ \lim\limits_{t\to 0^+}G_{\mu} (x,y_0;t )=G_{\mu} (x,y_0 )\quad\textrm{ and } \quad \lim\limits_{t\to 0^+} K_{\mu} (x,y_0;t)= \delta_{y_0}- \mu (C)^{-1} , $$
where $G_{\mu} (x,y_0 )$ is the $\mu$-admissable Green function of the diagonal $\Delta\subset C\times C$, and $\delta_{y_0}$ is the Dirac measure at $y_0 $. Hence $G_{\mu}(x,y;t)$ can be regarded as an extension of the Green function $G_{\mu}(x,y)$ to positive times. Furthermore, 
$$ \frac{\partial G_{\mu}}{\partial t} (x,y;t) = \Delta_{\mu ,x} G_{\mu} (x,y;t) = \Delta_{\mu ,y} G_{\mu} (x,y;t) = -K_{\mu} (x,y;t) .$$
Letting $t\to 0^+$, we obtain in the sense of distributions
$$ \Delta_{\mu ,x} G_{\mu} (x,y_0) = -\delta_{y_0}+ \mu (C)^{-1} . $$
In particular, by the Poincar{\'e}--Lelong formula \cite[Theorem B.2.19]{mamari1}, for any holomorphic coordinate $z$ centered at $y_0$, the Green function
$G_{\mu} (x,y_0)$ admits the local expansion
$$ G_{\mu} (x,y_0) = -\log |z(x)| + o(1) .$$
When $\mu=\mu_{\mathrm{Ar}}$, the Green function $G_{\mu_{\mathrm{Ar}}}(x,y)$ coincides with the Green function $G_a$ defined in \S\ref{section_ad}.

By the maximum principle for the heat equation, we obtain 
$$K_{\mu} (x,y;t) \geq - \mu (C)^{-1} ,$$ 
and consequently, 
$$  G_{\mu} (x,y;t) \leq  t\mu (C)^{-1}+G_{\mu} (x,y) ,\quad \textrm{ if } x\neq y.$$
See also \cite[Theorem 10.1]{peterli1} for further details.

When $\mu=\mu_{\mathrm{Ar}}$, $\mu_{\mathrm{hyp}}$ or $\mu_{\mathrm{KE}}$, we omit the $\mu$ in the subscript. For example, we write the corresponding Laplacians as $\Delta_{\mathrm{Ar}}$, $\Delta_{\mathrm{hyp}}$ and $\Delta_{\mathrm{KE}}$.

In this section, we establish estimates for the eigenvalues and eigenfunctions associated with $\Delta_{\mathrm{KE}}$ and $\Delta_{\mathrm{hyp}}$. Our approach is to control the larger eigenvalues and eigenfunctions via the heat kernel, while estimating the smaller ones individually. In particular, a result of \cite[Th\'eor\`eme 2]{otro1} shows that there are at most $2g-3$ positive eigenvalues $\leq \frac{1}{8\pi}$. Note that our normalization of the Laplacian differs from the classical convention by a factor of $2\pi$.

\subsection{An explicit upper bound of hyperbolic heat kernel}

Let $\rho_{\mathcal{H}} : \mathcal{H} \to C$ be the universal covering map, where $\mathcal{H}$ is the Poincar\'e upper half-plane equipped with the metric $\widetilde{\mu_{\mathrm{KE}}} =\frac{\mu_{\mathrm{Euc}}}{(\mathrm{Im}(\tau))^{2}}$, with $\mu_{\mathrm{Euc}} = \frac{idz\wedge d\bar{z}}{2} $ being the standard Euclidean metric on $\mathbb{C}$. See also Theorem \ref{thmpoincaremodel}. Then the deck transformation group $\pi_1(C)$ embeds into $\Aut(\mathcal{H}) \cong \mathrm{PSL}(2, \mathbb{R})$, and $C$ is isomorphic, as a K\"ahler curve, to the quotient $\mathcal{H} / \Gamma$, where $\Gamma$ is the image of $\pi_1(C)$. Moreover, $\rho_{\mathcal{H}}$ corresponds to the natural projection $\mathcal{H} \to \mathcal{H}/\Gamma$.

We now identify $C$ with the quotient $\mathcal{H} / \Gamma$ and compute the heat kernel on $C$ via the $\Gamma$-action on $\mathcal{H}$.

\begin{lem}
\label{lemheatkernelformulaquotientstructure}
Let $C$, $\rho_{\mathcal{H}}$, $\mathcal{H}$ and $\Gamma$ be as above. Then for any $z,w\in\mathcal{H} $ and $t>0$, we have
\begin{eqnarray*}
K_{\mathrm{hyp}} (\rho_{\mathcal{H}} (z),\rho_{\mathcal{H}} (w);t) & = & 4\pi (g-1) K_{\mathrm{KE}} (\rho_{\mathcal{H}} (z),\rho_{\mathcal{H}} (w);4\pi (g-1)t) \\
& = & -1+ 4\pi (g-1)\sum_{\sigma\in\Pi} K_{\mathcal{H}} (z,\sigma w;2 (g-1)t)  ,
\end{eqnarray*}
where 
$$K_{\mathcal{H}} (z,w; t) = \frac{\sqrt{2}e^{-\frac{t}{4}}}{(4\pi t)^{\frac{3}{2}}} \int_{\mathrm{dist}_{\mathcal{H}} (z,w)}^\infty \frac{\varsigma e^{-\frac{\varsigma^2}{4t}}}{(\cosh \varsigma - \cosh (\mathrm{dist}_{\mathcal{H}} (z,w)) )^{\frac{1}{2}}} d\varsigma ,$$
and $\mathrm{dist}_{\mathcal{H}}$ denote the hyperbolic distance on $\mathcal{H}$ determined by $\frac{\mu_{\mathrm{Euc}}}{(\mathrm{Im}(\tau))^{2}} $.
\end{lem}

\begin{proof}
The first equality follows directly from $\mu_{\mathrm{KE}}=4\pi (g-1)\mu_{\mathrm{hyp}}$, while the second one requires a more detailed argument.

Let $K_{\Gamma}(\rho_{\mathcal{H}} (z),\rho_{\mathcal{H}} (w);t) = \sum_{\sigma\in\Gamma} K_{\mathcal{H}} (z,\sigma w; t)$. By \cite[Theorem 7.5.11]{bu1}, $K_{\Gamma}$ is the fundamental solution of the heat equation $\frac{\partial u}{\partial t} = ( \star_{\mathrm{KE}} d )^2 u = 2\pi \Delta_{\mathrm{KE}} u $. Hence for any $y_0\in \mathcal{H}/\Gamma$, we have the weak limit $\lim\limits_{t\to 0^+} K_{\Gamma} (x,y_0;t)=\delta_{y_0}$. Define $K'_{\mathrm{KE}} (\rho_{\mathcal{H}} (z),\rho_{\mathcal{H}} (w); t) = K_{\Gamma}(\rho_{\mathcal{H}} (z),\rho_{\mathcal{H}} (w);\frac{t}{2\pi }) - \frac{1}{4\pi (g-1)} $. One can verify that for any $y_0$, $(x;t)\mapsto K'_{\mathrm{KE}}(x,y_0;t)$ satisfies the heat equation $\frac{\partial u}{\partial t} = \Delta_{\mathrm{KE}} u $ together with the initial condition $\lim\limits_{t\to 0^+} K'_{\mathrm{KE}} (x,y_0;t)= \delta_{y_0}- \frac{1}{4\pi (g-1)}$. By the uniqueness of the heat kernel \cite[Theorem 10.2]{peterli1}, we conclude that $K_{\mathrm{KE}} (x,y;t) = K'_{\mathrm{KE}} (x,y;t) $, which completes the proof.
\end{proof}

Before estimating the hyperbolic heat kernel $K_{\mathrm{hyp}} (x,x;t)$, we provide the following estimate for $K_{\mathcal{H}} $.

\begin{lem}
\label{lemestimatehyperbolicheatkernelonhalfplane}
Set $u (a,t)=\int_{a}^\infty \frac{\varsigma e^{-\frac{\varsigma^2}{4t}}}{(\cosh \varsigma - \cosh a )^{\frac{1}{2}}} d\varsigma $, $a> 0$, $t>0$. Then
$$ u (a,t) < \frac{3\sqrt{2}}{4} ae^{-\frac{a^2}{4t}} + \frac{5\sqrt{2\pi t}}{3} e^{-\frac{25a^2}{64t}} \leq \frac{3\sqrt{t}}{\sqrt{2e}} \cdot e^{-\frac{a^2}{8t}} + \frac{5\sqrt{2\pi t}}{3} e^{-\frac{25a^2}{64t}} $$
In particular, for $a\geq 2\arcsinh 1$ and $t\in (0,\frac{1}{20}]$, it holds that 
$$ \frac{1}{t^{\frac{3}{2}}}\int_{a}^\infty \frac{\varsigma e^{-\frac{\varsigma^2}{4t}}}{(\cosh \varsigma - \cosh a )^{\frac{1}{2}}} d\varsigma < \frac{3}{100000} .$$
\end{lem}

\begin{proof}
The proof follows the argument of \cite[Lemma 7.4.6]{bu1}.

Since $\cosh\varsigma - \cosh a=\int_{a}^\varsigma \sinh \vartheta d\vartheta > \int_{a}^\varsigma \vartheta d\vartheta = \frac{1}{2} (\varsigma^2 -a^2) $, we obtain
\begin{eqnarray*}
\int_{\frac{5a}{4}}^{\infty} \frac{\varsigma e^{-\frac{\varsigma^2}{4t}}}{(\cosh \varsigma - \cosh a )^{\frac{1}{2}}} d\varsigma & < & \sqrt{2}\int_{\frac{5a}{4}}^{\infty} \frac{\varsigma}{(\varsigma^2 -a^2)^{\frac{1}{2}} } \cdot e^{-\frac{\varsigma^2}{4t}} d\varsigma < \frac{5\sqrt{2}}{3} \int_{0}^{\infty} e^{-\frac{(\varsigma+\frac{5a}{4})^2}{4t}} d\varsigma \\
& < & \frac{5\sqrt{2}}{3}  e^{-\frac{25a^2}{64t}} \int_{0}^{\infty} e^{-\frac{\varsigma^2}{4t}} d\varsigma = \frac{5\sqrt{2\pi t}}{3} e^{-\frac{25a^2}{64t}} .
\end{eqnarray*}

Similarly,
\begin{equation*}
\int_{a}^{\frac{5a}{4}} \frac{\varsigma e^{-\frac{\varsigma^2}{4t}}}{(\cosh \varsigma - \cosh a )^{\frac{1}{2}}} d\varsigma < \sqrt{2} \int_{a}^{\frac{5a}{4}} \frac{\varsigma e^{-\frac{a^2}{4t}}}{(\varsigma^2 -a^2)^{\frac{1}{2}} } d\varsigma = \frac{3\sqrt{2}}{4} ae^{-\frac{a^2}{4t}} .
\end{equation*}

It follows that
$$ \int_{a}^\infty \frac{\varsigma e^{-\frac{\varsigma^2}{4t}}}{(\cosh \varsigma - \cosh a )^{\frac{1}{2}}} d\varsigma < \frac{3\sqrt{2}}{4} ae^{-\frac{a^2}{4t}} + \frac{5\sqrt{2\pi t}}{3} e^{-\frac{25a^2}{64t}} .$$

The second inequality in the statement of the lemma then follows directly from the elementary estimate $r e^{-\frac{r^2}{2}} \leq e^{-\frac{1}{2}} $.

When $a\geq 2\arcsinh 1 $ and $t\leq\frac{1}{20}$, one checks that the sum $\frac{3\sqrt{2}}{4t^{\frac{3}{2}}} ae^{-\frac{a^2}{4t}} + \frac{5\sqrt{2\pi }}{3t} e^{-\frac{25a^2}{64t}}$ is increasing in $t$. In particular, for $a= 2\arcsinh 1 $ and $t=\frac{1}{20}$, the sum is approximately $0.0000299 \cdots < \frac{3}{100000}$, which proves the lemma.
\end{proof}

In the following, we establish an upper bound for the hyperbolic heat kernel $K_{\mathrm{hyp}} (x,x;t)$ on the diagonal.

\begin{lem}
\label{lemestimatehyperbolicheatkernelondiagonal}
Let 
$$K_{\mathrm{hyp}} (x,y;t)= \sum_{l=1}^\infty e^{-t\lambda_{\mathrm{hyp},l} } \phi_{\mathrm{hyp} , l } (x) \phi_{\mathrm{hyp} , l }(y) $$
be the hyperbolic heat kernel on $C\cong \mathcal{H} / \Gamma$. Let $\mathrm{sys}(C)$ denote the systole of $(C,\mu_{\mathrm{KE}})$. 

Then for any $x\in C$ and $t\in \left( 0,\frac{1}{40(g-1)} \right] $, the following inequality holds:
$$ 1+ K_{\mathrm{hyp}} (x,x;t) \leq \frac{1}{2t} + \frac{6 \sqrt{g-1} }{ \sqrt{t}\cdot \mathrm{sys}(C)}  + \frac{g-1}{100} \cdot \left(1+\frac{4\arcsinh 1}{\mathrm{sys}(C)}\right) . $$
\end{lem}

\begin{proof}
Combining Lemma \ref{lemheatkernelformulaquotientstructure} and Lemma \ref{lemestimatehyperbolicheatkernelonhalfplane}, it follows that for every point $z \in\mathcal{H}$ and every $t>0$,
\begin{eqnarray*}
& & 1+K_{\mathrm{hyp}} (\rho_{\mathcal{H}} (z),\rho_{\mathcal{H}} (z);t) = 4\pi (g-1)\sum_{\sigma\in\Gamma} K_{\mathcal{H}} (z,\sigma z;2 (g-1)t) \\
& = & \sum_{\sigma\in\Gamma} \frac{4\sqrt{2} \pi (g-1)e^{-\frac{(g-1)t}{2}}}{(8\pi (g-1) t)^{\frac{3}{2}}} \int_{\mathrm{dist}_{\mathcal{H}} (z,\sigma z)}^\infty \frac{\varsigma e^{-\frac{\varsigma^2}{8(g-1)t}}}{(\cosh \varsigma - \cosh (\mathrm{dist}_{\mathcal{H}} (z,\sigma z)) )^{\frac{1}{2}}} d\varsigma .
\end{eqnarray*}

By selecting a maximal disjoint family of hyperbolic balls of radius $\arcsinh 1$, centered at points in $\Gamma z$, we obtain a sequence $\mathcal{S}_{\Gamma}=\{\sigma_j\}_{j=0}^\infty\subset \Gamma $ satisfying the following properties:
\begin{itemize}
    \item $\sigma_0=\mathrm{id}_{\mathcal{H}} $, and the sequence $\mathbf{d}_j =\mathrm{dist}_{\mathcal{H}} (\sigma_j z , z) $ is non-decreasing with respect to $j$.
    \item $\mathrm{dist}_{\mathcal{H}} (\sigma_j z ,\sigma_k z) \geq 2\arcsinh 1 $ for all $j\neq k$.
    \item For every $\sigma\in \Gamma $, there exists a $\sigma_j\in\mathcal{S}_{\Gamma} $ such that $\mathrm{dist}_{\mathcal{H}} (\sigma_j z ,\sigma z) < 2\arcsinh 1 $.
\end{itemize}

We now estimate $K_{\mathrm{hyp}} (\rho_{\mathcal{H}} (z),\rho_{\mathcal{H}} (z);t)$.

By definition, we have the following decomposition:
$$\Gamma z= \{z\} \bigcup \left( \Gamma z \cap \mathbf{B}_0 \setminus \{ z\} \right) \bigcup \bigcup_{j=1}^\infty \left( \Gamma z \cap \left( \mathbf{B}_j \setminus\mathbf{B}_0 \right) \right) ,$$
where $\mathbf{B}_j = B_{2\arcsinh1} (\sigma_j z)$, and $B_r (w)=\{ w'\in\mathcal{H} :\mathrm{dist}_{\mathcal{H}} (w,w') <r \} $.

Arguing as in Lemma \ref{lemestimatehyperbolicheatkernelonhalfplane}, we have
\begin{eqnarray*}
4\pi (g-1)  K_{\mathcal{H}} (z,z;2 (g-1)t) & = & \frac{4\sqrt{2} \pi (g-1)e^{-\frac{(g-1)t}{2}}}{(8\pi (g-1) t)^{\frac{3}{2}}} \int_{0}^\infty \frac{\varsigma e^{-\frac{\varsigma^2}{8(g-1)t}}}{(\cosh \varsigma - 1 )^{\frac{1}{2}}} d\varsigma \\
& < & \frac{4\sqrt{2} \pi (g-1)e^{-\frac{(g-1)t}{2}}}{(8\pi (g-1) t)^{\frac{3}{2}}} \cdot \sqrt{2} \cdot \int_{0}^\infty e^{-\frac{\varsigma^2}{8(g-1)t}} d\varsigma \\
& = & \frac{ e^{-\frac{(g-1)t}{2}}}{ 2 t} <\frac{1}{2t} ,
\end{eqnarray*}
and this gives the estimate for the first part.

For the second part, choose $\sigma_z\in\Gamma $ such that 
$$\mathrm{dist}_{\mathcal{H}} ( z,\sigma_z z) =\inf \left\{ \mathrm{dist}_{\mathcal{H}} ( z,\sigma z) ,\sigma\in\Gamma\setminus \left\{ \mathrm{id}_{\mathcal{H}} \right\}  \right\} .$$
By Lemma \ref{lemlocalstructureorbitcollartheorem}, if $\sigma'\in \mathbf{B}_j $, then $\sigma_j^{-1}\sigma'\in \langle\sigma_z\rangle$. Moreover, for any $k\geq 0$, it follows that 
$$ \mathrm{dist}_{\mathcal{H}} ( z,\sigma_z^{k+1} z) \geq \mathrm{dist}_{\mathcal{H}} ( z,\sigma_z^k z) + \mathrm{sys}(C) .$$
Hence there exists an integer $k_z\geq 0$ such that $\sigma_z^k z \in \mathbf{B}_0 $ if and only if $|k|\leq k_z $. If $\Gamma z \cap \mathbf{B}_0 \setminus\{z\} \neq \varnothing $, then $k_z\geq 1$, and
\begin{eqnarray*}
& & 4\pi (g-1) \sum_{\sigma z\in\Gamma z \cap \mathbf{B}_0 \setminus\{z\} } K_{\mathcal{H}} (z,\sigma z;2 (g-1)t) \\
& < & \frac{1}{t} \sum_{\sigma z\in\Gamma z \cap \mathbf{B}_0 \setminus\{z\} } \left( \frac{3 }{4\sqrt{e\pi }} \cdot e^{-\frac{\mathrm{dist}_{\mathcal{H}} (z,\sigma z)^2}{16(g-1)t}} + \frac{5 }{6} \cdot e^{-\frac{25\cdot\mathrm{dist}_{\mathcal{H}} (z,\sigma z)^2}{128(g-1)t}} \right) \\
& < & \frac{2}{t}\sum_{k=1}^{k_z} \frac{\int_{\mathrm{dist}_{\mathcal{H}} (z,\sigma_z^{k-1} z)}^{\mathrm{dist}_{\mathcal{H}} (z,\sigma_z^k z)} \left( \frac{3 }{4\sqrt{e\pi }} \cdot e^{-\frac{\vartheta^2}{16(g-1)t }} + \frac{5 }{6} \cdot e^{-\frac{25\cdot \vartheta^2}{128(g-1)t }} \right) d\vartheta }{\mathrm{dist}_{\mathcal{H}} (z,\sigma_z^{k} z) -\mathrm{dist}_{\mathcal{H}} (z,\sigma_z^{k-1} z)} \\
& < & \frac{2}{ \mathrm{sys}(C)\cdot t} \int_{0}^{\infty} \left( \frac{3 }{4\sqrt{e\pi }} \cdot e^{-\frac{\vartheta^2}{16(g-1)t }} + \frac{5 }{6} \cdot e^{-\frac{25\cdot \vartheta^2}{128(g-1)t }} \right) d\vartheta \\
& = & \frac{(9+4\sqrt{2e\pi })\sqrt{g-1}}{3\sqrt{e} \cdot \mathrm{sys}(C)\sqrt{t}} < \frac{6 \sqrt{g-1} }{ \sqrt{t}\cdot \mathrm{sys}(C)}  .
\end{eqnarray*}
Hence we obtain the claimed estimate for the second part.

It remains to consider the last part. By the construction, we have 
$$ \bigcup\limits_{j=0}^k B_{\arcsinh 1} (\sigma_j z) \subset B_{\mathbf{d}_k +\arcsinh 1} ( z) .$$
It follows that for any $k\geq 1$,
\begin{eqnarray*}
2(k+1)\pi \left( \sqrt{2} -1 \right) & = &\sum_{j=0}^{k} \mathrm{vol}_{\mathcal{H}} \left( B_{\arcsinh 1} (\sigma_j z) \right) \leq \mathrm{vol}_{\mathcal{H}} \left( B_{\mathbf{d}_k +\arcsinh 1} ( z) \right) \\
& = & 2\pi \left( \cosh\left(\mathbf{d}_k +\arcsinh 1 \right) -1 \right) \\
& \leq & \pi \left( \left(\sqrt{2} +1\right)e^{\mathbf{d}_k} +  \left( \sqrt{2} -1 \right)^3 -2 \right) .
\end{eqnarray*}
Hence for all $k\geq 1$, $ \mathbf{d}_k \geq \min\left\{ \log \left( 2k \right) -2\arcsinh1 ,2\arcsinh 1 \right\} $. Moreover, for any $\sigma z\in \Gamma z \cap \left( \mathbf{B}_k \setminus\mathbf{B}_0 \right)$, the distance
$$ \mathrm{dist}_{\mathcal{H}} (z,\sigma z) \geq \min\{ \log (2k)-4\arcsinh 1, 2\arcsinh 1 \} .$$
By Lemma \ref{lemlocalstructureorbitcollartheorem}, $\# (\Gamma z \cap \left( \mathbf{B}_j \setminus\mathbf{B}_0 \right)) \leq 1+\frac{4\arcsinh 1}{\mathrm{sys}(C)}$. Then Lemma \ref{lemestimatehyperbolicheatkernelonhalfplane} implies that for any $k\geq 1$, $t\in \left( 0,\frac{1}{40(g-1)} \right] $, and $\sigma z\in \Gamma z \cap \left( \mathbf{B}_k \setminus\mathbf{B}_0 \right)$,
\begin{equation*}
 \frac{1}{(2 (g-1) t)^{\frac{3}{2}}} \int_{\mathrm{dist}_{\mathcal{H}} (z,\sigma z)}^\infty \frac{\varsigma e^{-\frac{\varsigma^2}{8(g-1)t}}}{\sqrt{\cosh \varsigma - \cosh (\mathrm{dist}_{\mathcal{H}} (z,\sigma z)) }} d\varsigma \leq \frac{3}{100000},
\end{equation*}
and hence 
$$4\pi (g-1) K_{\mathcal{H}} (z,\sigma z;2 (g-1)t) < \frac{(g-1)e^{-\frac{(g-1)t}{2}}}{\sqrt{2\pi}} \cdot \frac{3}{100000} < \frac{g-1}{80000} .$$
It follows that
\begin{equation*}
4\pi (g-1)\sum_{\substack{
\sigma z\in \Gamma z \cap \left( \mathbf{B}_k \setminus\mathbf{B}_0 \right) \\
1 \le k \le 400}} K_{\mathcal{H}} (z,\sigma z;2 (g-1)t) < \frac{g-1}{200} \cdot \left(1+\frac{4\arcsinh 1}{\mathrm{sys}(C)}\right) .
\end{equation*}

Now we consider the case $k\geq 401$. By Lemma \ref{lemestimatehyperbolicheatkernelonhalfplane} again, for any $k\geq 401$, $t\in \left( 0,\frac{1}{40(g-1)} \right] $, and $\sigma z\in \Gamma z \cap \left( \mathbf{B}_k \setminus\mathbf{B}_0 \right)$,
\begin{eqnarray*}
& & \frac{1}{\sqrt{2(g-1)t}}\int_{\mathrm{dist}_{\mathcal{H}} (z,\sigma z)}^\infty \frac{\varsigma e^{-\frac{\varsigma^2}{8(g-1)t}}}{\sqrt{\cosh \varsigma - \cosh (\mathrm{dist}_{\mathcal{H}} (z,\sigma z)) }} d\varsigma \\
& \leq & \frac{3 }{\sqrt{2e}} \cdot e^{-\frac{\mathrm{dist}_{\mathcal{H}} (z,\sigma (z))^2}{16(g-1)t}} + \frac{5\sqrt{2\pi }}{3} e^{-\frac{25\mathrm{dist}_{\mathcal{H}} (z,\sigma (z))^2}{128(g-1)t}} \\
& \leq & \left( \frac{3 }{\sqrt{2e}} + \frac{5\sqrt{2\pi }}{3} e^{-\frac{85}{4}(\arcsinh 1)^2} \right) e^{-\frac{\mathrm{dist}_{\mathcal{H}} (z,\sigma (z))^2}{16(g-1)t}} < \frac{3}{2} \cdot e^{-\frac{\mathrm{dist}_{\mathcal{H}} (z,\sigma (z))^2}{16(g-1)t}} ,
\end{eqnarray*}
and hence
\begin{eqnarray*}
4\pi (g-1) K_{\mathcal{H}} (z,\sigma z;2 (g-1)t) & < & \frac{1}{2t\sqrt{2\pi}} \cdot \frac{3}{2} \cdot e^{-\frac{\mathrm{dist}_{\mathcal{H}} (z,\sigma (z))^2}{16(g-1)t}-\frac{(g-1)t}{2}} \\
& < & \frac{20(g-1)}{ \sqrt{2\pi}} \cdot \frac{3}{2} \cdot e^{-\frac{5\cdot\mathrm{dist}_{\mathcal{H}} (z,\sigma (z))^2}{2}} \\
& < & 20(g-1) \cdot e^{-\frac{5\cdot(\log (2k)-4\arcsinh 1)^2}{2}} .
\end{eqnarray*}
It follows that
\begin{eqnarray*}
& & 4\pi (g-1) \sum_{\substack{
\sigma z\in \Gamma z \cap \left( \mathbf{B}_k \setminus\mathbf{B}_0 \right) \\
 k \ge 401}} K_{\mathcal{H}} (z,\sigma z;2 (g-1)t) \\
 & < & \sum_{k=401}^\infty \left(1+\frac{4\arcsinh 1}{\mathrm{sys}(C)}\right) \cdot 20(g-1) \cdot e^{-\frac{5\cdot(\log (2k)-4\arcsinh 1)^2}{2}} \\
 & < & \left(1+\frac{4\arcsinh 1}{\mathrm{sys}(C)}\right) \cdot 20(g-1) \cdot \int_{400}^{\infty} e^{-\frac{5\cdot(\log (2\varsigma)-4\arcsinh 1)^2}{2}} d\varsigma \\
 & = & \left(1+\frac{4\arcsinh 1}{\mathrm{sys}(C)}\right) \cdot 20(g-1) \cdot \frac{e^{4\arcsinh 1 +\frac{1}{10}}}{ \sqrt{10}} \int_{\frac{\sqrt{10}}{2} (\log (800) -4\arcsinh 1-\frac{1}{10})}^{\infty} e^{-\varsigma^2} d\varsigma \\
 & < & \left(1+\frac{4\arcsinh 1}{\mathrm{sys}(C)}\right) \cdot 2\sqrt{10}(g-1) \cdot e^{4\arcsinh 1 +\frac{1}{10} }\int_0^{\infty} e^{-\frac{5}{2}\cdot (\log (800) -4\arcsinh 1-\frac{1}{10})^2-\varsigma^2} d\varsigma \\
 & < & \frac{g-1}{200} \cdot \left(1+\frac{4\arcsinh 1}{\mathrm{sys}(C)}\right).
\end{eqnarray*}
Here we used the formula 
$$\int_{\varsigma_0}^\infty e^{-a(\log (2\varsigma) -b)^2} d\varsigma = \frac{e^{b+\frac{1}{4a}}}{2\sqrt{a}} \int_{\sqrt{a}(\log (2\varsigma_0) -b-\frac{1}{2a})}^\infty e^{-\varsigma^2} d\varsigma ,\;\; \forall a,b,\varsigma_0 >0.$$
Now we can conclude that
\begin{equation*}
4\pi (g-1)\sum_{\substack{
\sigma z\in \Gamma z \cap \left( \mathbf{B}_k \setminus\mathbf{B}_0 \right) \\
 k \geq 1}} K_{\mathcal{H}} (z,\sigma z;2 (g-1)t) < \frac{g-1}{100} \cdot \left(1+\frac{4\arcsinh 1}{\mathrm{sys}(C)}\right) .
\end{equation*}

Combining the three parts, we conclude that
$$ 1+K_{\mathrm{hyp}} (\rho_{\mathcal{H}} (z),\rho_{\mathcal{H}} (z);t) \leq \frac{1}{2t} + \frac{6 \sqrt{g-1} }{ \sqrt{t}\cdot \mathrm{sys}(C)}  + \frac{g-1}{100} \cdot \left(1+\frac{4\arcsinh 1}{\mathrm{sys}(C)}\right) ,\;\; \forall z\in \mathcal{H} .$$
This completes the proof.
\end{proof}

\subsection{First eigenvalue and decompositions of hyperbolic surfaces}

Let $C $ be a connected compact Riemann surface of genus $g\geq 2$, and let $\mu_{\mathrm{KE}} = 4\pi (g-1) \mu_{\mathrm{hyp}} $ denote the K\"ahler metric on $C $ with constant curvature $-1 $. 

Recall that a separating multicurve on $C$ is a finite collection of disjoint simple closed geodesics whose union separates $C$ into at least two connected components. Let $\ell_{\Sigma} (C,\mu )$ denote the infimum of the total length of all geodesics in a separating multicurve. By definition, there exists a separating multicurve $\{\gamma_j \}_{j=1}^m$ such that $\ell_{\Sigma} (C,\mu ) = \sum\limits_{j=1}^m \ell (\gamma_j ) $. 

A classical result due to Schoen--Wolpert--Yau \cite{schwolyau1} states that the quotient $\frac{\lambda_{\mathrm{KE},1}}{\ell_{\Sigma}(C,\mu)}$ admits a positive lower bound depending only on the genus $g$; see also Dodziuk--Randol \cite{dora1} for a different proof. Recently, Wu--Xue \cite{wuxue1} proved that the lower bound can be taken to be the reciprocal of a quadratic polynomial in $g$. In this subsection, we present an improvement of this estimate, and our method is different from that of Wu–Xue. More precisely, our approach is based on combining conformal transformations with Cheeger’s inequality. Similar conformal techniques have been used in related estimates \cite{grolanach1}.

We start with Cheeger’s inequality for surfaces.

\begin{thm}[Cheeger's inequality for surfaces]
\label{thmcheegerinequality}
Let $ C $ be a connected compact Riemann surface equipped with a K\"ahler metric $\mu $, and let $\lambda^{\mathrm{Rm}}_{\mu ,1} >0$ be the first eigenvalue of the Riemannian Laplacian $\Delta^{\mathrm{Rm}}_{\mu}= {}\star_{\mu} d\star_{\mu} d $. Then we have
$$ \lambda^{\mathrm{Rm}}_{\mu ,1} \geq \frac{1}{4 } \cdot \inf_{U\subset C}\frac{\ell_\mu^2 (\partial U)}{\min\left\{\mu (U) , \mu (C\setminus U) \right\}^2} ,$$
where $U$ ranges over the set of all connected open subsets of $C$ with connected complement $C\setminus U$ and smooth boundary $\partial U$, $\ell_\mu (\partial U)$ is the length of $\partial U$, and $ \mu(U)$ is the measure of $U$.
\end{thm}

\begin{proof}
See \cite[Theorem 8.3.3, Lemma 8.3.6]{bu1}.
\end{proof}

The following weighted version of Cheeger's inequality is well known to experts. For completeness, we derive it from the classical Cheeger's inequality.

\begin{cor}
\label{cormodifiedcheegerinequality}
Let $ C $ be a connected compact Riemann surface equipped with a K\"ahler metric $\mu $, let $u:C\to [1,\infty )$ be a continuous function, let $\mu_u = u\cdot \mu$ be the weighted measure, and let $\lambda_{\mu ,1} >0$ denote the first eigenvalue of the Laplacian $\Delta_{\mu}= \frac{dd^c}{\mu} $. Then we have
$$ \lambda_{\mu ,1} \geq \frac{1}{8\pi } \cdot \inf_{U\subset C}\frac{\ell_{\mu_u}^2 (\partial U)}{\min\left\{\mu_u (U) , \mu_u (C\setminus U) \right\}^2} ,$$
where $U$ ranges over the set of all connected open subsets of $C$ with connected complement $C\setminus U$ and smooth boundary $\partial U$, $\ell_{\mu_u} (\partial U)$ is the weighted length of $\partial U$, and $ \mu_u (U)$ is the weighted measure of $U$.
\end{cor}

\begin{remark}
Note that $\Delta_{\mu} = \frac{1}{2\pi} \Delta^{\mathrm{Rm}}_{\mu} $, and hence $\lambda_{\mu ,1} = \frac{1}{2\pi} \lambda_{\mu ,1}^{\mathrm{Rm}}$.
\end{remark}

\begin{proof}
First, assume that $u$ is smooth. Then $\mu_u $ is also a K\"ahler metric, and hence the classical Cheeger's inequality implies that
$$ \lambda_{\mu_u ,1} \geq \frac{1}{8\pi } \cdot \inf_{U\subset C}\frac{\ell_{\mu_u}^2 (\partial U)}{\min\left\{\mu_u (U) , \mu_u (C\setminus U) \right\}^2} .$$

Now we will prove that $\lambda_{\mu ,1}\geq \lambda_{\mu_u ,1} $. Since $u\geq 0$, we have $\mu_u \geq \mu $. For any smooth function $\eta$, set $\eta_{\mathrm{aver} }=\mu (C)^{-1} \int_C \eta\, d\mu $ and $\eta_{\mathrm{aver},u}=\mu_u (C)^{-1} \int_C \eta\, d\mu_u$. By the Rayleigh principle for eigenvalues \cite[Theorem 8.2.1]{bu1}, we have
\begin{eqnarray*}
\lambda_{\mu_u ,1} & = & \inf_{\substack{\eta\in C^\infty (C;\mathbb{R}) \setminus \mathbb{R}}} \frac{-\int_C \eta\wedge dd^c \eta}{\int_C \left(\eta - \eta_{\mathrm{aver},u} \right)^2 d\mu_u } \leq \inf_{\substack{\eta\in C^\infty (C;\mathbb{R}) \setminus \mathbb{R}}} \frac{-\int_C \eta\wedge dd^c \eta}{\int_C \left(\eta - \eta_{\mathrm{aver},u} \right)^2 d\mu } \\
& \leq & \inf_{\substack{\eta\in C^\infty (C;\mathbb{R}) \setminus \mathbb{R}}} \frac{-\int_C \eta\wedge dd^c \eta}{\int_C \left(\eta - \eta_{\mathrm{aver}} \right)^2 d\mu } =\lambda_{\mu ,1} .
\end{eqnarray*}

In general, approximate the continuous function $u$ by a sequence of smooth functions $u_j$ converging uniformly to $u$ on $C$, with $u_j\geq 1$, $\forall j\in \mathbb{N}$. By the argument above, $\lambda_{\mu ,1}\geq \lambda_{\mu_{u_j} ,1} $ for all $ j\in\mathbb{N}$. Thus, for every $j\in\mathbb{N}$,
$$ \lambda_{\mu ,1} \geq \frac{1}{8\pi } \cdot \inf_{U\subset C}\frac{\ell_{\mu_{u_j}}^2 (\partial U)}{\min\left\{\mu_{u_j} (U) , \mu_{u_j} (C\setminus U) \right\}^2} .$$
Because $u_j \to u$ uniformly on $C$, we obtain
$$ \lim_{j\to\infty} \left( \inf_{U\subset C}\frac{\ell_{\mu_{u_j}}^2 (\partial U)}{\min\left\{\mu_{u_j} (U) , \mu_{u_j} (C\setminus U) \right\}^2}\right) = \inf_{U\subset C}\frac{\ell_{\mu_u}^2 (\partial U)}{\min\left\{\mu_u (U) , \mu_u (C\setminus U) \right\}^2} ,$$
which completes the proof.
\end{proof}

The following elementary isoperimetric inequalities are well known to experts.

\begin{lem}
\label{lemfirsteigenvalueestimatediscannulus}
Let $C $ be a connected smooth complex curve of genus $g\geq 2$, equipped with the unique K\"ahler metric $\mu_{\mathrm{KE}} $ of constant curvature $-1 $, and let $U\subset C$ be a connected open subset with smooth boundary $\partial U$. Assume that the complement $C\setminus U$ is also connected. If $\partial U$ has a null-homotopic connected component, or two connected components that are freely homotopic, then $\ell (\partial U) \geq \mu_{\mathrm{KE}} (U)$, where $\ell (\partial U)=\ell_{\mu_{\mathrm{KE}}} (\partial U)$ is the length of $\partial U$ in $(C,\mu_{\mathrm{KE}} )$. Combining this with Theorem \ref{thmcheegerinequality}, we obtain
$$ \lambda_{\mathrm{KE} ,1} \geq \min\left\{\frac{1}{8\pi} , \frac{(\ell_{\Sigma} (C,\mu_{\mathrm{KE}} ) )^2}{32\pi^3 (g-1)^2} \right\}. $$
\end{lem}

\begin{proof}
See \cite[(25)]{wuxue1} and \cite[Lemma 8, Lemma 9]{wuxue2}.
\end{proof}

Now assume that $\ell_{\Sigma} (C,\mu_{\mathrm{KE}} )\leq 1$. In this case, we construct a suitable conformal perturbation of $\mu_{\mathrm{KE}}$ and then apply Corollary \ref{cormodifiedcheegerinequality}.

Let $\{\gamma_j\}_{j=1}^m$ be the set of all simple closed geodesics on $(C,\mu_{\mathrm{KE}} )$ of length $ \leq 1$. Set $\epsilon_j = \max\{ \ell (\gamma_j), \ell_{\Sigma} (C,\mu_{\mathrm{KE}} ) \} $, where $\ell (\gamma_j) $ is the length of $\gamma_j$ in $(C,\mu_{\mathrm{KE}} )$. Now we use the structure of K\"ahler collars $\mathscr{C}(\gamma_j )$ to construct the weight function. By Lemma \ref{lemkahlercollar}, each $\gamma_j$ admits a K\"ahler collar $\mathscr{C}(\gamma_j)$ with holomorphic coordinate 
$$z_j : \mathscr{C} (\gamma_j ) \to \overline{\mathbb{D}_{e^{-\nu (\gamma_j)}} (0)} \setminus \mathbb{D}_{e^{\nu (\gamma_j ) -\frac{2\pi^2}{\ell (\gamma_j )} }} (0) $$
such that 
$$ z_j^* (\mu_{\mathrm{KE}} ) = \frac{ \ell^2 (\gamma_j ) \mu_{\mathrm{Euc}} }{ 4\pi^2 |z_j|^2 \sin^2 ( \frac{\ell (\gamma_j )}{2\pi} \log |z_j | ) } ,$$
where $\mu_{\mathrm{Euc}} = \frac{idz\wedge d\bar{z}}{2} $ is the standard Euclidean metric on $\mathbb{C}$, and $\nu (\gamma_j) \in (0, \frac{\pi^2}{2\ell (\gamma_j)}) $ is determined by $\nu (\gamma_j) \ell (\gamma_j) = {2\pi} \arcsin \left( \tanh \frac{\ell (\gamma_j )}{2} \right) $. On each collar $\mathscr{C}(\gamma_j )$, we write
\begin{eqnarray*}
  u_j (z_j ) & = & \max\left\{ 1,\frac{\sin^2 \left( \frac{\ell(\gamma_j)}{2\pi} \log |z_j| \right)}{\epsilon_j} \right\} \\
  & = & \left\{
\begin{aligned}
1, \;\;\;\;\;\;\;\;\;\;\;\;\;\;\; & \;\;\;\;\;\; \left|\log |z_j| + \frac{\pi^2}{\ell(\gamma_j)}\right| > \frac{\pi^2 - 2\pi \arcsin \sqrt{\epsilon_j}}{\ell(\gamma_j)} ;\\
\frac{\sin^2 \left( \frac{\ell(\gamma_j)}{2\pi} \log |z_j| \right)}{\epsilon_j}, \; & \;\;\;\;\;\; \left|\log |z_j| + \frac{\pi^2}{\ell(\gamma_j)}\right| \leq \frac{\pi^2 - 2\pi \arcsin \sqrt{\epsilon_j}}{\ell(\gamma_j)} ,
\end{aligned}
\right.
\end{eqnarray*}
and
\begin{eqnarray*}
  u (x) & = & \left\{
\begin{aligned}
u_j (z_j (x)), \;\;\;\;\; & \;\;\;\;\;\; x\in\mathscr{C} (\gamma_j) ;\\
1, \;\;\;\;\;\;\;\;\;\;\; & \;\;\;\;\;\; x\in C\setminus \bigcup_{j=1}^m \mathscr{C} (\gamma_j) . 
\end{aligned}
\right.
\end{eqnarray*}
Thus $u \geq 1$ on $C$, and $u$ is Lipschitz continuous with respect to $\mu_{\mathrm{KE}}$.

Let $\mu_{\mathrm{KE};u}=u\,\mu_{\mathrm{KE}}$ be the weighted measure. Then $\mu_{\mathrm{KE};u}$ gives a Lipschitz Riemannian metric on $C$. Now we list some basic properties of the weighted measure $\mu_{\mathrm{KE};u} $.

\begin{lem}
\label{lembasicpropertiesweightedmeasure}
Let $C$, $\mu_{\mathrm{KE}}$, $\{\gamma_j\}_{j=1}^m$, $u_j$, $u$, and $\mu_{\mathrm{KE};u} $ be as above. Then the following properties hold:
\begin{enumerate}[(1)]
    \item $\mu_{\mathrm{KE};u} \geq \mu_{\mathrm{KE}} $ and $\mu_{\mathrm{KE};u} |_{C\setminus \bigcup_{j=1}^m \mathscr{C} (\gamma_j)} = \mu_{\mathrm{KE}} |_{C\setminus \bigcup_{j=1}^m \mathscr{C} (\gamma_j)} $.
    \item $\mu_{\mathrm{KE};u} (\mathscr{C} (\gamma_j) ) < \mu_{\mathrm{KE}} (\mathscr{C} (\gamma_j)) +\pi $, $j=1,\cdots ,m$.
    \item The distance $\mathrm{dist}_{\mu_{\mathrm{KE}}} \left( \{ u\neq 1 \} , C\setminus \bigcup_{j=1}^m \mathscr{C} (\gamma_j) \right) \geq \log 2 $.
\end{enumerate}
\end{lem}

\begin{proof}
(1) is immediate from the definition of $u$. We now prove (2).

By the construction,
$$\{ u\neq 1 \}\cap \mathscr{C} (\gamma_j) =\left\{\left|\log |z_j| + \frac{\pi^2}{\ell(\gamma_j)}\right| \leq \frac{\pi^2 - 2\pi \arcsin \sqrt{\epsilon_j}}{\ell(\gamma_j)}\right\} ,$$
so that
\begin{eqnarray*}
\mu_{\mathrm{KE};u} (\mathscr{C} (\gamma_j) ) -\mu_{\mathrm{KE}} (\mathscr{C} (\gamma_j)) & = & \int_{\mathscr{C}(\gamma_j)} (u_j -1) \mu_{\mathrm{KE}} < \int_{\{ u\neq 1 \}\cap \mathscr{C} (\gamma_j) } u_j \mu_{\mathrm{KE}} \\
& = & \int_{\left\{\left|\log |z_j| + \frac{\pi^2}{\ell(\gamma_j)}\right| \leq \frac{\pi^2 - 2\pi \arcsin \sqrt{\epsilon_j}}{\ell(\gamma_j)}\right\}} \frac{\ell (\gamma_j)^2  \mu_{\mathrm{Euc}}}{4\pi^2 |z_j|^2 \epsilon_j} \\
& = & \frac{(\pi - 2 \arcsin \sqrt{\epsilon_j} )\cdot\ell (\gamma_j)}{\epsilon_j} <\pi .
\end{eqnarray*}

For (3), it suffices to check the distance $\mathrm{dist}_{\mu_{\mathrm{KE}}} \left( \{ u\neq 1 \}\cap \mathscr{C} (\gamma_j) , \partial \mathscr{C} (\gamma_j) \right) \geq \log 2 $. By Proposition \ref{propbasicpropertiesriemanniankahlercollars}, we have
\begin{eqnarray*}
\mathrm{dist}_{\mu_{\mathrm{KE}}} \left( \{ u\neq 1 \}\cap \mathscr{C} (\gamma_j) , \partial \mathscr{C} (\gamma_j) \right) & \geq & \log \left( \frac{2\pi \arcsin \sqrt{\epsilon_j }}{\ell (\gamma_j) \nu (\gamma_j)} \right) \\
& \geq & \log \left( \frac{ \arcsin \sqrt{\ell (\gamma_j ) }}{\arcsin \left( \tanh \left(  {\ell (\gamma_j )}/{2} \right) \right)} \right) \\
& > & \log \left( \frac{ 2}{ \sqrt{\ell (\gamma_j ) } } \right) \geq \log 2 .
\end{eqnarray*}
This proves (3), and hence the lemma.
\end{proof}

Building on the previous analysis, we now consider the isoperimetric inequalities on $(C,\mu_{\mathrm{KE};u})$.

\begin{lem}
\label{lemmodifiedfirsteigenvalueestimatediscannulus}
Let $C$, $\mu_{\mathrm{KE}}$, $\{\gamma_j\}_{j=1}^m$, $u_j$, $u$, and $\mu_{\mathrm{KE};u} $ be as in Lemma \ref{lembasicpropertiesweightedmeasure}. Let $U\subset C$ be an open subset of $C$ with smooth boundary $\partial U$. Assume that $\ell_{\Sigma} (C,\mu_{\mathrm{KE}} ) \leq 1 $, $\ell_{{\mathrm{KE};u}} (\partial U) \leq 2\log 2$, and that $U$ is homeomorphic to an open disk or a cylinder. Then $\ell_{\mathrm{KE};u} (\partial U) \geq \frac{2\sqrt{\ell_{\Sigma} (C,\mu_{\mathrm{KE}} )}}{\pi} \mu_{\mathrm{KE};u} (U)$, where $ \ell_{\mathrm{KE};u} (\partial U)$ is the length of $\partial U$ in $(C,\mu_{\mathrm{KE};u} )$. 
\end{lem}

\begin{proof}
We first assume that $\partial U \nsubseteq \mathscr{C} (\gamma_j) $ for any $\gamma_j$. Since $\arcsinh 1 >\log 2 $, Theorem \ref{thmcollarthm} shows that each connected component of $\partial U$ is null-homotopic. By Lemma \ref{lembasicpropertiesweightedmeasure}, this implies $\partial U \subset \{ u= 1 \} $. Hence $U\subset \{ u= 1 \}$. In this case, the inequality follows from Lemma \ref{lemfirsteigenvalueestimatediscannulus}.

Now assume that $\partial U \subset \mathscr{C} (\gamma_j) $ for some $\gamma_j$. Then either $ U \subset \mathscr{C} (\gamma_j) $ or $C\setminus U \subset \mathscr{C} (\gamma_j) $. Since both $\mathscr{C}(\gamma_j)$ and $U$ have first Betti number $\leq 1$, whereas $C$ has first Betti number $2g \geq 4$, the latter case is impossible. Hence $U\subset \mathscr{C}(\gamma_j)$. Set 
\begin{eqnarray*}
  \Psi (z_j ) & = & \left\{
\begin{aligned}
-\log \left( -\sin \left( \frac{\ell (\gamma_j)}{2\pi} \log |z_j | \right) \right), \;\;\;\;\;\;\;\;\;\;\;\; &  \;\;\;\;\;\; z_j\in \{ u_j = 1 \} ;\\
\frac{\ell^2 (\gamma_j) \tan \theta_j}{8\pi^2 \theta_j } \cdot \left(\log |z_j| + \frac{\pi^2}{\ell(\gamma_j)} \right)^2 + \vartheta_j, \; &  \;\;\;\;\;\; z_j\in \{ u_j \neq 1 \}  ,
\end{aligned}
\right.
\end{eqnarray*}
where $ \theta_{j} = \frac{\pi}{2} - \arcsin \sqrt{\epsilon_j} $, and $\vartheta_j = -\log \sqrt{\epsilon_j} - \frac{\theta_j \tan \theta_j}{2} $. Since
$$\{ u_j \neq 1 \} =\left\{ \left|\log |z_j| + \frac{\pi^2}{\ell(\gamma_j)}\right| < \frac{2\pi \theta_{j}}{\ell(\gamma_j)} \right\} ,$$
we have
\begin{eqnarray*}
 d \Psi (z_j ) & = & \left\{
\begin{aligned}
\frac{-\ell (\gamma_j )\left( \bar{z}_j dz_j + z_j d\bar{z}_j  \right)}{4\pi |z_j |^2\tan \left( \frac{\ell (\gamma_j)}{2\pi} \log |z_j | \right)}, \;\;\;\;\;\;\;\;\;\;\;\;\;\;\;\;\;\;\; &  \; z_j\in \{ u_j = 1 \} ;\\
\frac{\ell^2 (\gamma_j) \left(\log |z_j| + \frac{\pi^2}{\ell(\gamma_j)} \right) \tan \theta_j}{8\pi^2 \theta_j |z_j |^2 }  \cdot \left( \bar{z}_j dz_j + z_j d\bar{z}_j  \right), \; & \; z_j\in \{ u_j \neq 1 \}  .
\end{aligned}
\right.
\end{eqnarray*}
Thus $\Psi$ and $d\Psi $ are Lipschitz continuous on $\mathscr{C} (\gamma_j)$. Moreover,
\begin{eqnarray*}
 d\star d \Psi (z_j ) & = & \left\{
\begin{aligned}
\mu_{\mathrm{KE}}, \;\;\;\;\;\;\;\;\;\;\;\; & \; z_j\in \{ u_j = 1 \} ;\\
\frac{\ell^2 (\gamma_j) \tan \theta_j }{4\pi^2 \theta_j |z_j |^2 }   \mu_{\mathrm{Euc}}, \; & \; z_j\in \{ u_j \neq 1 \}  ,
\end{aligned}
\right.
\end{eqnarray*}
where $\star$ denotes the Hodge $\star$-operator. Note that on $1$-forms its action depends only on the complex structure.

Since $\cos \theta_j = \sqrt{\epsilon_j }$, for $z_j\in \{ u_j \neq 1 \} $, we obtain
\begin{eqnarray*}
d\star d \Psi (z_j ) & = & \frac{\ell^2 (\gamma_j) \tan \theta_j }{4\pi^2 \theta_j |z_j |^2 } \mu_{\mathrm{Euc}} = \frac{\epsilon_j \tan \theta_j}{\theta_j } \mu_{\mathrm{KE};u} \\
& = & \frac{\sqrt{\epsilon_j } \sin \theta_j }{\theta_j }  \mu_{\mathrm{KE};u} \geq \frac{2 \sqrt{\epsilon_j } }{\pi }  \mu_{\mathrm{KE};u} .
\end{eqnarray*}
Applying the argument in Lemma \ref{lemfirsteigenvalueestimatediscannulus} then yields
$$\ell_{\mathrm{KE};u} (\partial U) \geq \frac{2\sqrt{\epsilon_j } }{\pi} \mu_{\mathrm{KE};u} (U) \geq \frac{2\sqrt{\ell_{\Sigma} (C,\mu_{\mathrm{KE}} )}}{\pi} \mu_{\mathrm{KE};u} (U) , $$
which is our assertion.
\end{proof}

By combining the arguments above, we obtain the following estimate:

\begin{prop}
\label{propfirsteigenvalueestimate}
Let $C $ be a connected smooth complex curve of genus $g\geq 2$, and let $\mu_{\mathrm{KE}} $ denote the unique K\"ahler metric on $C $ with constant curvature $-1 $. Then at least one of the following holds:
\begin{enumerate}[(1)]
\item $\ell_{\Sigma} (C,\mu_{\mathrm{KE}} ) \geq 1 $ and $ \lambda_{\mathrm{KE} ,1} \geq \min\left\{\frac{1}{8\pi} , \frac{(\ell_{\Sigma} (C,\mu_{\mathrm{KE}} ) )^2}{32\pi^3 (g-1)^2} \right\} $.
\item $\ell_{\Sigma} (C,\mu_{\mathrm{KE}} ) \leq 1 $ and $ \lambda_{\mathrm{KE} ,1} \geq \min\left\{\frac{\ell_{\Sigma} (C,\mu_{\mathrm{KE}} ) }{2\pi^3} , \frac{1}{52\pi^3 (g-1)^2} \right\} $.
\item $\ell_{\Sigma} (C,\mu_{\mathrm{KE}} ) \leq 1 $ and there exists a minimal separating multicurve $\{\gamma_{j} \}_{k=1}^{m'}$ splitting $C$ in to two parts $C_1$, $C_2$, such that
$$ \lambda_{\mathrm{KE} ,1} \geq \frac{ \sum_{j=1}^{m'} \ell ( \gamma_j )}{8\pi^3} \cdot \frac{1}{\min\left\{ 7(g_1 -1) +5m' , 7(g_2 -1) +5m' ,\frac{7(g-1)}{2} \right\}^2 } , $$
where $g_1$, $g_2$ are the genera of $C_1$ and $C_2$, respectively.
\end{enumerate}
\end{prop}

\begin{remark}
As a direct consequence of Proposition \ref{propfirsteigenvalueestimate}, we obtain $\lambda_{\mathrm{KE} ,1}^{\mathrm{Rm}} \geq \frac{\ell_{\Sigma} (C,\mu_{\mathrm{KE}} )}{49\pi^2 (g-1)^2} $. Note that the classical estimate for Bers' constant \cite[Theorem 5.1.4]{bu1} yields the upper bound $\ell_{\Sigma} (C,\mu_{\mathrm{KE}} )\leq 78 (g-1)$ \cite[Lemma 6]{wuxue1}. Moreover, an improved bound of the form $\ell_{\Sigma} (C,\mu_{\mathrm{KE}} )\leq 38\log g$ was obtained in \cite[Proposition 5]{hewu1}.
\end{remark}

\begin{proof}
If $\ell_{\Sigma} (C,\mu_{\mathrm{KE}} ) \geq 1 $, then property (1) follows directly from Lemma \ref{lemfirsteigenvalueestimatediscannulus}. Now we assume that $\ell_{\Sigma} (C,\mu_{\mathrm{KE}} ) \leq 1 $. 

Let $\{\gamma_j\}_{j=1}^m$, $u_j$, $u$, and $\mu_{\mathrm{KE};u} $ be as in Lemma \ref{lembasicpropertiesweightedmeasure}. By Theorem \ref{thmcollarthm}, we have
$$ \mu_{\mathrm{KE};u} (C ) < \mu_{\mathrm{KE} } (C ) + m\pi \leq 7(g-1)\pi .$$

Let $U\subset C$ be a connected open subset with smooth boundary $\partial U$, such that $C\setminus U$ is also connected. 

First, consider the case $\ell_{\mu_{\mathrm{KE};u}} (\partial U) \geq 2\log 2$. In this case,
\begin{eqnarray*}
\frac{1}{8\pi} \cdot \frac{\ell_{\mu_{\mathrm{KE};u}}^2 (\partial U)}{\min\left\{\mu_{\mathrm{KE};u} (U) , \mu_{\mathrm{KE};u} (C\setminus U) \right\}^2} & \geq & \frac{(2\log 2)^2}{2\pi \cdot (7\pi (g-1))^2} \\
& = & \frac{(2\log 2)^2}{98\pi^3 (g-1)^2} > \frac{1}{52\pi^3 (g-1)^2} .
\end{eqnarray*}

Now suppose that $\ell_{\mu_{\mathrm{KE};u}} (\partial U) \leq 2\log 2$. If $\partial U$ has a null-homotopic connected component, or two connected components that are freely homotopic, then Lemma \ref{lemmodifiedfirsteigenvalueestimatediscannulus} implies
\begin{eqnarray*}
\frac{1}{8\pi} \cdot \frac{\ell_{\mu_{\mathrm{KE};u}}^2 (\partial U)}{\min\left\{\mu_{\mathrm{KE};u} (U) , \mu_{\mathrm{KE};u} (C\setminus U) \right\}^2} & \geq & \frac{1}{8\pi} \cdot \frac{\ell_{\mu_{\mathrm{KE};u}}^2 (\partial U)}{ \mu_{\mathrm{KE};u} (U)^2} \\
& \geq & \frac{1}{8\pi} \cdot \frac{4\ell_{\Sigma} (C,\mu_{\mathrm{KE}} ) }{\pi^2} = \frac{\ell_{\Sigma} (C,\mu_{\mathrm{KE}} ) }{2\pi^3} .
\end{eqnarray*}

Otherwise, every connected component of $\partial U$ is not null-homotopic, and no two connected components of $\partial U$ are homotopically equivalent.

Let $\{\mathbf{c}_j\}_{j=1}^{m'}$ denote the set of connected components of $\partial U$. Since $\ell_{\mu_{\mathrm{KE};u}} (\partial U) \leq 2\log 2<2\arcsinh 1$, it follows that, after reordering the short closed geodesics $\{\gamma_j\}_{j=1}^m$, for each $j=1,\cdots,m'$, we have $\mathbf{c}_j \subset \mathscr{C} (\gamma_j) $. Thus $\{\gamma_j \}_{j=1}^{m'} $ separates $C$ into two domains $C_1$ and $C_2$. Let $g_1$, $g_2$ denote the genera of $C_1$ and $C_2$, respectively. By the Gauss-Bonnet formula, $\mu_{\mathrm{KE}} (C_1) = 2\pi (2g_1 +m'-2) $, and $\mu_{\mathrm{KE}} (C_2) = 2\pi (2g_2 +m'-2) $. Moreover, by Theorem \ref{thmcollarthm}, we have $\# \{\gamma_j : \gamma_j\subset C_1 ,1\leq j \leq m\} \leq 3(g_1 -1) +m'$, and $\# \{\gamma_j : \gamma_j\subset C_2 ,1\leq j \leq m\} \leq 3(g_2 -1) +m'$. 

A straightforward computation shows that
\begin{eqnarray*}
\mu_{\mathrm{KE}} (\mathscr{C}(\gamma_j) \cap C_1 ) & = & \mu_{\mathrm{KE}} (\mathscr{C}(\gamma_j) \cap C_2 ) \\
& = & \int_{0}^{w(\gamma_j)} \ell (\gamma_j) \cosh rdr = \frac{\ell (\gamma_j)}{\sinh \frac{\ell (\gamma_j)}{2}} <2 ,
\end{eqnarray*}
where $w(\gamma_j ) >0$ satisfies $\sinh (w(\gamma_j ))\sinh \left( \frac{\ell (\gamma_j)}{2} \right) =1 $. It follows that
\begin{eqnarray*}
& & \min\left\{\mu_{\mathrm{KE};u} (U) , \mu_{\mathrm{KE};u} (C\setminus U) \right\} \\
& \leq & \min\left\{\mu_{\mathrm{KE};u} \left( C_1 \bigcup \left( \bigcup\limits_{j=1}^{m'} \mathscr{C} (\gamma_j) \right) \right) , \right.\\
& & \;\;\;\;\;\;\;\;\;\;\;\;\;\;\;\;\;\;\;\;\;\;\;\;\;\; \left. \mu_{\mathrm{KE};u} \left( C_2 \bigcup \left( \bigcup\limits_{j=1}^{m'} \mathscr{C} (\gamma_j) \right) \right) , \frac{1}{2}\mu_{\mathrm{KE};u} (C ) \right\} \\
& \leq & \min\left\{\mu_{\mathrm{KE}} \left( C_1 \right) + \sum_{j=1}^{m'} \mu_{\mathrm{KE}} (\mathscr{C}(\gamma_j) \cap C_2 ) + (3(g_1 -1) +2m')\pi , \right.\\
& & \;\; \left. \mu_{\mathrm{KE}} \left( C_2 \right) + \sum_{j=1}^{m'} \mu_{\mathrm{KE}} (\mathscr{C}(\gamma_j) \cap C_1 ) + (3(g_2 -1) +2m')\pi , \frac{1}{2}\mu_{\mathrm{KE};u} (C ) \right\} \\
& \leq & \min\left\{ (7 (g_1 -1) +5m' )\pi , (7(g_2 -1) +5m' )\pi ,\frac{7(g-1)\pi}{2} \right\} .
\end{eqnarray*}

We now consider $\ell_{\mu_{\mathrm{KE};u}} (\partial U) $. By the construction, we have $\sum_{j=1}^{m'} \ell (\mathbf{\gamma}_j ) \geq \ell_{\Sigma} (C,\mu_{\mathrm{KE}} ) $. Moreover, if $\ell (\gamma_j )\leq 1$, we have $ \ell_{\mu_{\mathrm{KE};u}} (\gamma_j ) = \frac{\ell (\gamma_j )}{\max\{\sqrt{\ell_{\Sigma} (C,\mu_{\mathrm{KE}} )} ,\;  \sqrt{\ell (\gamma_j )} \} }$. 

Suppose $\ell (\gamma_j )\leq \ell_{\Sigma} (C,\mu_{\mathrm{KE}} )$ for all $j=1,\cdots ,m'$. Then
\begin{eqnarray*}
\ell_{\mu_{\mathrm{KE};u}}^2 (\partial U) & = & \left( \sum_{j=1}^{m'} \ell_{\mu_{\mathrm{KE};u}} (\mathbf{c}_j ) \right)^2 \geq \left( \sum_{j=1}^{m'} \ell_{\mu_{\mathrm{KE};u}} (\gamma_j ) \right)^2 \\
& \geq & \left(  \frac{\sum_{j=1}^{m'} \ell (\mathbf{\gamma}_j )}{ \sqrt{\ell_{\Sigma} (C,\mu_{\mathrm{KE}} )} } \right)^2 \geq \left(  \frac{\sum_{j=1}^{m'} \ell (\mathbf{\gamma}_j )}{ \sqrt{\sum_{j=1}^{m'} \ell (\mathbf{\gamma}_j )}} \right)^2 = \sum_{j=1}^{m'} \ell (\mathbf{\gamma}_j ) .
\end{eqnarray*}

If there exists some $j$ with $\ell (\gamma_j ) > \ell_{\Sigma} (C,\mu_{\mathrm{KE}} )$, then we have $ \ell_{\mu_{\mathrm{KE};u}} (\gamma_j ) \geq \sqrt{\ell  (\gamma_j )} \geq \sqrt{\ell_{\Sigma} (C,\mu_{\mathrm{KE}} )}$. Without loss of generality, we may assume that $\ell (\gamma_j ) > \ell_{\Sigma} (C,\mu_{\mathrm{KE}} )$ if and only if $m''\leq j\leq m'$ for some index $m''$ with $1 \leq m'' \leq m'-1$. Hence
\begin{eqnarray*}
\ell_{\mu_{\mathrm{KE};u}}^2 (\partial U) & = & \left( \sum_{j=1}^{m'} \ell_{\mu_{\mathrm{KE};u}} (\mathbf{c}_j ) \right)^2 \geq \left( \sum_{j=1}^{m'} \ell_{\mu_{\mathrm{KE};u}} (\gamma_j ) \right)^2 \\
& \geq & \left(  \frac{\sum_{j=1}^{m''-1} \ell (\gamma_j )}{ \sqrt{\ell_{\Sigma} (C,\mu_{\mathrm{KE}} )} } + \sum_{j=m''}^{m' } \sqrt{\ell (\gamma_j )} \right)^2 \\
& > & \sqrt{\frac{\ell (\gamma_{m'})}{\ell_{\Sigma} (C,\mu_{\mathrm{KE}} )}} \cdot \sum_{j=1}^{m''-1} \ell (\gamma_j ) +  \sum_{j=m''}^{m' }  \ell (\gamma_j ) \geq \sum_{j=1}^{m'} \ell (\gamma_j ) .
\end{eqnarray*}
It follows that for any $\epsilon>0$, the estimates in this proposition hold if $\lambda_{\mathrm{KE} ,1}$ is replaced by $\lambda_{\mathrm{KE} ,1}+\epsilon$. Letting $\epsilon\to 0$, the estimates also hold for $\lambda_{\mathrm{KE} ,1}$. Note that there are only finitely many minimal separating multicurves on $C$. This completes the proof.
\end{proof}

\subsection{Eigenfunctions of small eigenvalues}

Let $C $ be a connected smooth complex curve of genus $g\geq 2$, and let $\mu_{\mathrm{KE}} $ denote the unique K\"ahler metric on $C $ with constant curvature $-1 $.

We consider the sup-norm of a $\mathbb{R}$-valued eigenfunction $\Phi$ associated with a small eigenvalue $\lambda = \epsilon \in (0,\frac{1}{10})$ of the Riemannian Laplacian 
$$\Delta^{\mathrm{Rm}}_{\mathrm{KE}} = {}\star_{\mathrm{KE}} d\star_{\mathrm{KE}} d = 2\pi \Delta_{\mathrm{KE}} .$$

The following lemma provides a local estimate for an eigenfunction on the universal covering, which in turn yields an estimate for $\Phi$ on the thick part.

\begin{lem}
\label{lemlocalsupnormthickpart}
Let $(\mathbb{D},\mu_{\mathbb{D}})$ denote the Poincar\'e disk (see Theorem \ref{thmpoincaremodel}), and $\Phi $ be a $\mathbb{R}$-valued function on $\mathbb{D}$ satisfying $\Delta^{\mathrm{Rm}}_{\mathrm{KE}} \Phi = -\epsilon \Phi $ for some $\epsilon \in (0,1)$. Then for any $r\in (0,\arcsinh 1]$, we have
$$ |\Phi (0)|^2 \leq \frac{1}{2\pi (\cosh r -1)} \cdot \frac{1}{ \left( 1-\epsilon \log \left( \frac{1+\cosh r}{2} \right) \right)^2} \int_{B_r (0)} |\Phi (z)|^2 d\mu_{\mathbb{D}} , $$
where $ B_r(0) = \left\{ z\in \mathbb{D} :\mathrm{dist}_{\mu_{\mathbb{D}}} (z,0) <r \right\} = \left\{ z\in \mathbb{D} :|z| <\tanh \frac{r}{2} \right\} .$
\end{lem}

\begin{proof}
Let $\Phi_0 (z) = \frac{1}{2\pi} \int_{0}^{2\pi} \Phi (ze^{i\theta}) d\theta $ denote the circular average of $\Phi$ around the origin. Then $\Phi_0 (0) = \Phi (0) $, $d\Phi_0 (0) =0$, $\Delta^{\mathrm{Rm}}_{\mathrm{KE}} \Phi_0 = -\epsilon \Phi_0 $, and $\Phi_0 (z) $ depends only on $|z|$, that is, $\Phi_0$ is radial. Write $\mathbf{r} (z) = \mathrm{dist}_{\mu_{\mathbb{D}}} (z,0) = 2 \mathrm{arctanh} (|z|) $. Then on the polar coordinates $(\mathbf{r},\vartheta)$, Theorem \ref{thmanalysispartspaceform} yields $\mu_{\mathbb{D}} = \sinh \mathbf{r} d\mathbf{r} \wedge d\vartheta $. Set $\Phi_1 (\mathbf{r})= \Phi_0 (\tanh \frac{\mathbf{r}}{2})$. Since $\Phi_0$ depends only on $\mathbf{r}$, the Laplacian reduces to the radial operator
$$ \Delta^{\mathrm{Rm}}_{\mathrm{KE}} \Phi_0 = \frac{\partial^2 \Phi_1}{\partial \mathbf{r}^2} + \frac{\cosh \mathbf{r}}{\sinh \mathbf{r}} \cdot \frac{\partial \Phi_1}{\partial \mathbf{r} } = -\epsilon \Phi_0 .$$
Without loss of generality, assume that $\Phi_0 (0) = \Phi (0) > 0 $. Then we have $\frac{\partial^2 \Phi_1}{\partial \mathbf{r}^2} (0) = -\epsilon \Phi_0 (0) <0 $, so there exists an open neighborhood $U$ of $0$ such that $\Phi_0 (z) < \Phi_0 (0)$, $\forall z\in U\setminus\{0\}$. For any $\mathbf{r}>0$, a straightforward calculation shows that 
$$\frac{\partial }{\partial \mathbf{r} } \left( \sinh \mathbf{r} \cdot \frac{\partial \Phi_1}{\partial \mathbf{r} } (\mathbf{r}) \right)  = -\epsilon \sinh \mathbf{r} \cdot \Phi_1 (\mathbf{r}) .$$ 
Hence, if $\Phi_1 (\mathbf{r} ) \geq 0 $ for $\mathbf{r} \in (0,\mathbf{r}_0)$, then $\frac{\partial \Phi_1}{\partial \mathbf{r} } \bigg|_{\mathbf{r}=\mathbf{r}_0} < 0 $, and consequently $\Phi_1 (\mathbf{r}_0) < \Phi_1 (0)$.

By the Cauchy-Schwarz inequality, for any $r>0$, we have
\begin{eqnarray*}
\int_{B_r (0)} |\Phi (z)|^2 d\mu_{\mathbb{D}} & = & \int_0^r \sinh \mathbf{r} \int_{0}^{2\pi} |\Phi (\tanh (\mathbf{r}) e^{i\vartheta} )|^2 d\vartheta d\mathbf{r} \\
& \geq & \int_0^r \sinh \mathbf{r} \int_{0}^{2\pi} |\Phi_0 (\tanh \mathbf{r} )|^2 d\vartheta d\mathbf{r} = \int_{B_r (0)} |\Phi_0 (z)|^2 d\mu_{\mathbb{D}} .
\end{eqnarray*}
It remains to prove that
$$ |\Phi_0 (0)|^2 \leq \frac{1}{2\pi (\cosh r -1)} \cdot \frac{1}{ \left( 1-\epsilon \log \left( \frac{1+\cosh r}{2} \right) \right)^2} \int_{B_r (0)} |\Phi_0 (z)|^2 d\mu_{\mathbb{D}} . $$

Set $u (t) = \log \left( \frac{1+\cosh t}{2} \right)$. Then $u(0)=u'(0)=0$, and a direct computation shows
\begin{equation*}
u'' (t) + \frac{\cosh t}{\sinh t} \cdot u' (t) = \frac{1}{1+\cosh t} + \frac{\cosh t}{1+\cosh t} =1 .
\end{equation*}

Write $\Phi_u (z) = (1-\epsilon u(\mathbf{r} (z))) \Phi_0 (0) $ and $\Phi_{u,1} (\mathbf{r})=\Phi_u (\tanh\frac{\mathbf{r}}{2})$. Then we have $\Phi_u (0) = \Phi_0 (0) $, $d\Phi_u (0) = d\Phi_0 (0) =0 $, and
$$ \Delta^{\mathrm{Rm}}_{\mathrm{KE}} \Phi_u = \frac{\partial^2 \Phi_{u,1}}{\partial \mathbf{r}^2} + \frac{\cosh \mathbf{r}}{\sinh \mathbf{r}} \cdot \frac{\partial \Phi_{u,1}}{\partial \mathbf{r} } = -\epsilon \Phi_0 (0) .$$

We now show that $\Phi_0 (0) \geq\Phi_0 (z) \geq\Phi_u (z) $, $\forall z\in B_r (0)$. Since $\epsilon \in (0,1)$ and $r\in (0,\arcsinh 1]$, we obtain for any $z\in B_r (0)$, 
$$\Phi_0 (0)\geq\Phi_u (z) = \Phi_0 (0) \cdot (1-\epsilon u(\mathbf{r} (z))) > \Phi_0 (0) \cdot (1- \log 2) >0 .$$

For any $\mathbf{r}>0$, a straightforward calculation gives
$$\frac{\partial }{\partial \mathbf{r} } \left( \sinh \mathbf{r} \cdot \frac{\partial \Phi_1}{\partial \mathbf{r} } - \sinh \mathbf{r} \cdot \frac{\partial \Phi_{u,1}}{\partial \mathbf{r} } \right) = \epsilon \sinh \mathbf{r} \cdot \left( \Phi_0 (0) - \Phi_0 \left( \tanh \frac{\mathbf{r}}{2} \right) \right) .$$
Thus, if $\Phi_0 (0 ) \geq \Phi_0 \left( \tanh \frac{\mathbf{r}}{2} \right) $ for $\mathbf{r} \in (0,\mathbf{r}_0)$, then $\frac{\partial (\Phi_0 - \Phi_u )}{\partial \mathbf{r} } (\mathbf{r}_0) > 0 $, which implies $\Phi_u \left( \tanh\frac{\mathbf{r}_0}{2} \right) < \Phi_0 \left( \tanh \frac{\mathbf{r}_0}{2} \right) $.

Let $\mathbf{r}_{\mathrm{max}} $ be the largest radius in $[0,\infty)$ such that $\Phi_0 (0) \geq\Phi_0 (z) \geq\Phi_u (z) $, $\forall z\in B_{\mathbf{r}_{\mathrm{max}}} (0)$. Since there exists an open neighborhood $U$ of $0$ such that $\Phi_0 (z) < \Phi_0 (0)$, $\forall z\in U\setminus\{0\}$, one can find $\mathbf{r}_{1}\in (0,r)$ such that $\Phi_u \left( \tanh \frac{\mathbf{r}}{2} \right) < \Phi_0 \left( \tanh \frac{\mathbf{r}}{2} \right) $, $\forall \mathbf{r}\in (0,\mathbf{r}_{1} )$. Hence $\mathbf{r}_{\mathrm{max}}\geq \mathbf{r}_1 >0$. Since $\Phi_0 >0$ on $ B_r (0)$, if $\mathbf{r}_{\mathrm{max}}< r$, then 
$$\frac{\partial (\Phi_1 - \Phi_{u,1} )}{\partial \mathbf{r} } \bigg|_{\mathbf{r}=\mathbf{r}_{\mathrm{max}}} > 0 > \frac{\partial \Phi_1}{\partial \mathbf{r} } \bigg|_{\mathbf{r}=\mathbf{r}_{\mathrm{max}}} $$ 
implies that $\Phi_0 (0) \geq\Phi_0 (z) \geq\Phi_u (z) $ on an open neighborhood of $\bar{B}_{\mathbf{r}_{\mathrm{max}}} (0)$, contradicting the definition of $\mathbf{r}_{\mathrm{max}}$. Hence $\mathbf{r}_{\mathrm{max}}\geq r$, and $\Phi_0 (0) \geq\Phi_0 (z) \geq\Phi_u (z) $, $\forall z\in B_r (0)$.

Finally, combining the estimates above gives
\begin{eqnarray*}
\int_{B_r (0)} |\Phi (z)|^2 d\mu_{\mathbb{D}} & \geq & \int_{B_r (0)} |\Phi_0 (z)|^2 d\mu_{\mathbb{D}} \geq \int_{B_r (0)} |\Phi_u (z)|^2 d\mu_{\mathbb{D}} \\
& = & \int_0^r 2\pi \sinh \varsigma \cdot |\Phi_0 (0)|^2 \cdot (1-\epsilon u(\varsigma))^2 d\varsigma \\
& \geq & 2\pi |\Phi_0 (0)|^2 \cdot \int_0^r \sinh \varsigma \cdot \left( 1-\epsilon \log \left( \frac{1+\cosh r}{2} \right) \right)^2 d\varsigma \\
& = & 2\pi \left( \cosh r - 1 \right) \cdot \left( 1-\epsilon \log \left( \frac{1+\cosh r}{2} \right) \right)^2   |\Phi (0)|^2 ,
\end{eqnarray*}
which completes the proof.
\end{proof}

Before estimating the sup-norm of local eigenfunctions on collars, we recall the following fundamental comparison principle for ordinary differential equations.

\begin{lem}
\label{lemODEcomparison}
Let $T\in (0,\infty)$, let $a_1 (t)\geq a_2 (t) \geq 0$ be smooth functions on $[0,T)$, and let $u_1 (t)$, $u_2 (t)$ be smooth functions on $[0,T )$ satisfying the following equations:
\begin{equation*}
\left\{
\begin{aligned}
u_j'' (t) = a_j (t) u_j (t), \; & \;\; \forall t\in [0,T ) ,\;\forall j\in\{ 1,2 \} ;\\
 u_1 (0) \geq u_2 (0)\geq 0, \; & \;\; u'_1 (0) \geq u'_2 (0) \geq 0  .
\end{aligned}
\right.
\end{equation*}
Then we have $u_1 (t) \geq u_2 (t)$, $u'_1 (t) \geq u'_2 (t)$, $\forall t\in [0,T )$.
\end{lem}

\begin{proof}
Although this result is well known to experts, we include a proof here for completeness. The argument essentially follows the standard proof of existence and uniqueness for solutions of ordinary differential equations.

For $j\in\{1,2 \}$, set $u_{j,0,1} (t) = u_j (0) $, and $u_{j,0,2} (t) = u'_j (0) $. We then define two sequences of functions $\{u_{j,k,1} \}_{k=1}^\infty $ and $\{u_{j,k,2} \}_{k=1}^\infty $ inductively by:
\begin{eqnarray*}
u_{j,k,1} (t) & = & u_j (0) + \int_0^t u_{j,k-1,2} (\varkappa ) d\varkappa \; ,\\
u_{j,k,2} (t) & = & u'_j (0) + \int_0^t a_j (\varkappa) u_{j,k-1,1} (\varkappa ) d\varkappa \; .
\end{eqnarray*}
By induction on $k$, one checks that $u_{1,k,1} \geq u_{2,k,1} $ and $u_{1,k,2} \geq u_{2,k,2} $, $\forall k\in\mathbb{N}$. Since $u_j (t)= \lim\limits_{k\to\infty} u_{j,k,1} (t) $, $j\in\{1,2 \}$, we obtain $u_1 (t) \geq u_2 (t)$, $\forall t\in [0,T )$. A similar argument yields $u'_1 (t) \geq u'_2 (t)$ on $[0,T)$, completing the proof.
\end{proof}

As an application of Lemma \ref{lemODEcomparison}, we have the following estimates.

\begin{lem}
\label{lemcollareigenvalueODEcomparison}
Let $\ell \in (0,\infty )$ and $\epsilon\in (0,\frac{1}{4})$ be constants. For any $k\in\mathbb{Z}$, $j\in\{1,2\}$, let $u_{k,j} (t)$ denote the unique solution to the initial value problem:
$$u_{k,j}'' (t) = \left( \frac{1}{4} -\epsilon + \frac{1}{4\cosh^2 t} + \frac{4k^2\pi^2}{\ell^2 \cosh^2 t} \right) \cdot u_{k,j} (t) $$
with initial conditions
\begin{equation*}
\left\{
\begin{aligned}
u_{k,1} (0) =1, \; &\;  u_{k,1}' (0) =0;\\
u_{k,2} (0) =0, \; &\;  u_{k,2}' (0) =1 .
\end{aligned}
\right.
\end{equation*}

Then the following properties hold:
\begin{enumerate}[(1)]
    \item The functions $u_{k,1} (t)$ are even functions, and $u_{k,2} (t)$ are odd functions.
    \item For any $T\in (0,\infty)$, the function $\frac{u_{k,j} (t)}{\cosh (\frac{2k\pi t}{\ell \cosh T})}$ is increasing on $[0,T]$.
    \item $\cosh \left( \frac{\sqrt{1-4\epsilon}}{2} t \right) \leq u_{0,1} (t) \leq \sqrt{\cosh t} $, $\forall t\in [0,\infty )$.
    \item $ \frac{2}{\sqrt{1-4\epsilon }} \sinh \left( \frac{\sqrt{1-4\epsilon}}{2} t \right) \leq u_{0,2} (t) \leq 2\left( \arctan e^t - \frac{\pi}{4} \right) \sqrt{\cosh t } $, $\forall t\in [0,\infty )$.
\end{enumerate}
\end{lem}

\begin{proof}
The property (1) follows immediately from 
the symmetry of the differential equation together with the given initial conditions and uniqueness of solutions.

The next step is to prove (2). By definition, we have 
$$\frac{1}{4} -\epsilon + \frac{1}{4\cosh^2 t} + \frac{4k^2\pi^2}{\ell^2 \cosh^2 t} > \left( \frac{2k\pi}{\ell \cosh T} \right)^2 ,\;\; \forall t\in [0,T] .$$
Then Lemma \ref{lemODEcomparison} implies that $u_{k,j} (t) >0 $ and $u'_{k,j} (t) >0 $, $\forall t\in (0,T]$. A straightforward calculation shows that for any $t\in (0,T]$,
\begin{eqnarray*}
& & \frac{d}{dt} \left(u'_{k,j} (t) \cosh \left(\frac{2k\pi t}{\ell \cosh T}\right) - \frac{2k\pi}{\ell \cosh T} \sinh \left( \frac{2k\pi t}{\ell \cosh T} \right) u_{k,j} (t) \right) \\
& = & \cosh \left(\frac{2k\pi t}{\ell \cosh T}\right) \left(u''_{k,j} (t) - \left( \frac{2k\pi}{\ell \cosh T} \right)^2 u_{k,j} (t) \right) >0 ,
\end{eqnarray*}
and hence
\begin{eqnarray*}
\frac{d}{dt} \left(\frac{u_{k,j} (t)}{\cosh (\frac{2k\pi t}{\ell \cosh T})} \right) & = & \frac{u'_{k,j} (t) \cosh (\frac{2k\pi t}{\ell \cosh T}) - \frac{2k\pi}{\ell \cosh T} \sinh (\frac{2k\pi t}{\ell \cosh T}) u_{k,j} (t)}{\cosh^2 (\frac{2k\pi t}{\ell \cosh T})} > u'_{k,j} (0) \geq 0.
\end{eqnarray*}
This proves (2).

The properties (3) and (4) follow directly from Lemma \ref{lemODEcomparison}.
\end{proof}

Now we are ready to consider the sup-norm of local eigenfunctions on collars.

\begin{lem}
\label{lmmlocalsupnormthinpart}
Let $\gamma\subset (C,\mu_{\mathrm{KE}})$ be a short simple closed geodesic such that the length $\ell (\gamma) = \ell\in (0,\frac{\arcsinh 1}{4}]$, and let
$$ (\mathscr{C} (\gamma ), ds^2)\cong \left( [-w , w ] \times (\mathbb{R}/\mathbb{Z}) ,dr^2 + \ell^2 \cosh^2 r dt^2\right) $$ 
be the collar of $\gamma$, where $w = \arcsinh \{ 1/\sinh ( \frac{\ell}{2} ) \}$. Let $\Phi $ be a $\mathbb{R}$-valued function on $\mathscr{C} (\gamma ) $ satisfying $\Delta^{\mathrm{Rm}}_{\mathrm{KE}} \Phi = -\epsilon \Phi $ for some $\epsilon \in (0,\frac{1}{10} )$. Then for any $r\in (0,w -\log 2 )$, $t\in \mathbb{R}/\mathbb{Z}$, we have
$$ |\Phi (r,t)|^2 \leq (2\ell^{-4\epsilon } + 1) \int_{\mathscr{C} (\gamma )} |\Phi |^2 d\mu_{\mathrm{KE}}, $$
\end{lem}

\begin{proof}
Without loss of generality, we can assume that $\int_{\mathscr{C} (\gamma )} |\Phi |^2 \mu_{\mathrm{KE}} =1 $.

For any $k\in\mathbb{Z}$ and $j\in\{ 1,2 \}$, set
\begin{equation*}
\Phi_{k,1} (r ) = \int_0^1 \Phi (r,t) \cos (2k\pi t) dt\quad \textrm{ and } \quad \Phi_{k,2} (r ) = \int_0^1 \Phi (r,t) \sin (2k\pi t) dt.
\end{equation*}
Then by the Fourier expansion of $\Phi $, we obtain
\begin{eqnarray*}
\Phi (r,t) & = & \Phi_0 (r) + 2 \sum_{k=1}^\infty \Phi_{k,1 } (r ) \cos (2k\pi t) + 2 \sum_{k=1}^\infty \Phi_{k,2 } (r ) \sin (2k\pi t) .
\end{eqnarray*}
Since $\Phi $ is smooth on $C$, it follows that the Fourier expansion converges in the $C^\infty$ sense. By definition,
$$ \Delta^{\mathrm{Rm}}_{\mathrm{KE}} \Phi = \frac{\partial^2 \Phi}{\partial r^2} + \tanh r \cdot \frac{\partial \Phi}{\partial r} + \frac{1}{\ell^2 \cosh^2 r} \cdot \frac{\partial^2 \Phi}{\partial t^2} .$$
Hence for any $k\in\mathbb{Z}$ and $j\in\{ 1,2 \}$, the function $\Phi_{k,j} (r ) $ satisfies the ordinary differential equation:
$$ \Phi_{k,j}'' (r ) + \tanh r \cdot \Phi_{k,j}' (r ) = \left( \frac{4\pi^2 k^2}{\ell^2 \cosh^2 r} - \epsilon \right) \cdot \Phi_{k,j} (r ) . $$
A straightforward calculation shows that the functions $ \Phi_{k} (r) \cdot \sqrt{\cosh r} $ satisfy
$$\frac{  \frac{d^2}{dr^2} \left( \Phi_{k,j} (r ) \sqrt{\cosh r}\right) }{ \Phi_{k,j} (r ) \sqrt{\cosh r} } = \left( \frac{1}{4} - \epsilon + \frac{1}{4\cosh^2 r} + \frac{4\pi^2 k^2}{\ell^2 \cosh^2 r}  \right) . $$

For each $k\in\mathbb{Z}$ and $m\in \{1,2\} $, let $u_{k,m}$ denote the unique solution of the initial value problem:
$$u_{k,m}'' (t) = \left( \frac{1}{4} -\epsilon + \frac{1}{4\cosh^2 t} + \frac{4k^2\pi^2}{\ell^2 \cosh^2 t} \right) \cdot u_{k,m} (t) $$
with initial conditions
\begin{equation*}
\left\{
\begin{aligned}
u_{k,1} (0) =1, \; &\;  u_{k,1}' (0) =0;\\
u_{k,2} (0) =0, \; &\;  u_{k,2}' (0) =1 .
\end{aligned}
\right.
\end{equation*}

Then for any $k\in\mathbb{Z}$ and $j,m\in\{ 1,2 \}$, there exists a constant $a_{k,j,m}$ such that 
$$\Phi_{k,j } (r) \sqrt{\cosh r} = a_{k,j,1} u_{k,1} (r) + a_{k,j,2} u_{k,2} (r) .$$
Hence we obtain the following $L^2$-orthogonal decomposition:
\begin{eqnarray*}
\Phi (r,t) \sqrt{\cosh r} & = & \sum_{m=1}^2 a_{0,1,m} u_{0,m} (r) + 2 \sum_{k=1}^\infty \sum_{m=1}^2 a_{k,1,m} u_{k,m} (r) \cos (2k\pi t) \\
& & + 2 \sum_{k=1}^\infty \sum_{m=1}^2 a_{k,2,m} u_{k,m} (r) \sin (2k\pi t) .
\end{eqnarray*}
Since $\int_C |\Phi |^2 \mu_{\mathrm{KE}} =1 $, we have
\begin{eqnarray*}
1 & = & \int_{-w}^w \left( \int_0^1 \left| \Phi (r,t) \right|^2 \ell \cosh r dt \right) dr \\
& = & \ell \sum_{m=1}^2 |a_{0,1,m}|^2 \int_{-w}^w |u_{0,m} (r)|^2 dr + 2\ell \sum_{k=1}^\infty \sum_{m=1}^2 |a_{k,1,m}|^2 \int_{-w}^w |u_{k,m} (r)|^2 dr \\
& & + 2\ell \sum_{k=1}^\infty \sum_{m=1}^2 |a_{k,2,m}|^2 \int_{-w}^w |u_{k,m} (r)|^2 dr .
\end{eqnarray*}
By the Cauchy-Schwarz inequality, it follows that
\begin{eqnarray*}
|\Phi (r,t)|^2 \cosh(r) & = & \left| \sum_{m=1}^2 a_{0,1,m} u_{0,m} (r) + 2 \sum_{k=1}^\infty \sum_{m=1}^2 a_{k,1,m} u_{k,m} (r) \cos (2k\pi t) \right.\\
& & \left. + 2 \sum_{k=1}^\infty \sum_{m=1}^2 a_{k,2,m} u_{k,m} (r) \sin (2k\pi t) \right|^2 \\
& \leq & \sum_{m=1}^2 \frac{|u_{0,m} (r)|^2}{\ell \int_{-w}^w |u_{0,m} (\varsigma)|^2 d\varsigma} + \sum_{k=1}^\infty \sum_{m=1}^2 \frac{2|u_{k,m} (r)|^2 \cos^2 (2k\pi t)}{\ell \int_{-w}^w |u_{k,m} (\varsigma)|^2 d\varsigma} \\
& & + \sum_{k=1}^\infty \sum_{m=1}^2 \frac{2|u_{k,m} (r)|^2 \sin^2 (2k\pi t)}{\ell \int_{-w}^w |u_{k,m} (\varsigma)|^2 d\varsigma} \\
& = & \sum_{m=1}^2 \frac{|u_{0,m} (r)|^2}{\ell \int_{-w}^w |u_{0,m} (\varsigma)|^2 d\varsigma} + \sum_{k=1}^\infty \sum_{m=1}^2 \frac{2|u_{k,m} (r)|^2 }{\ell \int_{-w}^w |u_{k,m} (\varsigma)|^2 d\varsigma} .
\end{eqnarray*}

Now we estimate the terms in the above inequality using Lemma \ref{lemcollareigenvalueODEcomparison}.

By Lemma \ref{lemcollareigenvalueODEcomparison}, we have $ \cosh \left( \frac{\sqrt{1-4\epsilon}}{2} r \right) \leq u_{0,1} (r) \leq \sqrt{\cosh r} $, and hence
\begin{eqnarray*}
\frac{|u_{0,1} (r)|^2}{\ell \cosh r \int_{-w}^w |u_{0,1} (\varsigma)|^2 d\varsigma} & \leq & \frac{1}{ \ell \int_{-w}^w \cosh^2 \left( \frac{\sqrt{1-4\epsilon}}{2} \varsigma \right) d\varsigma} \\
& = & \frac{2}{ \ell \int_{-w}^w \left( \cosh \left( \sqrt{1-4\epsilon}\cdot \varsigma \right) +1\right) d\varsigma} \\
& \leq & \frac{\sqrt{1-4\epsilon}}{ \ell \left( \sinh \left( \sqrt{1-4\epsilon}\cdot w \right) +\sqrt{1-4\epsilon}\cdot w\right) } .
\end{eqnarray*}

Since 
$$\sqrt{1-4\epsilon}\cdot w \geq \sqrt{\frac{3}{5}} \cdot \arcsinh \left\{ \frac{1}{\sinh ( \frac{\arcsinh 1}{8} )} \right\} \approx 2.246232\cdots  > 2 ,$$ 
we conclude that $2\sqrt{1-4\epsilon}\cdot w > e^{-\sqrt{1-4\epsilon}\cdot w}  $, and 
\begin{eqnarray*}
2\sqrt{1-4\epsilon} \cdot w e^{- \sqrt{1-4\epsilon}\cdot w} + e^{-2\sqrt{1-4\epsilon}\cdot w} \leq 0.486\cdots < \frac{1}{2}.
\end{eqnarray*}
Hence
\begin{eqnarray*}
\frac{|u_{0,1} (r)|^2}{\ell \cosh r \int_{-w}^w |u_{0,1} (\varsigma)|^2 d\varsigma} & \leq & \frac{\sqrt{1-4\epsilon}}{ \ell \left( \sinh \left( \sqrt{1-4\epsilon}\cdot w \right) +\sqrt{1-4\epsilon}\cdot w\right) } \\
& < & 2\ell^{-1} e^{-\sqrt{1-4\epsilon}\cdot w} .
\end{eqnarray*}

Similarly, by Lemma \ref{lemcollareigenvalueODEcomparison}, for any $r\in(0,w]$, we have
$$ \frac{2}{\sqrt{1-4\epsilon }} \sinh \left( \frac{\sqrt{1-4\epsilon}}{2} r \right) \leq u_{0,2} (r) \leq 2\left( \arctan e^r - \frac{\pi}{4} \right) \sqrt{\cosh (r )} .$$ 
Hence
\begin{eqnarray*}
\frac{|u_{0,2} (r)|^2}{\ell \cosh (r) \int_{-w}^w |u_{0,2} (t)|^2 dt} & \leq & \frac{\left( \arctan e^{|r|} - \frac{\pi}{4} \right)^2}{ \ell  \int_{-w}^w \frac{1}{1-4\epsilon} \sinh^2 \left( \frac{\sqrt{1-4\epsilon}}{2} \varsigma \right) d\varsigma} \\
& < & \frac{\pi^2 (1-4\epsilon )}{ 16\ell \int_{-w}^w \sinh^2 \left( \frac{\sqrt{1-4\epsilon}}{2} \varsigma \right) d\varsigma} \\
& = & \frac{\pi^2 (1-4\epsilon )^{\frac{3}{2}}}{ 16\ell \left( \sinh \left( \sqrt{1-4\epsilon}\cdot w \right) -\sqrt{1-4\epsilon}\cdot w\right)} \\
& \leq & \frac{\pi^2}{4\ell } \cdot e^{-\sqrt{1-4\epsilon}\cdot w} .
\end{eqnarray*}

It follows that
\begin{eqnarray*}
\sum_{m=1}^2 \frac{|u_{0,m} (r)|^2}{\ell \int_{-w}^w |u_{0,m} (\varsigma)|^2 d\varsigma}
& \leq & \left( 2+\frac{\pi^2}{4 } \right) \cdot \ell^{-1} e^{-\sqrt{1-4\epsilon}\cdot w} \\
& = & \left( 2+\frac{\pi^2}{4 } \right) \cdot \ell^{-1} \cdot \left( \tanh \left( \frac{\ell}{4} \right)\right)^{\sqrt{1-4\epsilon}} \\
& < & \left( 2+\frac{\pi^2}{4 } \right) \cdot 2^{-2\sqrt{\frac{3}{5}}} \cdot \ell^{-1+\sqrt{1-4\epsilon}} \\
& < & 2\ell^{-4\epsilon } .
\end{eqnarray*}

If $k\geq 1$, Lemma \ref{lemcollareigenvalueODEcomparison} implies that $\frac{u_{k,j} (r)}{\cosh (\frac{2k\pi r}{\ell \cosh w})}$ is increasing on $[0,w]$. It follows that
\begin{eqnarray*}
\int_{-w}^w |u_{k,m} (\varsigma)|^2 d\varsigma & \geq & 2 \int_{w-\log 2}^w |u_{k,m} (\varsigma)|^2 d\varsigma \\
& > & 2 \int_{w-\log 2}^w \frac{ |u_{k,m} (w-\log 2)|^2 \cosh^2 (3k\varsigma)}{\cosh^2 (3k (w-\log 2) )} d\varsigma \\
& \geq & 2 |u_{k,m} (w-\log 2)|^2 \int_{w-\log 2}^w \frac{e^{6k(\varsigma-w+\log 2)}}{(1+e^{-6(w-\log 2)})^2} d\varsigma \\
& \geq & \frac{2 |u_{k,m} (w-\log 2)|^2}{(1+e^{-6 })^2 } \cdot \frac{1}{6k} \cdot \left( 2^{6k} -1 \right) \\
& \geq & \frac{(1-2^{-6 })^3}{3k} \cdot 2^{6k} \cdot |u_{k,m} (w-\log 2)|^2 .
\end{eqnarray*}

Hence for any $r\in [-w+\log2,w-\log 2]$, 
\begin{eqnarray*}
\sum_{k=1}^\infty \sum_{m=1}^2 \frac{2|u_{k,m} (r)|^2 }{\ell\cosh r \int_{-w}^w |u_{k,m} (\varsigma)|^2 d\varsigma} & \leq & \sum_{k=1}^\infty \sum_{m=1}^2 \frac{2|u_{k,m} (w-\log 2)|^2 }{\ell \cosh (w-\log 2)} \cdot \frac{1}{ \int_{-w}^w |u_{k,m} (\varsigma)|^2 d\varsigma} \\
& \leq & \sum_{k=1}^\infty \sum_{m=1}^2 \frac{2 }{2^{-1} \ell \cosh w} \cdot \frac{3k \cdot 2^{-6k}}{ (1-2^{-6 })^3} \\
& = & \sum_{k=1}^\infty \sum_{m=1}^2 \frac{2\tanh \frac{\ell}{2} }{ \ell } \cdot \frac{6k \cdot 2^{-6k}}{ (1-2^{-6 })^3} \\
& \leq & \sum_{k=1}^\infty 24 k \cdot 2^{-6k} = \frac{3}{8} \cdot \frac{1}{(1-\frac{1}{64})^2} <\frac{1}{2} .
\end{eqnarray*}
Combining the estimates above, we conclude that
\begin{eqnarray*}
|\Phi (r,t)|^2 & \leq & \sum_{m=1}^2 \frac{|u_{0,m} (r)|^2 (\cosh r )^{-1}}{\ell \int_{-w}^w |u_{0,m} (\varsigma)|^2 d\varsigma} + \sum_{k=1}^\infty \sum_{m=1}^2 \frac{2|u_{k,m} (r)|^2 (\cosh r )^{-1} }{\ell \int_{-w}^w |u_{k,m} (\varsigma)|^2 d\varsigma} \\
& < & 2\ell^{-4\epsilon } + 1 .
\end{eqnarray*}
This completes the proof.
\end{proof}

Combining the estimates on the thick part and the collars, we obtain the following global estimate.

\begin{prop}
\label{propglobalsupnormeigenfunctionsmalleigenvalue}
Let $\gamma\subset (C,\mu_{\mathrm{KE}})$ be a shortest simple closed geodesic, and let $\Phi $ be a $\mathbb{R}$-valued function on $\mathscr{C} (\gamma ) $ satisfying $\Delta_{\mathrm{KE}} \Phi = -\epsilon \Phi $ for some $\epsilon \in (0,\frac{1}{20\pi} )$. Then we have
$$ \sup_{C}|\Phi |^2 \leq \max\left\{ 2\ell^{-8\pi \epsilon } +1, \frac{16}{5} \right\} \cdot \int_{C} |\Phi |^2 d\mu_{\mathrm{KE}}, $$
where $\ell=\ell(\gamma)$ denotes the length of $\gamma$.
\end{prop}

\begin{proof}
Let $\gamma_0 ,\gamma_1 ,\cdots ,\gamma_m$ be the simple closed geodesics of length $\ell (\gamma_j)\leq \frac{\arcsinh 1}{4}$, and let
$$ (\mathscr{C} (\gamma_j ), ds^2)\cong \left( [-w_j , w_j ] \times (\mathbb{R}/\mathbb{Z}) ,dr_j^2 + \ell(\gamma_j)^2 \cosh^2 r_j dt_j^2\right) $$ 
be the collar of $\gamma_j$, where $w_j = \arcsinh \{ 1/\sinh ( \ell (\gamma_j) /2 ) \}$.

Set $U=\bigcup\limits_{j=1}^{m} \left\{ x\in \mathcal{C} (\gamma_j) :|r_j (x)|\leq w-\log 2 \right\}$. By Lemma \ref{lmmlocalsupnormthinpart}, we have
$$\sup_{ U}|\Phi |^2 \leq (2\ell^{-8\pi\epsilon } + 1) \int_{C} |\Phi |^2 d\mu_{\mathrm{KE}} .$$
Note that $\Delta^{\mathrm{Rm}}_{\mathrm{KE}}=2\pi \Delta_{\mathrm{KE}}$.

Now consider the sup-norm of $\Phi$ on $C\setminus U$. For any $x\in C\setminus U$, take the covering map $\rho_{\mathbb{D}} : \mathbb{D}\to C$ from the Poincar\'e disk such that $\rho_{\mathbb{D}} (0)=x $. By Lemma \ref{lemlocalstructureorbitcollartheorem}, we have 
$$ \# \left(B_{\arcsinh 1} (0) \cap \rho_{\mathbb{D}}^{-1} (x') \right) \leq \max \left\{ \frac{2\arcsinh 1}{\frac{1}{4} \arcsinh 1} ,\, 2 \right\} =8,\quad \forall x'\in C ,$$
where $B_{r} (y)$ is the metric ball in $(\mathbb{D},\mu_{\mathbb{D}})$ of radius $r$ and center at $y$. Lemma \ref{lemlocalsupnormthickpart} then yields
\begin{eqnarray*}
\sup_{C\setminus U}|\Phi |^2 & \leq & 8 \cdot \frac{1}{2\pi (\sqrt{2} -1)} \cdot \frac{1}{ \left( 1-\frac{1}{10} \log \left( \frac{1+\sqrt{2}}{2} \right) \right)^2} \cdot \int_{C} |\Phi |^2 d\mu_{\mathrm{KE}} \\
& = & \frac{4}{ \pi (\sqrt{2} -1)} \cdot \frac{1}{ \left( 1-\frac{1}{10} \log \left( \frac{1+\sqrt{2}}{2} \right) \right)^2} \cdot \int_{C} |\Phi |^2 d\mu_{\mathrm{KE}} \\
& \leq & \frac{16}{5} \int_{C} |\Phi |^2 d\mu_{\mathrm{KE}} .
\end{eqnarray*}

Hence
$$ \sup_{C}|\Phi |^2 \leq \max\left\{ 2\ell^{-8\pi\epsilon } +1, \frac{16}{5} \right\} \int_{C} |\Phi |^2 \mu_{\mathrm{KE}} , $$
which proves this proposition.
\end{proof}

\section{Lower bounds of Zhang’s \texorpdfstring{$\varphi$}{Lg}-invariant}

In this section, we derive an explicit estimate for Zhang’s $\varphi$-invariant for curves over $\mathbb{C}$. Let $C$ be a curve of genus $g\geq 2$ over $\mathbb{C}$, and let $\omega_{C }$ denote its canonical bundle. Let $\mu_{\mathrm{KE}} $ be the unique K\"ahler metric on $C $ with constant curvature $-1 $. Let $\{ \alpha_j \}_{j=1}^g $ be an $L^2$-orthonormal basis of $\Gamma (C,\omega_C)$. The Arakelov K\"ahler metric is defined by 
$$\mu_{\mathrm{Ar}} = \frac{i}{2g} \sum_{k=1}^g \alpha_k \wedge \bar{\alpha}_k ,$$
and is independent of the choice of basis. The $\varphi$-invariant is defined as
$$\varphi(C)=-\int_{C^2}G_{\mathrm{Ar}}\, c_1( \Delta, G_{\mathrm{Ar}})^{\wedge 2} ,$$
where $G_{\mathrm{Ar}}:C^2\setminus \Delta \longrightarrow \mathbb{R}$ denotes the Arakelov Green function and $\Delta \subset C^2 $ is the diagonal.

\subsection{Equivalent formulations  of the invariant}
\label{sectionzhangphi1}

Before starting to estimate Zhang's $\varphi$-invariant, we first reformulate it into a form more suitable for analysis.

We begin by recalling the spectral decomposition of the Laplacian $\Delta_{\mathrm{Ar}}$ associated with the Arakelov K\"ahler metric $\mu_{\mathrm{Ar}}$. The eigenvalues of $\Delta_{\mathrm{Ar}}$, counted with multiplicities, form an increasing sequence
$$ 0=\lambda_{\mathrm{Ar} ,0 } < \lambda_{\mathrm{Ar} ,1 } \leq \lambda_{\mathrm{Ar} ,2 } \leq \cdots .$$
Corresponding to these eigenvalues, there exists an $L^2$-orthonormal basis $\{ \phi_{\mathrm{Ar} , l } \}_{l=0}^\infty $ of $L^2 (C,\mu_{\mathrm{Ar}} ;\mathbb{R} )$, the $L^2$-completion of $\mathscr{A}^{0}_{C,\mathbb{R}} (C) $, such that 
$$\Delta_{\mathrm{Ar}} \phi_{\mathrm{Ar} , l } + \lambda_{\mathrm{Ar} ,l } \phi_{\mathrm{Ar} , l } =0 .$$
The following formula is due to Zhang \cite{Zhang_phi}.

\begin{prop}[{\cite[Proposition 2.5.3]{Zhang_phi}}]
\label{propzhangphiinvariantcalculation}
With notation as above, the $\varphi$-invariant is given by
$$ \varphi (C) = \sum_{l=1}^\infty \sum_{j=1}^g \sum_{k=1}^g \frac{1}{2\lambda_{\mathrm{Ar} ,l} } \left| \int_C \phi_{\mathrm{Ar} ,l} \cdot \alpha_j \wedge \bar{\alpha}_k \right|^2 $$
\end{prop}

\begin{proof}
This follows from the $L^2$-expansion of $G_{\mathrm{Ar}}$ together with the identity
$$c_1( \Delta, G_{\mathrm{Ar}})=dd^c G_{\mathrm{Ar}} + \delta_{\Delta} = p_1^*\mu_{\mathrm{Ar}} + p_2^*\mu_{\mathrm{Ar}} - \frac{i}{2} \sum_{j=1}^g ( p_1^* \alpha_j \wedge p_2^* \bar{\alpha}_j + p_2^* \alpha_j \wedge p_1^* \bar{\alpha}_j ) ,$$ 
where $p_k$ denotes the projection $C\times C\to C$ onto the $k$-th component. Note that $\frac{i}{2} \int_C \alpha_j \wedge \bar{\alpha}_k = \delta_{jk} $. For more details, see \cite[Proposition 2.5.3]{Zhang_phi}.
\end{proof}

We now reformulate the above expression in preparation for estimating the $\varphi$-invariant.

By the classical theory of elliptic differential equations, for each pair of indices $j,k$, there exists a unique function $\mathfrak{u}_{j,k} \in  \mathscr{A}_{C,\mathbb{C}}^0 (C) $ satisfying $\int_C \mathfrak{u}_{j,k} \mu_{\mathrm{Ar}} =0$ and 
$$dd^c \mathfrak{u}_{j,k} = i\alpha_j \wedge\bar{\alpha}_k - \left( \int_{C} i\alpha_j \wedge\bar{\alpha}_k \right) \mu_{\mathrm{Ar}} .$$
Moreover, when $j=k$, $\mathfrak{u}_{j,k} $ is $\mathbb{R}$-valued, i.e., $\mathfrak{u}_{j,j} \in  \mathscr{A}_{C,\mathbb{R}}^0 (C) $. This also follows from the $\partial\bar{\partial}$-lemma \cite[Lemma 2.1]{tgbook1}.

These functions $\mathfrak{u}_{j,k}$ provide an alternative expression for the $\varphi$-invariant.

\begin{prop}
\label{propzhangphiinvariantpotential}
Let $C,\mu_{\mathrm{Ar}} $ and $\mathfrak{u}_{j,k} $ be as above. Then the $\varphi$-invariant of $C$ satisfies
\begin{equation*}
    \varphi (C) = \frac{1}{2}\sum_{j,k=1}^g \int_{C} d\mathfrak{u}_{j,k} \wedge d^c \bar{\mathfrak{u}}_{j,k} = \frac{1}{2\pi } \sum_{j,k=1}^g \left\Vert d {\mathfrak{u}}_{j,k} \right\Vert^2_{L^2 ;C} \; ,
\end{equation*}
where $d^c = \frac{1}{2\pi i} (\partial -\bar{\partial}) $, and $\left\Vert \cdot \right\Vert_{L^2 ;C}$ is the $L^2$-norm defined on the space $ \mathscr{A}^1_{C,\mathbb{C}}(C)$ of smooth $1$-forms on $C$ with coefficients in $\mathbb{C}$.
\end{prop}

\begin{proof}
Let $\{\lambda_{\mathrm{Ar} ,l}\}_{j=0}^\infty$, $\{\phi_{\mathrm{Ar} ,l}\}_{j=0}^\infty$ be the eigenvalues and eigenfunctions of $\Delta_{\mathrm{Ar}} $ with respect to $\mu_{\mathrm{Ar}}$. Then for any $j,k$, we can expand
$$ \mathfrak{u}_{j,k} = \sum_{l=0}^{\infty} a_{j,k,l} \phi_{\mathrm{Ar} ,l} ,$$
where $a_{j,k,l} = \int_C  \mathfrak{u}_{j,k} \phi_{\mathrm{Ar} ,l} \mu_{\mathrm{Ar}} \in \mathbb{C} $. Since $\int_C \mathfrak{u}_{j,k}\,\mu_{\mathrm{Ar}}=0$, we have $a_{j,k,0}=0$. By definition of $\mathfrak{u}_{j,k}$, for each $l\geq 1$ we obtain
$$\int_C \phi_{\mathrm{Ar} ,l} \alpha_j \wedge \bar{\alpha}_k =\int_{C} \phi_{\mathrm{Ar} ,l} dd^c \mathfrak{u}_{j,k} = \int_C \mathfrak{u}_{j,k} dd^c \phi_{\mathrm{Ar} ,l} = -\lambda_{\mathrm{Ar} ,l} a_{j,k,l} .$$
Since $\int_{C} \phi_{\mathrm{Ar} ,l } \mu_{\mathrm{Ar}} =0 $ for all $ l\geq 1$, we obtain
\begin{equation*}
\varphi (C) = \sum_{l=1}^\infty \sum_{j=1}^g \sum_{k=1}^g \frac{1}{2\lambda_{\mathrm{Ar} ,l} } \left| \int_C \phi_{\mathrm{Ar} ,l} \alpha_j \wedge \bar{\alpha}_k \right|^2 = \sum_{l=1}^\infty \sum_{j=1}^g \sum_{k=1}^g \frac{\lambda_{\mathrm{Ar} ,l} }{2} \left| a_{j,k,l} \right|^2 .
\end{equation*}

On the other hand,
\begin{equation*}
    \sum_{l=1}^\infty \lambda_{\mathrm{Ar} ,l} \left| a_{j,k,l} \right|^2 = -\int_{C} \mathfrak{u}_{j,k} dd^c \bar{\mathfrak{u}}_{j,k} = \int_{C} d \mathfrak{u}_{j,k} \wedge d^c \bar{\mathfrak{u}}_{j,k} .
\end{equation*}
Hence
\begin{equation*}
\varphi (C) = \frac{1}{2} \sum_{j,k=1}^g \int_{C} d\mathfrak{u}_{j,k} \wedge d^c \bar{\mathfrak{u}}_{j,k} = \frac{1}{2\pi } \sum_{j,k=1}^g \left\Vert d {\mathfrak{u}}_{j,k} \right\Vert^2_{L^2 ;C} ,
\end{equation*}
where the last equality follows from $\left\Vert d {\mathfrak{u}}_{j,k} \right\Vert^2_{L^2 ;C} = \frac{1}{2}\int_C d {\mathfrak{u}}_{j,k} \wedge \star d \mathfrak{u}_{j,k} $ and
$$ d^c \bar{\mathfrak{u}}_{j,k} = \frac{1}{2\pi i} (\partial \bar{\mathfrak{u}}_{j,k} -\bar{\partial} \bar{\mathfrak{u}}_{j,k}) = \frac{1}{2\pi} \star (\partial \mathfrak{u}_{j,k} + \bar{\partial} \mathfrak{u}_{j,k} ) = \frac{1}{2\pi} \star  d \mathfrak{u}_{j,k} . $$
This completes the proof.
\end{proof}

\subsection{On the thick part}

In this subsection, we derive an explicit estimate for Zhang’s $\varphi$-invariant for curves over $\mathbb{C}$ on the thick part of the moduli space. 
By the above preparation, we proceed to estimate the $\varphi$-invariant by considering the pull-backs of holomorphic $1$-forms on $C$ to the unit disk $\mathbb{D}$. For this purpose, we require the following auxiliary estimate on small disks.

\begin{lem}
\label{lemintegralgradientintegralsmalldisc}
Let $\mathbf{f}$ be a holomorphic function on $\mathbb{D}_{\sqrt{2} -1} (0) $ with $|\mathbf{f}(0)|=1$, and let $\mathbf{a} \in [0,2\sqrt{2}-2] $ be a non-negative constant. Suppose $\mathbf{u}$ is a smooth $\mathbb{R}$-valued function on $\mathbb{D}_{\sqrt{2} -1} (0) $ satisfying
$$dd^c \mathbf{u} \geq \left( |\mathbf{f} (z)|^2 - \frac{\mathbf{a}}{(1-|z|^2)^2} \right) \mu_{\mathrm{Euc}} ,$$
where $d^c = \frac{1}{2\pi i} (\partial -\bar{\partial}) $, and $\mu_{\mathrm{Euc}} = \frac{idz\wedge d\bar{z}}{2} $ is the standard Euclidean metric. Then
\begin{eqnarray*}
\int_{\mathbb{D}_{\sqrt{2} -1} (0)} d\mathbf{u} \wedge \star d \mathbf{u} & \geq & 2\pi^3 \int_0^{\sqrt{2} -1} t^3 \left(  1 - \frac{\mathbf{a} }{1-t^2} \right)^2 dt = 2\pi^3 \mathbf{I}_1 - 4\pi^3 \mathbf{a}\mathbf{I}_2 +2\pi^3 \mathbf{a}^2 \mathbf{I}_3  , 
\end{eqnarray*}
where
\begin{eqnarray*}
\mathbf{I}_1 & = & \int_0^{\sqrt{2} -1} t^3 dt = \frac{(\sqrt{2} -1 )^4}{4} \approx 0.0073593\cdots , \\
\mathbf{I}_2 & = & \int_0^{\sqrt{2} -1} \frac{t^3 }{1-t^2} dt = \frac{1}{2}\left( 2\sqrt{2} -3 + \log \left( \frac{\sqrt{2}+1}{2} \right) \right) \approx 0.0083267\cdots ,\\
\mathbf{I}_3 & = & \int_0^{\sqrt{2} -1} \frac{t^3 }{(1-t^2)^2} dt = \frac{1}{4}\left( \sqrt{2} - 1 - 2 \log \left( \frac{\sqrt{2}+1}{2} \right) \right) \approx 0.0094402 \cdots .
\end{eqnarray*}
\end{lem}

\begin{proof}
Let $\mathbf{u}_0 (z) = \frac{1}{2\pi} \int_{0}^{2\pi} \mathbf{u}(ze^{i\theta_1}) d\theta_1 $ denote the circular average of $\mathbf{u}$ around the origin. Then $\mathbf{u}_0$ is a smooth, $\mathbb{R}$-valued, radially symmetric function on $\mathbb{D}_{\sqrt{2} -1}(0)$. Since $\mu_{\mathrm{Euc}}$ is invariant under the rotations $z\mapsto ze^{i\theta_1}$, $\forall \theta_1\in \mathbb{R}$, we obtain 
$$dd^c \mathbf{u}_0 \geq \left( \frac{1}{2\pi} \int_{0}^{2\pi} |\mathbf{f} (ze^{i\theta_1})|^2 d\theta_1 - \frac{\mathbf{a}}{(1-|z|^2)^2} \right) \mu_{\mathrm{Euc}} .$$
By the Cauchy-Schwarz inequality and the holomorphicity of $\mathbf{f}$,
$$\frac{1}{2\pi} \int_{0}^{2\pi} |\mathbf{f} (ze^{i\theta_1})|^2 d\theta_1 \geq \left( \frac{1}{2\pi} \int_{0}^{2\pi} \mathbf{f} (ze^{i\theta_1})  d\theta_1 \right)^2 = \mathbf{f} (0) =1 .$$
Hence
$$dd^c \mathbf{u}_0 \geq \left( 1 - \frac{\mathbf{a}}{(1-|z|^2)^2} \right) \mu_{\mathrm{Euc}} .$$

Now switch to polar coordinates $z=\varsigma e^{i\theta}$ with $\varsigma=|z|$. Then $\star d\varsigma = \varsigma d\theta $, and $\star d\theta = -\frac{1}{\varsigma}d\varsigma $. Therefore,
\begin{eqnarray*}
 \int_{\mathbb{D}_{\sqrt{2} -1} (0)} d\mathbf{u} \wedge \star d \mathbf{u} & = & \int_{\mathbb{D}_{\sqrt{2} -1} (0)} \left( \varsigma \left|\frac{\partial \mathbf{u}}{\partial \varsigma} \right|^2 + \frac{1}{\varsigma}\left|\frac{\partial \mathbf{u}}{\partial \theta} \right|^2 \right) d\varsigma \wedge d\theta \\
 & \geq & \int_{0}^{\sqrt{2}-1} \left(\varsigma \int_{0}^{2\pi}  \left|\frac{\partial \mathbf{u}}{\partial \varsigma} (\varsigma e^{i\theta }) \right|^2 d\theta \right) d\varsigma \\
 & \geq & \int_{0}^{\sqrt{2}-1} \left(\frac{\varsigma}{2\pi} \left| \int_{0}^{2\pi} \frac{\partial \mathbf{u}}{\partial \varsigma} (\varsigma e^{i\theta }) d\theta \right|^2 \right) d\varsigma \\
 & = & \int_{\mathbb{D}_{\sqrt{2} -1} (0)} \left( \varsigma \left|\frac{1}{2\pi} \int_{0}^{2\pi} \frac{\partial \mathbf{u}}{\partial \varsigma} (\varsigma e^{i\theta +i\theta_1}) d\theta_1 \right|^2 \right) d\varsigma \wedge d\theta \\
 & = & \int_{\mathbb{D}_{\sqrt{2} -1} (0)} \varsigma \left|\frac{\partial \mathbf{u}_0 }{\partial \varsigma} \right|^2 d\varsigma \wedge d\theta =\int_{\mathbb{D}_{\sqrt{2} -1} (0)} d\mathbf{u}_0 \wedge \star d \mathbf{u}_0 .
\end{eqnarray*}

Since $\mathbf{u}_0 $ is radial, we have $ \lim\limits_{\varsigma \to 0^+}\frac{\partial \mathbf{u}_0 }{\partial \varsigma} (\varsigma ) =0 $. A direct computation gives
\begin{eqnarray*}
dd^c \mathbf{u}_0 (\varsigma e^{i\theta} ) & = & \frac{1}{2\pi \varsigma } \frac{\partial}{\partial \varsigma} \left( \varsigma \frac{\partial \mathbf{u}_0 }{\partial \varsigma } \right) (\varsigma e^{i\theta}) \mu_{\mathrm{Euc}} ,\;\; \forall \varsigma \in (0,\sqrt{2} -1) ,\; \theta\in [0,2\pi ).
\end{eqnarray*}
Hence 
$$ \frac{\partial}{\partial \varsigma} \left( \varsigma \frac{\partial \mathbf{u}_0 }{\partial \varsigma } \right) (\varsigma e^{i\theta}) \geq 2\pi \left( \varsigma - \frac{\mathbf{a} \varsigma}{(1-\varsigma^2)^2} \right) ,\;\; \forall \varsigma \in (0,\sqrt{2} -1) ,\; \theta\in [0,2\pi ). $$
Integrating in $\varsigma$ yields
$$ \varsigma \frac{\partial \mathbf{u}_0 }{\partial \varsigma }  (\varsigma e^{i\theta}) \geq 2\pi \int_0^\varsigma t - \frac{\mathbf{a} t}{(1-t^2)^2} dt = \pi \varsigma^2 \left( 1 - \frac{\mathbf{a} }{1-\varsigma^2} \right) \geq \pi \varsigma^2 \left( 1 - \frac{\mathbf{a} }{2(\sqrt{2}-1)} \right) \geq 0 , $$
and thus
$$ \left|\frac{\partial \mathbf{u}_0 }{\partial \varsigma }  (\varsigma e^{i\theta}) \right| \geq \pi \varsigma \left| 1 - \frac{\mathbf{a} }{1-\varsigma^2} \right| . $$

Substituting back into the integral, we find
\begin{eqnarray*}
 \int_{\mathbb{D}_{\sqrt{2} -1} (0)} d\mathbf{u}_0 \wedge \star d \mathbf{u}_0 & = & \int_{\mathbb{D}_{\sqrt{2} -1} (0)} \varsigma \left|\frac{\partial \mathbf{u}_0 }{\partial \varsigma} \right|^2 d\varsigma \wedge d\theta \geq 2\pi^3 \int_0^{\sqrt{2} -1} \varsigma^3 \left(  1 - \frac{\mathbf{a} }{1-\varsigma^2} \right)^2 d\varsigma .
\end{eqnarray*}
This proves the lemma.
\end{proof}

We now proceed to establish estimates for the $\varphi$-invariant. In the following lemma, two cases are considered according to whether the $L^\infty$ norm of the Bergman kernel is large or small.

\begin{lem}
\label{lemphiinvariantthickpart}
Let $C$ be a compact Riemann surface of genus $g \geq 2$, equipped with the hyperbolic metric $\mu_{\mathrm{KE}}$ of constant curvature $-1$ and total volume $4\pi(g-1)$. Denote by $\mathrm{sys}(C)$ the systole of $(C,\mu_{\mathrm{KE}})$, and by $\mathbf{b}_0$ the Bergman kernel function associated with $(\Gamma(C,\omega_C) , \mu_{\mathrm{KE}} )$. Then the following estimates hold:
\begin{enumerate}[(1)]
    \item For any $g\geq 2$,
    \begin{eqnarray*}
    \left(1+\frac{2\arcsinh 1}{\mathrm{sys}(C)}\right) \varphi (C) & \geq & \frac{32\pi^2 g \Vert \mathbf{b}_0 \Vert_{L^\infty}^2 }{g-1} \left( \mathbf{I}_1 - \frac{2}{g} \mathbf{I}_2 +\frac{1}{g^2} \mathbf{I}_3 \right) ,
    \end{eqnarray*}
    where $\mathbf{I}_1 $, $ \mathbf{I}_2 $ and $ \mathbf{I}_3$ are the constants defined in Lemma \ref{lemintegralgradientintegralsmalldisc}.
    \item For any $g\geq 3$ and $k\in \{ 1,\cdots ,g-1 \}$, if $\Vert \mathbf{b}_0 \Vert_{L^\infty} \leq \frac{(\sqrt{2} -1)g (g-k) }{2\pi (g-1)} $, then
    \begin{eqnarray*}
    \left(1+\frac{2\arcsinh 1}{\mathrm{sys}(C)}\right) \varphi (C) & \geq & \frac{k (k+1)(2k+1) }{ 3(g-1)^2 } \mathbf{I}_1 .
    \end{eqnarray*}
\end{enumerate}
\end{lem}

\begin{proof}
    Let $\mathbf{X}_0 = \Gamma (C,\omega_C ) \cong \mathbb{C}^g $. We construct, by induction on $j$, a sequence of data $(x_j, \alpha_j, \mathbf{X}_j)$ for $j = 1, \dots, g$, satisfying:
    \begin{itemize}
        \item $x_j$ is a peak point of $(\mathbf{X}_{j-1} , \mu_{\mathrm{KE}} )$.
        \item $\alpha_j \in \mathbf{X}_{j-1} $ is a peak section in $\mathbf{X}_{j-1}$ at $x_j$.
        \item $\mathbf{X}_{j}$ is the $L^2$-orthogonal complement of $\mathbb{C}\alpha_j $ in $\mathbf{X}_{j-1}$.
    \end{itemize}
    For each $j$, let $\mathfrak{u}_j \in \mathscr{A}^0_{C,\mathbb{R}}(C)$ be the unique solution to
    $$ dd^c \mathfrak{u}_{j} = i\alpha_j \wedge\bar{\alpha}_j - 2 \mu_{\mathrm{Ar}} ,\;\;\;\; \int_C \mathfrak{u}_{j} \mu_{\mathrm{Ar}} =0 . $$

    Next, we lift the data to the unit disk $\mathbb{D}$ to obtain explicit estimates. Let $\rho_{\mathbb{D} ,j} : \mathbb{D} \to C $ be a universal covering mapping with $\rho_{\mathbb{D} ,j} (0) =x_j $. Set $\tilde{\alpha}_j = \rho^*_{\mathbb{D} ,j} \alpha_j $ and $\tilde{\mathfrak{u}}_j = \mathfrak{u}_j \circ \rho_{\mathbb{D} ,j} $. Then there are holomorphic functions $\mathbf{f}_j$ on $\mathbb{D}$ such that the pullback $\tilde{\alpha}_j = \rho^*_{\mathbb{D} ,j} \alpha_j = \mathbf{f}_j dz $. Let $\mathbf{b}_j$ denote the Bergman kernel function associated with $(\mathbf{X}_{j} , \mu_{\mathrm{KE}} )$. By Lemma \ref{lembasicpropertiesBergmankernelfunction}, we have
    \begin{equation*}
    dd^c \tilde{\mathfrak{u}}_{j} = 2|\mathbf{f}_j|^2 \mu_{\mathrm{Euc}} - \frac{2}{g} (\mathbf{b}_0 \circ \rho_{\mathbb{D} ,j} ) \rho^*_{\mathbb{D} ,j} \mu_{\mathrm{KE}} \geq \left( 2|\mathbf{f}_j|^2 - \frac{8 \Vert \mathbf{b}_0 \Vert_{L^\infty} }{g(1-|z|^2)^2} \right) \mu_{\mathrm{Euc}} ,
    \end{equation*}
    and $ |\mathbf{f}_j (0) |^2 = 4\mathbf{b}_{j-1} (x_j) = 4\Vert \mathbf{b}_{j-1} \Vert_{L^\infty} $. To estimate the $L^2$-norm of $d\mathfrak{u}_j$, Lemma \ref{lemlocalstructureorbitcollartheorem} gives
    \begin{eqnarray*}
    \int_C d\mathfrak{u}_{j } \wedge \star d \mathfrak{u}_{j } & \geq & \int_{\rho_{\mathbb{D} ,j} ( \mathbb{D}_{\sqrt{2}-1} (0) ) } d\mathfrak{u}_{j } \wedge \star d \mathfrak{u}_{j } \\
    & \geq & \frac{\mathrm{sys}(C)}{2\arcsinh 1 +\mathrm{sys}(C)} \int_{ \mathbb{D}_{\sqrt{2}-1} (0) } d\tilde{\mathfrak{u}}_{j} \wedge \star d \tilde{\mathfrak{u}}_{j} .
    \end{eqnarray*}

    We are now ready to establish the claimed estimates.

    We begin by proving (1). For $j=1$, we have
    $$ dd^c \tilde{\mathfrak{u}}_{1} \geq 8 \Vert \mathbf{b}_0 \Vert_{L^\infty}\left( \frac{ \left| \mathbf{f}_1 \right|^2}{4\Vert \mathbf{b}_0 \Vert_{L^\infty}} - \frac{ 1 }{g(1-|z|^2)^2} \right) \mu_{\mathrm{Euc}} ,\;\;\;\;\; \frac{ \left| \mathbf{f}_1 (0) \right|^2}{4\Vert \mathbf{b}_0 \Vert_{L^\infty}} = 1 . $$
    Then Lemma \ref{lemintegralgradientintegralsmalldisc} yields
    \begin{eqnarray*}
    \int_{ \mathbb{D}_{\sqrt{2}-1} (0) } d\tilde{\mathfrak{u}}_{1} \wedge \star d \tilde{\mathfrak{u}}_{1} & \geq & 128\pi^2 \Vert \mathbf{b}_0 \Vert_{L^\infty}^2 \int_0^{\sqrt{2} -1} t^3 \left(  1 - \frac{1 }{g(1-t^2)} \right)^2 dt .
    \end{eqnarray*}
    Since $dd^c ( \sum_{j=1}^g \mathfrak{u}_j ) = 0$, we have $\sum_{j=1}^g \mathfrak{u}_j = 0 $, and by Proposition \ref{propzhangphiinvariantpotential}:
    \begin{eqnarray*}
    \varphi (C) & \geq & \frac{1}{4\pi} \int_C d\mathfrak{u}_{1 } \wedge \star d \mathfrak{u}_{1 } + \frac{1}{4\pi (g-1)} \int_C d \left( \sum_{j=2}^g \mathfrak{u}_{j } \right) \wedge \star d \left( \sum_{j=2}^g \mathfrak{u}_{j } \right) \\
    & = & \frac{g}{4\pi (g-1)} \int_C d\mathfrak{u}_{1 } \wedge \star d \mathfrak{u}_{1 } \\
    & \geq & \frac{32\pi^2 g\Vert \mathbf{b}_0 \Vert_{L^\infty}^2 \cdot \mathrm{sys}(C) }{(g-1)(2\arcsinh 1 +\mathrm{sys}(C))} \left( \mathbf{I}_1 - \frac{2}{g} \mathbf{I}_2 +\frac{1}{g^2} \mathbf{I}_3 \right) ,
    \end{eqnarray*}
    where $\mathbf{I}_1 $, $ \mathbf{I}_2 $ and $ \mathbf{I}_3$ are the constants defined in Lemma \ref{lemintegralgradientintegralsmalldisc}. This proves (1).

    Now we consider the estimate (2).

    Since $\Vert \mathbf{b}_{j-1} \Vert_{L^\infty} \geq \frac{\dim \mathbf{X}_{j-1}}{\mu_{\mathrm{KE}} (C)} = \frac{g-j+1}{4\pi (g-1)} $, we have
    \begin{eqnarray*}
    dd^c \tilde{\mathfrak{u}}_{j} & \geq & 8 \Vert \mathbf{b}_{j-1} \Vert_{L^\infty}\left( \frac{ \left| \mathbf{f}_j \right|^2}{4\Vert \mathbf{b}_{j-1} \Vert_{L^\infty}} - \frac{ \Vert \mathbf{b}_{0} \Vert_{L^\infty} }{g\Vert \mathbf{b}_{j-1} \Vert_{L^\infty}(1-|z|^2)^2} \right) \mu_{\mathrm{Euc}}  \\
    & \geq & \frac{2(g-j+1)}{\pi (g-1)} \left( \frac{ \left| \mathbf{f}_j \right|^2}{4\Vert \mathbf{b}_{j-1} \Vert_{L^\infty}} - \frac{ 2(\sqrt{2} -1 )(g-k) }{(g-j+1) (1-|z|^2)^2} \right) \mu_{\mathrm{Euc}}
    \end{eqnarray*}
    and $\frac{ \left| \mathbf{f}_j (0) \right|^2}{4\Vert \mathbf{b}_{j-1} \Vert_{L^\infty}} = 1 $.  The same argument as in (1) gives
    \begin{eqnarray*}
    \varphi (C) & \geq & \frac{ \mathrm{sys}(C)}{4\pi (2\arcsinh 1 +\mathrm{sys}(C))} \sum_{j=1}^{k+1} \int_{ \mathbb{D}_{\sqrt{2}-1} (0) } d\tilde{\mathfrak{u}}_{j} \wedge \star d \tilde{\mathfrak{u}}_{j} \\
    & \geq & \frac{32\pi^2 \cdot \mathrm{sys}(C) }{ 2\arcsinh 1 +\mathrm{sys}(C) } \sum_{j=0}^{k } \int_0^{\sqrt{2} -1} t^3 \left(  \frac{g-j}{4\pi (g-1)} - \frac{g-k}{4\pi (g-1)} \right)^2 dt \\
    & = & \frac{ k(k+1)(2k+1) \mathrm{sys}(C) }{ 3 (g-1)^2 ( 2\arcsinh 1 +\mathrm{sys}(C) ) } \mathbf{I}_1 ,
    \end{eqnarray*}
    where $\mathbf{I}_1 $ is the constant defined in Lemma \ref{lemintegralgradientintegralsmalldisc}. This is our assertion.
\end{proof}

For any $g\geq 3$ and $k\in \{ 1,\cdots ,g-1 \}$, let $\mathbf{I}_4 (g,k)$ denote
$$\mathbf{I}_4 (g,k) = \min \left\{ \frac{k (k+1)(2k+1) }{ 3(g-1)^2 } \mathbf{I}_1 , \frac{ 8(\sqrt{2}-1)^2 (g-k)^2 g^3 }{(g-1)^3} \left( \mathbf{I}_1 - \frac{2}{g} \mathbf{I}_2 +\frac{1}{g^2} \mathbf{I}_3 \right)  \right\} ,$$ 
where $\mathbf{I}_1 $, $ \mathbf{I}_2 $ and $ \mathbf{I}_3$ are the constants defined in Lemma \ref{lemintegralgradientintegralsmalldisc}. As a corollary of Lemma \ref{lemphiinvariantthickpart}, we obtain the following estimate.

\begin{cor}
\label{corphiinvariantthickpart}
Let $C$ be a compact Riemann surface of genus $g \geq 2$, equipped with the hyperbolic metric $\mu_{\mathrm{KE}}$ of constant curvature $-1$. Let $\mathrm{sys}(C)$ denote the systole of $(C,\mu_{\mathrm{KE}})$. Then, for any $k\in \{ 1,\cdots ,g-1 \}$, we have
\begin{equation*}
\left(1+\frac{2\arcsinh 1}{\mathrm{sys}(C)}\right) \varphi (C) \geq \max \left\{ \mathbf{I}_4 (g,k) ,\frac{ 2\cdot g^3 }{(g-1)^3} \left( \mathbf{I}_1 - \frac{2}{g} \mathbf{I}_2 +\frac{1}{g^2} \mathbf{I}_3 \right) \right\} .
\end{equation*}
\end{cor}

\begin{proof}
Since $\Vert \mathbf{b}_{0} \Vert_{L^\infty} \geq \frac{\dim \mathbf{X}_{0}}{\mu_{\mathrm{KE}} (C)} = \frac{g }{4\pi (g-1)} $, we can apply Lemma \ref{lemintegralgradientintegralsmalldisc} to deduce that
\begin{equation*}
\left(1+\frac{2\arcsinh 1}{\mathrm{sys}(C)}\right) \varphi (C) \geq \frac{ 2\cdot g^3 }{(g-1)^3} \left( \mathbf{I}_1 - \frac{2}{g} \mathbf{I}_2 +\frac{1}{g^2} \mathbf{I}_3 \right) .
\end{equation*}
Moreover, if $\Vert \mathbf{b}_0 \Vert_{L^\infty} \geq \frac{(\sqrt{2} -1)g (g-k) }{2\pi (g-1)} $, then Lemma \ref{lemintegralgradientintegralsmalldisc} further yields
\begin{eqnarray*}
\left(1+\frac{2\arcsinh 1}{\mathrm{sys}(C)}\right) \varphi (C) & \geq & \frac{ 8(\sqrt{2}-1)^2 (g-k)^2 g^3 }{(g-1)^3} \left( \mathbf{I}_1 - \frac{2}{g} \mathbf{I}_2 +\frac{1}{g^2} \mathbf{I}_3 \right) .
\end{eqnarray*}

On the other hand, if $\Vert \mathbf{b}_0 \Vert_{L^\infty} \leq \frac{(\sqrt{2} -1)g (g-k) }{2\pi (g-1)} $, then Lemma \ref{lemintegralgradientintegralsmalldisc} implies
\begin{eqnarray*}
\left(1+\frac{2\arcsinh 1}{\mathrm{sys}(C)}\right) \varphi (C) & \geq & \frac{k (k+1)(2k+1) }{ 3(g-1)^2 } \mathbf{I}_1 .
\end{eqnarray*}
This completes the proof.
\end{proof}

Now we are ready to prove the following estimate for the $\varphi$-invariant.

\begin{prop}
\label{propexplicitlowerboundphiinvthickpart}
Let $C$ be a compact Riemann surface of genus $g \geq 2$, equipped with the hyperbolic metric $\mu_{\mathrm{KE}}$ of constant curvature $-1$. Denote by $\mathrm{sys}(C)$ the systole of $(C,\mu_{\mathrm{KE}})$. Then the $\varphi$-invariant satisfies
$$ \left(1+\frac{2\arcsinh 1}{\mathrm{sys}(C)}\right) \varphi (C) \geq \mathbf{I}_5 (g) ,$$
where $\mathbf{I}_5 (g) $ is defined as follows:
$$
\begin{array}{|c|c|c|c|c|c|c|c|c|c|c|c|c|}
\hline
	g		& 2 & 3 	& 4 & 5 		& 6 	& 7 	& 8 		& 9    & 10&11&12&\geq 13\\\hline
    \mathbf{I}_5(g)	& \frac{11}{495} & \frac{1}{52} 	& \frac{1}{56} & \frac{3}{175}	& \frac{7}{400} & \frac{7}{312}	& \frac{17}{625}	& \frac{4}{125} & \frac{1}{27} & \frac{1}{24} & \frac{2}{47} &\frac{g}{400}\\
\hline
\end{array}.
$$
\end{prop}

\begin{proof}
By Lemma \ref{lemphiinvariantthickpart}, for any $g\geq 2$,
\begin{eqnarray*}
    \left(1+\frac{2\arcsinh 1}{\mathrm{sys}(C)}\right) \varphi (C) & \geq & \frac{2g^3 }{ (g-1)^{3} }  \left( \mathbf{I}_1 - \frac{2}{g} \mathbf{I}_2 +\frac{1}{g^2} \mathbf{I}_3 \right) ,
\end{eqnarray*}
where $\mathbf{I}_1 $, $ \mathbf{I}_2 $ and $ \mathbf{I}_3$ are the constants defined in Lemma \ref{lemintegralgradientintegralsmalldisc}. For $g\in\{2,3,4,5\}$, the right-hand side is given by:
\begin{center}
\begin{tabular}{|c|c|c|}
\hline
$g$ & expression & numerical value \\ \hline
$2$ & $16\left( \mathbf{I}_1 - \mathbf{I}_2 +\tfrac{1}{4} \mathbf{I}_3 \right)$ & $0.02228\cdots > \tfrac{11}{495}$ \\ \hline
$3$ & $\tfrac{27}{4}\left( \mathbf{I}_1 - \tfrac{2}{3} \mathbf{I}_2 +\tfrac{1}{9} \mathbf{I}_3 \right)$ & $0.01928\cdots > \tfrac{1}{52}$ \\ \hline
$4$ & $\tfrac{128}{27}\left( \mathbf{I}_1 - \tfrac{1}{2} \mathbf{I}_2 +\tfrac{1}{16} \mathbf{I}_3 \right)$ & $0.01794\cdots > \tfrac{1}{56}$ \\ \hline
$5$ & $\tfrac{125}{32}\left( \mathbf{I}_1 - \tfrac{2}{5} \mathbf{I}_2 +\tfrac{1}{25} \mathbf{I}_3 \right)$ & $0.01721\cdots > \tfrac{3}{175}$ \\
\hline
\end{tabular}
\end{center}
 
By Lemma \ref{lemphiinvariantthickpart}, for any $g\geq 2$ and $k\in \{ 1,\cdots ,g-1 \}$,
\begin{eqnarray*}
    \left(1+\frac{2\arcsinh 1}{\mathrm{sys}(C)}\right) \varphi (C) & \geq & \mathbf{I}_4 (g,k) ,
\end{eqnarray*}
where 
$$\mathbf{I}_4 (g,k) = \min \left\{ \frac{k (k+1)(2k+1) }{ 3(g-1)^2 } \mathbf{I}_1 , \frac{ 8(\sqrt{2}-1)^2 (g-k)^2 g^3 }{(g-1)^3} \left( \mathbf{I}_1 - \frac{2}{g} \mathbf{I}_2 +\frac{1}{g^2} \mathbf{I}_3 \right)  \right\} .$$

For $g=6,\cdots , 12$, define $k(g)$ as follows:
$$
\begin{array}{|c|c|c|c|c|c|c|c|}
\hline
	g &6 	& 7 	& 8 		& 9    & 10&11&12 \\\hline
     k(g)	& 4 & 5 	& 6 & 7	& 8 & 9	& 10	\\
\hline
\end{array}.
$$
A direct calculation shows that $\mathbf{I}_4(g,k(g)) > \mathbf{I}_5 (g)$ for all $6 \leq g \leq 12$.

We now turn to the case $g \geq 13$. Let $k(g)$ be the unique integer in 
$$\left[ g-\frac{(\sqrt{2}+1)\sqrt{g}}{2\sqrt{3}} - \frac{1}{3},g-\frac{(\sqrt{2}+1)\sqrt{g}}{2\sqrt{3}} + \frac{2}{3} \right).$$
Define $\mathbf{p} (x) = 1-\frac{\sqrt{2}+1}{2\sqrt{3x}} +\frac{2}{3x} $. Since $\mathbf{p}' (x) = \frac{\sqrt{2}+1}{4x\sqrt{3x}} -\frac{2}{3x^2} >0 $ for $x\geq 4$, it follows that $\mathbf{p}(g)$ is an increasing function of $g$. Similarly, for $g \geq 13$, the quantity $\frac{ g^3 }{(g-1)^3} \left( \mathbf{I}_1 - \frac{2}{g} \mathbf{I}_2 +\frac{1}{g^2} \mathbf{I}_3 \right) $ is decreasing in $g$, while $\frac{(\sqrt{2}+1)}{2\sqrt{3}} - \frac{2}{3\sqrt{g}}$ is increasing in $g$. 

Hence
\begin{eqnarray*}
\frac{k (k+1)(2k+1) }{ 3g(g-1)^2 } \mathbf{I}_1 & \geq & \frac{2\mathbf{I}_1}{3} \cdot \mathbf{p} (13) \left( \mathbf{p} (13) -\frac{1}{26} \right)\left( \mathbf{p} (13) -\frac{1}{13} \right) \\
& \approx & 0.0026945 \cdots >\frac{1}{400} ,
\end{eqnarray*}
and
\begin{eqnarray*}
 \frac{ 8(\sqrt{2}-1)^2 (g-k)^2 g^3 }{g(g-1)^3} \left( \mathbf{I}_1 - \frac{2}{g} \mathbf{I}_2 +\frac{1}{g^2} \mathbf{I}_3 \right) & \geq & \frac{8(\sqrt{2}-1)^2 (g-k)^2 \mathbf{I}_1}{g} \\
& \geq & 8(\sqrt{2}-1)^2 \cdot \left( \frac{(\sqrt{2}+1)}{2\sqrt{3}} - \frac{2}{3\sqrt{13}} \right)^2 \cdot \frac{(\sqrt{2}-1)^4}{4} \\
& \approx & 0.0026482\cdots >\frac{1}{400} .
\end{eqnarray*}
Therefore $\mathbf{I}_4 (g,k(g))> \mathbf{I}_5 (g)$ for all $g\geq 13$, which proves the proposition.
\end{proof}

\subsection{Separating case}
\label{sectionzhangphi2}

As above, let $C$ be a curve of genus $g \geq 2$ over $\mathbb{C}$, and let $\omega_C$ denote its canonical bundle. Write $\mu_{\mathrm{KE}}$ for the unique Kähler metric on $C$ with constant curvature $-1$.

In this subsection we estimate the $\varphi$-invariant in the case where a finite collection of short geodesics separates $C$ into two components. A finite collection of disjoint simple closed geodesics is called a separating multicurve if it divides $C$ into at least two connected components. Such a multicurve is minimal if no proper subset of it is separating.

As a corollary of Poincar\'e duality, a collection $\{\gamma_j\}_{j=1}^m$ of disjoint simple closed geodesics is separating if and only if the $\mathbb{C}$-span $ \sum_{j=1}^m \mathbb{C} [\gamma_j ] \subset H_1 (C;\mathbb{C}) $ has dimension at most $m-1$. Note also that for any open subset $U\subset C$ with smooth boundary $\partial U$, the chain $\partial U$ is a nontrivial linear combination of its connected components, yet its homology class in $H_1(C;\mathbb{C})$ vanishes.

We begin with the following topological lemma.

\begin{lem}
\label{lemtopologyseparatingmulticurvehomologicaldescripton}
Let $C$ be a compact Riemann surface of genus $g\geq 2$, and let $\gamma_1,\dots,\gamma_m$ be pairwise disjoint simple closed geodesics on $C$. Suppose $\{\gamma_j\}_{j=1}^m$ forms a minimal separating multicurve. Then $C\setminus \bigcup_{j=1}^m \gamma_j$ has exactly two connected components, say $C_1$ and $C_2$. Their closures $\overline{C}_1,\overline{C}_2$ are compact surfaces with common boundary $\bigcup_{j=1}^m \gamma_j$, and if $g_k$ denotes the genus of $\overline{C}_k$, then $g_1+g_2+m-1=g$.

Moreover, there exists a collection of smooth loops $\{\mathbf{c}_j\}_{j=1}^{2(g_1+g_2)}$ on $C$ such that:
\begin{enumerate}[(1)]
    \item $\mathbf{c}_j \subset C_1$ for $1\leq j\leq 2g_1$, and $\mathbf{c}_j \subset C_2$ for $2g_1+1\leq j\leq 2(g_1+g_2)$;
    \item $\{\mathbf{c}_j\}_{j=1}^{g_1}\cup \{\mathbf{c}_j\}_{j=2g_1+1}^{2g_1+g_2}$ consists of pairwise disjoint loops;
    \item $\{[\mathbf{c}_j]\}_{j=1}^{2(g_1+g_2)}$ are $\mathbb{C}$-linearly independent in $H_1(C;\mathbb{C})$;
    \item $\bigcup_{j=1}^m \gamma_j$ lies in a single connected component of both 
    $\overline{C}_1 \setminus \bigcup_{j=1}^{2g_1}\mathbf{c}_j$ and $\overline{C}_2 \setminus \bigcup_{j=2g_1+1}^{2(g_1+g_2)}\mathbf{c}_j$.
\end{enumerate}
\end{lem}

\begin{proof}
This is a purely topological fact and follows directly from the classification theorem of compact surfaces with boundary (see \cite[Chapter~I, 41A]{alhsait1}). We therefore omit the proof.
\end{proof}

To apply Tian's peak section method, we need the following lemma, which may be regarded as an analogue of H\"ormander's $L^2$ estimate appearing in Tian's original paper \cite{tg1}. Recall that harmonic $1$-forms are precisely the peak sections along cycles defined earlier.

\begin{lem}
\label{lemconstructionoflocalizationoperator}
Let $C $ be a compact Riemann surface of genus $g\geq 2$, and let $\{\gamma_j \}_{j=1}^m $ be a minimal separating multicurve on $C$ with total length $ \sum_{j=1}^m \ell (\gamma_j) \leq \frac{1}{4}$. For any piecewise smooth loop $\mathbf{c}:\mathbb{R}/\mathbb{Z} \to C $, denote by $\mathbf{T}_{\mathbf{c}} $ the linear functional on $\mathscr{Z}_{C,\mathbb{C}}^1 (C ) $ given by $\beta \mapsto \int_{\mathbf{c}} \beta $. Let $K\in\{\mathbb{R},\;\mathbb{C}\} $ and set $\mathbf{X}\subset \mathscr{H}_{C,K}^{1} (C) $ to be the subspace of harmonic $1$-forms $\beta$ satisfying $\mathbf{T}_{\mathbf{\gamma_j}} (\beta) = \mathbf{T}_{\mathbf{\gamma_j}} (\star \beta) =0 $, $j=1,\cdots ,m$. Let $C_1$ be one of the two connected components of $C\setminus \left( \bigcup_{j=1}^m \gamma_j \right) $, and let $\{\mathbf{c}_j\}_{j=1}^{2(g_1 +g_2)}$ be the system of smooth loops constructed in Lemma \ref{lemtopologyseparatingmulticurvehomologicaldescripton}. 

Then there exists a linear operator
$$ \Psi_{ C_1, \mathrm{loc}} : \mathbf{X} \longrightarrow \mathscr{Z}^1_{C;K} (C) ,$$
satisfying:
\begin{enumerate}[(1)]
    \item The image of $\Psi_{ C_1, \mathrm{loc}}$ satisfies 
    $$\Psi_{ C_1, \mathrm{loc}} (\mathbf{X} ) \subset \left( \bigcap\limits_{j=2 g_1 +1}^{2(g_1 +g_2 )} \mathrm{ker} \mathbf{T}_{\mathbf{c}_j} \right) \bigcap \left( \bigcap\limits_{j= 1}^{m} \mathrm{ker} \mathbf{T}_{\gamma_j} \right) .$$
    \item The deviation from identity satisfies
    $$ \left( \Psi_{ C_1, \mathrm{loc}} - \mathrm{id}_\mathbf{X} \right) (\mathbf{X} ) \subset \bigcap\limits_{j= 1}^{2 g_1 } \mathrm{ker} \mathbf{T}_{\mathbf{c}_j} ,$$ 
    where $\mathrm{id}_\mathbf{X} $ is the identity map of $\mathbf{X}$.
    \item For any $\beta\in\mathbf{X} $,
    $$\Vert \Psi_{ C_1, \mathrm{loc}} (\beta) \Vert^2_{L^2 ;C } \leq \Vert  \beta \Vert^2_{L^2 ;C_1 } + 200 e^{ - \frac{\pi^2}{\sum_{j=1}^m \ell (\gamma_j)}} \Vert  \beta \Vert^2_{L^2 ;C } , $$
    where $\Vert  \beta \Vert^2_{L^2 ;C_1 } = \frac{1}{2} \int_{C_1} \beta \wedge \star \beta $.
\end{enumerate}
\end{lem}

\begin{proof}
Without loss of generality, we assume $K = \mathbb{R}$. We will use the local structure of the K\"ahler collars $\mathscr{C}(\gamma_j)$ to construct the linear operator $\Psi_{C_1, \mathrm{loc}}$. By Lemma \ref{lemkahlercollar}, there exist holomorphic coordinates 
$$z_j : \mathscr{C} (\gamma_j ) \to \overline{\mathbb{D}_{e^{-\nu (\gamma_j)}} (0) }\setminus \mathbb{D}_{e^{\nu (\gamma_j ) -\frac{2\pi^2}{\ell (\gamma_j )} }} (0) $$
such that $z_j \left( \mathscr{C} (\gamma_j ) \cap C_1 \right) = \overline{\mathbb{D}_{e^{-\nu (\gamma_j)}} (0)} \setminus \mathbb{D}_{ e^{-\frac{\pi^2}{\ell (\gamma_j )} }} (0) $, where $\nu (\gamma_j) \in (0, \frac{\pi}{2\ell (\gamma_j)}) $ satisfies
$$\nu (\gamma_j) \ell (\gamma_j) = {2\pi} \arcsin \left( \tanh  \frac{\ell (\gamma_j )}{2} \right) .$$

Fix $x_j\in\gamma_j $ for $j=1,\cdots ,m$. Then $|z_j (x_j)|=e^{-\frac{\pi^2}{\ell(\gamma_j)}}$. Let $\beta $ be a harmonic $1$-form on $C$ satisfying $\mathbf{T}_{\mathbf{\gamma_j}} (\beta) = \mathbf{T}_{\mathbf{\gamma_j}} (\star \beta) =0 $ for all $j$. For each $j$, there exists a unique function $u_{\beta,j}$ on $\mathscr{C}(\gamma_j )$ such that $du_{\beta,j} = \beta + i\star \beta $ and $u_{\beta,j} (x_j)=0$. Note that $\frac{\pi^2}{\ell (\gamma_j)} \geq 4\pi^2 > 6\nu (\gamma_j) $. Since $\beta + i\star \beta $ is holomorphic, Proposition \ref{propexplicitestimateskahlercollar} implies that for any $x\in \mathscr{C} (\gamma_j ) \cap C_1 $ with $|z_j(x)|\leq e^{-2\nu(\gamma_j)}$, we have
\begin{eqnarray*}
|u_{\beta,j} (x) |^2 & \leq & \frac{ e^{2\nu (\gamma_j)} \Vert\beta + i\star \beta\Vert^2_{L^2 ;\mathscr{C}(\gamma_j)} }{\pi (1-e^{-3\pi})} \cdot \frac{|z_j (x) - z_j (x_j )|^2}{\left( 1-|z_j (x)|^2 e^{2\nu (\gamma_j )} \right)^2} \\
& & + \frac{ e^{2\nu (\gamma_j)} \Vert\beta + i\star \beta\Vert^2_{L^2 ;\mathscr{C}(\gamma_j)} }{\pi (1-e^{-3\pi})} \cdot \frac{|z_j (x)^{-1} - z_j (x_j )^{-1} |^{2} e^{-\frac{4\pi^2}{\ell (\gamma_j)}} }{\left(1- e^{2\left( \nu (\gamma_j) - \frac{\pi^2}{\ell (\gamma_j)} \right)} \right)^2} \\
& \leq & \frac{ 2 e^{2\nu (\gamma_j)} \Vert\beta \Vert^2_{L^2 ;\mathscr{C}(\gamma_j)} }{\pi (1-e^{-3\pi}) } \left( \frac{\left(|z_j (x)| + e^{-\frac{\pi^2}{\ell(\gamma_j)}}\right)^2}{\left( 1-|z_j (x)|^2 e^{2\nu (\gamma_j )} \right)^2} + \frac{4 e^{-\frac{2\pi^2}{\ell (\gamma_j)}}}{(1-e^{-4\pi^2})^2} \right) .
\end{eqnarray*}
Here we used
$$\Vert\beta + i\star \beta\Vert^2_{L^2 ;\mathscr{C}(\gamma_j)} = \frac{1}{2} \int_{\mathscr{C}(\gamma_j)} (\beta + i\star \beta) \wedge \star \overline{(\beta + i\star \beta)} = \int_{\mathscr{C}(\gamma_j)} \beta \wedge \star \beta = 2 \Vert\beta \Vert^2_{L^2 ;\mathscr{C}(\gamma_j)} . $$

There are two subcases to consider.

If $|z_j (x)|^2 e^{2\nu (\gamma_j )} \leq e^{-3\pi} $, then $1-|z_j (x)|^2 e^{2\nu (\gamma_j )} \geq 1-e^{-3\pi } $. By the Cauchy--Schwarz inequality, we obtain
$$ |u_{\beta,j} (x) |^2 \leq \frac{ 4 e^{2\nu (\gamma_j)} \Vert\beta \Vert^2_{L^2 ;\mathscr{C}(\gamma_j)} }{\pi (1-e^{-3\pi})^3} \left( |z_j (x)|^2 + 3e^{- \frac{2\pi^2}{\ell (\gamma_j)}} \right) .$$

If $e^{-3\pi} \leq |z_j (x)|^2 e^{2\nu (\gamma_j )} \leq e^{-2\nu (\gamma_j )} $, then
\begin{eqnarray*}
\frac{\left(|z_j (x)| + e^{-\frac{\pi^2}{\ell(\gamma_j)}}\right)^2}{\left( 1-|z_j (x)|^2 e^{2\nu (\gamma_j )} \right)^2} & \leq & \frac{\left(1 + e^{\frac{3\pi}{2}-\left(\frac{\pi^2}{\ell(\gamma_j)}-\nu (\gamma_j )\right)}\right)^2}{\left( 1- e^{-2\nu (\gamma_j )} \right)^2} \cdot |z_j (x)|^2 \\
& \leq & \frac{(1+e^{-2\pi^2})^2}{(1-e^{-4})^2} \cdot |z_j (x)|^2 < 2 |z_j (x)|^2 .
\end{eqnarray*}
Consequently, for any $x\in z_j^{-1} \left( \overline{\mathbb{D}_{e^{-2\nu (\gamma_j)}} (0)} \setminus \mathbb{D}_{e^{ -\frac{ \pi^2}{\ell (\gamma_j )} }} (0) \right) \subset \mathscr{C} (\gamma_j ) \cap C_1 $, we obtain the bound
$$ |u_{\beta,j} (x) |^2 \leq \frac{ 4 e^{2\nu (\gamma_j)} \Vert\beta \Vert^2_{L^2 ;\mathscr{C}(\gamma_j)} }{\pi (1-e^{-3\pi})^3} \left( |z_j (x)|^2 + 3e^{- \frac{2\pi^2}{\ell (\gamma_j)}} \right) .$$

Next, let $\eta_{j,0} \in C^{0,1}(\mathbb{R}; [0,1])$ be a cut-off function defined by
\begin{eqnarray*}
\eta_{j,0} (r) = \left\{
\begin{aligned}
0, \;\;\;\;\;\;\;\;\;\;\;\;\;\;\;\;\;\;\;\;\;\;\;\; &\;\;\;  \;\;\;\;\;\;\;\;\;\;\;\;\;\;\; r< e^{-\frac{\pi^2}{\ell (\gamma_j )}} ;\\
 \frac{2\log2 -\log \left( 1+3r^{-2} e^{-\frac{2\pi^2}{\ell (\gamma_j )}} \right)}{2\log2 -\log \left( 1+3 e^{4\nu (\gamma_j )-\frac{2\pi^2}{\ell (\gamma_j )}} \right)}, \; &\;\;\;  e^{-\frac{\pi^2}{\ell (\gamma_j )}} \leq r\leq e^{-2\nu (\gamma_j)} ;\\
 1,\;\;\;\;\;\;\;\;\;\;\;\;\;\;\;\;\;\;\;\;\;\;\;\; &\;\;\;  \;\;\;\;\;\;\;\;\;\;\;\;\;\;\; r\geq e^{-2\nu (\gamma_j )}.
\end{aligned}
\right.
\end{eqnarray*}
It is straightforward to check that $\eta_{j,0}'\in L^\infty (\mathbb{R} ;[0,\infty)) $. Then $ u_{\beta,j} d( \eta_{j,0} ( | z_j|))$ is an $1$-form with $L^\infty$-coefficients on $\mathscr{C}(\gamma_j) $, whose $L^2$-norm can be estimated using Proposition \ref{propexplicitestimateskahlercollar}:
\begin{eqnarray*}
& & \left\Vert \re(u_{\beta,j}) d( \eta_{j,0} (| z_j|)) \right\Vert^2_{L^2 ;\mathscr{C} (\gamma_j)} \leq \frac{1}{2}\int_{e^{-\frac{\pi^2}{\ell (\gamma_j)}} \leq|z_j|\leq e^{-2\nu (\gamma_j)}}  |u_{\beta,j} |^2 |\eta_{j,0}' (|z_j|)|^2 d\mu_{\mathrm{Euc}} \\
& \leq & \frac{ 4 e^{2\nu (\gamma_j)} \Vert\beta \Vert^2_{L^2 ;\mathscr{C}(\gamma_j)} }{ (1-e^{-3\pi})^3 } \int_{e^{-\frac{\pi^2}{\ell (\gamma_j )}}}^{e^{-2\nu (\gamma_j)} } t \left( t^2 + 3e^{- \frac{2\pi^2}{\ell (\gamma_j)}} \right) |\eta_{j,0}' (t)|^2 dt \\
& \leq & \frac{ 4 e^{2\nu (\gamma_j)} \Vert\beta \Vert^2_{L^2 ;\mathscr{C}(\gamma_j)} }{ (1-e^{-3\pi})^3 } \cdot \frac{6e^{-\frac{2\pi^2}{\ell (\gamma_j)}}}{2\log 2 -3e^{-6\pi }} \leq \frac{240 e^{2\pi- \frac{2\pi^2}{\ell (\gamma_j)}} }{13} \Vert\beta \Vert^2_{L^2 ;\mathscr{C}(\gamma_j)} .
\end{eqnarray*}
Note that $\left(1-e^{-3\pi}\right)^{4}>1-\frac{1}{1000}$ and $2\log 2\approx 1.386\cdots >\frac{13}{10} \cdot (1+\frac{3}{50}) .$
Hence the sum can be bounded by
\begin{eqnarray*}
 & &\left\Vert \re(u_{\beta,j} ) d( \eta_{j,0} (| z_j|)) \right\Vert^2_{L^2 ;\mathscr{C} (\gamma_j)} + 2\Vert\beta \Vert_{L^2 ;\mathscr{C}(\gamma_j)} \left\Vert \re(u_{\beta,j} ) d( \eta_{j,0} (| z_j|)) \right\Vert_{L^2 ;\mathscr{C} (\gamma_j)} \\
& \leq & \left(\frac{240 e^{2\pi- \frac{2\pi^2}{\ell (\gamma_j)}} }{13} + 2\sqrt{\frac{240 e^{2\pi- \frac{2\pi^2}{\ell (\gamma_j)}} }{13}} \right) \Vert\beta \Vert^2_{L^2 ;\mathscr{C}(\gamma_j)} \\
& \leq & \left(\frac{240 e^{2\pi-4\pi^2} }{13} + 2\sqrt{\frac{240 e^{2\pi } }{13}} \right)e^{- \frac{\pi^2}{\ell (\gamma_j)}} \Vert\beta \Vert^2_{L^2 ;\mathscr{C}(\gamma_j)} \leq 199  e^{- \frac{\pi^2}{\ell (\gamma_j)}} \Vert\beta \Vert^2_{L^2 ;\mathscr{C}(\gamma_j)} .
\end{eqnarray*}
Note that $\left(\frac{240 e^{2\pi-4\pi^2} }{13} + 2\sqrt{\frac{240 e^{2\pi } }{13}} \right) \approx 198.8567\cdots <199 .$ Therefore,
\begin{eqnarray*}
\left\Vert d \left( ( \eta_{j,0} (| z_j|))\re(u_{\beta,j} ) \right) \right\Vert^2_{L^2 ;\mathscr{C} (\gamma_j)} & \leq &\left\Vert \re(u_{\beta,j} ) d( \eta_{j,0} (| z_j|)) \right\Vert^2_{L^2 ;\mathscr{C} (\gamma_j)} \\
& & + 2\Vert\beta \Vert_{L^2 ;\mathscr{C}(\gamma_j)} \left\Vert \re(u_{\beta,j} ) d( \eta_{j,0} (| z_j|)) \right\Vert_{L^2 ;\mathscr{C} (\gamma_j)} + \Vert\beta \Vert^2_{L^2 ;\mathscr{C}(\gamma_j) \cap C_1} \\
& \leq & \Vert\beta \Vert^2_{L^2 ;\mathscr{C}(\gamma_j) \cap C_1} + 199  e^{- \frac{\pi^2}{\ell (\gamma_j)}} \Vert\beta \Vert^2_{L^2 ;\mathscr{C}(\gamma_j)} .
\end{eqnarray*}

To obtain smooth cut-off functions, we mollify $\eta_{j,0}$. For any $\epsilon > 0$, there exists $\eta_{j,\epsilon} \in C^\infty(\mathbb{R}; [0,1])$ such that
$$ \eta_{j,\epsilon} =0 \textrm{ on } \big( -\infty , e^{-\frac{\pi^2}{\ell(\gamma_j)}} \big] ,\;\;\;\;\;\; \eta_{j,\epsilon} =1 \textrm{ on } \big[ e^{-2\nu (\gamma_j )},\infty \big) , $$
and
$$ \sup | \eta_{j,\epsilon} - \eta_{j,0} |^2 + \int_{e^{-\nu (\gamma_j)-\frac{\pi^2}{\ell(\gamma_j)}}}^{e^{-\nu (\gamma_j)}} | \eta'_{j,\epsilon} (t) - \eta'_{j,0} (t) |^2 dt \leq \epsilon . $$
It follows that
$$ \left\Vert d\bigg( ( \eta_{j,0} (| z_j|) - \eta_{j,\epsilon} (| z_j |)) \re(u_{\beta,j}) \bigg) \right\Vert^2_{L^2 ;\mathscr{C} (\gamma_j)} \leq 100\epsilon \Vert\beta\Vert_{L^2;\mathscr{C} (\gamma_j)}^2 .$$

Then for any $\epsilon >0$, we can define a linear operator
$$ \Psi_{ C_1, \mathrm{loc} ; \epsilon} : \mathbf{X} \longrightarrow \mathscr{Z}^1_{C;K} (C) $$
by
\begin{eqnarray*}
\Psi_{ C_1, \mathrm{loc} ; \epsilon} (\beta ) = \left\{
\begin{aligned}
 0 , \;\;\;\;\;\;\;\;\;\;\;\;\;\;\;\; &\;\;\; \textrm{ on $C\setminus \left( C_1 \bigcup \left(\bigcup_{j=1}^m \mathscr{C} (\gamma_j) \right) \right) $} ;\\
 d\bigg( \eta_{j,\epsilon} ( |z_j| ) \re(u_{\beta,j}) \bigg), \; &\;\;\;  \textrm{ on $ \mathscr{C} (\gamma_j) $} ;\\
 \beta ,\;\;\;\;\;\;\;\;\;\;\;\;\;\;\;\; &\;\;\; \textrm{ on $ C_1 \setminus \left(\bigcup_{j=1}^m \mathscr{C} (\gamma_j) \right) $}.
\end{aligned}
\right.
\end{eqnarray*}

By the construction, one can easily see that for any $\epsilon >0$, $\Psi_{ C_1, \mathrm{loc} ; \epsilon}$ satisfies (1) and (2). Now we consider the estimate (3). By the above argument, we can conclude that there exists a small constant $0<\epsilon\leq e^{-\sum_{j=1}^m \frac{100\pi^2}{\ell (\gamma_j)} } $ such that for any $\beta\in\mathbf{X}$,
\begin{eqnarray*}
 \Vert \Psi_{ C_1, \mathrm{loc} ; \epsilon} (\beta ) \Vert_{L^2;C} & = & \frac{1}{2}\int_{C_1 \setminus \left(\bigcup_{j=1}^m \mathscr{C} (\gamma_j) \right)} \beta \wedge \star \beta + \sum_{j=1}^m  \left\Vert d\left( \eta_{j,\epsilon} (| z_j |) \re(u_{\beta,j} ) \right) \right\Vert^2_{L^2 ;\mathscr{C} (\gamma_j)} \\
& \leq & \frac{1}{2}\int_{C_1 \setminus \left(\bigcup_{j=1}^m \mathscr{C} (\gamma_j) \right)} \beta \wedge \star \beta \\
& & + \sum_{j=1}^m \left( \left\Vert d\left( \eta_{j,0} (| z_j |) \re(u_{\beta,j} ) \right) \right\Vert_{L^2 ;\mathscr{C} (\gamma_j)} + e^{- \frac{50\pi^2}{\ell (\gamma_j)}} \Vert\beta \Vert_{L^2 ;\mathscr{C}(\gamma_j)} \right)^2 \\
&\leq & \frac{1}{2}\int_{C_1 \setminus \left(\bigcup_{j=1}^m \mathscr{C} (\gamma_j) \right)} \beta \wedge \star \beta + \sum_{j=1}^m \left( \Vert\beta \Vert^2_{L^2 ;\mathscr{C}(\gamma_j) \cap C_1} + 200  e^{- \frac{\pi^2}{\ell (\gamma_j)}} \Vert\beta \Vert^2_{L^2 ;\mathscr{C}(\gamma_j)} \right) \\
& \leq & \Vert  \beta \Vert^2_{L^2 ;C_1 } + 200 e^{ - \frac{\pi^2}{\sum_{j=1}^m \ell (\gamma_j)}} \Vert  \beta \Vert^2_{L^2 ;C } .
\end{eqnarray*}

This completes the proof.
\end{proof}

By combining the above lemma with Tian's peak section method, we obtain the following lemma, which may be regarded as a version of the geometric localization principle for certain holomorphic $1$-forms. For related results concerning the geometric localization principle, see \cite{mamari1, sxz1}.

\begin{lem}
\label{lemlocalizationpeaksectioncutsurfaces}
Let $C $ be a compact Riemann surface of genus $g\geq 2$, and let $\{\gamma_j \}_{j=1}^m $ be a minimal separating multicurve on $C$ such that $ \sum_{j=1}^m \ell (\gamma_j) \leq \frac{1}{4} $. Let $C_1$, $C_2 $ be the two connected components of $C\setminus \left( \bigcup_{j=1}^m \gamma_j \right) $, and let $g_1$ and $g_2$ denote the genera of $\overline{C}_1$ and $\overline{C}_2$, respectively. Set $g_1' = \max\{g_1 - (m-1) ,0\}$ and $g_2' = \max\{g_2 - (m-1) ,0\}$. Then there exists an $L^2$ orthonormal basis of $\Gamma (C,\omega_C)$, $\{\alpha_j\}_{j=1}^g$, such that:
\begin{enumerate}[(1).]
    \item $\Vert \alpha_j \Vert^2_{L^2 ;C_2} \leq 200 e^{ - \frac{\pi^2}{\sum_{j=1}^m \ell (\gamma_j)}} $, $1\leq j\leq g_1' $.
    \item $\Vert \alpha_j \Vert^2_{L^2 ;C_1} \leq 2000 e^{ - \frac{\pi^2}{\sum_{j=1}^m \ell (\gamma_j)}} $, $g'_1 +1\leq j\leq g_1' +g_2' $.
\end{enumerate}
\end{lem}

\begin{proof}
Let $\{\mathbf{c}_j\}_{j=1}^{2(g_1 +g_2)}$ be the system of smooth loops constructed in Lemma \ref{lemtopologyseparatingmulticurvehomologicaldescripton}. For each $j$, let $\beta_j\in \mathscr{H}^1_{C,\mathbb{R}} (C) $ be the unique harmonic $1$-form whose cohomology class is the Poincar\'e dual of $[\mathbf{c}_j] \in H_{1} (C;\mathbb{R})$, and let $\beta_{\gamma_j}\in \mathscr{H}^1_{C,\mathbb{R}} (C) $ be the unique harmonic $1$-form whose cohomology class is the Poincar\'e dual of $\gamma_j $. Define the following $\mathbb{R}$-subspaces of $\mathscr{H}^1_{C,\mathbb{R}}(C)$:
$$ \mathbf{X}_1=\mathrm{span}_{\mathbb{R}} \{ \beta_j \}_{j=1}^{g_1} ,\;\; \mathbf{X}_2=\mathrm{span}_{\mathbb{R}} \{ \beta_j \}_{j=2g_1 +1}^{2g_1+g_2} ,\;\; \mathbf{X}_3=\mathrm{span}_{\mathbb{R}} \{\beta_{\gamma_j}\}_{j=1}^{m-1} .$$
Next, consider the subspaces
\begin{eqnarray*}
\mathbf{X}'_1 & = & \left( \mathbf{X}_1 + \mathbf{X}_3 \right) \bigcap \left(\bigcap\limits_{k= 1}^{m} \star \mathrm{ker} \mathbf{T}_{\gamma_k} \right) \subset \mathscr{H}^1_{C,\mathbb{R}} (C) ,\\
\mathbf{X}'_2 & = & \left( \mathbf{X}_2 + \mathbf{X}_3 \right) \bigcap \left(\bigcap\limits_{k= 1}^{m} \star \mathrm{ker} \mathbf{T}_{\gamma_k} \right) \subset \mathscr{H}^1_{C,\mathbb{R}} (C) ,
\end{eqnarray*}
where $\star$ is the Hodge $\star$-operator, and $\mathbf{T}_{\mathbf{c}} $ is the linear functional on $\mathscr{Z}_{C,\mathbb{C}}^1 (C ) $ defined by $\beta \mapsto \int_{\mathbf{c}} \beta $. Note that the Poincar\'e duality implies that $\mathbf{X}_3 \bigcap \left(\bigcap\limits_{k= 1}^{m} \star \mathrm{ker} \mathbf{T}_{\gamma_k} \right) =0$, and hence $\mathbf{X}'_1\cap \mathbf{X}'_2=0$. Since the subspace of $H_1 (C;\mathbb{R})$ generated by the homology classes $\{[\gamma_j]\}_{j=1}^m$ is isomorphic to $\mathbb{R}^{m-1} $, it follows that $\mathrm{dim}_{\mathbb{R}} \mathbf{X}'_1 \geq g_1 $ and $\mathrm{dim}_{\mathbb{R}} \mathbf{X}'_2 \geq g_2 $. Let $ \Psi_{ C_1, \mathrm{loc}} : \mathbf{X}'_1 \longrightarrow \mathscr{Z}^1_{C;\mathbb{R}} (C) $ be the linear operator constructed in Lemma \ref{lemconstructionoflocalizationoperator}. Since $\mathbf{X}_1 + \mathbf{X}_3 \subset \left( \bigcap\limits_{k=1}^{m} \mathrm{ker} \mathbf{T}_{\gamma_k} \right) \bigcap \left( \bigcap\limits_{j=2 g_1 +1}^{2(g_1 +g_2 )} \mathrm{ker} \mathbf{T}_{\mathbf{c}_j} \right) $, we have
$$ \left( \Psi_{ C_1, \mathrm{loc}} - \mathrm{id}_{\mathbf{X}'_1} \right) (\mathbf{X}'_1 ) \subset \left( \bigcap\limits_{k=1}^{m} \mathrm{ker} \mathbf{T}_{\gamma_k} \right) \bigcap \left( \bigcap\limits_{j= 1}^{2(g_1 +g_2 )} \mathrm{ker} \mathbf{T}_{\mathbf{c}_j} \right) ,$$ 
where $\mathrm{id}_{\mathbf{X}'_1} $ is the identity map of $\mathbf{X}'_1 $. It follows that
\begin{eqnarray*}
& & \frac{\left( \Psi_{ C_1, \mathrm{loc}} - \mathrm{id}_{\mathbf{X}_1} \right) (\mathbf{X}'_1 )}{\left( \Psi_{ C_1, \mathrm{loc}} - \mathrm{id}_{\mathbf{X}'_1} \right) (\mathbf{X}'_1 ) \cap d\mathscr{A}^0_{C,\mathbb{R}} (C)}\\
& \subset & \frac{\left( \bigcap\limits_{k=1}^{m} \mathrm{ker} \mathbf{T}_{\gamma_k} \right) \bigcap \left( \bigcap\limits_{j=1}^{2(g_1 +g_2 )} \mathrm{ker} \mathbf{T}_{\mathbf{c}_j} \right)}{\left( \bigcap\limits_{k=1}^{m} \mathrm{ker} \mathbf{T}_{\gamma_k} \right) \bigcap \left( \bigcap\limits_{j=1}^{2(g_1 +g_2 )} \mathrm{ker} \mathbf{T}_{\mathbf{c}_j} \right) \bigcap d\mathscr{A}^0_{C,\mathbb{R}} (C)} \cong \mathbb{R}^{m-1} ,
\end{eqnarray*}
and hence
\begin{eqnarray*}
\left( \Psi_{ C_1, \mathrm{loc}} - \mathrm{id}_{\mathbf{X}'_1} \right)^{-1} (d\mathscr{A}^0_{C,\mathbb{R}} (C) ) & \geq & \max\{ \mathrm{dim}_{\mathbb{R}} \mathbf{X}'_1 - (m-1) ,0\} \\
& \geq & \max\{ g_1 -(m-1) ,0 \} =g_1' .
\end{eqnarray*}

Set $\mathbf{X}''_1 = \left( \Psi_{ C_1, \mathrm{loc}} - \mathrm{id}_{\mathbf{X}_1} \right)^{-1} (d\mathscr{A}^0_{C,\mathbb{R}} (C) ) $. Then for any $\beta\in \mathbf{X}''_1 $, we have
$$\Vert \beta \Vert^2_{L^2 ;C } = \Vert \Psi_{ C_1, \mathrm{loc}} (\beta) \Vert^2_{L^2 ;C } \leq \Vert \beta \Vert^2_{L^2 ;C_1 } + 200 e^{ - \frac{\pi^2}{\sum_{j=1}^m \ell (\gamma_j)}} \Vert  \beta \Vert^2_{L^2 ;C } ,$$
and hence $ \Vert \beta \Vert^2_{L^2 ;C_2 } \leq 200 e^{ - \frac{\pi^2}{\sum_{j=1}^m \ell (\gamma_j)}} \Vert  \beta \Vert^2_{L^2 ;C } $. Moreover, for any pair of $1$-forms $ \beta ,\beta'\in \mathbf{X}''_1 $, we have $\int_C \beta \wedge \beta' =0 $. Let $\beta'_1 ,\cdots ,\beta'_{\mathrm{dim}\mathbf{X}''_1} $ be an $L^2$-orthonormal basis of $\mathbf{X}''_1 $, and $\alpha_j = \frac{\sqrt{2}}{2} \left( \beta'_j + i\star \beta'_j \right) $. Then we have
\begin{equation*}
\left\langle \alpha_j ,\alpha_k \right\rangle_{L^2 ;C} = \left\langle \beta'_j , \beta'_k  \right\rangle_{L^2 ;C} + \frac{i}{2} \int_C \beta'_j \wedge \beta'_k = \left\langle \beta'_j , \beta'_k  \right\rangle_{L^2 ;C} ,
\end{equation*}
and hence $\alpha_1 ,\cdots ,\alpha_{\mathrm{dim}\mathbf{X}''_1} $ forms an $L^2$-orthonormal subset of $\Gamma (C,\omega_C )$. Furthermore,
\begin{equation*}
\Vert \alpha_j \Vert^2_{L^2 ;C_2} = \left\langle \beta'_j , \beta'_j  \right\rangle_{L^2 ;C_2} + \frac{i}{2} \int_{C_2} \beta'_j \wedge \beta'_j = \left\langle \beta'_j , \beta'_j  \right\rangle_{L^2 ;C_2} \leq 200 e^{ - \frac{\pi^2}{\sum_{j=1}^m \ell (\gamma_j)}} .
\end{equation*}

Similarly, by considering the linear operator $ \Psi_{ C_2, \mathrm{loc}} : \mathbf{X}'_2 \longrightarrow \mathscr{Z}^1_{C;\mathbb{R}} (C) $ and the space
$$\mathbf{X}''_2 = \left( \Psi_{ C_2, \mathrm{loc}} - \mathrm{id}_{\mathbf{X}_2} \right)^{-1} (d\mathscr{A}^0_{C,\mathbb{R}} (C) ) ,$$
we can conclude that $\mathrm{dim}\mathbf{X}''_2\geq g'_2$, and for any $\beta\in\mathbf{X}''_2 $, $ \Vert \beta \Vert^2_{L^2 ;C_1 } \leq 200 e^{ - \frac{\pi^2}{\sum_{j=1}^m \ell (\gamma_j)}} \Vert  \beta \Vert^2_{L^2 ;C } $.

Let $\mathbf{X}''=\mathbf{X}''_1 + \mathbf{X}''_2$, and let $(\mathbf{X}''_1)^{\perp} $ be the $L^2$-orthogonal completion of $\mathbf{X}''_1$ in $\mathbf{X}''$. Then the projection $\mathbf{X}''\to(\mathbf{X}''_1)^{\perp}$ induces a vector space isomorphism $ \mathbf{X}''_2 \to (\mathbf{X}''_1)^{\perp} $.

Now, let $\beta'' \in (\mathbf{X}''_1)^{\perp} \subset \mathbf{X}'' $. Then there exists a unique $1$-form $\beta_1''\in \mathbf{X}''_1 $ such that $\beta'' +\beta_1''\in \mathbf{X}''_2 $. We then estimate:
\begin{eqnarray*}
\Vert \beta_1'' \Vert_{L^2;C}^2 & = & \left\langle \beta'' + \beta_1'' , \beta_1''   \right\rangle_{L^2 ;C_1} + \left\langle \beta'' + \beta_1'' , \beta_1''   \right\rangle_{L^2 ;C_2} \\
& \leq & 20\sqrt{2} e^{ - \frac{\pi^2}{2\sum_{j=1}^m \ell (\gamma_j)}} \Vert \beta'' + \beta''_1 \Vert_{L^2;C} \Vert \beta_1'' \Vert_{L^2;C} \\
& \leq & 20\sqrt{2} e^{ - \frac{\pi^2}{2\sum_{j=1}^m \ell (\gamma_j)}} \left( \Vert \beta'' \Vert_{L^2;C}+\Vert \beta''_1 \Vert_{L^2;C} \right) \Vert \beta_1'' \Vert_{L^2;C} ,
\end{eqnarray*}
which implies $\Vert \beta_1'' \Vert_{L^2;C} \leq \frac{20\sqrt{2} e^{ - \frac{\pi^2}{2\sum_{j=1}^m \ell (\gamma_j)}}}{1-20\sqrt{2} e^{ - \frac{\pi^2}{2\sum_{j=1}^m \ell (\gamma_j)}}} \Vert \beta'' \Vert_{L^2;C} $. 

Given that $\sum_{j=1}^m \ell(\gamma_j) \leq \frac{1}{4}$, we further obtain:
\begin{eqnarray*}
\Vert \beta'' \Vert_{L^2;C_1} & \leq & \Vert \beta'' + \beta''_1 \Vert_{L^2;C_1} + \Vert \beta_1'' \Vert_{L^2;C_1} \\
& \leq & 10\sqrt{2} e^{ - \frac{\pi^2}{2\sum_{j=1}^m \ell (\gamma_j)}} \Vert \beta'' \Vert_{L^2 ;C} + \left( 1+10\sqrt{2} e^{ - \frac{\pi^2}{2\sum_{j=1}^m \ell (\gamma_j)}} \right) \Vert \beta''_1 \Vert_{L^2;C} \\
& \leq & \frac{30\sqrt{2} e^{ - \frac{\pi^2}{2\sum_{j=1}^m \ell (\gamma_j)}}}{\left( 1-20\sqrt{2} e^{ - \frac{\pi^2}{2\sum_{j=1}^m \ell (\gamma_j)}} \right)^2 } \Vert \beta'' \Vert_{L^2;C} \\
& \leq & \frac{30\sqrt{2}e^{ - \frac{\pi^2}{2\sum_{j=1}^m \ell (\gamma_j)}}}{\left(1-20\sqrt{2}e^{-2\pi^{2}}\right)^{2}}  \Vert \beta'' \Vert_{L^2;C}\leq 20\sqrt{5} e^{ - \frac{\pi^2}{2\sum_{j=1}^m \ell (\gamma_j)}} \Vert \beta'' \Vert_{L^2;C} .
\end{eqnarray*}
Let $\beta'''_1,\cdots ,\beta'''_{\mathrm{dim}\mathbf{X}''_2}$ be an $L^2$-orthonormal basis of $(\mathbf{X}''_1)^{\perp}$. Via complexification of $(\mathbf{X}''_1)^{\perp}$, and following a construction analogous to that of the $\alpha_j$, we obtain an $L^2$-orthonormal system $\{ \alpha'_1 ,\cdots ,\alpha'_{\mathrm{dim}\mathbf{X}''_2} \} \subset \Gamma (C,\omega_C ) $, satisfying
$$\Vert \alpha'_j \Vert^2_{L^2 ;C_1} \leq 2000 e^{ - \frac{\pi^2}{\sum_{j=1}^m \ell (\gamma_j)}} ,$$
which completes the proof of the lemma.
\end{proof}

We are now ready to estimate the $\varphi$-invariant.

\begin{prop}
\label{propositionanalysispartphiinvariantcutsurfaces}
Let $C$ be a compact Riemann surface of genus $g \geq 2$, equipped with the hyperbolic metric $\mu_{\mathrm{KE}}$ of constant curvature $-1$. Suppose $\{\gamma_k \}_{k=1}^m $ is a minimal separating multicurve on $C$ satisfying $ \sum_{k=1}^m \ell (\gamma_k) \leq \frac{1}{6 \log g} $. Let $C_1$, $C_2 $ denote the two connected components of $C\setminus \left( \bigcup_{k=1}^m \gamma_k \right) $, and let $g_1$ and $g_2$ denote the genera of $\overline{C}_1$ and $\overline{C}_2$, respectively. Set $g_1' = \max\{g_1 -(m-1) ,0\}$, $g_2' = \max\{g_2 -(m-1) ,0\}$ and $g'=g_1'+g_2'$. Then the $\varphi$-invariant is bounded below by
$$\varphi(C) \geq \frac{\pi^2\cdot g_1'g_2' }{ \sum_{k=1}^m  \ell (\gamma_k) }\left(\frac{ g'}{g^2}-\frac{2 }{g } \cdot e^{15-\frac{\pi^2}{ \sum_{k=1}^m \ell (\gamma_k)}} \right) . $$
\end{prop}

\begin{proof}
By Lemma \ref{lemlocalizationpeaksectioncutsurfaces}, one can choose an $L^2$-orthonormal basis $\{\alpha_j\}_{j=1}^g$ of $\Gamma (C,\omega_C )$
\begin{itemize}
\item $\int_{\gamma_k} \alpha_j =0 $ for $j=1,\cdots ,g'$ and $k=1,\cdots ,m$.
\item $\Vert \alpha_j \Vert^2_{L^2 ;C_2} \leq 200 e^{ - \frac{\pi^2}{\sum_{k=1}^m \ell (\gamma_k)}} $ for $1\leq j\leq g_1' $.
\item $\Vert \alpha_j \Vert^2_{L^2 ;C_1} \leq 2000 e^{ - \frac{\pi^2}{\sum_{k=1}^m \ell (\gamma_k)}} $ for $g'_1 +1\leq j\leq g' $.
\end{itemize}

Let $\mathfrak{u}_{j} \in  \mathscr{A}_{C,\mathbb{R}}^0 (C) $ be the unique solution to
$$ dd^c \mathfrak{u}_{j} = i\alpha_j \wedge\bar{\alpha}_j - 2 \mu_{\mathrm{Ar}} ,\;\;\;\; \int_C \mathfrak{u}_{j} \mu_{\mathrm{Ar}} =0 . $$

For each $j$, let $\eta_j$ be a $\mathbb{R}$-valued smooth function to be specified later. Then Proposition \ref{propzhangphiinvariantpotential} implies
\begin{eqnarray*}
\varphi (C) & \geq & \frac{1}{2\pi } \sum_{j=1}^g \left\Vert d {\mathfrak{u}}_{j} \right\Vert^2_{L^2 ;C} \geq \frac{1}{2\pi} \sum_{j=1}^g \frac{\left| \left\langle d {\mathfrak{u}}_{j} ,d\eta_j \right\rangle_{L^2 ;C} \right|^2}{\left\Vert d {\eta}_{j} \right\Vert^2_{L^2 ;C}} \\
& = & \frac{\pi}{2} \sum_{j=1}^g \frac{1}{\left\Vert d {\eta}_{j} \right\Vert^2_{L^2 ;C}} \left| \int_C \eta_j \left(i\alpha_j \wedge\bar{\alpha}_j - 2 \mu_{\mathrm{Ar}}\right) \right|^2 .
\end{eqnarray*}

Our next step is to construct the functions $\eta_j$ and estimate the $\varphi$-invariant. As in the proof of Lemma \ref{lemconstructionoflocalizationoperator}, a standard mollification argument shows that it suffices to construct $C^{0,1}$ functions $\eta_{j,0}$ for the present estimate.

We make use of the local structure of the Kähler collars $\mathscr{C}(\gamma_k)$ to construct the functions $\eta_j$. By Lemma \ref{lemkahlercollar}, there exist holomorphic coordinates 
$$z_k : \mathscr{C} (\gamma_k ) \to \overline{\mathbb{D}_{e^{-\nu (\gamma_k)}} (0)} \setminus \mathbb{D}_{e^{\nu (\gamma_k ) -\frac{2\pi^2}{\ell (\gamma_k )} }} (0) $$
such that $z_k \left( \mathscr{C} (\gamma_k ) \cap C_1 \right) = \overline{\mathbb{D}_{e^{-\nu (\gamma_k)}} (0)} \setminus \mathbb{D}_{e^{-\frac{\pi^2}{\ell (\gamma_k )} }} (0) $, where $\nu (\gamma_k) \in (0, \frac{\pi}{2\ell (\gamma_k)}) $ is determined by $\nu (\gamma_k) \ell (\gamma_k) = {2\pi} \arcsin \left( \tanh \frac{\ell (\gamma_k )}{2} \right) $. Let $\eta_{\mathrm{cut}} \in C^{0,1} (\mathbb{R} ;[0,1]) $ be the cut-off function defined by:
\begin{eqnarray*}
\eta_{\mathrm{cut}} = \left\{
\begin{aligned}
 0 ,\;\;\;\;\;\;\;\;\;\;\;\;\;\;\;\; &\;\;\;  \textrm{ on $C\setminus \left( C_1 \bigcup \left(\bigcup_{k=1}^m \mathscr{C} (\gamma_k) \right) \right) $} ;\\
 \eta_{[0,1]} \left( 2+\frac{2\ell (\gamma_k)}{\pi^2} \log |z_k| \right) ,\; &\;\;\;  \textrm{ on $ \mathscr{C} (\gamma_k) $} ;\\
 1 ,\;\;\;\;\;\;\;\;\;\;\;\;\;\;\;\; &\;\;\;  \textrm{ on $ C_1 \setminus \left(\bigcup_{k=1}^m \mathscr{C} (\gamma_k) \right) $},
\end{aligned}
\right.
\end{eqnarray*}
where $\eta_{[0,1]} (t) = \max\{ 0,\min \{t,1\}\}$ denotes the clipping function. Clearly, $\eta_{\mathrm{cut}} =0$ on $C\setminus C_1$. For $1\leq j\leq g_1'$, set $\eta_j=\eta_{\mathrm{cut}}$. For such $j$, we compute
\begin{eqnarray*}
\left\Vert d {\eta}_{j } \right\Vert^2_{L^2 ;C} & = & \frac{1}{2} \sum_{k=1}^m \int_{\mathscr{C} (\gamma_k)} d {\eta}_{j } \wedge \star d {\eta}_{j } =\frac{1}{2 } \sum_{k=1}^m \int_{ e^{-\frac{\pi^2}{\ell (\gamma_k)}} }^{e^{-\frac{\pi^2}{2\ell (\gamma_k)}} } 2\pi \cdot \frac{4 \ell^2 (\gamma_k)}{\pi^4 } \frac{ d\varsigma_k }{ \varsigma_k } \\
& = & \sum_{k=1}^m \frac{4 \ell^2 (\gamma_k)}{\pi^3 } \left(  \frac{\pi^2}{\ell (\gamma_k)} - \frac{\pi^2}{2\ell (\gamma_k)} \right) = \frac{2}{\pi} \sum_{k=1}^m  \ell (\gamma_k) ,
\end{eqnarray*}
where $z_k=\varsigma_k e^{i\theta_k} $ denotes polar coordinate. Moreover, for $1\leq j\leq g_1'$, it follows from Proposition \ref{propexplicitestimateskahlercollar} and Lemma \ref{lemlocalizationpeaksectioncutsurfaces} that
\begin{eqnarray*}
i\int_C {\eta}_{j } \alpha_j \wedge \bar{\alpha}_j & \geq & 2\Vert\alpha_j \Vert_{L^2 ;C_1}^2 - \sum_{k=1}^m \Vert\alpha_j \Vert_{L^2 ;\mathscr{C}(\gamma_k)}^2 \int_{\left\{ -1 \leq \frac{ \ell (\gamma_k) \log |z_k |}{\pi^2} \leq -\frac{1}{2} \right\} } \frac{4e^{2\pi}}{\pi} \mu_{\mathrm{Euc}} \\
& \geq & 2\Vert\alpha_j \Vert_{L^2 ;C_1}^2 - 4e^{2\pi} e^{-\frac{\pi^2}{ \sum_{k=1}^m \ell (\gamma_k)}} \Vert\alpha_j \Vert_{L^2 ;C }^2 \\
& \geq & 2-\left(400 + 4e^{2\pi} \right) e^{-\frac{\pi^2}{ \sum_{k=1}^m \ell (\gamma_k)}} \geq 2-e^{10-\frac{\pi^2}{ \sum_{k=1}^m \ell (\gamma_k)}}.
\end{eqnarray*}

Let $\mu'_{\mathrm{Ar}} = \frac{i}{2g} \sum_{j=1}^{g'} \alpha_j \wedge \bar{\alpha}_j $. Then $0\leq \mu'_{\mathrm{Ar}} \leq \mu_{\mathrm{Ar}} $, and by orthonormality we have $\mu'_{\mathrm{Ar}}(C) = \frac{g'}{g} $. By Lemma \ref{lemlocalizationpeaksectioncutsurfaces},
\begin{eqnarray*}
\left| \mu'_{\mathrm{Ar}} (C_1) - \frac{g_1'}{g} \right| & \leq & \frac{1}{g} \sum_{j=1}^{g_1'} \left| \Vert\alpha_j \Vert_{L^2 ;C_1}^2 -1 \right| + \frac{1}{g} \sum_{j=g_1'+1}^{g' } \Vert\alpha_j \Vert_{L^2 ;C_1}^2 \\
& \leq & \frac{2000g'}{g} e^{ -\frac{\pi^2}{ \sum_{k=1}^m \ell (\gamma_k)}} ,
\end{eqnarray*}
and hence
\begin{eqnarray*}
\left| \mu'_{\mathrm{Ar}} (C_2) - \frac{g_2'}{g} \right| & = & \left| \frac{g'}{g} - \mu'_{\mathrm{Ar}} (C_1) - \frac{g' -g_1'}{g} \right|\\
& = & \left| \mu'_{\mathrm{Ar}} (C_1) - \frac{g_1'}{g} \right| \leq \frac{2000g'}{g} e^{ -\frac{\pi^2}{ \sum_{k=1}^m \ell (\gamma_k)}} .
\end{eqnarray*}

Suppose that $g_2'\ge 1$. For any $1\leq j\leq g_1'$ we have
\begin{eqnarray*}
 & & \frac{1}{\left\Vert d {\eta}_{j } \right\Vert^2_{L^2 ;C}} \left| \int_C \eta_{j } \left(i\alpha_j \wedge\bar{\alpha}_j - 2 \mu_{\mathrm{Ar}}\right) \right|^2 \\
 & \geq & \frac{\pi}{2 \sum_{k=1}^m  \ell (\gamma_k) } \cdot \left( \left| \int_C \eta_j \alpha_j \wedge\bar{\alpha}_j \right|- 2 \left| \int_{C_1} 1\mu_{\mathrm{Ar}} \right| \right)^2 \\
 & \geq & \frac{\pi}{2 \sum_{k=1}^m  \ell (\gamma_k) } \cdot \left( 2-e^{10-\frac{\pi^2}{ \sum_{k=1}^m \ell (\gamma_k)}}-2(\mu_{\mathrm{Ar}} (C)-\mu_{\mathrm{Ar}}' (C ))- 2\mu_{\mathrm{Ar}}' (C_1) \right)^2 \\
 & = & \frac{\pi}{2 \sum_{k=1}^m  \ell (\gamma_k) } \cdot \left( 2\mu_{\mathrm{Ar}}' (C_2) -e^{10-\frac{\pi^2}{ \sum_{k=1}^m \ell (\gamma_k)}} \right)^2 \\
 & \geq & \frac{\pi}{2 \sum_{k=1}^m  \ell (\gamma_k) } \cdot \left( \frac{2g_2' }{g } -e^{10-\frac{\pi^2}{ \sum_{k=1}^m \ell (\gamma_k)}}- \frac{4000g'}{g} e^{ -\frac{\pi^2}{ \sum_{k=1}^m \ell (\gamma_k)}} \right)^2 .
\end{eqnarray*}

Since $e^{30-\frac{\pi^2}{ \sum_{k=1}^m \ell (\gamma_k)}} \leq e^{30 } g^{-6\pi^2} < \frac{1}{100g} $, we obtain
\begin{eqnarray*}
 & & \frac{1}{\left\Vert d {\eta}_{j } \right\Vert^2_{L^2 ;C}} \left| \int_C \eta_{j } \left(i\alpha_j \wedge\bar{\alpha}_j - 2 \mu_{\mathrm{Ar}}\right) \right|^2 \\
 & \geq & \frac{\pi}{2 \sum_{k=1}^m  \ell (\gamma_k) } \cdot \left( \frac{2g_2' }{g } -e^{10-\frac{\pi^2}{ \sum_{k=1}^m \ell (\gamma_k)}}- \frac{4000g'}{g} e^{ -\frac{\pi^2}{ \sum_{k=1}^m \ell (\gamma_k)}} \right)^2 \\
 & \geq & \frac{2\pi}{ \sum_{k=1}^m  \ell (\gamma_k) } \cdot \left( \frac{g_2'^2}{g^2}-\frac{( e^{10 } + 4000)g_2' }{g }  e^{-\frac{\pi^2}{ \sum_{k=1}^m \ell (\gamma_k)}} -e^{25-\frac{2\pi^2}{ \sum_{k=1}^m \ell (\gamma_k)}}  \right) \\
 & \geq & \frac{2\pi}{ \sum_{k=1}^m \ell (\gamma_k) } \cdot \left( \frac{g_2'^2}{g^2}-\frac{g_2' }{g } \cdot e^{15-\frac{\pi^2}{ \sum_{k=1}^m \ell (\gamma_k)}} \right) .
\end{eqnarray*}
Hence
\begin{eqnarray*}
\frac{1}{2\pi } \sum_{j=1}^{g_1'} \left\Vert d {\mathfrak{u}}_{j} \right\Vert^2_{L^2 ;C} & \geq & \frac{\pi}{2} \sum_{j=1}^{g_1'} \frac{1}{\left\Vert d {\eta}_{j } \right\Vert^2_{L^2 ;C}} \left| \int_C \eta_{j } \left(i\alpha_j \wedge\bar{\alpha}_j - 2 \mu_{\mathrm{Ar}}\right) \right|^2 \\
& \geq & \frac{\pi^2 \cdot g_1'}{ \sum_{k=1}^m  \ell (\gamma_k) } \cdot \left( \frac{g_2'^2}{g^2}-\frac{g_2' }{g } \cdot e^{15-\frac{\pi^2}{ \sum_{k=1}^m \ell (\gamma_k)}} \right) .
\end{eqnarray*}
If instead $g_2'=0$, the right-hand side vanishes and the inequality is trivial.

For any $g_1' +1\leq j\leq g'$, let $\eta_j=1-\eta_{\mathrm{cut}}$. By \ref{propexplicitestimateskahlercollar}, for any $g_1' +1\leq j\leq g'$, we obtain
\begin{eqnarray*}
\left|\int_{C_1} {\eta}_{j} \mu'_{\mathrm{Ar}} \right| & \leq & \frac{1}{2g} \sum_{j=1}^{g'} \sum_{k=1}^m \left|\int_{ \mathscr{C} (\gamma_k) \cap C_1 } (1-{\eta}_{\mathrm{cut}}) \alpha_j \wedge \bar{\alpha}_j \right| \\
& \leq & \frac{1}{2g} \sum_{j=1}^{g'} \sum_{k=1}^m \Vert\alpha_j \Vert_{L^2 ;\mathscr{C}(\gamma_k)}^2 \int_{\left\{ -1 \leq \frac{ \ell (\gamma_k) \log |z_k |}{\pi^2} \leq -\frac{1}{2} \right\} } \frac{4e^{2\pi}}{\pi} \mu_{\mathrm{Euc}} \\
& \leq & \frac{1}{2g} \sum_{j=1}^{g'} \sum_{k=1}^m \Vert\alpha_j \Vert_{L^2 ;\mathscr{C}(\gamma_k)}^2 \cdot 4e^{2\pi} \cdot e^{-\frac{\pi^2}{\ell (\gamma_k )}} \\
& \leq & \frac{ g'}{g} \cdot 2e^{2\pi} e^{ -\frac{\pi^2}{ \sum_{k=1}^m \ell (\gamma_k)}} \leq \frac{ g'}{g} \cdot e^{ 8 -\frac{\pi^2}{ \sum_{k=1}^m \ell (\gamma_k)}} .
\end{eqnarray*}
Arguing as above, for $g_1' +1\leq j\leq g'$ we deduce
\begin{eqnarray*}
\frac{1}{2\pi } \sum_{j=g'_1+1}^{g'} \left\Vert d {\mathfrak{u}}_{j} \right\Vert^2_{L^2 ;C} & \geq & \frac{\pi^2 \cdot g_2'}{ \sum_{k=1}^m  \ell (\gamma_k) } \cdot \left( \frac{g_1'^2}{g^2}-\frac{g_1' }{g } \cdot e^{15-\frac{\pi^2}{ \sum_{k=1}^m \ell (\gamma_k)}} \right) .
\end{eqnarray*}
Therefore,
\begin{eqnarray*}
\varphi (C) & \geq & \frac{1}{2\pi } \sum_{j=1}^{g'} \left\Vert d {\mathfrak{u}}_{j} \right\Vert^2_{L^2 ;C} \geq \frac{1}{2\pi} \sum_{j=1}^{g'} \frac{\left| \left\langle d {\mathfrak{u}}_{j} ,d\eta_j \right\rangle_{L^2 ;C} \right|^2}{\left\Vert d {\eta}_{j} \right\Vert^2_{L^2 ;C}} \\
& \geq & \frac{\pi^2 }{ \sum_{k=1}^m  \ell (\gamma_k) }\left(\frac{g_1'g_2'g'}{g^2}-\frac{2g_1'g_2' }{g } \cdot e^{15-\frac{\pi^2}{ \sum_{k=1}^m \ell (\gamma_k)}} \right) .
\end{eqnarray*}
This completes the proof.
\end{proof}

\subsection{Non-separating case}
\label{sectionzhangphi3}

As above, let $C$ be a curve of genus $g \geq 2$ over $\mathbb{C}$, and let $\omega_C$ denote its canonical bundle. Let $\mu_{\mathrm{KE}} $ be the unique K\"ahler metric on $C $ with constant curvature $-1 $.

Our goal in this subsection is to estimate the $\varphi$-invariant in the case where $C$ contains a short non-separating closed geodesic. A loop on $C$ is called non-separating if its homology class is nonzero. Our approach combines Tian's peak section method with the fundamental relationship between loops and harmonic $1$-forms developed by Buser--Sarnak \cite{busa1}.

We next consider the relationship between the $\varphi$-invariant and the peak section along a short closed geodesic.

\begin{lem}
\label{lemW12integralpeaksectionalongshortnonseparatedgeodesic}
Let $C$ be a compact Riemann surface of genus $g \geq 2$, equipped with the hyperbolic metric $\mu_{\mathrm{KE}}$ of constant curvature $-1$. Let $\gamma $ be a non-separating closed geodesic on $(C,\mu_{\mathrm{KE}})$ with length $\ell (\gamma) =\ell <2\arcsinh 1 $, and let $\alpha_1 ,\cdots ,\alpha_g $ be an $L^2 $-orthonormal basis of $\Gamma (C,\omega_C)$ such that $\alpha_1$ is the peak section in $\Gamma (C,\omega_C )$ along $\gamma$. Assume that $|\int_\gamma \alpha_{1} | =2\pi \mathbf{t}_{\gamma } $. Let $\mathfrak{u}_{j} \in  \mathscr{A}_{C,\mathbb{R}}^0 (C) $ satisfy $\int_C \mathfrak{u}_{j} \mu_{\mathrm{Ar}} =0$ and $dd^c \mathfrak{u}_{j} = i\alpha_j \wedge\bar{\alpha}_j - 2 \mu_{\mathrm{Ar}} $. Then
$$ \sum_{j=1}^g \int_{\mathscr{C} (\gamma)} d\mathfrak{u}_{j} \wedge \star d\mathfrak{u}_{j} \geq \frac{32\pi^2 (g-1) \mathbf{t}_{\gamma }^2}{ g } \left( \frac{\pi^2}{\ell} -\nu \right) \left( \frac{2\pi \mathbf{t}_{\gamma }^2}{3} \left( \frac{\pi^2}{\ell} -\nu \right)^2 -1 \right) , $$
where $\mathscr{C} (\gamma) $ denotes the collar of $\gamma$ as in Theorem \ref{thmcollarthm}. 
\end{lem}

\begin{proof}
By Lemma \ref{lemkahlercollar}, we have
$$(\mathscr{C} (\gamma) ,\mu_{\mathrm{KE}} |_{\mathscr{C} (\gamma)} ) \cong \left( \overline{\mathbb{D}_{e^{-\nu }} (0)} \setminus {\mathbb{D}}_{e^{ \nu  - \frac{2\pi^2}{\ell } }} (0) , \frac{ \ell^2 \mu_{\mathrm{Euc}} }{ 4\pi^2 |z|^2 \sin^2 ( \frac{\ell }{2\pi} \log |z| ) } \right) ,$$ 
where $\mu_{\mathrm{Euc}} = \frac{idz\wedge d\bar{z}}{2} $ denotes the standard Euclidean metric on $ \mathbb{C} $, and $\nu \in (0, \frac{\pi}{2\ell }) $ satisfies $\nu \ell = {2\pi} \arcsin \left( \tanh \frac{\ell }{2} \right) $.

Define $\mathfrak{u}_{\gamma } (z) = \frac{1}{2\pi} \int_{0}^{2\pi} \mathfrak{u}_{1} ( ze^{i\theta} ) d\theta $, $\forall z\in \overline{\mathbb{D}_{e^{-\nu }} (0)} \setminus {\mathbb{D}}_{e^{ \nu - \frac{2\pi^2}{\ell } }} (0) $. By a similar argument as in the proof of Lemma \ref{lemintegralgradientintegralsmalldisc}, we obtain
$$ \int_{\mathscr{C} (\gamma)} d\mathfrak{u}_{1} \wedge \star d\mathfrak{u}_{1} \geq \int_{\mathscr{C} (\gamma)} d\mathfrak{u}_{\gamma} \wedge \star d\mathfrak{u}_{\gamma} .$$

Write $\alpha_j = \sum_{k\in\mathbb{Z} } a_k \mathbf{t}_{j,k } z^k dz $, where $a_k\in\mathbb{C} $ satisfy
\begin{eqnarray*}
|a_k |^2 = \left\{
\begin{aligned}
\frac{1}{4\pi} \left( \frac{\pi^2}{\ell } -\nu \right)^{-1} ,\;\;\;\;\;\;\;\;\;\;\;\;\;\;\;\; & \; k=-1 ;\\
 \frac{k+1}{\pi} \left( e^{-2(k+1) \nu } - e^{2(k+1) (\nu - \frac{2\pi^2}{\ell } ) } \right)^{-1}, \; & \; k\neq -1 . 
\end{aligned}
\right.
\end{eqnarray*}
By Lemma \ref{lemlocalorthonormalbasis}, the set $\{ a_k z^k \}_{k\in\mathbb{Z}}$ forms a local $L^2$-orthonormal basis on $\mathscr{C} (\gamma)$. Hence, we have $\sum_{k\in\mathbb{Z} } | \mathbf{t}_{j,k } |^2 \leq 1 $, $j=1,\cdots ,g $. Since $|\int_\gamma \alpha_{1} | =2\pi \mathbf{t}_{\gamma } $, it follows that $a_{-1} | \mathbf{t}_{1,-1 } | = \mathbf{t}_{\gamma } $.

By definition, on $\mathscr{C} (\gamma)\cong \overline{\mathbb{D}_{e^{-\nu }} (0)} \setminus {\mathbb{D}}_{e^{ \nu  - \frac{2\pi^2}{\ell } }} (0) $, we have
\begin{equation*}
dd^c \mathfrak{u}_{1} = i\alpha_1 \wedge\bar{\alpha}_1 - 2 \mu_{\mathrm{Ar}} = \frac{(g-1)i}{g} \alpha_1 \wedge\bar{\alpha}_1 - \frac{i}{g}\sum_{j=2}^g \alpha_j \wedge\bar{\alpha}_j .
\end{equation*}
Hence
\begin{eqnarray*}
dd^c \mathfrak{u}_{\gamma} & = & \left( \frac{(g-1)i}{g} \sum_{k\in\mathbb{Z} } |a_k|^2 |\mathbf{t}_{1,k }|^2 |z|^{2k} - \frac{i}{g} \sum_{j=2}^g \sum_{k\in\mathbb{Z} } |a_k|^2 |\mathbf{t}_{j,k }|^2 |z|^{2k} \right) dz \wedge d\bar{z} .
\end{eqnarray*}

We express $z$ in polar coordinates by setting $z=\varsigma e^{i\theta}$, where $\varsigma = |z| $ and $\theta\in [0,2\pi)$. Then $\star d\varsigma = \varsigma d\theta $, and $\star d\theta = -\frac{1}{\varsigma}d\varsigma $. Hence
$$ dd^c \mathfrak{u}_{\gamma} =  \frac{i}{4\pi} \left( \frac{\partial^2 \mathfrak{u}_{\gamma}}{\partial\varsigma^2 } + \frac{1}{\varsigma} \frac{\partial\mathfrak{u}_{\gamma}}{\partial\varsigma} \right) dz\wedge d\bar{z} .$$
Observing that $\mathbf{t}_{j,-1} = 0$ for $j \geq 2$, we obtain
\begin{equation*}
\frac{1}{4\pi\varsigma } \frac{\partial}{\partial\varsigma} \left( \varsigma \frac{\partial\mathfrak{u}_{\gamma}}{\partial\varsigma} \right)  = \frac{1}{4\pi} \left( \frac{\partial^2 \mathfrak{u}_{\gamma}}{\partial\varsigma^2 } + \frac{1}{\varsigma} \frac{\partial\mathfrak{u}_{\gamma}}{\partial\varsigma} \right) = \frac{g-1}{g} \mathbf{t}_{\gamma }^2 \varsigma^{-2} + \Psi_0 (\varsigma ) ,
\end{equation*}
where $\Psi_0 (\varsigma ) = \sum\limits_{k\neq -1 } |a_k|^2\left( \frac{g-1}{g} |\mathbf{t}_{1,k }|^2 - \frac{1}{g} \sum\limits_{j=2}^g  |\mathbf{t}_{j,k }|^2 \right)\varsigma^{2k} $. Integrating, we have
\begin{eqnarray*}
\frac{\varsigma}{4\pi} \frac{\partial\mathfrak{u}_{\gamma}}{\partial\varsigma}  & = & \int_{e^{ -\frac{\pi^2}{\ell}}}^{\varsigma} \left( \frac{g-1}{g} \mathbf{t}_{\gamma }^2 \varsigma^{-1} + \varsigma \Psi_0 (\varsigma ) \right) d\varsigma + \frac{e^{ -\frac{\pi^2}{\ell}}}{4\pi} \frac{\partial\mathfrak{u}_{\gamma}}{\partial\varsigma} (e^{ -\frac{\pi^2}{\ell}}) \\
& = & \frac{g-1}{g} \mathbf{t}_{\gamma }^2 \left( \log \varsigma +\frac{\pi^2}{\ell} \right) + \Psi_1 (\varsigma ) - \Psi_1 (e^{\frac{-\pi^2}{\ell}} ) + \frac{e^{ -\frac{\pi^2}{\ell}}}{4\pi} \frac{\partial\mathfrak{u}_{\gamma}}{\partial\varsigma} (e^{ -\frac{\pi^2}{\ell}}) ,
\end{eqnarray*}
where 
$$\Psi_1 (\varsigma ) = \sum_{k\neq -1} \frac{|a_k|^2}{2k+2}\left( \frac{g-1}{g} |\mathbf{t}_{1,k }|^2 - \frac{1}{g} \sum\limits_{j=2}^g  |\mathbf{t}_{j,k }|^2 \right) \varsigma^{2k+2} .$$

Set $ \vartheta_0 =  \frac{e^{ -\frac{\pi^2}{\ell}}}{4\pi} \frac{\partial\mathfrak{u}_{\gamma}}{\partial\varsigma} (e^{ -\frac{\pi^2}{\ell}}) - \Psi_1 (e^{\frac{-\pi^2}{\ell}} ) $. Then 
$$ \frac{\varsigma}{4\pi} \frac{\partial\mathfrak{u}_{\gamma}}{\partial\varsigma} = \frac{g-1}{g} \mathbf{t}_{\gamma }^2 \left( \log \varsigma +\frac{\pi^2}{\ell} \right)  + \Psi_1 (\varsigma ) + \vartheta_0 ,$$
and consequently,
\begin{eqnarray*}
 \int_{\mathscr{C} (\gamma)} d\mathfrak{u}_{\gamma} \wedge \star d\mathfrak{u}_{\gamma} & = & \int_{\mathscr{C} (\gamma)} \frac{16\pi^2 (g-1)^2 \mathbf{t}_{\gamma }^4 }{g^2 \varsigma^2 } \left( \log \varsigma +\frac{\pi^2}{\ell} \right)^2 \mu_{\mathrm{Euc}} \\
& & + \int_{\mathscr{C} (\gamma)} \frac{16\pi^2}{\varsigma^2 } \left( \Psi_1 (\varsigma ) + \vartheta_0 \right)^2 \mu_{\mathrm{Euc}} \\
& & + \int_{\mathscr{C} (\gamma)} \frac{32\pi^2 (g-1)}{g \varsigma^2 } \mathbf{t}_{\gamma }^2 \left( \log \varsigma +\frac{\pi^2}{\ell} \right) \left( \Psi_1 (\varsigma ) + \vartheta_0 \right) \mu_{\mathrm{Euc}} .
\end{eqnarray*}

By a straightforward calculation, we have 
$$ \int_{e^{\nu -\frac{2\pi^2}{\ell}}}^{e^{-\nu}} \frac{1}{ \varsigma } \left( \log \varsigma +\frac{\pi^2}{\ell} \right) d\varsigma =0 ,\;\; \int_{e^{\nu -\frac{2\pi^2}{\ell}}}^{e^{-\nu}} \frac{1}{ \varsigma } \left( \log \varsigma +\frac{\pi^2}{\ell} \right)^2 d\varsigma =\frac{2}{3} \left( -\nu +\frac{\pi^2}{\ell} \right)^3 . $$

It follows that
\begin{eqnarray*}
\int_{\mathscr{C} (\gamma)} d\mathfrak{u}_{\gamma} \wedge \star d\mathfrak{u}_{\gamma} & \geq & \frac{32\pi^3 (g-1)^2 \mathbf{t}_{\gamma }^4 }{g^2 } \int_{e^{\nu -\frac{2\pi^2}{\ell}}}^{e^{-\nu}} \frac{1 }{ \varsigma } \left( \log \varsigma +\frac{\pi^2}{\ell} \right)^2 d\varsigma \\
& & + \frac{64\pi^3 (g-1) \mathbf{t}_{\gamma }^2}{g } \int_{e^{\nu -\frac{2\pi^2}{\ell}}}^{e^{-\nu}} \frac{\Psi_1 (\varsigma) +\vartheta_0 }{ \varsigma } \left( \log \varsigma +\frac{\pi^2}{\ell} \right) d\varsigma \\
& \geq & \frac{64\pi^3 (g-1)^2 \mathbf{t}_{\gamma }^4 }{3g^2 } \left( -\nu +\frac{\pi^2}{\ell} \right)^3 \\
& & + \frac{64\pi^3 (g-1) \mathbf{t}_{\gamma }^2}{g } \int_{e^{\nu -\frac{2\pi^2}{\ell}}}^{e^{-\nu}} \frac{\Psi_1 (\varsigma) }{ \varsigma } \left( \log \varsigma +\frac{\pi^2}{\ell} \right) d\varsigma .
\end{eqnarray*}

Since $\Psi_1 (\varsigma ) = \sum_{k\neq -1} \frac{|a_k|^2}{2k+2}\left( \frac{g-1}{g} |\mathbf{t}_{1,k }|^2 - \frac{1}{g} \sum\limits_{j=2}^g  |\mathbf{t}_{j,k }|^2 \right) \varsigma^{2k+2} $, we obtain
\begin{eqnarray*}
& & \int_{e^{\nu -\frac{2\pi^2}{\ell}}}^{e^{-\nu}} \frac{\Psi_1 (\varsigma) }{ \varsigma } \left( \log \varsigma +\frac{\pi^2}{\ell} \right) d\varsigma \\
& = & \sum_{k\neq -1} \frac{|a_k|^2}{2k+2}\left( \frac{g-1}{g} |\mathbf{t}_{1,k }|^2 - \frac{1}{g} \sum\limits_{j=2}^g  |\mathbf{t}_{j,k }|^2 \right) \int_{e^{\nu -\frac{2\pi^2}{\ell}}}^{e^{-\nu}} \varsigma^{2k+1} \left( \log \varsigma +\frac{\pi^2}{\ell} \right) d\varsigma \\
& = & \sum_{k\neq -1} \frac{\frac{\pi^2}{\ell} -\nu }{4\pi (k+1)}\left( \frac{g-1}{g} |\mathbf{t}_{1,k }|^2 - \frac{1}{g} \sum\limits_{j=2}^g  |\mathbf{t}_{j,k }|^2 \right) \frac{ 1 + e^{4(k+1)(\nu - \frac{\pi^2 }{\ell} )} }{1 - e^{4(k+1)(\nu - \frac{\pi^2 }{\ell} )}} \\
& & - \sum_{k\neq -1} \frac{1 }{8\pi (k+1)^2}\left( \frac{g-1}{g} |\mathbf{t}_{1,k }|^2 - \frac{1}{g} \sum\limits_{j=2}^g  |\mathbf{t}_{j,k }|^2 \right) .
\end{eqnarray*}

Since $\sum_{k\in\mathbb{Z} } | \mathbf{t}_{j,k } |^2 \leq 1 $, $j=1,\cdots ,g $, we have
$$ \sum_{k\neq -1} \left| \frac{g-1}{g} |\mathbf{t}_{1,k }|^2 - \frac{1}{g} \sum\limits_{j=2}^g  |\mathbf{t}_{j,k }|^2 \right| \leq \frac{g-1}{g} + \frac{1}{g} \sum\limits_{j=2}^g  1 = \frac{2(g-1)}{g} .$$
It follows that
\begin{eqnarray*}
& & \left| \int_{e^{\nu -\frac{2\pi^2}{\ell}}}^{e^{-\nu}} \frac{\Psi_1 (\varsigma) }{ \varsigma } \left( \log \varsigma +\frac{\pi^2}{\ell} \right) d\varsigma \right| \\
& = & \left| \sum_{k\neq -1} \frac{\frac{\pi^2}{\ell} -\nu }{4\pi (k+1)}\left( \frac{g-1}{g} |\mathbf{t}_{1,k }|^2 - \frac{1}{g} \sum\limits_{j=2}^g  |\mathbf{t}_{j,k }|^2 \right) \frac{ 1 + e^{4(k+1)(\nu - \frac{\pi^2 }{\ell} )} }{1 - e^{4(k+1)(\nu - \frac{\pi^2 }{\ell} )}} \right. \\
& & \left. - \sum_{k\neq -1} \frac{1 }{8\pi (k+1)^2}\left( \frac{g-1}{g} |\mathbf{t}_{1,k }|^2 - \frac{1}{g} \sum\limits_{j=2}^g  |\mathbf{t}_{j,k }|^2 \right) \right| \\
& \leq & \frac{(g-1)\left( \frac{\pi^2}{\ell} -\nu \right)}{2\pi g} \cdot \sup_{k\neq -1} \left| \frac{ 1 + e^{4(k+1)(\nu - \frac{\pi^2 }{\ell} )} }{(k+1) \left(1 - e^{4(k+1)(\nu - \frac{\pi^2 }{\ell} )}\right)} - \frac{1}{2(k+1)^2\left( \frac{\pi^2}{\ell} -\nu \right)} \right| \\
& = & \frac{(g-1)\left( \frac{\pi^2}{\ell} -\nu \right)}{2\pi g} \cdot \sup_{k\geq 0} \left( \left| \frac{1}{k+1} \right|\left| 1+ \frac{ 2 }{  e^{4(k+1)( \frac{\pi^2 }{\ell} -\nu )} -1 } - \frac{1}{2(k+1) \left( \frac{\pi^2}{\ell} -\nu \right)} \right| \right)\\
& \leq & \frac{(g-1)\left( \frac{\pi^2}{\ell} -\nu \right)}{2\pi g} .
\end{eqnarray*}

Therefore, we obtain the lower bound
\begin{eqnarray*}
\int_{\mathscr{C} (\gamma)} d\mathfrak{u}_{\gamma} \wedge \star d\mathfrak{u}_{\gamma} & \geq & \frac{64\pi^3 (g-1)^2 \mathbf{t}_{\gamma }^4 }{3g^2 } \left( -\nu +\frac{\pi^2}{\ell} \right)^3 \\
& & + \frac{64\pi^3 (g-1) \mathbf{t}_{\gamma }^2}{g } \int_{e^{\nu -\frac{2\pi^2}{\ell}}}^{e^{-\nu}} \frac{\Psi_1 (\varsigma) }{ \varsigma } \left( \log \varsigma +\frac{\pi^2}{\ell} \right) d\varsigma \\
& \geq & \frac{32\pi^2 (g-1)^2 \mathbf{t}_{\gamma }^2}{ g^2 } \left( \frac{\pi^2}{\ell} -\nu \right) \left( \frac{2\pi \mathbf{t}_{\gamma }^2}{3} \left( \frac{\pi^2}{\ell} -\nu \right)^2 -1 \right) .
\end{eqnarray*}

By definition, $\sum_{j=1}^g \mathfrak{u}_{j} = 0 $, and using an argument analogous to that in Lemma \ref{lemphiinvariantthickpart}, we conclude that
\begin{eqnarray*}
\sum_{j=1}^g \int_{\mathscr{C} (\gamma)} d\mathfrak{u}_{j} \wedge \star d\mathfrak{u}_{j} & \geq & \frac{g}{g-1} \int_{\mathscr{C} (\gamma)} d\mathfrak{u}_{1} \wedge \star d\mathfrak{u}_{1} .
\end{eqnarray*}
This completes the proof.
\end{proof}

We now consider the shortest non-separating simple closed geodesic $\gamma_0$, together with the shortest loop in the free loop space whose homology class has nonzero intersection number with $[\gamma_0]$.

\begin{lem}
\label{lemshortgeodesicintersectionanothergeodesicloop}
Let $C$ be a compact Riemann surface of genus $g \geq 2$, equipped with the hyperbolic metric $\mu_{\mathrm{KE}}$ of constant curvature $-1$. Let $\gamma_0 $ be the shortest non-separating simple closed geodesic on $(C,\mu_{\mathrm{KE}})$, and define
$$ \mathcal{S}_{\gamma_0} = \left\{ \mathbf{c}: \mathbb{S}^1 = \mathbb{R}/\mathbb{Z} \to C \; \big| \; \mathbf{c}\; \textrm{ piecewise smooth, and } \mathbf{c}_* (e_{ \mathbb{S}^1}) \cap [\gamma_0] \neq 0 \right\} ,$$
where $e_{ \mathbb{S}^1} $ denotes the generator of $H_1 (\mathbb{S}^1;\mathbb{Z})\cong\mathbb{Z}$. Assume that $\ell (\gamma_0) \leq 2\arcsinh 1 $. Then there exists a closed geodesic $\mathbf{c}_0 \in \mathcal{S}_{\gamma_0} $ satisfying the following properties:

\begin{enumerate}[(1)]
    \item $\ell (\mathbf{c}_0 ) \leq \ell (\mathbf{c} ) $ for all $ \mathbf{c} \in \mathcal{S}_{\gamma_0} $.
    \item For any $t_0\in (0,1) $ and $t\in\mathbb{R}$, we have
    $$ \mathrm{dist}_{\mu_{\mathrm{KE}}} ( \mathbf{c}_0 (t ) , \mathbf{c}_0 (t + t_0 ) ) = \min \{ \ell (\mathbf{c}_0) t_0 ,\ell (\mathbf{c}_0) (1 -t_0) \} .$$
    \item For any simple closed geodesic $\gamma' $ on $C$ with $\ell (\gamma')\leq 2\arcsinh 1$, we have $\# (\gamma' \cap \mathbf{c}_0 ) \leq 1 $.
\end{enumerate}
\end{lem}

\begin{proof}
The argument is standard. Since the injectivity radius of any compact surface is strictly positive, every sufficiently short loop must be homotopically trivial. Hence the infimum of lengths among loops in $\mathcal{S}_{\gamma_0}$ is strictly positive. Consider a minimizing sequence of loops in $\mathcal{S}_{\gamma_0}$ whose lengths converge to this infimum. By compactness of the space of loops with bounded length, a subsequence converges to a closed geodesic loop realizing the minimum. This limit loop is the desired $\mathbf{c}_0$, and the lemma follows.
\end{proof}

To construct a suitable curve near the given geodesic loop, we invoke the segment inequality of Cheeger--Colding, specialized to the setting of hyperbolic surfaces.

\begin{thm}[Cheeger--Colding’s segment inequality (hyperbolic surface version)]
\label{thmcheegercoldingsegmentinequality}
Let $C$ be a Riemann surface equipped with a complete hyperbolic metric $\mu_{\mathrm{KE}}$ of constant curvature $-1$. Let $A,B,W\subset C$ be Borel measurable subsets, and let $u:C\to [0,\infty )$ be a nonnegative Borel measurable function. For any pair $(x,y)\in A\times B$, choose a minimal geodesic $\mathbf{c}_{x,y}:[0,1]\to C$ joining $x$ and $y$. Assume that the geodesics lie entirely within $W$, i.e., $\mathbf{c}_{x,y}\subset W$. Let $d_{W}=\mathrm{diam}_{\mu_{\mathrm{KE}}} (W) $ denote the diameter of $W$. Then
\begin{eqnarray*}
\int_{A\times B} \int_0^1 u(\mathbf{c}_{x,y} (t)) dt\; \mu_{\mathrm{KE}} (x)\wedge \mu_{\mathrm{KE}} (y) & \leq & \cosh \left( \frac{d_{W}}{2} \right) \left( \mu_{\mathrm{KE}} (A) \cdot \sup_{x\in A,t\in[\frac{1}{2},1]}\int_{\mathbf{c}_{x,B} (t)} u\mu_{\mathrm{KE}} \right. \\
& & \quad\quad\quad\quad\quad\quad\quad\quad+ \mu_{\mathrm{KE}} (B) \cdot \sup_{y\in B,t\in[0,\frac{1}{2}]}\int_{\mathbf{c}_{A,y} (t)} u\mu_{\mathrm{KE}} ) \\
& \leq & \cosh \left( \frac{d_{W}}{2} \right) ( \mu_{\mathrm{KE}} (A) + \mu_{\mathrm{KE}} (B) ) \int_{W} u\mu_{\mathrm{KE}},
\end{eqnarray*}
where $\mathbf{c}_{x,B} (t) = \{\mathbf{c}_{x,y} (t) :y\in B \textrm{ and the minimal geodesic from $x$ to $y$ is unique}\} $, and $\mathbf{c}_{A,y}(t)$ is defined analogously.
\end{thm}

\begin{remark}
For any point $x\in C$, the set of points in $C$ that are connected to $x$ by more than one minimal geodesic is contained in the ``cut locus" of $x$, and hence has measure zero \cite{jahe1, itta1}. Therefore, the integral appearing in Cheeger--Colding's segment inequality is unchanged if one restricts to pairs of points connected by a unique minimal geodesic.
\end{remark}

\begin{proof}
By the Ricci comparison (see \cite[Lemma 7.1.2]{pp1}), for any $x\in A$ and $t\in (0,1)$,
$$\int_{B} u(\mathbf{c}_{x,y} (t)) \mu_{\mathrm{KE}} (y) = \int_{\mathbf{c}_{x,B} (1)} u(\mathbf{c}_{x,y} (t)) \mu_{\mathrm{KE}} (y) \leq \frac{\sinh d_{W}}{\sinh (t d_{W} ) } \int_{\mathbf{c}_{x,B} (t)} u\mu_{\mathrm{KE}} .$$
Note that the measure $\mu_{\mathrm{KE}}(B\setminus \mathbf{c}_{x,B}(1))=0$. It follows that
$$\int_A \int_B \int_{\frac{1}{2}}^1 u(\mathbf{c}_{x,y} (t)) dt\wedge \mu_{\mathrm{KE}} (y) \wedge \mu_{\mathrm{KE}} (x) \leq \frac{\sinh d_{W}}{2\sinh\frac{d_{W}}{2}} \mu_{\mathrm{KE}} (A) \cdot \sup_{x\in A,t\in[\frac{1}{2},1]}\int_{\mathbf{c}_{x,B} (t)} u\mu_{\mathrm{KE}} .$$ 

Similarly,
$$\int_B \int_A \int_0^{\frac{1}{2}} u(\mathbf{c}_{x,y} (t)) dt\wedge\mu_{\mathrm{KE}} (x)\wedge\mu_{\mathrm{KE}} (y) \leq \frac{\sinh d_{W}}{2\sinh\frac{d_{W}}{2}} \mu_{\mathrm{KE}} (B) \cdot \sup_{y\in B,t\in[0,\frac{1}{2}]}\int_{\mathbf{c}_{A,y} (t)} u\mu_{\mathrm{KE}} .$$ 
Combining the two estimates yields the stated inequalities. For more details, see \cite[Theorem 2.11]{chco1} or \cite[Theorem 7.1.10]{pp1}.
\end{proof}

To proceed, we investigate the behavior of peak sections associated with short geodesics intersecting the loop $\mathbf{c}_0$ constructed in Lemma \ref{lemshortgeodesicintersectionanothergeodesicloop}. In particular, we apply Cheeger–Colding’s segment inequality to derive lower bounds on their periods along the corresponding geodesics.

\begin{lem}
\label{lempeaksectionalongshortgeodesicintersectionanothergeodesicloop}
Let $C$ be a compact Riemann surface of genus $g \geq 2$, equipped with the hyperbolic metric $\mu_{\mathrm{KE}}$ of constant curvature $-1$. Let $\gamma_0 $ be the shortest non-separating simple closed geodesic on $(C,\mu_{\mathrm{KE}})$, and let $\mathbf{c}_0 $ be the closed geodesic defined in Lemma \ref{lemshortgeodesicintersectionanothergeodesicloop}. Assume that $\ell (\gamma_0) \leq 2\arcsinh \frac{1}{8} $.

Let $\gamma_0 ,\gamma_1 ,\cdots ,\gamma_m$ be the simple closed geodesics with $\ell (\gamma_j)\leq 2\arcsinh \frac{1}{2}$ that intersect $\mathbf{c}_0$. For each $l=1,\dots,m$, let $\alpha_{\gamma_l}$ denote the peak section in $\Gamma (C,\omega_C)$ along $\gamma_l$. Let $J_0 \subset \{0, \dots, m\}$ denote the set of indices $j$ such that $\gamma_j$ is non-separating.

Then for any $\gamma_l$ with $\ell (\gamma_l) \leq 2\arcsinh \frac{1}{3} $, we have
$$ \left| \int_{\gamma_l} \alpha_{\gamma_l} \right| \geq \frac{\sqrt{\pi}}{240 \sqrt{g-1} + \left( \sum_{j\in J_0} \left| \frac{ \pi^2}{\ell (\gamma_j)} - \nu (\gamma_j) \right| \right)^{\frac{1}{2}}} .$$
\end{lem}

\begin{proof}
Let $\beta_{\gamma_l}$ be the harmonic form representing the Poincar\'e dual of $\gamma_l $. Then the peak section in $\Gamma (C,\omega_C)$ along $\gamma_l$ is given by $\alpha_{\gamma_l} = \frac{- \sqrt{2}}{2 \Vert \beta_{\gamma_l} \Vert_{L^2;C}} (\star \beta_{\gamma_l} -i\beta_{\gamma_l} ) $, and hence
\begin{eqnarray*}
 \int_{\gamma_l} \alpha_{\gamma_l} & = & \int_C \alpha_{\gamma_l} \wedge \beta_{\gamma_l} = \int_C \frac{- \sqrt{2}}{2 \Vert \beta_{\gamma_l} \Vert_{L^2;C}} (\star \beta_{\gamma_l} -i\beta_{\gamma_l} ) \wedge \beta_{\gamma_l} = \sqrt{2}\Vert \beta_{\gamma_l} \Vert_{L^2;C} .
\end{eqnarray*}
We now proceed to estimate $\Vert \beta_{\gamma_l} \Vert_{L^2;C}$.

Let $\rho_{\mathcal{H}}:\mathcal{H}\to C$ be a universal covering map, let $\tilde{\mathbf{c}}_1 :\mathbb{R}\to \mathcal{H}$ be a unit-speed parametrization lifting $\mathbf{c}_0 $, and set $\mathbf{c}_1=\rho_{\mathcal{H}}\circ\tilde{\mathbf{c}}_1$. For convenience, we write $\gamma_{j+m+1}=\gamma_{j}$. After reparameterizing $\tilde{\mathbf{c}}_1$, we may choose constants $\{\mathbf{r}_j\}_{j\in\mathbb{Z}}$ satisfying $\mathbf{r}_0 =0$, $\mathbf{r}_{2j} < \mathbf{r}_{2j+1}$, $\mathbf{r}_{2j-1} \leq \mathbf{r}_{2j}$, $\mathbf{r}_{2m+2+j}=\ell (\mathbf{c}_0) +\mathbf{r}_j$, and $\mathbf{c}_1 ([\mathbf{r}_{2j} , \mathbf{r}_{2j+1}])$ are the connected components of $\mathbf{c}_0 \cap \mathscr{C} (\gamma_j) $ containing $\mathbf{c}_0 \cap \gamma_j $, where $\mathscr{C} (\gamma_j) $ is the collar of $\gamma_j$. By Lemma \ref{lemkahlercollar}, there exist holomorphic coordinates 
$$z_j : \mathscr{C} (\gamma_j ) \to \overline{\mathbb{D}_{e^{-\nu (\gamma_j)}} (0)} \setminus \mathbb{D}_{e^{\nu (\gamma_j ) -\frac{2\pi^2}{\ell (\gamma_j )} }} (0) $$
such that $|z_j(\mathbf{c}_1 (\mathbf{r}))| $ is strictly increasing on $[\mathbf{r}_{2j} , \mathbf{r}_{2j+1}] $, where $\nu (\gamma_j) \in (0, \frac{\pi}{2\ell (\gamma_j)}) $ satisfies $\nu (\gamma_j) \ell (\gamma_j) = {2\pi} \arcsin \left( \tanh \frac{\ell (\gamma_j )}{2} \right) $. See also Corollary \ref{coranothergeodesicscrossorawayfromshortgeodesics}.

We now choose points $\mathbf{r}'_{2j} , \mathbf{r}'_{2j+1} \in ( \mathbf{r}_{2j} , \mathbf{r}_{2j+1} ) $ such that 
$ |z_j (\mathbf{r}'_{2j}) | = e^{2\nu (\gamma_j ) -\frac{2\pi^2}{\ell (\gamma_j)} } $ and $|z_j (\mathbf{r}'_{2j+1}) | = e^{-2\nu (\gamma_j ) } $. Let $\mathbf{r}\in\mathbb{R}$. By the above construction and Corollary \ref{coranothergeodesicscrossorawayfromshortgeodesics}, we conclude that, if $\mathbf{r}\notin (\mathbf{r}'_{2j} , \mathbf{r}'_{2j+1}) $ for all $j\in\mathbb{Z}$, then $ \mathrm{inj}_{\mu_{\mathrm{KE}}} (\mathbf{c}_1 (\mathbf{r})) \geq \arcsinh \frac{1}{2}$.

For any $j\in\mathbb{Z}$, choose $N_j +1$ points $\mathbf{r}'_{2j-1} =\mathbf{r}''_{j ,0} <\mathbf{r}''_{j ,1}<\cdots < \mathbf{r}''_{j ,N_j} = \mathbf{r}'_{2j } $ so that the interval $[ \mathbf{r}'_{2j-1} , \mathbf{r}'_{2j } ]$ is divided into $N_j$ segments of equal length, satisfying $\mathbf{c}_1 (\mathbf{r}''_{j ,k}) = \mathbf{c}_1 (\mathbf{r}''_{j +m+1 ,k}) $ and $\frac{\mathbf{r}''_{j ,k+1} - \mathbf{r}''_{j ,k}}{\arcsinh \frac{1}{2}} \in [1,\frac{3}{2}] $. This is possible since Proposition \ref{propbasicpropertiesriemanniankahlercollars} implies that $ |\mathbf{r}'_{2j }-\mathbf{r}'_{2j-1 }| \geq 2\log 2 > 2\arcsinh \frac{1}{2} $.

Let $\mathbf{r}'''_{j ,k} =\frac{\mathbf{r}''_{j ,k} + \mathbf{r}''_{j ,k+1} }{2}$ be the midpoint of $[\mathbf{r}''_{j ,k} , \mathbf{r}''_{j ,k+1}]$. Set
$$\widetilde{B}_{j,k} = \left\{ \tau\in \mathcal{H}: \mathrm{dist}_{\mu_{\mathcal{H}}} (\tau,\tilde{\mathbf{c}}_1 (\mathbf{r}''_{j ,k})) <\arcsinh \frac{1}{2} \right\} ,$$
and 
$$\widetilde{W}_{j,k} = \left\{ \tau\in \mathcal{H}: \mathrm{dist}_{\mu_{\mathcal{H}}} \left(\tau,\tilde{\mathbf{c}}_1 \left(\mathbf{r}'''_{j ,k}\right)\right) <\arcsinh \frac{1}{2} + \frac{\mathbf{r}''_{j ,k+1} - \mathbf{r}''_{j ,k} }{2} \right\} ,$$
where $\mu_{\mathcal{H}}=\rho_{\mathcal{H}}^*\mu_{\mathrm{KE}}$ is the pullback of $\mu_{\mathrm{KE}}$. Write $B_{j,k}=\rho_{\mathcal{H}} (\widetilde{B}_{j,k})$ and $W_{j,k}=\rho_{\mathcal{H}} (\widetilde{W}_{j,k})$. Since $ \mathrm{inj}_{\mu_{\mathrm{KE}}} (\mathbf{c}_1 (\mathbf{r}''_{j,k})) \geq \arcsinh \frac{1}{2}$, $\rho_{\mathcal{H}}|_{\widetilde{B}_{j,k}} : \widetilde{B}_{j,k} \to B_{j,k}$ is a diffeomorphism. By a straightforward calculation,
$$\mu_{\mathrm{KE}} (B_{j,k} ) = \mu_{\mathcal{H}} (\widetilde{B}_{j,k}) = 2\pi \int_{0}^{\arcsinh\frac{1}{2}} \sinh t dt = \pi (\sqrt{5} -2) ,$$
and
$$ \mathrm{diam}_{\mu_{\mathrm{KE}}} (W_{j,k} ) = \mathrm{diam}_{\mu_{\mathcal{H}}} (\widetilde{W}_{j,k} ) = 2\arcsinh \frac{1}{2} + \mathbf{r}''_{j ,k+1} - \mathbf{r}''_{j ,k} \leq \frac{7}{2} \arcsinh \frac{1}{2} .$$
For each $0\leq k<N_j$ and any pair $(\tilde{x}_{j,k} ,\tilde{x}_{j,k+1}) \in \widetilde{B}_{j,k} \times \widetilde{B}_{j,k+1} $, there exists a unique geodesic segment $\tilde{\mathbf{c}}_{\tilde{x}_{j,k} ,\tilde{x}_{j,k+1}} : [0,1]\to \mathcal{H}$ joining $\tilde{x}_{j,k} $ and $\tilde{x}_{j,k+1}$, which lies entirely within $\widetilde{W}_{j,k}$. Moreover, by the hyperbolic triangle formula (see \cite[Theorem 2.2.2]{bu1}), for any $t\in (0,1)$,
\begin{eqnarray*}
\mathrm{dist}_{\mu_{\mathcal{H}}} ( \tilde{\mathbf{c}}_{\tilde{x}_{j,k} ,\tilde{x}_{j,k+1}} (t) ,\tilde{\mathbf{c}}_{\tilde{\mathbf{c}}_1 (\mathbf{r}''_{j ,k}) ,\tilde{x}_{j,k+1}} (t) ) & < & \arcsinh \frac{1}{2},\\ 
\mathrm{dist}_{\mu_{\mathcal{H}}} ( \tilde{\mathbf{c}}_{\tilde{x}_{j,k} ,\tilde{x}_{j,k+1}} (t) , \tilde{\mathbf{c}}_{\tilde{x}_{j,k} ,\tilde{\mathbf{c}}_1 (\mathbf{r}''_{j ,k+1})} (t) ) & < & \arcsinh \frac{1}{2} ,\\
\mathrm{dist}_{\mu_{\mathcal{H}}} ( \tilde{\mathbf{c}}_1 (t\mathbf{r}''_{j ,k+1}+(1-t)\mathbf{r}''_{j ,k}) ,\tilde{\mathbf{c}}_{\tilde{\mathbf{c}}_1 (\mathbf{r}''_{j ,k}) ,\tilde{x}_{j,k+1}} (t) ) & < & \arcsinh \frac{1}{2},\\
\mathrm{dist}_{\mu_{\mathcal{H}}} ( \tilde{\mathbf{c}}_1 (t\mathbf{r}''_{j ,k+1}+(1-t)\mathbf{r}''_{j ,k}) , \tilde{\mathbf{c}}_{\tilde{x}_{j,k} ,\tilde{\mathbf{c}}_1 (\mathbf{r}''_{j ,k+1})} (t) ) & < & \arcsinh \frac{1}{2} .
\end{eqnarray*}

We now apply Cheeger--Colding's segment inequality to the balls $\widetilde{B}_{j,k} $. 

Write $x_{j,k}=\rho_{\mathcal{H}}(\tilde{x}_{j,k})$ and $\mathbf{c}_{x_{j,k} ,x_{j,k+1}}=\rho_{\mathcal{H}}\circ\tilde{\mathbf{c}}_{\tilde{x}_{j,k} ,\tilde{x}_{j,k+1}} $. By definition, we have
\begin{eqnarray*}
\left| \int_{\mathbf{c}_{x_{j,k} ,x_{j,k+1}}} \beta_{\gamma_l} \right| & \leq & \frac{\sqrt{2}}{2} \ell (\mathbf{c}_{x_{j,k} ,x_{j,k+1}}) \int_{0}^{1} \Vert \beta_{\gamma_l} \Vert_{\mathrm{hyp}} (\mathbf{c}_{x_{j,k} ,x_{j,k+1}} (t) ) dt \\
& \leq & \frac{7 \sqrt{2}\cdot \arcsinh \frac{1}{2} }{4} \int_{0}^{1} \Vert \beta_{\gamma_l} \Vert_{\mathrm{hyp}} (\mathbf{c}_{x_{j,k} ,x_{j,k+1}} (t) ) dt ,
\end{eqnarray*}
where the norm $\Vert \cdot \Vert_{\mathrm{hyp}}$ is defined by $\beta \wedge \star \beta = \frac{\Vert \beta \Vert^2_{\mathrm{hyp}}}{2} \mu_{\mathrm{KE}} $. By the above argument, for any $t\in [0,1]$, $\mathrm{inj}_{\mu_{\mathrm{KE}}} (\mathbf{c}_1 (t\mathbf{r}''_{j ,k+1}+(1-t)\mathbf{r}''_{j ,k})) \geq \arcsinh \frac{1}{2} $. Hence, by combining Corollary \ref{coranothergeodesicscrossorawayfromshortgeodesics} with Theorem \ref{thmcheegercoldingsegmentinequality}, we obtain
\begin{eqnarray*}
& & \int_{B_{j,k} \times B_{j,k+1}} \int_{0}^{1} \Vert \beta_{\gamma_l} \Vert_{\mathrm{hyp}} (\mathbf{c}_{x_{j,k} ,x_{j,k+1}} (t) ) dt\wedge\mu_{\mathrm{KE}} (x_{j,k}) \wedge\mu_{\mathrm{KE}} (x_{j,k+1}) \\
& = & \int_{\widetilde{B}_{j,k} \times \widetilde{B}_{j,k+1}} \int_{0}^{1} \Vert\beta_{\gamma_l} \Vert_{\mathrm{hyp}} ( \rho_{\mathcal{H}}(\tilde{\mathbf{c}}_{\tilde{x}_{j,k} ,\tilde{x}_{j,k+1}} (t) )) dt\wedge\mu_{\mathcal{H}} (\tilde{x}_{j,k}) \wedge\mu_{\mathcal{H}} (\tilde{x}_{j,k+1}) \\
& \leq & \cosh \left( \frac{7}{4} \arcsinh \frac{1}{2} \right) \left( \mu_{\mathcal{H}} (\widetilde{B}_{j,k}) \cdot \sup_{\tilde{x}_{j,k}\in \widetilde{B}_{j,k},t\in[\frac{1}{2},1]}\int_{\tilde{\mathbf{c}}_{\tilde{x}_{j,k},\widetilde{B}_{j,k+1}} (t)} (\Vert\beta_{\gamma_l} \Vert_{\mathrm{hyp}} \circ \rho_{\mathcal{H}})\mu_{\mathcal{H}} \right. \\
& & \left.\quad\quad\quad\quad\quad\quad\quad\quad+ \mu_{\mathcal{H}} (\widetilde{B}_{j,k+1}) \cdot \sup_{\tilde{x}_{j,k+1}\in \widetilde{B}_{j,k+1},t\in[0,\frac{1}{2}]}\int_{\tilde{\mathbf{c}}_{\widetilde{B}_{j,k},\tilde{x}_{j,k+1}} (t)} (\Vert\beta_{\gamma_l} \Vert_{\mathrm{hyp}} \circ \rho_{\mathcal{H}})\mu_{\mathcal{H}} \right) \\
& \leq & 2\pi (\sqrt{5}-2) \cosh \left( \frac{7}{4} \arcsinh \frac{1}{2} \right) \left(  \sup_{\tilde{x}_{j,k}\in \widetilde{B}_{j,k},t\in[\frac{1}{2},1]}\int_{\rho_{\mathcal{H}}(\tilde{\mathbf{c}}_{\tilde{x}_{j,k},\widetilde{B}_{j,k+1}} (t))} \Vert\beta_{\gamma_l} \Vert_{\mathrm{hyp}} \mu_{\mathrm{KE}} \right. \\
& & \left.\quad\quad\quad\quad\quad\quad\quad\quad+  \sup_{\tilde{x}_{j,k+1}\in \widetilde{B}_{j,k+1},t\in[0,\frac{1}{2}]}\int_{\rho_{\mathcal{H}}(\tilde{\mathbf{c}}_{\widetilde{B}_{j,k},\tilde{x}_{j,k+1})} (t)} \Vert\beta_{\gamma_l} \Vert_{\mathrm{hyp}} \mu_{\mathrm{KE}} \right) \\
& \leq & 4\pi (\sqrt{5}-2) \cosh \left( \frac{7}{4} \arcsinh \frac{1}{2} \right) \int_{W_{j,k} } \Vert \beta_{\gamma_l} \Vert_{\mathrm{hyp}}  \; \mu_{\mathrm{KE}},
\end{eqnarray*}
where $\tilde{\mathbf{c}}_{\widetilde{B}_{j,k},\tilde{x}_{j,k+1}} (t)$ and $\tilde{\mathbf{c}}_{\tilde{x}_{j,k},\widetilde{B}_{j,k+1}} (t)$ are defined as in Theorem  \ref{thmcheegercoldingsegmentinequality}.

It follows that
\begin{eqnarray*}
& & \fint_{\prod\limits_{k=0}^{N_j} B_{j,k} } \left| \sum_{k=0}^{N_{j} -1} \int_{\mathbf{c}_{x_{j,k} ,x_{j,k+1}}} \beta_{\gamma_l} \right|\; \bigwedge_{k=0}^{N_j} \mu_{\mathrm{KE}}(x_{j,k}) \\
& \leq & \frac{7 \sqrt{2}\cdot \arcsinh \frac{1}{2} }{4}  \fint_{\prod\limits_{k=0}^{N_j} B_{j,k} }\sum_{k=0}^{N_{j} -1} \int_{0}^{1} \Vert \beta_{\gamma_l} \Vert_{\mathrm{hyp}} (\mathbf{c}_{x_{j,k} ,x_{j,k+1}} (t) ) dt\wedge \bigwedge_{k=0}^{N_j} \mu_{\mathrm{KE}}(x_{j,k}) \\
& = & \sum_{k=0}^{N_{j} -1} \frac{7 \sqrt{2}\cdot \arcsinh \frac{1}{2} }{4 \pi^2 (\sqrt{5} -2)^2} \cdot \int_{B_{j,k} \times B_{j,k+1}} \int_{0}^{1} \Vert \beta_{\gamma_l} \Vert_{\mathrm{hyp}} (\mathbf{c}_{x_{j,k} ,x_{j,k+1}} (t) ) dt\wedge\mu_{\mathrm{KE}} (x_{j,k})\wedge\mu_{\mathrm{KE}} (x_{j,k+1})\\
& \leq & \frac{7 \sqrt{2}\cdot \arcsinh \frac{1}{2} \cosh \left( \frac{7}{4} \arcsinh \frac{1}{2}\right)}{\pi (\sqrt{5} -2) } \cdot \sum_{k=0}^{N_{j} -1} \int_{W_{j,k} } \Vert \beta_{\gamma_l} \Vert_{\mathrm{hyp}}  \; \mu_{\mathrm{KE}} ,
\end{eqnarray*}
where $\fint$ denotes the integral average. Hence, for any $\epsilon >0$, there exists a choice of points $x_{j,k}\in B_{j,k} $ such that
\begin{eqnarray*}
\left| \sum_{k=0}^{N_{j} -1} \int_{\mathbf{c}_{x_{j,k} ,x_{j,k+1}}} \beta_{\gamma_l} \right| & \leq & \epsilon + \frac{7 \sqrt{2}\cdot \arcsinh \frac{1}{2} \cosh \left( \frac{7}{4} \arcsinh \frac{1}{2}\right)}{\pi (\sqrt{5} -2) } \cdot \sum_{k=0}^{N_{j} -1} \int_{W_{j,k} } \Vert \beta_{\gamma_l} \Vert_{\mathrm{hyp}}  \; d\mu_{\mathrm{KE}} .
\end{eqnarray*}

Since $|z_j (\mathbf{r}'_{2j+1})| = e^{-2\nu (\gamma_j )} $, it follows from Proposition \ref{propbasicpropertiesriemanniankahlercollars} that
\begin{eqnarray*}
    \mathrm{dist}_{\mu_{\mathrm{KE}}} \left( \mathbf{c}_1 (\mathbf{r}''_{j+1 ,0}) , \left\{ x\in\mathscr{C} (\gamma_j) : \left| \frac{\log |z_{j} (x)|}{\nu (\gamma_j)} + \sqrt{5} \right| =1 \right\} \right)  & > &  \arcsinh \frac{1}{2} .
\end{eqnarray*}
Therefore,
$$ B_{j+1,0} \subset \left\{ x\in\mathscr{C} (\gamma_j) : \left| \frac{\log |z_{j} (x)|}{\nu (\gamma_j)} + \sqrt{5} \right| < 1 \right\} .$$
Similarly,
$$ B_{j,N_j} \subset \left\{ x\in\mathscr{C} (\gamma_j) : \left| \frac{\log |z_{j} (x)| + \frac{2\pi^2}{\ell (\gamma_j)}}{\nu (\gamma_j)} - \sqrt{5} \right| < 1 \right\} .$$

By Theorem \ref{thmcollarthm}, $\gamma_j\cap \gamma_l =\varnothing$, and thus $\int_{\gamma_j} \beta_{\gamma_l} =0$. It follows that 
$$ \alpha_{\gamma_l}'= \beta_{\gamma_l} + i\star \beta_{\gamma_{l}} - \left( \frac{1}{2\pi} \int_{\gamma_j} \star \beta_{\gamma_l} \right) \frac{dz_j}{z_j} $$
defines a holomorphic $1$-form on $\mathscr{C} (\gamma_j)$ satisfying $\int_{\gamma_j} \alpha_{\gamma_l}' =0 $. Applying Proposition \ref{propexplicitestimateskahlercollar}, we see that for any pair $(x_{j,N_j} ,x_{j+1 ,0})\in B_{j,N_j}\times B_{j+1,0} $, there exists a piecewise smooth curve $\mathbf{c}_{x_{j,N_j} ,x_{j+1 ,0}}$ joining $ x_{j,N_j} $ and $x_{j+1 ,0} $, such that
\begin{eqnarray*}
 \left| \int_{\mathbf{c}_{x_{j,N_j} ,x_{j+1 ,0}}} \beta_{\gamma_l} \right| & = & \left| \int_{\mathbf{c}_{x_{j,N_j} ,x_{j+1 ,0}}} \re(\beta_{\gamma_l} + i\star \beta_{\gamma_{l}}) \right| \\
  & \leq & \left(  \frac{-4\log \left(  1-e^{ -2(\sqrt{5} -2) \nu (\gamma_j ) } \right) }{\pi (1-e^{-3\pi} ) \left( 1+e^{-\frac{2\pi^2}{\ell (\gamma_j )} +2(\sqrt{5}+1)\nu (\gamma_j) } \right)^{-2}} \right)^{\frac{1}{2}} \Vert \beta_{\gamma_l} + i\star \beta_{\gamma_{l}} \Vert_{L^2 ;\mathscr{C} (\gamma_j)} \\
  & & + \frac{1}{2\pi}\left| \int_{\gamma_j} \star \beta_{\gamma_{l}} \right| \left| \frac{2\pi^2}{\ell (\gamma_j)} - 2(\sqrt{5} -1) \nu (\gamma_j) \right| ,
\end{eqnarray*}
where $\nu (\gamma_j ) \in (0, \frac{\pi}{2\ell (\gamma_j )}) $ satisfies $\nu (\gamma_j) \ell (\gamma_j) = 2\pi \arcsin \left( \tanh \frac{\ell (\gamma_j )}{2} \right) $. Moreover, Proposition \ref{propbasicpropertiesriemanniankahlercollars} implies that $\nu (\gamma_j) \geq \frac{ \pi \arcsin \frac{\sqrt{5}}{5} }{ \arcsinh \frac{1}{2}} \approx 2.83\cdots >\frac{5}{2} $. Consequently,
\begin{eqnarray*}
  & & \frac{-2\log \left(  1-e^{ -2(\sqrt{5} -2) \nu (\gamma_j ) } \right) }{\pi (1-e^{-3\pi} ) \left( 1+e^{-\frac{2\pi^2}{\ell (\gamma_j )} +2(\sqrt{5}+1)\nu (\gamma_j) } \right)^{-2}} \\ 
  & < & \frac{-2\log\left( 1-e^{ -5(\sqrt{5} -2) } \right)}{\pi (1-e^{-3\pi} ) \left( 1+e^{-\frac{\pi^2}{\arcsinh \left( \frac{1}{2} \right)} +2(\sqrt{5}+1)\pi } \right)^{-2} } \approx 0.789\cdots <1.
\end{eqnarray*}
Thus 
\begin{eqnarray*}
 \left| \int_{\mathbf{c}_{x_{j,N_j} ,x_{j+1 ,0}}} \beta_{\gamma_l} \right| & < & \Vert \beta_{\gamma_l} + i\star \beta_{\gamma_{l}} \Vert_{L^2 ;\mathscr{C} (\gamma_j)}  + \frac{1}{ \pi}\left| \int_{\gamma_j} \star \beta_{\gamma_{l}} \right| \left| \frac{ \pi^2}{\ell (\gamma_j)} - \nu (\gamma_j) \right| \\
 & \leq & 2 \Vert \beta_{\gamma_l} \Vert_{L^2 ;\mathscr{C} (\gamma_j)}  + \frac{1}{ \pi}\left| \int_{\gamma_j} \star \beta_{\gamma_{l}} \right| \left| \frac{ \pi^2}{\ell (\gamma_j)} - \nu (\gamma_j) \right| .
\end{eqnarray*}

Set $\mathbf{c}_{\gamma_l}$ to be the piecewise smooth loop on $C$ given by
$$ \mathbf{c}_{\gamma_l} = \left( \bigcup_{j=0}^m \bigcup_{k=0}^{N_j -1} \mathbf{c}_{x_{j,k} ,x_{j,k+1}} \right) \bigcup \left( \bigcup_{j=0}^m \mathbf{c}_{x_{j,N_j} ,x_{j+1 ,0}} \right) .$$

Our next step is to show that $|[\mathbf{c}_{\gamma_l}] \cap [\gamma_l]| \geq 1$. By construction, the difference of the homology classes satisfies
$$[\mathbf{c}_{\gamma_l}]-[\mathbf{c}_{0}] \in \sum_{j=0}^m \mathbb{Z} [\gamma_j] \subset H_1 (C;\mathbb{Z}) ,$$
where $[\gamma_j]\in H_1 (C;\mathbb{Z})$ denotes the homology class represented by $\gamma_j$. By Theorem \ref{thmcollarthm}, the short closed geodesics $\gamma_j$ are pairwise disjoint, and hence
$$ [\mathbf{c}_{\gamma_l}] \cap [\gamma_l]=[\mathbf{c}_{0}] \cap [\gamma_l] \neq 0 .$$
Since $[\mathbf{c}_{\gamma_l}] \cap [\gamma_l]\in\mathbb{Z}$, it follows that $|[\mathbf{c}_{\gamma_l}] \cap [\gamma_l]| \geq 1$. Consequently, $\left|\int_{\mathbf{c}_{\gamma_l}} \beta_{\gamma_l} \right| = |[\mathbf{c}_{\gamma_l}] \cap [\gamma_l]| \geq 1 $.
See also \cite[(6.31)]{bot1}. Hence, for any $\epsilon >0$,
\begin{eqnarray*}
 1 & \leq & \sum_{j=0}^m \left| \sum_{k=0}^{N_{j} -1} \int_{\mathbf{c}_{x_{j,k} ,x_{j,k+1}}} \beta_{\gamma_l} \right| + \sum_{j=0}^m \left| \int_{\mathbf{c}_{x_{j,N_j} ,x_{j+1 ,0}}} \beta_{\gamma_l} \right| \\
 & \leq & \epsilon + \frac{7 \sqrt{2}\cdot \arcsinh \frac{1}{2} \cosh \left( \frac{7}{4} \arcsinh \frac{1}{2}\right)}{\pi (\sqrt{5} -2) } \cdot \sum_{j=0}^m \sum_{k=0}^{N_{j} -1} \int_{W_{j,k} } \Vert \beta_{\gamma_l} \Vert_{\mathrm{hyp}} \mu_{\mathrm{KE}} \\
 &  & +  2 \sum_{j=0}^m \Vert \beta_{\gamma_l} \Vert_{L^2 ;\mathscr{C} (\gamma_j)}  + \frac{1}{ \pi}\sum_{j=0}^m \left| \int_{\gamma_j} \star \beta_{\gamma_{l}} \right| \left| \frac{ \pi^2}{\ell (\gamma_j)} - \nu (\gamma_j) \right| .
\end{eqnarray*}
Letting $\epsilon\to 0$, the first term vanishes. We next estimate the remaining three terms.

Let $y\in C$. We now estimate $\# \{ W_{j,k} : y\in W_{j,k} \} $. If $y\in W_{j,k}\cap W_{j',k'}$, then
\begin{eqnarray*}
\mathrm{dist}_{\mu_{\mathrm{KE}}} (\mathbf{c}_1 (\mathbf{r}'''_{j ,k}) , \mathbf{c}_1 (\mathbf{r}'''_{j' ,k'})) & \leq & \mathrm{dist}_{\mu_{\mathrm{KE}}} (y , \mathbf{c}_1 (\mathbf{r}'''_{j ,k})) + \mathrm{dist}_{\mu_{\mathrm{KE}}} (y , \mathbf{c}_1 (\mathbf{r}'''_{j' ,k'})) \\
& < & \frac{  |\mathbf{r}''_{j ,0} -\mathbf{r}''_{j ,1}| + |\mathbf{r}''_{j' ,0} -\mathbf{r}''_{j' ,1}|}{2} + 2\arcsinh \frac{1}{2} .
\end{eqnarray*}
By construction, we have $\frac{  |\mathbf{r}''_{j ,0} -\mathbf{r}''_{j ,1}| + |\mathbf{r}''_{j' ,0} -\mathbf{r}''_{j' ,1}|}{2} \leq \frac{3}{2} \arcsinh \frac{1}{2} $ and 
$$\ell (\mathbf{c}_0 ) \geq 2\arcsinh 8 \approx 5.55294 \cdots > \frac{21}{2} \arcsinh \frac{1}{2} \approx 5.05272 \cdots . $$
Then Lemma \ref{lemshortgeodesicintersectionanothergeodesicloop} implies that there exists an integer $n$ such that
\begin{eqnarray*}
|\mathbf{r}'''_{j ,k} -\mathbf{r}'''_{j' ,k'} -n\ell (\mathbf{c}_0 )| & = & \mathrm{dist}_{\mu_{\mathrm{KE}}} (\mathbf{c}_1 (\mathbf{r}'''_{j ,k}) , \mathbf{c}_1 (\mathbf{r}'''_{j' ,k'})) \\
& < & \frac{  |\mathbf{r}''_{j ,0} -\mathbf{r}''_{j ,1}| + |\mathbf{r}''_{j' ,0} -\mathbf{r}''_{j' ,1}|}{2} + 2\arcsinh \frac{1}{2} \\
& \leq & \frac{7}{2}\arcsinh \frac{1}{2} .
\end{eqnarray*}
Without loss of generality, assume that there exists $W_{j_0 ,k_0} $ containing $y$, such that for any $W_{j,k}$ containing $y$, we have
$$|\mathbf{r}'''_{j ,k} -\mathbf{r}'''_{j_0 ,k_0} | < \frac{  |\mathbf{r}''_{j ,0} -\mathbf{r}''_{j ,1}| + |\mathbf{r}''_{j_0 ,0} -\mathbf{r}''_{j_0 ,1}|}{2} + 2\arcsinh \frac{1}{2} .$$

For any pair $W_{j_0 ,k_1} , W_{j_0 ,k_2} \in \{ W_{j,k} : y\in W_{j,k} \} $, we have
$$|\mathbf{r}'''_{j_0 ,k_1} -\mathbf{r}'''_{j_0 ,k_2} | = |k_1 -k_2| |\mathbf{r}''_{j_0 ,1} -\mathbf{r}''_{j_0 ,0} | \leq 2\arcsinh \frac{1}{2} + |\mathbf{r}''_{j_0 ,0} -\mathbf{r}''_{j_0 ,1}| ,$$
and hence $|k_1 -k_2| \leq 2$. It follows that $\# \{ W_{j_0 ,k} : y\in W_{j_0 ,k} \} \leq 3 $.

Suppose that there exist $ W_{j_0 ,k_1 } , W_{j_1 ,k_2 } \in \{ W_{j,k} : y\in W_{j,k} \}  $ with $m+1 \nmid j_1 - j_0$. Without loss of generality we may assume $j_0+m\geq j_1\geq j_0 +1$. By Proposition \ref{propbasicpropertiesriemanniankahlercollars} we obtain
\begin{eqnarray*}
|\mathbf{r}'''_{j_0 ,k_1} -\mathbf{r}'''_{j_1 ,k_2} | & \geq & |\mathbf{r}''_{j_0 ,N_{j_0}} -\mathbf{r}''_{j_0 +1 ,0} | + \frac{  |\mathbf{r}''_{j_0+1 ,0} -\mathbf{r}''_{j_0+1 ,1}| + |\mathbf{r}''_{j_0 ,N_{j_0}} -\mathbf{r}''_{j_0 ,N_{j_0}-1}|}{2} \\
& \geq & 2\arcsinh 2-\log \left( \frac{\tan (\frac{\ell (\gamma_{j_0})\nu (\gamma_{j_0})}{2\pi})}{\tan (\frac{\ell (\gamma_{j_0})\nu (\gamma_{j_0})}{4\pi})} \right)-\log \left( \frac{\tan (\frac{\ell (\gamma_{j_0+1})\nu (\gamma_{j_0+1})}{2\pi})}{\tan (\frac{\ell (\gamma_{j_0+1})\nu (\gamma_{j_0+1})}{4\pi})} \right) + \arcsinh \frac{1}{2} \\
& \geq & 2\log (\sqrt{5}+2)-2\log \left( \frac{\tan (\arcsin \frac{\sqrt{5}}{5})}{\tan (\frac{1}{2}\arcsin \frac{\sqrt{5}}{5})} \right) + \arcsinh \frac{1}{2} \\
& = & 2\log 2 + \arcsinh \frac{1}{2} > \frac{7}{2}\arcsinh \frac{1}{2} ,
\end{eqnarray*}
a contradiction. Hence $\# \{ W_{j,k} : y\in W_{j,k} \} \leq 3 $, $\forall y\in C$. Then
\begin{eqnarray*}
\sum_{j=0}^m \sum_{k=0}^{N_{j} -1} \int_{W_{j,k} } \Vert \beta_{\gamma_l} \Vert_{\mathrm{hyp}} \mu_{\mathrm{KE}} & \leq & 3 \int_{C } \Vert \beta_{\gamma_l} \Vert_{\mathrm{hyp}} \mu_{\mathrm{KE}} \\
& \leq &  3 \left( \int_{C } \Vert \beta_{\gamma_l} \Vert^2_{\mathrm{hyp}} d\mu_{\mathrm{KE}} \right)^{\frac{1}{2}} \left( \mu_{\mathrm{KE}} (C) \right)^{\frac{1}{2}} \\
& = & 12 \sqrt{\pi (g-1)} \Vert \beta_{\gamma_l} \Vert_{L^2;C} .
\end{eqnarray*}

Since the collars $\mathscr{C} (\gamma_j)$ are pairwise disjoint, we have
\begin{eqnarray*}
\sum_{j=0}^m \Vert \beta_{\gamma_l} \Vert_{L^2 ;\mathscr{C} (\gamma_j)} & \leq & \left( (m+1)\sum_{j=0}^m \Vert \beta_{\gamma_l} \Vert^2_{L^2 ;\mathscr{C} (\gamma_j)} \right)^{\frac{1}{2}} \\
& \leq & \left( (m+1) \Vert \beta_{\gamma_l} \Vert^2_{L^2 ;C} \right)^{\frac{1}{2}}  \leq \sqrt{3(g-1)} \Vert \beta_{\gamma_l} \Vert_{L^2 ;C} ,
\end{eqnarray*}
where the last inequality uses Theorem \ref{thmcollarthm}, which implies $m+1\leq 3(g-1)$.

We now consider the K\"ahler collar $\mathscr{C}(\gamma_j)$. Denote by $z_j=\varsigma_j e^{i\theta_j}$ the polar coordinates on $\mathscr{C}(\gamma_j)$. On $\mathscr{C}(\gamma_j)$, we observe that the decomposition
$$\star\beta_{\gamma_l}= \left( \frac{1}{2\pi} \int_{\gamma_j} \star\beta_{\gamma_l} \right) d\theta_j + \left( \star\beta_{\gamma_l}- \left( \frac{1}{2\pi} \int_{\gamma_j} \star\beta_{\gamma_l} \right) d\theta_j \right) $$
is $L^2$-orthogonal. It follows that
\begin{eqnarray*}
 2\Vert \beta_{\gamma_l} \Vert^2_{L^2 ;\mathscr{C} (\gamma_j)} & = & \int_{\mathscr{C} (\gamma_{j})} \beta_{\gamma_l} \wedge \star \beta_{\gamma_l} \geq \left| \frac{1}{2\pi} \int_{\gamma_j} \star\beta_{\gamma_l} \right|^2 \int_{\mathscr{C} (\gamma_{j})} d\theta_j \wedge \star d\theta_j \\
 & = & \frac{1}{ \pi } \left| \int_{\gamma_j} \star\beta_{\gamma_l} \right|^2 \left(\frac{ \pi^2}{\ell (\gamma_j)} -\nu (\gamma_j)\right) .
\end{eqnarray*}

Applying the Cauchy–Schwarz inequality, we obtain
\begin{equation*}
\left( \sum_{j=0}^m \left| \int_{\gamma_j} \star \beta_{\gamma_{l}} \right| \left| \frac{ \pi^2}{\ell (\gamma_j)} - \nu (\gamma_j) \right| \right)^2 \leq 2\pi \Vert \beta_{\gamma_l} \Vert^2_{L^2 ;C} \left( \sum_{j\in J_0} \left| \frac{ \pi^2}{\ell (\gamma_j)} - \nu (\gamma_j) \right| \right) .
\end{equation*}
Note that $ \gamma_j$ is non-separating if and only if $j\leq J_0 \in \{0,\cdots ,m\} $.

Combining the above estimates, a direct computation yields
\begin{eqnarray*}
 1 & \leq & \frac{7 \sqrt{2}\cdot \arcsinh \frac{1}{2} \cosh \left( \frac{7}{4} \arcsinh \frac{1}{2}\right)}{\pi (\sqrt{5} -2) } \cdot \sum_{j=0}^m \sum_{k=0}^{N_{j} -1} \int_{W_{j,k} } \Vert \beta_{\gamma_l} \Vert_{\mathrm{hyp}} \mu_{\mathrm{KE}} \\
 &  & +  2 \sum_{j=0}^m \Vert \beta_{\gamma_l} \Vert_{L^2 ;\mathscr{C} (\gamma_j)}  + \frac{1}{ \pi}\sum_{j=0}^m \left| \int_{\gamma_j} \star \beta_{\gamma_{l}} \right| \left| \frac{ \pi^2}{\ell (\gamma_j)} - \nu (\gamma_j) \right| \\
 & \leq & \frac{84 \sqrt{2\pi}\cdot \arcsinh \frac{1}{2} \cosh \left( \frac{7}{4} \arcsinh \frac{1}{2}\right) \sqrt{g-1}}{\pi (\sqrt{5} -2) } \cdot \Vert \beta_{\gamma_l} \Vert_{L^2 ;C} \\
 &  & + 2\sqrt{3(g-1)} \Vert \beta_{\gamma_l} \Vert_{L^2 ;C} + \frac{\sqrt{2\pi}}{\pi} \left( \sum_{j\in J_0} \left| \frac{ \pi^2}{\ell (\gamma_j)} - \nu (\gamma_j) \right| \right)^{\frac{1}{2}} \Vert \beta_{\gamma_l} \Vert_{L^2 ;C} \\
 & \leq & \sqrt{\frac{2}{\pi}} \left( 240 \sqrt{g-1} + \left( \sum_{j\in J_0} \left| \frac{ \pi^2}{\ell (\gamma_j)} - \nu (\gamma_j) \right| \right)^{\frac{1}{2}}  \right) \Vert \beta_{\gamma_l}\Vert_{L^2 ;C} .
\end{eqnarray*}
Note that
$$ \frac{84\cdot \arcsinh \frac{1}{2} \cosh \left( \frac{7}{4} \arcsinh \frac{1}{2}\right) }{ \sqrt{5} -2 } + \sqrt{6\pi} \approx 239.960239024 \cdots . $$

Hence
\begin{eqnarray*}
 \left|\int_{\gamma_l} \alpha_{\gamma_l}\right| & = & \sqrt{2}\Vert \beta_{\gamma_l} \Vert_{L^2;C} \geq \frac{\sqrt{\pi}}{240 \sqrt{g-1} + \left( \sum_{j\in J_0} \left| \frac{ \pi^2}{\ell (\gamma_j)} - \nu (\gamma_j) \right| \right)^{\frac{1}{2}}} ,
\end{eqnarray*}
and the proof is complete.
\end{proof}

Before estimate the $\varphi$-invariant, we need the following elementary lemma from calculus.

\begin{lem}
\label{lemexplicitestimatecalculussumofcollarnonseparated}
Set $u (k,t) = \frac{1+kt^3}{(1+kt)^2} $, $k\in\mathbb{N}$, $t\in [0,1]$. Then the following hold:
\begin{enumerate}[(1)]
    \item For any given $t\in [0,1]$, $u (k,t) \geq u (k+1,t) $, $\forall k\in\mathbb{N}$. 
    \item We have $u (3k-2 ,t)> \mathbf{I}_6 (k)$, where $\mathbf{I}_6 (k) $ is defined as follows:
    $$
    \begin{array}{|c|c|c|c|c|}
    \hline
	k		& 2 & 3 	& 4 & \geq 5 \\\hline
    \mathbf{I}_6(k)	& \frac{1}{4} & \frac{1}{8} 	& \frac{3}{40} & \frac{7}{20}\cdot k^{-\frac{4}{3}}  \\
    \hline
    \end{array}.
    $$
\end{enumerate}
\end{lem}

\begin{proof}
We begin with proving (1).

Extend $u(k,t)$ to a function of a real variable $x\geq 1$ by setting $ u (x,t) = \frac{1+xt^3}{(1+xt)^2} $, $\forall x\in[1,\infty )$, $t\in [0,1]$. A direct computation gives 
$$ \frac{\partial u (x,t)}{\partial x} = \frac{t( t^2 -2-xt^3 ) }{(1+xt)^3} \leq 0 ,$$
and hence $u (k+1 ,t)\leq u (k,t)  $, $\forall k\in\mathbb{N}$. This proves (1).

We now turn to the estimate (2).

For any $k\in\mathbb{N}$, let $t_k $ be the unique number in $(0,1)$ that satisfies $kt_k^3 +3t_k^2 =2 $. Clearly, $t_{k+1}<t_k$, $\forall k\in\mathbb{N}$. Then $ \frac{\partial u (k,t)}{\partial t} = \frac{k( kt^3 +3t^2 -2 ) }{(1+kt)^3} $ shows that $u (k,t) \geq u(k,t_k) $. Moreover, $kt_k^3 +3t_k^2 =2 $ implies that $ u(k,t_k) = \frac{3t_k^4}{4 (1-t_k^2)} $.

Now we consider the case $k\in \{ 2,3,4 \}$.

Since $2\cdot (\frac{67}{100})^3 +3\cdot (\frac{67}{100})^2 -2 = -\frac{25887}{500000} <0 $, we have $t_2 > \frac{67}{100}$, and hence
$$ u(2,t_2) = \frac{3t_2^4}{4 (1-t_2^2)} > \frac{3}{4} \cdot\frac{\left(\frac{67}{100}\right)^{4}}{1-\left(\frac{67}{100}\right)^{2}} > \frac{1}{4} .$$
Similarly, one finds $t_5 > \frac{29}{50} $ and $t_8 > \frac{13}{25} $, which yield $u(5,t_5) >\frac{1}{8}$ and $u(8,t_8)> \frac{3}{40}$.

Since $ kt_k^3 +3t_k^2 =2 $, we have
\begin{equation*}
(k+1)^{\frac{4}{3}} u(k,t_k ) =  \frac{3 }{4 } \cdot (1-t_k )^{\frac{1}{3}} (1+t_k )^{-1} (2 + 2t_k - t_k^2 )^{\frac{4}{3}} .
\end{equation*}

Let $u_1 (t) = (1-t ) (1+t )^{-3} (2 + 2t - t^2 )^{4} $, $\forall t\in (0,1)$. Then $(k+1)^{\frac{4}{3}} u(k,t ) \geq \frac{3 }{4 }u_1 (t_k )^{\frac{1}{3}} $. Since 
$$ \frac{du_1 (t)}{dt} =  \frac{6t(t^2 -2)\cdot u_1 (t)}{ (1-t^2) (2+2t-t^2)} <0 ,$$
the function $u_1 (t) $ is decreasing on $t$. Moreover, since $\frac{11}{8} + \frac{3}{4} = \frac{17}{8}>2$, we have $t_{11}<\frac{1}{2}$. Hence for all $k\geq 11$, 
\begin{equation*}
(k+1)^{\frac{4}{3}} u(k,t ) \geq \frac{3 }{4 } u_1 (t_k )^{\frac{1}{3}} \geq \frac{3 }{4 } u_1 \left( \frac{1}{2} \right)^{\frac{1}{3}} = \frac{(11)^{\frac{4}{3}} }{16 } > \frac{7\cdot  3^{\frac{4}{3}}}{20} .
\end{equation*}
This completes the lemma.
\end{proof}

Combining Lemma \ref{lemW12integralpeaksectionalongshortnonseparatedgeodesic}, Lemma \ref{lempeaksectionalongshortgeodesicintersectionanothergeodesicloop} and Lemma \ref{lemexplicitestimatecalculussumofcollarnonseparated}, we obtain the following estimate for Zhang's $\varphi$-invariant.

\begin{lem}
\label{lemzhangphiinvariantsumofquotient}
Let $C$ be a compact Riemann surface of genus $g \geq 2$, equipped with the hyperbolic metric $\mu_{\mathrm{KE}}$ of constant curvature $-1$. Let $\gamma_0 $ be the shortest non-separating simple closed geodesic on $(C,\mu_{\mathrm{KE}})$, and let $\mathbf{c}_0 $ be the geodesic loop defined in Lemma \ref{lemshortgeodesicintersectionanothergeodesicloop}.

Let $\gamma_0 ,\gamma_1 ,\cdots ,\gamma_m$ be the simple closed geodesics intersecting $\mathbf{c}_0$ and satisfying $\ell (\gamma_j)\leq 2\arcsinh \frac{1}{2}$. Assume that $\ell (\gamma_0) \leq 2\arcsinh \tfrac{1}{8}$. Let $J_0 \subset \{0, \dots, m\}$ denote the set of indices $j$ such that $\gamma_j$ is non-separating. Define $\varkappa_j  =\frac{\pi^2}{\ell (\gamma_j)} -\nu (\gamma_j) $ and $\varkappa_{\Sigma} = \sum_{j\in J_0} \varkappa_j$. Then for any $\epsilon\in (0,\frac{1}{2})$, if
$$ \varkappa_0 \geq \max\left\{ 240 (1+\epsilon)\cdot \sqrt{7(g-1)} ,\;\; \frac{7\cdot (1+\epsilon)^2}{\mathbf{I}_6 (g) \cdot \epsilon^2} \right\} ,$$
we have
$$ \varphi (C) \geq \min\left\{ \frac{ \varkappa_0 }{3\cdot (240)^2 (1+\epsilon)^2 g } ,\;\; \frac{\epsilon^4 (g-1) \mathbf{I}_6 (g) \cdot \varkappa_0 }{21 (1+\epsilon)^4 g} \right\} ,$$
where $\mathbf{I}_6 $ is the constant defined in Lemma \ref{lemexplicitestimatecalculussumofcollarnonseparated}.
\end{lem}

\begin{proof}
let $\alpha_j$ denote the peak section in $\Gamma (C,\omega_C)$ along $\gamma_j$, and let $\mathbf{t}_j =\left| \int_{\gamma_j} \alpha_j \right|$. By Lemma \ref{lempeaksectionalongshortgeodesicintersectionanothergeodesicloop}, if $\ell (\gamma_j ) \leq 2\arcsinh  \frac{1}{3} $, then $ \mathbf{t}_{j } \geq \frac{\sqrt{\pi}}{240 \sqrt{g-1} + \sqrt{\varkappa_{\Sigma}}} $. Set $J'_0=\{j\in J_0: \ell (\gamma_j ) \leq 2\arcsinh \frac{1}{3} \}$.

Combining Proposition \ref{propzhangphiinvariantpotential} with Lemma \ref{lemW12integralpeaksectionalongshortnonseparatedgeodesic}, we deduce
\begin{equation*}
    \varphi (C) \geq \sum_{j=0}^m \max\left\{0,  \frac{2 (g-1) \mathbf{t}_{j }^2 \varkappa_j }{ \pi g }  \left( \frac{\mathbf{t}_{j }^2 \varkappa_j^2 }{6\pi } -1 \right) \right\} .
\end{equation*}

Let $\epsilon\in (0,\frac{1}{2})$. 

If $\sqrt{\varkappa_{\Sigma}} \leq 240 \epsilon \sqrt{g-1} $, then $ \mathbf{t}_0 \geq \frac{\sqrt{\pi}}{240 (1+\epsilon ) \sqrt{g-1} } $, which implies
\begin{eqnarray*}
\varphi (C) & \geq & \frac{ 2\varkappa_0 }{ (240)^2 (1+\epsilon)^2 g } \left( \frac{ \varkappa_0^2 }{6\cdot (240)^2 (1+\epsilon)^2 (g-1) } -1 \right) \\
& \geq & \frac{ 2\varkappa_0 }{ (240)^2 (1+\epsilon)^2 g } \left( \frac{7}{6} -1 \right) = \frac{ \varkappa_0 }{ 3\cdot (240)^2 (1+\epsilon)^2 g } .
\end{eqnarray*}

If instead $\sqrt{\varkappa_{\Sigma}} \geq 240 \epsilon \sqrt{g-1} $, then for any $j\in J'_0$, $ \mathbf{t}_j \geq \frac{\epsilon\sqrt{\pi}}{ (1+\epsilon ) \sqrt{\varkappa_{\Sigma}} } $. Using this bound and enlarging the summation from $J'_0$ to $J_0$, we obtain
$$ \varphi (C) \geq \frac{ 2\epsilon^2 (g-1) }{ (1+\epsilon)^2 g } \left( \frac{ \epsilon^2 \sum_{j\in J'_0} \varkappa_j^3 }{6 (1+\epsilon )^2 \varkappa_\Sigma^2 } -\frac{\sum_{j\in J'_0} \varkappa_j}{\varkappa_\Sigma} \right) \geq \frac{ 2\epsilon^2 (g-1) }{ (1+\epsilon)^2 g } \left( \frac{ \epsilon^2 \sum_{j\in J_0} \varkappa_j^3 }{6 (1+\epsilon )^2 \varkappa_\Sigma^2 } -1 \right) ,$$
where we used that $\ell (\gamma_j ) \geq 2\arcsinh \frac{1}{3} $ implies $ \varkappa_j^2 < 6 (1+\frac{1}{\epsilon})^2 \varkappa_\Sigma $. By Lemma \ref{lemexplicitestimatecalculussumofcollarnonseparated} and our assumption on $\kappa_0$, we can conclude that
$$ \varphi (C) \geq \frac{ 2\epsilon^2 (g-1) }{ (1+\epsilon)^2 g } \left( \frac{ \epsilon^2 \mathbf{I}_6 (g) \cdot \varkappa_0 }{6 (1+\epsilon )^2 } -\frac{ \epsilon^2 \mathbf{I}_6 (g) \cdot \varkappa_0 }{7 (1+\epsilon )^2 } \right) = \frac{ \epsilon^4 (g-1) \mathbf{I}_6 (g) \cdot \varkappa_0 }{ 21(1+\epsilon)^4 g } . $$
This completes the proof.
\end{proof}

We are now in a position to estimate the $\varphi$-invariant of hyperbolic Riemann surfaces that admit a short non-separating geodesic.

\begin{prop}
\label{propositionanalysispartphiinvariantshortnonseparatedgeodesic}
Let $C$ be a compact Riemann surface of genus $g \geq 2$, equipped with the hyperbolic metric $\mu_{\mathrm{KE}}$ of constant curvature $-1$. Let $\gamma_0 $ be the shortest non-separating simple closed geodesic on $(C,\mu_{\mathrm{KE}})$. Assume that $\ell (\gamma_0) \leq 2\arcsinh \tfrac{1}{8}$.

If $\ell (\gamma_0)\leq \mathbf{I}_7 (g)$, then $\varphi (C) \geq \frac{\mathbf{I}_8 (g)}{\ell (\gamma_0)}$, where the constants $\mathbf{I}_7 (g) $, $\mathbf{I}_8 (g) $ are given in the following table:
$$
\begin{array}{|c|c|c|c|c|}
\hline
	g		& 2 & 3 	& 4 & \geq 5 \\\hline
    \mathbf{I}_7 (g)	&\frac{1}{80} & \frac{17}{2000} 	& \frac{4}{625} & \frac{1}{22} \cdot g^{-\frac{4}{3}}  \\\hline
    \mathbf{I}_8 (g)	& \frac{1}{100000} & \frac{1}{120000} & \frac{1}{140000} & \frac{1}{30000}\cdot g^{-\frac{4}{3}}  \\
\hline
\end{array}.
$$
\end{prop}

\begin{proof}
Let $\mathbf{c}_0 $ be the geodesic loop defined in Lemma \ref{lemshortgeodesicintersectionanothergeodesicloop}. Let $\gamma_0 ,\gamma_1 ,\cdots ,\gamma_m$ be the short simple closed geodesics with $\ell (\gamma_j)\leq 2\arcsinh\left( \frac{1}{2}\right)$ and $\gamma_j\cap\mathbf{c}_0 \neq \varnothing $. For $j=0,\cdots ,m$, let $\alpha_{\gamma_j}$ be a peak section in $\Gamma (C,\omega_C)$ along $\gamma_j $, and let $\mathscr{C} (\gamma_j)$ be the collar of $\gamma_j $. Denote $\frac{\pi^2}{\ell (\gamma_j)} -\nu (\gamma_j) $ and $\left| \int_{\gamma_j} \alpha_j \right|$ by $\varkappa_j $ and $\mathbf{t}_j $, respectively. Let $J_0 \subset \{0,\dots,m\}$ be the set of indices $j$ for which $\gamma_j$ is non-separating, and set $\varkappa_{\Sigma} = \sum_{j\in J_0} \varkappa_j$.

For any $g\geq 2$, define $\epsilon_1(g)$ by
$$
\begin{array}{|c|c|c|c|c|}
\hline
	g		& 2 & 3 	& 4 & \geq 5 \\\hline
    \epsilon_1 (g)	& \frac{47}{200} & \frac{113}{400} 	& \frac{42}{125} & \frac{11}{25}   \\
\hline
\end{array}.
$$
Applying Lemma \ref{lemzhangphiinvariantsumofquotient} with the choice $\epsilon=\epsilon_1(g)$, we conclude that if $\varkappa_0\geq \vartheta_1 (g) $, then $\varphi (C) \geq \epsilon_2 (g)\varkappa_0$, where $\vartheta_1 (g)$ and $\epsilon_2 (g)$ are constants defined by
$$
\begin{array}{|c|c|c|c|c|}
\hline
	g		& 2 & 3 	& 4 & \geq 5 \\\hline
    \vartheta_1 (g)	& 785 & 1155 	& 1476 & \frac{429}{2} \cdot g^{\frac{4}{3}}  \\\hline
    \epsilon_2 (g)	& \frac{1}{600000} & \frac{1}{1000000} 	& \frac{1}{1250000} & \frac{1}{250000}\cdot g^{-\frac{4}{3}}  \\
\hline
\end{array}.
$$

Since $\varkappa_0 =\frac{\pi^2}{\ell (\gamma_0)} -\nu (\gamma_0) $ and $\nu (\gamma_0) <\pi$, we can conclude that if $\ell (\gamma_0) \leq \frac{\pi^2}{\vartheta_1 (g)+\pi} $, then $\varphi (C) \geq \frac{\epsilon_2 (g) \left(\pi^2 -\nu (\gamma_0)\ell (\gamma_0)\right)}{ \ell (\gamma_0)}$. It follows that if $\ell (\gamma_0)\leq \mathbf{I}_7 (g)$, then $\varphi (C) \geq \frac{\mathbf{I}_8 (g)}{\ell (\gamma_0)}$, where the constants $\mathbf{I}_7 (g) $, $\mathbf{I}_8 (g) $ are as given in the table above. This completes the proof.
\end{proof}

Combining Proposition \ref{propexplicitlowerboundphiinvthickpart}, Proposition \ref{propositionanalysispartphiinvariantcutsurfaces} and Proposition \ref{propositionanalysispartphiinvariantshortnonseparatedgeodesic}, we obtain the following global estimate for Zhang's $\varphi$-invariant.

\begin{thm}
\label{thmglobalestimatezhangphiinvariant}
Let $C$ be a compact Riemann surface of genus $g \geq 2$, equipped with the hyperbolic metric $\mu_{\mathrm{KE}}$ of constant curvature $-1$. Denote by $\mathrm{sys}(C)$ the systole of $(C,\mu_{\mathrm{KE}})$, and by $\varphi (C) $ the $\varphi$-invariant of $C$. Then 
$$\varphi (C) \geq \max \left\{ \xi_1 (g) , \frac{\xi_2 (g)}{\mathrm{sys}(C)} \right\} ,$$
where the constants $\xi_1 (g) $, $\xi_2 (g) $ are given in the following table:
$$
\begin{array}{|c|c|c|c|c|}
\hline
	g		& 2 & 3 & 4 & \geq 5 \\\hline
    \xi_1 (g)	&\frac{1}{6400} & \frac{1}{10900} 	& \frac{1}{15500} & \frac{1}{15625} \cdot g^{-\frac{1}{3}}  \\\hline
    \xi_2 (g)	& \frac{1}{512000} & \frac{17}{21800000} & \frac{1}{2421875} & \frac{1}{343750} \cdot g^{-\frac{5}{3}}  \\
\hline
\end{array}.
$$
\end{thm}

\begin{proof}
Let $\gamma_0$ be the shortest closed geodesic on $(C,\mu_{\mathrm{KE}})$. Then $\ell (\gamma_0)=\mathrm{sys}(C)$. Let $\mathbf{I}_7 (g)$ be the constant defined in Proposition \ref{propositionanalysispartphiinvariantshortnonseparatedgeodesic}.

If $\mathrm{sys}(C)\geq \mathbf{I}_7 (g)$, then Proposition \ref{propexplicitlowerboundphiinvthickpart} yields
$$ \varphi (C)\geq \frac{\mathbf{I}_5 (g)\cdot \mathbf{I}_7 (g)}{\mathbf{I}_7 (g)+2\arcsinh1} ,$$
where $\mathbf{I}_5 (g)$ is the constant defined in Proposition \ref{propexplicitlowerboundphiinvthickpart}. By a direct computation, when $g=2,3,4$, we obtain $\frac{\mathbf{I}_5 (g)\cdot \mathbf{I}_7 (g)}{\mathbf{I}_7 (g)+2\arcsinh1} \geq \xi_1 (g) $. Now assume $g\geq 5$. By definition, one checks that $\mathbf{I}_5 (g)\geq \frac{g}{400}$, and hence for any $g\geq 5$,
\begin{equation*}
    \varphi (C) \geq \frac{\mathbf{I}_5 (g)\cdot \mathbf{I}_7 (g)}{\mathbf{I}_7 (g)+2\arcsinh1} \geq \frac{g}{400} \cdot \frac{g^{-\frac{4}{3}}}{22} \cdot \frac{1}{2\arcsinh1 + \frac{g^{-\frac{4}{3}}}{22} } \geq \frac{g^{-\frac{1}{3}}}{15625} = \xi_1 (g).
\end{equation*}
Moreover, by definition we have $\xi_2 (g)=\mathbf{I}_7 (g) \cdot\xi_1 (g)$. Therefore, 
$$\varphi (C) \geq \max \left\{ \xi_1 (g) , \frac{\xi_2 (g)}{\mathrm{sys}(C)} \right\},\;\;\;\textrm{ if }\; \mathrm{sys}(C)\geq \mathbf{I}_7 (g) .$$

Next, consider the case $\mathrm{sys}(C)\leq \mathbf{I}_7 (g)$. In this case, $\max \left\{ \xi_1 (g) , \frac{\xi_2 (g)}{\mathrm{sys}(C)} \right\} = \frac{\xi_2 (g)}{\mathrm{sys}(C)} $. 

Assume first that $\gamma_0$ is separated. Let $C_1$, $C_2 $ denote the two connected components of $C\setminus \left( \bigcup_{j=1}^m \gamma_j \right) $, and let $g_1$ and $g_2$ denote the genera of $\overline{C}_1$ and $\overline{C}_2$, respectively. Since $e^{15-\frac{\pi^2}{\mathbf{I}_7 (g)}} \leq \frac{1}{4} $, Proposition \ref{propositionanalysispartphiinvariantcutsurfaces} gives
\begin{eqnarray*}
\varphi(C) & \geq & \frac{\pi^2\cdot g_1g_2 }{\mathrm{sys}(C) }\left(\frac{1}{g}-\frac{2 }{g } \cdot e^{15-\frac{\pi^2}{\mathrm{sys}(C)}} \right) \geq \frac{\pi^2\cdot g_1 (g-g_1) }{\mathrm{sys}(C)\cdot g }\left(1-\frac{1}{2} \right) \geq \frac{\pi^2}{4 \mathrm{sys}(C) } > \frac{\xi_2 (g)}{\mathrm{sys}(C)} .
\end{eqnarray*}

If $\gamma_0$ is non-separated, then combining $\mathbf{I}_8 (g)\geq \xi_2 (g)$ with Proposition \ref{propositionanalysispartphiinvariantshortnonseparatedgeodesic} yields $\varphi (C) \geq \frac{\xi_2 (g)}{\mathrm{sys}(C)} $. This proves the theorem.
\end{proof}

\section{Lower bounds of the Faltings--Elkies invariant}

Let $C$ be a curve of genus $g \geq 2$ over $\mathbb{C}$, and let $\omega_C$ denote its canonical bundle. Let $\mu_{\mathrm{hyp}}$ be the unique K\"ahler metric on $C$ with constant curvature $4\pi(1-g)$, normalized to have volume $1$. Then $\mu_{\mathrm{KE}} = 4\pi(g-1)\,\mu_{\mathrm{hyp}}$ is the unique K\"ahler metric on $C$ with constant curvature $-1$. Since $\dim C = 1$, we identify the Kähler metric $\mu_{\mathrm{hyp}}$ with the measure it induces. Let $\{\alpha_j\}_{j=1}^g$ be an $L^2$-orthonormal basis of $\Gamma(C,\omega_C)$. The Arakelov K\"ahler metric is defined by
$$
\mu_{\mathrm{Ar}} = \frac{i}{2g} \sum_{k=1}^g \alpha_k \wedge \bar{\alpha}_k .
$$
This definition is independent of the choice of basis.

Let $G_{\mathrm{Ar}}$ and $G_{\mathrm{hyp}}$ denote the Green functions associated with $\mu_{\mathrm{Ar}}$ and $\mu_{\mathrm{hyp}}$, respectively. Our goal in this section is to derive explicit lower bounds for the Faltings--Elkies invariant
$$
\mathrm{FE}_{v}(n,C)=\inf_{x_1,\dots, x_n\in C} \sum_{1\leq j<k\leq n} G_{\mathrm{Ar}} (x_j,x_k),
$$ 
most notably an explicit bound in terms of Zhang’s $\varphi$-invariant.

The key idea is to transform the estimate for the Arakelov Green function into one involving hyperbolic geometric data. The first step is to observe that the difference between the two Green functions can be expressed through a single-variable function.

\begin{thm}
\label{thmdifferencearakelovgreenhyperbolicgreen}
There exists a unique $\mathbb{R}$-valued smooth function $\psi_{\mathrm{Ar}}$ on $C$, such that
$$G_{\mathrm{Ar}} (x,y) - G_{\mathrm{hyp}} (x,y) = \psi_{\mathrm{Ar}} (x) +\psi_{\mathrm{Ar}} (y) .$$

In particular, $dd^c \psi_{\mathrm{Ar}}  = \mu_{\mathrm{Ar}} - \mu_{\mathrm{hyp}}$ and $\int_C \psi_{\mathrm{Ar}}  \mu_{\mathrm{Ar}} + \int_{C} \psi_{\mathrm{Ar}} \mu_{\mathrm{hyp}} =0 $.
\end{thm}

\begin{proof}
It suffices to check that the sum $G_{\mathrm{hyp}} (x,y) + \psi_{\mathrm{Ar}} (x) +\psi_{\mathrm{Ar}} (y)$ satisfies the defining properties of $G_{\mathrm{Ar}} (x,y)$, and then apply the uniqueness of the Green function. See \cite[Theorem 3.8]{jorkra1} for more details.
\end{proof}

As a preliminary step toward the explicit lower bound estimates, we establish the following lemma. Its proof is based on an argument similar to that of \cite[Proposition 3.2.5]{autissier1}, which generalizes the original estimate of Faltings and Elkies \cite[Theorem 5.1]{lang1}.

\begin{lem}
\label{lemarakelovgreenhyperbolicgreen}
Let $x_1,\cdots ,x_n $ be pairwise distinct points on $C$. Then for any $t>0$,
$$ \sum_{1\leq j< k\leq n} G_{\mathrm{Ar}} (x_j,x_k) \geq - \frac{1}{2} \sum_{j=1}^n G_{\mathrm{hyp}} (x_j,x_j;t ) -\frac{n(n-1)t}{2} \left( 1 +\sup_C \left| \frac{\mu_{\mathrm{Ar}} }{\mu_{\mathrm{hyp}}} -1 \right| \right) .$$
\end{lem}

\begin{proof}
Let $\psi_{\mathrm{Ar}}$ be the function defined in Theorem \ref{thmdifferencearakelovgreenhyperbolicgreen}, and $ \psi_{t} (x)$ is the solution of the hyperbolic heat equation $\frac{du}{dt} = \Delta_{\mathrm{hyp},x} u $ with initial data $u (x ;0)=\psi_{\mathrm{Ar}}  (x) $. Set $a_{l}=\int_{C} (\psi_{\mathrm{Ar}} \phi_{\mathrm{hyp} , l }) \mu_{\mathrm{hyp}} $, $\forall l\geq 0$. Then $\psi_t (x ) = \sum_{l=0}^\infty a_le^{-t\lambda_{\mathrm{hyp} , l }}\phi_{\mathrm{hyp} , l } (x) $. Since $\frac{d\psi_t (x)}{dt} = \Delta_{\mathrm{hyp},x} \psi_t (x ) $ also satisfies the heat equation, we can apply the maximum principle to show that 
\begin{eqnarray*}
|\psi_t (x ) - \psi_{\mathrm{Ar}}  (x )| & = & \left|\int_0^t \frac{d}{d\varsigma} \psi_\varsigma (x )d\varsigma\right| = \left|\int_0^t\Delta_{\mathrm{hyp},x} \psi_{\varsigma} (x )d\varsigma\right| \\
& \leq & t\sup_C \left|\Delta_{\mathrm{hyp},x} \psi_{\mathrm{Ar}}  (x )\right|=t\sup_C \left| \frac{\mu_{\mathrm{Ar}}-\mu_{\mathrm{hyp}}}{\mu_{\mathrm{hyp}}} \right|  .
\end{eqnarray*}

By the maximum principle of the heat equations, we have $  G_{\mathrm{hyp}} (x,y;t) \leq t+G_{\mathrm{hyp}} (x,y) $, $\forall x\neq y$. See also \cite[Lemma 5.2]{lang1}. Then we have
\begin{eqnarray*}
\sum_{1\leq j< k\leq n} G_{\mathrm{Ar}} (x_j,x_k) & = & \sum_{1\leq j< k\leq n} \left( G_{\mathrm{hyp}} (x_j,x_k) + \psi_{\mathrm{Ar}}  (x_j) + \psi_{\mathrm{Ar}}  (x_k) \right) \\
& \geq & \sum_{1\leq j< k\leq n} G_{\mathrm{hyp} } (x_j,x_k ;t) -\frac{n(n-1)t}{2} \\
& & + (n-1) \sum_{j=1}^n \psi_{\frac{t}{2}} \left( x_j \right) - \frac{n(n-1)t}{2}\sup_C \left| \frac{\mu_{\mathrm{Ar}} }{\mu_{\mathrm{hyp}}} -1 \right| .
\end{eqnarray*}

For any $l\geq 1$, set $\Phi_l (t) = \sum_{j=1}^n e^{-\frac{t\lambda_{\mathrm{hyp},l}}{2} } \phi_{\mathrm{hyp} , l }(x_j) $. Hence we have
\begin{eqnarray*}
\sum_{j=1}^n \psi_{\frac{t}{2}} \left( x_j \right)
& = & \sum_{j=1}^n \sum_{l=0}^\infty e^{-\frac{t\lambda_{\mathrm{hyp},l}}{2} } a_l \phi_{\mathrm{hyp} , l }(x_j) = n a_0 + \sum_{l=1}^\infty a_l \Phi_l (t) ,
\end{eqnarray*}
and
\begin{eqnarray*}
\sum_{1\leq j< k\leq n} G_{\mathrm{hyp} } (x_j,x_k ;t) & = & \sum_{l=1}^\infty \sum_{1\leq j< k\leq n} \frac{e^{-t\lambda_{\mathrm{hyp},l} }}{\lambda_{\mathrm{hyp},l}} \phi_{\mathrm{hyp} , l } (x_j) \phi_{\mathrm{hyp} , l }(x_k) \\
& = & \sum_{l=1}^\infty\frac{\Phi_l^2 (t) }{2 \lambda_{\mathrm{hyp},l}} - \sum_{l=1}^\infty\sum_{j=1}^n \frac{e^{-t\lambda_{\mathrm{hyp},l} }}{2\lambda_{\mathrm{hyp},l}} (\phi_{\mathrm{hyp} , l } (x_j))^2 \\
& = & \sum_{l=1}^\infty \frac{\Phi_l^2(t) }{2 \lambda_{\mathrm{hyp},l}} - \frac{1}{2} \sum_{j=1}^n G_{\mathrm{hyp};t} (x_j,x_j ) .
\end{eqnarray*}

It follows that
\begin{eqnarray*}
\sum_{1\leq j< k\leq n} G_{\mathrm{Ar}} (x_j,x_k) & \geq & \sum_{l=1}^\infty \frac{\Phi_l^2 (t) }{2 \lambda_{\mathrm{hyp},l}} + (n-1)  \sum_{l=1}^\infty  a_l \Phi_l (t) - \frac{1}{2} \sum_{j=1}^n G_{\mathrm{hyp};t} (x_j,x_j ) \\
& & + n(n-1)a_0 -\frac{n(n-1)t}{2} \left( 1+ \left\Vert \mathbf{b}_0 -1 \right\Vert_{\sup} \right) \\
& \geq & - \sum_{l=1}^\infty \frac{(n-1)^2 a_l^2 \lambda_{\mathrm{hyp},l} }{2} - \frac{1}{2} \sum_{j=1}^n G_{\mathrm{hyp};t} (x_j,x_j ) \\
& & + n(n-1)a_0 -\frac{n(n-1)t}{2} \left( 1+ \sup_C \left| \frac{\mu_{\mathrm{Ar}} }{\mu_{\mathrm{hyp}}} -1 \right| \right) .
\end{eqnarray*}

By definition, we have $a_0 = \int_C \psi_{\mathrm{Ar}}  d\mu_{\mathrm{hyp}} $ and $$\sum_{l=1}^\infty a_l^2 \lambda_{\mathrm{hyp},l} = -\int_C \psi_{\mathrm{Ar}}  dd^c\psi_{\mathrm{Ar}}  = \int_C \psi_{\mathrm{Ar}}  (\mu_{\mathrm{hyp}} - \mu_{\mathrm{Ar}}) = 2\int_C \psi_{\mathrm{Ar}}  \mu_{\mathrm{hyp}}=2a_0 .$$
Hence $a_0 >0$, and
\begin{eqnarray*}
n(n-1)a_0 - \sum_{l=1}^\infty \frac{(n-1)^2 a_l^2 \lambda_{\mathrm{hyp},l} }{2} & = & n(n-1)a_0 - (n-1)^2 a_0 \\
& = & (n-1)a_0 \geq 0 .
\end{eqnarray*}
This completes the proof.
\end{proof}

Applying Lemma \ref{lemarakelovgreenhyperbolicgreen} with $t= \frac{1}{40n(g-1)} $, we obtain the following estimate.

\begin{lem}
\label{lemarakelovgreenhyperbolicdatathreeparts}
Let $x_1,\cdots ,x_n $ be pairwise distinct points on $C$. Then
\begin{eqnarray*}
\sum_{1\leq j< k\leq n} G_{\mathrm{Ar}} (x_j,x_k) & \geq & -\frac{ n-1 }{80(g-1)} \left( 1 + \sup_{ C} \left| \frac{\mathbf{b}_0 }{g} -1 \right| \right)  \\
& & - \frac{5n}{2(g-1)} \sup_{x\in C} \left| K_{\mathrm{hyp}} \left( x,x;\frac{1}{40(g-1)}\right) \right| \\
& & -\frac{n}{2}\int_{\frac{1}{40n(g-1)}}^{\frac{1}{40(g-1)}} \sup_{x\in C} K_{\mathrm{hyp}} (x,x;\varsigma) d\varsigma \\
& & - \frac{n}{2} \sum_{0<\lambda_{\mathrm{KE},l} \leq \frac{1}{20\pi}} \frac{ \sup\limits_{C} |\phi_{\mathrm{KE} , l } |^2}{ \lambda_{\mathrm{KE},l}},
\end{eqnarray*}
where $\mathbf{b}_0(x) $ is the Bergman kernel function associated with $(\Gamma(C,\omega_C),\mu_{\mathrm{hyp}})$.
\end{lem}

\begin{proof}
From Lemma \ref{lemarakelovgreenhyperbolicgreen} with $t=\tfrac{1}{40n(g-1)}$, we deduce that
$$\frac{n(n-1)t}{2} \left( 1 +\sup_C \left| \frac{\mu_{\mathrm{Ar}} }{\mu_{\mathrm{hyp}}} -1 \right| \right)=\frac{ n-1 }{80(g-1)} \left( 1 + \frac{1}{g}\sup_{x\in C}\left| \mathbf{b}_0 (x) -g \right| \right).$$
We are thus reduced to proving that, for each $j=1,\dots,n$,
\begin{eqnarray*}
 G_{\mathrm{hyp}} \left( x_j,x_j;\frac{1}{40n(g-1)} \right) & \leq & \sum_{\lambda_{\mathrm{KE},l} \leq \frac{1}{20\pi}} \frac{ \sup\limits_{C} |\phi_{\mathrm{KE} , l }|^2}{ \lambda_{\mathrm{KE},l}} + \int_{\frac{1}{40n(g-1)}}^{\frac{1}{40(g-1)}} \sup_{x\in C} \left| K_{\mathrm{hyp}} (x,x;\varsigma) \right| d\varsigma \\
 & & + \frac{5}{g-1} \sup_{x\in C} \left| K_{\mathrm{hyp}} \left( x,x;\frac{1}{40(g-1)}\right) \right|.
\end{eqnarray*}

By definition, we have
\begin{eqnarray*}
G_{\mathrm{hyp}} \left( x_j,x_j;\frac{1}{40n(g-1)} \right) & = & \sum_{l=1}^\infty \frac{e^{-\frac{\lambda_{\mathrm{hyp},l}}{40n(g-1)} }}{ \lambda_{\mathrm{hyp},l}} |\phi_{\mathrm{hyp} , l } (x_j )|^2 \\
& \leq & \sum_{\lambda_{\mathrm{hyp},l} \geq \frac{g-1}{5}} \frac{e^{-\frac{\lambda_{\mathrm{hyp},l}}{40n(g-1)} }}{ \lambda_{\mathrm{hyp},l}} |\phi_{\mathrm{hyp} , l } (x_j )|^2 \\
& & + \sum_{0<\lambda_{\mathrm{hyp},l} \leq \frac{g-1}{5}} \frac{ |\phi_{\mathrm{hyp} , l } (x_j )|^2}{ \lambda_{\mathrm{hyp},l}} .
\end{eqnarray*}

It is enough to prove that
\begin{eqnarray*}
 \sum_{\lambda_{\mathrm{hyp},l} \geq \frac{g-1}{5}} \frac{e^{-\frac{\lambda_{\mathrm{hyp},l}}{40n(g-1)} }}{ \lambda_{\mathrm{hyp},l}} |\phi_{\mathrm{hyp} , l } (x_j )|^2 & \leq & \int_{\frac{1}{40n(g-1)}}^{\frac{1}{40(g-1)}} \sup_{x\in C} K_{\mathrm{hyp}} (x,x;\varsigma) d\varsigma \\
 & & + \frac{5}{g-1} \sup_{x\in C} \left| K_{\mathrm{hyp}} \left( x,x;\frac{1}{40(g-1)}\right) \right|.
\end{eqnarray*}

Indeed, by definition of the heat kernel,
\begin{eqnarray*}
\int_{\frac{1}{40n(g-1)}}^{\frac{1}{40(g-1)}} K_{\mathrm{hyp}} (x,x;\varsigma ) d\varsigma & = & \int_{\frac{1}{40n(g-1)}}^{\frac{1}{40(g-1)}} \sum_{l=1}^\infty e^{-\varsigma \lambda_{\mathrm{hyp},l} } |\phi_{\mathrm{hyp} , l } (x)|^2 d\varsigma \\
& \geq & \sum_{\lambda_{\mathrm{hyp},l} \geq \frac{g-1}{5}} \frac{e^{-\frac{\lambda_{\mathrm{hyp},l}}{40n(g-1)}}}{\lambda_{\mathrm{hyp},l}} \left| \phi_{\mathrm{hyp} , l } (x) \right|^2 - \sum_{\lambda_{\mathrm{hyp},l} \geq \frac{g-1}{5}} \frac{e^{-\frac{\lambda_{\mathrm{hyp},l}}{40(g-1)}}}{\lambda_{\mathrm{hyp},l}} \left| \phi_{\mathrm{hyp} , l } (x) \right|^2 \\
& \geq &  \sum_{\lambda_{\mathrm{hyp},l} \geq \frac{g-1}{5}} \frac{e^{-\frac{\lambda_{\mathrm{hyp},l}}{40n(g-1)}}}{\lambda_{\mathrm{hyp},l}} \left| \phi_{\mathrm{hyp} , l } (x) \right|^2 - \frac{5}{g-1} \sum_{\lambda_{\mathrm{hyp},l} \geq \frac{g-1}{5}}e^{-\frac{\lambda_{\mathrm{hyp},l}}{40(g-1)}} \left| \phi_{\mathrm{hyp} , l } (x) \right|^2 \\
& \geq & \sum_{\lambda_{\mathrm{hyp},l} \geq \frac{g-1}{5}} \frac{e^{-\frac{\lambda_{\mathrm{hyp},l}}{40n(g-1)}}}{\lambda_{\mathrm{hyp},l}} \left| \phi_{\mathrm{hyp} , l } (x) \right|^2 - \frac{5}{g-1} \sup_{x\in C} \left| K_{\mathrm{hyp}} \left( x,x;\frac{1}{40(g-1)}\right) \right|,
\end{eqnarray*}
and the lemma follows.
\end{proof}

We now estimate the three contributions in the right-hand side of
Lemma \ref{lemarakelovgreenhyperbolicdatathreeparts}. The first term involves the Bergman kernel, for which we have the following bound.

\begin{lem}
\label{lmmbergmankernelfunctionupperbound}
Let $\varphi (C)$ denote Zhang's $\varphi$-invariant of $C$, and let $\mathbf{b}_0 $ be the Bergman kernel function associated with $(\Gamma(C,\omega_C),\mu_{\mathrm{hyp}})$. The following chain of inequalities holds:
$$ \frac{ n-1 }{80(g-1)} \left( 1 + \sup_{C} \left| \frac{\mathbf{b}_0}{g} -1 \right| \right) <\frac{ n}{8g} \cdot \left( 1 + \frac{2}{\mathrm{sys}(C)} \right) \leq\frac{ n}{8g} \cdot \left( \frac{1}{\xi_1} + \frac{2}{\xi_2} \right) \cdot \varphi (C) ,$$
where $\mathrm{sys}(C)$ denotes the systole of $(C,\mu_{\mathrm{KE}})$, and $\xi_1(g)$, $\xi_2 (g)$ are the positive constants defined in Theorem \ref{thmglobalestimatezhangphiinvariant}.
\end{lem}

\begin{remark}
The first inequality in Lemma \ref{lmmbergmankernelfunctionupperbound} improves upon an estimate of Jorgensen--Kramer \cite[Proposition 4.4]{jorkra2}. More precisely, they proved that for any unramified covering $C\to C_0$,
$$ \sup_{C} \frac{\mu_{\mathrm{Ar}}}{\mu_{\mathrm{hyp}}} \leq \frac{1.2\cdot 10^3 e^{\mathrm{sys}(C_0)/2}}{(1-e^{-\mathrm{sys}(C_0)/4})^{5/2}} .$$
Although our lemma is stated only in the case $C_0 = C$, the general case considered in \cite{jorkra2} can be deduced immediately. Indeed, one applies the lemma together with the inequality $\mathrm{sys}(C)\geq \mathrm{sys}(C_0)$. In particular, our estimate improves the exponent in the denominator from $\frac{5}{2}$ to $1$.
\end{remark}

\begin{proof}
Let $\mathbf{b}_{0,\mu_{\mathrm{KE}} }$ be the Bergman kernel function associated with $(\Gamma(C,\omega_C),\mu_{\mathrm{KE}})$. Since $\mu_{\mathrm{KE}} = 4\pi (g-1) \mu_{\mathrm{hyp}} $, it follows that $\mathbf{b}_0=4\pi (g-1)\mathbf{b}_{0,\mu_{\mathrm{KE}} } $. Because $C$ is compact, there exists a point $x\in C$ with $\mathbf{b}_{0,\mu_{\mathrm{KE}} } (x) = \sup_C \left| \mathbf{b}_{0,\mu_{\mathrm{KE}} } \right| $. To estimate this quantity, we pass to the covering map $\rho_{\mathbb{D}} : \mathbb{D}\to C$ from the Poincar\'e disk such that $\rho_{\mathbb{D}} (0)=x $. Let $\alpha\in \Gamma (C,\omega_C)$ be the peak section at $x$. Then $\frac{i}{2}\int_{C} \alpha \wedge \bar{\alpha}=1$. Write $\rho_{\mathbb{D}}^* \alpha = \mathbf{f} (z)dz$. Since $\rho_{\mathbb{D}}^* \mu_{\mathrm{KE}} (0) = 2i dz \wedge d\bar{z}$, we obtain $\mathbf{b}_{0,\mu_{\mathrm{KE}} } (x) = \frac{1}{4} |\mathbf{f} (0) |^2 $.

By Lemma \ref{lemlocalstructureorbitcollartheorem}, we have
$$\# \left( \rho^{-1}_{\mathbb{D}} (y) \cap \mathbb{D}_{\sqrt{2} -1} (0)  \right) < 1 + \frac{2\arcsinh 1}{\mathrm{sys}(C)} ,\;\; \forall y\in C ,$$
where $\mathrm{sys}(C)$ denotes the systole of $(C,\mu_{\mathrm{KE}})$, i.e., the length of the shortest closed geodesic on $(C,\mu_{\mathrm{KE}})$. Hence
\begin{eqnarray*}
\pi \left( \sqrt{2} -1 \right)^2 |\mathbf{f} (0) |^2 & \leq & \int_{\mathbb{D}_{\sqrt{2} -1} (0)} |\mathbf{f} (z) |^2 \mu_{\mathrm{Euc}} \\
& = & \frac{i}{2}\int_{\mathbb{D}_{\sqrt{2} -1} (0)} \rho_{\mathbb{D}}^* \alpha \wedge \rho_{\mathbb{D}}^* \bar{\alpha} \\
& \leq & 1 + \frac{2\arcsinh 1}{\mathrm{sys}(C)} .
\end{eqnarray*}

Consequently,
\begin{eqnarray*}
\sup_{x\in C} \left| \mathbf{b}_{0 } (x)\right| & = 4\pi (g-1) \mathbf{b}_{0,\mu_{\mathrm{KE}} } (x) = \pi (g-1) |\mathbf{f} (0) |^2 \\
& \leq \left( \sqrt{2} +1 \right)^2 (g-1) \left(1 + \frac{2\arcsinh 1}{\mathrm{sys}(C)} \right) ,
\end{eqnarray*}
and hence
$$1 + \sup_{x\in C} \left| \frac{\mathbf{b}_0 (x)}{g} -1 \right| \leq \frac{\left( \sqrt{2} +1 \right)^2 (g-1)}{g} \left(1 + \frac{2\arcsinh 1}{\mathrm{sys}(C)} \right) .$$

By Theorem \ref{thmglobalestimatezhangphiinvariant}, we have
\begin{eqnarray*}
\frac{ n-1 }{80(g-1)} \left( 1 + \sup_{x\in C} \left| \frac{\mathbf{b}_0 (x)}{g} -1 \right| \right) & \leq & \frac{\left( \sqrt{2} +1 \right)^2 (n-1) }{80 g} \left(1 + \frac{2\arcsinh 1}{\mathrm{sys}(C)} \right) \\
& < & \frac{9n}{80g} \left( 1 + \frac{2}{\mathrm{sys}(C)} \right) \\
& < & \frac{ n}{8g} \cdot \left( \frac{1}{\xi_1 (g)} + \frac{2}{\xi_2 (g)} \right) \cdot \varphi (C) .
\end{eqnarray*}
This completes the proof.
\end{proof}

Our next step is to estimate the sup norm of the hyperbolic heat kernel.

\begin{lem}
\label{lemsupnormhyperbolicheatkernel}
Let $x\in C$, and $n\geq 2 $. Then we have
\begin{eqnarray*}
\frac{5n}{2(g-1)} \sup_{x\in C} \left| K_{\mathrm{hyp}} \left( x,x;\frac{1}{40(g-1)}\right) \right| & < & \left(51 + \frac{100}{\mathrm{sys}(C)} \right)n \\
& \leq & \left(\frac{51}{\xi_1 (g)} + \frac{100}{\xi_2 (g)} \right)n \varphi (C) ,
\end{eqnarray*}
where $\xi_1(g)$, $\xi_2 (g)$ are the positive constants defined in Theorem \ref{thmglobalestimatezhangphiinvariant}.
\end{lem}

\begin{proof}
By Lemma \ref{lemestimatehyperbolicheatkernelondiagonal}, we have
\begin{eqnarray*}
\sup_{x\in C} \left| K_{\mathrm{hyp}} \left( x,x;\frac{1}{40(g-1)}\right) \right| & < & \frac{40(g-1)}{2} + \frac{6 \sqrt{g-1} }{\mathrm{sys}(C)} \cdot \sqrt{40(g-1)} + \frac{g-1}{100} \cdot \left(1+\frac{4\arcsinh 1}{\mathrm{sys}(C)}\right) \\
 & = & \frac{2001(g-1)}{100} + \left(6\sqrt{40}+\frac{\arcsinh1}{25}\right) \cdot \frac{g-1}{\mathrm{sys}(C)} \\
 & < & \frac{2001(g-1)}{100} + \frac{40(g-1)}{\mathrm{sys}(C)} .
\end{eqnarray*}

Hence
\begin{eqnarray*}
\frac{5n}{2(g-1)}\sup_{x\in C} \left| K_{\mathrm{hyp}} \left( x,x;\frac{1}{40(g-1)}\right) \right| & < & 51n + \frac{100n}{\mathrm{sys}(C)} \\
& \leq & \left(\frac{51}{\xi_1 (g)} + \frac{100}{\xi_2 (g)} \right)n \varphi (C) ,
\end{eqnarray*}
which completes the proof.
\end{proof}

Next we consider the integral involving the hyperbolic heat kernel.

\begin{lem}
\label{lemhyperbolicgreenlargeeigenvalue}
Let $x\in C$, and $n\geq 2 $. Then we have
\begin{eqnarray*}
\frac{n}{2}\int_{\frac{1}{40n(g-1)}}^{\frac{1}{40(g-1)}} K_{\mathrm{hyp}} (x,x;t)  dt & < & \frac{n\log n }{4} + \left(\frac{1}{8000} + \frac{1}{\mathrm{sys}(C)} \right)n \\
& \leq & \frac{n\log n}{4\xi_1 (g)} \varphi (C) + \left(\frac{1 }{8000\xi_1 (g)} + \frac{1}{\xi_2 (g)} \right)n \varphi (C) ,
\end{eqnarray*}
where $\xi_1(g)$, $\xi_2 (g)$ are the positive constants defined in Theorem \ref{thmglobalestimatezhangphiinvariant}.
\end{lem}

\begin{proof}
Applying Lemma \ref{lemestimatehyperbolicheatkernelondiagonal}, we obtain
\begin{eqnarray*}
 \int_{\frac{1}{40n(g-1)}}^{\frac{1}{40(g-1)}} K_{\mathrm{hyp}} (x,x;t) dt & < &  \int_{\frac{1}{40n(g-1)}}^{\frac{1}{40(g-1)}} \frac{1}{2t} + \frac{6 \sqrt{g-1} }{ \sqrt{t}\cdot \mathrm{sys}(C)}  + \frac{g-1}{100} \cdot \left(1+\frac{4\arcsinh 1}{\mathrm{sys}(C)}\right) dt \\
 & = & \frac{\log n}{2} + \frac{12 \sqrt{g-1} }{ \mathrm{sys}(C)} \left( \sqrt{\frac{1}{40(g-1)}} - \sqrt{\frac{1}{40n(g-1)}} \right) \\
 & & + \frac{g-1}{100} \cdot \left(1+\frac{4\arcsinh 1}{\mathrm{sys}(C)}\right) \cdot \left( \frac{1}{40(g-1)} - \frac{1}{40n(g-1)} \right) \\
 & < & \frac{\log n}{2} + \frac{1}{4000} + \left(\frac{3\sqrt{10}}{5 }+\frac{\arcsinh 1}{1000 }\right) \frac{1}{\mathrm{sys}(C)} .
\end{eqnarray*}

It follows that
\begin{eqnarray*}
  \frac{n}{2}\int_{\frac{1}{40n(g-1)}}^{\frac{1}{40(g-1)}} K_{\mathrm{hyp}} (x,x;t)  dt & < & \frac{n}{2} \left( \frac{\log n}{2} + \frac{1}{4000} + \left(\frac{3\sqrt{10}}{5 }+\frac{\arcsinh 1}{1000 }\right) \frac{1}{\mathrm{sys}(C)} \right) \\
 & < & \frac{n\log n}{4} + \left(\frac{1}{8000} + \frac{1}{\mathrm{sys}(C)} \right)n.
\end{eqnarray*}
By Theorem \ref{thmglobalestimatezhangphiinvariant}, we conclude that
\begin{eqnarray*}
\frac{n}{2}\int_{\frac{1}{40n(g-1)}}^{\frac{1}{40(g-1)}} K_{\mathrm{hyp}} (x,x;t)  dt & < & \frac{n\log n}{4\xi_1 (g)} \varphi(C) + \left(\frac{1 }{8000\xi_1 (g)} + \frac{1}{\xi_2 (g)} \right)n \varphi(C),
\end{eqnarray*}
where $\xi_1(g)$, $\xi_2 (g)$ are the positive constants defined in Theorem \ref{thmglobalestimatezhangphiinvariant}. This completes the proof.
\end{proof}

The remaining term arises from small eigenvalues. To deal with this part, we first establish a quantitative lower bound for the first nonzero eigenvalue in terms of Zhang’s $\varphi$-invariant. This is provided by the following lemma.

\begin{lem}
\label{lemfirsteigenvalueestimatephiinvariant}
Let $C $ be a connected smooth complex curve of genus $g\geq 2$, let $\mu_{\mathrm{KE}} $ be the unique K\"ahler metric on $C $ with constant curvature $-1 $, and let $\varphi (C)$ be the $\varphi$-invariant of $C$. Then the following estimate holds:
\begin{equation*}
\frac{1}{\lambda_{\mathrm{KE} ,1}} < \min\left\{ \frac{98\pi^3 (g-1)^2}{\mathrm{sys}(C)} , \frac{385\pi^3(g-1)  \cdot \varphi (C)}{\xi_2 (g)} \right\} ,
\end{equation*}
where $\mathrm{sys}(C)$ denotes the systole of $(C,\mu_{\mathrm{KE}})$, and the constant $\xi_2 (g) $ is the positive constant defined in Theorem \ref{thmglobalestimatezhangphiinvariant}.
\end{lem}

\begin{proof}
The first inequality $\frac{1}{\lambda_{\mathrm{KE} ,1}} <  \frac{98\pi^3 (g-1)^2}{\mathrm{sys}(C)}$ follows directly from Proposition \ref{propfirsteigenvalueestimate}. See also the remark following that proposition. We now turn to the second inequality.

According to Proposition \ref{propfirsteigenvalueestimate}, one of the alternatives (1)–(3) in its statement must hold. We first consider the case where (1) or (2) applies. In this situation, we obtain
\begin{eqnarray*}
\lambda_{\mathrm{KE} ,1} & \geq & \min\left\{ \frac{1}{8\pi} , \frac{1}{32\pi^3 (g-1)^2} , \frac{\ell_{\Sigma} (C,\mu_{\mathrm{KE}} ) }{2\pi^3} , \frac{1}{52\pi^3 (g-1)^2} \right\} \\
& \geq & \min\left\{ \frac{\ell_{\Sigma} (C,\mu_{\mathrm{KE}} ) }{2\pi^3} , \frac{1}{52\pi^3 (g-1)^2} \right\} \geq \min\left\{ \frac{\mathrm{sys}(C) }{2\pi^3} , \frac{1}{52\pi^3 (g-1)^2} \right\},
\end{eqnarray*}
where $\ell_{\Sigma}(C,\mu)$ denotes the infimum of the total length of all geodesics in a separating multicurve, and $\mathrm{sys}(C)$ denotes the systole of $(C,\mu_{\mathrm{KE}})$. By Theorem \ref{thmglobalestimatezhangphiinvariant}, we conclude that
\begin{eqnarray*}
\frac{1}{\lambda_{\mathrm{KE} ,1}} & \leq & \max\left\{ \frac{2\pi^3}{\mathrm{sys}(C)} , 52\pi^3 (g-1)^2 \right\} \\
& \leq & \max\left\{ \frac{2\pi^3}{\xi_2 (g)} , \frac{52\pi^3 (g-1)^2}{\xi_1 (g)} \right\} \cdot \varphi (C) < \frac{385\pi^3(g-1)  \cdot \varphi (C)}{\xi_2 (g)} .
\end{eqnarray*}
Here $\xi_1(g)$, $\xi_2 (g)$ are the positive constants defined in Theorem \ref{thmglobalestimatezhangphiinvariant}, and we used the fact that $\xi_1(g) > 20 g^{\frac{4}{3}}\,\xi_2(g)$.

Now we assume that (iii) in Proposition \ref{propfirsteigenvalueestimate} holds. Then there exists a minimal separating multicurve $\{\gamma_j \}_{j=1}^{m} $ splitting $C$ into two parts $C_1$, $C_2$, such that
$$ \lambda_{\mathrm{KE} ,1} \geq \frac{ \sum_{j=1}^{m} \ell ( \gamma_j )}{8\pi^3} \cdot \frac{1}{\min\left\{ 7(g_1 -1) +5m , 7(g_2 -1) +5m ,\frac{7(g-1)}{2} \right\}^2 } , $$
where $g_1$, $g_2$ are the genera of $C_1$ and $C_2$, respectively.

We now consider three subcases. If $\sum_{j=1}^{m} \ell ( \gamma_j ) > \frac{1}{6\log g} $, then
\begin{eqnarray*}
\frac{1}{\lambda_{\mathrm{KE} ,1}} & \leq & \frac{8\pi^3}{ \sum_{j=1}^{m} \ell ( \gamma_j )} \cdot \min\left\{ 7(g_1 -1) +5m , 7(g_2 -1) +5m ,\frac{7(g-1)}{2} \right\}^2 \\
& \leq & 8\pi^3 \cdot (6\log g) \cdot \frac{49(g-1)^2}{4} = 588\pi^3 (g-1)^2 \log g \\
& \leq & \frac{588\pi^3 (g-1)^2 \log g}{\xi_1 (g)} \cdot \varphi (C) < \frac{90\pi^3 (g-1) }{ \xi_2 (g)} \cdot \varphi (C) .
\end{eqnarray*}

On the other hand, if $\sum_{j=1}^{m} \ell ( \gamma_j ) \leq \frac{1}{6\log g} $ and $m\leq \frac{4}{5} \cdot \min \{ g_1 ,g_2 \}$, then $g_1+g_2+m-1=g$ implies $m\le\frac{2g+2}{7}$. In this situation, Proposition \ref{propositionanalysispartphiinvariantcutsurfaces} yields
\begin{eqnarray*}
\varphi (C) & \geq & \frac{\pi^2\cdot (g_1 -m+1)(g_2 -m+1)}{\sum_{j=1}^{m} \ell ( \gamma_j )} \cdot \left( \frac{g-3(m-1)}{g^2} - \frac{2}{g} \cdot e^{15-\frac{\pi^2}{\sum_{j=1}^{m} \ell ( \gamma_j )}} \right) \\
& \geq & \frac{\pi^2\cdot g_1 g_2}{25\sum_{j=1}^{m} \ell ( \gamma_j )} \cdot \left( \frac{g+15}{7g^2} - \frac{2}{g} \cdot e^{15-59 \log g} \right) \\
& > & \frac{\pi^2\cdot (g_1 +g_2 -1)}{25\sum_{j=1}^{m} \ell ( \gamma_j )} \cdot \frac{1}{7g } = \frac{\pi^2\cdot (g-m)}{25\sum_{j=1}^{m} \ell ( \gamma_j )} \cdot \frac{1}{7g } \\
& > & \frac{\pi^2 }{25\sum_{j=1}^{m} \ell ( \gamma_j )} \cdot \frac{5g-2}{7} \cdot \frac{1}{7g } > \frac{\pi^2}{400} \cdot \frac{1}{\sum_{j=1}^{m} \ell ( \gamma_j )}.
\end{eqnarray*}

It follows that
\begin{eqnarray*}
\frac{1}{\lambda_{\mathrm{KE} ,1}} & \leq & \frac{8\pi^3}{ \sum_{j=1}^{m} \ell ( \gamma_j )} \cdot \min\left\{ 7(g_1 -1) +5m , 7(g_2 -1) +5m ,\frac{7(g-1)}{2} \right\}^2 \\
& \leq & 98\pi^3 (g-1)^2 \cdot \frac{1}{ \sum_{j=1}^{m} \ell ( \gamma_j ) } < 98\pi^3 (g-1)^2 \cdot \frac{400}{\pi^2} \varphi (C) \\
& = & 39200\pi \cdot (g-1)^2 \cdot \varphi (C) < \frac{g-1}{\xi_2 (g)}\cdot \varphi (C) .
\end{eqnarray*}

Finally, in the remaining situation where $\sum\limits_{j=1}^{m} \ell(\gamma_j) \leq \frac{1}{6\log g}$ and $m \geq \frac{4}{5}\min\{g_1,g_2\}$, we may assume without loss of generality that $g_1 \leq g_2$. Then $7g_1 +5m\leq \frac{55}{4} m$, and hence
\begin{eqnarray*}
\frac{1}{\lambda_{\mathrm{KE} ,1}} & \leq & \frac{8\pi^3}{ \sum_{j=1}^{m} \ell ( \gamma_j )} \cdot \min\left\{ 7(g_1 -1) +5m , 7(g_2 -1) +5m ,\frac{7(g-1)}{2} \right\}^2 \\
& \leq & \frac{8\pi^3\cdot\min\left\{ \frac{55}{4} m ,\frac{7(g-1)}{2} \right\}^2}{m\cdot \mathrm{sys}(C)} = \frac{8\pi^3}{\mathrm{sys}(C)} \cdot \min\left\{ \frac{55}{4} \sqrt{m} ,\frac{7(g-1)}{2\sqrt{m}} \right\}^2 \\
& \leq & \frac{8\pi^3}{ \mathrm{sys}(C)} \cdot \frac{55}{4} \cdot \frac{7 (g-1)}{2} = \frac{385\pi^3 \cdot (g-1)}{\mathrm{sys}(C)} \leq \frac{385\pi^3(g-1)  \cdot \varphi (C)}{\xi_2 (g)} .
\end{eqnarray*}
This completes the proof.
\end{proof}

We are now in a position to estimate the contribution of the small eigenvalues.

\begin{lem}
\label{lemsumsmalleigenvalueestimatephiinvariant}
Let $C $ be a connected smooth complex curve of genus $g\geq 2$, let $\mu_{\mathrm{KE}} $ be the unique K\"ahler metric on $C $ with constant curvature $-1 $, and let $\varphi (C)$ be the $\varphi$-invariant of $C$. Then the estimate below holds:
\begin{equation*}
\frac{n}{2} \sum_{0<\lambda_{\mathrm{KE},l} \leq \frac{1}{20\pi}} \frac{ \sup\limits_{ C} |\phi_{\mathrm{KE} , l } |^2}{ \lambda_{\mathrm{KE},l}} < 616\pi^3 (g-1) (2g-3)n \cdot\min\left\{ \frac{\varphi (C)}{\xi_2 (g)} , \frac{14(g-1)\cdot \max\{1,\mathrm{sys}(C)\}}{55\cdot \mathrm{sys}(C)} \right\} ,
\end{equation*}
where $\mathrm{sys}(C)$ denotes the systole of $(C,\mu_{\mathrm{KE}})$, and the constant $\xi_2 (g) $ is the positive constant defined in Theorem \ref{thmglobalestimatezhangphiinvariant}.
\end{lem}

\begin{proof}
Since $\Delta_{\mathrm{KE}}= \frac{1}{2\pi} \star_{\mathrm{KE}}d \star_{\mathrm{KE}}d$, it follows that $\lambda_{\mathrm{KE},l} \leq \frac{1}{20\pi}$ if and only if $\lambda^{\mathrm{Rm}}_{\mathrm{KE},l} \leq \frac{1}{10}$. Recall from \cite[Th\'eor\`eme 2]{otro1} that if 
$\lambda^{\mathrm{Rm}}_{\mathrm{KE},l} \leq \frac{1}{4} $, then $l \leq 2g-3$. Thus it remains to show that for $0<\lambda_{\mathrm{KE},l} \leq \frac{1}{20\pi}$,
\begin{equation*}
\frac{ \sup\limits_{ C} |\phi_{\mathrm{KE} , l } |^2}{ \lambda_{\mathrm{KE},l}} < 1232\pi^3 (g-1) \cdot\min\left\{ \frac{\varphi (C)}{\xi_2 (g)} , \frac{14(g-1)\cdot \max\{1,\mathrm{sys}(C)\}}{55\cdot \mathrm{sys}(C)} \right\} .
\end{equation*}

By Proposition \ref{propglobalsupnormeigenfunctionsmalleigenvalue}, we obtain in this case
\begin{equation*}
\sup\limits_{C} |\phi_{\mathrm{KE} , l } |^2 \leq \max\left\{ 2\cdot \mathrm{sys}(C)^{-8\pi \lambda_{\mathrm{KE},l}} +1 , \frac{16}{5} \right\} .
\end{equation*}

If $\sup\limits_{C} |\phi_{\mathrm{KE} , l } |^2 \leq \frac{16}{5}  $, then by Lemma \ref{lemfirsteigenvalueestimatephiinvariant}, we have
\begin{equation*}
\frac{ \sup\limits_{ C} |\phi_{\mathrm{KE} , l } |^2}{ \lambda_{\mathrm{KE},l}} \leq \frac{16}{5}\cdot \frac{1}{\lambda_{\mathrm{KE},1}} < 1232\pi^3 (g-1) \cdot\min\left\{ \frac{\varphi (C)}{\xi_2 (g)} , \frac{14(g-1)\cdot \max\{1,\mathrm{sys}(C)\}}{55\cdot \mathrm{sys}(C)} \right\} .
\end{equation*}

Otherwise, suppose $\sup\limits_{C} |\phi_{\mathrm{KE} , l } |^2 \leq 2\cdot \mathrm{sys}(C)^{-8\pi \lambda_{\mathrm{KE},l}} +1 $. By Theorem \ref{thmglobalestimatezhangphiinvariant}, we have $\frac{1}{\mathrm{sys}(C)}\leq \frac{\varphi (C)}{\xi_2 (g)}$. Hence
\begin{equation*}
\frac{ \sup\limits_{ C} |\phi_{\mathrm{KE} , l } |^2}{ \lambda_{\mathrm{KE},l}} \leq \frac{2\cdot \mathrm{sys}(C)^{-8\pi \lambda_{\mathrm{KE},l}} +1}{\lambda_{\mathrm{KE},l}} < \frac{2 }{\lambda_{\mathrm{KE},l}} \cdot \left(\frac{\varphi (C)}{\xi_2 (g)}\right)^{8\pi \lambda_{\mathrm{KE},l}} + \frac{1}{\lambda_{\mathrm{KE},l}} .
\end{equation*}
Set $\lambda_{\min}=\frac{\xi_2 (g)}{385\pi^3(g-1)  \cdot \varphi (C)}$ and $\lambda_{\max}=\frac{1}{20\pi}$. By convexity,
\begin{equation*}
 \frac{1 }{\lambda_{\mathrm{KE},l}} \cdot \left(\frac{\varphi (C)}{\xi_2 (g)}\right)^{8\pi \lambda_{\mathrm{KE},l}} < \max\left\{ \frac{1 }{\lambda_{\max}} \cdot \left(\frac{\varphi (C)}{\xi_2 (g)}\right)^{8\pi \lambda_{\max}}, \frac{1 }{\lambda_{\min}} \cdot \left(\frac{\varphi (C)}{\xi_2 (g)}\right)^{8\pi \lambda_{\min}}  \right\}.
\end{equation*}

Evaluating at $\lambda_{\max}$ and $\lambda_{\min}$ gives
\begin{equation*}
\frac{1 }{\lambda_{\max}} \cdot \left(\frac{\varphi (C)}{\xi_2 (g)}\right)^{8\pi \lambda_{\max}} = 20\pi \cdot \left(\frac{\varphi (C)}{\xi_2 (g)}\right)^{\frac{2}{5}} \leq \frac{20\pi}{\xi_2(g)}\cdot\varphi (C) ,
\end{equation*}
and
\begin{eqnarray*}
\frac{1 }{\lambda_{\min}} \cdot \left(\frac{\varphi (C)}{\xi_2 (g)}\right)^{8\pi \lambda_{\min}} & = & \frac{385\pi^3(g-1)  \cdot \varphi (C)}{\xi_2 (g)} \cdot \left(\frac{\varphi (C)}{\xi_2 (g)}\right)^{\frac{8 \xi_2 (g)}{385\pi^2(g-1)  \cdot \varphi (C)}} \\
& \leq & \frac{385\pi^3(g-1)  \cdot \varphi (C)}{\xi_2 (g)} \cdot \sup_{t>0} e^{\frac{8\log t }{385\pi^2(g-1)  \cdot t}} \\
& \leq & \frac{385\pi^3(g-1)  \cdot \varphi (C)}{\xi_2 (g)} \cdot  e^{\frac{8 }{385e\pi^2 }}.
\end{eqnarray*}

It follows that
\begin{equation*}
\frac{ \sup\limits_{ C} |\phi_{\mathrm{KE} , l } |^2}{ \lambda_{\mathrm{KE},l}} < \frac{385\pi^3(g-1)  \cdot \varphi (C)}{\xi_2 (g)} \cdot \left( 2e^{\frac{8 }{385e\pi^2 }} +1\right) < \frac{1232\pi^3 (g-1) }{\xi_2 (g)} \cdot \varphi (C) .
\end{equation*}

Set $\lambda'_{\min}=\frac{\mathrm{sys}(C)}{98\pi^3(g-1)^2 }$. Similarly to above, we have
\begin{eqnarray*}
 \frac{\mathrm{sys}(C)^{-8\pi \lambda_{\mathrm{KE},l}} }{\lambda_{\mathrm{KE},l}} & < & \max\left\{ \frac{\mathrm{sys}(C)^{-8\pi \lambda_{\max}}  }{\lambda_{\max}} , \frac{\mathrm{sys}(C)^{-8\pi \lambda'_{\min}}  }{\lambda'_{\min}} \right\} \\
 & \leq & \max\left\{ 20\pi \cdot \mathrm{sys}(C)^{-\frac{2}{5}} , \frac{98\pi^3 (g-1)^2 \mathrm{sys}(C)^{-\frac{4\cdot \mathrm{sys}(C)}{49\pi^2}} }{\mathrm{sys}(C)} \right\} \\
 & \leq & 98\pi^3 e^{\frac{4}{49e\pi^2}} (g-1)^2 \cdot \frac{\max\{1,\mathrm{sys}(C)\}}{\mathrm{sys}(C)} .
\end{eqnarray*}
Hence
\begin{eqnarray*}
 \frac{ \sup\limits_{ C} |\phi_{\mathrm{KE} , l } |^2}{ \lambda_{\mathrm{KE},l}} & \leq & \frac{2\cdot \mathrm{sys}(C)^{-8\pi \lambda_{\mathrm{KE},l}} +1}{\lambda_{\mathrm{KE},l}} \leq 98\pi^3 (2e^{\frac{4}{49e\pi^2}} +1) (g-1)^2 \cdot \frac{\max\{1,\mathrm{sys}(C)\}}{\mathrm{sys}(C)} \\
 & < & \frac{1568\pi^3 (g-1)^2 \cdot \max\{1,\mathrm{sys}(C)\}}{5\cdot \mathrm{sys}(C)} ,
\end{eqnarray*}
which completes the proof.
\end{proof}

By combining Lemma \ref{lemarakelovgreenhyperbolicdatathreeparts}, Lemma \ref{lmmbergmankernelfunctionupperbound}, Lemma \ref{lemsupnormhyperbolicheatkernel}, Lemma \ref{lemhyperbolicgreenlargeeigenvalue} and Lemma \ref{lemsumsmalleigenvalueestimatephiinvariant}, we first obtain the following Faltings–-Elkies type estimate in terms of the systole.

\begin{prop}
\label{propfaltingselkiessystole}
Let $x_1,\cdots ,x_n $ be pairwise distinct points on $C$. Then we have
$$ \sum_{1\leq j< k\leq n} G_{\mathrm{Ar}} (x_j,x_k) \geq -\frac{1}{4}  n\log n - \frac{400\pi^3(g-1)^3\max\left\{ 1,\mathrm{sys}(C) \right\}}{\mathrm{sys}(C)} n ,$$
where $\mathrm{sys}(C)$ is the systole of $(C,\mu_{\mathrm{KE}})$.
\end{prop}

\begin{proof}
By combining Lemma \ref{lemarakelovgreenhyperbolicdatathreeparts}, Lemma \ref{lmmbergmankernelfunctionupperbound}, Lemma \ref{lemsupnormhyperbolicheatkernel}, Lemma \ref{lemhyperbolicgreenlargeeigenvalue} and Lemma \ref{lemsumsmalleigenvalueestimatephiinvariant}, we obtain
\begin{eqnarray*}
\sum_{1\leq j< k\leq n} G_{\mathrm{Ar}} (x_j,x_k) & \geq & -\frac{ n}{8g} \cdot \left( 1 + \frac{2}{\mathrm{sys}(C)} \right) - \left(51+\frac{100}{\mathrm{sys}(C)}\right) n \\
& & -\frac{ n\log n}{4 } - \left(\frac{1 }{8000 } + \frac{1 }{\mathrm{sys}(C)} \right)n \\
& & - \frac{784\pi^3 (g-1)^2 (2g-3)n}{5} \cdot \max\left\{1,\frac{1}{\mathrm{sys}(C)}\right\}\\
& > & -\frac{n\log n}{4} - 400\pi^3 (g-1)^3 n \cdot \max\left\{1,\frac{1}{\mathrm{sys}(C)}\right\} ,
\end{eqnarray*}
where $\mathrm{sys}(C)$ is the systole of $(C,\mu_{\mathrm{KE}})$. This completes the proof.
\end{proof}

The above result gives a systole-based lower bound. In contrast, the next theorem provides a similar estimate in terms of Zhang’s $\varphi$-invariant, which is of primary importance for our later arguments.

\begin{thm}[Theorem \ref{thmfaltingselkieszhangphiinv part I}]
\label{thmfaltingselkieszhangphiinv}
Let $x_1,\cdots ,x_n $ be pairwise distinct points on $C$. Then we have
$$ \sum_{1\leq j< k\leq n} G_{\mathrm{Ar}} (x_j,x_k) \geq - \left(\xi_3 (g) n\log n+ \xi_4 (g) n \right) \cdot \varphi (C)
\geq - \left(4\cdot 10^3 g^{\frac{1}{3}} n\log n+ 1.32\cdot 10^{10} g^{\frac{11}{3}} n \right) \cdot \varphi (C),$$
where the constants $\xi_3 (g) $, $\xi_4 (g) $ are given in the following table:
$$
\begin{array}{|c|c|c|}
\hline
	g		& \xi_3 (g) & \xi_4 (g) \\\hline
    2	& 1600 & 19558263230 \\\hline
    3	& 2725 & 195942159193 \\\hline
    4	& 3875 & 832634802404 \\\hline
    \geq 5	& 4000\cdot g^{\frac{1}{3}} & 13131158175 \cdot g^{\frac{11}{3}}\\
\hline
\end{array}.
$$
\end{thm}

\begin{proof}
Again, combining the five lemmas yields
\begin{eqnarray*}
\sum_{1\leq j< k\leq n} G_{\mathrm{Ar}} (x_j,x_k) & \geq & -\frac{ n}{8g} \cdot \left( \frac{1}{\xi_1} + \frac{2}{\xi_2} \right) \cdot \varphi (C) - n \cdot \left( \frac{51}{\xi_1} + \frac{100}{\xi_2} \right) \cdot \varphi (C) \\
& & -\frac{n\log n}{4\xi_1 (g)} \cdot\varphi (C)- \left(\frac{1 }{8000\xi_1 (g)} + \frac{1 }{\xi_2 (g)} \right)n \varphi (C) \\
& & - \frac{616\pi^3 (g-1) (2g-3)n}{\xi_2 (g)} \cdot \varphi (C)\\
& \geq & -\frac{n\log n}{4\xi_1 (g)} \cdot\varphi (C)- \frac{1232\pi^3 (g-1)^2 n}{\xi_2 (g)} \varphi (C) ,
\end{eqnarray*}
where the constants $\xi_1 (g)$, $\xi_2 (g) $ are the positive constants defined in Theorem \ref{thmglobalestimatezhangphiinvariant}.

Finally, one checks that
\begin{equation*}
\frac{1 }{4\xi_1 (g)} \leq \xi_3 (g) \qquad \textrm{ and } \qquad \frac{1232\pi^3 (g-1)^2 }{\xi_2 (g) }\leq \xi_4 (g) ,
\end{equation*}
which establishes the theorem.
\end{proof}

\section{Arakelov metric versus hyperbolic metric}

Let $C$ be a curve of genus $g \geq 2$ over $\mathbb{C}$, and let $\omega_C$ denote its canonical bundle. Let $\mu_{\mathrm{hyp}}$ denote the unique K\"ahler metric on $C $ with constant curvature $4\pi (1-g) $ and volume $1$, and $\mu_{\mathrm{KE}} = 4\pi (g-1) \mu_{\mathrm{hyp}} $ denote the unique K\"ahler metric on $C $ with constant curvature $-1 $. Since $\dim C =1 $, we identify a K\"ahler metric $\mu$ with its associated volume form and measure. Let $\{ \alpha_j \}_{j=1}^g $ be an $L^2$-orthonormal basis of $\Gamma (C,\omega_C)$. The Arakelov K\"ahler metric can be defined by 
$$\mu_{\mathrm{Ar}} = \frac{i}{2g} \sum_{k=1}^g \alpha_k \wedge \bar{\alpha}_k .$$

Let $G_{\mathrm{Ar}}$ be the Green function associated with $\mu_{\mathrm{Ar}}$. Recall that the Arakelov hermitian metric $\Vert \cdot \Vert_{\mathrm{Ar}}$ on the canonical bundle $\omega_C$ is defined as follows.

For any point $x\in C$ and any $\alpha \in \omega_C|_x$, choose a local coordinate $z:U_x \to \mathbb{C}$ centered at $x$ (i.e. $z(x)=0$) such that $\alpha$ coincides with $dz$ at $x$. Then the norm $\Vert \alpha \Vert_{\mathrm{Ar}} (x)$ can be determined by
$$\log \Vert \alpha \Vert_{\mathrm{Ar}} (x) = \lim_{y\to x} \left( G_{\mathrm{Ar}} (x,y) + \log |z(y)| \right) .$$

Equivalently, one can describe this construction via the diagonal embedding $C\hookrightarrow C^2 $. Indeed, the canonical bundle $\omega_C$ is the restriction of $\mathcal{O}(-\Delta)$ to the diagonal $\Delta\subset C^2 $. On $C^2$, let $\mathbf{1}_{\Delta}\in \Gamma(C^2,\mathcal{O}(\Delta))$ be the canonical section of $\mathcal{O}(\Delta)$. By declaring its hermitian norm to satisfy $-\log \Vert\mathbf{1}_{\Delta}\Vert_{\Delta} (x,y)=G_{\mathrm{Ar}} (x,y)$, we obtain a hermitian metric $\Vert\cdot\Vert_{\Delta}$ on $\mathcal{O}(\Delta)$. The Arakelov metric on $\omega_C$ is then recovered as the dual metric of $(\mathcal{O}(\Delta),\Vert\cdot\Vert_{\Delta})$ via the canonical residue isomorphism $\omega_C\cong\mathcal{O} (-\Delta)|_{\Delta}$. 

Another important hermitian metric on $\omega_C$ is the hyperbolic hermitian metric $\Vert \alpha \Vert_{\mathrm{hyp}}$, defined by the relation $i\alpha \wedge \bar{\alpha} = \frac{ \Vert \alpha \Vert_{\mathrm{hyp}}^2}{2} \mu_{\mathrm{KE}} $. 

In this section, we compare the Arakelov hermitian metric $\Vert\cdot\Vert_{\mathrm{Ar}}$ and the hyperbolic hermitian metric $\Vert\cdot \Vert_{\mathrm{hyp}}$. It is worth noting that our arguments also apply to the off-diagonal setting $(\mathcal{O}(\Delta),\Vert\cdot\Vert_{\Delta})$

In order to compare the Arakelov hermitian metric $\Vert\cdot\Vert_{\mathrm{Ar}}$ with the hyperbolic hermitian metric $\Vert\cdot \Vert_{\mathrm{hyp}}$, certain auxiliary constructions are required to connect local data with global geometry. By the thick-thin decomposition (Theorem \ref{thmcollarthm}), the analysis is naturally divided into two cases: points contained in a collar around a short geodesic and points contained in the thick part. The discussion starts with the thick part.

\begin{lem}
\label{lemconstructionlocalresidue_thick}
Let $x\in C$ be a point such that the injectivity radius $\mathrm{inj}_{\mu_{\mathrm{KE}}} (x)\geq\epsilon $ for some $\epsilon\in (0,\mathrm{arcsinh} \frac{1}{2} ]$. Let $z_{x,\epsilon}:B_{\epsilon} (x) \to \mathbb{D}$ be the local coordinate given by the inverse of the universal covering map $\rho_{\mathbb{D}} : \mathbb{D} \to C$ with $\rho_{\mathbb{D}} (0)=x$, where $B_r (x)=\{ y\in C:\dist_{\mu_{\mathrm{KE}}} (y,x)<r \}$. Then the following hold:
\begin{enumerate}[(1)]
    \item The map $z_{x,\epsilon}$ restricts to a biholomorphism 
    $$B_{\epsilon} (x)\longrightarrow \mathbb{D}_{\tanh\frac{\epsilon}{2}} (0)=\left\{z\in\mathbb{D} : |z|<\tanh\frac{\epsilon}{2} \right\} .$$
    \item $z_{x,\epsilon} (x)=0$, and the norm $\left| \log \left( \Vert dz_{x,\epsilon} \Vert_{\mathrm{hyp}} \right) \right| \leq 1 $.
\end{enumerate}

Moreover, there exists a smooth function $ \digamma_{x,\epsilon} : C\setminus \{x\} \to \mathbb{R}$ on $C\setminus \{x\}$, such that 
\begin{enumerate}[(1)]
\addtocounter{enumi}{2}
    \item  $\digamma_{x,\epsilon} (y)  = 0$ for all $y\in C \setminus B_\epsilon (x)$.
    \item The function $ \digamma_{x,\epsilon} (y) +\log | z_{x,\epsilon}(y) | $ extends to a smooth function on $B_\epsilon (x)$, and
    $$ \left| \digamma_{x,\epsilon} (y) +\log | z_{x,\epsilon}(y) | \right| \leq \log \left( \frac{3}{\tanh\frac{\epsilon}{2}} \right) . $$
    \item The integral satisfies $\left| \int_C \digamma_{x,\epsilon} (y) \mu_{\mathrm{KE}} (y) \right| \leq 4 $.
    \item In the sense of currents,
    $$-\frac{135}{32\pi | \mathrm{tanh} \frac{\epsilon}{2} |^2} \cdot\log \left( \frac{3}{\mathrm{tanh} \frac{\epsilon}{2}} \right)\cdot\mu_{\mathrm{KE}} \leq dd^c \digamma_{x,\epsilon} +\delta_x \leq \frac{135}{32\pi | \mathrm{tanh} \frac{\epsilon}{2} |^2} \cdot\log \left( \frac{3}{\mathrm{tanh} \frac{\epsilon}{2}} \right)\cdot \mu_{\mathrm{KE}} ,$$
    where $\delta_x$ is the Dirac measure at $x$.
\end{enumerate}
\end{lem}

\begin{proof}
On the Poincar\'e disk $(\mathbb{D},\mu_{\mathbb{D}})$ (see Theorem \ref{thmpoincaremodel}), we have $\mathrm{dist}_{\mu_{\mathbb{D}}} (z,0) = 2 \mathrm{arctanh} (|z|)$, $\forall z\in\mathbb{D}$, which immediately yields (1).

For (2), note that $\epsilon\leq \mathrm{arcsinh} \frac{1}{2} $ implies $\mathrm{tanh} \frac{\epsilon}{2}\leq \sqrt{5}-2$. Hence
$$\left|\log \Vert dz_{x,\epsilon}\Vert_{\mathrm{hyp}  }\right| = \left|\log \left(1-|z_{x,\epsilon}|^2\right)\right| < \log \frac{4}{3} <1 ,$$
proving (2).

We now turn to the construction of the function $\digamma_{x,\epsilon}$ required in (3)–(6).

Let $\eta :\mathbb{R}\to\mathbb{R}_{\geq 0} $ be a smooth cutoff function such that
$$\eta(\varsigma)=1 \ \text{ for }\varsigma\leq \frac{\tanh\frac{\epsilon}{2}}{3}, \qquad \eta(\varsigma)=0 \ \text{ for }\varsigma\geq \tanh\frac{\epsilon}{2},$$
and
$$|\eta'|\leq \frac{9}{2\tanh\frac{\epsilon}{2}}, \qquad |\eta''|\leq \frac{81}{4(\tanh\frac{\epsilon}{2})^2} . $$
Set $\digamma_{x,\epsilon} (y)=-\eta (|z_{x,\epsilon}(y)|) \cdot \log |z_{x,\epsilon}(y)|$ on $B_\epsilon (x)$, and set $\digamma_{x,\epsilon} =0$ on $C\setminus B_\epsilon (x)$. Then $ \digamma_{x,\epsilon} : C\setminus \{x\} \to \mathbb{R}$ is a smooth function on $C\setminus \{x\}$, and (3) follows.

Moreover, $ \digamma_{x,\epsilon} (y) +\log | z_{x,\epsilon}(y) | $ extends smoothly to $B_\epsilon (x)$, Since $\eta(\varsigma)=1$ for $\varsigma\leq \frac{\tanh\frac{\epsilon}{2}}{3}$, we obtain
$$ \left| \digamma_{x,\epsilon} (y) +\log | z_{x,\epsilon}(y) | \right| = \left| 1-\eta (|z_{x,\epsilon} (y)|) \right|\cdot \left| \log | z_{x,\epsilon}(y) | \right| \leq \log \left( \frac{3}{\tanh\frac{\epsilon}{2}} \right) , $$
establishing (4).

For (5), observe that $|\digamma_{x,\epsilon} (y)|\leq\left| \log | z_{x,\epsilon}(y) | \right| $. By Theorem \ref{thmpoincaremodel},
\begin{eqnarray*}
\left| \int_C \digamma_{x,\epsilon} (y) \mu_{\mathrm{KE}} (y) \right| & \leq & \int_{\mathbb{D}_{\tanh(\frac{\epsilon}{2})} (0)} \left| \log |z| \right| \mu_{\mathbb{D}} (z) \leq \int_{\mathbb{D}_{\sqrt{5}-2} (0)} \frac{4 \left| \log |z| \right|\mu_{\mathrm{Euc}}}{(1-|z|^2)^2} \\
& < & \frac{8\pi}{(1-(\sqrt{5}-2)^2)^2} \int_0^{\sqrt{5}-2} |\varsigma\log \varsigma | d\varsigma \\
& < & \frac{8\pi\left(\sqrt{5}-2\right)^2\left(1-\log (\sqrt{5}-2) \right)}{(1-(\sqrt{5}-2)^2)^2} \approx 3.838\cdots <4,
\end{eqnarray*}
which proves (5).

Finally, we establish (6). A direct computation on $B_\epsilon(x)$ gives
\begin{eqnarray*}
dd^c \digamma_{x,\epsilon} +\delta_x & = & - \log |z_{x,\epsilon}| \cdot dd^c \eta(|z_{x,\epsilon}|) - d\eta (|z_{x,\epsilon}|) \wedge d^c \log |z_{x,\epsilon}| -d\log |z_{x,\epsilon}| \wedge d^c \eta (|z_{x,\epsilon}|) \\
& = & \frac{1}{2\pi} \cdot\left( - \log |z_{x,\epsilon}| \cdot \eta''(|z_{x,\epsilon}|) - \frac{2+\log |z_{x,\epsilon}|}{|z_{x,\epsilon}|}\eta' (|z_{x,\epsilon}|) \right)\cdot \frac{(1-|z_{x,\epsilon}|^2)^2}{4} \cdot \mu_{\mathrm{KE}} \\
& \leq & \frac{1}{8\pi} \cdot\left( - \log |z_{x,\epsilon}| \cdot \eta''(|z_{x,\epsilon}|) - \frac{2+\log |z_{x,\epsilon}|}{|z_{x,\epsilon}|}\eta' (|z_{x,\epsilon}|) \right)\cdot \mu_{\mathrm{KE}} .
\end{eqnarray*}

From the bounds on $\eta'$ and $\eta''$,
\begin{equation*}
 |\log |z_{x,\epsilon}| \cdot \eta''(|z_{x,\epsilon}|)| \leq \frac{81}{4(\tanh\frac{\epsilon}{2})^2} \cdot \log \left( \frac{3}{\tanh\frac{\epsilon}{2}} \right) ,
\end{equation*}
and
\begin{eqnarray*}
\left| \frac{2+\log |z_{x,\epsilon}|}{|z_{x,\epsilon}|}\eta' (|z_{x,\epsilon}|) \right| & \leq & \sup_{\frac{\tanh\frac{\epsilon}{2}}{3}\leq|\varsigma|\leq \tanh\frac{\epsilon}{2}} \left| 2+\log \varsigma  \right| \cdot \frac{3}{\tanh\frac{\epsilon}{2}} \cdot \frac{9}{2 \tanh\frac{\epsilon}{2}} \\
& \leq & \frac{27}{2(\tanh\frac{\epsilon}{2})^2} \cdot \log \left( \frac{3}{\tanh\frac{\epsilon}{2}} \right) .
\end{eqnarray*}
Thus
\begin{eqnarray*}
\left| \log |z_{x,\epsilon}| \cdot \eta''(|z_{x,\epsilon}|) + \frac{2+\log |z_{x,\epsilon}|}{|z_{x,\epsilon}|}\eta' (|z_{x,\epsilon}|) \right| & \leq & \left( \frac{81}{4} + \frac{27}{2} \right) \cdot \frac{1}{(\tanh\frac{\epsilon}{2})^2} \cdot \log \left( \frac{3}{\tanh\frac{\epsilon}{2}} \right) \\
& = & \frac{135}{4} \cdot \frac{1}{(\tanh\frac{\epsilon}{2})^2} \cdot \log \left( \frac{3}{\tanh\frac{\epsilon}{2}} \right) ,
\end{eqnarray*}
which establishes (6) and completes the proof.
\end{proof}

We now turn to the thin part. In this case we work in a collar $\mathscr{C}(\gamma)$ around a short geodesic $\gamma$, endowed with the coordinate from Lemma \ref{lemkahlercollar}. The following lemma provides the analogue of Lemma \ref{lemconstructionlocalresidue_thick} in this setting.

\begin{lem}
\label{lemconstructionlocalresidue_thin}
Let $\gamma $ be a simple closed geodesic on $C$ of length $\ell \leq 2\arcsinh \frac{1}{2}$. Let $\mathscr{C}(\gamma)$ and $\nu(\gamma)$ be as in Lemma \ref{lemkahlercollar}, and let 
$$z_\gamma : \mathscr{C} (\gamma ) \to \overline{\mathbb{D}_{e^{-\nu (\gamma )}} (0)} \setminus \mathbb{D}_{e^{\nu (\gamma ) -\frac{2\pi^2}{\ell } }} (0) $$
be the holomorphic coordinate given there. Then
$ \left|\frac{\pi^2}{\ell} + \log \left( \Vert dz_\gamma \Vert_{\mathrm{hyp}} \right) \right| < \frac{\pi^2}{\ell} $.

Moreover, for any $x\in \mathscr{C} (\gamma) $ satisfying $|\log |z_\gamma (x)| + \frac{\pi^2}{\ell} |\leq \frac{\pi^2}{\ell} - 2\nu(\gamma)$, there exists a smooth function $ \digamma_{x,\gamma} : C\setminus \{x\} \to \mathbb{R}$ on $C\setminus \{x\}$, such that 
\begin{enumerate}[(1)]
    \item  $\digamma_{x,\gamma} (y)  = 0$ for all $y\in C \setminus \mathscr{C} (\gamma)$.
    \item The function $ \digamma_{x,\gamma} (y) +\log | z_\gamma (y) -z_\gamma (x) | $ extends smoothly to $\mathscr{C} (\gamma)$, and
    $$ \left| \digamma_{x,\gamma} (y) +\log | z_\gamma (y) - z_\gamma (x) | \right| \leq \frac{3\pi^2}{\ell} . $$
    \item The integral satisfies $\left| \int_C \digamma_{x,\gamma} (y) \mu_{\mathrm{KE}} (y) \right| \leq \frac{8\pi^2}{\ell} $.
    \item In the sense of currents,
    $$-600\cdot\mu_{\mathrm{KE}} \leq dd^c \digamma_{x,\gamma} +\delta_x \leq 600 \cdot \mu_{\mathrm{KE}} ,$$
    where $\delta_x$ is the Dirac measure at $x$.
\end{enumerate}
\end{lem}

\begin{proof}
We first estimate the norm $ \Vert\cdot\Vert_{\mathrm{hyp}} $ using the explicit expression of the metric $\mu_{\mathrm{KE}}$ in the collar $\mathscr{C}(\gamma)$. By Lemma \ref{lemkahlercollar}, we have
$$ \mu_{\mathrm{KE}} = \frac{ \ell^2 z_\gamma^* (\mu_{\mathrm{Euc}}) }{ 4\pi^2 |z_\gamma|^2 \sin^2 ( \frac{\ell }{2\pi} \log |z_\gamma | ) } ,$$
where $\mu_{\mathrm{Euc}} = \frac{idz\wedge d\bar{z}}{2} $ is the standard Euclidean metric on $\mathbb{C}$.

Set $\varsigma=\frac{\pi^2}{\ell}+\log|z_\gamma|$. Then $|\varsigma|\leq \frac{\pi^2}{\ell}-\nu(\gamma)$, and the relation between $dz_\gamma$ and $\mu_{\mathrm{KE}}$ yields
\begin{equation*}
\frac{\pi^2}{\ell} + \log(\Vert dz_\gamma \Vert_{\mathrm{hyp}}) = \log\frac{4\pi}{\ell} +\varsigma + \log \left( \cos \frac{\ell \varsigma }{2\pi} \right) .
\end{equation*}
A direct calculation shows
\begin{equation*}
\frac{d}{d\varsigma} \left( \varsigma + \log \left( \cos \frac{\ell \varsigma }{2\pi} \right) \right) = 1 - \frac{\ell}{2\pi} \tan \frac{\ell \varsigma }{2\pi} \geq 1-\frac{\ell}{2\pi} \cot \frac{\ell\nu(\gamma)}{2\pi} = 1-\frac{\ell}{2\pi\sinh \frac{\ell}{2}} >0 . 
\end{equation*}
and hence the function is strictly increasing in the admissible range of $\varsigma$.

Since $\ell\leq 2\arcsinh\frac{1}{2}<1$, $\frac{2\tanh\frac{\ell}{2}}{\ell} > 2\tanh\frac{1}{2} \approx 0.924\cdots >\frac{9}{10} $. It follows from Proposition \ref{propbasicpropertiesriemanniankahlercollars} and the calculations above that
\begin{eqnarray*}
\frac{\pi^2}{\ell} + \log(\Vert dz_\gamma \Vert_{\mathrm{hyp}}) & \leq & \log\frac{4\pi}{\ell} + \frac{\pi^2}{\ell} - \nu(\gamma) + \log \left( \sin \frac{\ell \nu (\gamma) }{2\pi} \right) \\
& < & \log 2\pi - \nu(\gamma) +\frac{\pi^2}{\ell} < \frac{\pi^2}{\ell} ,
\end{eqnarray*}
and
\begin{eqnarray*}
\frac{\pi^2}{\ell} + \log(\Vert dz_\gamma \Vert_{\mathrm{hyp}}) & \geq & \log\frac{4\pi}{\ell} - \frac{\pi^2}{\ell} + \nu(\gamma) + \log \left( \sin \frac{\ell \nu (\gamma) }{2\pi} \right) \\
& > & \log \frac{9\pi}{5} + \nu(\gamma) -\frac{\pi^2}{\ell} > -\frac{\pi^2}{\ell} .
\end{eqnarray*}
Hence $ \left|\frac{\pi^2}{\ell} + \log \left( \Vert dz_\gamma \Vert_{\mathrm{hyp}} \right) \right| < \frac{\pi^2}{\ell} $, which establishes the desired bound on $ \Vert\cdot\Vert_{\mathrm{hyp}} $.

We now turn to the construction of the function $\digamma_{x,\gamma}$. For any $x\in \mathscr{C} (\gamma) $ satisfying $|\log |z_\gamma (x)| + \frac{\pi^2}{\ell} |\leq \frac{\pi^2}{\ell} - 2\nu(\gamma)$, set $\digamma_{x,\mathrm{loc}} (y) = \log |z_\gamma (y)-z_\gamma (x)| + \frac{\ell  \log |z_\gamma (x)| }{2\pi^2} \log |z_\gamma (y)| $ on $\mathscr{C} (\gamma )\setminus \{x\}$. To extend $\digamma_{x,\mathrm{loc}}$ globally, we choose a cut-off function $\eta $ supported on $\mathscr{C}(\gamma )$ as follows. Let $\eta_{1 }:\mathbb{R}\to\mathbb{R}_{\geq 0} $ be a smooth function with 
$$\eta_{1 }(\varsigma)=1 \ \text{ for }\varsigma\leq e^{-\frac{3\nu(\gamma )}{2}}, \qquad \eta_{1 }(\varsigma)=0 \ \text{ for }\varsigma\geq e^{-\nu(\gamma ) },$$
and derivatives bounded by $|\eta'_{1 }| \leq 3e^{\nu(\gamma)} $, $|\eta''_{1 }| \leq 15e^{2\nu(\gamma_j)} $. Note that $e^{\frac{\nu (\gamma )}{2}}> e^{\frac{\pi^2}{8\arcsinh 1}}\approx 4.05\cdots >4$.

Set $\eta (y)= \eta_{1 } (|z_\gamma (y)|) \cdot \eta_{1 } \left( e^{-\frac{2\pi^2}{\ell }}|z_\gamma (y)|^{-1} \right)$, and define
\begin{eqnarray*}
\digamma_{x,\gamma} (y) & = & \left\{
\begin{aligned}
-\eta (y) \cdot \digamma_{x,\mathrm{loc}} (y), &\;\;\;  y\in \mathscr{C} (\gamma )\setminus \{x\} ;\\
0, \;\;\;\;\;\;\;\;\;\; &\;\;\; y\in C\setminus \mathscr{C} (\gamma ) .
\end{aligned}
\right.
\end{eqnarray*}
Clearly, the support of $d \eta $ is contained in
$$ \left\{ y\in\mathscr{C} (\gamma ): \frac{\pi^2}{\ell } - \frac{3\nu (\gamma )}{2} \leq \left| \frac{\pi^2}{\ell } +\log |z_\gamma (y)| \right| \leq \frac{\pi^2}{\ell } - \nu (\gamma ) \right\} .$$

The property (1) follows immediately from the support of $\eta$.  
For (2), note that by construction the function $ \digamma_{x,\gamma} (y) +\log | z_\gamma (y) -z_\gamma (x) | $ extends smoothly to $\mathscr{C}(\gamma)$. Moreover, a direct estimate yields
\begin{eqnarray*}
\left| \digamma_{x,\gamma} (y) +\log | z_\gamma (y) - z_\gamma (x) | \right| & \leq & \left| \digamma_{x,\mathrm{loc}} (y) +\log | z_\gamma (y) - z_\gamma (x) | \right| + |\eta (y) -1| \left| \digamma_{x,\mathrm{loc}} (y) \right| \\
& < & \frac{2\pi^2}{\ell} + |\eta (y) -1| \left| \digamma_{x,\mathrm{loc}} (y) \right| .
\end{eqnarray*}
When $y\in\mathscr{C} (\gamma)$ and $ -\frac{3 \nu (\gamma )}{2} \leq \log |z_\gamma (y)| \leq - \nu (\gamma ) $, we have
\begin{eqnarray*}
|\digamma_{x,\mathrm{loc}} (y)| & \leq & |\log |z_\gamma (y)||+\left|\log \left|1-\frac{z_\gamma (x)}{z_\gamma (y)}\right|\right| + \left|\frac{\ell \log |z_\gamma (x)| }{2\pi^2} \log |z_\gamma (y)| \right|\\
& \leq & \frac{3 \nu (\gamma )}{2} - \log \left(1-e^{-\frac{\nu (\gamma )}{2}}\right) + \frac{3\nu (\gamma )}{2} < 3\pi +\log\frac{4}{3} .
\end{eqnarray*}
Using the involution $(x,y)\mapsto (x',y')$ with $z_\gamma(x')z_\gamma(x)=z_\gamma(y')z_\gamma(y)=e^{-\frac{2\pi^2}{\ell}}$, we further obtain $\digamma_{x',\mathrm{loc}} (y')=\digamma_{x,\mathrm{loc}} (y)$, and thus $|\eta (y) -1| \left| \digamma_{x,\mathrm{loc}} (y) \right|\leq 3\pi +\log\frac{4}{3} $. Hence
\begin{equation*}
\left| \digamma_{x,\gamma} (y) +\log | z_\gamma (y) - z_\gamma (x) | \right|
< \frac{2\pi^2}{\ell} + 3\pi +\log\frac{4}{3} <\frac{3\pi^2}{\ell} ,
\end{equation*}
which establishes property (2).

We now proceed to estimate the integral in (3). Since $|\digamma_{x,\gamma }(y)|\leq |\digamma_{x,\mathrm{loc}}(y)| $, it suffices to perform the same type of estimate for $\digamma_{x,\mathrm{loc}}(y)$. By definition, we can decompose the integral as
\begin{eqnarray*}
\left| \int_{\mathscr{C} (\gamma)} \digamma_{x,\gamma} (y) d\mu_{\mathrm{KE}} (y) \right|& \leq & \left| \int_{\mathscr{C} (\gamma )} \frac{\ell \log |z_\gamma (x)| }{2\pi^2} \log |z_\gamma (y)| d\mu_{\mathrm{KE}} (y) \right| \\
& & + \left| \int_{ \left\{|z_\gamma (y)-z_\gamma (x)| \leq \frac{|z_\gamma(x)|}{10} \right\}} \log |z_\gamma (y)-z_\gamma (x)| d\mu_{\mathrm{KE}} (y) \right| \\
& & + \left| \int_{\left\{|z_\gamma (y)-z_\gamma (x)| \geq \frac{|z_\gamma(x)|}{10} \right\}} \log |z_\gamma (y)-z_\gamma (x)| d\mu_{\mathrm{KE}} (y) \right| .
\end{eqnarray*}

By symmetry we may assume $\log|z_\gamma (x)|\geq -\frac{\pi^2}{\ell}$. For the first term, a straightforward calculation gives
\begin{eqnarray*}
\left| \int_{\mathscr{C} (\gamma )} \frac{\ell \log |z_\gamma (x)| }{2\pi^2} \log |z_\gamma (y)| d\mu_{\mathrm{KE}} (y) \right| & \leq & \left| \int_{e^{\nu (\gamma ) - \frac{2\pi^2}{\ell }}}^{e^{-\nu (\gamma )}} \frac{\pi\varsigma\cdot \ell^2 \log\varsigma d\varsigma }{4\pi^2 \varsigma^2 \sin^2 (\frac{\ell }{2\pi} \log \varsigma ) } \right| \\
& = & \frac{\pi^2}{\sinh \frac{\ell}{2}} < \frac{2\pi^2}{\ell} .
\end{eqnarray*}

If $|z_\gamma (y)-z_\gamma (x)| \leq \frac{|z_\gamma(x)|}{10}$, then $|z_\gamma (y)| \geq \frac{9}{10}e^{-\frac{\pi^2}{\ell}}$, and it follows that $|\sin ( \frac{\ell }{2\pi} \log |z_\gamma (y)| )| \geq | \frac{\ell }{8\pi} \log |z_\gamma (x)| | $. A direct estimate then yields
\begin{eqnarray*}
 & & \left| \int_{ \left\{|z_\gamma (y)-z_\gamma (x)| \leq \frac{|z_\gamma(x)|}{10} \right\}} \log |z_\gamma (y)-z_\gamma (x)| d\mu_{\mathrm{KE}} (y) \right| \\
 & \leq & 16\cdot\frac{100}{81} \left| \int_{ \left\{|z_\gamma (y)-z_\gamma (x)| \leq \frac{|z_\gamma(x)|}{10} \right\}} \frac{ \log |z_\gamma (y)-z_\gamma (x)|  z_\gamma^* (\mu_{\mathrm{Euc}}) (y) }{ |z_\gamma (x)|^2 |\log |z_\gamma (x)||^2 } \right| \\
 & < & \frac{-20}{ |z_\gamma (x)|^2 |\log |z_\gamma (x)||^2 } \int_0^{\frac{|z_\gamma(x)|}{10}} 2\pi \varsigma \log\varsigma d\varsigma < 1 .
\end{eqnarray*}

For the last term, since $\mu_{\mathrm{KE}} (\mathscr{C}(\gamma)) = \frac{2\ell}{\sinh \frac{\ell}{2}} <4 $, we have
\begin{eqnarray*}
\left| \int_{ \left\{|z_\gamma (y)-z_\gamma (x)| \geq \frac{|z_\gamma(x)|}{10} \right\}} \log |z_\gamma (y)-z_\gamma (x)| d\mu_{\mathrm{KE}} (y) \right| & < & \left|\log \frac{|z_\gamma (x)|}{10} \right| \cdot\mu_{\mathrm{KE}} (\mathscr{C}(\gamma)) \\
& < & 4\left( \frac{\pi^2}{\ell}+\log 10 \right) .
\end{eqnarray*}

Combining these estimates, we obtain
\begin{equation*}
\left| \int_{\mathscr{C} (\gamma)} \digamma_{x,\gamma} (y) d\mu_{\mathrm{KE}} (y) \right| < \frac{2\pi^2}{\ell} + 1+4\left( \frac{\pi^2}{\ell}+\log 10 \right) < \frac{8\pi^2}{\ell} ,
\end{equation*}
which proves (3).

Finally, we establish property (4) for $\digamma_{x,\gamma}$. Recall that
$$ \mu_{\mathrm{KE}} = \frac{ \ell^2 z_\gamma^* (\mu_{\mathrm{Euc}}) }{ 4\pi^2 |z_\gamma|^2 \sin^2 ( \frac{\ell }{2\pi} \log |z_\gamma | ) } \geq \frac{ z_\gamma^* (\mu_{\mathrm{Euc}}) }{  |z_\gamma|^2 |\log |z_\gamma ||^2 } .$$
A straightforward calculation shows that on $\mathrm{supp} \left( d\left(\eta_1 (|z_\gamma|)\right) \right)$,
\begin{eqnarray*}
& & d \digamma_{x,\mathrm{loc}} (y) \wedge d^c \digamma_{x,\mathrm{loc}} (y) \\
& = & \frac{i}{\pi} \left|\frac{1}{2(z_\gamma (y)-z_\gamma (x))} + \frac{\ell \log |z_\gamma (x)| }{4\pi^2 z_\gamma (y)}\right|^2 dz_\gamma (y) \wedge d\bar{z}_\gamma (y) \\
& \leq & \frac{2}{\pi} \cdot \left|\frac{1}{2\left(1-e^{-\frac{\nu (\gamma)}{2}}\right) |z_\gamma (y)|} + \frac{1 }{2 |z_\gamma (y)|}\right|^2 z_\gamma^* (\mu_{\mathrm{Euc}}) \\
& \leq & \frac{49}{18\pi\cdot |z_\gamma (y)|^2} \cdot z_\gamma^* (\mu_{\mathrm{Euc}}) \leq \frac{49}{18\pi } \cdot \frac{9\nu^2 (\gamma)}{4} \cdot \mu_{\mathrm{KE}} < \frac{49\pi }{8} \cdot \mu_{\mathrm{KE}}.
\end{eqnarray*}

Next, we estimate $d\eta$ and $dd^c \eta$. By symmetry, it suffices to consider $\eta_{1} (|z_\gamma|)$.

Since $|\eta'_{1}|\leq 3e^{\nu (\gamma)} $ and $|\eta''_{1}|\leq 15 e^{2\nu (\gamma)} $, we have
\begin{eqnarray*}
d\eta_{1} (|z_j (y)|)\wedge d^c \eta_{1} (|z_j (y)|) & = & |\eta'_{1}(|z_j (y)|)|^2 \frac{idz_j (y) \wedge d\bar{z}_j (y)}{4\pi} \\
& < & 9e^{2\nu (\gamma_j)} \cdot \frac{|z_\gamma|^2 |\log |z_\gamma ||^2 }{2\pi} \cdot \mu_{\mathrm{KE}} < \frac{9\pi}{2} \cdot \mu_{\mathrm{KE}} ,
\end{eqnarray*}
and
\begin{eqnarray*}
dd^c\eta_{1,1} (|z_j (y)|) & = & \left( \eta''_{1,1}(|z_j (y)|) + \frac{\eta'_{1,1}(|z_j (y)|)}{|z_j (y)|} \right) \frac{idz_j (y) \wedge d\bar{z}_j (y)}{4\pi} \\
& \leq & \left( 15 e^{2\nu (\gamma_j)} + \frac{3e^{\nu (\gamma_j)}}{|z_j (y)|} \right) \cdot \frac{z_\gamma^* (\mu_{\mathrm{Euc}})}{2\pi} \\
& < & \frac{1}{2\pi} \cdot \left( 15 \pi^2 + 3\cdot\frac{9\pi^2}{4} \right) \cdot \mu_{\mathrm{KE}} = \frac{87\pi}{8} \cdot \mu_{\mathrm{KE}} .
\end{eqnarray*}

Combining these estimates for $d\digamma_{x,\mathrm{loc}}$, $d\eta_1$, and $dd^c \eta_1$, we obtain
\begin{eqnarray*}
dd^c \digamma_{x,\gamma} +\delta_x & = & - \digamma_{x,\mathrm{loc}} \cdot dd^c \eta_1 - d\eta_1 \wedge d^c \digamma_{x,\mathrm{loc}} -d\digamma_{x,\mathrm{loc}} \wedge d^c \eta_1 \\
& \leq & \left( 3\pi+\log\frac{4}{3} \right)\cdot \frac{87\pi}{8} \cdot \mu_{\mathrm{KE}} + 2\cdot \sqrt{\frac{49\pi}{8} \cdot\frac{9\pi }{2}} \cdot \mu_{\mathrm{KE}} < 600\cdot \mu_{\mathrm{KE}} . 
\end{eqnarray*}
By the same reasoning, the lower bound also holds: $dd^c \digamma_{x,\gamma} +\delta_x > -600\cdot \mu_{\mathrm{KE}}$. This completes the proof.
\end{proof}

Combining Lemma \ref{lemconstructionlocalresidue_thick} and Lemma \ref{lemconstructionlocalresidue_thin} yields a global construction on $C\times C$.

\begin{lem}
\label{lemconstructionlocalresidue}
Let $C$ be a curve of genus $g\geq 2$ over $\mathbb{C}$, and let $\mu_{\mathrm{KE}}$ denote the unique K\"ahler metric on $C$ with constant curvature $-1$. Let $\mathrm{sys}(C)$ denote the systole of $(C,\mu_{\mathrm{KE}})$. Then there exists a Borel measurable $\mathbb{R}$-valued function $\digamma$ on $C\times C$, such that for every $x\in C$ we can choose a local coordinate
$$ z_x:U_x\to\mathbb{C},\qquad U_x=B_{\frac{1}{2}\arcsinh \frac{1}{2}}(x)=\left\{y\in C:\dist_{\mu_{\mathrm{KE}}} (x,y)<\frac{1}{2}\arcsinh \frac{1}{2}\right\} ,$$
with the following properties:
\begin{enumerate}[(1)]
    \item For any $x\in C$, $\left| \int_C \digamma_{x} (y) \mu_{\mathrm{KE}} (y) \right| \leq \frac{8\pi^2\cdot \max\{1,\mathrm{sys}(C)\}}{\mathrm{sys}(C)} $, where $\digamma_{x} (y)=\digamma (x,y)$.
    \item For any $x\in C$, in the sense of currents,
    $$-600\cdot\mu_{\mathrm{KE}} \leq dd^c \digamma_{x} +\delta_x \leq 600 \cdot \mu_{\mathrm{KE}} ,$$
    where $\delta_x$ is the Dirac measure at $x$.
    \item For any $x\in C$, the function $ \digamma_{x} (y) +\log | z_x (y) -z_x (x) | $ extends smoothly to $U_x$.
    \item If $\dist_{\mu_{KE}} (x,y)\geq \frac{1}{100}$, then $|\digamma (x ,y)  | \leq \frac{6\pi^2\cdot \max\{1,\mathrm{sys}(C)\}}{\mathrm{sys}(C)} $.
    \item For any $x\in C$ and $(y_1 ,y_2)\in U_x\times U_x$, $|\digamma (y_1 ,y_2) +\log | z_x (y_1) -z_x (y_2) || \leq \frac{8\pi^2\cdot \max\{1,\mathrm{sys}(C)\}}{\mathrm{sys}(C)} $.
\end{enumerate}
\end{lem}

\begin{proof}
Let $\{\gamma_j\}_{j=1}^m$ be the set of all simple closed geodesics on $(C,\mu_{\mathrm{KE}} )$ of length $ \leq 2\arcsinh \frac{1}{2} $. By Lemma \ref{lemkahlercollar}, there exist holomorphic coordinates
$$z_j : \mathscr{C} (\gamma_j ) \to \overline{\mathbb{D}_{e^{-\nu (\gamma_j)}} (0)} \setminus \mathbb{D}_{e^{\nu (\gamma_j ) -\frac{2\pi^2}{\ell (\gamma_j )} }} (0) $$
such that 
$$ \mu_{\mathrm{KE}} = \frac{ \ell^2 (\gamma_j ) z_j^* (\mu_{\mathrm{Euc}}) }{ 4\pi^2 |z_j|^2 \sin^2 ( \frac{\ell (\gamma_j )}{2\pi} \log |z_j| ) } ,$$
where $\mu_{\mathrm{Euc}} = \frac{idz\wedge d\bar{z}}{2} $ is the standard Euclidean metric on $\mathbb{C}$, and $\nu (\gamma_j) \in (0, \frac{\pi^2}{2\ell (\gamma_j)}) $ is chosen such that $\nu (\gamma_j) \ell (\gamma_j) = {2\pi} \arcsin \left( \tanh \frac{\ell (\gamma_j )}{2} \right) $. 

For any $x\in C$ with $\mathrm{inj}_{\mu_{\mathrm{KE}}} (x)\leq \arcsinh \frac{1}{2}$, Proposition \ref{propbasicpropertiesriemanniankahlercollars} shows that $x\in\mathscr{C} (\gamma_j ) $ for some $\gamma_j$, and $ \left| \frac{\pi^2}{\ell (\gamma_j )} +\log |z_j (x)| \right| \leq \frac{\pi^2}{\ell (\gamma_j )} - 2 \nu (\gamma_j ) $. Note that the collars $\mathscr{C}(\gamma_j)$ are disjoint, so the index $j$ is unique when $x\in\bigcup_{j=1}^m\mathscr{C}(\gamma_j)$. We then define
\begin{eqnarray*}
\digamma (x,y) & = & \left\{
\begin{aligned}
\digamma_{x,\gamma_j} (y) &,\;\;\; \text{if } x\in\mathscr{C}(\gamma_j)\ \text{ and } \left| \frac{\pi^2}{\ell (\gamma_j )} +\log |z_j (x)| \right| \leq \frac{\pi^2}{\ell (\gamma_j )} - 2 \nu (\gamma_j ),\\
\digamma_{x,\arcsinh \frac{1}{2}} (y) &,\;\;\; \textrm{ otherwise,} 
\end{aligned}
\right.
\end{eqnarray*}
where $\digamma_{x,\epsilon}$ and $\digamma_{x,\gamma}$ are as in Lemma \ref{lemconstructionlocalresidue_thick} and Lemma \ref{lemconstructionlocalresidue_thin}, respectively. The local coordinate $z_x$ is defined in the same way, by taking $z_\gamma$ in the first case and $z_{x,\arcsinh \frac{1}{2}}$ in the second, where $z_{x,\epsilon}$ and $z_\gamma$ are the coordinates in Lemma \ref{lemconstructionlocalresidue_thick} and Lemma \ref{lemconstructionlocalresidue_thin}, respectively.

By construction, (1)-(3) follow from Lemma \ref{lemconstructionlocalresidue_thick} and Lemma \ref{lemconstructionlocalresidue_thin}. 

We now turn to (4). The argument for (4) splits into two cases. First, suppose $\digamma_x (y) = \digamma_{x,\arcsinh \frac{1}{2}} (y)$. By Lemma \ref{lemconstructionlocalresidue_thick}, for any $y\in B_{\arcsinh \frac{1}{2}} (x)$, we have
$$\left| \digamma_{x} (y) +\log | z_{x,\arcsinh \frac{1}{2}}(y) | \right| \leq \log \left( \frac{3}{\tanh \frac{\arcsinh \frac{1}{2}}{2} } \right),$$
and $\digamma_{x} (y)=0$ otherwise. Hence, if $\dist_{\mu_{KE}} (x,y)\geq \frac{1}{100}$, then
\begin{eqnarray*}
\left| \digamma_{x} (y) \right| & \leq & \left| \log | z_{x,\arcsinh \frac{1}{2}}(y) | \right| + \log \left( \frac{3}{\tanh \frac{\arcsinh \frac{1}{2}}{2} } \right) \\
& \leq & \log \left( \frac{1}{\tanh \frac{1}{200} } \right) +\log \left( \frac{3}{\tanh \frac{\arcsinh \frac{1}{2}}{2} } \right) <10 .
\end{eqnarray*}

Next, suppose $\digamma_x (y) = \digamma_{x,\gamma_j} (y)$ for some $j$. By Lemma \ref{lemconstructionlocalresidue_thin}, for any $y\in \mathscr{C}(\gamma_j)$,
$$\left| \digamma_{x} (y) +\log | z_{\gamma_j}(y) -z_{\gamma_j}(x) | \right| \leq \frac{3\pi^2}{\ell (\gamma_j)},$$
and $\digamma_{x} (y)=0$ otherwise. Recall that on $\mathscr{C}(\gamma_j)$,
$$ \mu_{\mathrm{KE}} = \frac{ \ell^2 (\gamma_j) z_{\gamma_j}^* (\mu_{\mathrm{Euc}}) }{ 4\pi^2 |z_{\gamma_j}|^2 \sin^2 ( \frac{\ell (\gamma_j) }{2\pi} \log |z_{\gamma_j} | ) } \leq \frac{ \ell^2 (\gamma_j) z_{\gamma_j}^* (\mu_{\mathrm{Euc}}) }{ 4\pi^2 |z_{\gamma_j}|^2 \tanh^2 ( \frac{\ell (\gamma_j)}{2} ) } \leq e^{\frac{4\pi^2}{\ell(\gamma_j)}} z_{\gamma_j}^* (\mu_{\mathrm{Euc}}) .$$
Thus $\dist_{\mu_{KE}} (x,y)\geq \frac{1}{100}$ implies $| z_{\gamma_j}(y) - z_{\gamma_j}(x) |\geq \frac{1}{100} e^{-\frac{2\pi^2}{\ell (\gamma_j)}}$, and hence
\begin{equation*}
\left| \digamma_{x} (y) \right| \leq \left| \log | z_{\gamma_j}(y) - z_{\gamma_j}(x) | \right| + \frac{3\pi^2}{\ell (\gamma_j)} \leq \log 100+ \frac{2\pi^2}{\ell (\gamma_j)} +\frac{3\pi^2}{\ell (\gamma_j)} < \frac{6\pi^2}{\ell (\gamma_j)} .
\end{equation*}
Combining the two cases, we conclude that whenever $\dist_{\mu_{KE}} (x,y)\geq \frac{1}{100}$, we have $|\digamma (x ,y)  | \leq \frac{6\pi^2\cdot \max\{1,\mathrm{sys}(C)\}}{\mathrm{sys}(C)} $. Thus, property (4) follows.

It remains to establish (5). Similar to the argument in (4), for any $x\in C$ and $(y_1,y_2)\in U_x\times U_x$, we have
$$\left| \digamma(y_1 ,y_2) +\log | z_{y_1} (y_1) - z_{y_1} (y_2) | \right| \leq \frac{3\pi^2\cdot \max\{1,\mathrm{sys}(C)\}}{\mathrm{sys}(C)} ,$$
as follows directly from Lemma \ref{lemconstructionlocalresidue_thick} and Lemma
\ref{lemconstructionlocalresidue_thin}. What is left is to show that 
$$\left| \log | z_{x} (y_1) - z_{x} (y_2) |-  \log | z_{y_1} (y_1) - z_{y_1} (y_2) | \right| \leq \frac{5\pi^2\cdot \max\{1,\mathrm{sys}(C)\}}{\mathrm{sys}(C)} .$$

First consider the case $z_x = z_{x,\arcsinh \frac{1}{2}} $. In this setting the geodesic segment joining $y_1,y_2$ lies in $U_x$. By Theorem \ref{thmpoincaremodel}, on $z_x (U_x)$, we have $4\mu_{\mathrm{Euc}}\leq\mu_{\mathbb{D}}\leq 8\mu_{\mathrm{Euc}}$. Hence
$$\frac{\dist_{\mu_{KE}} (y_1 ,y_2)}{2\sqrt{2}} = \frac{\dist_{\mu_{\mathbb{D}}} (z_x(y_1) ,z_x(y_2))}{2\sqrt{2}} \leq |z_x(y_1)-z_x(y_2)| \leq \frac{\dist_{\mu_{KE}} (y_1 ,y_2)}{2} ,$$
which gives
$$\left| \log | z_{x} (y_1) - z_{x} (y_2) |-  \log \dist_{\mu_{KE}} (y_1 ,y_2) \right| \leq \frac32\log 2 .$$

On the other hand, if $z_x =z_{\gamma_j} $ for some $j$, then the geodesic segment joining $y_1,y_2$ lies in $\mathscr{C}(\gamma_j)$. A direct computation yields
$$ e^{\frac{4\pi^2}{\ell(\gamma_j)}} z_{\gamma_j}^* (\mu_{\mathrm{Euc}}) \geq \mu_{\mathrm{KE}} = \frac{ \ell^2 (\gamma_j) z_{\gamma_j}^* (\mu_{\mathrm{Euc}}) }{ 4\pi^2 |z_{\gamma_j}|^2 \sin^2 ( \frac{\ell (\gamma_j) }{2\pi} \log |z_{\gamma_j} | ) } \geq \frac{ \ell^2 (\gamma_j) z_{\gamma_j}^* (\mu_{\mathrm{Euc}}) }{ 4\pi^2 e^{-2\nu (\gamma_j)} } \geq e^{-\frac{\pi^2}{\ell(\gamma_j)}} z_{\gamma_j}^* (\mu_{\mathrm{Euc}}) ,$$
so that
$$e^{-\frac{\pi^2}{\ell(\gamma_j)}} \cdot\dist_{\mu_{KE}} (y_1 ,y_2) \leq |z_x(y_1)-z_x(y_2)| \leq e^{\frac{2\pi^2}{\ell(\gamma_j)}}\cdot \dist_{\mu_{KE}} (y_1 ,y_2) .$$
It follows that
$$\left| \log | z_{x} (y_1) - z_{x} (y_2) |-  \log \dist_{\mu_{KE}} (y_1 ,y_2) \right| \leq \frac{2\pi^2}{\ell (\gamma_j)} .$$

We now return to (5).

If both coordinate charts are thin–part charts, i.e. $z_x = z_{\gamma_{j_1}} $ and $z_{y_1} =z_{\gamma_{j_2}}$ for some $j_1,j_2 $, then the assumption $\mathrm{inj}_{\mu_{\mathrm{KE}}} (x) \leq \arcsinh \frac{1}{2}$ together with $y_1\in U_x$ implies that $j_1=j_2$, and hence
$$\left| \log | z_{x} (y_1) - z_{x} (y_2) |-  \log | z_{y_1} (y_1) - z_{y_1} (y_2) | \right|=0 .$$

Otherwise, at least one of the coordinate charts $z_x$ or $z_{y_1}$ 
comes from the thick part, that is, $z_x = z_{x,\arcsinh \frac{1}{2}} $ or $z_{y_1} = z_{y_1,\arcsinh \frac{1}{2}} $ (and possibly both). Combining the two estimates above, we obtain
\begin{eqnarray*}
\left| \log | z_{x} (y_1) - z_{x} (y_2) |-  \log | z_{y_1} (y_1) - z_{y_1} (y_2) | \right| & \leq & \left| \log | z_{x} (y_1) - z_{x} (y_2) |-  \log \dist_{\mu_{KE}} (y_1 ,y_2) \right| \\
& & + \left| \log | z_{y_1} (y_1) - z_{y_1} (y_2) |-  \log \dist_{\mu_{KE}} (y_1 ,y_2) \right| \\
& \leq & \max\left\{ \frac{2\pi^2}{\mathrm{sys}(C)} + \frac32\log2 ,3\log2 \right\} \\
& < & \frac{5\pi^2\cdot \max\{1,\mathrm{sys}(C)\}}{\mathrm{sys}(C)} .
\end{eqnarray*}
Thus property (5) is verified, and hence the lemma is proved.
\end{proof}

Before comparing the Arakelov hermitian metric $\Vert\cdot\Vert_{\mathrm{Ar}}$ with the hyperbolic hermitian metric $\Vert\cdot \Vert_{\mathrm{hyp}}$, we need the following estimate.

\begin{lem}
\label{lemsupnormpotential}
Let $u:C\to \mathbb{R}$ be a smooth function such that $\Delta_{\mathrm{hyp}} u \geq -\upsilon_1$ and $|\Delta_{\mathrm{hyp}} u| \leq \upsilon_2 $, for some constants $\upsilon_1,\upsilon_2>0$. Then the following inequalities hold:
\begin{eqnarray*}
 \int_C du\wedge d^c u & \leq & 5\upsilon_1 \cdot \upsilon_2 + \upsilon_1^2 \cdot 1232\pi^3 (g-1) (2g-3) \cdot\min\left\{ \frac{\varphi (C)}{\xi_2 (g)} , \frac{14(g-1)\cdot \max\{1,\mathrm{sys}(C)\}}{55\cdot \mathrm{sys}(C)} \right\} \\
 & \leq & \upsilon_1 \cdot \left( \frac{1232\pi^3 (g-1) (2g-3)}{\xi_2 (g)} \cdot \upsilon_1 \varphi (C) + 5 \cdot \upsilon_2 \right) ,
\end{eqnarray*}
and
\begin{eqnarray*}
\sup_C\left| u-\int_C u \mu_{\mathrm{hyp}} \right| & \leq & \sqrt{\upsilon_1\upsilon_2 \cdot \left(\frac{1}{2}+ \frac{1}{\mathrm{sys}(C)} \right) \cdot \frac{2000}{g-1} } \\
& & + \upsilon_1 \cdot 1232\pi^3 (g-1) (2g-3) \cdot\min\left\{ \frac{\varphi (C)}{\xi_2 (g)} , \frac{14(g-1)\cdot \max\{1,\mathrm{sys}(C)\}}{55\cdot \mathrm{sys}(C)} \right\} \\
& \leq & \sqrt{\left(\frac{1000}{\xi_1 (g)}+ \frac{2000}{\xi_2 (g)} \right) \cdot\frac{\upsilon_1\upsilon_2 \varphi (C)}{g-1}} + \frac{1232\pi^3 (g-1) (2g-3)}{\xi_2 (g)} \cdot \upsilon_1 \varphi (C) ,
\end{eqnarray*}
where $\xi_1(g)$, $\xi_2 (g)$ are the positive constants defined in Theorem \ref{thmglobalestimatezhangphiinvariant}.
\end{lem}

\begin{proof}
We expand $u$ into its $L^2$-orthonormal eigenfunction decomposition:
$$u= a_0 + \sum_{j=1}^{\infty} a_j \phi_{\mathrm{hyp} ,j} ,$$ 
where $a_j$ are constants, and $\phi_{\mathrm{hyp} ,j}$ are eigenfunctions of the hyperbolic Laplacian with eigenvalues $\lambda_{\mathrm{hyp},j}$.

By definition,
\begin{eqnarray*}
\sum_{j=1}^\infty |a_j|^2 \lambda^2_{\mathrm{hyp} ,j} & = & \int_C |\Delta_{\mathrm{hyp}} u|^2 d\mu_{\mathrm{hyp}} = \int_C \Delta_{\mathrm{hyp}} u \left( dd^c u + \upsilon_1\mu_{\mathrm{hyp}} \right) \\
& \leq & \sup_{C} |\Delta_{\mathrm{hyp}} u| \int_C  \left( dd^c u + \upsilon_1 \mu_{\mathrm{hyp}} \right) \leq \upsilon_1 \cdot \upsilon_2 .
\end{eqnarray*}
Similarly, for each $j$ we have
$$ |a_j|\lambda_{\mathrm{hyp} ,j} = \left| \int_C \phi_{\mathrm{hyp} ,j} dd^c u \right| \leq \upsilon_1 \sup_{C} \left| \phi_{\mathrm{hyp} ,j} \right| .$$

Following Lemma \ref{lemsumsmalleigenvalueestimatephiinvariant}, we obtain
\begin{eqnarray*}
    \sum_{0<\lambda_{\mathrm{hyp},l}\leq \frac{g-1}{5}} \frac{\sup_{C} \left| \phi_{\mathrm{hyp} ,j} \right|^2 }{\lambda_{\mathrm{hyp} ,j}} & = & \sum_{0<\lambda_{\mathrm{KE},l} \leq \frac{1}{20\pi}} \frac{ \sup\limits_{ C} |\phi_{\mathrm{KE} , l } |^2}{ \lambda_{\mathrm{KE},l}} \\
    & < & 1232\pi^3 (g-1) (2g-3) \cdot\min\left\{ \frac{\varphi (C)}{\xi_2 (g)} , \frac{14(g-1)\cdot \max\{1,\mathrm{sys}(C)\}}{55\cdot \mathrm{sys}(C)} \right\} .
\end{eqnarray*}
Therefore,
\begin{eqnarray*}
\int_C du\wedge d^c u & \leq & \frac{5}{g-1}\sum_{\lambda_{\mathrm{hyp},l}\geq \frac{g-1}{5}} |a_j|^2 \lambda^2_{\mathrm{hyp} ,j} + \sum_{0<\lambda_{\mathrm{hyp},l}\leq \frac{g-1}{5}} |a_j|^2 \lambda_{\mathrm{hyp} ,j} \\
& \leq & 5\upsilon_1 \cdot \upsilon_2 + \upsilon_1^2 \cdot \sum_{0<\lambda_{\mathrm{hyp},l}\leq \frac{g-1}{5}} \frac{\sup_{C} \left| \phi_{\mathrm{hyp} ,j} \right|^2 }{\lambda_{\mathrm{hyp} ,j}} \\
& \leq & 5\upsilon_1 \cdot \upsilon_2 + \upsilon_1^2 \cdot 1232\pi^3 (g-1) (2g-3) \cdot\min\left\{ \frac{\varphi (C)}{\xi_2 (g)} , \frac{14(g-1)\cdot \max\{1,\mathrm{sys}(C)\}}{55\cdot \mathrm{sys}(C)} \right\} .
\end{eqnarray*}

For the sup-norm of $u$, we have
\begin{eqnarray*}
\sup_C\left| u-\int_C u \mu_{\mathrm{hyp}} \right| & \leq & \left| \sum_{\lambda_{\mathrm{hyp},l}\geq \frac{g-1}{5}} a_j \phi_{\mathrm{hyp} ,j} \right| + \sum_{0<\lambda_{\mathrm{hyp},l}\leq \frac{g-1}{5}} \left| a_j \phi_{\mathrm{hyp} ,j} \right| \\
& \leq & \sqrt{\left( \sum_{\lambda_{\mathrm{hyp},l}\geq \frac{g-1}{5}} |a_j|^2 \lambda^2_{\mathrm{hyp} ,j} \right)\left( \sum_{\lambda_{\mathrm{hyp},l}\geq \frac{g-1}{5}} \lambda^{-2}_{\mathrm{hyp} ,j} \phi^2_{\mathrm{hyp} ,j} \right)} \\
& & + \upsilon_1 \cdot \sum_{0<\lambda_{\mathrm{hyp},l}\leq \frac{g-1}{5}} \frac{\sup_{C} \left| \phi_{\mathrm{hyp} ,j} \right|^2 }{\lambda_{\mathrm{hyp} ,j}} \\
& \leq & \sqrt{\upsilon_1\upsilon_2 \cdot \left( \sum_{\lambda_{\mathrm{hyp},l}\geq \frac{g-1}{5}} \lambda^{-2}_{\mathrm{hyp} ,j} \phi^2_{\mathrm{hyp} ,j} \right)} \\
& & + \upsilon_1 \cdot 1232\pi^3 (g-1) (2g-3) \cdot\min\left\{ \frac{\varphi (C)}{\xi_2 (g)} , \frac{14(g-1)\cdot \max\{1,\mathrm{sys}(C)\}}{55\cdot \mathrm{sys}(C)} \right\} .
\end{eqnarray*}

It remains to estimate the sum $\sum\limits_{\lambda_{\mathrm{hyp},l}\geq \frac{g-1}{5}} \lambda^{-2}_{\mathrm{hyp} ,j} \phi^2_{\mathrm{hyp} ,j}$. By a straightforward calculation,
\begin{eqnarray*}
\sum_{\lambda_{\mathrm{hyp},l} \geq \frac{g-1}{5}} \frac{\left| \phi_{\mathrm{hyp} , l } (x) \right|^2}{\lambda^2_{\mathrm{hyp},l}} & = & \int_{0}^{\frac{1}{80(g-1)}}\sum_{\lambda_{\mathrm{hyp},l} \geq \frac{g-1}{5}} \frac{e^{-t\lambda_{\mathrm{hyp},l} }}{\lambda_{\mathrm{hyp},l}\big( 1-e^{-\frac{\lambda_{\mathrm{hyp},l}}{80(g-1)}} \big)}\left| \phi_{\mathrm{hyp} , l } (x) \right|^2 dt\\
& \leq & \frac{1}{1-e^{-\frac{1}{400}}} \int_{0}^{\frac{1}{80(g-1)}} \sum_{\lambda_{\mathrm{hyp},l} \geq \frac{g-1}{5}} \frac{e^{-t\lambda_{\mathrm{hyp},l} }}{\lambda_{\mathrm{hyp},l}} \left| \phi_{\mathrm{hyp} , l } (x) \right|^2 dt .
\end{eqnarray*}
As in the proofs of Lemma \ref{lemarakelovgreenhyperbolicdatathreeparts}, Lemma \ref{lemsupnormhyperbolicheatkernel} and Lemma \ref{lemhyperbolicgreenlargeeigenvalue}, for any $t\in \left( 0,\frac{1}{80(g-1)} \right] $, we have
\begin{eqnarray*}
  \sum_{\lambda_{\mathrm{hyp},l} \geq \frac{g-1}{5}} \frac{e^{-t\lambda_{\mathrm{hyp},l} }}{\lambda_{\mathrm{hyp},l}} \left| \phi_{\mathrm{hyp} , l } (x) \right|^2 & \leq & -\frac{ \log \big(40(g-1)t \big)}{2 } + 101 + \frac{200 }{\mathrm{sys}(C)} .
\end{eqnarray*}
Hence
\begin{eqnarray*}
\sum_{\lambda_{\mathrm{hyp},l} \geq \frac{g-1}{5}} \frac{\left| \phi_{\mathrm{hyp} , l } (x) \right|^2}{\lambda^2_{\mathrm{hyp},l}} & \leq & \frac{1}{(g-1)\big(1-e^{-\frac{1}{400}}\big)} \left( \frac{ \log 2+1}{160 } + \frac{101}{80 }+ \frac{5 }{2\cdot\mathrm{sys}(C)} \right) \\
& \leq & \left(\frac{1}{2}+ \frac{1}{\mathrm{sys}(C)} \right) \cdot \frac{2000}{g-1} .
\end{eqnarray*}
This completes the proof.
\end{proof}

Consequently, we can estimate the difference between the Arakelov Green function and $\digamma_{x}$.

\begin{cor}
\label{corollarydiatanceArakelovGreen_Log}
Let $G_{\mathrm{Ar}}$ be the Arakelov Green function on $C$. Let $\digamma_x$ be the function constructed in Lemma \ref{lemconstructionlocalresidue}. Then we have
\begin{eqnarray*}
\sup_{y\in C} |G_{\mathrm{Ar}}(x,y)-\digamma_x (y)| \leq \frac{(2958032\pi^4 + 3696\pi^3) (g-1)^2 (2g-\frac{5}{2})}{\xi_2 (g)} \cdot \varphi (C) ,
\end{eqnarray*}
where $\xi_2 (g)$ is the constant from Theorem \ref{thmglobalestimatezhangphiinvariant}.
\end{cor}

\begin{proof}
By Lemma \ref{lemconstructionlocalresidue},
$$ -2400\pi \cdot (g-1) \cdot\mu_{\mathrm{hyp}} < dd^c \digamma_x +\delta_x (y) < 2400\pi \cdot (g-1) \cdot\mu_{\mathrm{hyp}} .$$
Note that $\mu_{\mathrm{KE}} = 4\pi (g-1)\mu_{\mathrm{hyp}}$.

Set $\mathbf{F}_x (y) = G_{\mathrm{hyp}} (x,y) - \digamma_x (y) $. By construction, $\mathbf{F}_x (y)$ is smooth and satisfies
$$ -2401\pi \cdot (g-1) \cdot\mu_{\mathrm{hyp}} < dd^c \mathbf{F}_x (y) < 2401\pi \cdot (g-1) \cdot\mu_{\mathrm{hyp}} .$$

Applying \ref{lemsupnormpotential} and Theorem \ref{thmglobalestimatezhangphiinvariant}, we obtain
\begin{eqnarray*}
 \sup_C\left| \mathbf{F}_x \right| & \leq & \sup_C\left| \mathbf{F}_x-\int_C \mathbf{F}_x d\mu_{\mathrm{hyp}} \right| + \left| \int_C G_{\mathrm{hyp}} (x,y) d\mu_{\mathrm{hyp}} (y) \right| + \left| \int_C \digamma_x d\mu_{\mathrm{hyp}} \right| \\
& \leq & 2401\pi \cdot\varphi (C)\cdot \sqrt{\left(\frac{1000}{\xi_1 (g)}+ \frac{2000}{\xi_2 (g)} \right)\cdot\frac{g-1}{ \xi_1 (g)}}  \\
& & + \frac{2958032\pi^4 (g-1)^2 (2g-3)}{\xi_2 (g)} \cdot \varphi (C) + 8\pi^2 \left( \frac{1}{\xi_1(g)} + \frac{1}{\xi_2(g)} \right) \cdot \varphi (C) \\
& < & \frac{2958032\pi^4 (g-1)^2 (2g-3) +30000\pi}{\xi_2 (g)} \cdot \varphi (C) ,
\end{eqnarray*}
where $\xi_1(g)$ and $\xi_2(g)$ are the constants from 
Theorem \ref{thmglobalestimatezhangphiinvariant}, and we have used the inequalities $ \frac{\varphi (C)}{\xi_1 (g)}\geq 1$ and $\xi_1 \geq 20g\cdot \xi_2$.

Next, let $\psi_{\mathrm{Ar}}$ be the function defined in Theorem \ref{thmdifferencearakelovgreenhyperbolicgreen}, so that $G_{\mathrm{Ar}}(x,y)-G_{\mathrm{hyp}}(x,y)=\psi_{\mathrm{Ar}}(x)+\psi_{\mathrm{Ar}}(y)$. By Lemma \ref{lmmbergmankernelfunctionupperbound},
$$ \frac{\mu_{\mathrm{Ar}}}{\mu_{\mathrm{hyp}}} \leq (\sqrt{2}+1)^2 \left( 1+\frac{2\arcsinh 1}{\mathrm{sys}(C)}\right) <11\left( 1+\frac{1}{\mathrm{sys}(C)}\right) ,$$
where $\mathrm{sys}(C)$ is the systole of $(C,\mu_{\mathrm{KE}})$. Since $ - \mu_{\mathrm{hyp}} < dd^c \psi_{\mathrm{Ar}} < \mu_{Ar} $ and $2\int_C \psi_{\mathrm{Ar}} \mu_{\mathrm{hyp}} = \int_C d \psi_{\mathrm{Ar}} \wedge d^c \psi_{\mathrm{Ar}} $, it follows that
\begin{eqnarray*}
& & \sup_C\left| \psi_{\mathrm{Ar}} \right| \leq \sup_C\left| \psi_{\mathrm{Ar}}-\int_C \psi_{\mathrm{Ar}} d\mu_{\mathrm{hyp}} \right| +\frac{1}{2} \left| \int_C d \psi_{\mathrm{Ar}} \wedge d^c \psi_{\mathrm{Ar}} \right| \\
& \leq & \sqrt{ \left(\frac{1000}{\xi_1 (g)}+ \frac{2000}{\xi_2 (g)} \right) \cdot\frac{11+11\mathrm{sys}(C)}{(g-1)\mathrm{sys}(C)}\cdot \varphi (C) } \\
& & +\frac{1232\pi^3 (g-1)^2 (2g-3)}{\xi_2 (g)} \cdot \varphi (C) + \frac{616\pi^3 (g-1)^2 (2g-3)}{\xi_2 (g)} \cdot \varphi (C) \\
& \leq & \frac{1848\pi^3 (g-1)^2 (2g-3) + 200}{\xi_2 (g)} \cdot \varphi (C) ,
\end{eqnarray*}
where we used $1\leq \frac{\varphi (C)}{\xi_1 (g)}$, $\frac{1}{\mathrm{sys}(C)}\leq \frac{\varphi (C)}{\xi_2 (g)}$ and $\xi_1 \geq 20g\cdot \xi_2$.

Finally, combining the bounds for $\mathbf{F}_x$ and $\psi_{\mathrm{Ar}}$, we deduce
\begin{eqnarray*}
|G_{\mathrm{Ar}}(x,y)-\digamma_x (y)| & \leq & |\psi_{\mathrm{Ar}} (x)| + |\psi_{\mathrm{Ar}} (y)| + | \mathbf{F}_x (y)| \\
& \leq & \frac{3696\pi^3 (g-1)^2 (2g-3)+400}{\xi_2 (g)} \cdot \varphi (C) \\
& & + \frac{2958032\pi^4 (g-1)^2 (2g-3)+30000\pi}{\xi_2 (g)} \cdot \varphi (C) \\
& \leq & \frac{(2958032\pi^4 + 3696\pi^3) (g-1)^2 (2g-\frac{5}{2})}{\xi_2 (g)} \cdot \varphi (C) .
\end{eqnarray*}
This completes the proof.
\end{proof}

In addition to the previous corollary, Lemma \ref{lemsupnormpotential} also yields an upper bound of Zhang's $\varphi$-invariant, which provides an inequality in the opposite direction to Theorem \ref{thmglobalestimatezhangphiinvariant}. Although this result will not be used later, it provides a complementary perspective on the inequality established earlier.

\begin{cor}
\label{coro_phiinv_upperbound}
Let $C$ be a compact Riemann surface of genus $g \geq 2$, equipped with the hyperbolic metric $\mu_{\mathrm{KE}}$ of constant curvature $-1$. Denote by $\mathrm{sys}(C)$ the systole of $(C,\mu_{\mathrm{KE}})$, and by $\varphi (C) $ the $\varphi$-invariant of $C$. Then 
$$\varphi (C) \leq 10^{6} g^{5} \cdot\max \left\{ 1 , \frac{1}{\mathrm{sys}(C)} \right\} ,$$
\end{cor}

\begin{proof}
For any $\alpha\in \Gamma (C,\omega_C)$, let $\mathfrak{u}_{\alpha}\in \mathscr{A}^0_{C,\mathbb{R}} (C)$ be the unique function satisfying $\int_C \mathfrak{u}_{\alpha} \mu_{\mathrm{hyp}} =0$ and
\begin{equation*}
    dd^c \mathfrak{u}_{\alpha} = i\alpha\wedge\bar{\alpha} - \left(\int_C i\alpha\wedge\bar{\alpha} \right) \mu_{\mathrm{hyp}} .
\end{equation*}

By the argument used in the proof of Lemma \ref{lmmbergmankernelfunctionupperbound}, we have
\begin{equation*}
    \frac{i}{2}\alpha\wedge\bar{\alpha} \leq (\sqrt{2}+1)^2 (g-1) \left( 1+\frac{2\arcsinh 1}{\mathrm{sys}(C)}\right)\Vert\alpha\Vert^2_{L^2;C} \cdot \mu_{\mathrm{hyp}} ,
\end{equation*}
where $\Vert\alpha\Vert^2_{L^2;C} = \frac{i}{2}\int_C\alpha\wedge\bar{\alpha}$. Hence
\begin{equation*}
    -2\Vert\alpha\Vert^2_{L^2;C}\leq dd^c \mathfrak{u}_{\alpha} \leq 2(\sqrt{2}+1)^2 (g-1) \left( 1+\frac{2\arcsinh 1}{\mathrm{sys}(C)}\right)\Vert\alpha\Vert^2_{L^2;C} .
\end{equation*}

Lemma \ref{lemsupnormpotential} then gives
\begin{eqnarray*}
 \int_C d \mathfrak{u}_{\alpha} \wedge d^c \mathfrak{u}_{\alpha} & \leq & 20(\sqrt{2}+1)^2 (g-1) \left( 1+\frac{2\arcsinh 1}{\mathrm{sys}(C)}\right)\Vert\alpha\Vert^4_{L^2;C} \\
 & & + \frac{6272\pi^3}{5} \cdot (g-1)^2 (2g-3) \cdot\frac{\max\{1,\mathrm{sys}(C)\}}{\mathrm{sys}(C)}\cdot\Vert\alpha\Vert^4_{L^2;C} \\
 & \leq & \frac{12544\pi^3}{5} \cdot (g-1)^3 \cdot\frac{\max\{1,\mathrm{sys}(C)\}}{\mathrm{sys}(C)}\cdot\Vert\alpha\Vert^4_{L^2;C} .
\end{eqnarray*}

Next, let $\{\alpha_j\}_{j=1}^g\subset \Gamma(C,\omega_C)$ be an $L^2$ orthonormal basis. For each pair $j,k$, let $\mathfrak{u}_{j,k} \in  \mathscr{A}_{C,\mathbb{C}}^0 (C) $ be the unique function satisfying $\int_C \mathfrak{u}_{j,k} \mu_{\mathrm{Ar}} =0$ and 
$$dd^c \mathfrak{u}_{j,k} = i\alpha_j \wedge\bar{\alpha}_k - \left( \int_{C} i\alpha_j \wedge\bar{\alpha}_k \right) \mu_{\mathrm{Ar}} .$$

By Proposition \ref{propzhangphiinvariantpotential}, the $\varphi$-invariant of $C$ satisfies
\begin{equation*}
    \varphi (C) = \frac{1}{2}\sum_{j,k=1}^g \int_{C} d\mathfrak{u}_{j,k} \wedge d^c \bar{\mathfrak{u}}_{j,k} = \frac{1}{2\pi } \sum_{j,k=1}^g \left\Vert d {\mathfrak{u}}_{j,k} \right\Vert^2_{L^2 ;C} .
\end{equation*}

We split the sum into two parts according to whether $j=k$ or $j\neq k$.

We first consider the case $j=k$. Since
$$dd^c \mathfrak{u}_{j,j} = i\alpha_j \wedge\bar{\alpha}_j - 2 \mu_{\mathrm{Ar}} = \frac{g-1}{g} dd^c \mathfrak{u}_{\alpha_j} - \frac{1}{g}\sum_{m\neq j} dd^c \mathfrak{u}_{\alpha_m} ,$$
the function $ \mathfrak{u}_{j,j} -\frac{g-1}{g}  \mathfrak{u}_{\alpha_j} + \frac{1}{g}\sum_{m\neq j}  \mathfrak{u}_{\alpha_m} $ is a constant. It follows that
\begin{eqnarray*}
\int_{C} d\mathfrak{u}_{j,j} \wedge d^c  \mathfrak{u}_{j,j} & = & \frac{1}{ \pi } \left\Vert d {\mathfrak{u}}_{j,j} \right\Vert^2_{L^2 ;C} = \frac{1}{ \pi } \left\Vert d\left(\frac{g-1}{g}  \mathfrak{u}_{\alpha_j} - \frac{1}{g}\sum_{m\neq j}  \mathfrak{u}_{\alpha_m}\right) \right\Vert^2_{L^2 ;C} \\
& \leq & \frac{1}{\pi}\left(\frac{g-1}{g}  \left\Vert d \mathfrak{u}_{\alpha_j} \right\Vert_{L^2 ;C} + \frac{1}{g}\sum_{m\neq j} \left\Vert d \mathfrak{u}_{\alpha_m} \right\Vert_{L^2 ;C}\right)^2 \\
& \leq & \frac{4(g-1)^2}{g^2} \cdot \frac{12544\pi^3}{5} \cdot (g-1)^3 \cdot\frac{\max\{1,\mathrm{sys}(C)\}}{\mathrm{sys}(C)} .
\end{eqnarray*}

Next, assume that $j\neq k$. A straightforward computation gives
\begin{eqnarray*}
dd^c \mathfrak{u}_{j,k} & = & \frac{i}{2} (\alpha_j+\alpha_k) \wedge (\bar{\alpha}_j+\bar{\alpha}_k) - \frac{i}{2} (\alpha_j-\alpha_k) \wedge (\bar{\alpha}_j-\bar{\alpha}_k) \\
& & + \frac{1}{2} (\alpha_j-i\alpha_k) \wedge (\bar{\alpha}_j+i\bar{\alpha}_k) - \frac{1}{2} (\alpha_j+i\alpha_k) \wedge (\bar{\alpha}_j-i\bar{\alpha}_k) \\
& = & \frac{1}{2} dd^c ( \mathfrak{u}_{\alpha_j+\alpha_k} -\mathfrak{u}_{\alpha_j-\alpha_k} +i\mathfrak{u}_{\alpha_j+i\alpha_k}-i\mathfrak{u}_{\alpha_j-i\alpha_k} ) .
\end{eqnarray*}

Applying the previous $L^2$-estimate to each term and using the triangle inequality, we obtain
\begin{eqnarray*}
\int_{C} d\mathfrak{u}_{j,k} \wedge d^c  \bar{\mathfrak{u}}_{j,k} & \leq & \frac{1}{4\pi} \left( \left\Vert d \mathfrak{u}_{\alpha_j+\alpha_k} \right\Vert_{L^2 ;C} + \left\Vert d \mathfrak{u}_{\alpha_j-\alpha_k} \right\Vert_{L^2 ;C}+\left\Vert d \mathfrak{u}_{\alpha_j+i\alpha_k} \right\Vert_{L^2 ;C} + \left\Vert d \mathfrak{u}_{\alpha_j-i\alpha_k} \right\Vert_{L^2 ;C}\right)^2 \\
& \leq & 16 \cdot \frac{12544\pi^3}{5} \cdot (g-1)^3 \cdot\frac{\max\{1,\mathrm{sys}(C)\}}{\mathrm{sys}(C)} .
\end{eqnarray*}

Combining the two cases,
\begin{eqnarray*}
\varphi (C)& \leq &\frac{1}{2}\sum_{j=1}^g \int_{C} d\mathfrak{u}_{j,j} \wedge d^c  \mathfrak{u}_{j,j} +\frac{1}{2}\sum_{j\neq k}\int_{C} d\mathfrak{u}_{j,k} \wedge d^c  \bar{\mathfrak{u}}_{j,k} \\
& \leq & \frac{g}{2}\cdot \frac{4(g-1)^2}{g^2} \cdot \frac{12544\pi^3}{5} \cdot (g-1)^3 \cdot\frac{\max\{1,\mathrm{sys}(C)\}}{\mathrm{sys}(C)} \\
& & +  \frac{g(g-1)}{2} \cdot 16 \cdot \frac{12544\pi^3}{5} \cdot (g-1)^3 \cdot\frac{\max\{1,\mathrm{sys}(C)\}}{\mathrm{sys}(C)} \\
& \leq & 10^6 \cdot g^5 \cdot\frac{\max\{1,\mathrm{sys}(C)\}}{\mathrm{sys}(C)},
\end{eqnarray*}
as claimed.
\end{proof}

We now establish an estimate for the hermitian line bundle $(\mathcal{O}(\Delta),\Vert\cdot\Vert_{\Delta})$, and derive a quantitative comparison between the Arakelov hermitian metric $\Vert\cdot\Vert_{\mathrm{Ar}}$ and hyperbolic hermitian metric $\Vert\cdot \Vert_{\mathrm{hyp}}$ on $\omega_C$.

\begin{thm}[Theorem \ref{propestimateArakelovofdiagonal_compareTwoMetrics part I}]
\label{propestimateArakelovofdiagonal_compareTwoMetrics}
Let $\Vert\cdot\Vert_{\Delta}$ denote the Arakelov hermitian metric on $\mathcal{O}(\Delta)\to C\times C$, and let $\varphi (C)$ be the $\varphi$-invariant of $C$. For any $(x_1,x_2)\in C\times C$, let $\rho=(\rho_1,\rho_2):\mathbb{D}\times \mathbb{D}\to C\times C$ be a universal covering map with $\rho_1(0)=x_1$ and $\rho_2(0)=x_2$. Then there exists a local holomorphic section $s$ of $(\rho^*\mathcal{O}(\Delta),\rho^*\Vert \cdot\Vert_{\Delta})$ on $\mathbb{D}_{\frac{1}{20}}(0)\times\mathbb{D}_{\frac{1}{20}}(0)$, such that
$$ \big| \log \Vert s\Vert \big| \leq \frac{(5916064\pi^4 + 7392\pi^3) (g-1)^3}{\xi_2 (g)} \cdot \varphi (C) < 2\cdot 10^{14}g^{\frac{14}{3}} \cdot \varphi (C) ,$$
where $\xi_2(g)$ is the constant from 
Theorem \ref{thmglobalestimatezhangphiinvariant}.

In particular, the Arakelov hermitian metric $\Vert\cdot\Vert_{\mathrm{Ar}}$ and the hyperbolic hermitian metric $\Vert\cdot \Vert_{\mathrm{hyp}}$ on $\omega_C\cong\mathcal{O} (-\Delta)|_{\Delta}$ satisfy
$$ \left| \log \frac{\Vert\cdot\Vert_{\mathrm{Ar}}}{\Vert\cdot\Vert_{\mathrm{hyp}}} \right| \leq \frac{(5916064\pi^4 + 7392\pi^3) (g-1)^3}{\xi_2 (g)} \cdot \varphi (C) < 2\cdot 10^{14}g^{\frac{14}{3}} \cdot \varphi (C) .$$
\end{thm}

\begin{proof}
The proof of the first statement divides naturally into two cases, depending on the distance between $x_1$ and $x_2$ with respect to the unique K\"ahler metric $\mu_{\mathrm{KE}}$ with constant curvature $-1$.

If $\dist_{\mu_{\mathrm{KE}}} (x_1,x_2)<\frac{1}{2}\arcsinh\frac{1}{2}$, there exists a point $x'\in B_{\frac{1}{4}\arcsinh\frac{1}{2}} (x_1)\cap B_{\frac{1}{4}\arcsinh\frac{1}{2}} (x_2)$. Let $U_{x'}$, $z_{x'}$ and $\digamma $ be as in Lemma \ref{lemconstructionlocalresidue}. Then $U_{x'}$ contains both $B_{\frac{1}{4}\arcsinh\frac{1}{2}} (x_1)$ and $ B_{\frac{1}{4}\arcsinh\frac{1}{2}} (x_2)$, so for any $(y_1,y_2)\in B_{\frac{1}{4}\arcsinh\frac{1}{2}} (x_1)\times B_{\frac{1}{4}\arcsinh\frac{1}{2}} (x_2)$, we have
\begin{equation*}
|G_{\mathrm{Ar}} (y_1 ,y_2)-\log |z_{x'}(y_1)-z_{x'}(y_2)|| \leq |G_{\mathrm{Ar}} (y_1 ,y_2)-\digamma_{y_1} (y_2)| + |\digamma_{y_1} (y_2)-\log |z_{x'}(y_1)-z_{x'}(y_2)|| .
\end{equation*}
Applying Theorem \ref{thmglobalestimatezhangphiinvariant}, Lemma \ref{lemconstructionlocalresidue} and Corollary \ref{corollarydiatanceArakelovGreen_Log}, we deduce
\begin{eqnarray*}
|G_{\mathrm{Ar}} (y_1 ,y_2)-\log |z_{x'}(y_1)-z_{x'}(y_2)|| & \leq & \frac{(2958032\pi^4 + 3696\pi^3) (g-1)^2 (2g-\frac{5}{2})}{\xi_2 (g)} \cdot \varphi (C) \\
& & + \frac{8\pi^2\cdot \max\{1,\mathrm{sys}(C)\}}{\mathrm{sys}(C)} \\
& \leq & \frac{(2958032\pi^4 + 3696\pi^3) (g-1)^2 (2g-\frac{5}{2}) }{\xi_2 (g)} \cdot \varphi (C) \\
& & + 8\pi^2 \cdot \max\left\{ \frac{1}{\xi_1 (g)} ,\frac{1}{\xi_2(g)} \right\} \cdot \varphi (C) \\
& < & \frac{2\cdot(2958032\pi^4 + 3696\pi^3) (g-1)^3  }{\xi_2 (g)} \cdot \varphi (C) ,
\end{eqnarray*}
where $\xi_1 (g)$ is the constant from 
Theorem \ref{thmglobalestimatezhangphiinvariant}, and we used that $1\leq \frac{\varphi (C)}{\xi_1 (g)}$, $\frac{1}{\mathrm{sys}(C)}\leq \frac{\varphi (C)}{\xi_2 (g)}$, and $\xi_1 \geq 20g\cdot \xi_2$. Let $\mathbf{1}_{\Delta}$ denote the canonical section of $\mathcal{O}(\Delta)$ and set $s'(y_1 ,y_2)= \frac{1}{z_{x'}(y_1)-z_{x'}(y_2)}\mathbf{1}_{\Delta}$ on $ B_{\frac{1}{4}\arcsinh\frac{1}{2}} (x_1)\times B_{\frac{1}{4}\arcsinh\frac{1}{2}} (x_2)$. Its Arakelov norm satisfies 
$$ \left| \log \Vert s'\Vert_{\Delta} \right| \leq \frac{(5916064\pi^4 + 7392\pi^3) (g-1)^3}{\xi_2 (g)} \cdot \varphi (C) < 2\cdot 10^{14}g^{\frac{14}{3}} \cdot \varphi (C) .$$
Since $\tanh \frac{\arcsinh\frac{1}{2}}{8} \approx 0.06\cdots >\frac{1}{20} $, $\rho^{-1}(B_{\frac{1}{4}\arcsinh\frac{1}{2}} (x_1)\times B_{\frac{1}{4}\arcsinh\frac{1}{2}} (x_2))$ contains $\mathbb{D}_{\frac{1}{20}}(0)\times\mathbb{D}_{\frac{1}{20}}(0)$. We define $s=\rho^* s'$, which then serves as the desired local holomorphic section on $\mathbb{D}_{\frac{1}{20}}(0)\times\mathbb{D}_{\frac{1}{20}}(0)$.

On the other hand, if $\dist_{\mu_{\mathrm{KE}}} (x_1,x_2)\geq \frac{1}{2}\arcsinh\frac{1}{2}$, then for any $(y_1,y_2)\in B_{\frac{2}{9}\arcsinh\frac{1}{2}} (x_1)\times B_{\frac{2}{9}\arcsinh\frac{1}{2}} (x_2)$, we have
$$\dist_{\mu_{\mathrm{KE}}} (y_1 ,y_2) \geq \dist_{\mu_{\mathrm{KE}}} (x_1 ,x_2)-\dist_{\mu_{\mathrm{KE}}} (x_1 ,y_1)-\dist_{\mu_{\mathrm{KE}}} (x_2 ,y_2) >\frac{\arcsinh\frac{1}{2}}{18} >\frac{1}{100} .$$
Lemma \ref{lemconstructionlocalresidue} then gives $|\digamma (y_1,y_2)|\leq \frac{6\pi^2\cdot \max\{1,\mathrm{sys}(C)\}}{\mathrm{sys}(C)} $, and thus
\begin{eqnarray*}
|G_{\mathrm{Ar}} (y_1 ,y_2) | & \leq & |G_{\mathrm{Ar}} (y_1 ,y_2)-\digamma (y_1 ,y_2)| + |\digamma (y_1 ,y_2)| \\& \leq & \frac{(2958032\pi^4 + 3696\pi^3) (g-1)^2 (2g-\frac{5}{2})+6\pi^2}{\xi_2 (g)} \cdot \varphi (C) \\
& < & \frac{2\cdot(2958032\pi^4 + 3696\pi^3) (g-1)^3  }{\xi_2 (g)} \cdot \varphi (C) .
\end{eqnarray*}
Then on $B_{\frac{2}{9}\arcsinh\frac{1}{2}} (x_1)\times B_{\frac{2}{9}\arcsinh\frac{1}{2}} (x_2)$,
$$ \left| \log \Vert \mathbf{1}_{\Delta}\Vert_{\Delta} \right| \leq \frac{(5916064\pi^4 + 7392\pi^3) (g-1)^3}{\xi_2 (g)} \cdot \varphi (C) < 2\cdot 10^{14}g^{\frac{14}{3}} \cdot \varphi (C) .$$
Since $\tanh \frac{\arcsinh\frac{1}{2}}{9} \approx 0.053\cdots >\frac{1}{20} $, $\rho^{-1}(B_{\frac{2}{9}\arcsinh\frac{1}{2}} (x_1)\times B_{\frac{2}{9}\arcsinh\frac{1}{2}} (x_2))$ contains $\mathbb{D}_{\frac{1}{20}}(0)\times\mathbb{D}_{\frac{1}{20}}(0)$. We take the pullback of the canonical section $\mathbf{1}_\Delta$ to $\mathbb{D}_{\frac{1}{20}}(0)\times\mathbb{D}_{\frac{1}{20}}(0)$ as $s$, which serves as the desired local holomorphic section. 

For the second statement, for any $x$, let $z_x$ and $ \digamma_x : C\setminus \{x\} \to \mathbb{R}$ be as in Lemma \ref{lemconstructionlocalresidue}. By the definition of the Arakelov hermitian metric,
\begin{eqnarray*}
\left| \log \frac{\Vert\cdot\Vert_{\mathrm{Ar}}}{\Vert\cdot\Vert_{\mathrm{hyp}}} (x) \right| & = & \left| \log \frac{\Vert dz_x\Vert_{\mathrm{Ar}}}{\Vert dz_x\Vert_{\mathrm{hyp}}} (x) \right| \\
& \leq & \left| \lim\limits_{y\to x} \left( \digamma_x (y) +\log |z_x (x)-z_x(y) | \right) \right| + \left| \log \Vert dz_x\Vert_{\mathrm{hyp}} (x) \right| \\
& & + \left| \lim\limits_{y\to x} \left( G_{\mathrm{Ar}} (x,y) - \digamma_x (y) \right) \right|  .
\end{eqnarray*}

Using the local estimates for $\digamma_x (y) +\log |z_x (x)-z_x(y) |$ and $\Vert dz_x\Vert_{\mathrm{hyp}}$ in terms of the systole $\mathrm{sys}(C)$ and the global bound from Theorem \ref{thmglobalestimatezhangphiinvariant}, we obtain
\begin{eqnarray*}
\left| \log \frac{\Vert\cdot\Vert_{\mathrm{Ar}}}{\Vert\cdot\Vert_{\mathrm{hyp}}} (x) \right| & \leq & \frac{12\pi^2 \cdot \max\left\{ 1 ,\mathrm{sys}(C)\right\}}{\mathrm{sys}(C)} + \frac{(2958032\pi^4 + 3696\pi^3) (g-1)^2 (2g-\frac{5}{2})}{\xi_2 (g)} \cdot \varphi (C) ,
\end{eqnarray*}

Using again that $1\leq \frac{\varphi (C)}{\xi_1 (g)}$, $\frac{1}{\mathrm{sys}(C)}\leq \frac{\varphi (C)}{\xi_2 (g)}$ and $\xi_1 \geq 20g\cdot \xi_2$, we deduce
\begin{eqnarray*}
\left| \log \frac{\Vert\cdot\Vert_{\mathrm{Ar}}}{\Vert\cdot\Vert_{\mathrm{hyp}}} (x) \right|
& \leq & \frac{(2958032\pi^4 + 3696\pi^3) (g-1)^2 (2g-\frac{5}{2}) +12\pi^2}{\xi_2 (g)} \cdot \varphi (C) \\
& < & \frac{2\cdot(2958032\pi^4 + 3696\pi^3) (g-1)^3  }{\xi_2 (g)} \cdot \varphi (C) < 2\cdot 10^{14}g^{\frac{14}{3}} \cdot \varphi (C) ,
\end{eqnarray*}
which completes the proof.
\end{proof}


\

{\footnotesize
\noindent Jiawei Yu

\noindent Address: \emph{School of Mathematical Sciences, Peking University, Haidian District, Beijing 100871, China}

\noindent  Email: \emph{yujiawei@pku.edu.cn}

\

\noindent Xinyi Yuan

\noindent Address: \emph{BICMR, Peking University, Haidian District, Beijing 100871, China}

\noindent Email: \emph{yxy@bicmr.pku.edu.cn}

\

\noindent Shengxuan Zhou

\noindent Address: \emph{Institut de Math{\'e}matiques de Toulouse, 118 route de Narbonne, F-31062 Toulouse Cedex 9, France}

\noindent  Email: \emph{zhoushx98@outlook.com}

\
}

\end{document}